\documentclass[5p]{elsarticle}

\usepackage{amssymb}
\usepackage{amsmath}

\usepackage{algorithm,algorithmicx}
\usepackage{algpseudocode}
\usepackage{siunitx}
\usepackage{multicol}
\usepackage{diagbox}
\usepackage{graphicx}
\usepackage{tikz-cd}
\usepackage{url}
\usepackage{bm}

\journal{Nuclear Physics B}

\begin{document}

\begin{frontmatter}



\title{Reconstructing Flash-Evoked Potentials in the Human Brain and Tracking Cortical Activity Using Spatiotemporal Filtering: Numerical Simulation Study}


\author[tam]{Qiong Tao\corref{cor}}
\author[tam]{Joonas Lahtinen}
\author[tam]{Sampsa Pursiainen}

\affiliation[tam]{organization={Tampere University},
            addressline={Korkeakoulunkatu 1}, 
            city={Tampere},
            postcode={33720}, 
            state={Pirkanmaa},
            country={Finland}\\ \today }

\cortext[cor]{Corresponding author}

\begin{abstract}

Electroencephalographic source reconstruction of visual evoked activity is challenging because distributed, temporally evolving, and overlapping cortical sources generate scalp potentials. We present a controlled simulation framework for flash visual evoked potentials (FVEPs) based on a coupled Jansen--Rit neural mass network and a realistic finite-element EEG forward model. The framework generates a sequential activation pattern across early visual, parietal, and frontal regions and enables direct comparison of reconstructed activity with known ground truth. Standardized Kalman filtering (SKF), change-rate SKF (CrSKF), conventional Kalman filtering (KF), and linearly constrained minimum variance (LCMV) beamforming are compared using parcel-wise time courses, temporal lags, and assessments of spatial focality and Earth Mover's Distance. In the inverse comparison, SKF and CrSKF preserve the principal temporal progression more consistently than KF and LCMV, while CrSKF produces smoother temporal estimates. The results also demonstrate parcel-wise leakage and reduced separation of closely spaced components, emphasizing the need to interpret temporal ordering together with spatial reconstruction quality. The proposed framework provides a reproducible test environment for evaluating dynamic EEG inverse methods under controlled noise conditions.
\end{abstract}




\begin{highlights}
\item A coupled Jansen--Rit network generates a controlled synthetic FVEP with known spatiotemporal source progression.
\item SKF and CrSKF preserve the principal temporal progression more consistently than conventional KF and LCMV in the present simulation.
\item Parcel-wise leakage and source overlap remain important limitations even when the underlying source sequence is known.
\end{highlights}


\end{frontmatter}


\section{Introduction}
Understanding how distributed neural activity gives rise to measured brain signals remains a central problem in computational neuroscience. Electroencephalography (EEG) provides millisecond-scale temporal resolution, but scalp recordings contain superposed contributions from multiple neuronal populations that have been spatially smoothed by volume conduction. Recovering the underlying cortical activity is therefore a highly ill-posed inverse problem, and different reconstruction methods can yield substantially different source estimates because they impose different spatial and temporal assumptions.

In particular, the Jansen–Rit neural mass model \cite{JansenRit1995} offers a biophysically interpretable representation of cortical population dynamics that can reproduce realistic EEG rhythms and evoked responses. Jansen and Rit's neural mass model describes mesoscopic activity of the neural population. In this model, the complex variety of neurons within a cortical column is simplified to three main groups: pyramidal neurons, excitatory interneurons, and inhibitory interneurons. The pyramidal neurons stimulate both types of interneurons, and these interneurons subsequently send feedback to the pyramidal cells. Because of this, the Jansen-Rit neural mass model strikes a balance between realism and computational complexity \cite{Escuain2018JRUKF}. Coupling neural mass models to an EEG forward model therefore enables a controlled test environment for inverse methods while retaining characteristic features of the evoked activity, including transient responses, polarity changes, multiple components, and propagation between cortical regions.

This study was inspired by Qiong Tao's master's thesis \cite{QiongTao2026thesis}, in which she conducted preliminary research on this topic. 
In this study, we develop a synthetic flash visual evoked potential (FVEP) generation framework based on a coupled Jansen--Rit network. FVEPs provide a particularly demanding spatiotemporal test because their scalp topography changes markedly across successive components. The model is designed to reproduce a sequential activation pattern across early visual, parietal, and frontal regions and to generate a sensor-level waveform with the principal polarity pattern and relative component structure of representative FVEP recordings. Jansen--Rit dynamics have previously been shown to reproduce FVEP-like responses \cite{Jansen1993FVEP}.

The simulated FVEP is used to compare dynamic and spatially filtered EEG inverse methods. Standardized Kal\-man filtering (SKF) \cite{Lahtinen2024SKF,Dilshanie2026SKF} and its change-rate extension (CrSKF) \cite{Lahtinen2026CRSKF,piispa2025tracking} have previously been evaluated with simulated and measured somatosensory evoked potentials. Here, both methods are tested with the richer spatial variability of FVEPs and compared with conventional Kalman filtering \cite{GalkaAndreas2004KalmanEEG} and a linearly constrained minimum variance (LCMV) beamformer \cite{VeenB1992Lois}. LCMV is widely used in neuroelectromagnetic source analysis and studies of neural interactions and dynamical connectivity \cite{HadjipapasAvgis2005LCMV,HuiHuaBrian2010LCMV,CoelliStefania2023LCMV}.
The computational burden of Kalman filters mentioned in previous studies of spatiotemporal source estimation \cite{GalkaAndreas2004KalmanEEG, OuHamalainen2009} does not reflect the current situation. Advanced techniques of multithreading and efficient memory management, combined with low-level coding, enable efficient matrix multiplication and inversion. 

The methodological contributions are threefold: (i) a coupled neural-mass model for generating a known multiregion FVEP source sequence, (ii) projection of that activity through a realistic finite-element EEG forward model, and (iii) a common evaluation of four inverse methods using parcel-wise time courses, spatial reconstructions, temporal lags, and Earth Mover's Distance (EMD). The aim is both to assess whether the synthetic data reproduce salient FVEP characteristics and to use the resulting known ground truth as a benchmark for spatiotemporal source reconstruction.

\section{Background}
\subsection{Identifying a dynamical system from VEP}
Visual evoked potentials (VEPs) provide a non-invasive way to characterize the temporal dynamics of neural responses to visual stimulation. Conventional VEPs are typically obtained by presenting discrete visual stimuli and averaging the resulting electroencephalographic (EEG) responses. Although this approach has been widely used to investigate visual processing, the requirement for repeated stimulus presentations can make measurements relatively time-consuming. Lalor et al. introduced visual-evoked spread spectrum analysis (VESPA) as an alternative approach in which the luminance of a visual stimulus is continuously varied according to a stochastic process under high SNR conditions [1]. Because the stimulus has a rich statistical structure, the relationship between the stimulus sequence and the simultaneously recorded EEG can be used to estimate the impulse response of the visual system. The resulting VESPA response was shown to resemble the conventional VEP and could be estimated rapidly, providing an alternative framework for measuring visually evoked activity \cite{LalorEdmund2006VESPA}. Subsequent comparison of VESPA and conventional VEP responses provides evidence that at least some components of the two measures may reflect common cortical generators, while also demonstrating that similarities in waveform morphology should not by themselves be taken as evidence of equivalent physiological origins \cite{MurphyJeremy2012IsolatingEarlyCorticalGenerators}. This distinction is important when interpreting system-identification-derived visual responses and relating them to established measures of visual-evoked activity.

The principle underlying VESPA is related to a broader class of system-identification approaches in which the response of a biological system is estimated from its response to a temporally structured input.  These approaches establish a framework in which continuously or pseudo-random\-ly varying visual stimulation can be combined with system-identification methods to estimate temporally resolved neural responses. The multifocal m-sequence work shows how structured stimulation and kernel analysis can characterize dynamic visual responses \cite{SutterErich2000MultifocalBinaryKernels,SutterErich2001MultifocalMsequence}. 

By modeling the evoked response as the convolution of an underlying response function with the stimulus sequence, deconvolution can estimate the temporal response associated with individual stimuli despite overlap between consecutive VEPs. This approach provides a means of recovering temporally resolved evoked responses without requiring sufficiently long inter-stimulus intervals for individual responses to return to baseline. It therefore complements system-identification methods that characterize the relationship between visual stimulation and EEG responses \cite{BohorquezJ2013ModelingOverlappingVEP}.

\subsection{VEP/VEF source estimation}
The study \cite{Ponasso2025} aimed to reconstruct the cortical sources of visual evoked potentials elicited by intermittent photic stimulation. In particular, the authors investigated whet\-her a high-resolution, anatomically realistic head model could localize the neural generators of the P100/P1 component of the visual evoked response. EEG data from 12 healthy volunteers were combined with individual MRI-derived he\-ad models. Principal component analysis was applied to the reciprocal basis functions to remove the dominant component associated with bias toward the reference electrode. Reconstructed activity was predominantly localized to the expected visual cortical regions, with negligible activation in lateral cortical regions. However, anatomical agreement alone does not prove that the estimated source distribution is correct. 

The study by Hagler et al. \cite{Hagler2009}, aimed to reconstruct the time courses of neural activity in multiple visual cortical areas, particularly V1, V2, and V3. The central challenge was distinguishing the contributions of nearby, simultaneously active sources from MEG and EEG recordings. The authors therefore combined electrophysiological recordings with individual structural MRI and fMRI-derived retinotopic maps. Dipolar activity was assumed, and dipole locations and orientations were specified using MRI-derived cortical surfaces and fMRI retinotopic maps, rather than estimating them freely from the electrophysiological data; e.g., the dipole locations were constrained inside V1--V3 regions. The study examined the waveforms of each stimulus location. The results show that distinguishing regional waveforms is difficult because they tend to be highly correlated across nearby anatomical regions. Although, observing the time delays is plausible. The authors projected estimated sources back into sensor space and calculated the normalized residual error between predicted and measured data. This quantified how well each inverse solution explained the observations. Importantly, the authors found that stronger constraints could produce more plausible and stable source waveforms even when the residual error increased. Thus, a smaller residual error alone was not sufficient to establish that a reconstruction was physiologically meaningful. 

Using a simulation of visual evoked potentials could provide the missing information about the reconstruction method's specific aspects that can limit source and waveform estimation accuracy.

\section{Methods}

\subsection{Jansen--Rit-based observation model}

The Jansen--Rit neural mass model \cite{JansenRit1995} describes cortical activity at the mesoscopic scale using three interacting population types: pyramidal neurons, excitatory interneurons, and inhibitory interneurons. The pyramidal population excites both interneuron groups, which in turn project back to the pyramidal cells. 
The mean rate of incoming action potentials is mapped to an average postsynaptic membrane potential through second -order synaptic operators. For excitatory input \noindent$f^{(e)}(t)$

\begin{equation}
    \ddot{u}^{(e)}(t)+2a\dot{u}^{(e)}(t)+a^2u^{(e)}(t)=aAf^{(e)}(t),
\end{equation}
and for inhibitory input,
\begin{equation}
    \ddot{u}^{(i)}(t)+2b\dot{u}^{(i)}(t)+b^2u^{(i)}(t)=bBf^{(i)}(t),
\end{equation}
where $a$ and $b$ are the excitatory and inhibitory inverse time constants and $A$ and $B$ are the corresponding postsy\-nap\-tic-potential gains.

Following \cite{Babajani2010MultiJR,Escuain2018JRUKF}, the coupled dynamics of column $k$ are written as
\begin{equation}
\begin{aligned}
&\ddot{u}^{(0)}_k(t)+2a\dot{u}^{(0)}_k(t)+a^2u^{(0)}_k(t)\\
&\qquad =Aa\,\mathrm{Sigm}\!\left(
u_k^{(e)}(t)-u_k^{(i)}(t)\right),\\[0.35em]
&\ddot{u}^{(e)}_k(t)+2a\dot{u}^{(e)}_k(t)+a^2u^{(e)}_k(t)\\
&\qquad =Aa\Bigg[
\sum_{j=1}^{N}K_{kj}\,
\mathrm{Sigm}\!\left(
u_j^{(e)}(t-\tau_{kj})
-u_j^{(i)}(t-\tau_{kj})
\right)\\
&\hspace{5.7em}
+C_2\,\mathrm{Sigm}\!\left(C_1u_k^{(0)}(t)\right)
+p_k(t)
\Bigg],\\[0.35em]
&\ddot{u}^{(i)}_k(t)+2b\dot{u}^{(i)}_k(t)+b^2u^{(i)}_k(t)\\
&\qquad =bBC_4\,\mathrm{Sigm}\!\left(C_3u_k^{(0)}(t)\right).
\end{aligned}
\end{equation}
Here, $K$ is the directed inter-column adjacency matrix, $\tau_{kj}$ is the propagation delay from column $j$ to column $k$, $C_1,\ldots,C_4$ are intrinsic connectivity constants, and $p_k(t)$ is the external afferent input. In the present visual-network implementation, nonzero entries of $K$ correspond to the directed connections shown in the demonstration setup below. A common inter-column coupling weight is used, while the relative amplitudes of the regional responses are controlled separately by the output gains specified below.

The firing-rate nonlinearity is
\begin{equation}
    \mathrm{Sigm}(\nu)=\frac{2e_0}{1+e^{r(\nu_0-\nu)}},
\end{equation}
where $e_0$ is the half-maximum firing-rate parameter, $r$ is the sigmoid slope, and $\nu_0$ is the membrane potential at half activation. The simulations use the classical Jansen--Rit parameterization
$A=3.25$ mV, $B=22.0$ mV, $a=100$ s$^{-1}$, $b=50$ s$^{-1}$,
$e_0=5$ Hz, $\nu_0=6$ mV, and $r=0.56$ mV$^{-1}$. The intrinsic
connectivity coefficients are $C_1=135$, $C_2=108$, $C_3=33.75$, and
$C_4=33.75$, corresponding to $C_2=0.8C_1$ and
$C_3=C_4=0.25C_1$.

The model output is the pyramidal membrane-potential difference
$V_k=u_k^{(e)}-u_k^{(i)}$. This quantity is interpreted as a source-level intracranial signal representing aggregate postsynaptic activity rather than as a scalp EEG voltage \cite{Sanchez-TodoRoser2023JRintepretation}. The extracranial signal is obtained by forward projection. Collecting the source outputs in ${\bf V}_t$, the sensor data are modeled as
\begin{equation}
    {\bf y}_t=\alpha L{\bf V}_t,
\end{equation}
where $L$ is the EEG lead-field matrix and $\alpha$ converts the relative source amplitudes to dipole-strength units. The scaling is calibrated so that a dipole of approximately 10 nAm produces a scalp potential of approximately 1~$\mu$V, providing a common amplitude scale for forward simulation and inversion.

\subsection{Source reconstruction methods}
Galka {\em et al.} introduced Kalman filtering (KF) to the distributed dynamical EEG inverse problem \cite{GalkaAndreas2004KalmanEEG}. The state vector contains the three Cartesian dipole components at every source location, so that ${\bf x}_t\in\mathbb{R}^{3n}$ for $n$ source locations. A random-walk model, $A_t=I$, is used for KF and SKF. The initial state is zero-mean with isotropic covariance $P_{0|0}=\eta_0 I$. Additive sensor noise is modeled as independent Gaussian noise with covariance $R_t=\sigma_\varepsilon^2 I$, where $\sigma_\varepsilon$ is determined from the prescribed measurement SNR. The process covariance is likewise isotropic, $Q_t=qI$, with $q$ scaled relative to $\eta_0$ so that the accumulated process variance remains independent of the number of time samples. The same SNR-based covariance calibration is used for KF and SKF, ensuring that their differences arise from standardization rather than different noise assumptions.

SKF is a spatiotemporal extension of the sLORETA principle \cite{PascualMarqui2002}. The filtered distributed estimate ${\bf x}_t$ is transformed to a standardized field ${\bf z}_t$ after each Kalman update. We use the standard-deviation normalization $p=1/2$ in the weighting operator:

\begin{algorithm}[H]
\caption{Standardized Kalman Filter algorithm}
\begin{algorithmic}[1]
\State \textbf{Input:} Initial state estimate ${\bf x}_{0\mid 0}$, source covariance $P_{0\mid 0}$, state-transition matrix $A_t$, lead field $L$, evolution covariance $Q_t$, observation noise covariance $R_t$, time steps $t = 1, \ldots, T$
\For{$t = 1$ to $T$}
    \State \textbf{Kalman filtering algorithm}
    \State {\em Prediction Step:}
    \State ${\bf x}_{t\mid t-1} = A_t {\bf x}_{t-1\mid t-1}$
    \State $P_{t\mid t-1} = A_t P_{t-1\mid t-1} A_t^\mathrm{T} + Q_t$
    
    \State $S_t = L P_{t\mid t-1} L^\mathrm{T} + R_t$
    \State $K_t = P_{t\mid t-1} L^\mathrm{T} S_t^{-1}$
    
    \State {\em Update Step:}
    \State ${\bf x}_{t\mid t} = {\bf x}_{t\mid t-1} + K_t({\bf y}_t - L {\bf x}_{t\mid t-1})$
    \State $P_{t\mid t} = P_{t\mid t-1} - K_t S_t K_t^\mathrm{T}$
    
    \State \textbf{Standardization: }${\bf x}_{t\mid t}\mapsto {\bf z}_{t\mid t}$
    \State $W_t = \mathrm{Diag}\left({P_{\bf x}}_{t\mid t-1}^{-1/2} K_t S_t K_t^\mathrm{T} {P_{\bf x}}_{t\mid t-1}^{-1/2} \right)^{-p} {P_{\bf x}}_{t\mid t-1}^{-1/2}$
    \State ${\bf z}_{t\mid t} = W_t {\bf x}_{t\mid t}$

   \EndFor
\end{algorithmic}\label{Algo:KalmanFilterStandardized}
\end{algorithm}

Change-rate Standardized Kalman filter (CrSKF) \cite{Lahtinen2026CRSKF} takes a step away from the random walk evolution model by introducing a hidden change-rate state variable ${\bf v}_t$ representing the change rate of ${\bf x}_t$, yielding the evolution model
\begin{align}
    {\bf x}_{t+1}&={\bf x}_t+{\bf v}_t\Delta t +{\bf q}_{t+1},\\
    {\bf v}_{t+1}&={\bf v}_t+{\bf c}_{t+1},
\end{align}
where ${\bf q}_{t+1}$ and ${\bf c}_{t+1}$ are independent zero-mean Gaussian evolution noises. Their diagonal covariances are scaled from the same source and measurement-noise calibration used for KF and SKF, with the change-rate covariance chosen so that the augmented state has comparable total process variance. Furthermore, the second-order backward difference equation is utilized to model the dynamical change of the observations as a piecewise quadratic polynomial behavior 
\begin{equation}
    \frac{3}{2}{\bf y}_t-2{\bf y}_{t-1}+\frac{1}{2}{\bf y}_{t-2}=\Delta t L{\bf v}_{t},
\end{equation}
which smooths out the drastic changes assumed to arise from measurement noise. Algorithmically, instead of tracking ${\bf x}_t$ by Kalman filter, we track $\left[{\bf x}_t^T\; {\bf v}_t^T\right]^T$ and instead of using ${\bf y}_t$ as the observations, we have extended observation vector $\left[{\bf y}_t^T\; {\bf f}_t^T\right]^T$, where ${\bf f}_t =\frac{3}{2}{\bf y}_t-2{\bf y}_{t-1}+\frac{1}{2}{\bf y}_{t-2}$.

In this paper, we use a vectorized version of the linearly constrained minimum variance (LCMV) beamformer \cite{VeenB1992Lois} provided by Zeffiro Interface \cite{he2018zeffiro,Rezaei2021SEPRAMUS}. The sensor covariance is estimated over the full simulated epoch and regularized by diagonal loading, $\Sigma_\lambda=\Sigma_y+\lambda I$, with $\lambda$ set to the sensor-noise variance implied by the imposed SNR. The filter is then derived by solving the following minimization problem
\begin{equation}
    \min_{W_i\in\mathbb{R}^{m\times 3}} \mathrm{tr}\left(W_i^T\Sigma W_i\right) \quad\mathrm{subject\: to}\quad W_i^TL_i=I,
\end{equation}
for $i=1,\cdots,n$, where $\Sigma$ is the error covariance matrix. By using the Lagrangian method of multipliers, we have the explicit form for the filtering matrix:
\begin{equation}
    W_i=\Sigma ^{-1}L_i\left[L_i^T\Sigma ^{-1}L_i\right]^{-1}.
\end{equation}
The solution is obtained by filtering the observations $\hat{x}_{t;i}=W_i^T{\bf y}_t$ for each source location $i=1,\cdots, n$.

\subsection{Demonstration setup}
Neural activity responsible for flash visually evoked potentials (FVEP) is primitive and therefore suitable for detecting abnormalities in optic neuropathy \cite{Mohammed2021UseOfFlshVEP}. It is also a useful testing mechanism for children and subjects with short attention spans or low cooperation \cite{SinghRamji2016VEPreview}. The synthetic FVEP is represented by a simplified visual-processing network with sources in V1, ventral V2 and V3, dorsal V2 and V3, lateral intraparietal cortex (LIP), and frontal eye field (FEF). The branching structure reflects the conventional distinction between ventral and dorsal visual processing streams; here, it serves as a modeling scaffold rather than a claim that an experimental FVEP follows a unique deterministic pathway. The directed sequence is
\begin{center}
\begin{tikzcd}%
        \textnormal{V1} \arrow[r] \arrow[rd] & \textnormal{V2v} \arrow[r] & \textnormal{V3v} \arrow[r]  &  \textnormal{LIP} \arrow[r]  & \textnormal{FEF} \\  
         & \textnormal{V2d} \arrow[r]  & \textnormal{V3d} 
    \end{tikzcd}%
\end{center}

\begin{figure}[h!]
    \centering
    \includegraphics[trim={5cm 2.5cm 9cm 5.5cm},clip,width=0.4\linewidth]{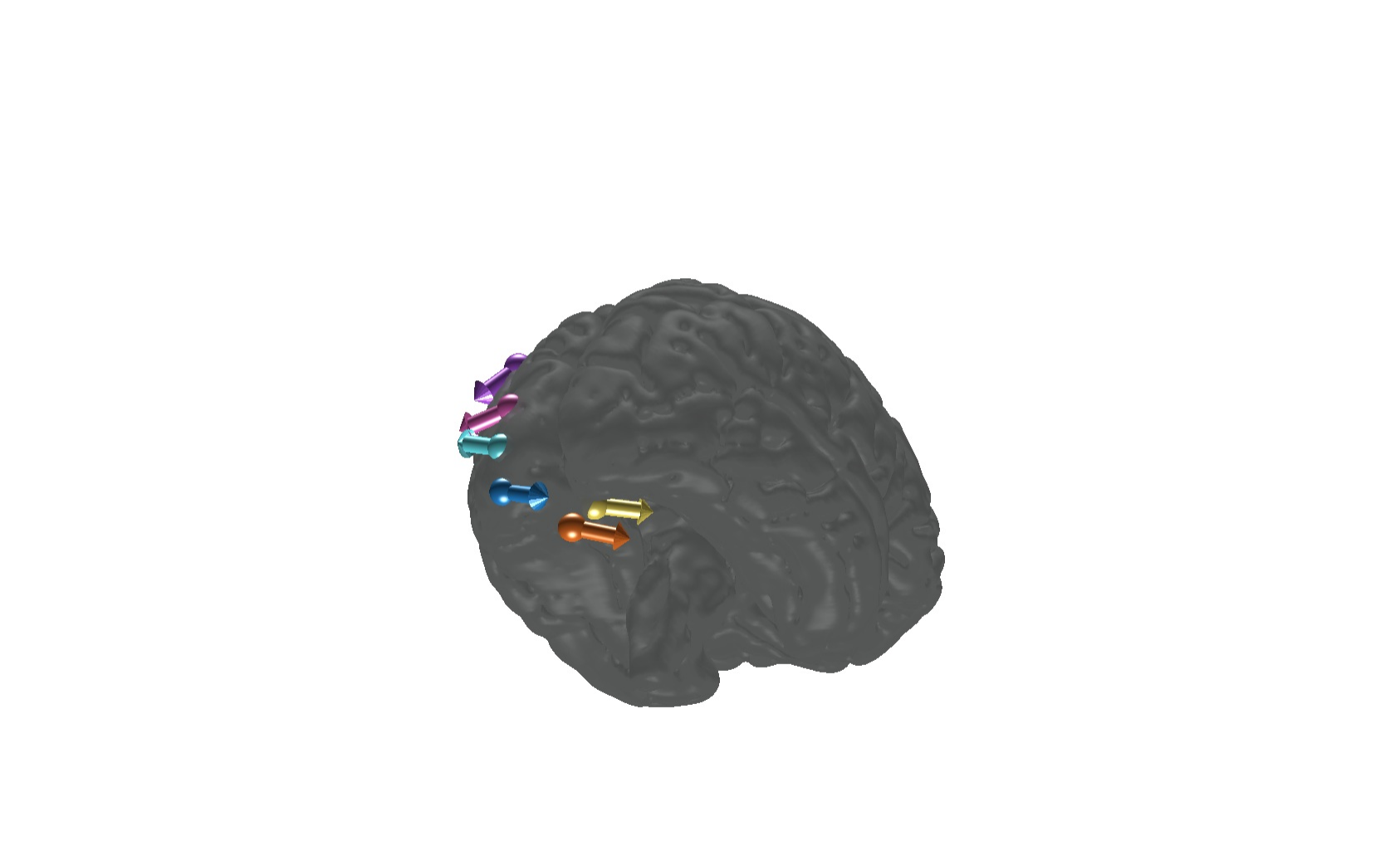}\hspace{0.5cm}\includegraphics[trim={5cm 2.5cm 9cm 5.5cm},clip,width=0.4\linewidth]{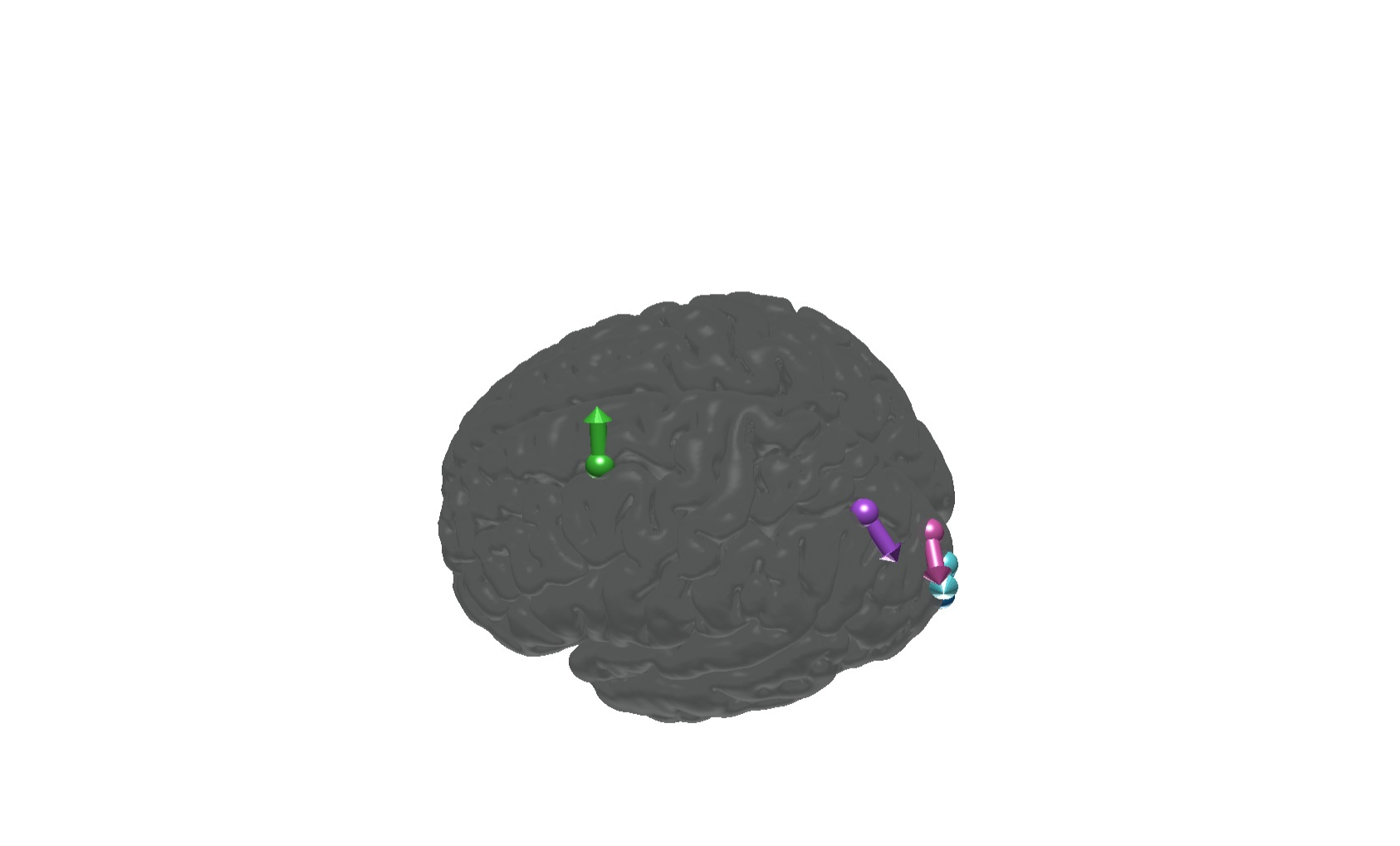}
    \caption{Dipolar source locations used to generate the synthetic FVEP. The sources represent the modeled visual-processing regions V1, V2v, V3v, V2d, V3d, LIP, and FEF.}
    \label{fig:sourceconfig}
\end{figure}

\begin{figure}[h!]
    \centering
    \includegraphics[width=0.75\linewidth]{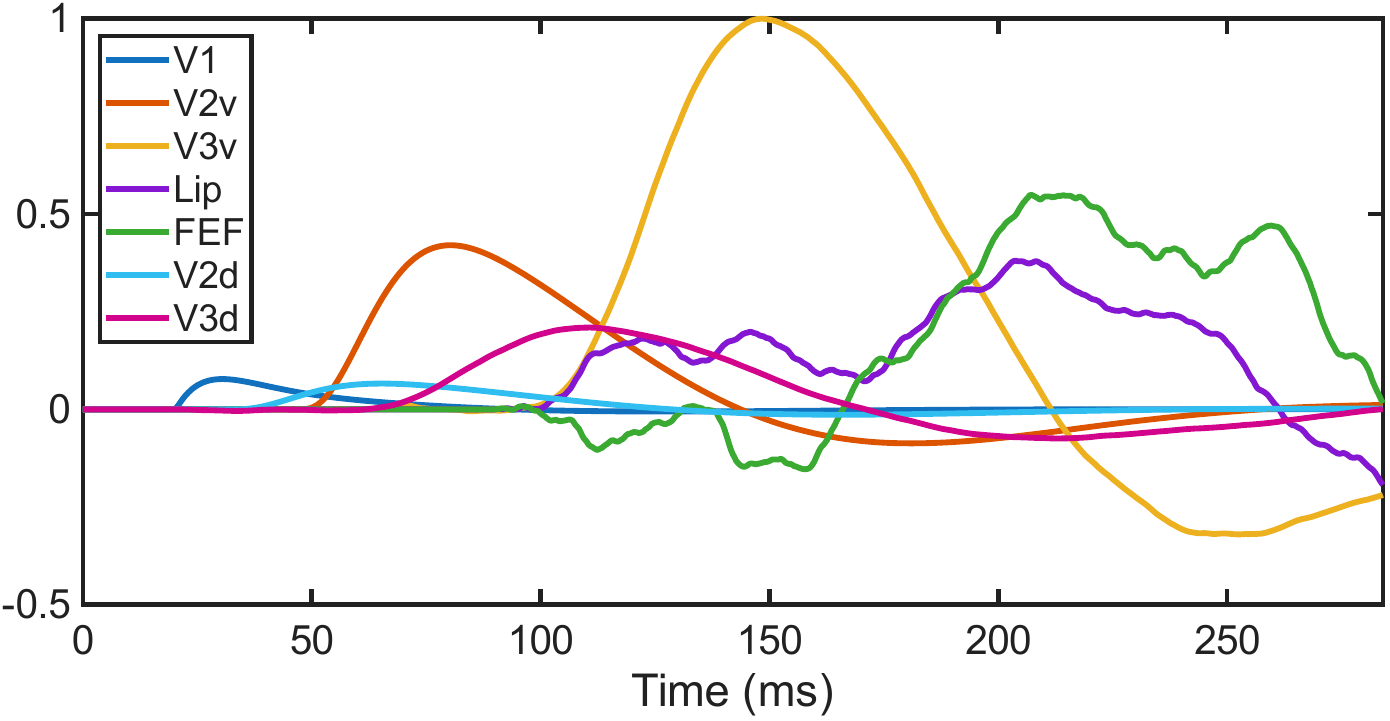}
    \caption{Source-level temporal evolution of the simulated FVEP components generated with the coupled Jansen--Rit neural mass model. Curves show the relative amplitudes of the modeled visual-network sources.}
    \label{fig:simulationcurves}
\end{figure}

Representative dipole locations were manually selected within the corresponding anatomical regions of the head model. The source orientations are the free three-dimensi\-onal dipole orientations represented by the lead field at those locations. Figure \ref{fig:sourceconfig} shows the selected sources, and Figure \ref{fig:simulationcurves} shows their relative time courses.

Propagation between successive source nodes is represented by a 10 ms inter-source delay. Relative output gains are 20 dB (V1), 47 dB (V2v), 36 dB (V3v), 23 dB (LIP), 25 dB (FEF), 45 dB (V2d), and 21 dB (V3d). The flash input is a Blackman--Harris pulse with 80 mV amplitude, 30 ms onset, and 3 ms duration. The source-generation simulation has duration 0.5 s and sampling rate 20 kHz and includes additive Gaussian noise with relative standard deviation $10^{-3}$. This source-generation noise is distinct from the 20 and 25 dB sensor-noise conditions used below to test inverse reconstruction.

\begin{table}[h!]
\centering
\begin{footnotesize}
\caption{Naming convention and functional interpretation of the cortical regions included in the simulated visual-processing network.}\label{tab:components}
\begin{tabular}{lll}
\hline
Naming convention & Functional area or interpretation \\
\hline
V1 &  Primary visual cortex \\
V2v &  Ventral visual pathway \\
V3v &  Higher ventral visual area \\
LIP &  Parietal integration area \\
FEF &  Frontal eye field \\
V2d &  Dorsal visual pathway \\
V3d &  Higher dorsal visual area \\
\hline
\end{tabular}
\end{footnotesize}
\end{table}

\begin{figure}[h!]
    \centering
    \includegraphics[trim={2cm 4cm 1cm 0},clip,width=0.4\linewidth]{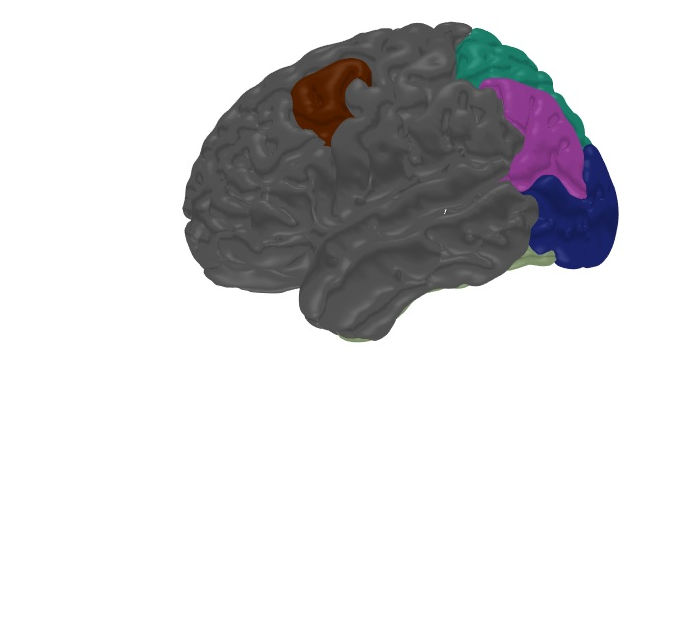}\hspace{0.5cm}\includegraphics[trim={2cm 4cm 1cm 0},clip,width=0.4\linewidth]{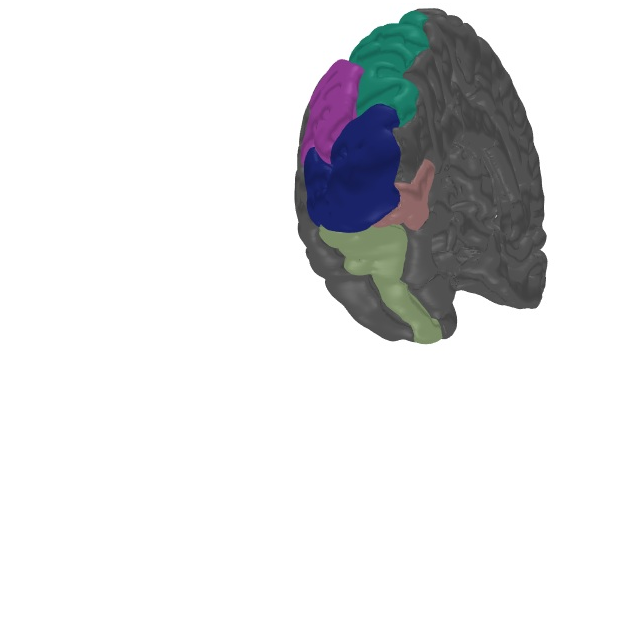}
    \caption{Desikan--Killiany regions of interest used for parcel-wise evaluation of the reconstructed source time series: lateral occipital (deep blue), lingual (beige), superior parietal (green), caudal middle frontal (brown), inferior parietal (violet), and fusiform (lime).}
    \label{fig:parcelROIs}
\end{figure}

We reconstruct the time series shown in Figure \ref{fig:simulationcurves} with the inverse methods described above. At each source location, the reconstructed magnitude is computed as the Euclidean norm of the three Cartesian dipole components. Parcel-wise activity is obtained by aggregating these magnitudes over the source locations belonging to each parcel and normalizing each resulting curve to its peak value for temporal comparison. The regions of interest are defined with the FreeSurfer parcellation\footnote{\url{https://surfer.nmr.mgh.harvard.edu/}} through SPM12 \citep{SPM12} and the Desikan--Killiany atlas \cite{DKatlas2006}. The selected regions are shown in Figure \ref{fig:parcelROIs}.

We use the population-average ICBM152 2009a multicompartment head model \cite{fonov2009,fonov2011,icbm-website}. The volume conductor includes scalp, skull, cerebrospinal fluid, gray matter, white matter, and the principal deeper compartments represented in the segmentation; conductivity values follow the tissue parameterization used by Dannhauer {\em et al.} \cite{dannhauer2011}. The geometry is discretized with a boundary-fitted tetrahedral mesh, and the EEG lead field is computed from the quasi-static Maxwell equations using the finite element method. The source space contains 10,000 fixed locations with free three-dimensional orientation. Electrodes follow the international 10--10 placement system. The same lead field is used for synthetic data generation and inverse reconstruction, so the present benchmark intentionally excludes geometric forward-model mismatch and isolates differences between inverse methods. The finite element solution is computed using the Local subtraction approach \cite{hoeltershinken2025LS}, where the total potential field is split into a known analytical singularity (the primary source) and an unknown correction potential. The singularity is analytically subtracted from the governing equation to make the remaining problem smooth and regular.

In this study, we assess localization accuracy and focality qualitatively by projecting and interpolating cortical activity on the brain model, and quantitatively using Earth Mover's Distance (EMD). EMD compares the reconstructed and ground-truth source distributions after normalization to unit mass. It is the minimum transport cost required to transform one spatial distribution into the other, and therefore reflects both localization error and distributional mismatch.

\section{Results}

\subsection{Assessment of the synthetic FVEP data}

\begin{figure}[h!]
    \centering
    \includegraphics[width=0.5\linewidth]{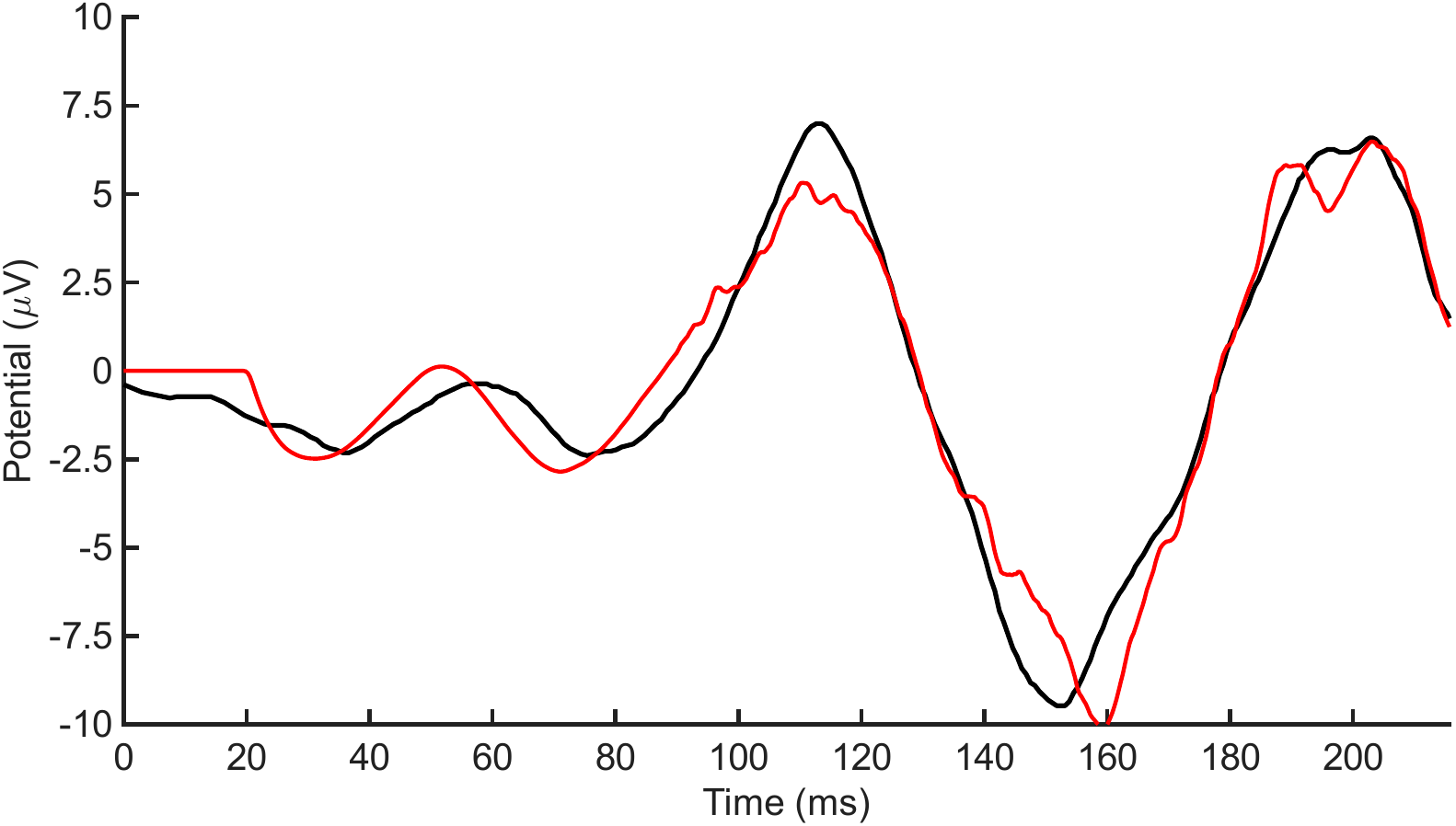}
    \caption{Simulated FVEP at electrode PO7 relative to the reference electrode (red), compared with the representative FVEP waveform (black) presented in \cite{odom2010TrueVEPPO7,AfrasiabiS2012TrueFVEPPO7}.}
    \label{fig:PO7waveform}
\end{figure}

To assess whether the simulated sensor-level response has a plausible FVEP waveform of trace, the PO7 signal relative to the reference electrode was compared with a representative FVEP waveform \cite{odom2010TrueVEPPO7}. Figure \ref{fig:PO7waveform} shows agreement in the principal polarity sequence and relative component structure. Peak times do not coincide exactly, as expected, because experimentally reported component latencies vary substantially across subjects and stimulation protocols. We therefore use this comparison to qualitatively validate the waveform rather than as a quantitative goodness-of-fit test. The synthetic waveform is intended primarily to provide physiologically plausible temporal structure with known source-level ground truth.

\subsection{Parcellation-wise time series reconstructions}

\begin{figure}[h!]
    \centering
    \begin{minipage}{0.02\linewidth}
        \rotatebox{90}{\scriptsize \hfill}
    \end{minipage}\begin{minipage}{0.196\linewidth}
        \centering
        {\bf Ground truth}
    \end{minipage}\begin{minipage}{0.196\linewidth}
        \centering
        {\bf CrSKF}
    \end{minipage}\begin{minipage}{0.196\linewidth}
        \centering
        {\bf SKF}
    \end{minipage}\begin{minipage}{0.196\linewidth}
        \centering
        {\bf KF}
    \end{minipage}\begin{minipage}{0.196\linewidth}
        \centering
        {\bf LCMV}
    \end{minipage}

    \begin{minipage}{0.02\linewidth}
        \rotatebox{90}{\scriptsize L.O.}
    \end{minipage}\begin{minipage}{0.196\linewidth}
        \centering
        \includegraphics[width=0.99\linewidth]{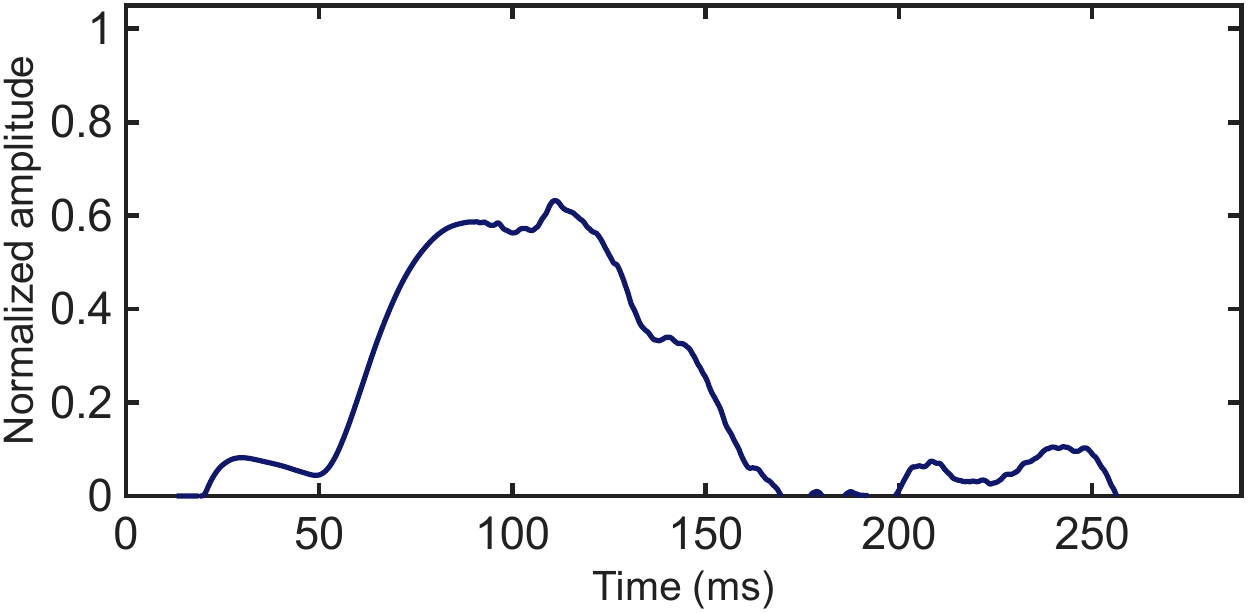}
    \end{minipage}\begin{minipage}{0.2\linewidth}
        \centering
        \includegraphics[width=0.99\linewidth]{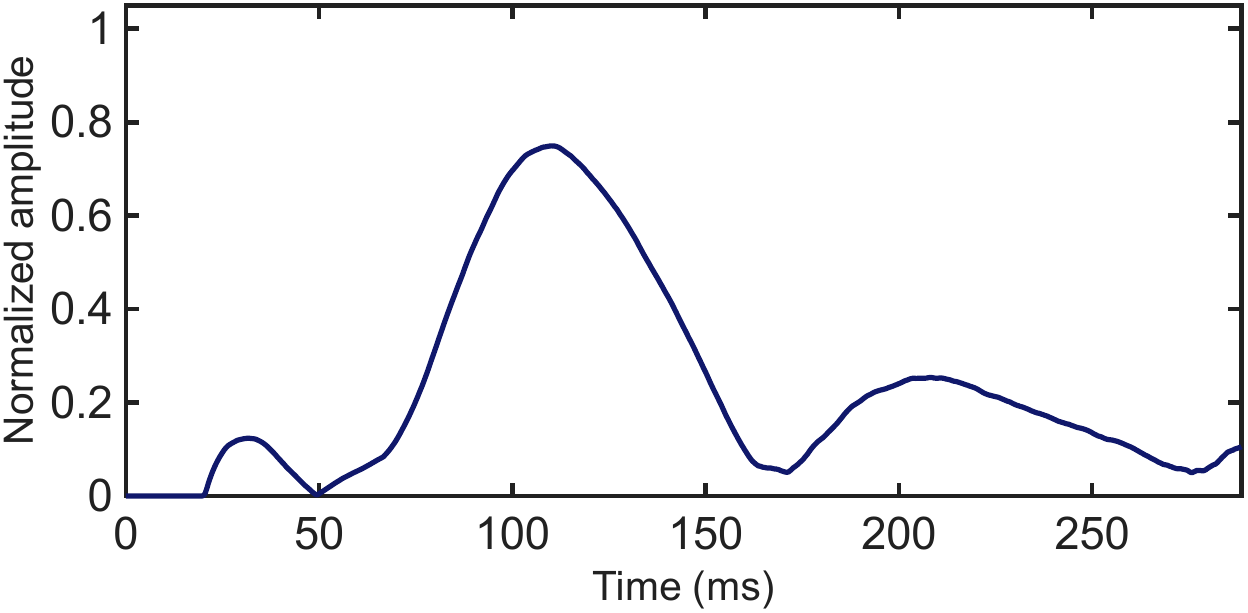}
    \end{minipage}\begin{minipage}{0.196\linewidth}
        \centering
        \includegraphics[width=0.99\linewidth]{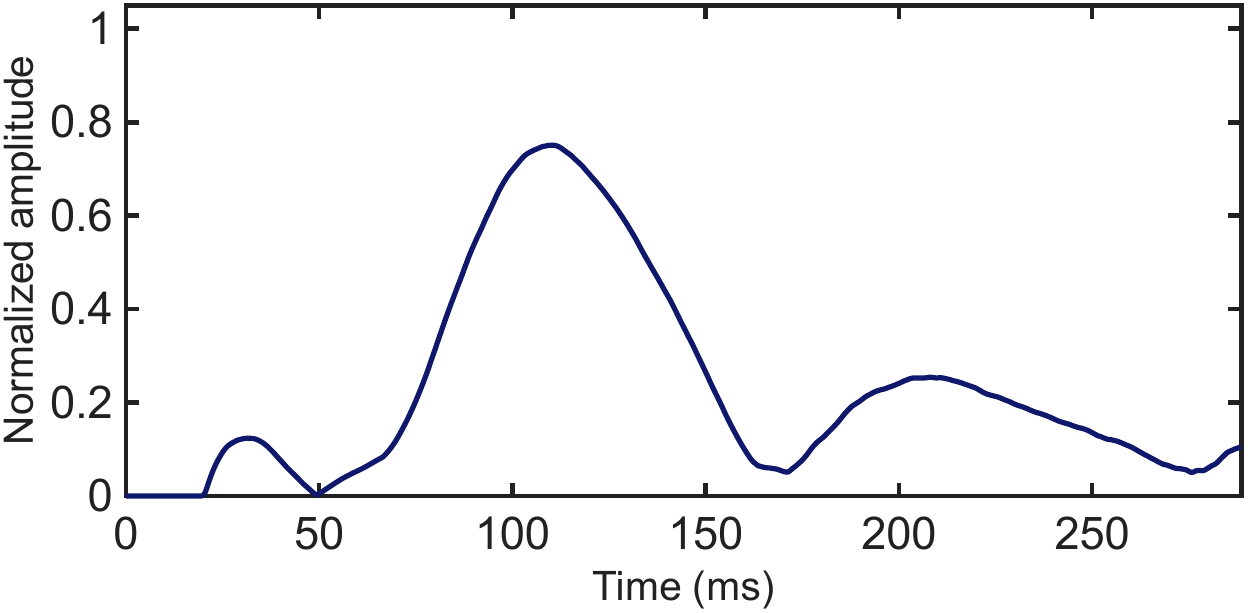}
    \end{minipage}\begin{minipage}{0.196\linewidth}
        \centering
        \includegraphics[width=0.99\linewidth]{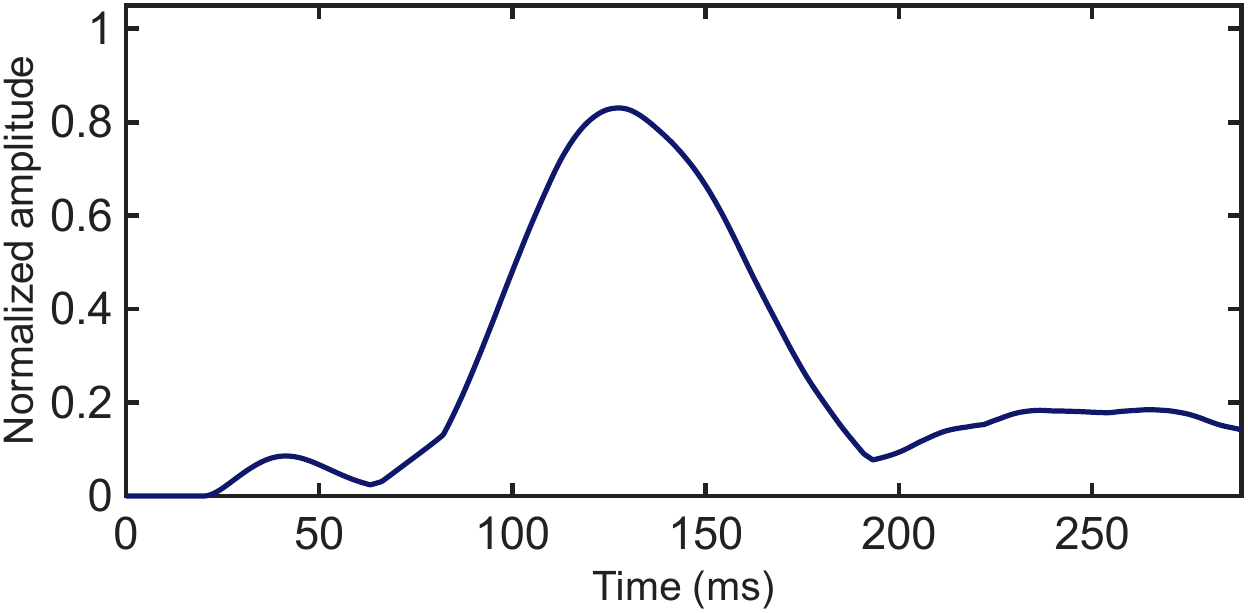}
    \end{minipage}\begin{minipage}{0.196\linewidth}
        \centering
        \includegraphics[width=0.99\linewidth]{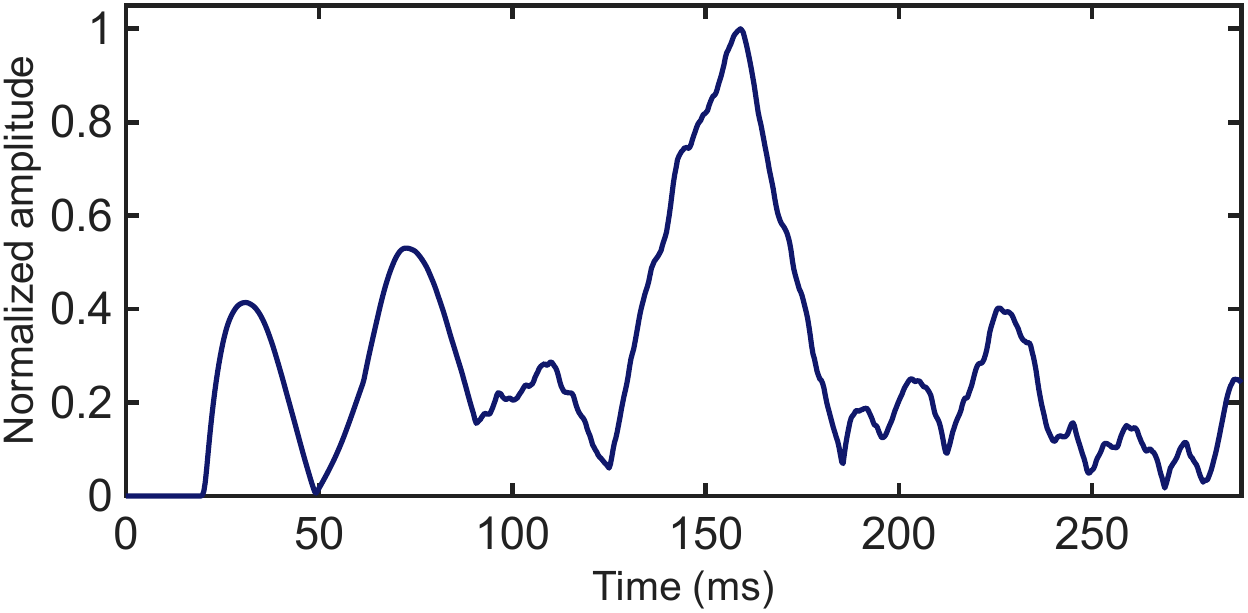}
    \end{minipage}

    \begin{minipage}{0.02\linewidth}
        \rotatebox{90}{\scriptsize Ling.}
    \end{minipage}\begin{minipage}{0.196\linewidth}
        \centering
        \includegraphics[width=0.99\linewidth]{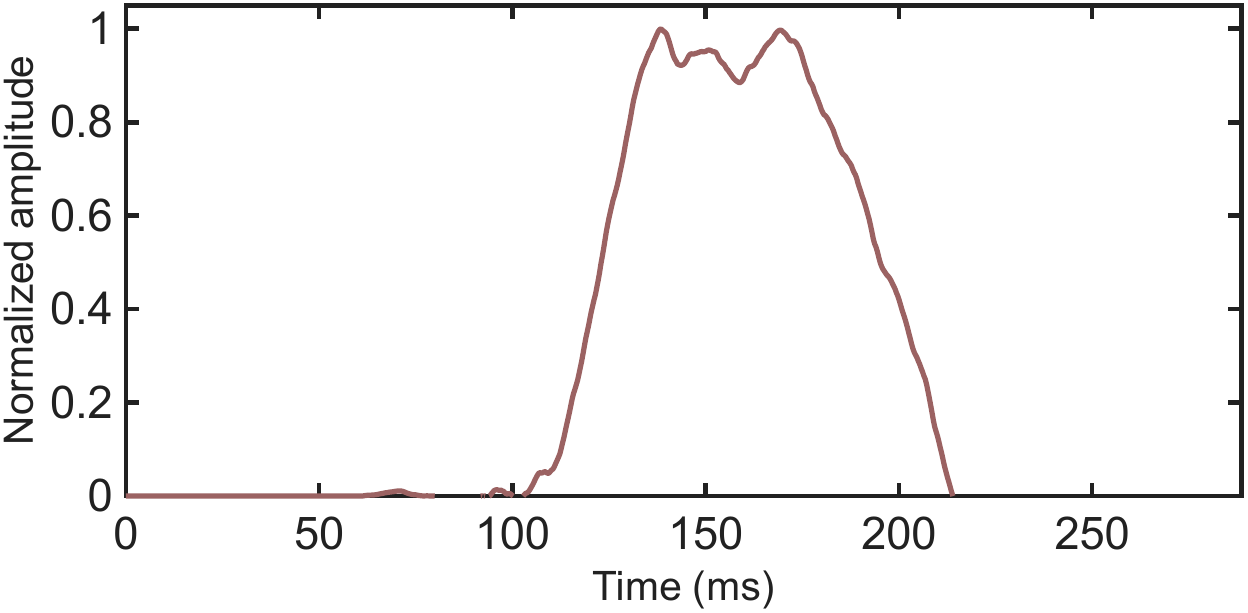}
    \end{minipage}\begin{minipage}{0.2\linewidth}
        \centering
        \includegraphics[width=0.99\linewidth]{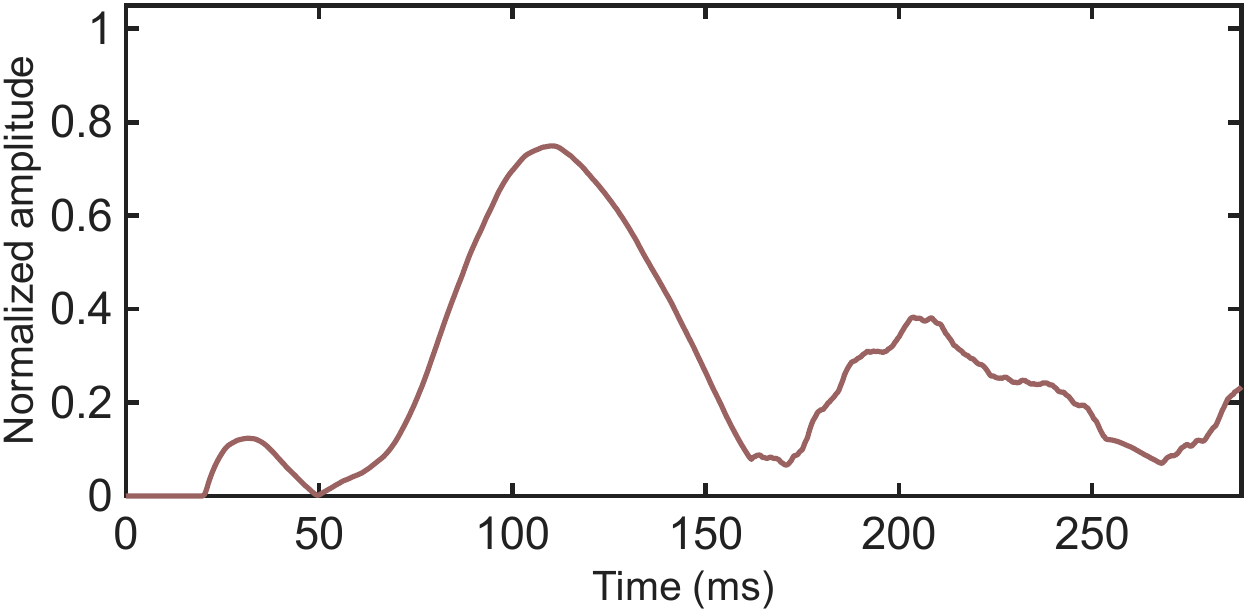}
    \end{minipage}\begin{minipage}{0.196\linewidth}
        \centering
        \includegraphics[width=0.99\linewidth]{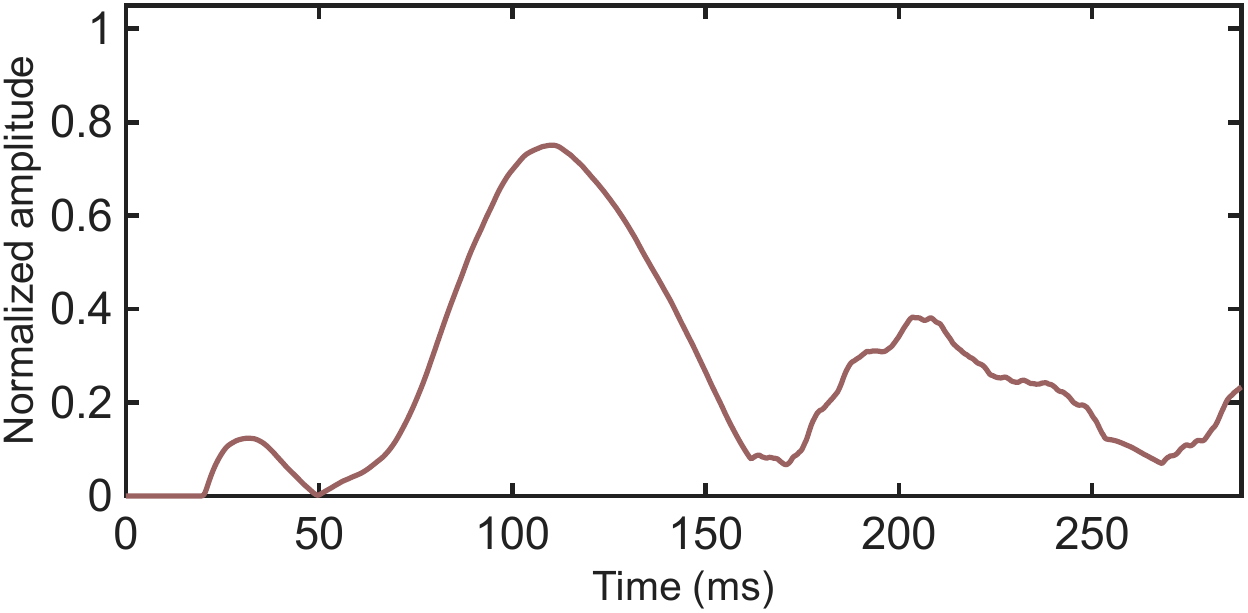}
    \end{minipage}\begin{minipage}{0.196\linewidth}
        \centering
        \includegraphics[width=0.99\linewidth]{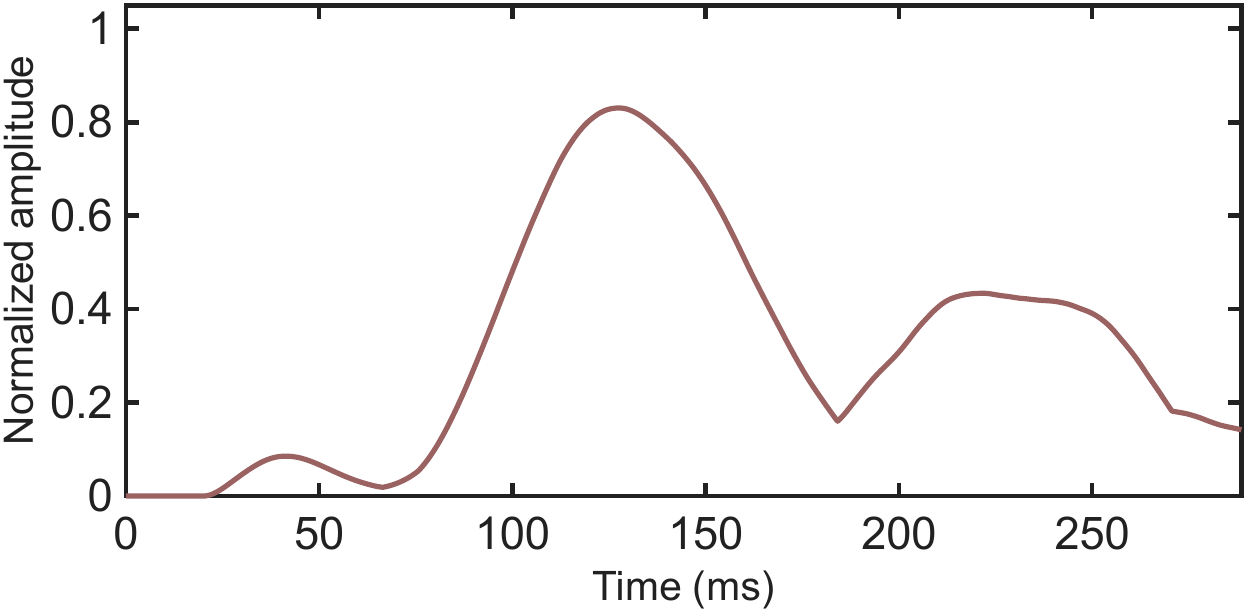}
    \end{minipage}\begin{minipage}{0.196\linewidth}
        \centering
        \includegraphics[width=0.99\linewidth]{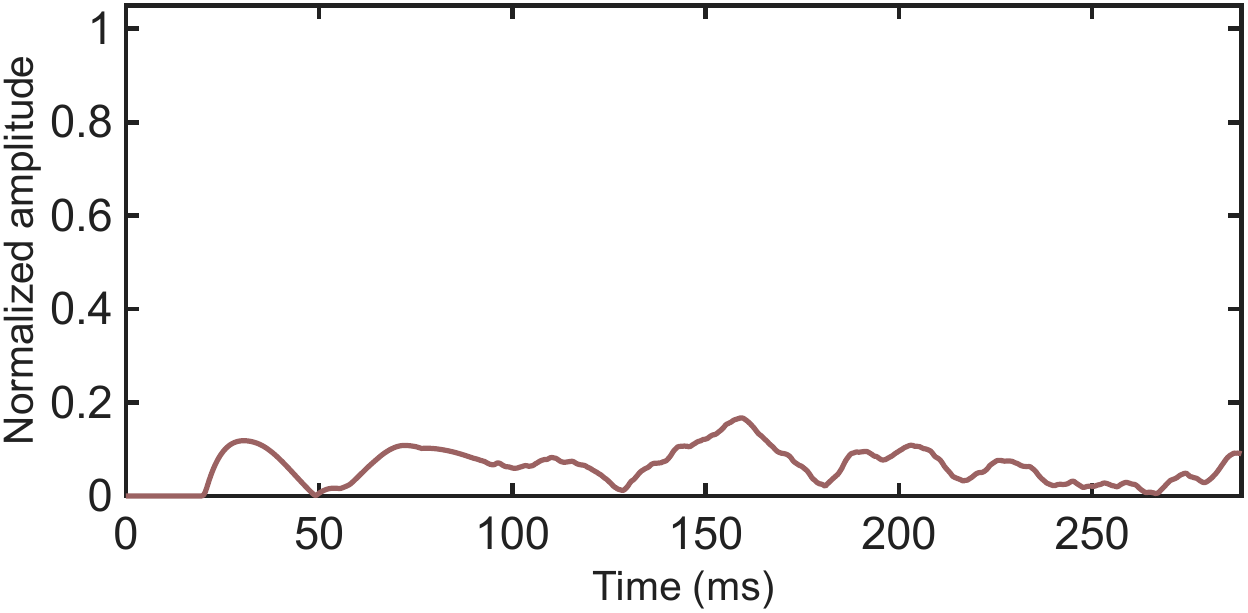}
    \end{minipage}

    \begin{minipage}{0.02\linewidth}
        \rotatebox{90}{\scriptsize S.P.}
    \end{minipage}\begin{minipage}{0.196\linewidth}
        \centering
        \includegraphics[width=0.99\linewidth]{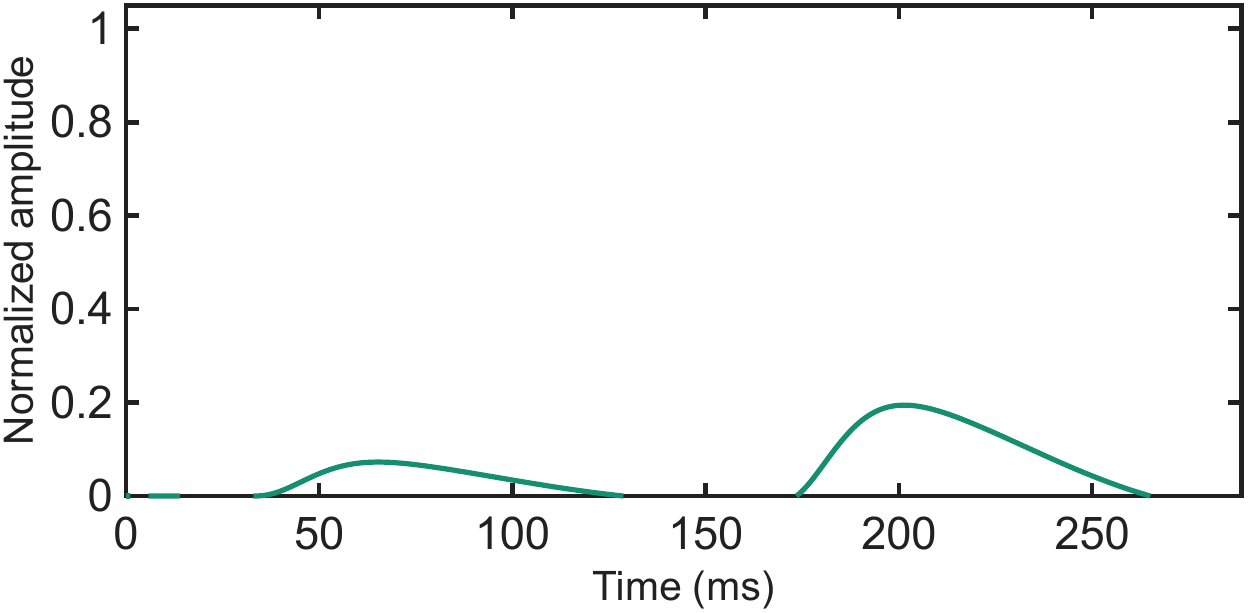}
    \end{minipage}\begin{minipage}{0.2\linewidth}
        \centering
        \includegraphics[width=0.99\linewidth]{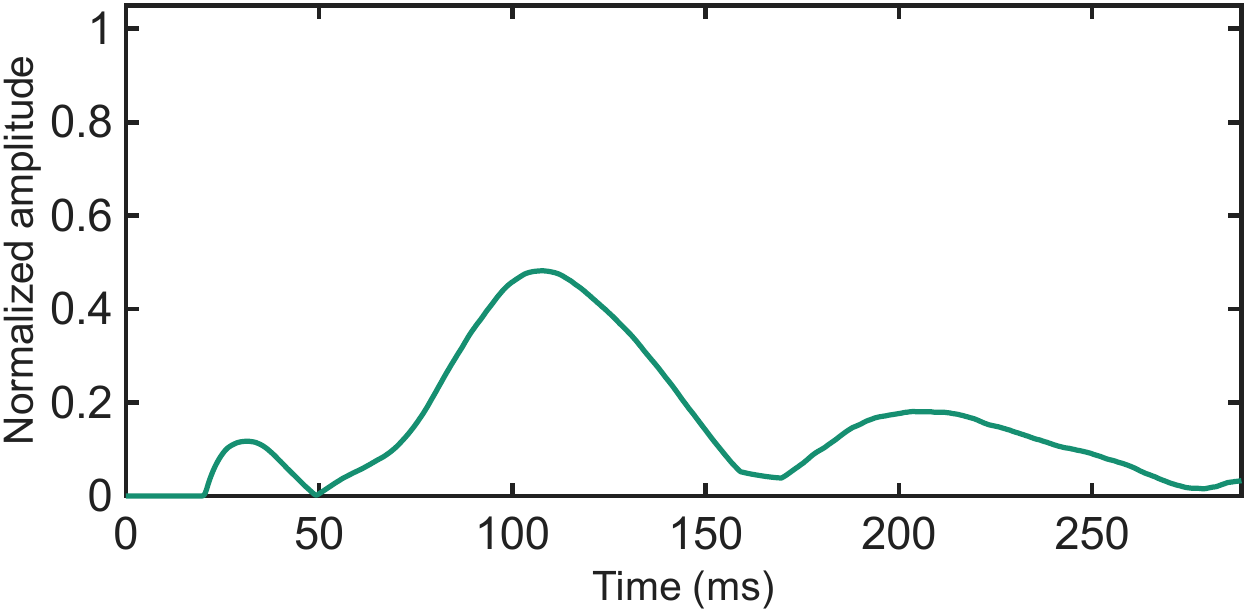}
    \end{minipage}\begin{minipage}{0.196\linewidth}
        \centering
        \includegraphics[width=0.99\linewidth]{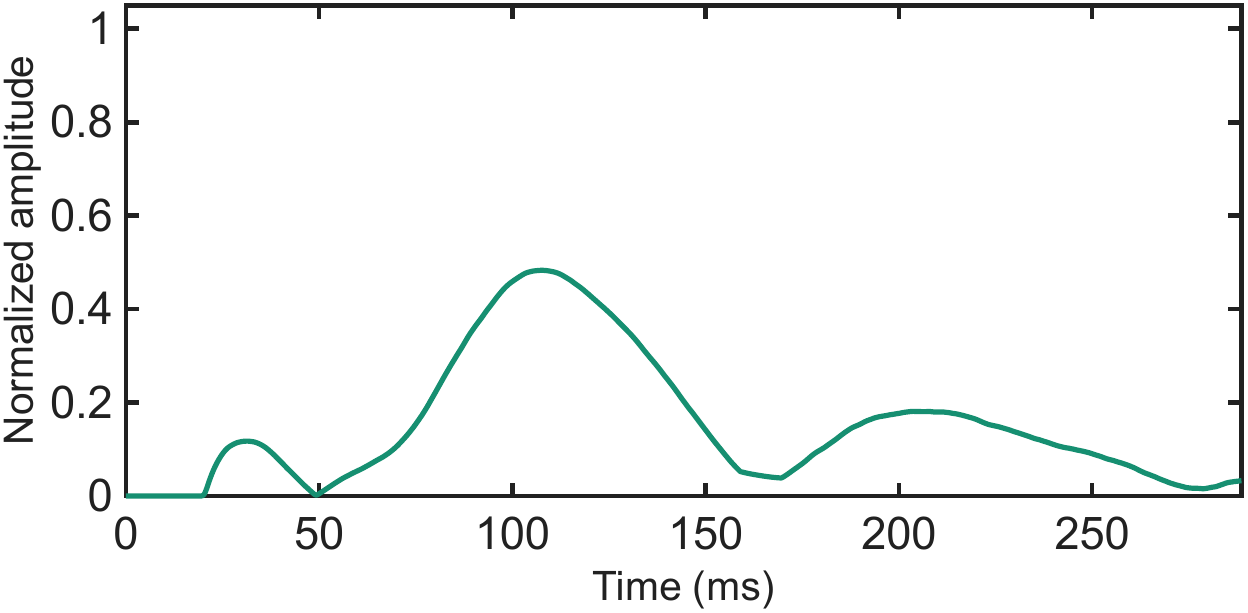}
    \end{minipage}\begin{minipage}{0.196\linewidth}
        \centering
        \includegraphics[width=0.99\linewidth]{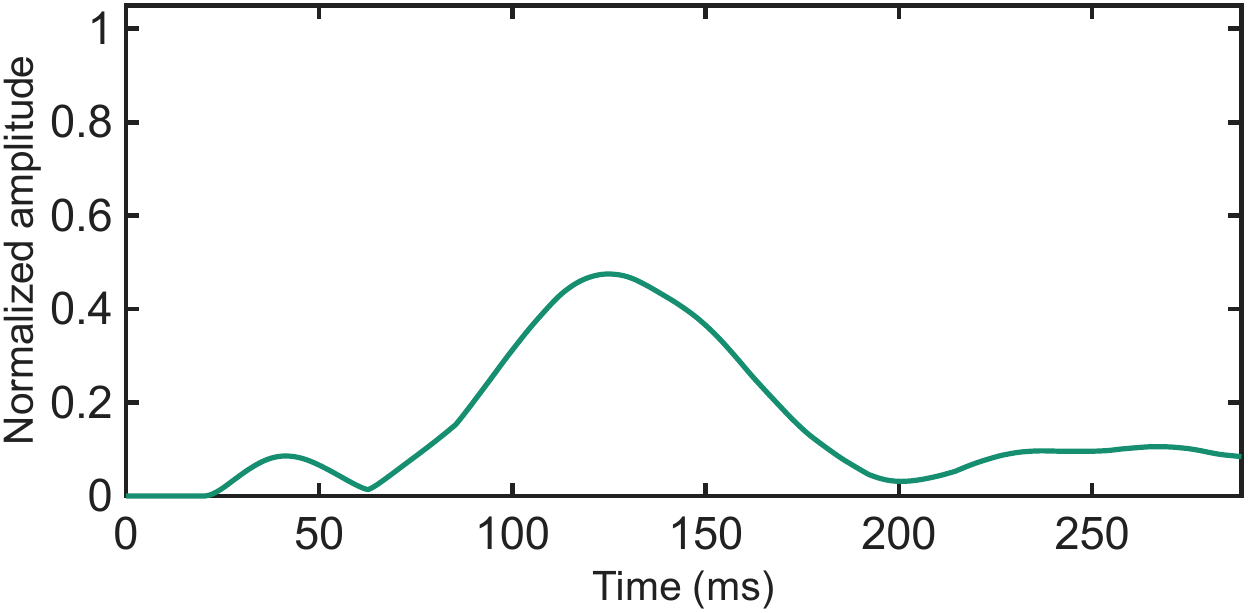}
    \end{minipage}\begin{minipage}{0.196\linewidth}
        \centering
        \includegraphics[width=0.99\linewidth]{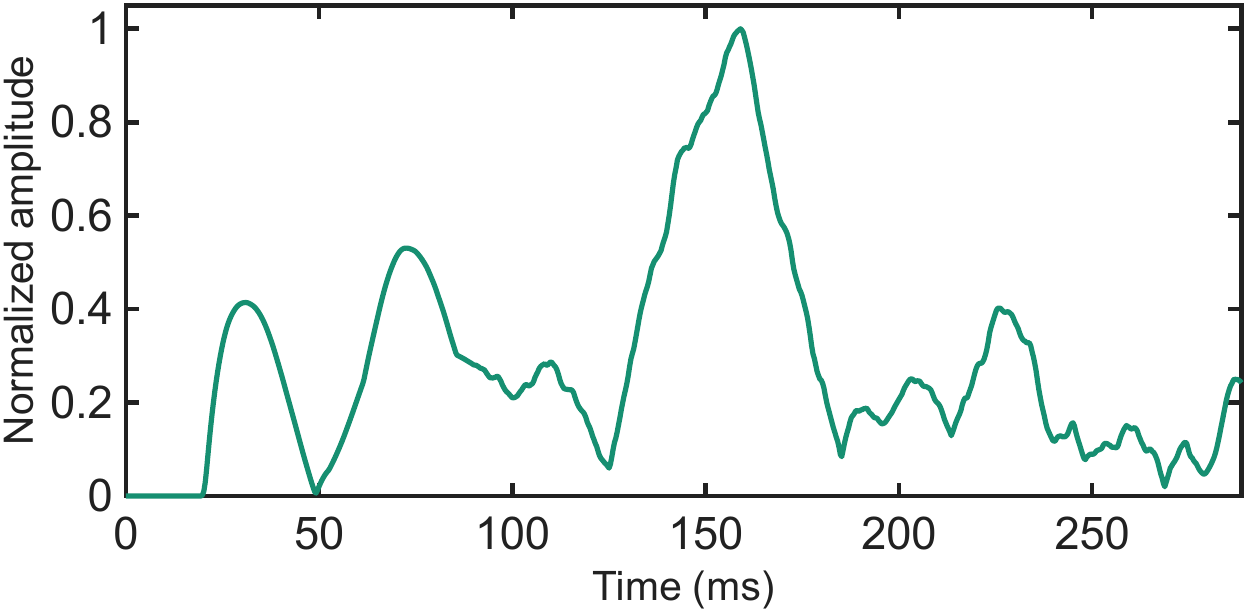}
    \end{minipage}

    \begin{minipage}{0.02\linewidth}
        \rotatebox{90}{\scriptsize C.M.F.}
    \end{minipage}\begin{minipage}{0.196\linewidth}
        \centering
        \includegraphics[width=0.99\linewidth]{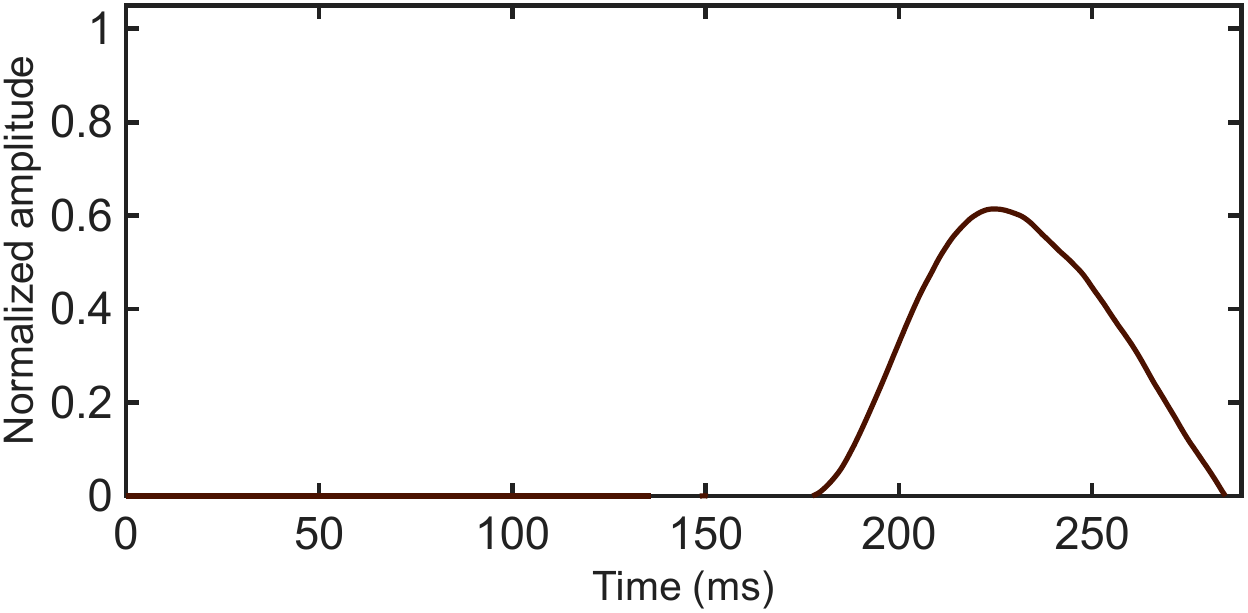}
    \end{minipage}\begin{minipage}{0.2\linewidth}
        \centering
        \includegraphics[width=0.99\linewidth]{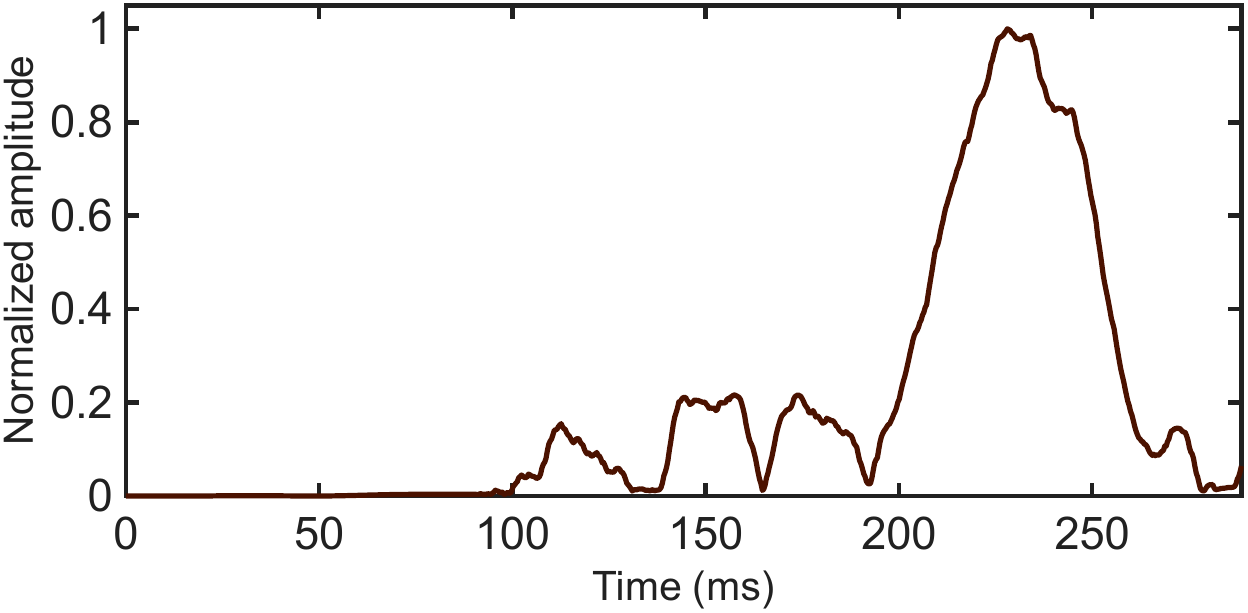}
    \end{minipage}\begin{minipage}{0.196\linewidth}
        \centering
        \includegraphics[width=0.99\linewidth]{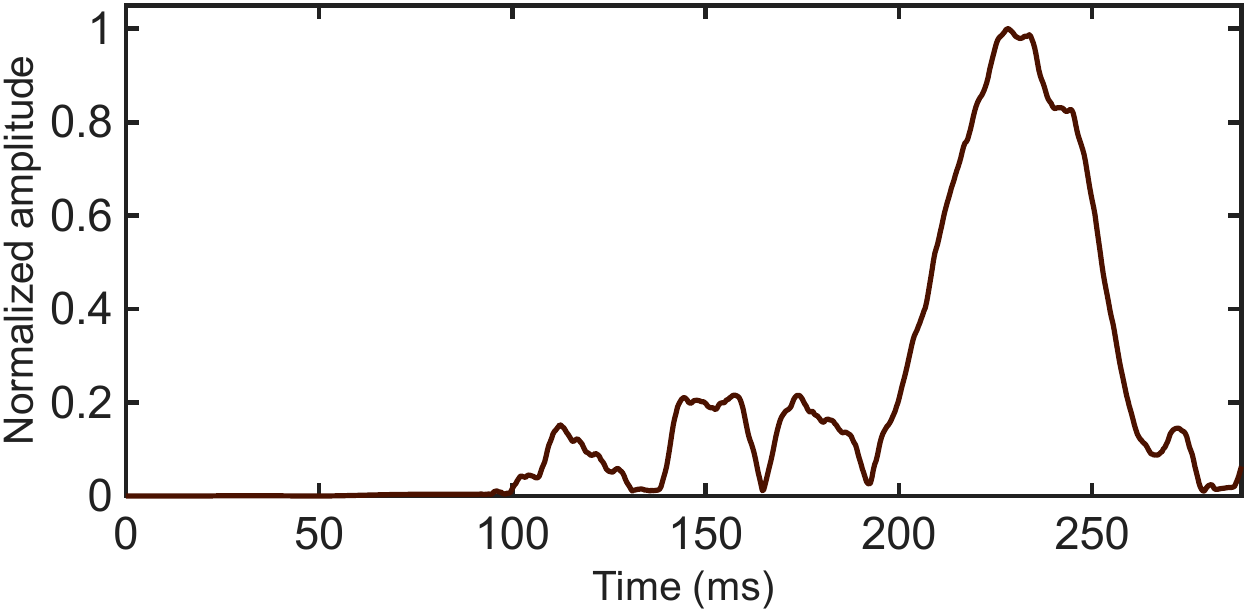}
    \end{minipage}\begin{minipage}{0.196\linewidth}
        \centering
        \includegraphics[width=0.99\linewidth]{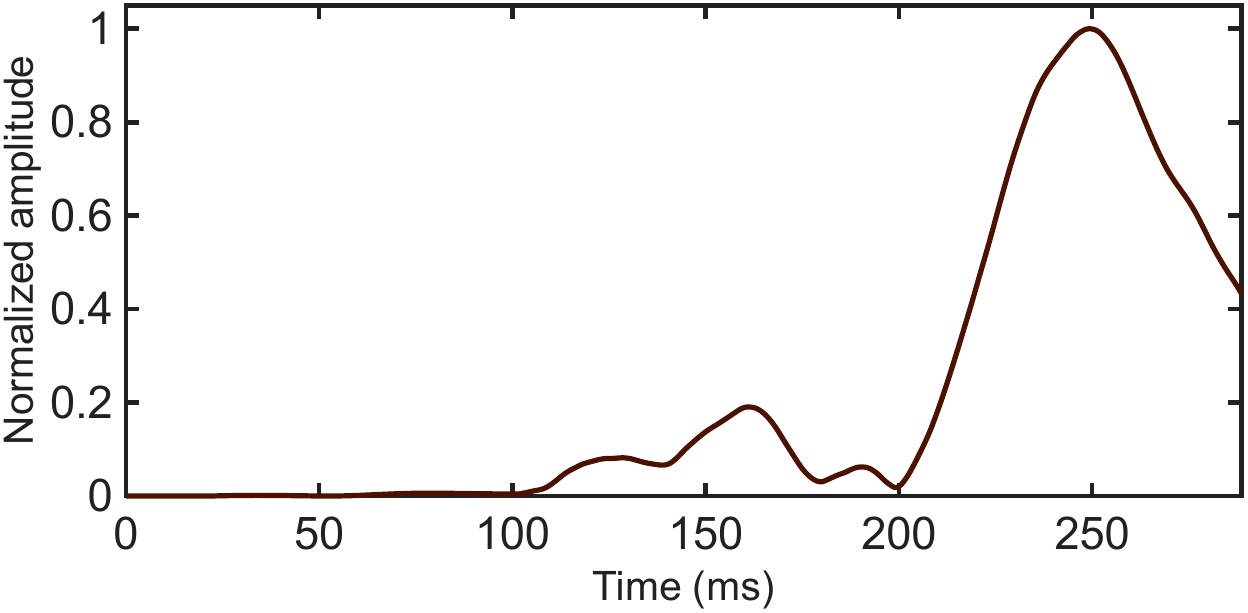}
    \end{minipage}\begin{minipage}{0.196\linewidth}
        \centering
        \includegraphics[width=0.99\linewidth]{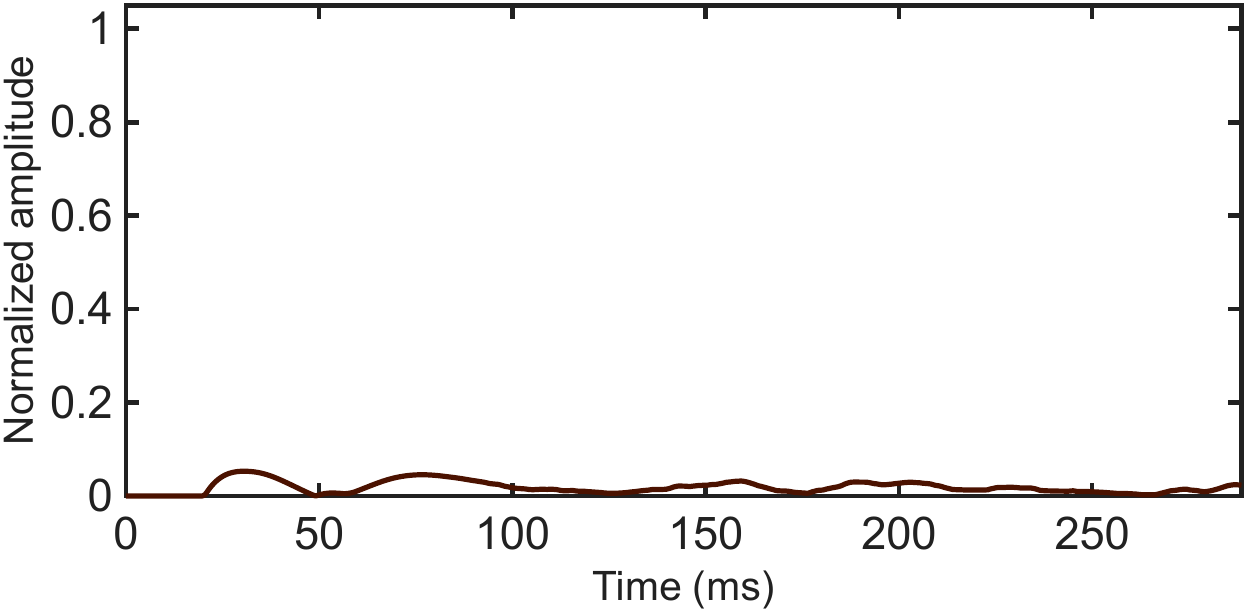}
    \end{minipage}

    \begin{minipage}{0.02\linewidth}
        \rotatebox{90}{\scriptsize I.P.}
    \end{minipage}\begin{minipage}{0.196\linewidth}
        \centering
        \hfill
        \end{minipage}\begin{minipage}{0.2\linewidth}
        \centering
        \includegraphics[width=0.99\linewidth]{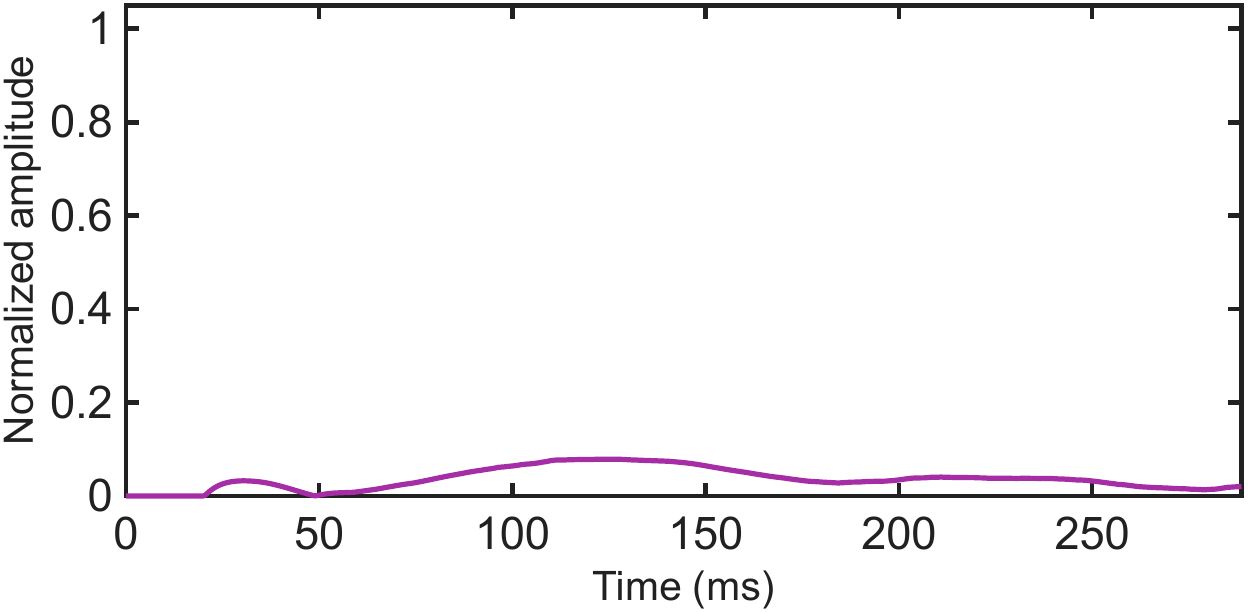}
    \end{minipage}\begin{minipage}{0.196\linewidth}
        \centering
        \includegraphics[width=0.99\linewidth]{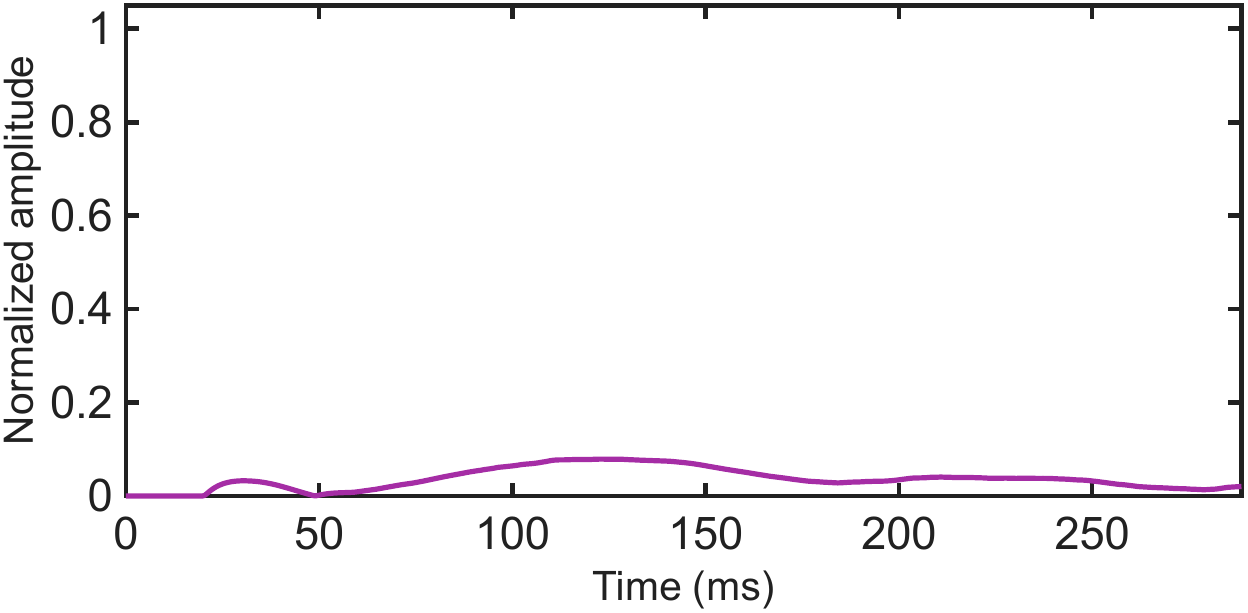}
    \end{minipage}\begin{minipage}{0.196\linewidth}
        \centering
        \includegraphics[width=0.99\linewidth]{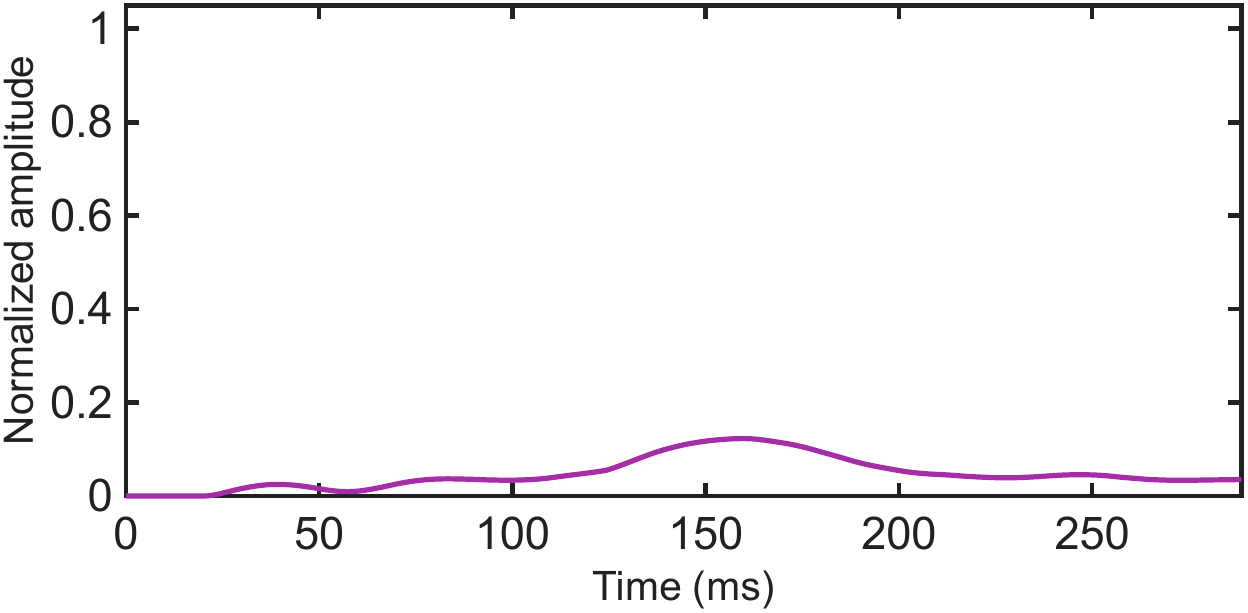}
    \end{minipage}\begin{minipage}{0.196\linewidth}
        \centering
        \includegraphics[width=0.99\linewidth]{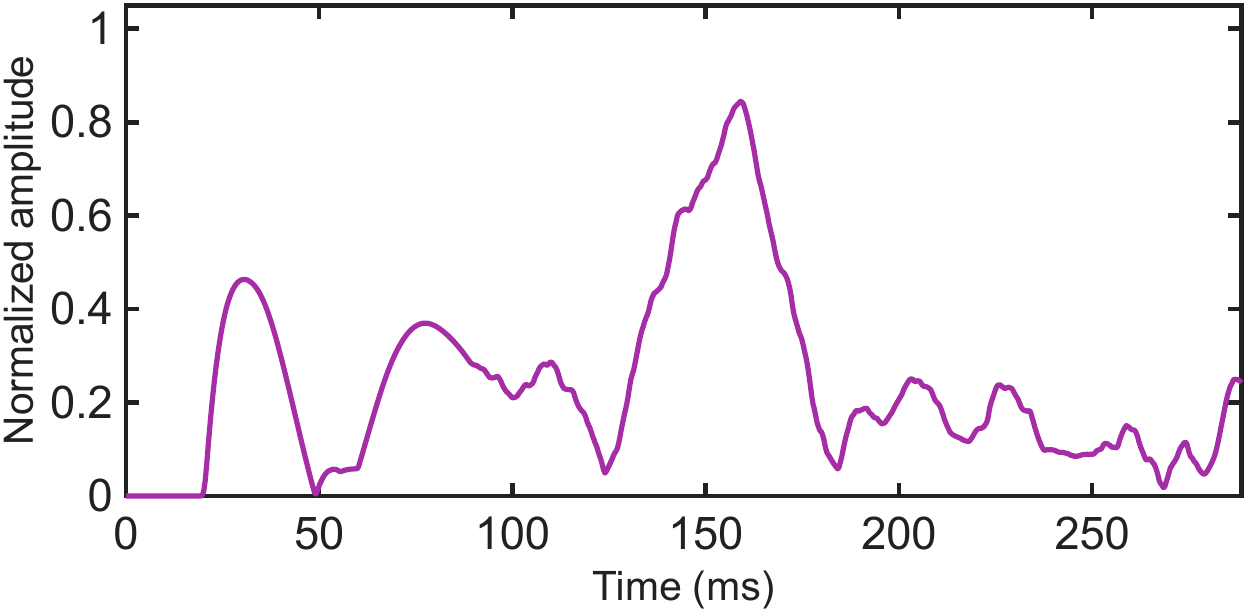}
    \end{minipage}

    \begin{minipage}{0.02\linewidth}
        \rotatebox{90}{\scriptsize Fusif.}
    \end{minipage}\begin{minipage}{0.196\linewidth}
        \centering
        \hfill
    \end{minipage}\begin{minipage}{0.2\linewidth}
        \centering
        \includegraphics[width=0.99\linewidth]{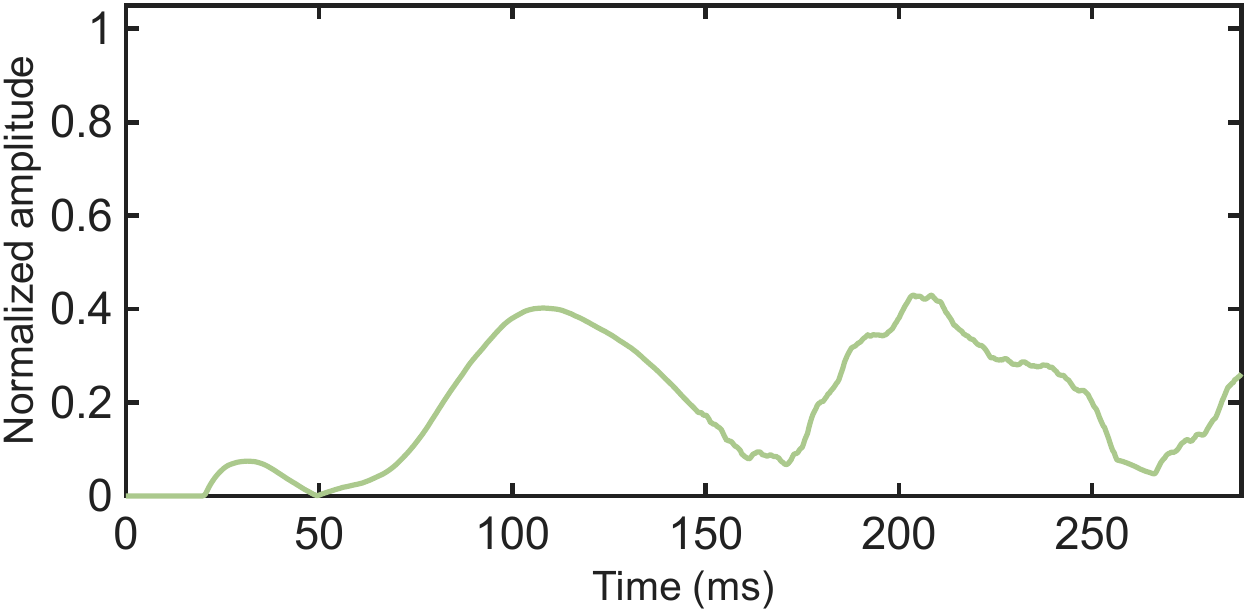}
    \end{minipage}\begin{minipage}{0.196\linewidth}
        \centering
        \includegraphics[width=0.99\linewidth]{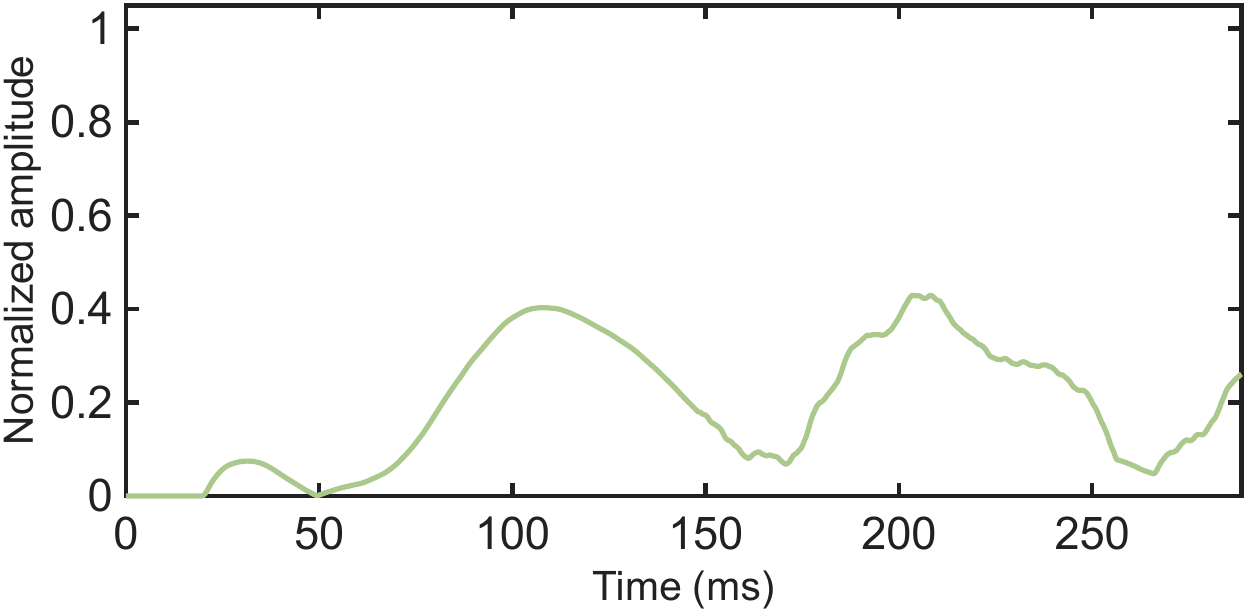}
    \end{minipage}\begin{minipage}{0.196\linewidth}
        \centering
        \includegraphics[width=0.99\linewidth]{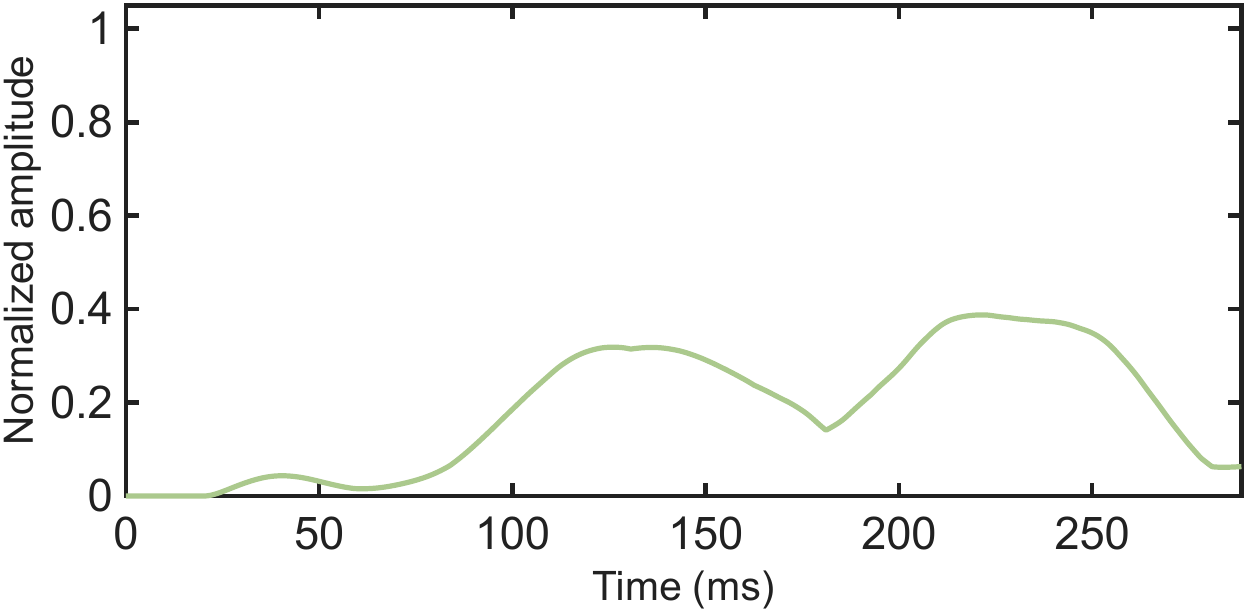}
    \end{minipage}\begin{minipage}{0.196\linewidth}
        \centering
        \includegraphics[width=0.99\linewidth]{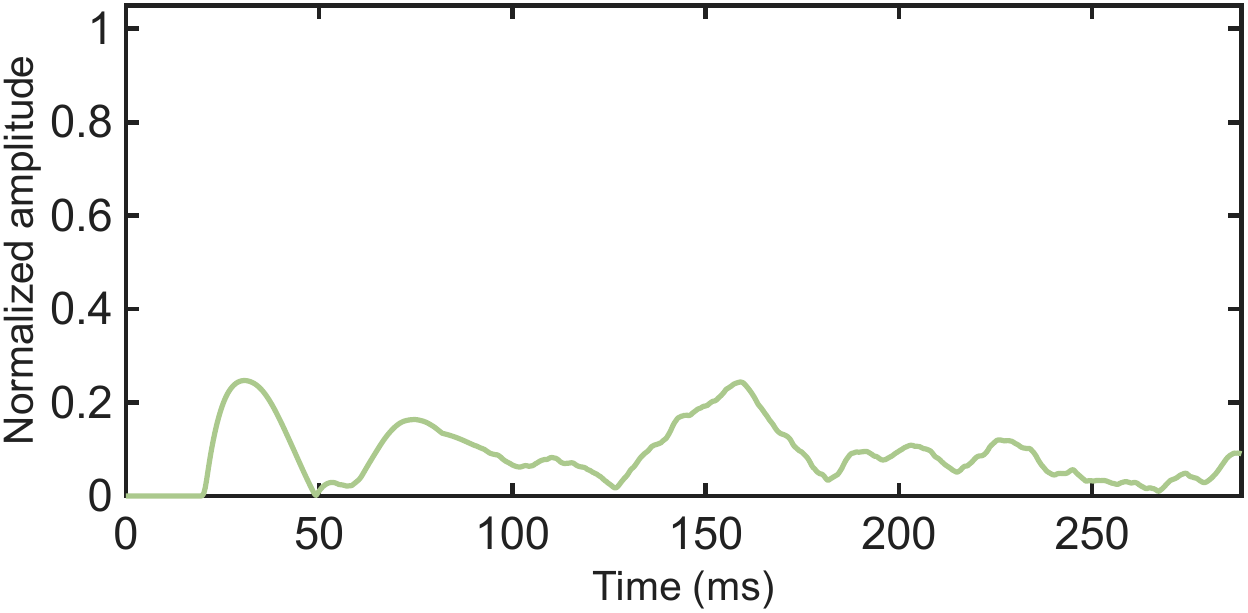}
    \end{minipage}
    \caption{Parcel-wise temporal activity under the noiseless condition. Rows correspond to the selected cortical parcels, and columns show the simulated ground-truth time series followed by CrSKF, SKF, KF, and LCMV reconstructions.}
    \label{fig:timeseries0dB}
\end{figure}

\begin{figure}[h!]
    \centering
    \begin{minipage}{0.02\linewidth}
        \rotatebox{90}{\scriptsize \hfill}
    \end{minipage}\begin{minipage}{0.196\linewidth}
        \centering
        {\bf Ground truth}
    \end{minipage}\begin{minipage}{0.196\linewidth}
        \centering
        {\bf CrSKF}
    \end{minipage}\begin{minipage}{0.196\linewidth}
        \centering
        {\bf SKF}
    \end{minipage}\begin{minipage}{0.196\linewidth}
        \centering
        {\bf KF}
    \end{minipage}\begin{minipage}{0.196\linewidth}
        \centering
        {\bf LCMV}
    \end{minipage}

    \begin{minipage}{0.02\linewidth}
        \rotatebox{90}{\scriptsize L.O.}
    \end{minipage}\begin{minipage}{0.196\linewidth}
        \centering
        \includegraphics[width=0.99\linewidth]{pictures/timeseries/true_timeseries_lateraloccipital.pdf}
    \end{minipage}\begin{minipage}{0.2\linewidth}
        \centering
        \includegraphics[width=0.99\linewidth]{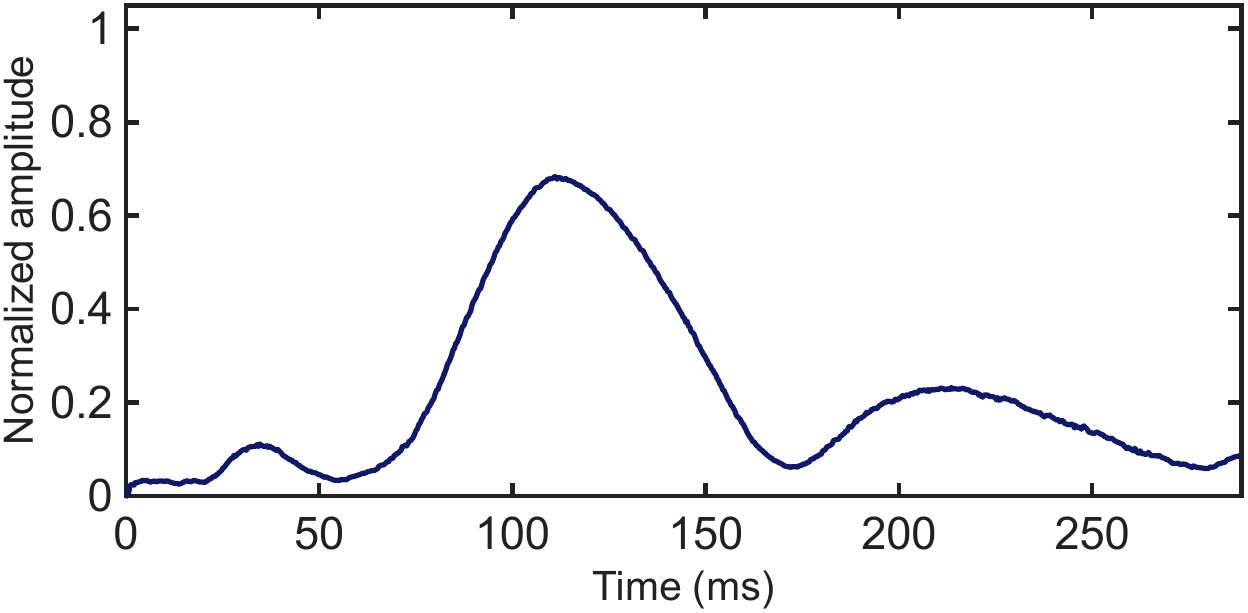}
    \end{minipage}\begin{minipage}{0.196\linewidth}
        \centering
        \includegraphics[width=0.99\linewidth]{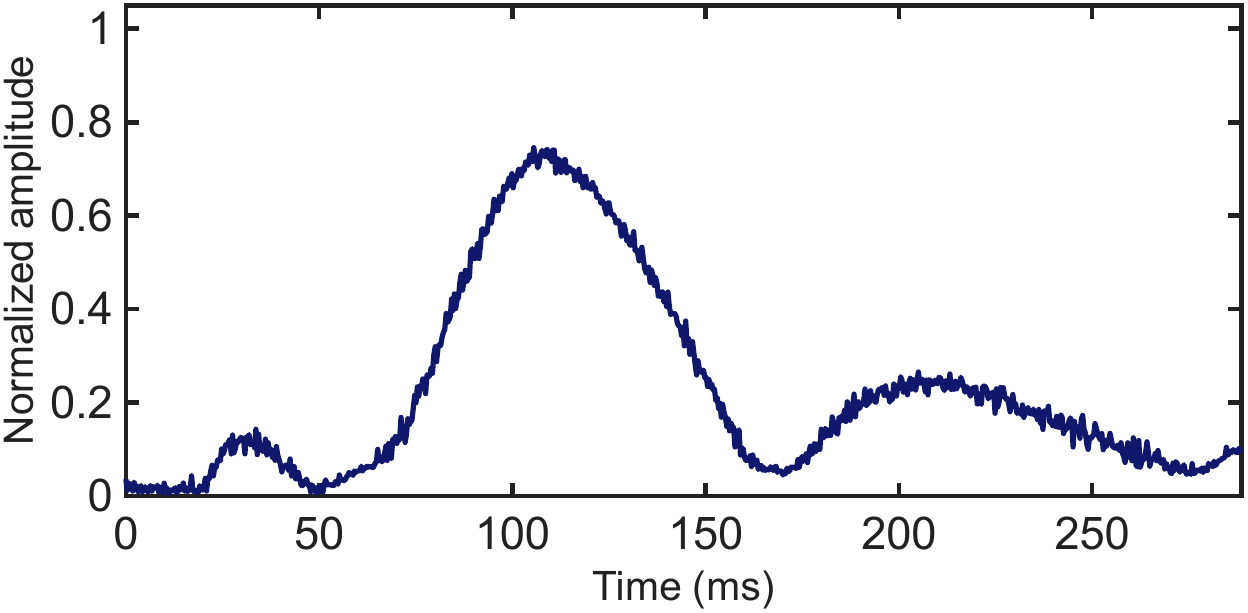}
    \end{minipage}\begin{minipage}{0.196\linewidth}
        \centering
        \includegraphics[width=0.99\linewidth]{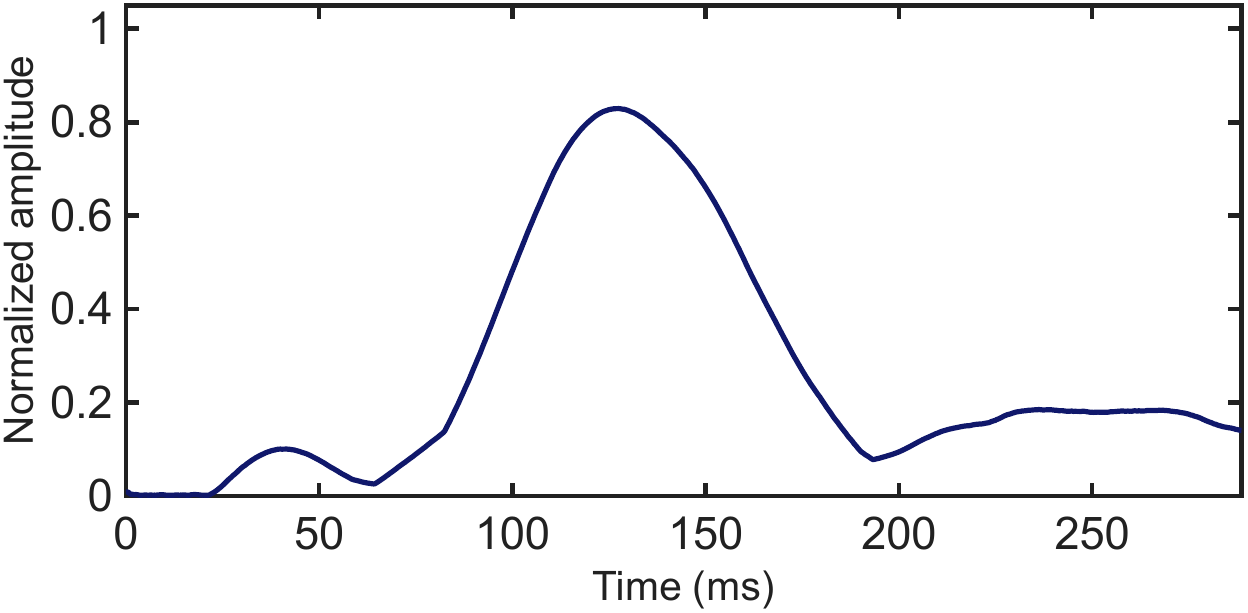}
    \end{minipage}\begin{minipage}{0.196\linewidth}
        \centering
        \includegraphics[width=0.99\linewidth]{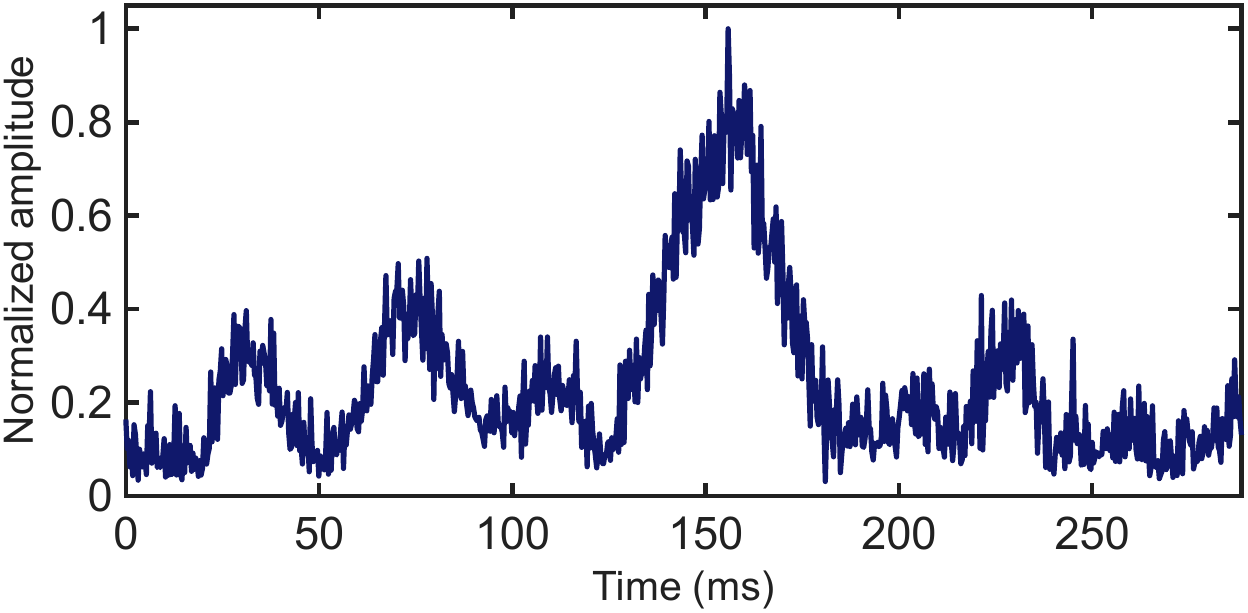}
    \end{minipage}

    \begin{minipage}{0.02\linewidth}
        \rotatebox{90}{\scriptsize Ling.}
    \end{minipage}\begin{minipage}{0.196\linewidth}
        \centering
        \includegraphics[width=0.99\linewidth]{pictures/timeseries/true_timeseries_lingual.pdf}
    \end{minipage}\begin{minipage}{0.2\linewidth}
        \centering
        \includegraphics[width=0.99\linewidth]{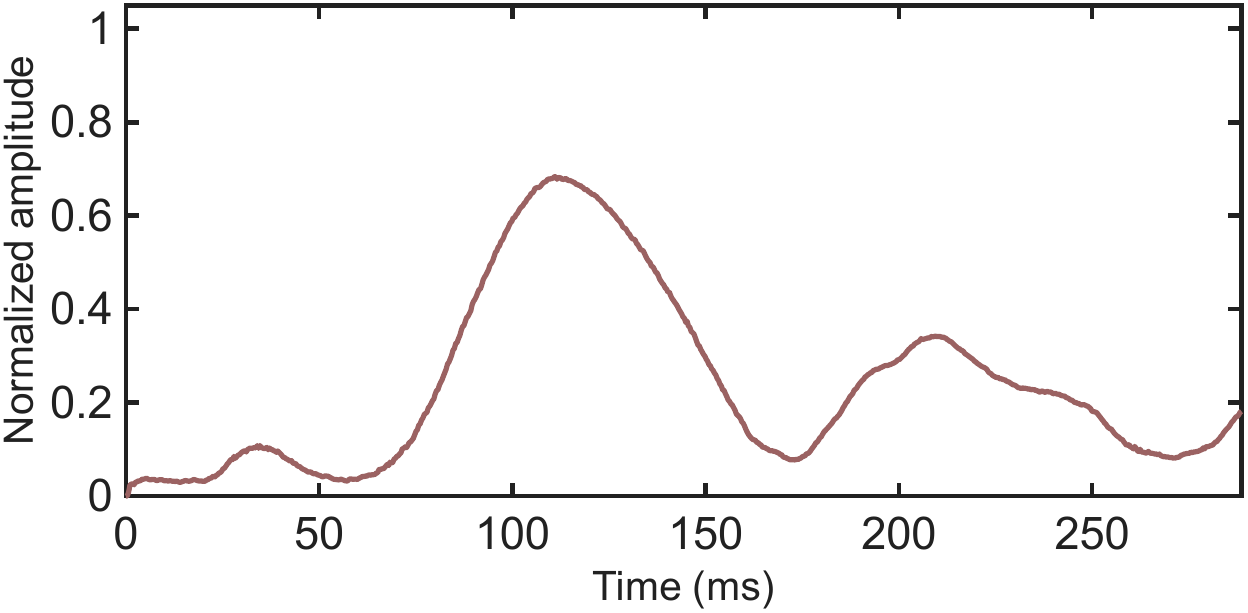}
    \end{minipage}\begin{minipage}{0.196\linewidth}
        \centering
        \includegraphics[width=0.99\linewidth]{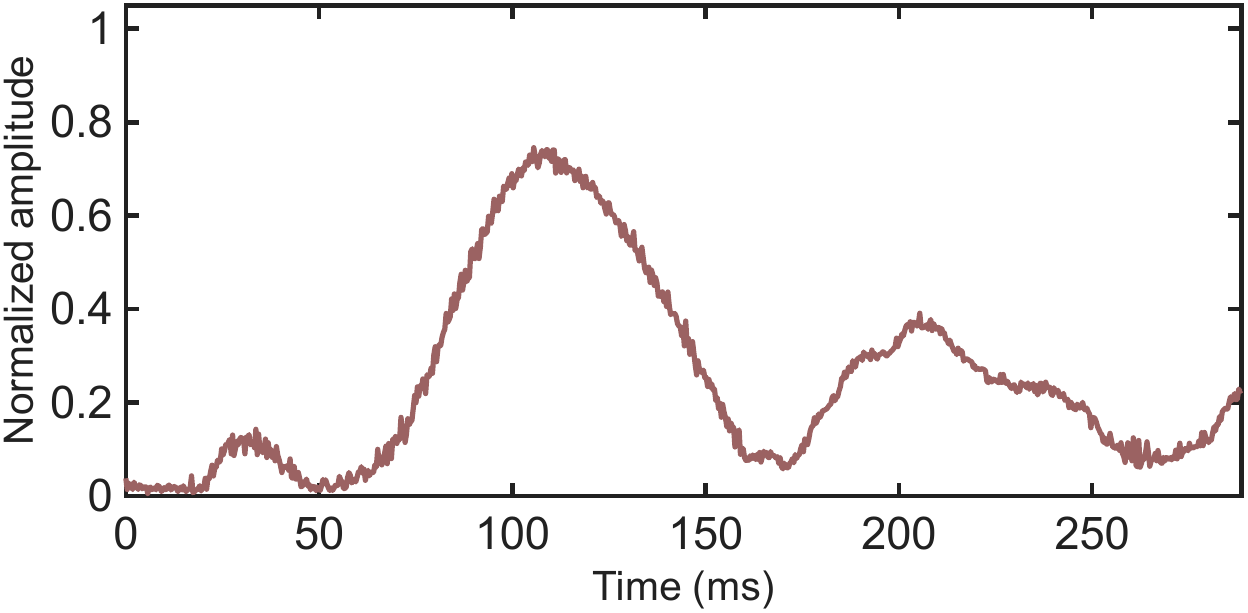}
    \end{minipage}\begin{minipage}{0.196\linewidth}
        \centering
        \includegraphics[width=0.99\linewidth]{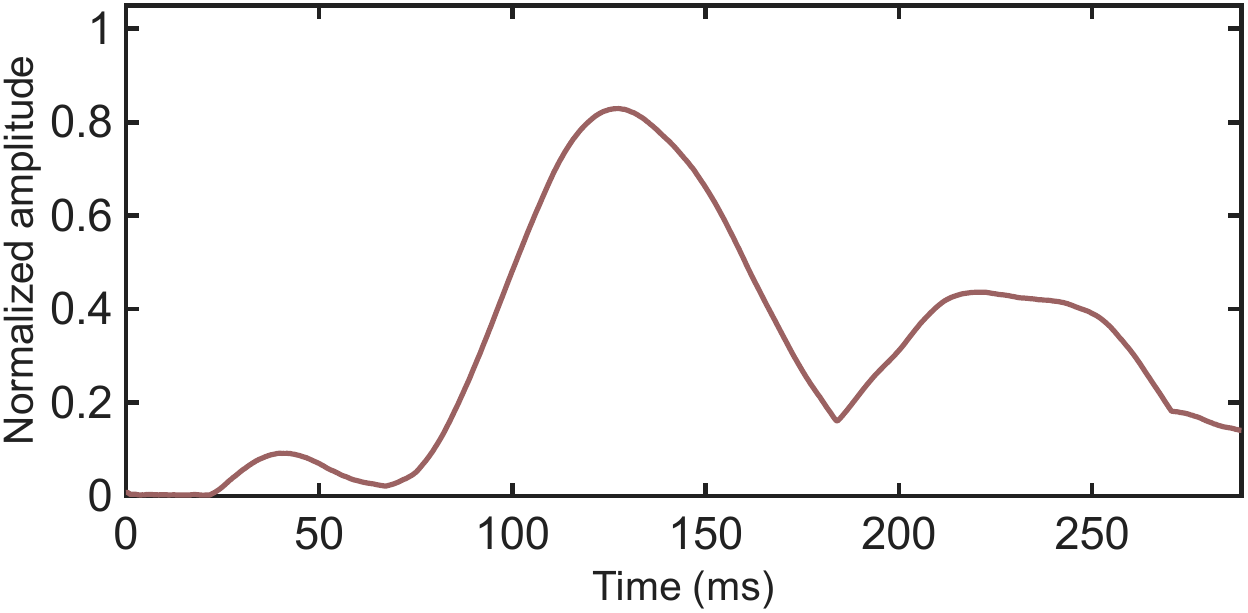}
    \end{minipage}\begin{minipage}{0.196\linewidth}
        \centering
        \includegraphics[width=0.99\linewidth]{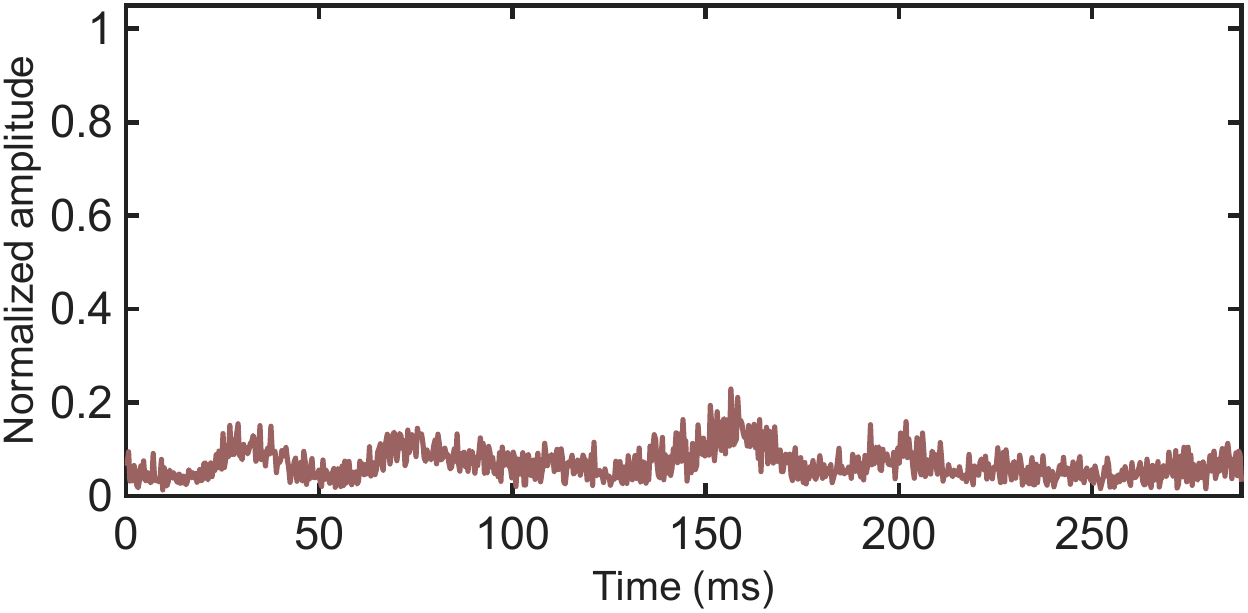}
    \end{minipage}

    \begin{minipage}{0.02\linewidth}
        \rotatebox{90}{\scriptsize S.P.}
    \end{minipage}\begin{minipage}{0.196\linewidth}
        \centering
        \includegraphics[width=0.99\linewidth]{pictures/timeseries/true_timeseries_superiorparietal.pdf}
    \end{minipage}\begin{minipage}{0.2\linewidth}
        \centering
        \includegraphics[width=0.99\linewidth]{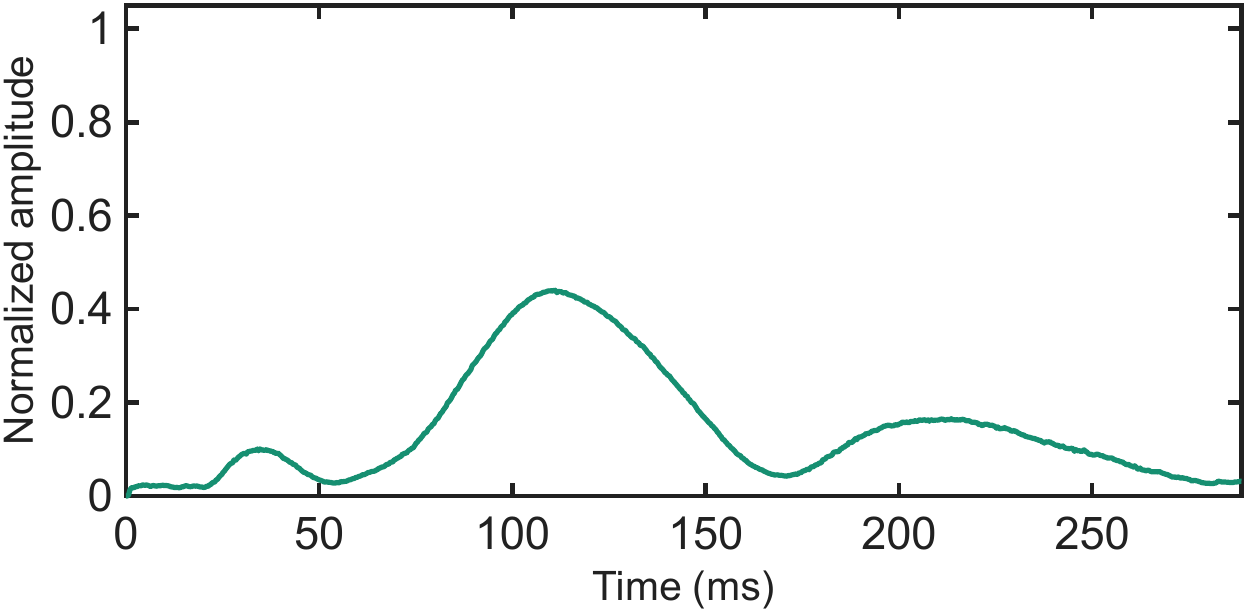}
    \end{minipage}\begin{minipage}{0.196\linewidth}
        \centering
        \includegraphics[width=0.99\linewidth]{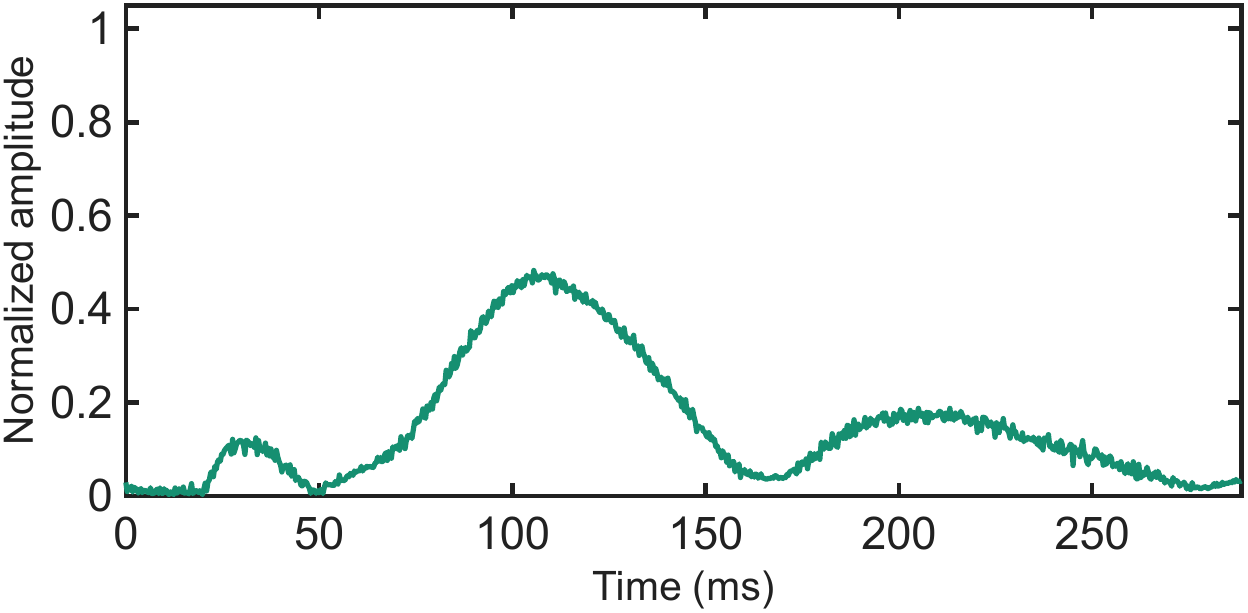}
    \end{minipage}\begin{minipage}{0.196\linewidth}
        \centering
        \includegraphics[width=0.99\linewidth]{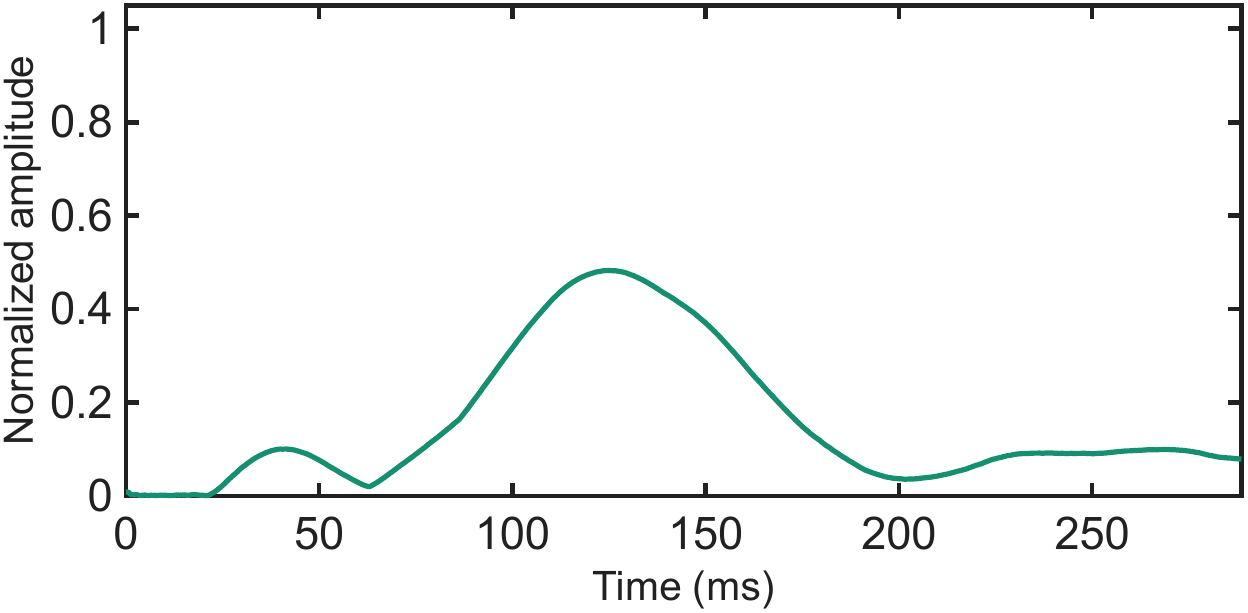}
    \end{minipage}\begin{minipage}{0.196\linewidth}
        \centering
        \includegraphics[width=0.99\linewidth]{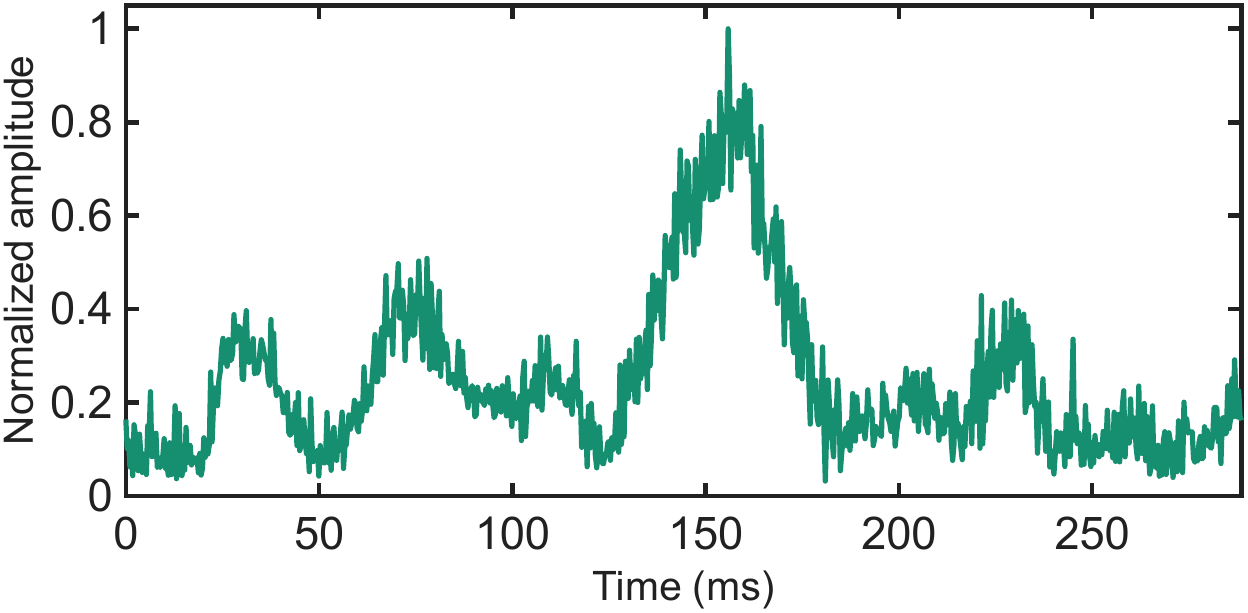}
    \end{minipage}

    \begin{minipage}{0.02\linewidth}
        \rotatebox{90}{\scriptsize C.M.F.}
    \end{minipage}\begin{minipage}{0.196\linewidth}
        \centering
        \includegraphics[width=0.99\linewidth]{pictures/timeseries/true_timeseries_caudalmiddlefrontal.pdf}
    \end{minipage}\begin{minipage}{0.2\linewidth}
        \centering
        \includegraphics[width=0.99\linewidth]{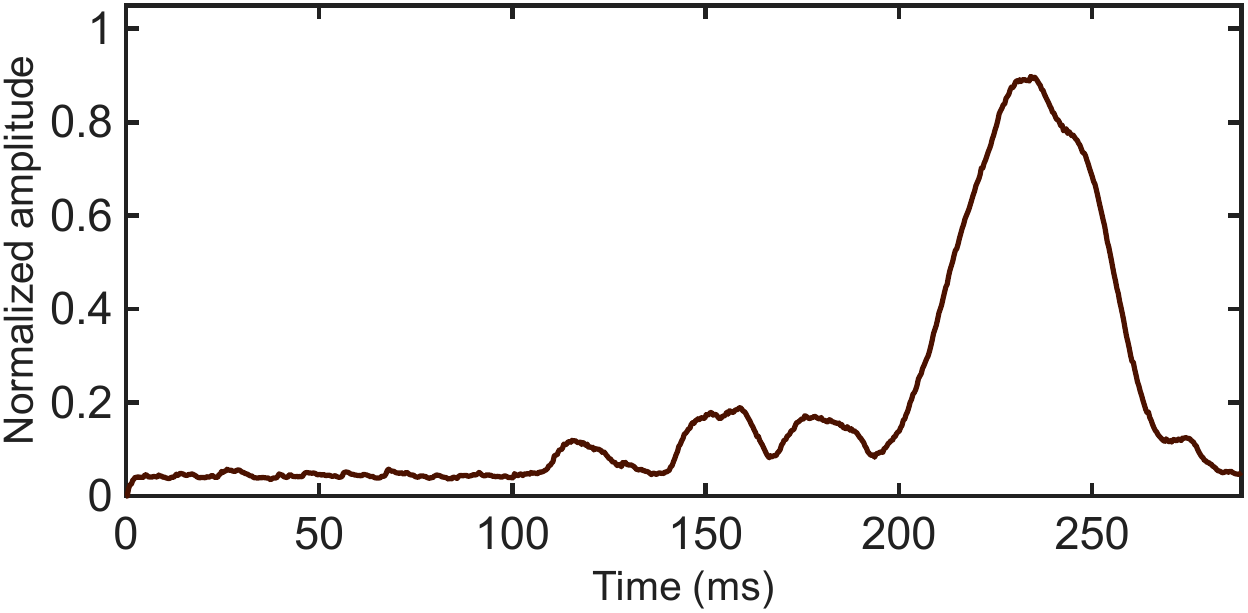}
    \end{minipage}\begin{minipage}{0.196\linewidth}
        \centering
        \includegraphics[width=0.99\linewidth]{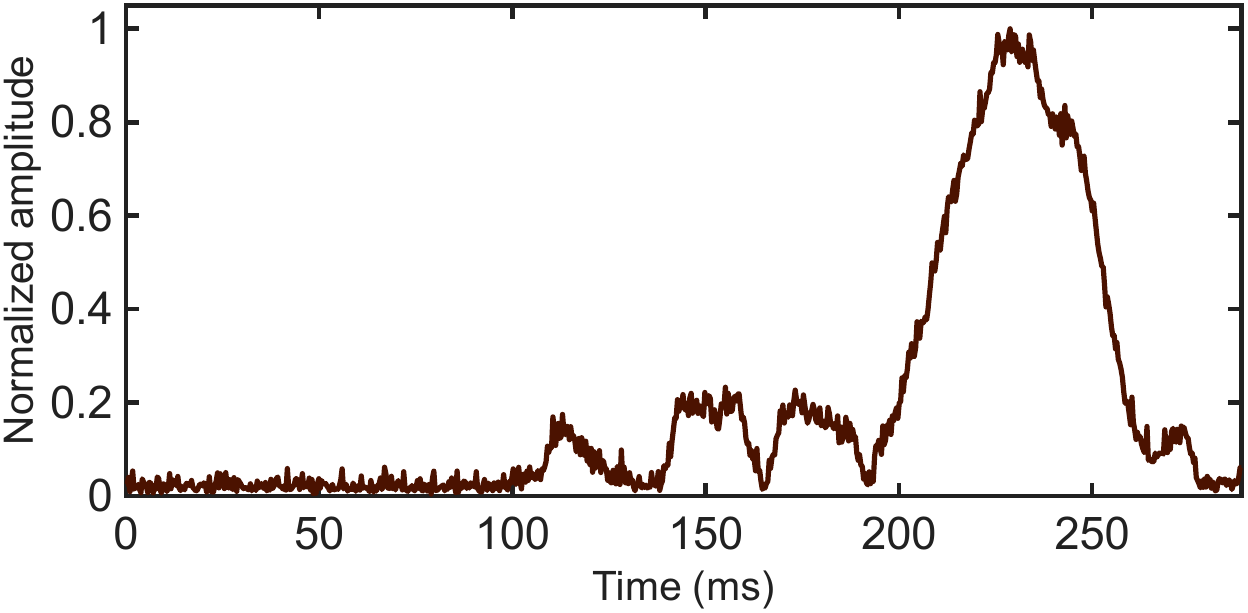}
    \end{minipage}\begin{minipage}{0.196\linewidth}
        \centering
        \includegraphics[width=0.99\linewidth]{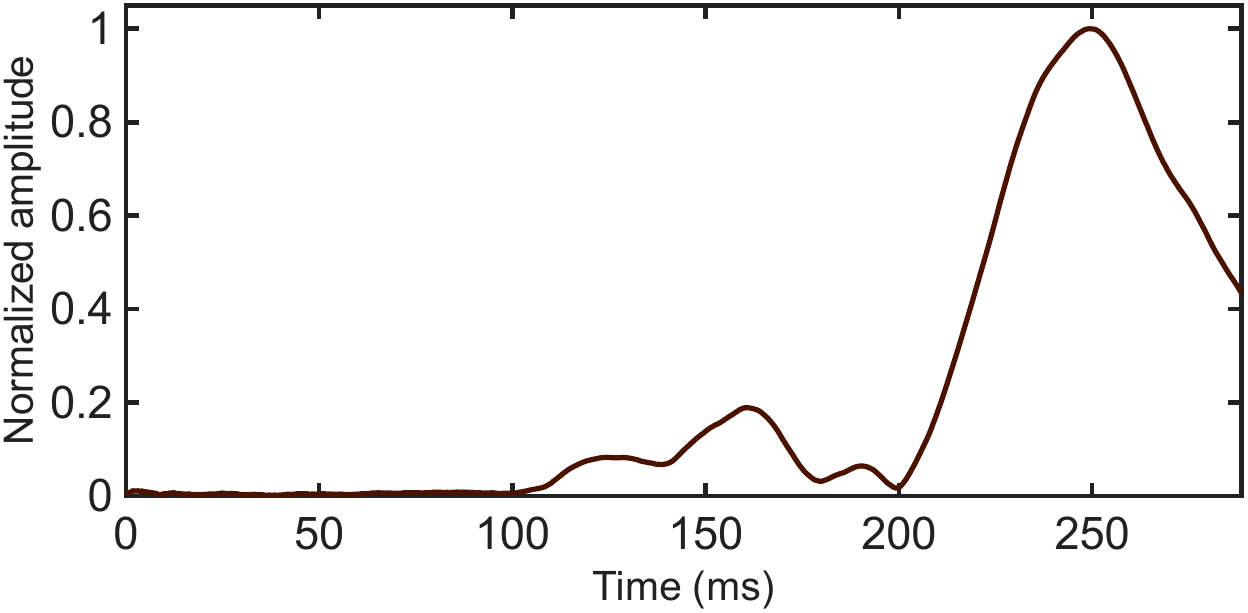}
    \end{minipage}\begin{minipage}{0.196\linewidth}
        \centering
        \includegraphics[width=0.99\linewidth]{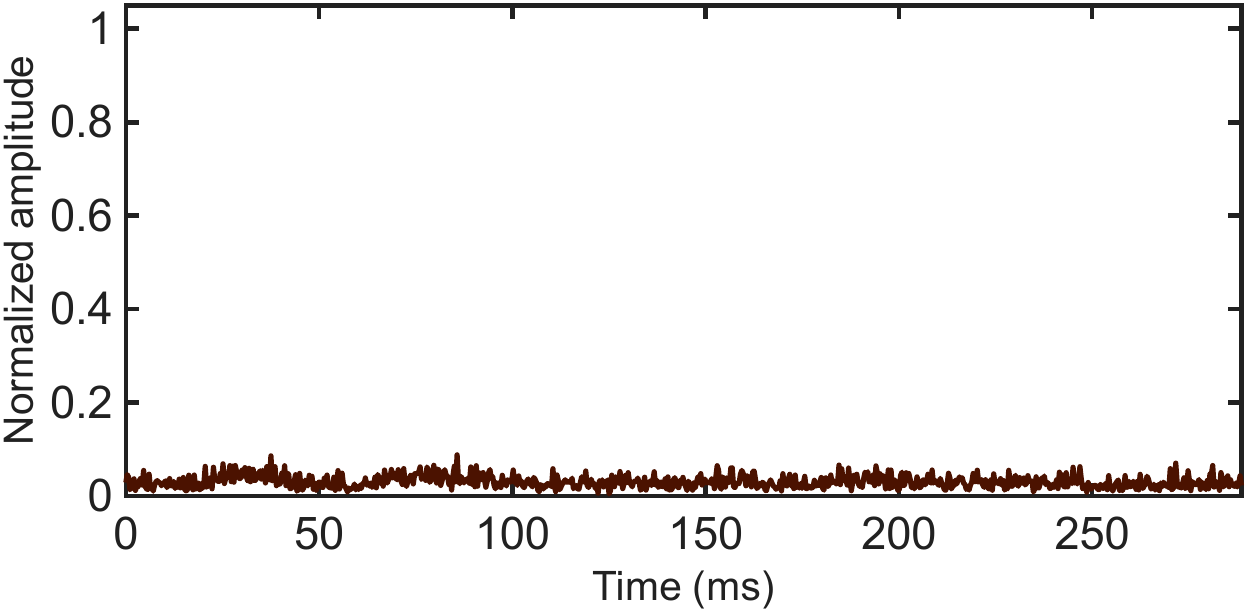}
    \end{minipage}

    \begin{minipage}{0.02\linewidth}
        \rotatebox{90}{\scriptsize I.P.}
    \end{minipage}\begin{minipage}{0.196\linewidth}
        \centering
        \hfill
        \end{minipage}\begin{minipage}{0.2\linewidth}
        \centering
        \includegraphics[width=0.99\linewidth]{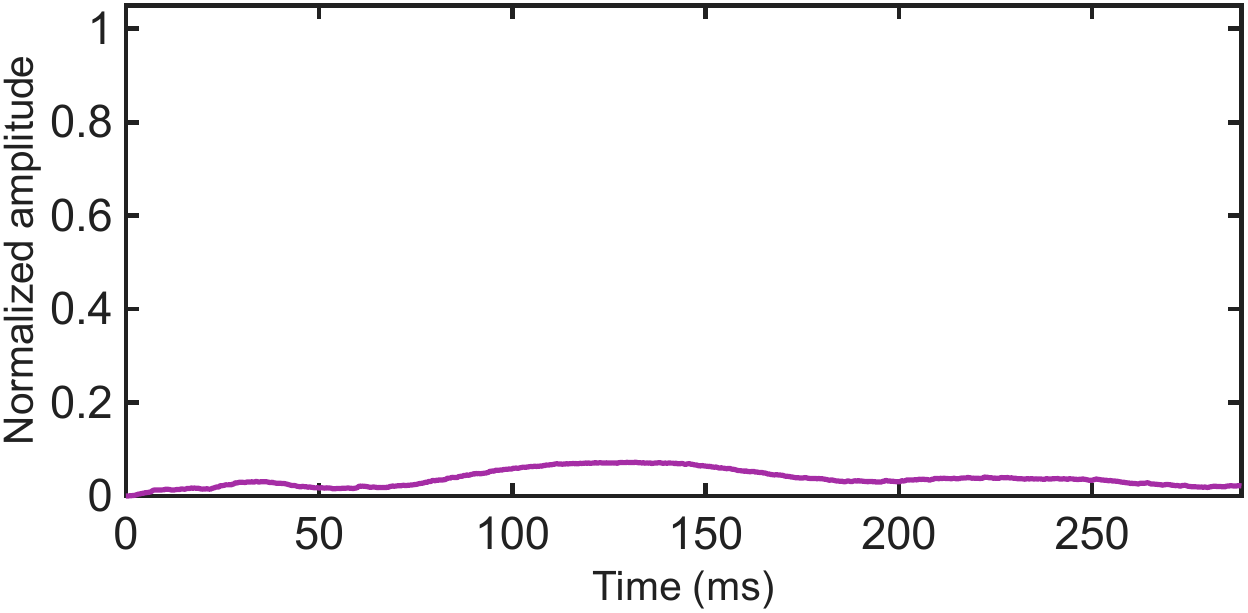}
    \end{minipage}\begin{minipage}{0.196\linewidth}
        \centering
        \includegraphics[width=0.99\linewidth]{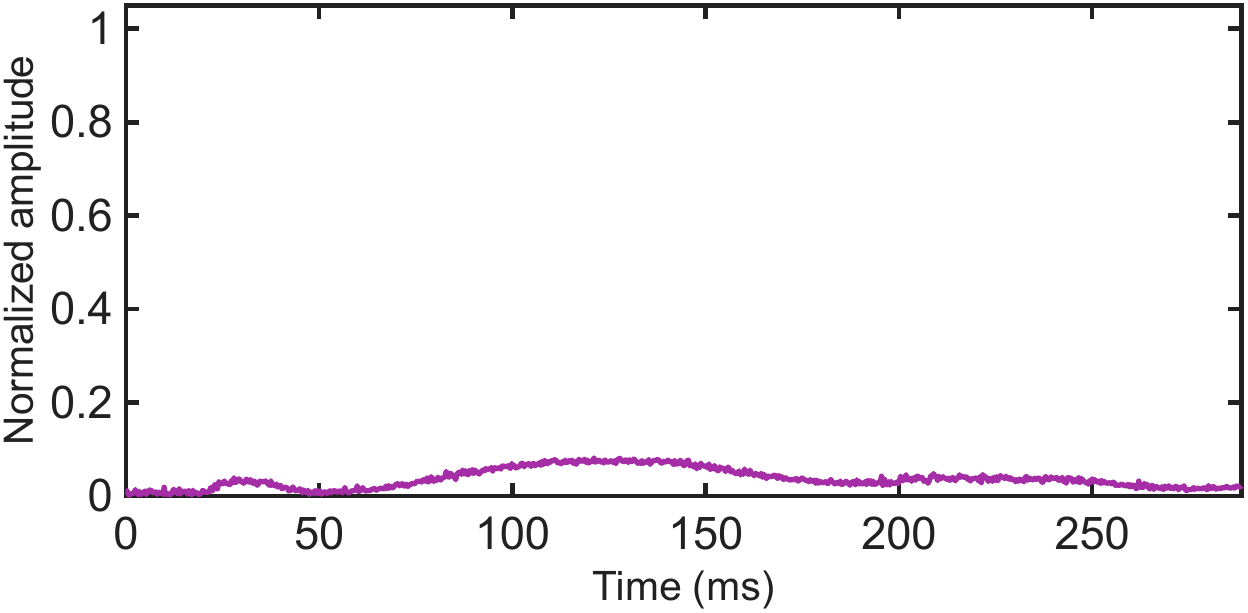}
    \end{minipage}\begin{minipage}{0.196\linewidth}
        \centering
        \includegraphics[width=0.99\linewidth]{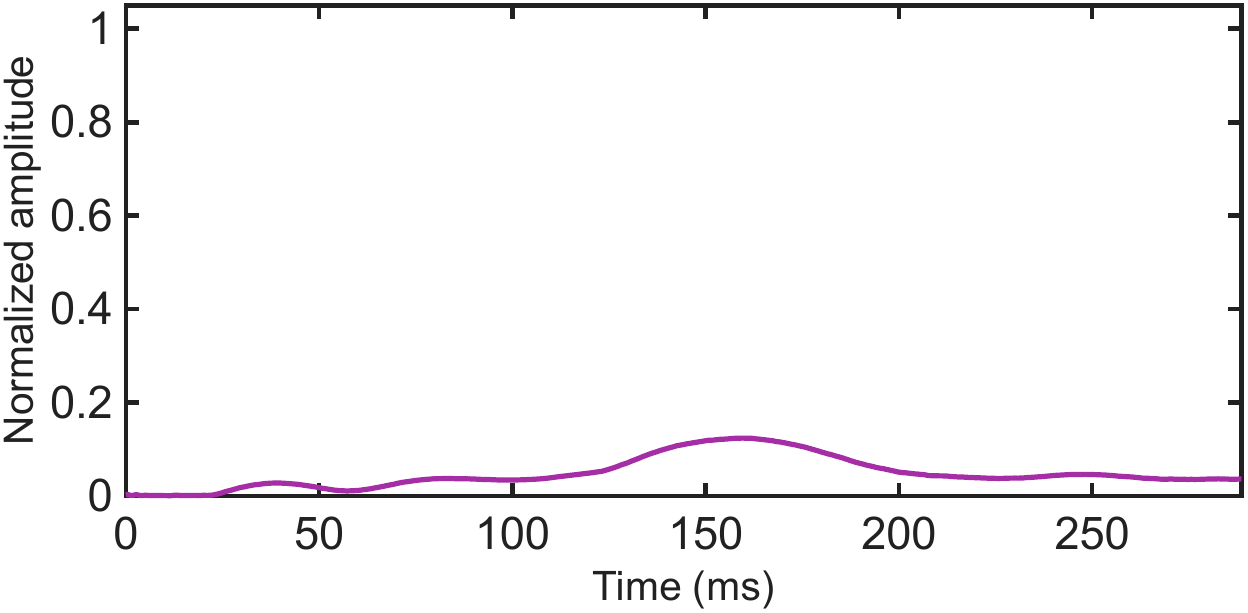}
    \end{minipage}\begin{minipage}{0.196\linewidth}
        \centering
        \includegraphics[width=0.99\linewidth]{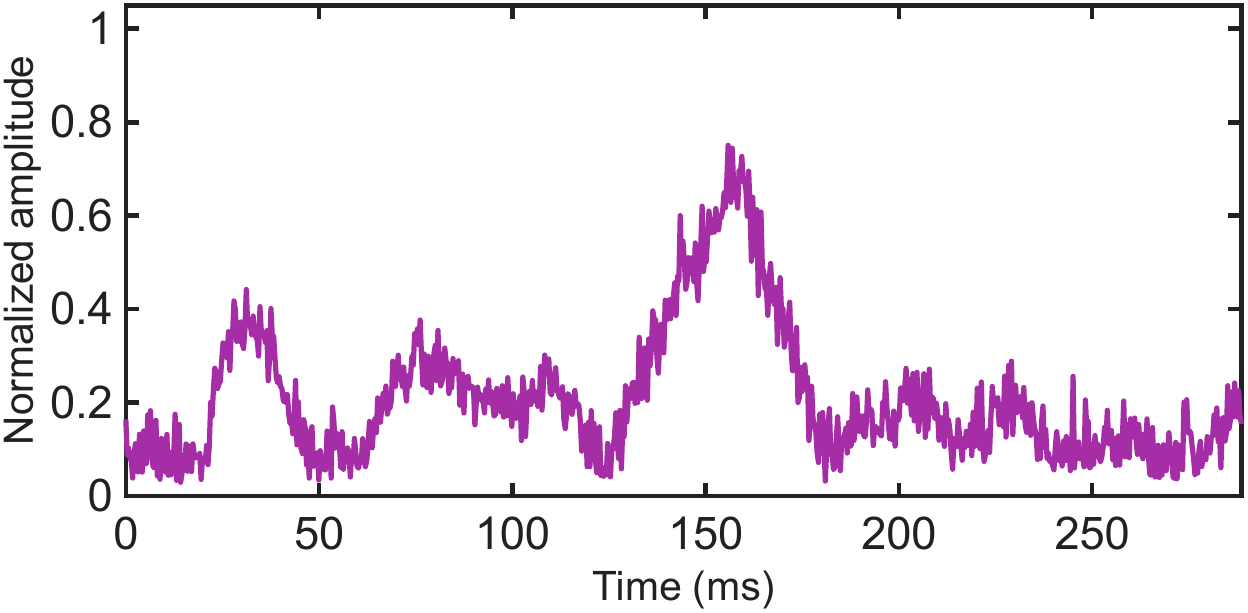}
    \end{minipage}

    \begin{minipage}{0.02\linewidth}
        \rotatebox{90}{\scriptsize Fusif.}
    \end{minipage}\begin{minipage}{0.196\linewidth}
        \centering
        \hfill
    \end{minipage}\begin{minipage}{0.2\linewidth}
        \centering
        \includegraphics[width=0.99\linewidth]{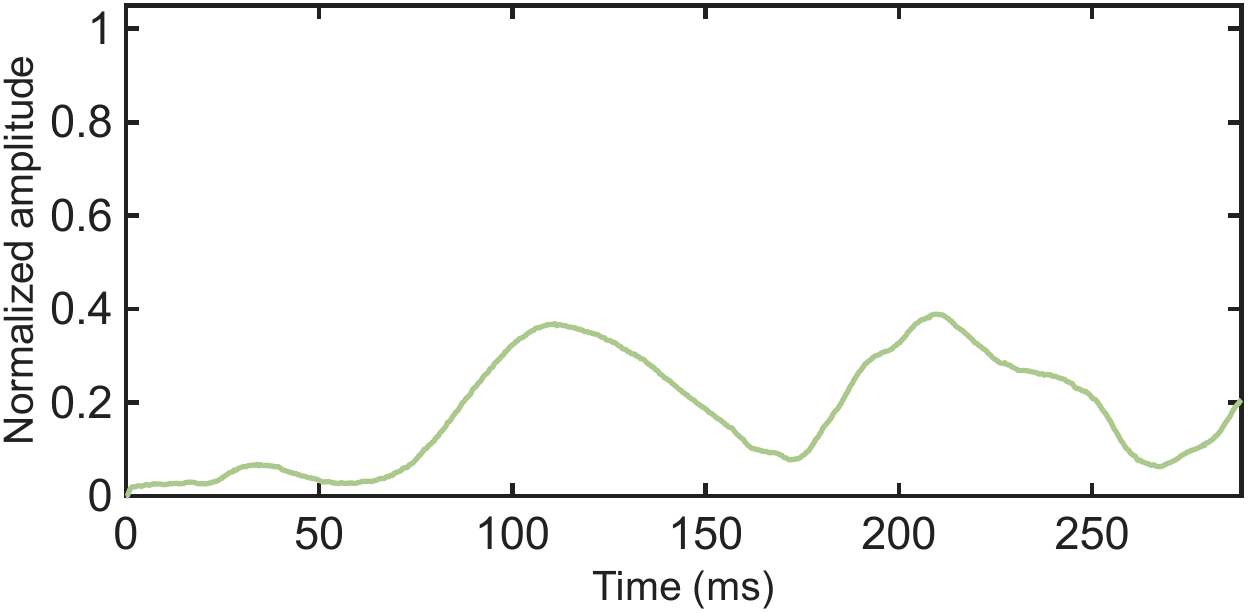}
    \end{minipage}\begin{minipage}{0.196\linewidth}
        \centering
        \includegraphics[width=0.99\linewidth]{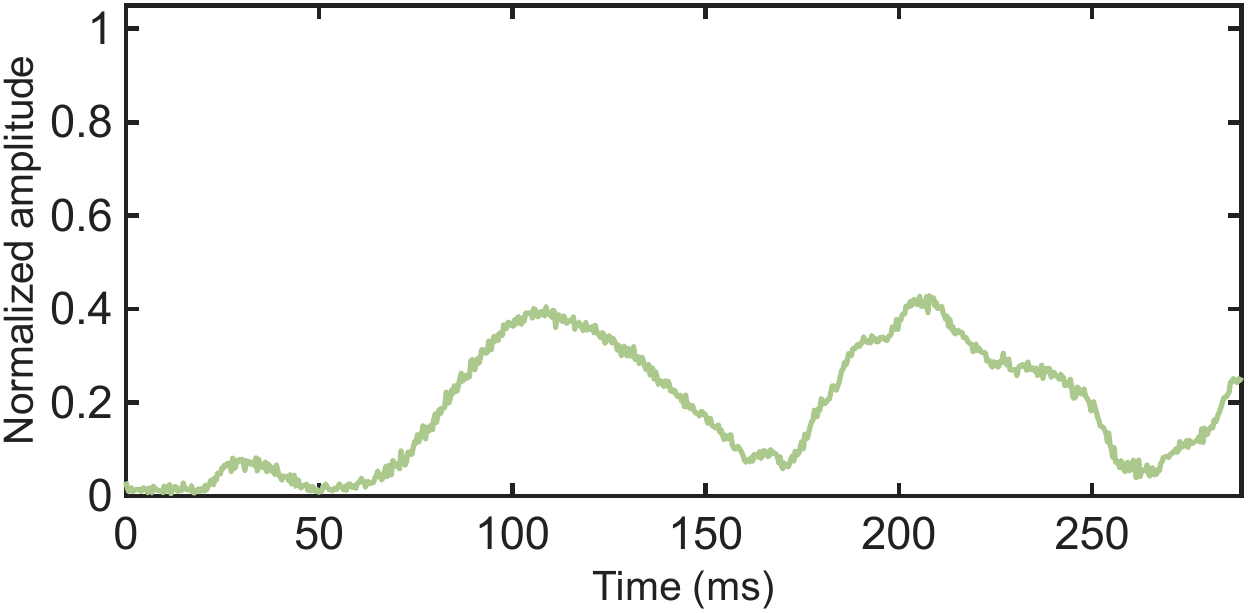}
    \end{minipage}\begin{minipage}{0.196\linewidth}
        \centering
        \includegraphics[width=0.99\linewidth]{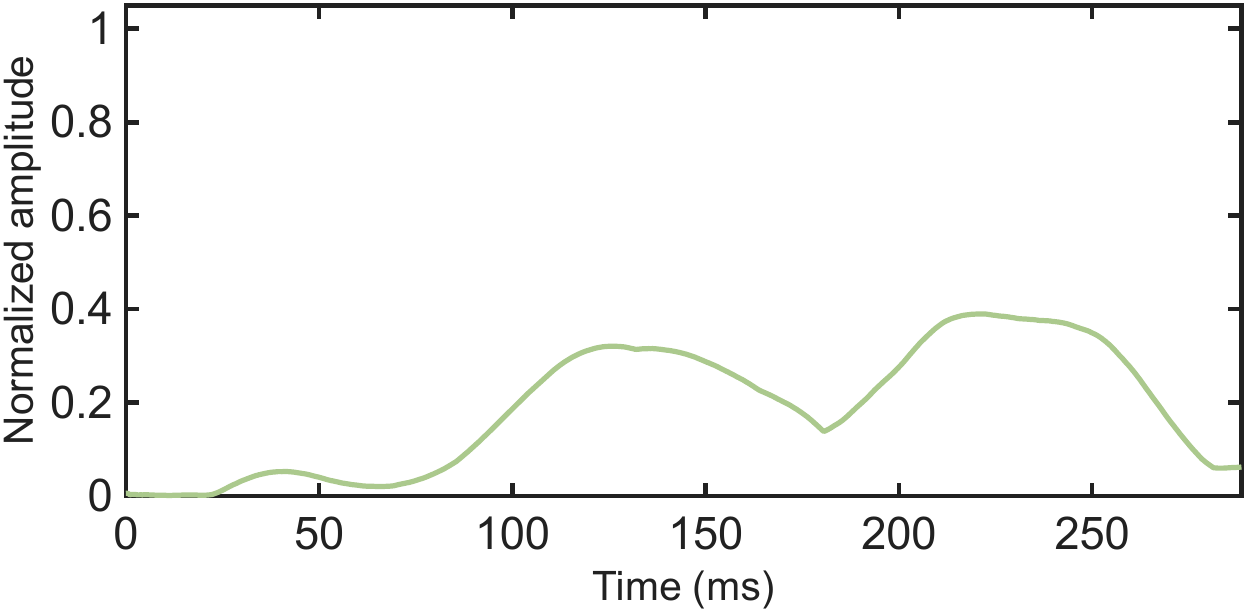}
    \end{minipage}\begin{minipage}{0.196\linewidth}
        \centering
        \includegraphics[width=0.99\linewidth]{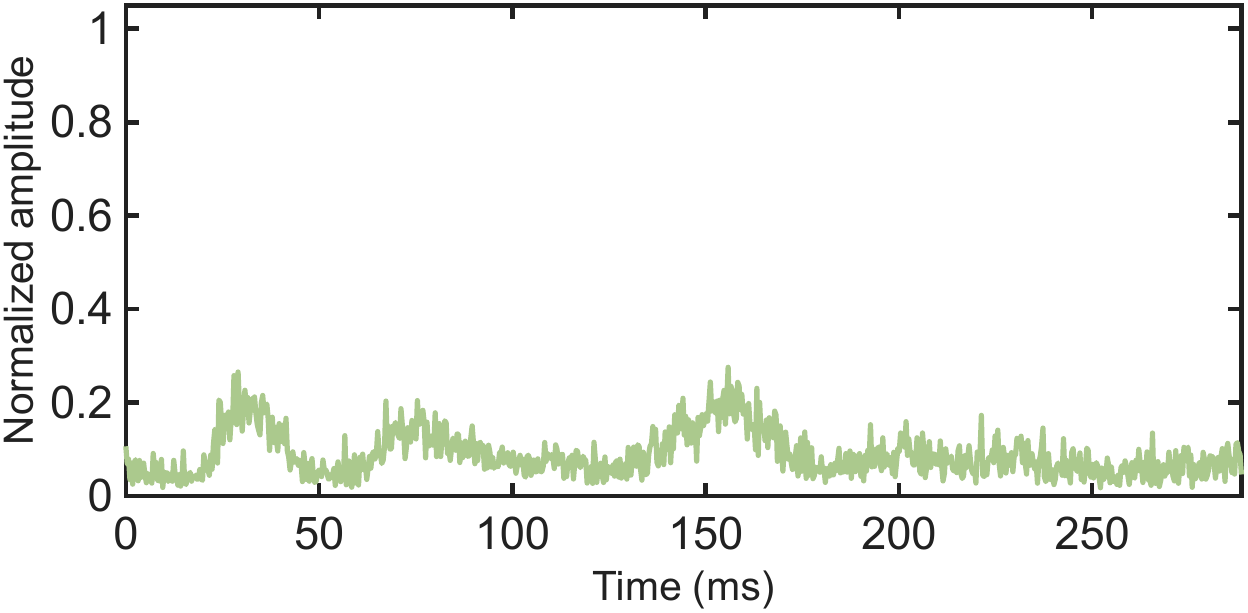}
    \end{minipage}
    \caption{Parcel-wise temporal activity at 25 dB SNR. Rows correspond to the selected cortical parcels, and columns show the simulated ground-truth time series followed by CrSKF, SKF, KF, and LCMV reconstructions.}
    \label{fig:timeseries25dB}
\end{figure}

\begin{figure}[h!]
    \centering
    \begin{minipage}{0.02\linewidth}
        \rotatebox{90}{\scriptsize \hfill}
    \end{minipage}\begin{minipage}{0.196\linewidth}
        \centering
        {\bf Ground truth}
    \end{minipage}\begin{minipage}{0.196\linewidth}
        \centering
        {\bf CrSKF}
    \end{minipage}\begin{minipage}{0.196\linewidth}
        \centering
        {\bf SKF}
    \end{minipage}\begin{minipage}{0.196\linewidth}
        \centering
        {\bf KF}
    \end{minipage}\begin{minipage}{0.196\linewidth}
        \centering
        {\bf LCMV}
    \end{minipage}

    \begin{minipage}{0.02\linewidth}
        \rotatebox{90}{\scriptsize L.O.}
    \end{minipage}\begin{minipage}{0.196\linewidth}
        \centering
        \includegraphics[width=0.99\linewidth]{pictures/timeseries/true_timeseries_lateraloccipital.pdf}
    \end{minipage}\begin{minipage}{0.2\linewidth}
        \centering
        \includegraphics[width=0.99\linewidth]{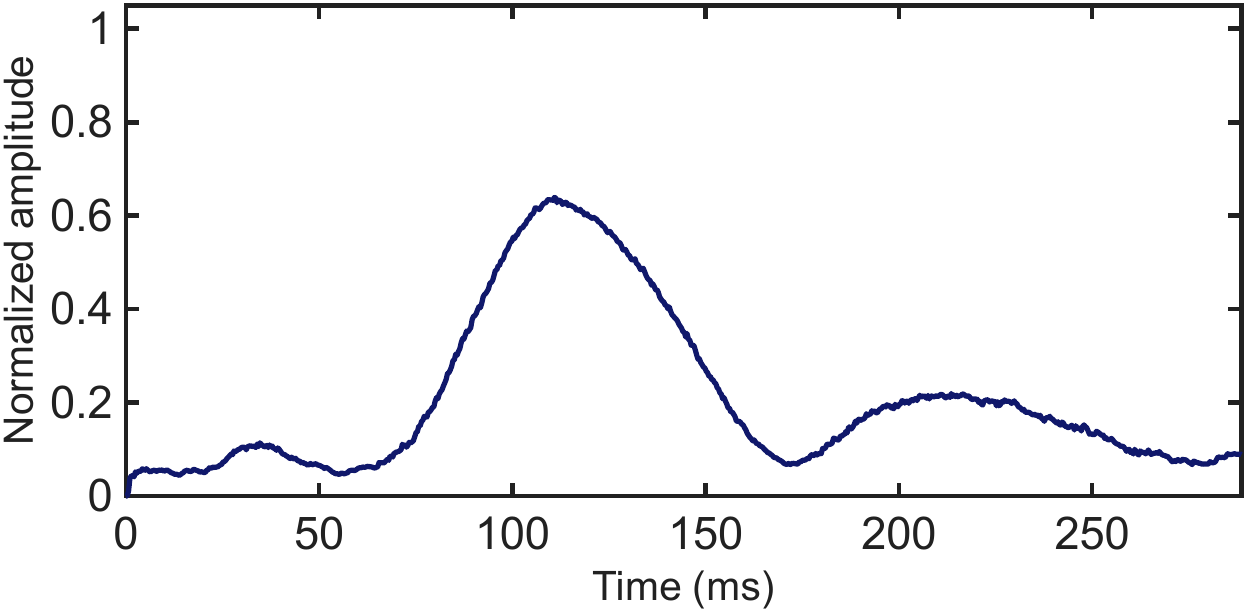}
    \end{minipage}\begin{minipage}{0.196\linewidth}
        \centering
        \includegraphics[width=0.99\linewidth]{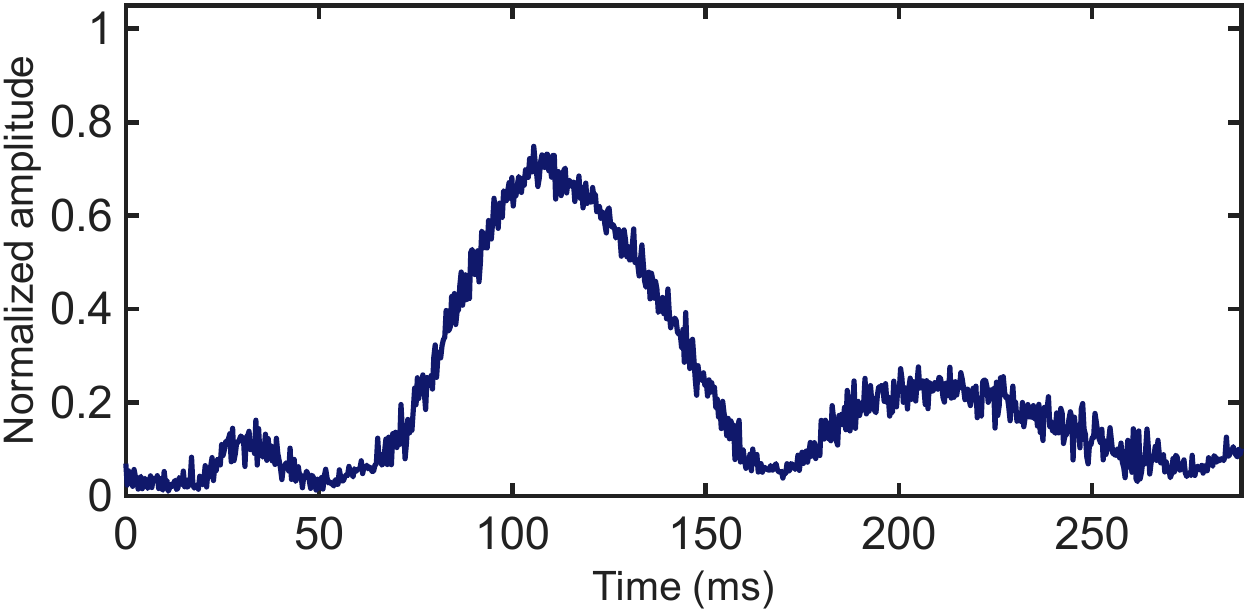}
    \end{minipage}\begin{minipage}{0.196\linewidth}
        \centering
        \includegraphics[width=0.99\linewidth]{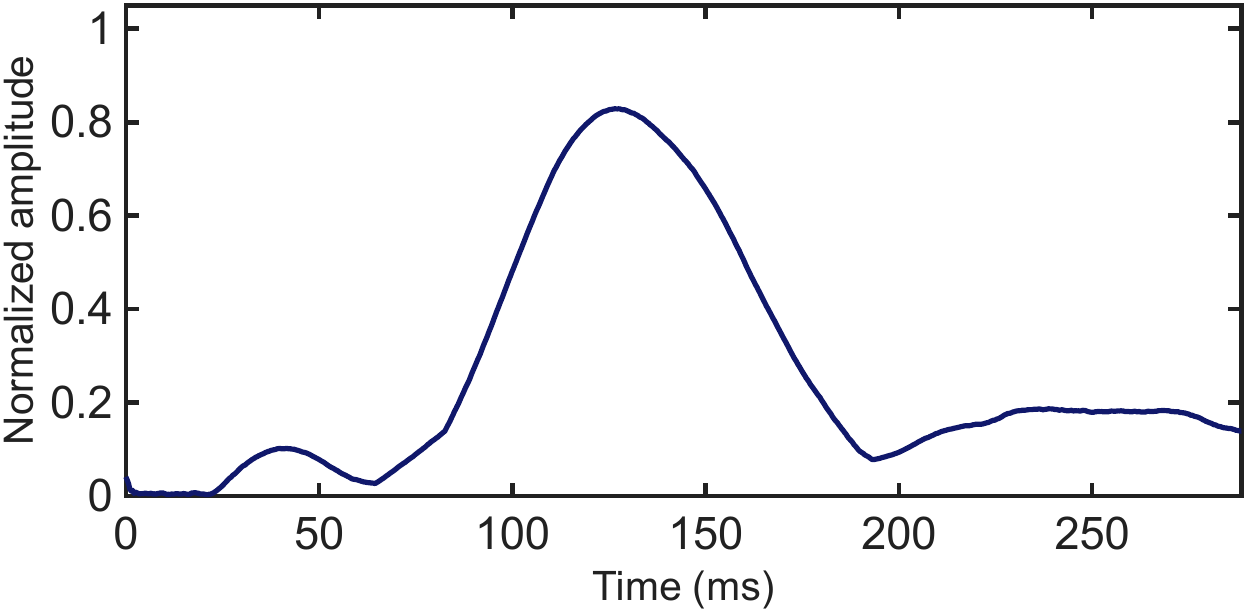}
    \end{minipage}\begin{minipage}{0.196\linewidth}
        \centering
        \includegraphics[width=0.99\linewidth]{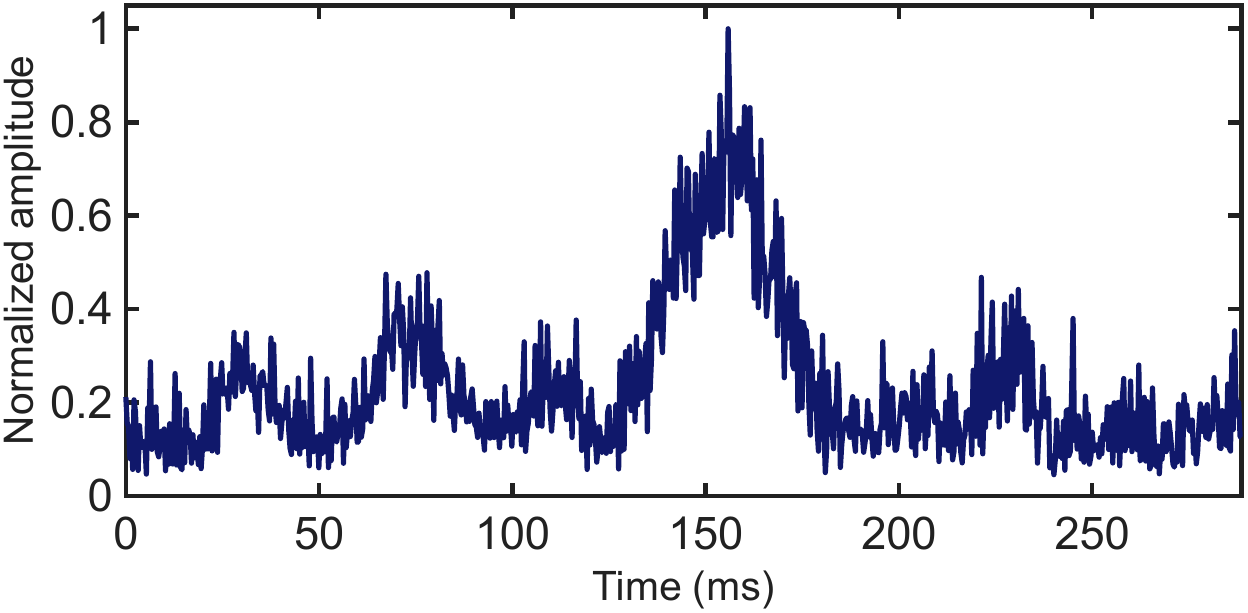}
    \end{minipage}

    \begin{minipage}{0.02\linewidth}
        \rotatebox{90}{\scriptsize Ling.}
    \end{minipage}\begin{minipage}{0.196\linewidth}
        \centering
        \includegraphics[width=0.99\linewidth]{pictures/timeseries/true_timeseries_lingual.pdf}
    \end{minipage}\begin{minipage}{0.2\linewidth}
        \centering
        \includegraphics[width=0.99\linewidth]{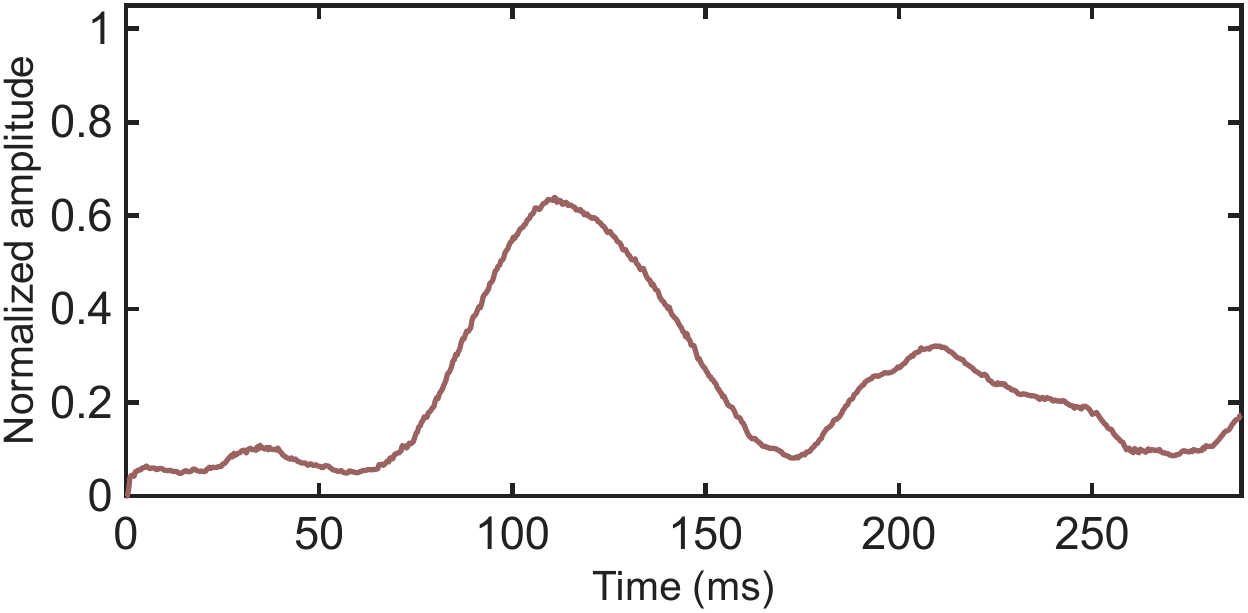}
    \end{minipage}\begin{minipage}{0.196\linewidth}
        \centering
        \includegraphics[width=0.99\linewidth]{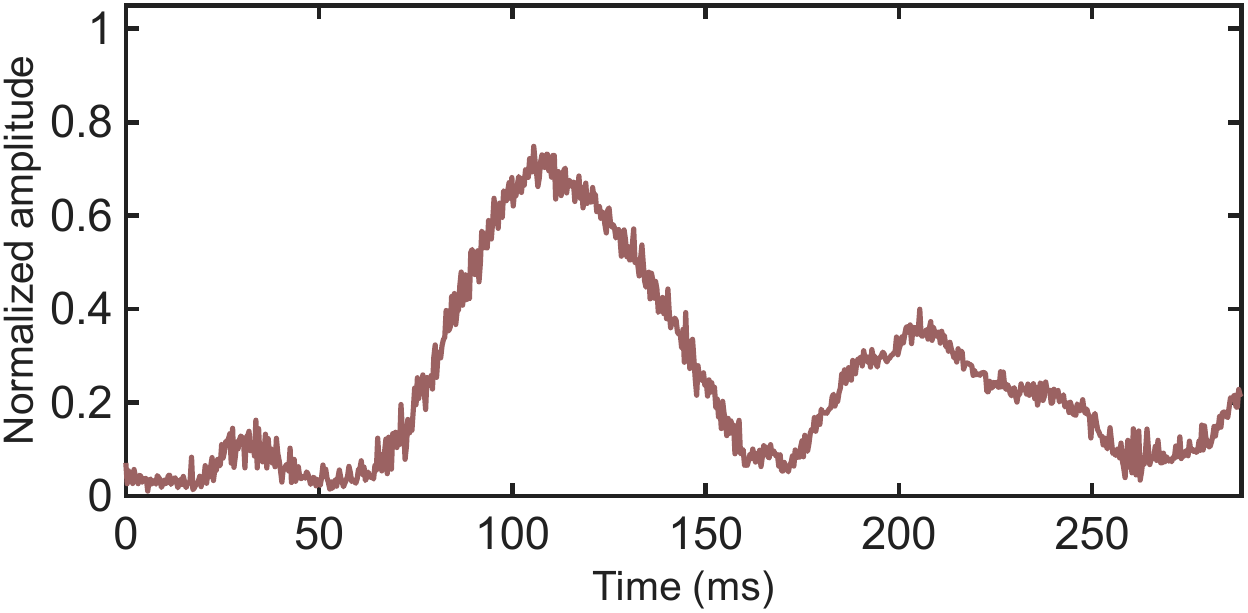}
    \end{minipage}\begin{minipage}{0.196\linewidth}
        \centering
        \includegraphics[width=0.99\linewidth]{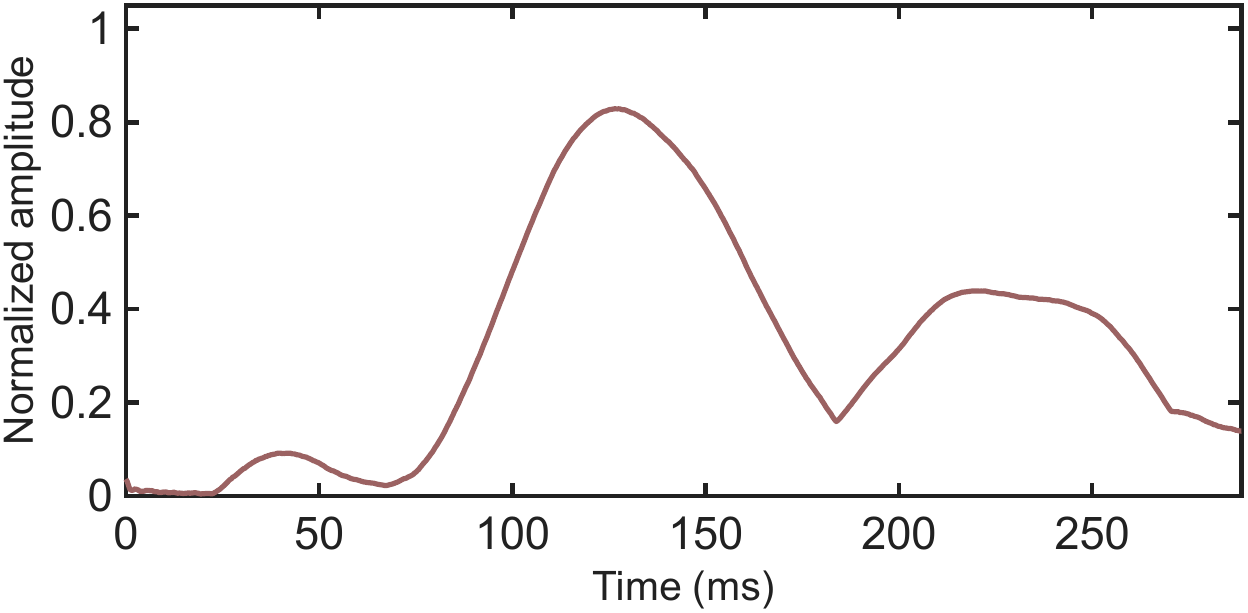}
    \end{minipage}\begin{minipage}{0.196\linewidth}
        \centering
        \includegraphics[width=0.99\linewidth]{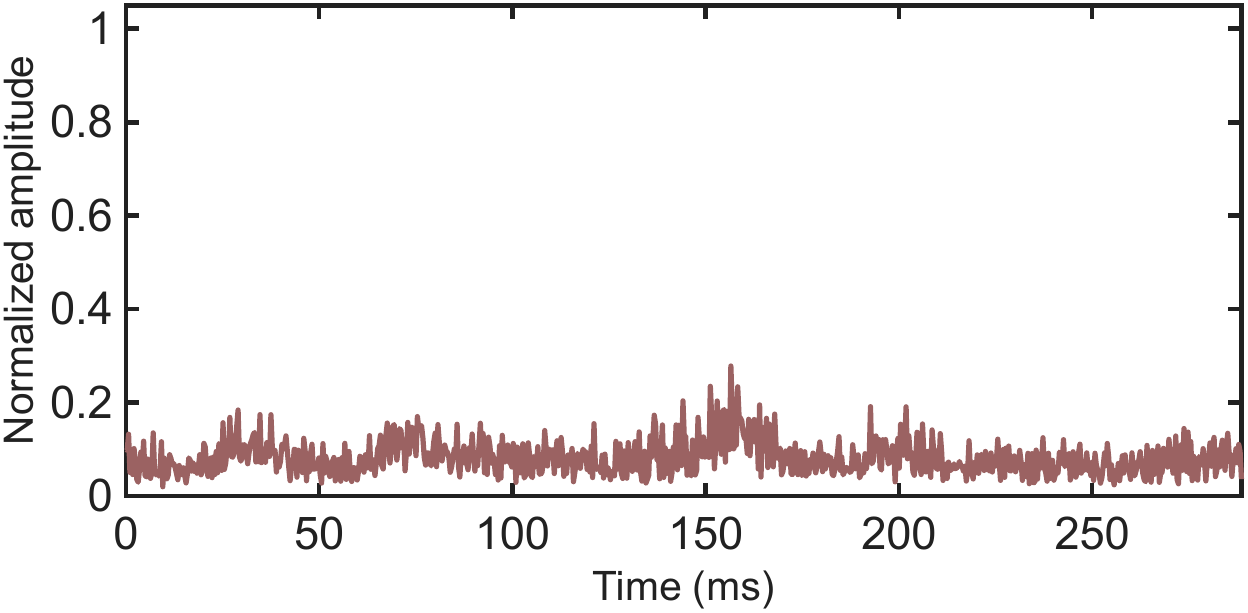}
    \end{minipage}

    \begin{minipage}{0.02\linewidth}
        \rotatebox{90}{\scriptsize S.P.}
    \end{minipage}\begin{minipage}{0.196\linewidth}
        \centering
        \includegraphics[width=0.99\linewidth]{pictures/timeseries/true_timeseries_superiorparietal.pdf}
    \end{minipage}\begin{minipage}{0.2\linewidth}
        \centering
        \includegraphics[width=0.99\linewidth]{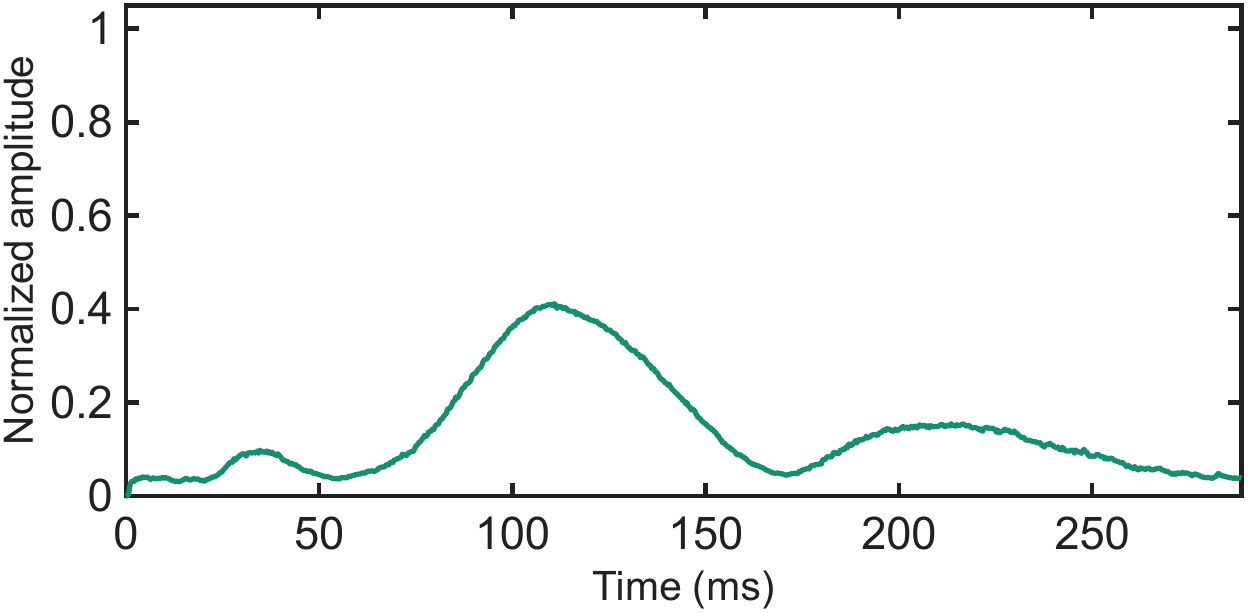}
    \end{minipage}\begin{minipage}{0.196\linewidth}
        \centering
        \includegraphics[width=0.99\linewidth]{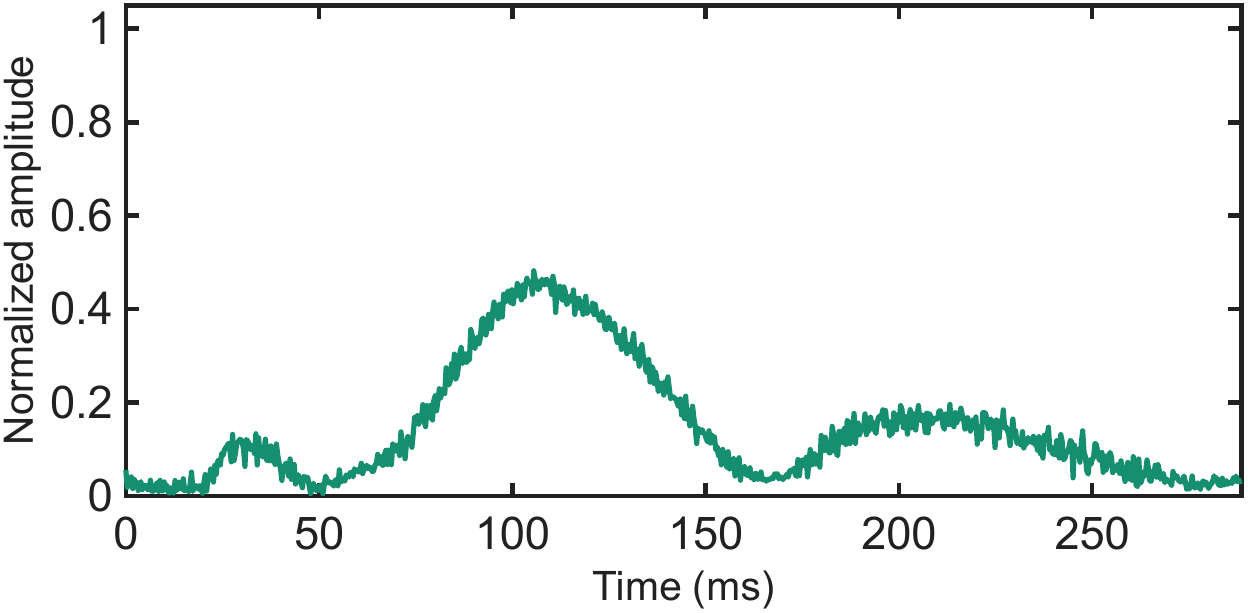}
    \end{minipage}\begin{minipage}{0.196\linewidth}
        \centering
        \includegraphics[width=0.99\linewidth]{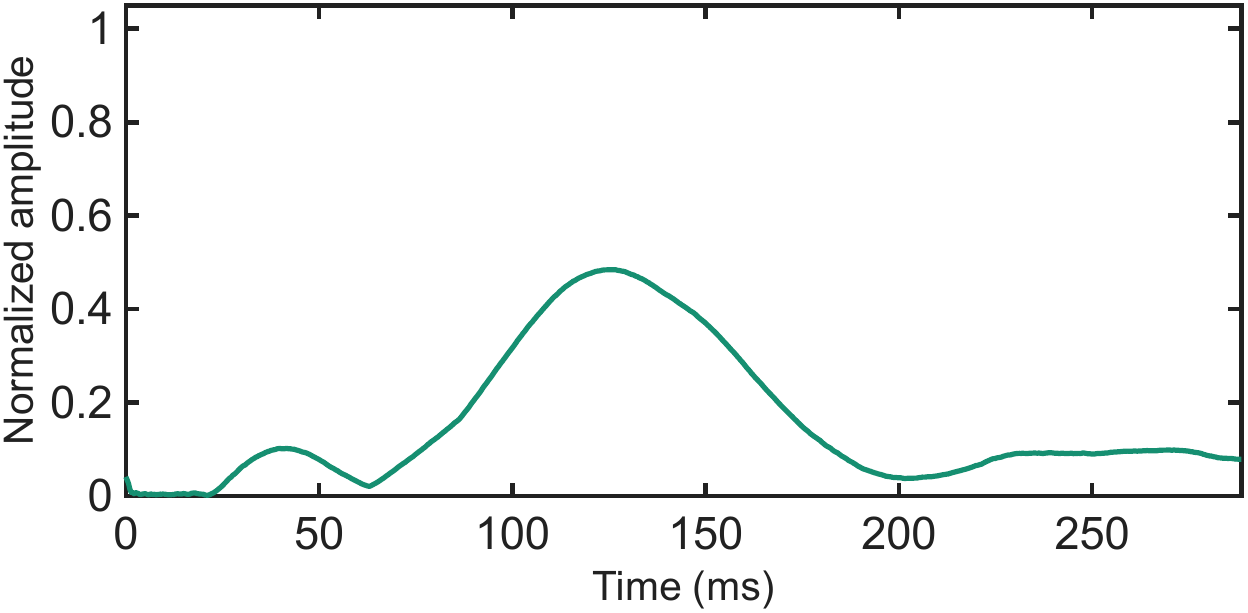}
    \end{minipage}\begin{minipage}{0.196\linewidth}
        \centering
        \includegraphics[width=0.99\linewidth]{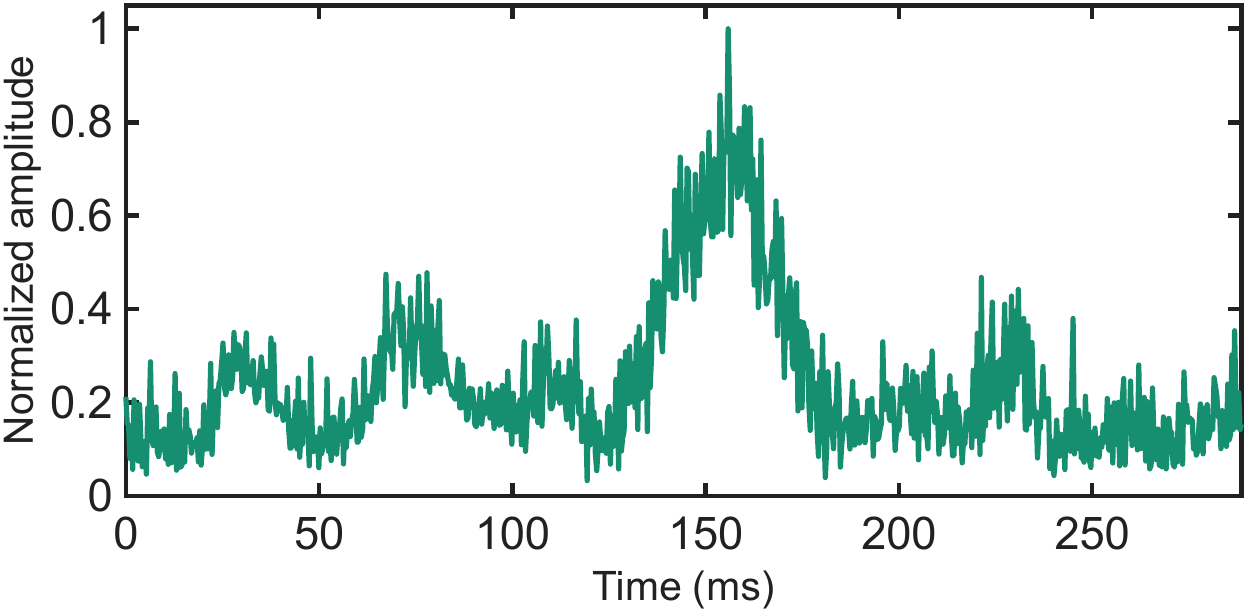}
    \end{minipage}

    \begin{minipage}{0.02\linewidth}
        \rotatebox{90}{\scriptsize C.M.F.}
    \end{minipage}\begin{minipage}{0.196\linewidth}
        \centering
        \includegraphics[width=0.99\linewidth]{pictures/timeseries/true_timeseries_caudalmiddlefrontal.pdf}
    \end{minipage}\begin{minipage}{0.2\linewidth}
        \centering
        \includegraphics[width=0.99\linewidth]{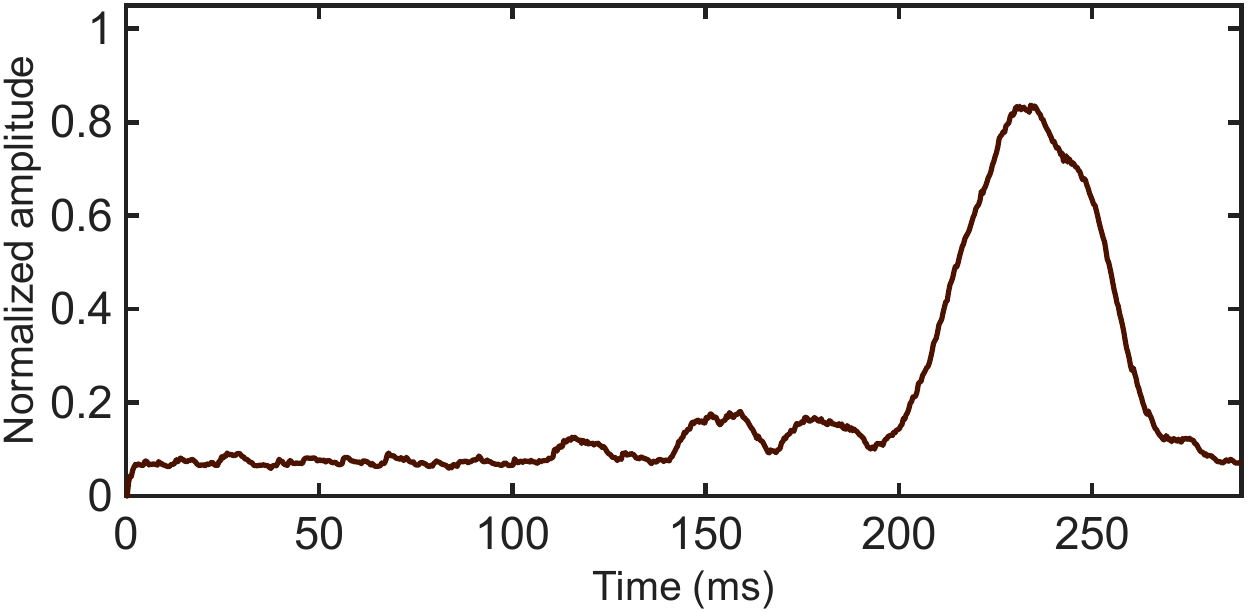}
    \end{minipage}\begin{minipage}{0.196\linewidth}
        \centering
        \includegraphics[width=0.99\linewidth]{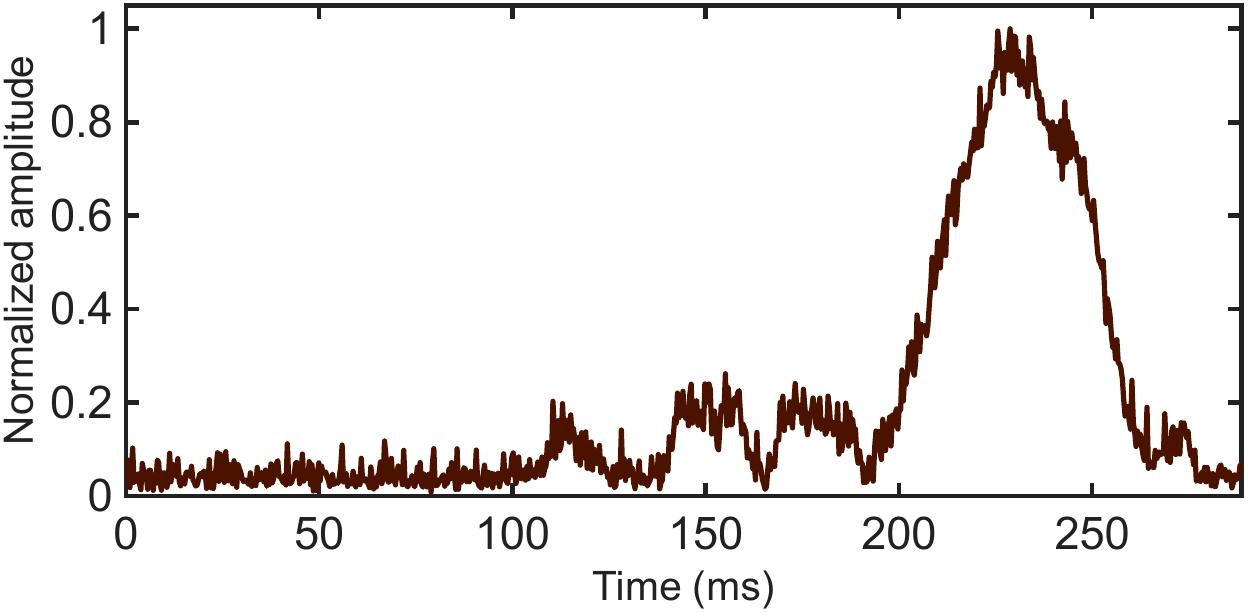}
    \end{minipage}\begin{minipage}{0.196\linewidth}
        \centering
        \includegraphics[width=0.99\linewidth]{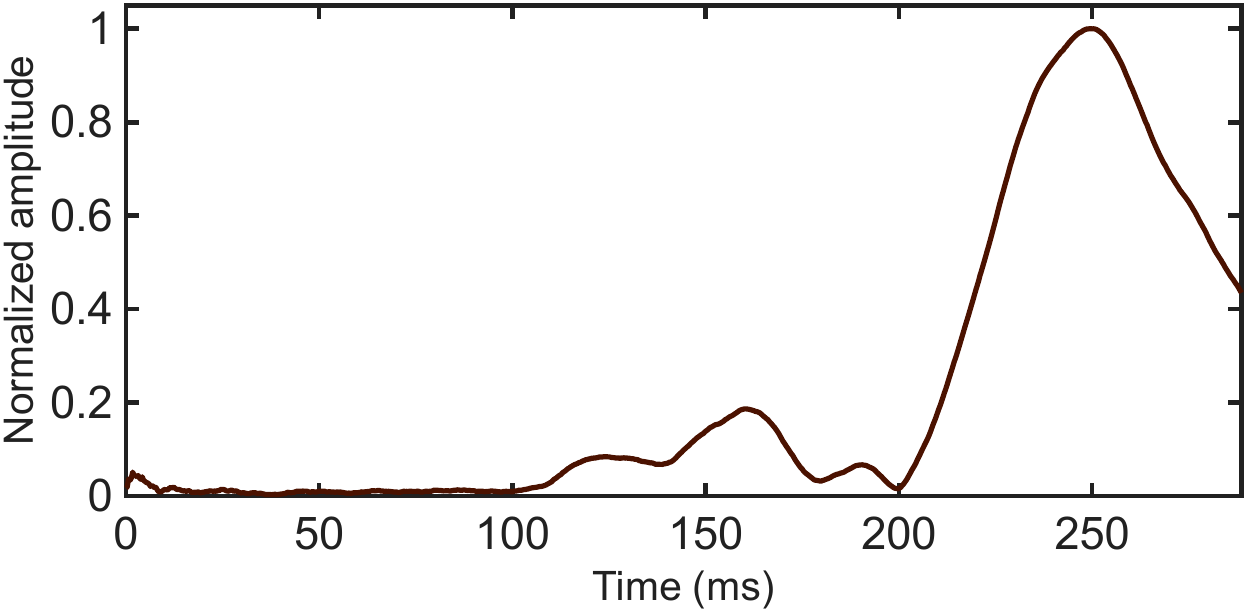}
    \end{minipage}\begin{minipage}{0.196\linewidth}
        \centering
        \includegraphics[width=0.99\linewidth]{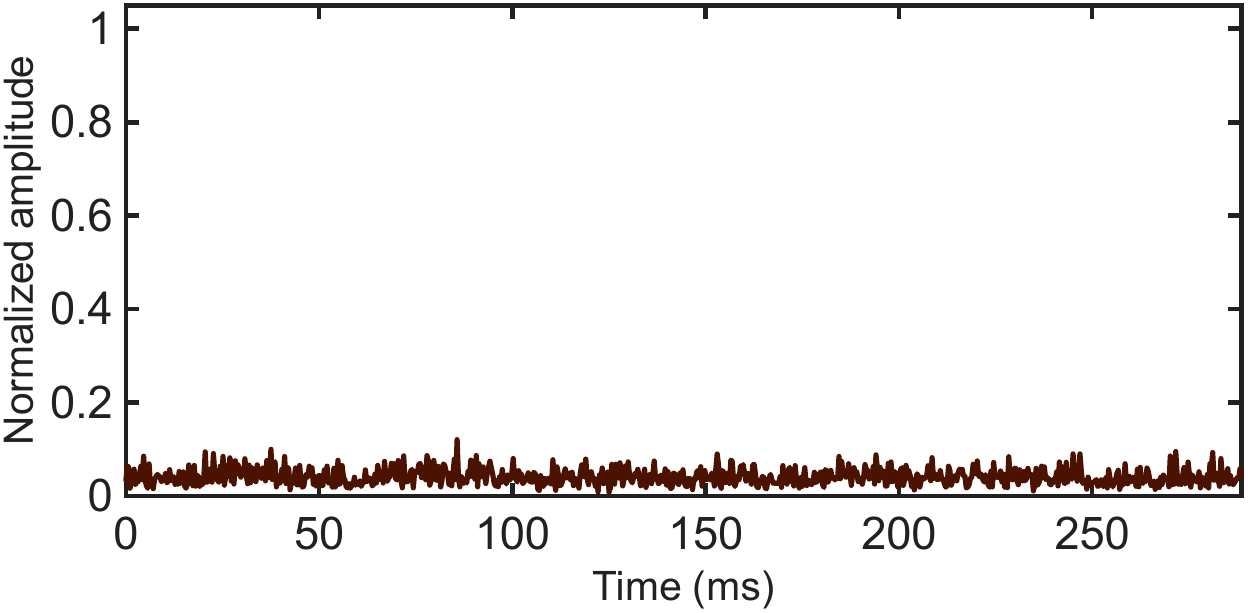}
    \end{minipage}

    \begin{minipage}{0.02\linewidth}
        \rotatebox{90}{\scriptsize I.P.}
    \end{minipage}\begin{minipage}{0.196\linewidth}
        \centering
        \hfill
        \end{minipage}\begin{minipage}{0.2\linewidth}
        \centering
        \includegraphics[width=0.99\linewidth]{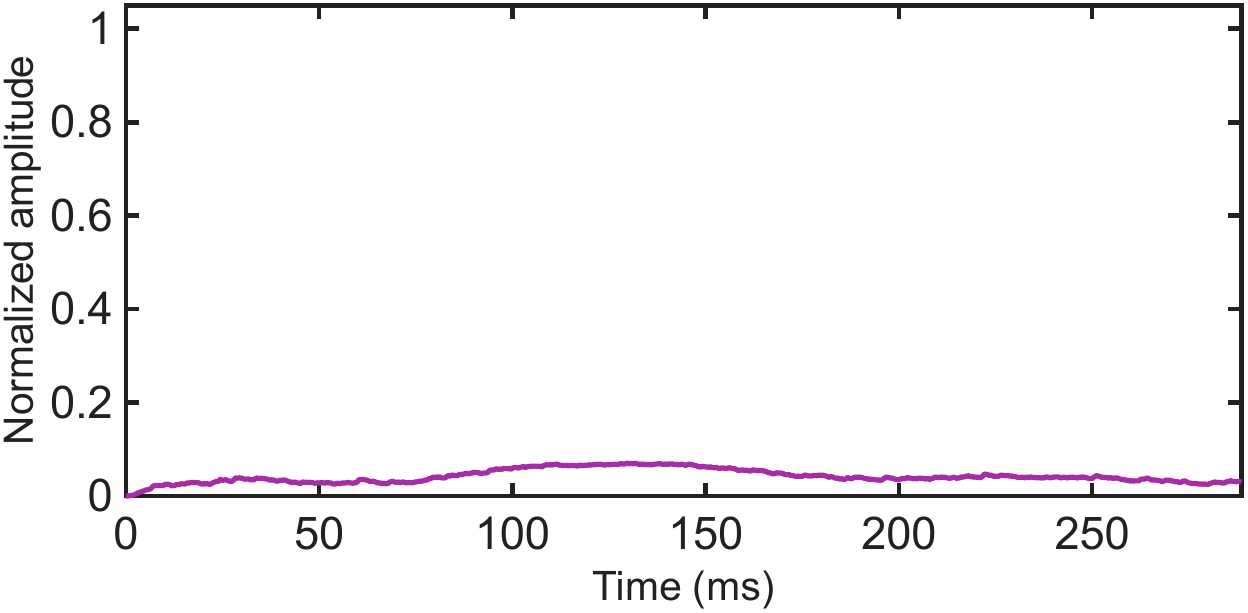}
    \end{minipage}\begin{minipage}{0.196\linewidth}
        \centering
        \includegraphics[width=0.99\linewidth]{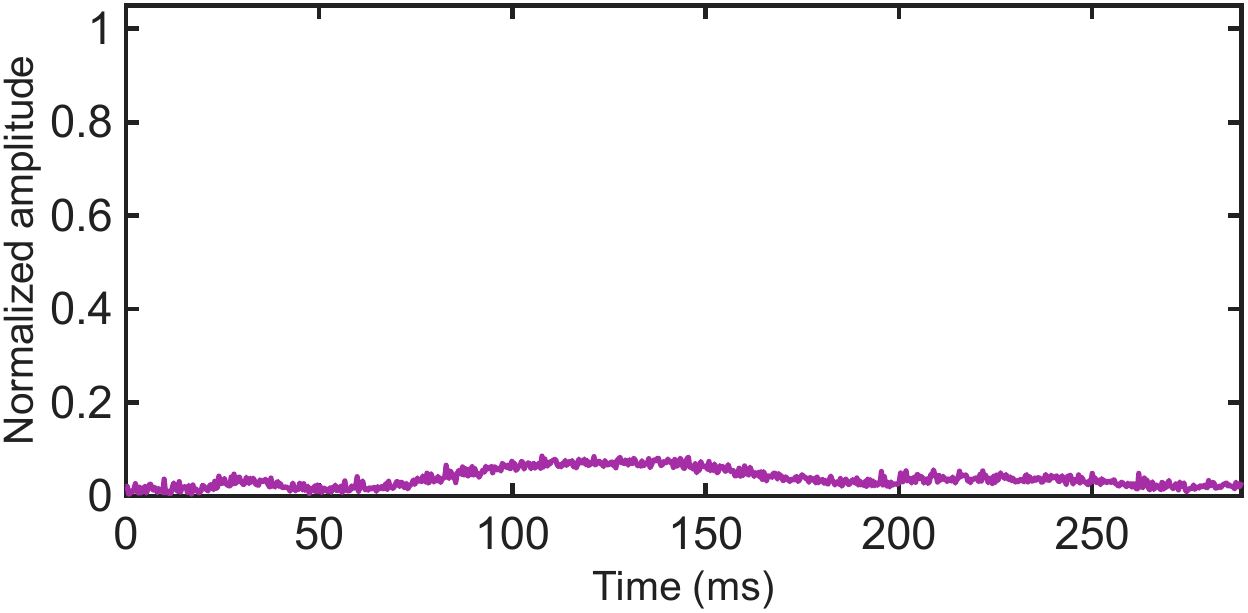}
    \end{minipage}\begin{minipage}{0.196\linewidth}
        \centering
        \includegraphics[width=0.99\linewidth]{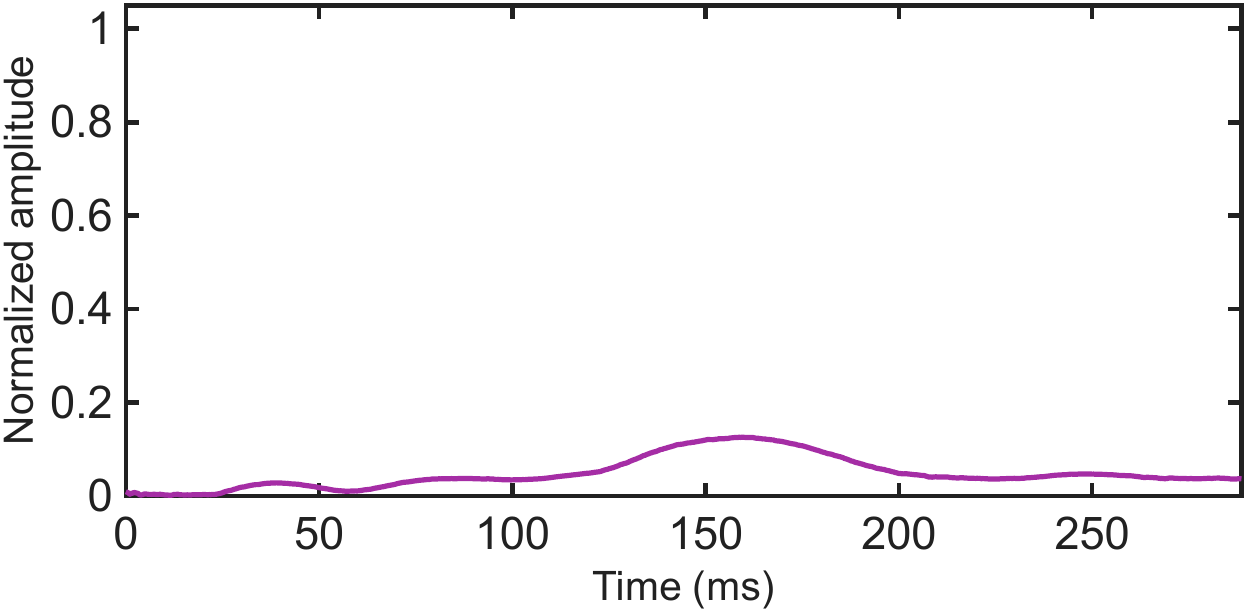}
    \end{minipage}\begin{minipage}{0.196\linewidth}
        \centering
        \includegraphics[width=0.99\linewidth]{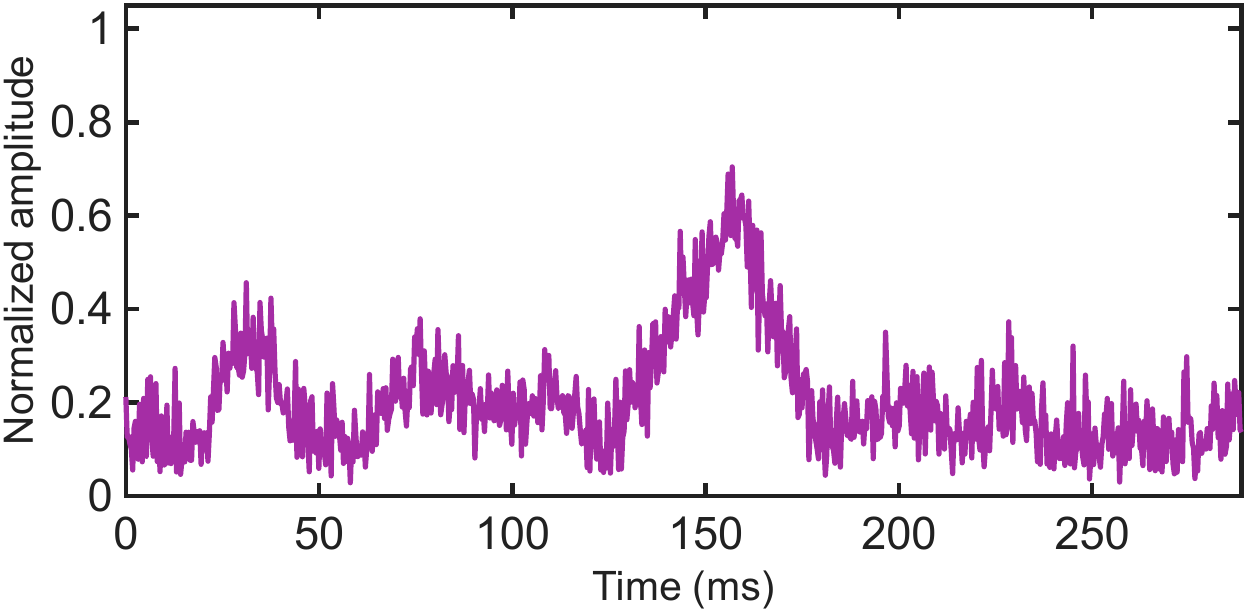}
    \end{minipage}

    \begin{minipage}{0.02\linewidth}
        \rotatebox{90}{\scriptsize Fusif.}
    \end{minipage}\begin{minipage}{0.196\linewidth}
        \centering
        \hfill
    \end{minipage}\begin{minipage}{0.2\linewidth}
        \centering
        \includegraphics[width=0.99\linewidth]{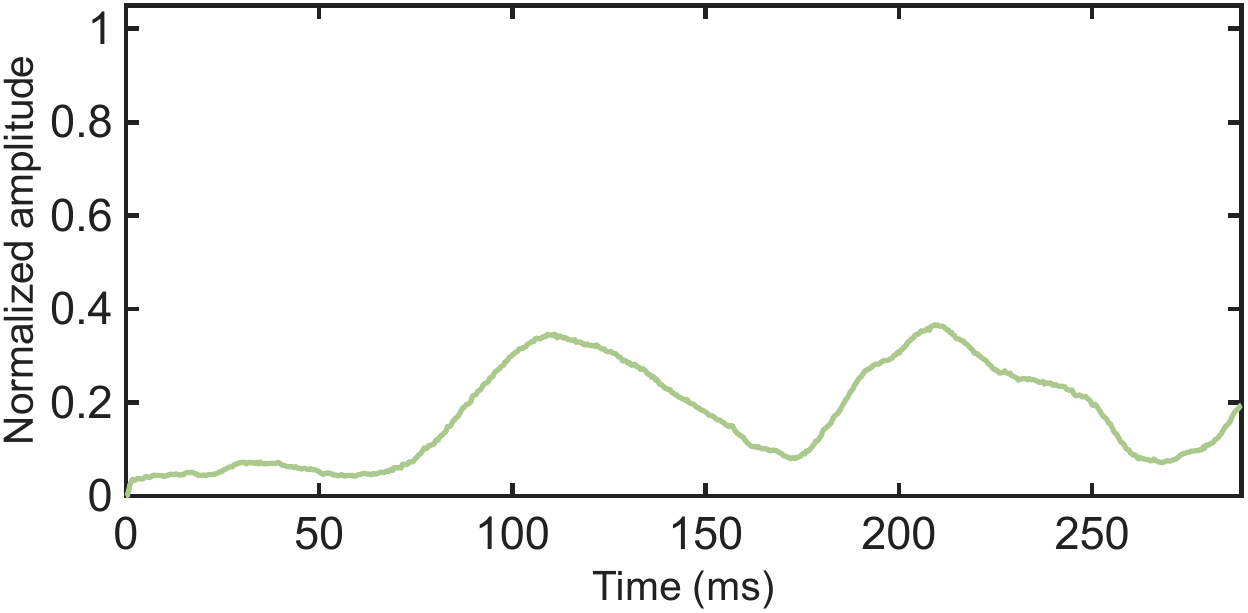}
    \end{minipage}\begin{minipage}{0.196\linewidth}
        \centering
        \includegraphics[width=0.99\linewidth]{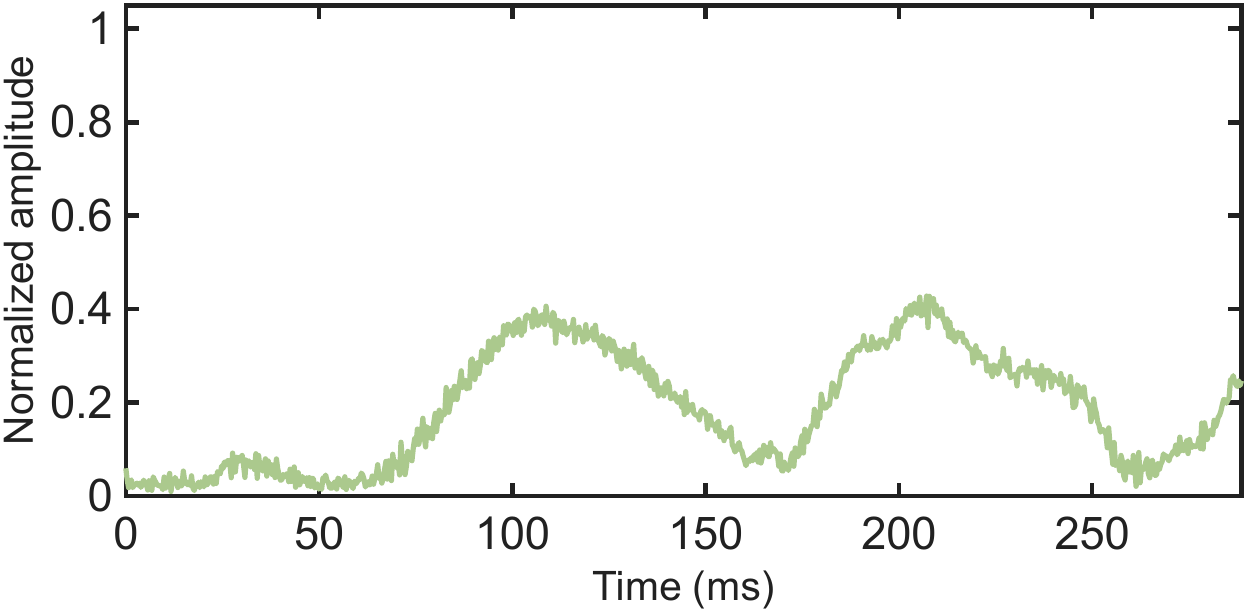}
    \end{minipage}\begin{minipage}{0.196\linewidth}
        \centering
        \includegraphics[width=0.99\linewidth]{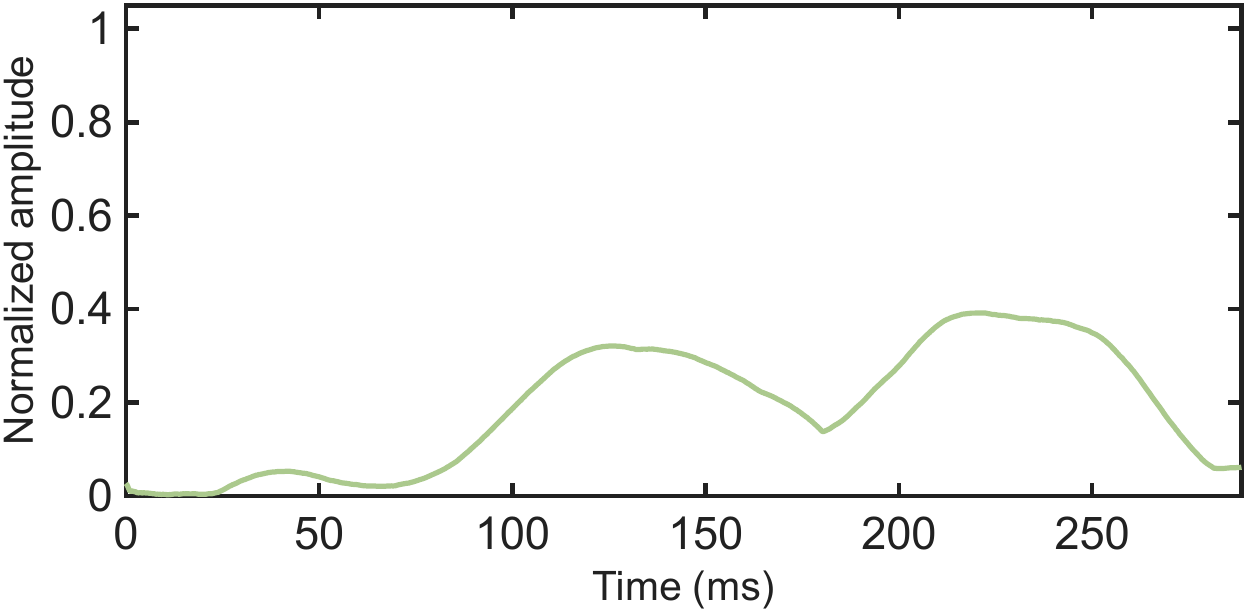}
    \end{minipage}\begin{minipage}{0.196\linewidth}
        \centering
        \includegraphics[width=0.99\linewidth]{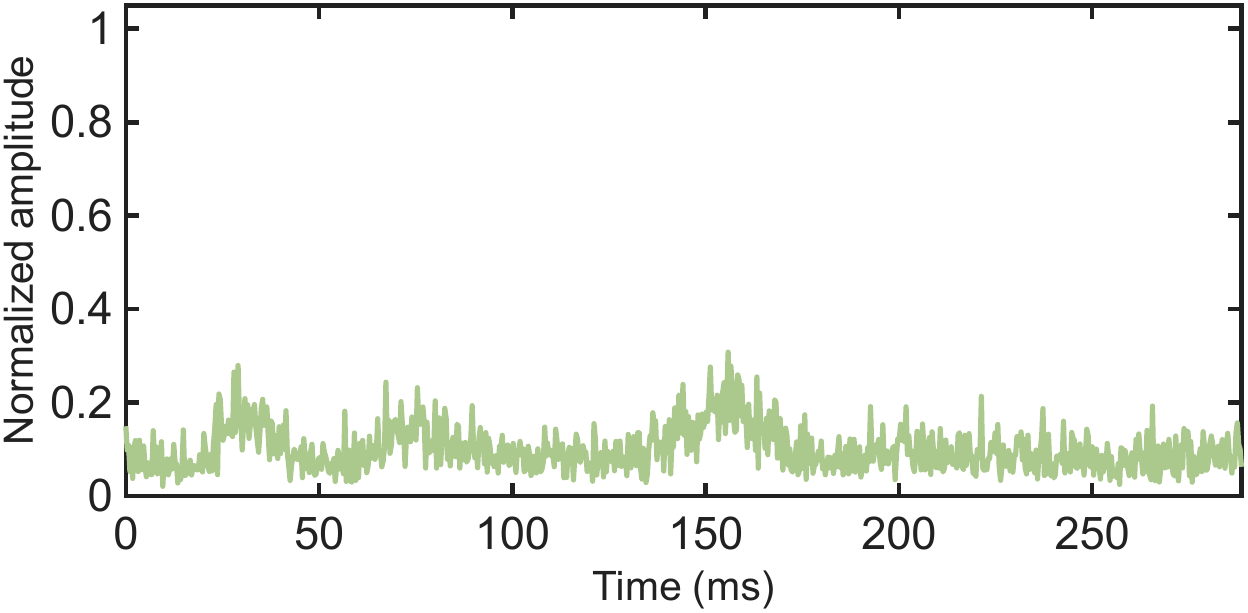}
    \end{minipage}
    \caption{Parcel-wise temporal activity at 20 dB SNR. Rows correspond to the selected cortical parcels, and columns show the simulated ground-truth time series followed by Change-rate Standardized Kalman filter (CrSKF), Standardized Kalman filter (SKF), Kalman filter (KF), and Linearly Constrained Minimum Variance (LCMV) beamformer reconstructions.}
    \label{fig:timeseries20dB}
\end{figure}

For robustness testing, zero-mean Gaussian sensor noise was added independently to each electrode. The noise level was set from
$\mathrm{SNR}=20\log_{10}(\|{\bf y}\|_2/\|{\boldsymbol\varepsilon}\|_2)$,
with noiseless, 25 dB, and 20 dB conditions. A representative realization at each noise level is shown in Figures \ref{fig:timeseries0dB}--\ref{fig:timeseries20dB}; the comparisons should therefore be interpreted as illustrative rather than as Monte Carlo statistics.

Both standardized Kalman filters recover the principal lateral-occipital activity at approximately the correct time frames. Activity in the lingual and superior-parietal parcels is less well separated and partly follows the stronger lateral-occipital response. A spurious superior-parietal peak is visible around 110 ms, while the later frontal activity is recovered in the caudal middle frontal parcel. CrSKF produces the smoothest temporal trajectories, consistent with its explicit change-rate state and observation-difference constraint; however, this smoothing can reduce the apparent separation between nearby components.

Conventional KF also yields smooth trajectories, but the lateral-occipital, lingual, and superior-parietal responses are more strongly mixed and some reconstructed peaks are delayed. LCMV produces simultaneous activity in several posterior parcels and suppresses part of the lingual and frontal activity. Its parcel-wise trajectories change more markedly as noise is added. Because beamformer performance depends on covariance estimation and diagonal loading, this observation is specific to the present full-epoch covariance estimate and noise-level regularization and should not be interpreted as a general lack of noise robustness.

Temporal relationships were quantified by cross-correla\-tion between the reconstructed lateral-occipital time series and each comparison parcel. For a lag $\ell$, we use
$C_{xy}(\ell)=\sum_t x(t)y(t+\ell)$; positive lag means that the comparison parcel follows the lateral-occipital reference, whereas negative lag means that it precedes it. Table \ref{tab:timeshifttable} reports the lag maximizing the cross-correlation. The near-zero lag between the lateral-occipital and superior-parietal estimates indicates that the two reconstructed responses are not temporally well separated, whereas the caudal middle frontal activity is shifted substantially later. These lags describe temporal ordering and association only; they do not by themselves establish causality.


\begin{table}[h!]
    \centering
    \footnotesize{
    \begin{tabular}{c|c|c|c|c|c}\hline
       \diagbox{Method}{Parcel}  & ling. & S.P. & C.M.F. & I.P. & fusif. \\ \hline
       \multicolumn{6}{c}{noiseless} \\ \hline
        CrSKF & 0 & -1 & 119 & 9 & 0\\ \hline
        SKF & 0 & -1 & 119 & 9 & 0 \\ \hline
        KF & 0 & -2 & 119 & 28 & 1 \\ \hline
        LCMV & 0 & 0 & 0 & 0 & 0 \\ \hline
        \multicolumn{6}{c}{25 \unit{\decibel} noise} \\ \hline
        CrSKF & 0 & 0 & 118 & 9 & 0\\ \hline
        SKF & 0 & 0 & 118 & 10 & 0 \\ \hline
        KF & 0 & -1 & 119 & 27 & 1 \\ \hline
        LCMV & 0 & 0 & 0 & 0 & 0 \\ \hline
        \multicolumn{6}{c}{20 \unit{\decibel} noise} \\ \hline
        CrSKF & 0 & 0 & 118 & 0 & 0\\ \hline
        SKF & 0 & 0 & 119 & 9 & 0 \\ \hline
        KF & 0 & -1 & 119 & 27 & 1 \\ \hline
        LCMV & 0 & 0 & 0 & 0 & 0 \\ \hline
    \end{tabular}}

    \caption{Time lags (ms) maximizing the cross-correlation between the reconstructed lateral-occipital time series and the reconstructed time series in the lingual (ling.), superior parietal (S.P.), caudal middle frontal (C.M.F.), inferior parietal (I.P.), and fusiform (fusif.) parcels. Positive lag means that the comparison parcel follows the lateral-occipital reference; negative lag means that it precedes it. Results are shown for the noiseless, 25 dB, and 20 dB conditions.}
    \label{tab:timeshifttable}
\end{table}

\begin{figure}[h!]
    \centering
    \begin{minipage}{0.02\linewidth}
        \hfill
    \end{minipage}\begin{minipage}{0.98\linewidth}
        \begin{minipage}{0.16\linewidth}
            \centering
            \scriptsize{V1}
        \end{minipage}\begin{minipage}{0.16\linewidth}
            \centering
            \scriptsize{V2d}
        \end{minipage}\begin{minipage}{0.16\linewidth}
            \centering
            \scriptsize{V2v+V3d}
        \end{minipage}\begin{minipage}{0.16\linewidth}
            \centering
            \scriptsize{V3v}
        \end{minipage}\begin{minipage}{0.16\linewidth}
        \centering
            \scriptsize{LIP}
        \end{minipage}\begin{minipage}{0.16\linewidth}
        \centering
            \scriptsize{FEF}
        \end{minipage}\end{minipage}
        
    \begin{minipage}{0.02\linewidth}
        \rotatebox{90}{\scriptsize CrSKF}
    \end{minipage}\begin{minipage}{0.98\linewidth}
        \begin{minipage}{0.16\linewidth}
            \includegraphics[trim={5cm 1cm 14cm 8cm},clip,width=\linewidth]{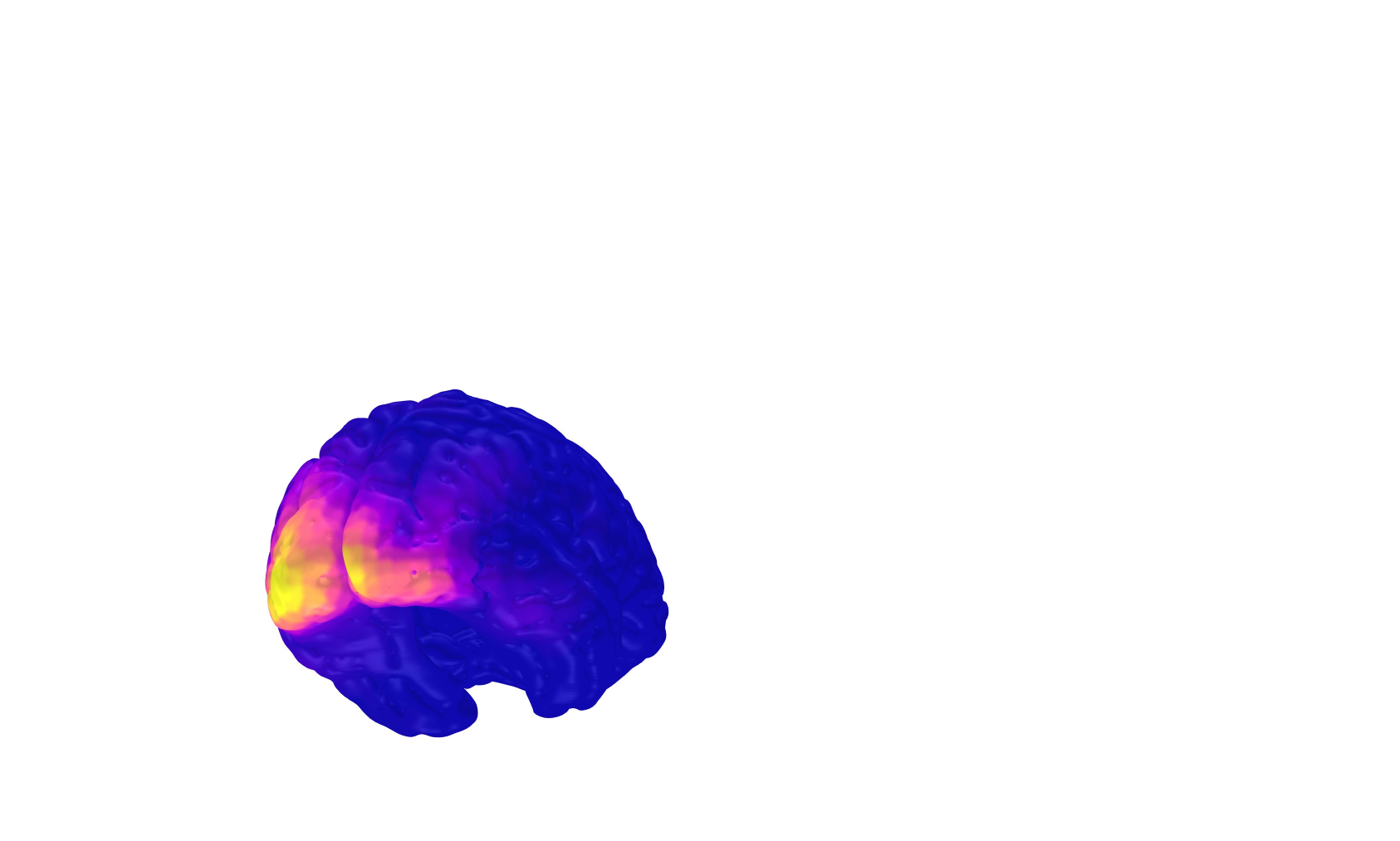}
        \end{minipage}\begin{minipage}{0.16\linewidth}
            \includegraphics[trim={5cm 1cm 14cm 8cm},clip,width=\linewidth]{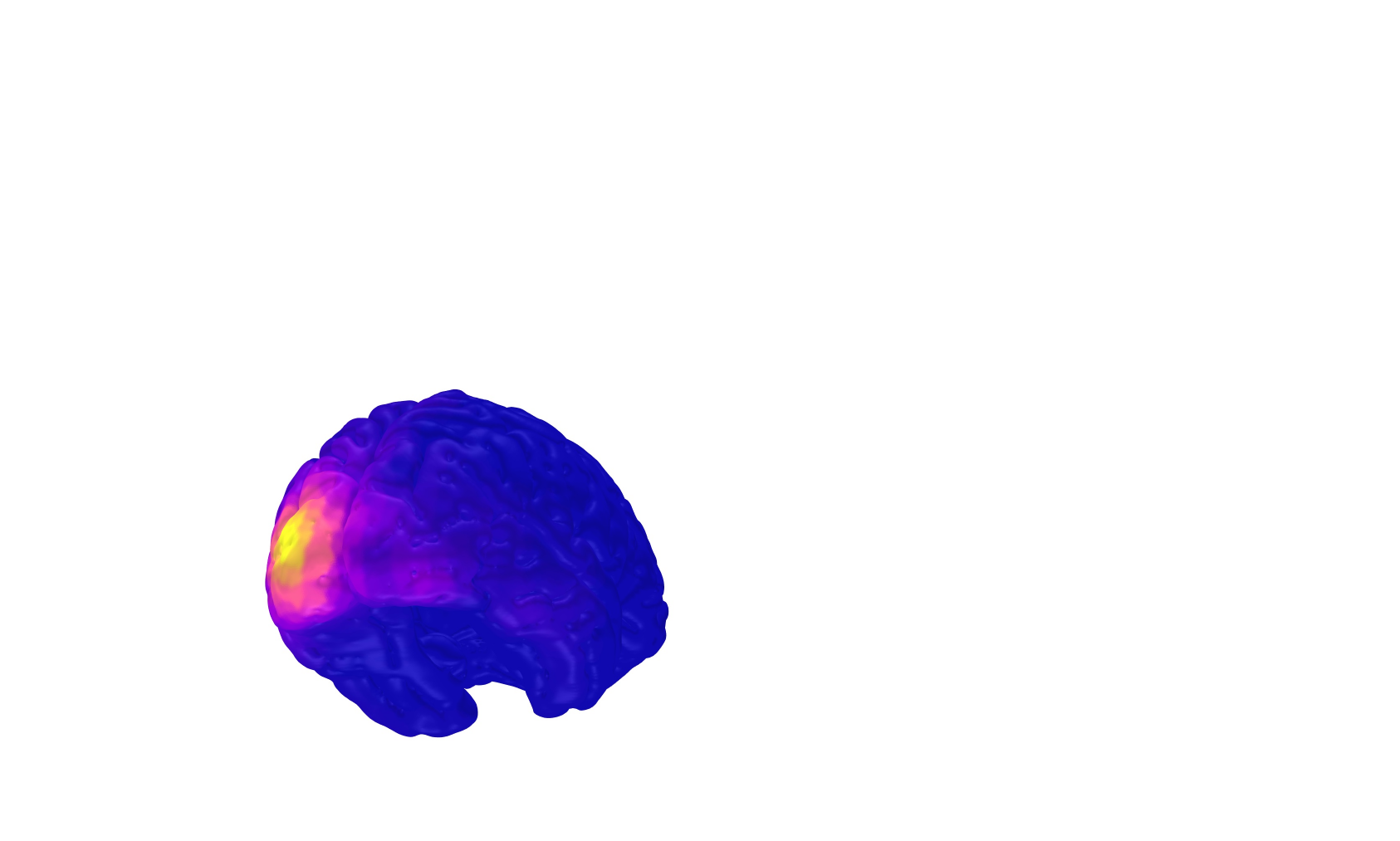}
        \end{minipage}\begin{minipage}{0.16\linewidth}
            \includegraphics[trim={5cm 1cm 14cm 8cm},clip,width=\linewidth]{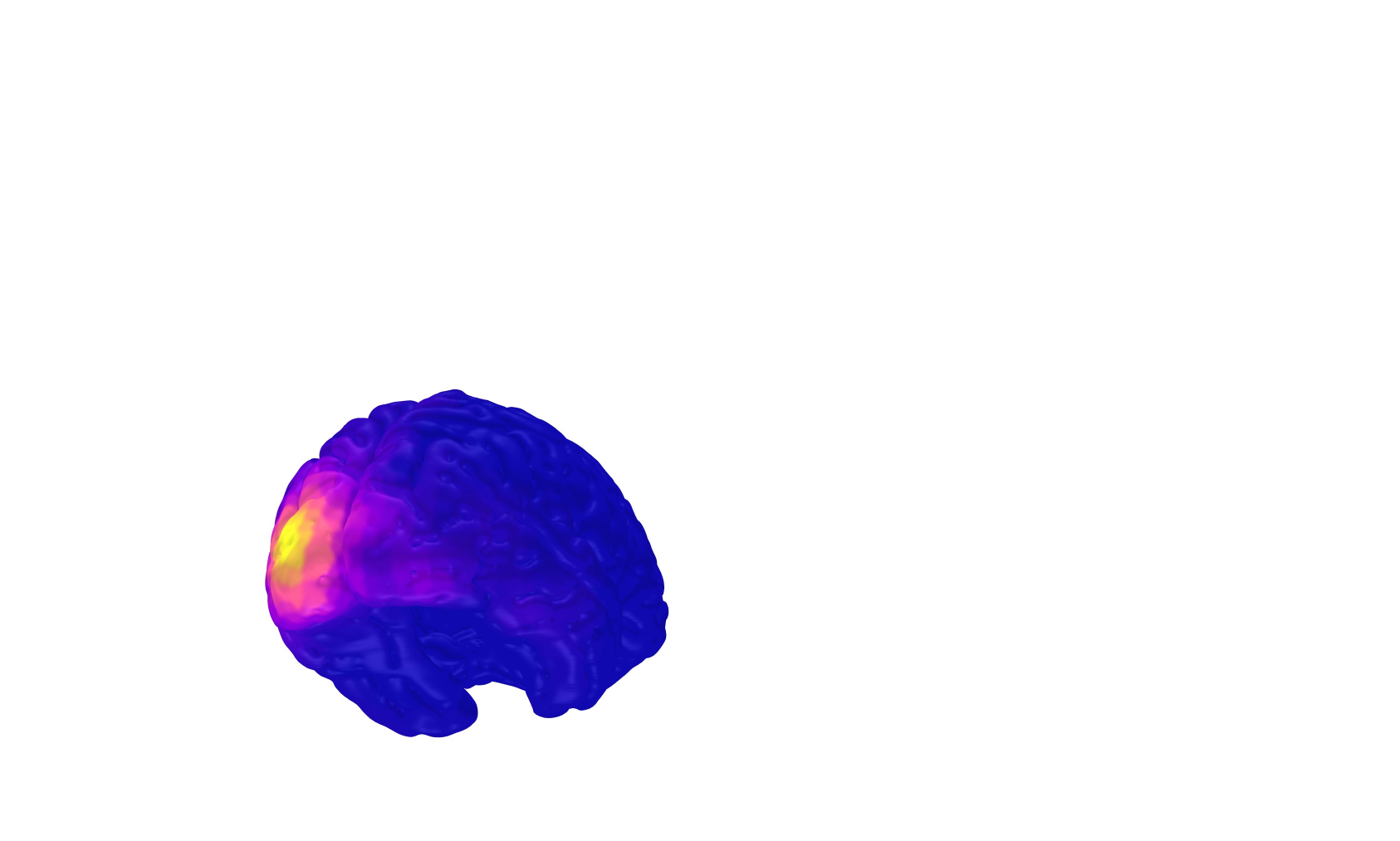}
        \end{minipage}\begin{minipage}{0.16\linewidth}
            \includegraphics[trim={5cm 1cm 14cm 8cm},clip,width=\linewidth]{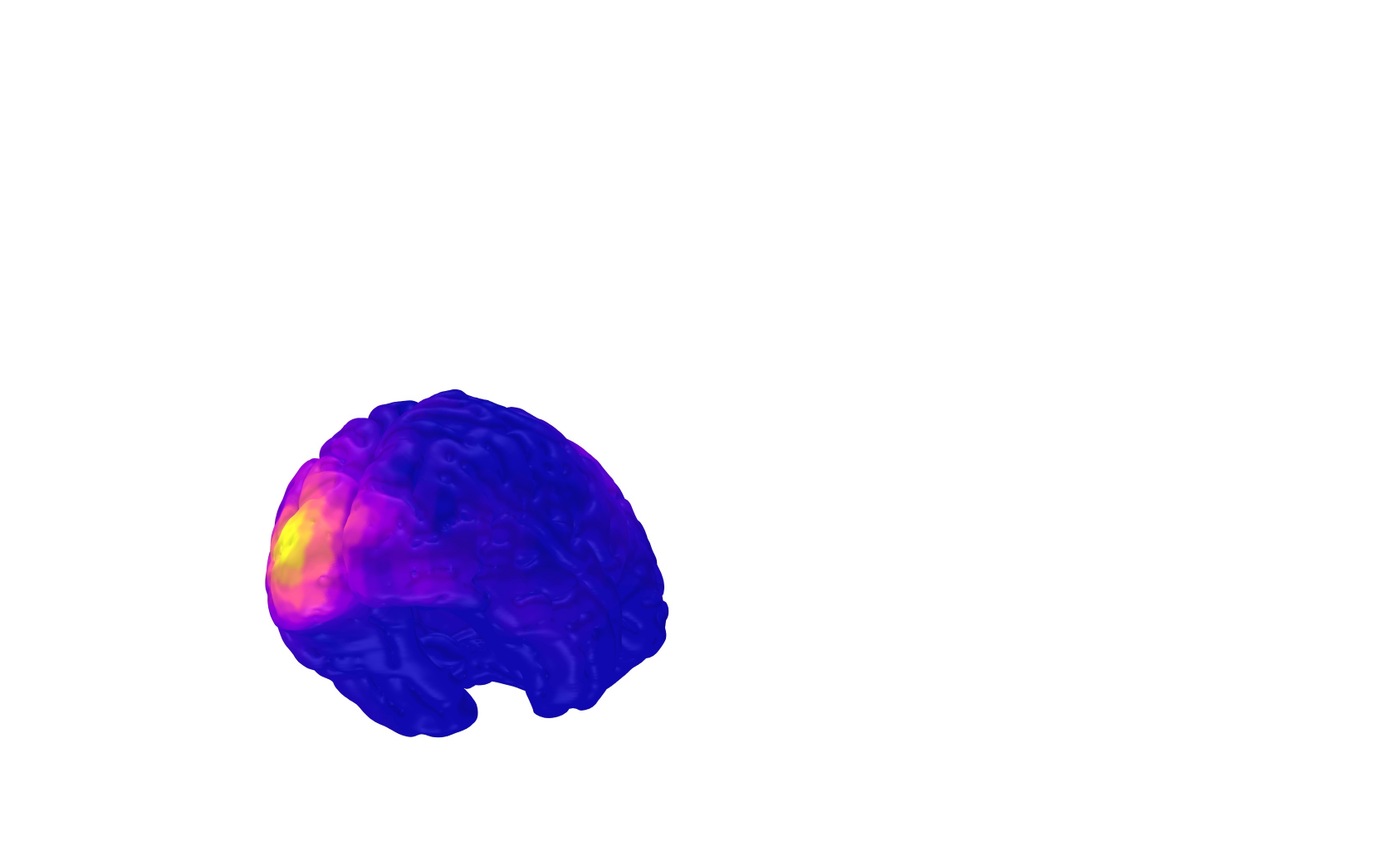}
        \end{minipage}\begin{minipage}{0.16\linewidth}
            \includegraphics[trim={5cm 1cm 14cm 8cm},clip,width=\linewidth]{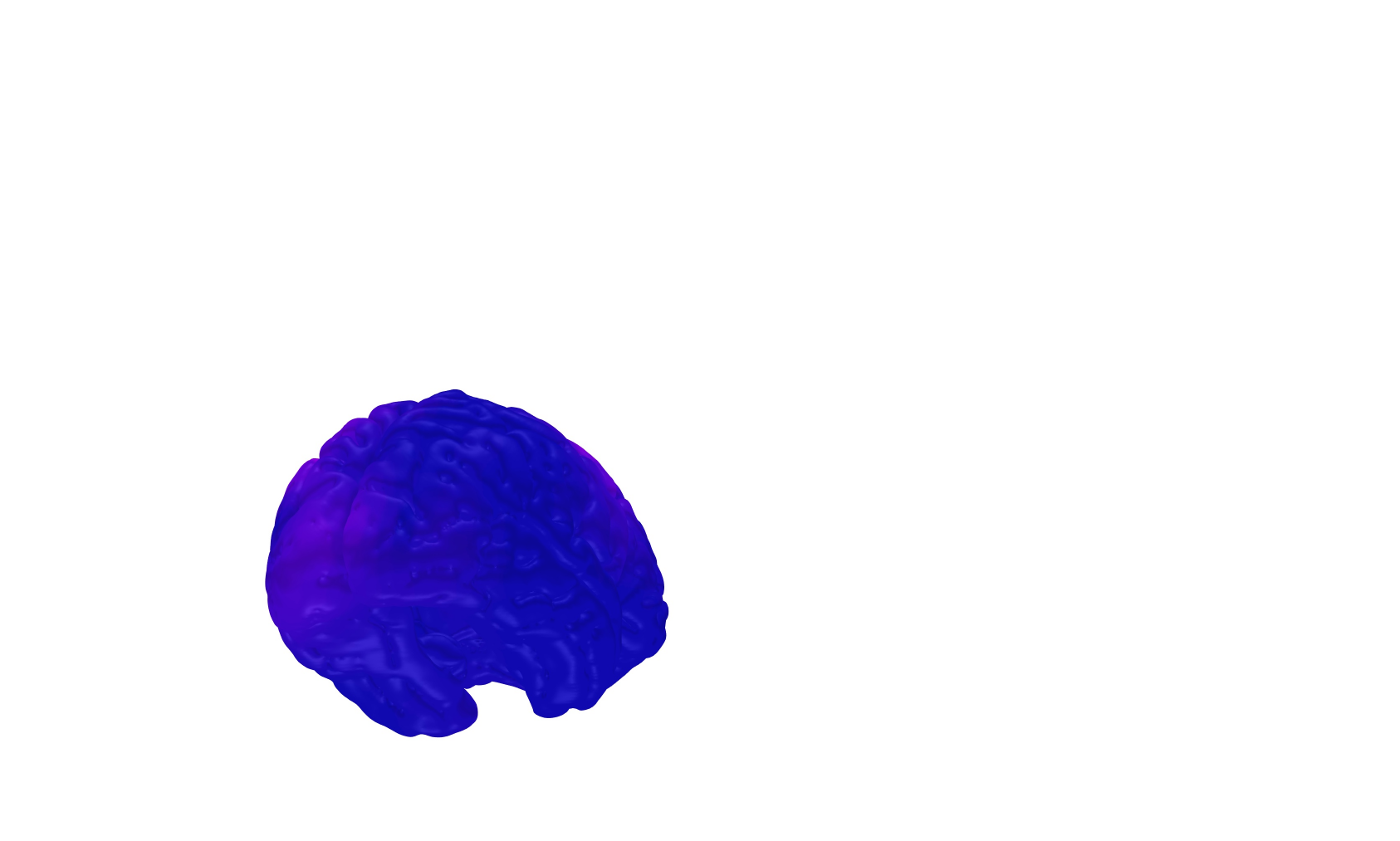}
        \end{minipage}\begin{minipage}{0.16\linewidth}
            \includegraphics[trim={5cm 1cm 14cm 8cm},clip,width=\linewidth]{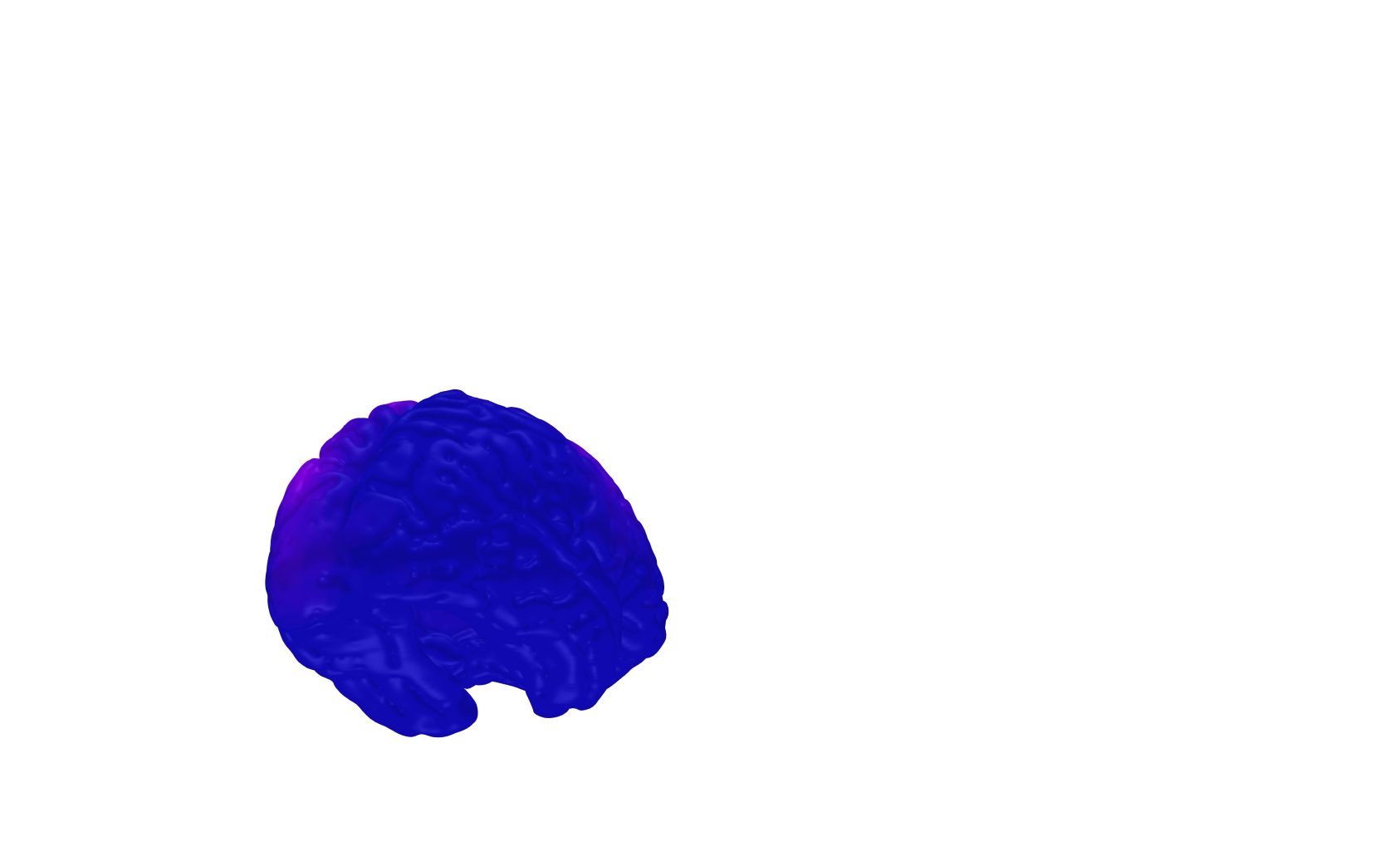}
        \end{minipage}

        \begin{minipage}{0.16\linewidth}
            \includegraphics[trim={5cm 1cm 14cm 8cm},clip,width=\linewidth]{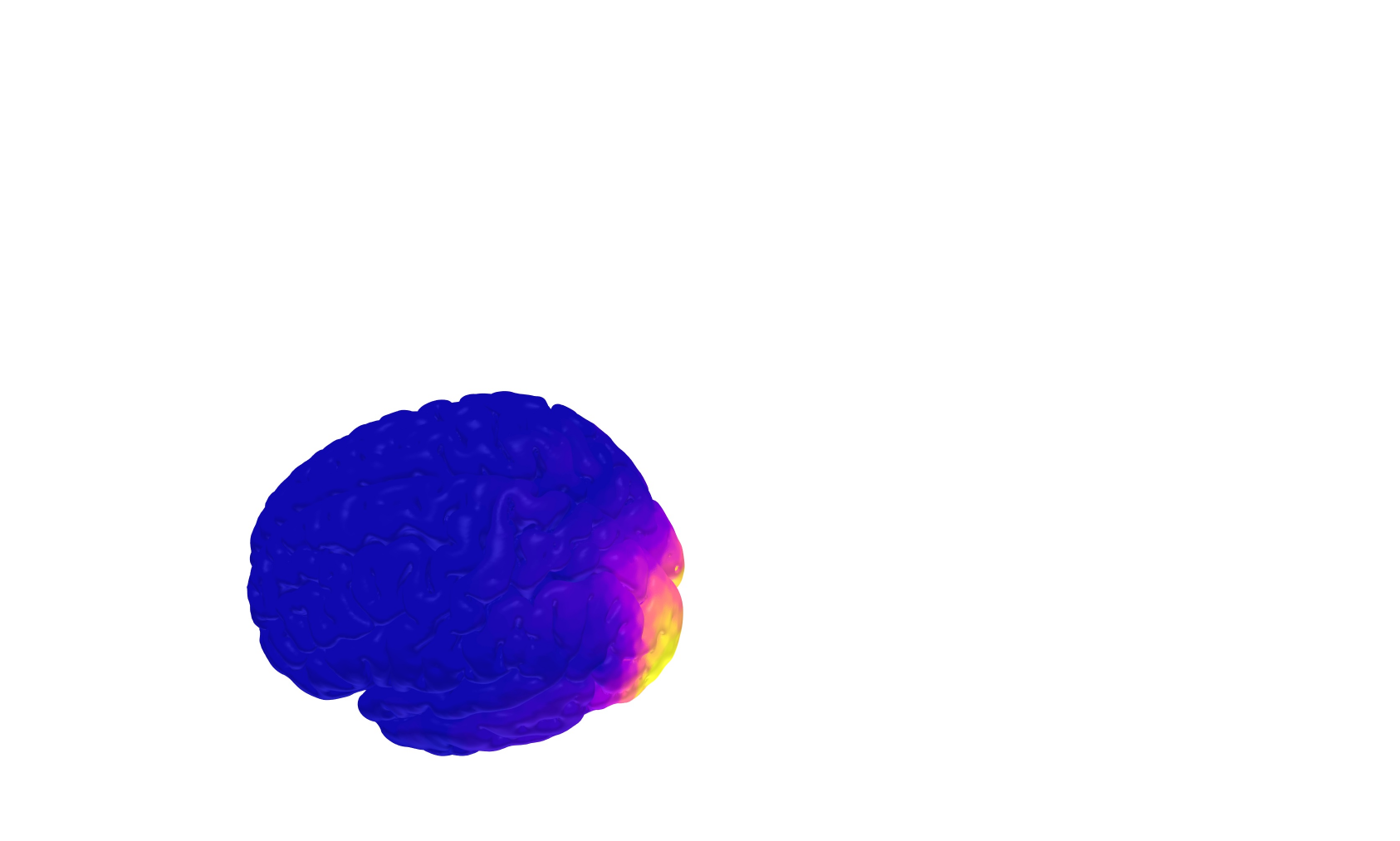}
        \end{minipage}\begin{minipage}{0.16\linewidth}
            \includegraphics[trim={5cm 1cm 14cm 8cm},clip,width=\linewidth]{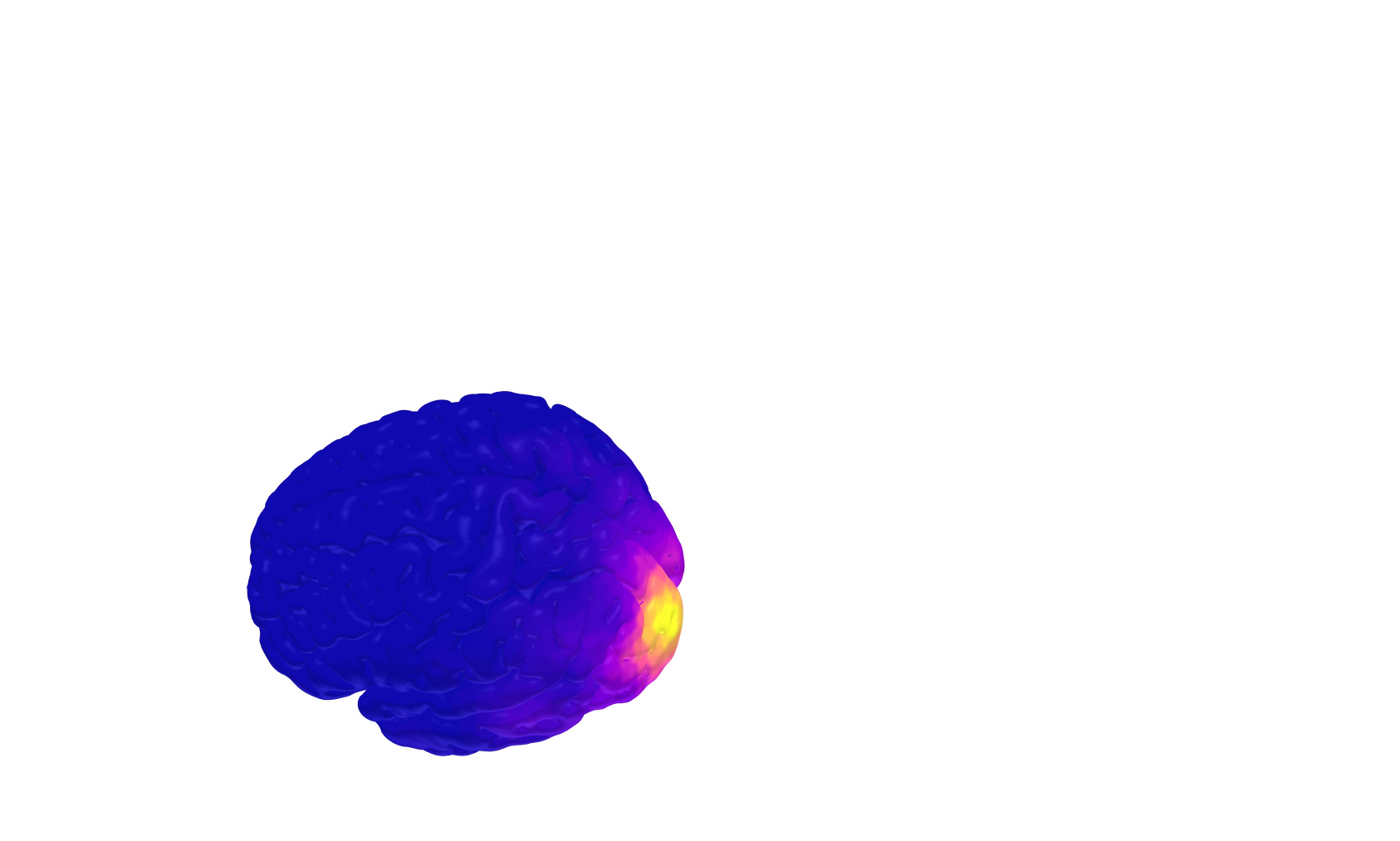}
        \end{minipage}\begin{minipage}{0.16\linewidth}
            \includegraphics[trim={5cm 1cm 14cm 8cm},clip,width=\linewidth]{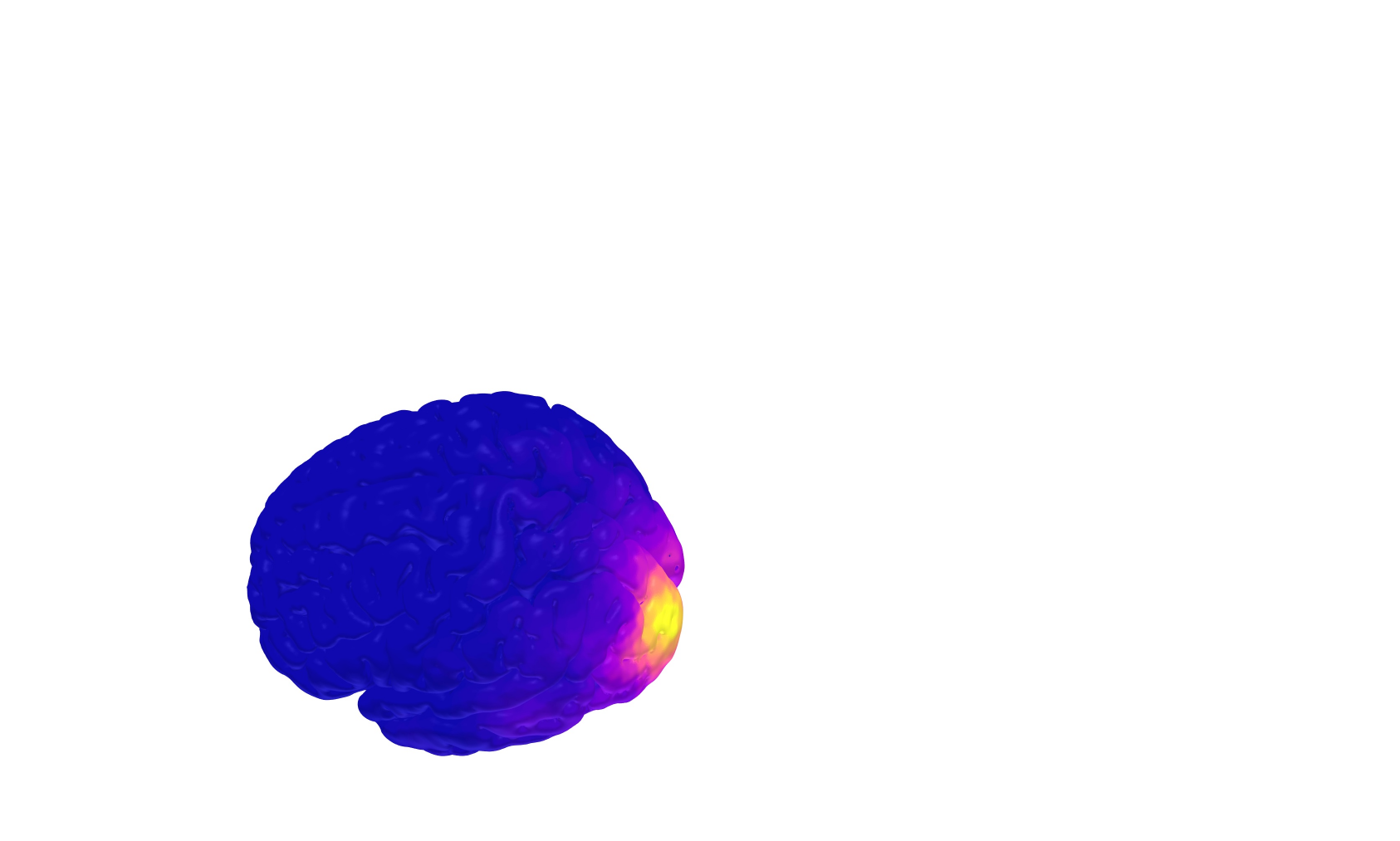}
        \end{minipage}\begin{minipage}{0.16\linewidth}
            \includegraphics[trim={5cm 1cm 14cm 8cm},clip,width=\linewidth]{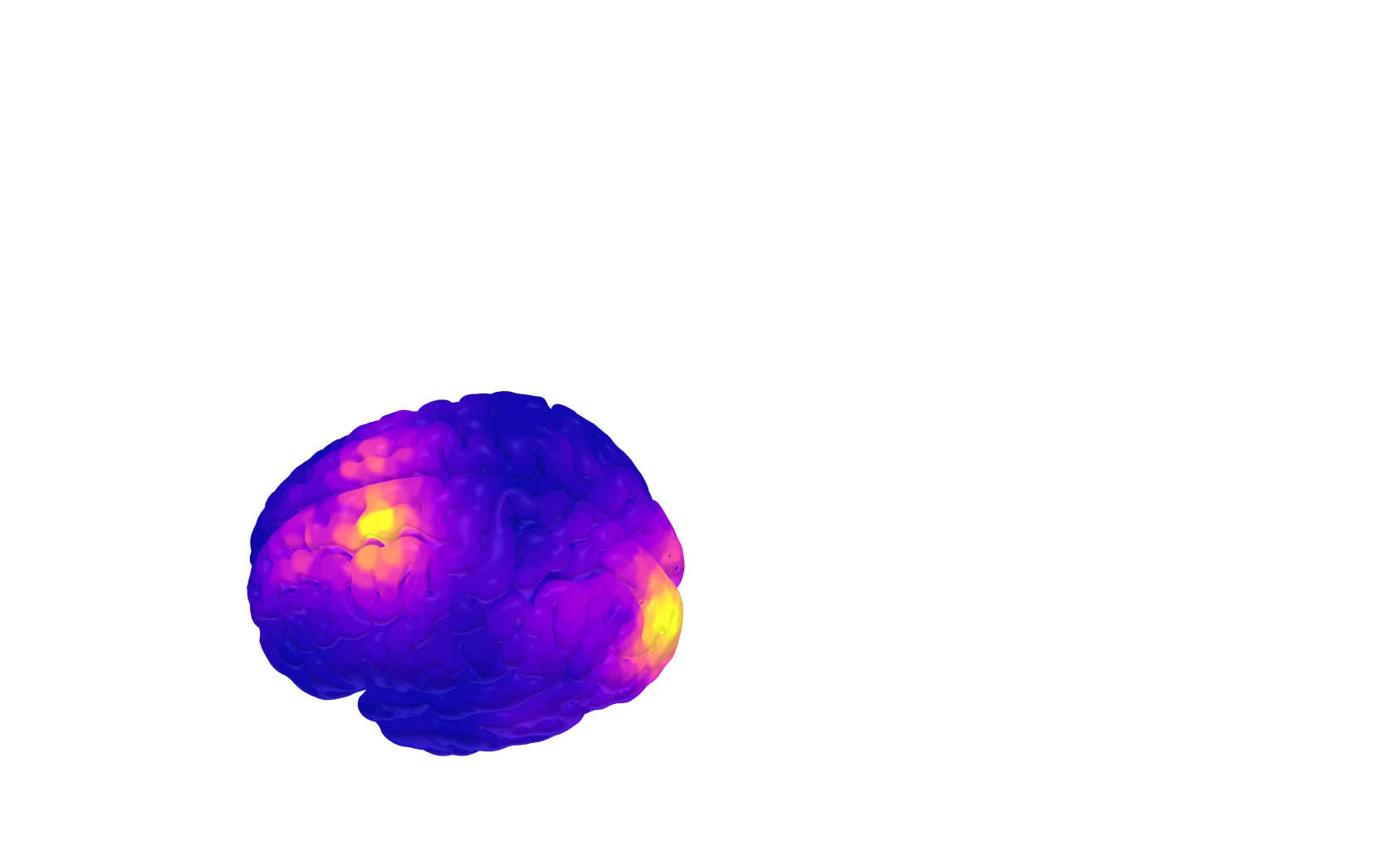}
        \end{minipage}\begin{minipage}{0.16\linewidth}
            \includegraphics[trim={5cm 1cm 14cm 8cm},clip,width=\linewidth]{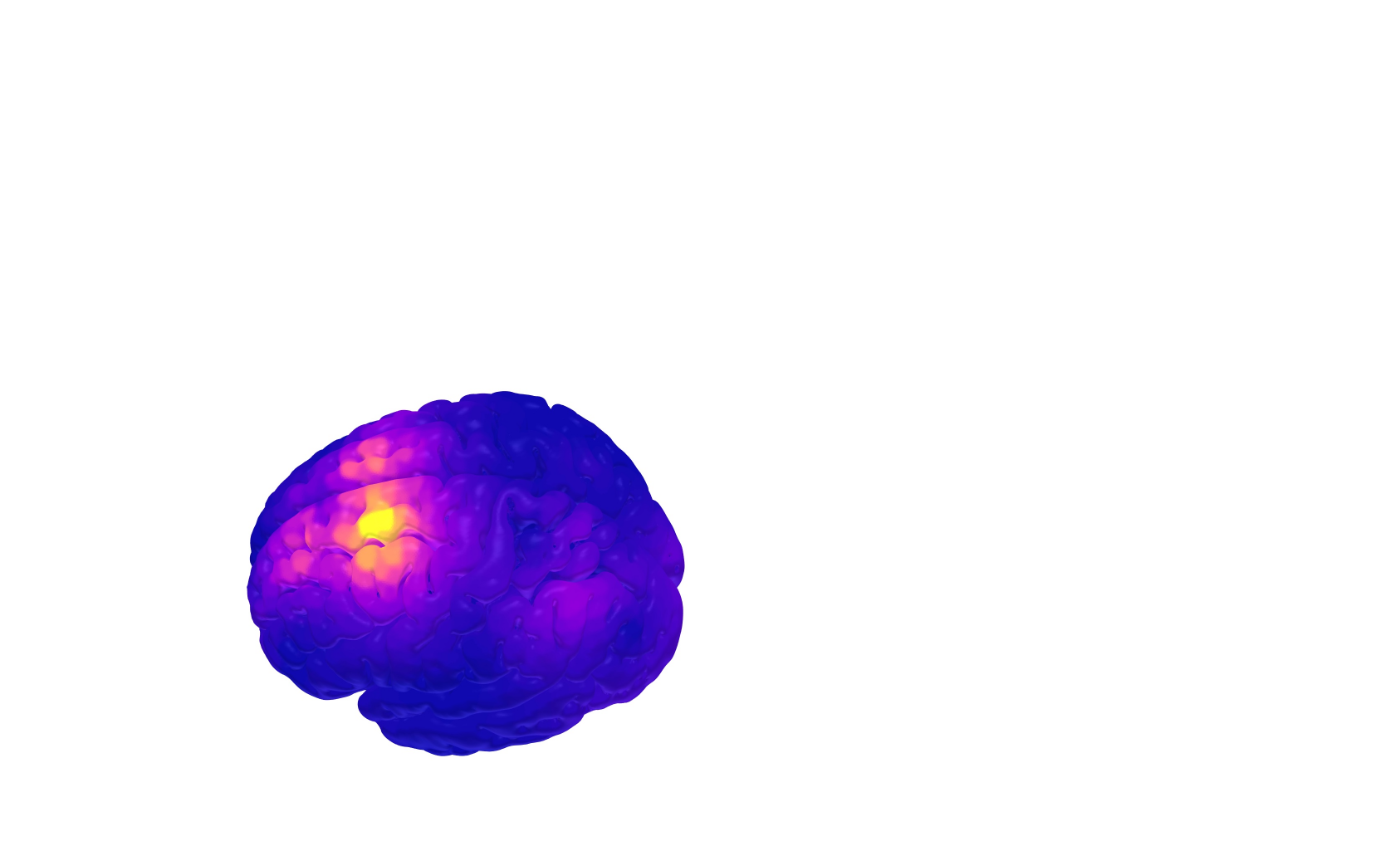}
        \end{minipage}\begin{minipage}{0.16\linewidth}
            \includegraphics[trim={5cm 1cm 14cm 8cm},clip,width=\linewidth]{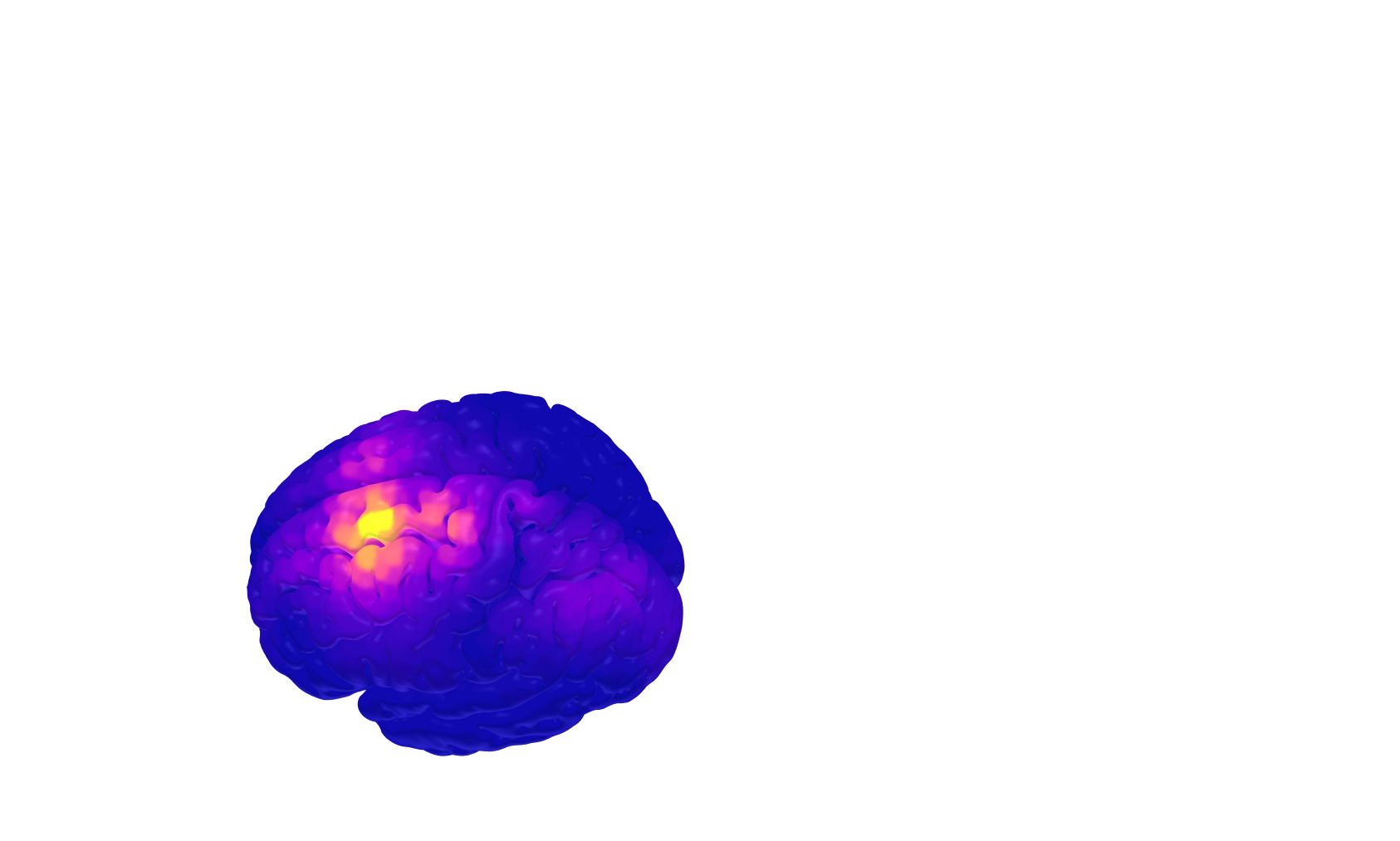}
        \end{minipage}
    \end{minipage}
    \noindent\rule[7pt]{\linewidth}{0.4pt}

    \begin{minipage}{0.02\linewidth}
        \rotatebox{90}{\scriptsize SKF}
    \end{minipage}\begin{minipage}{0.98\linewidth}
        \begin{minipage}{0.16\linewidth}
            \includegraphics[trim={5cm 1cm 14cm 8cm},clip,width=\linewidth]{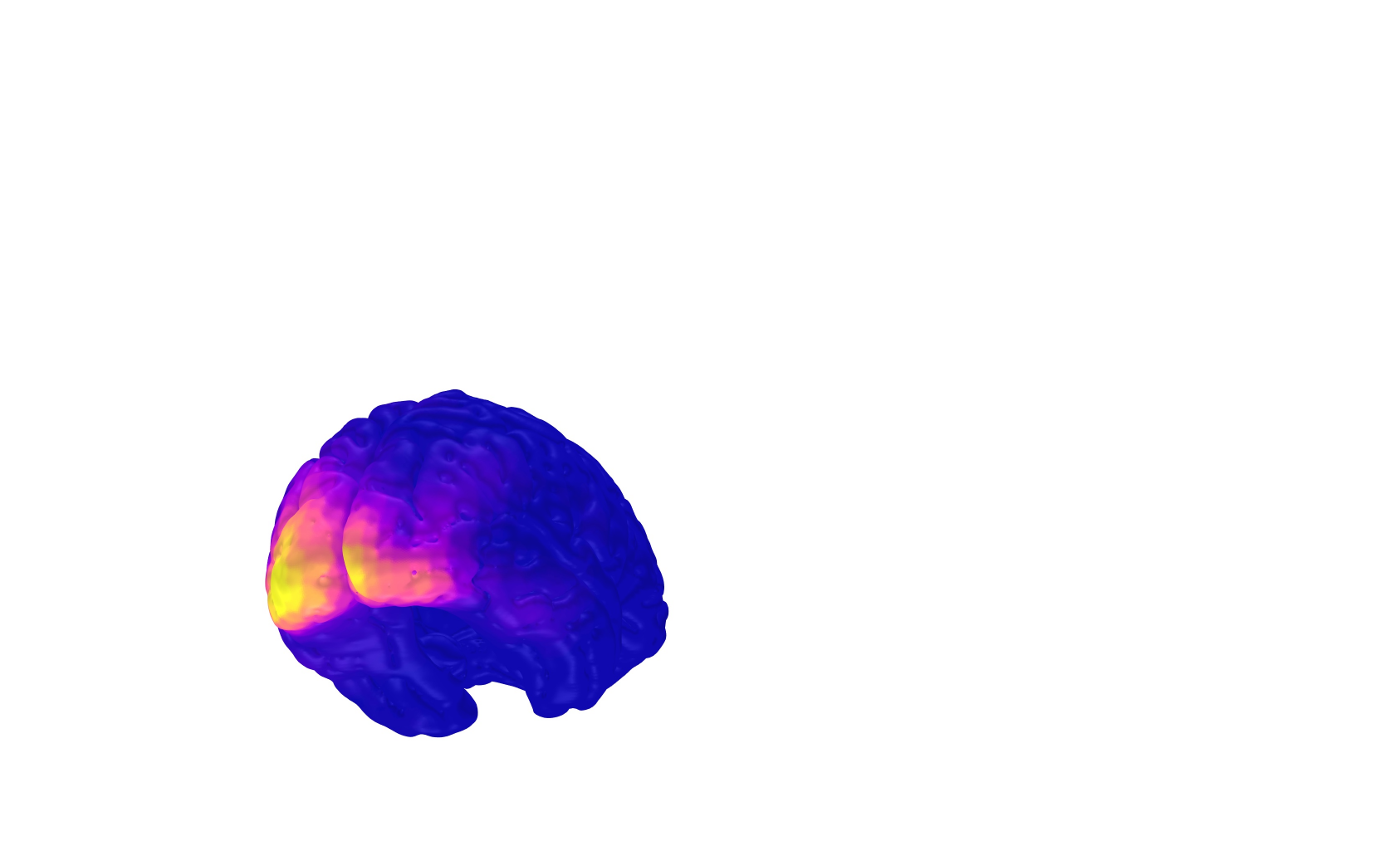}
        \end{minipage}\begin{minipage}{0.16\linewidth}
            \includegraphics[trim={5cm 1cm 14cm 8cm},clip,width=\linewidth]{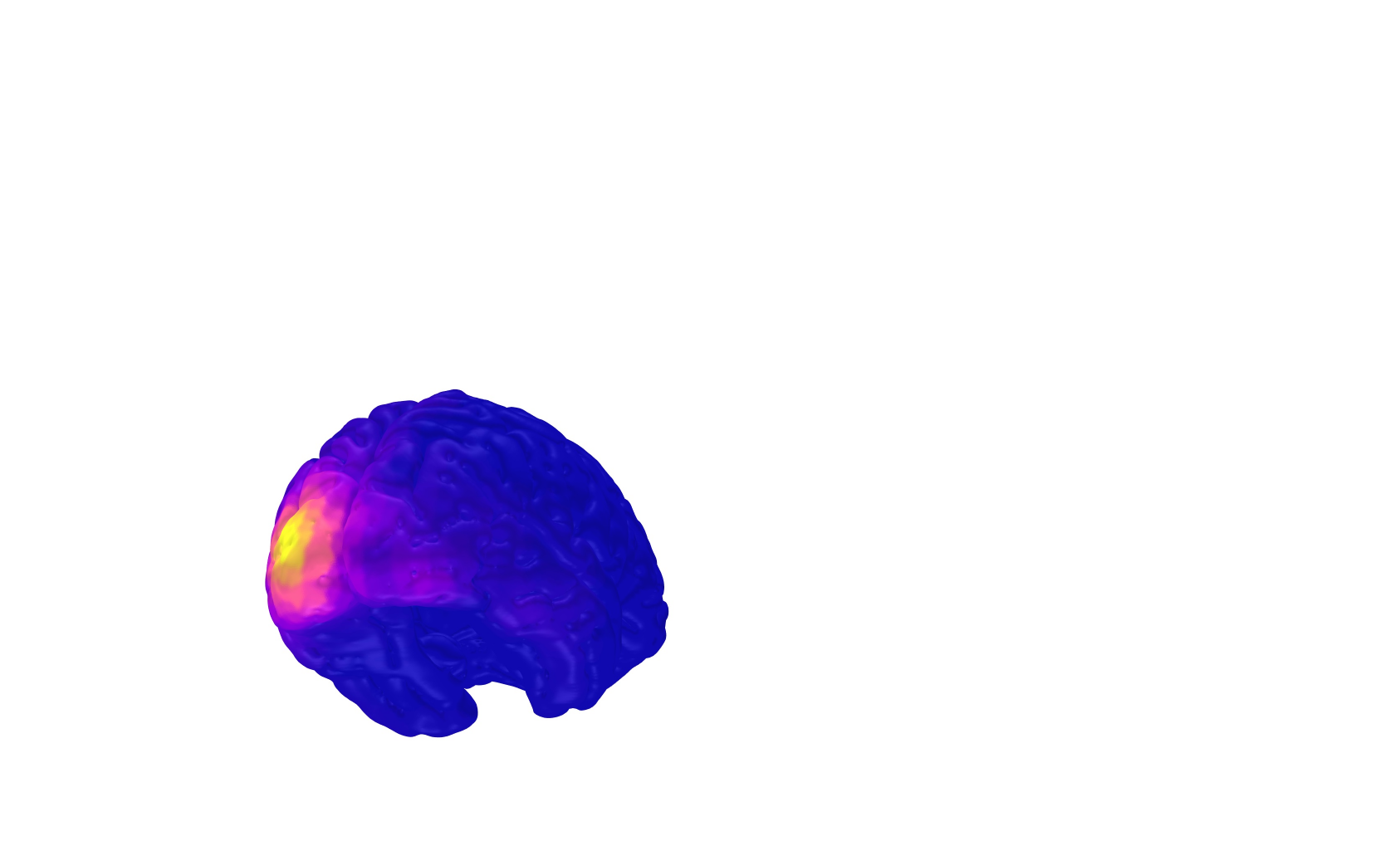}
        \end{minipage}\begin{minipage}{0.16\linewidth}
            \includegraphics[trim={5cm 1cm 14cm 8cm},clip,width=\linewidth]{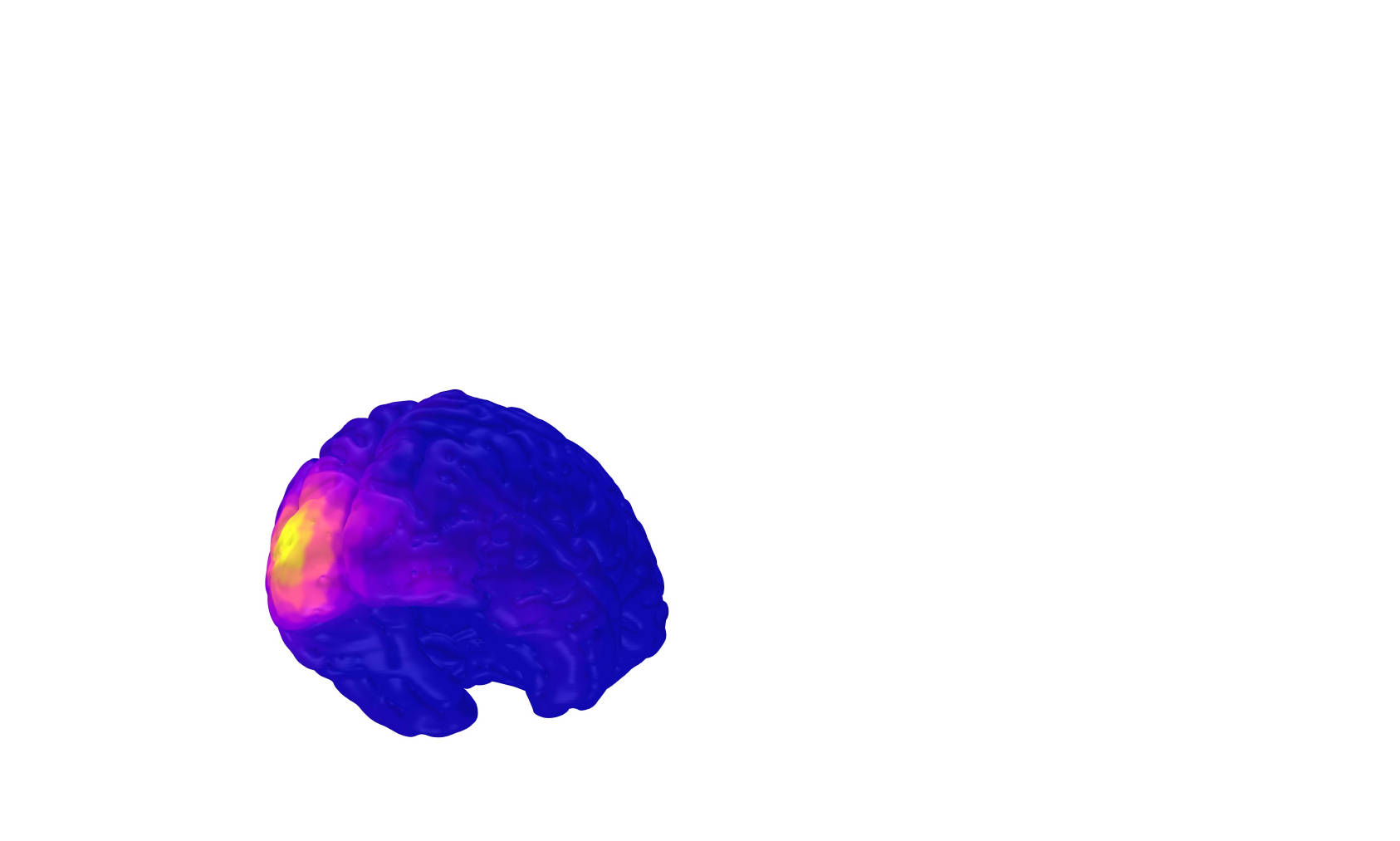}
        \end{minipage}\begin{minipage}{0.16\linewidth}
            \includegraphics[trim={5cm 1cm 14cm 8cm},clip,width=\linewidth]{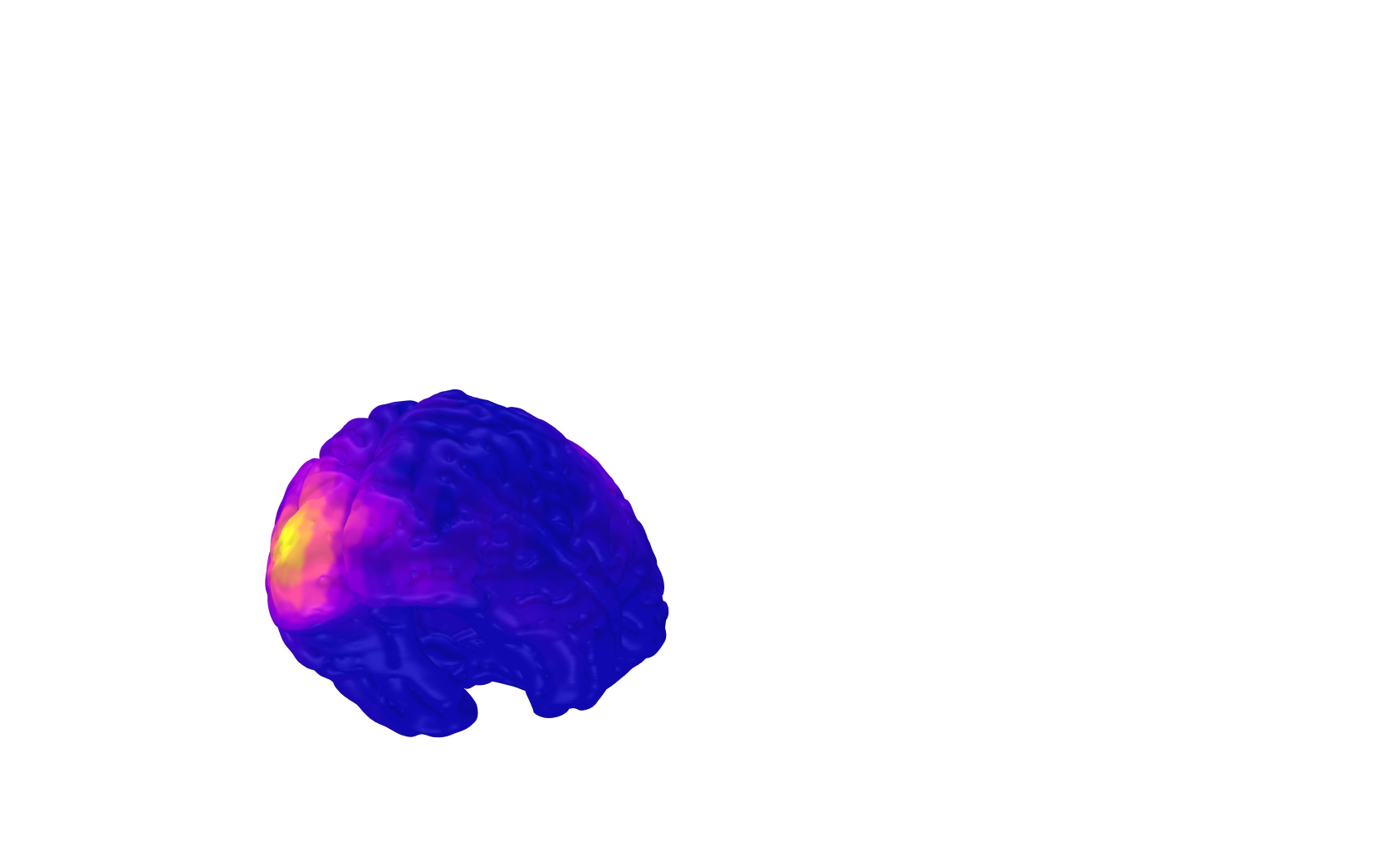}
        \end{minipage}\begin{minipage}{0.16\linewidth}
            \includegraphics[trim={5cm 1cm 14cm 8cm},clip,width=\linewidth]{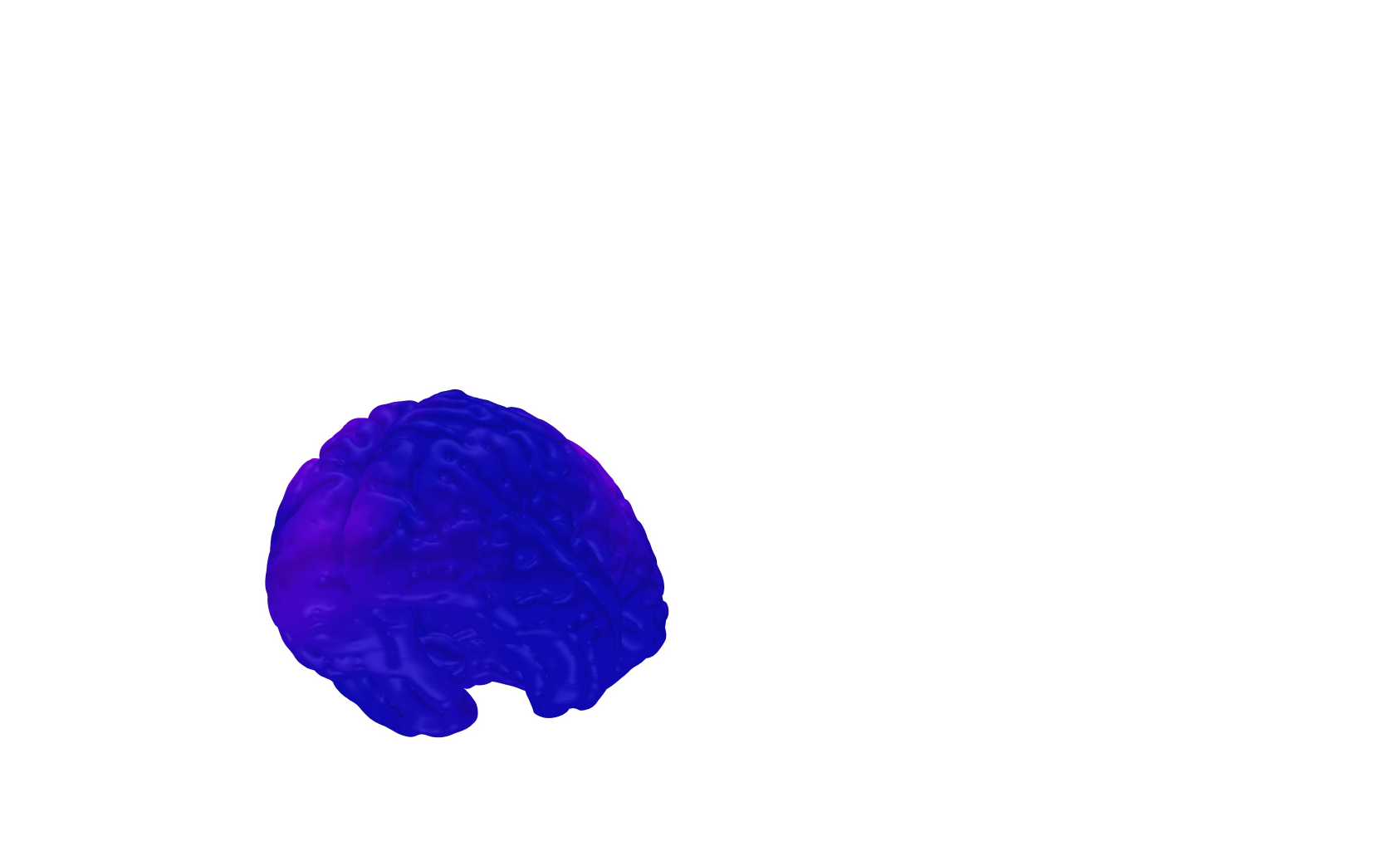}
        \end{minipage}\begin{minipage}{0.16\linewidth}
            \includegraphics[trim={5cm 1cm 14cm 8cm},clip,width=\linewidth]{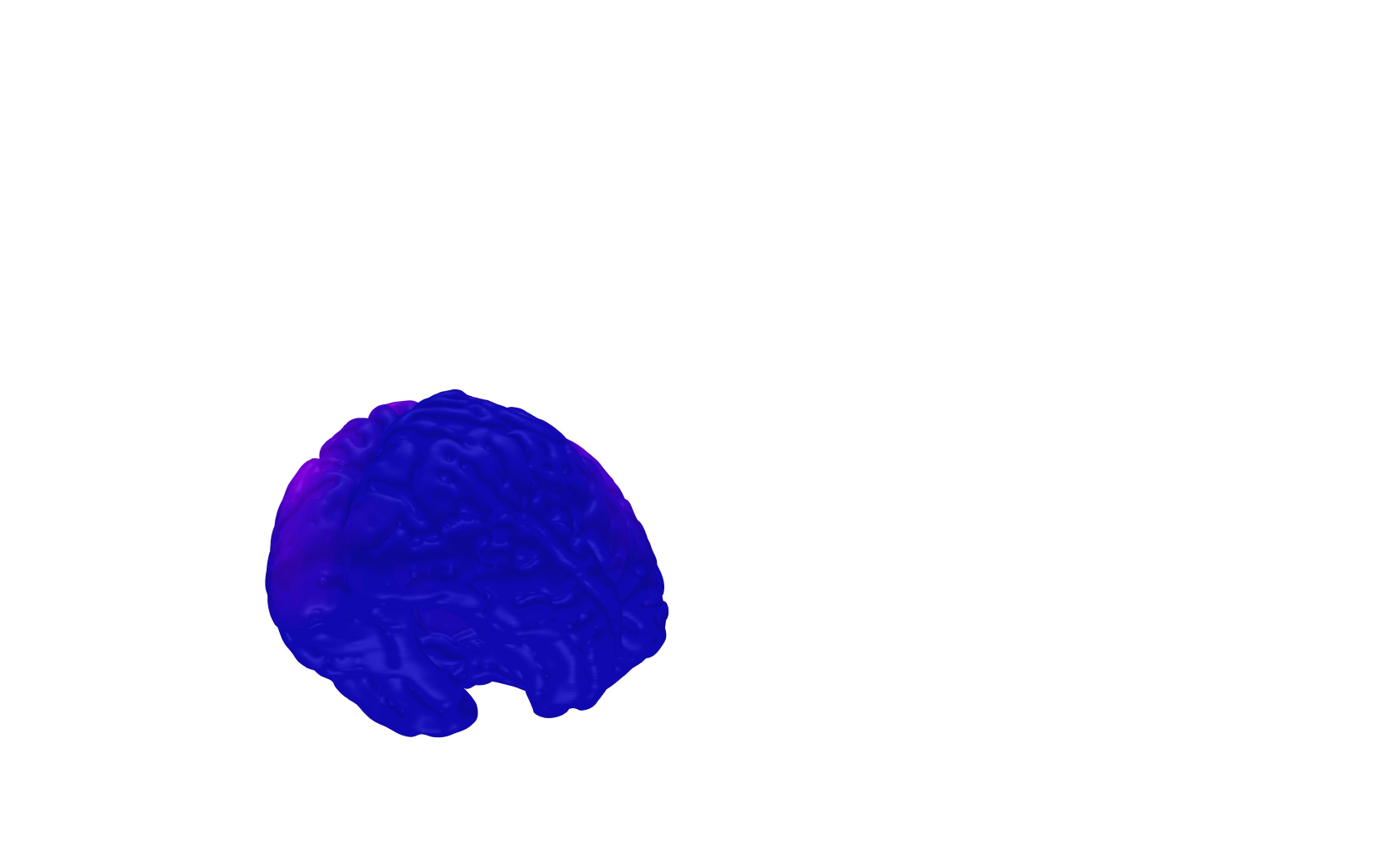}
        \end{minipage}

        \begin{minipage}{0.16\linewidth}
            \includegraphics[trim={5cm 1cm 14cm 8cm},clip,width=\linewidth]{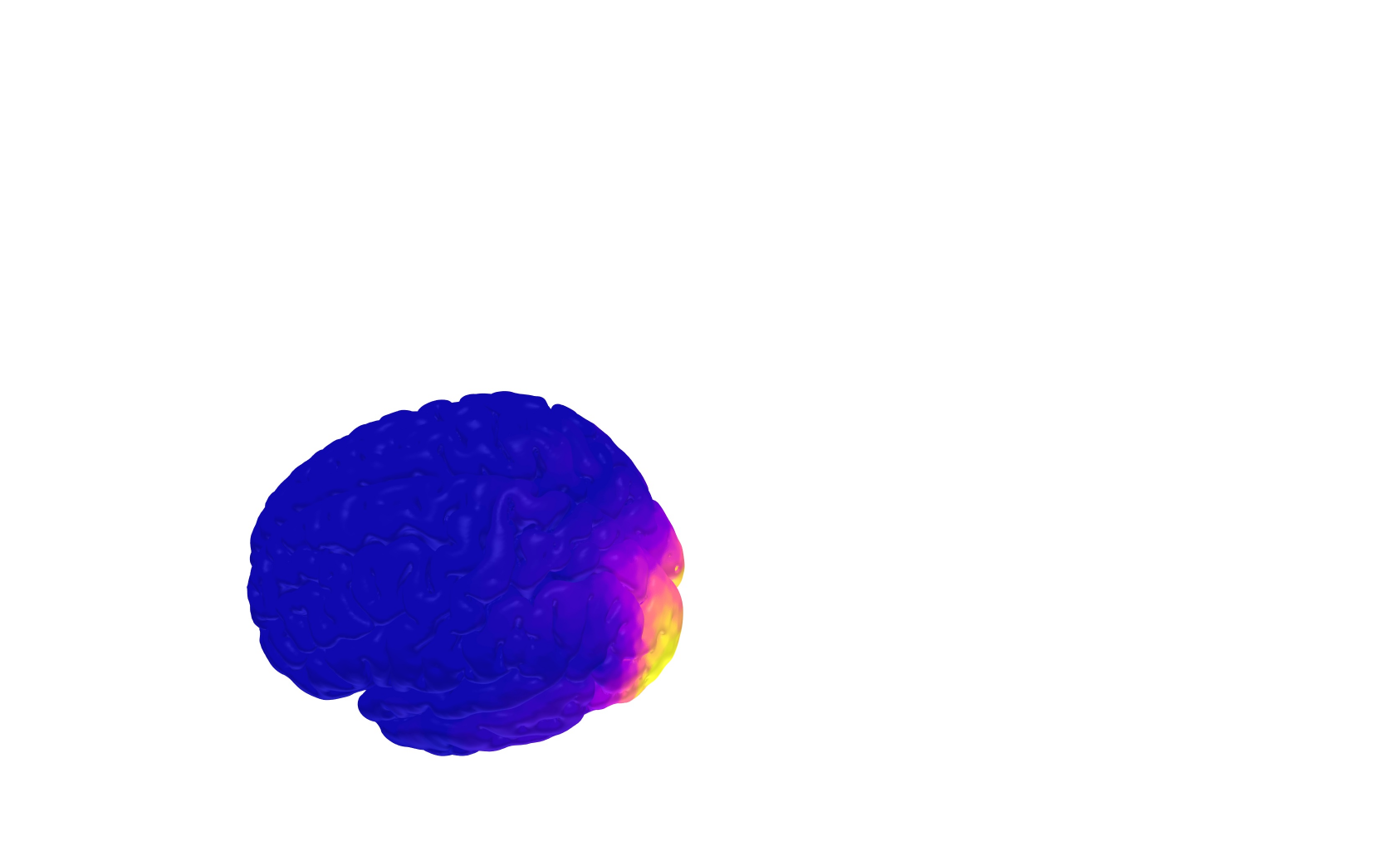}
        \end{minipage}\begin{minipage}{0.16\linewidth}
            \includegraphics[trim={5cm 1cm 14cm 8cm},clip,width=\linewidth]{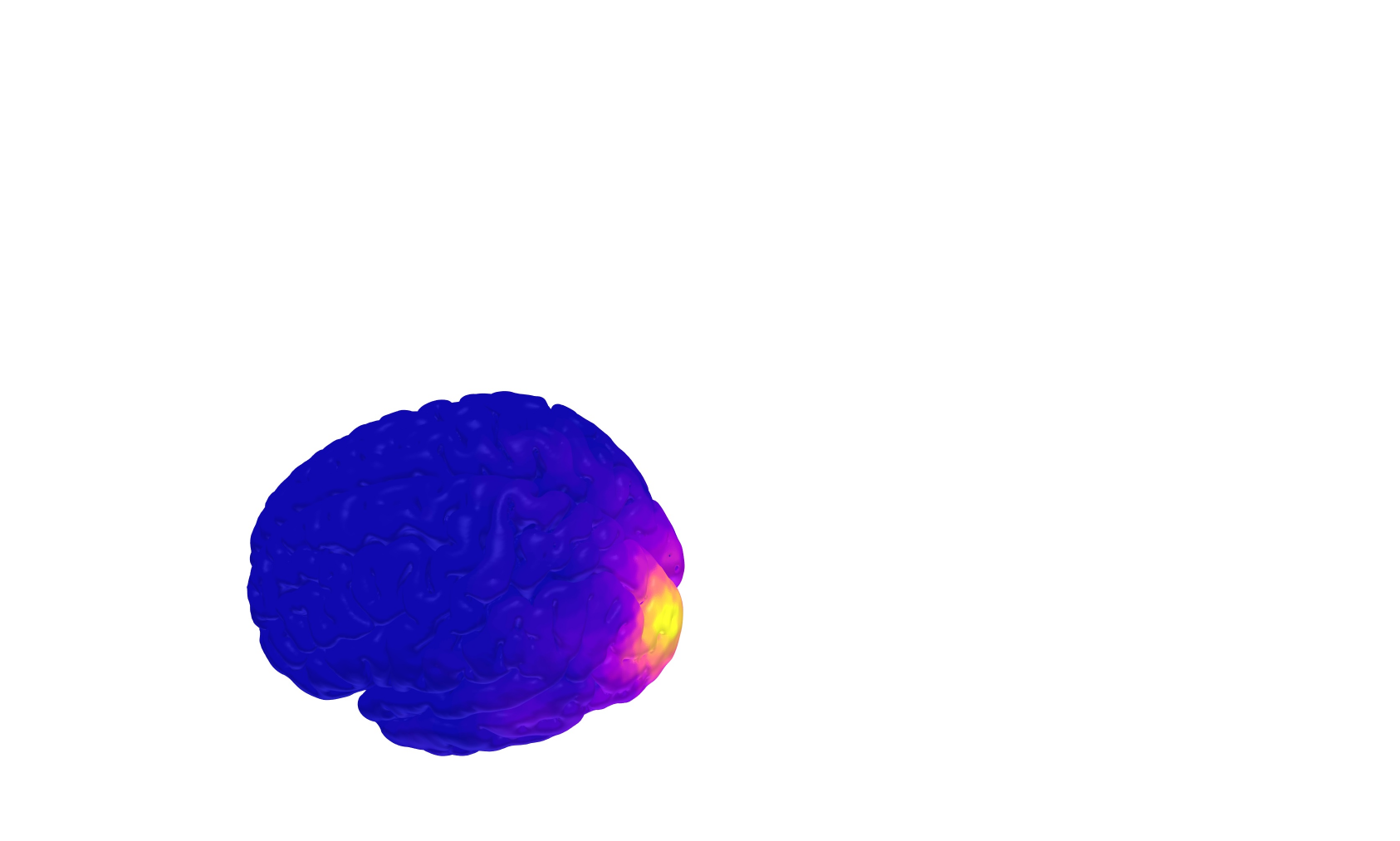}
        \end{minipage}\begin{minipage}{0.16\linewidth}
            \includegraphics[trim={5cm 1cm 14cm 8cm},clip,width=\linewidth]{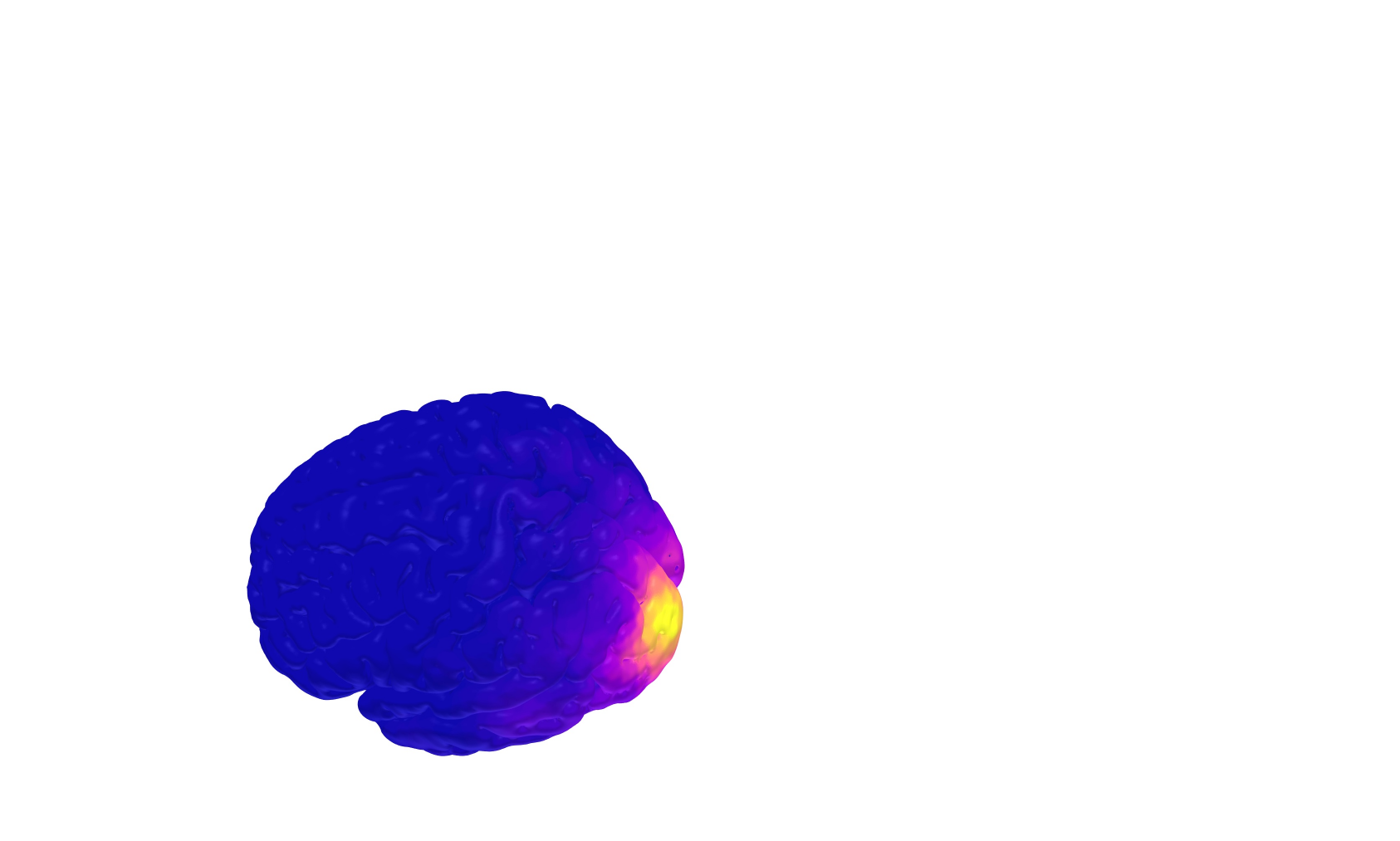}
        \end{minipage}\begin{minipage}{0.16\linewidth}
            \includegraphics[trim={5cm 1cm 14cm 8cm},clip,width=\linewidth]{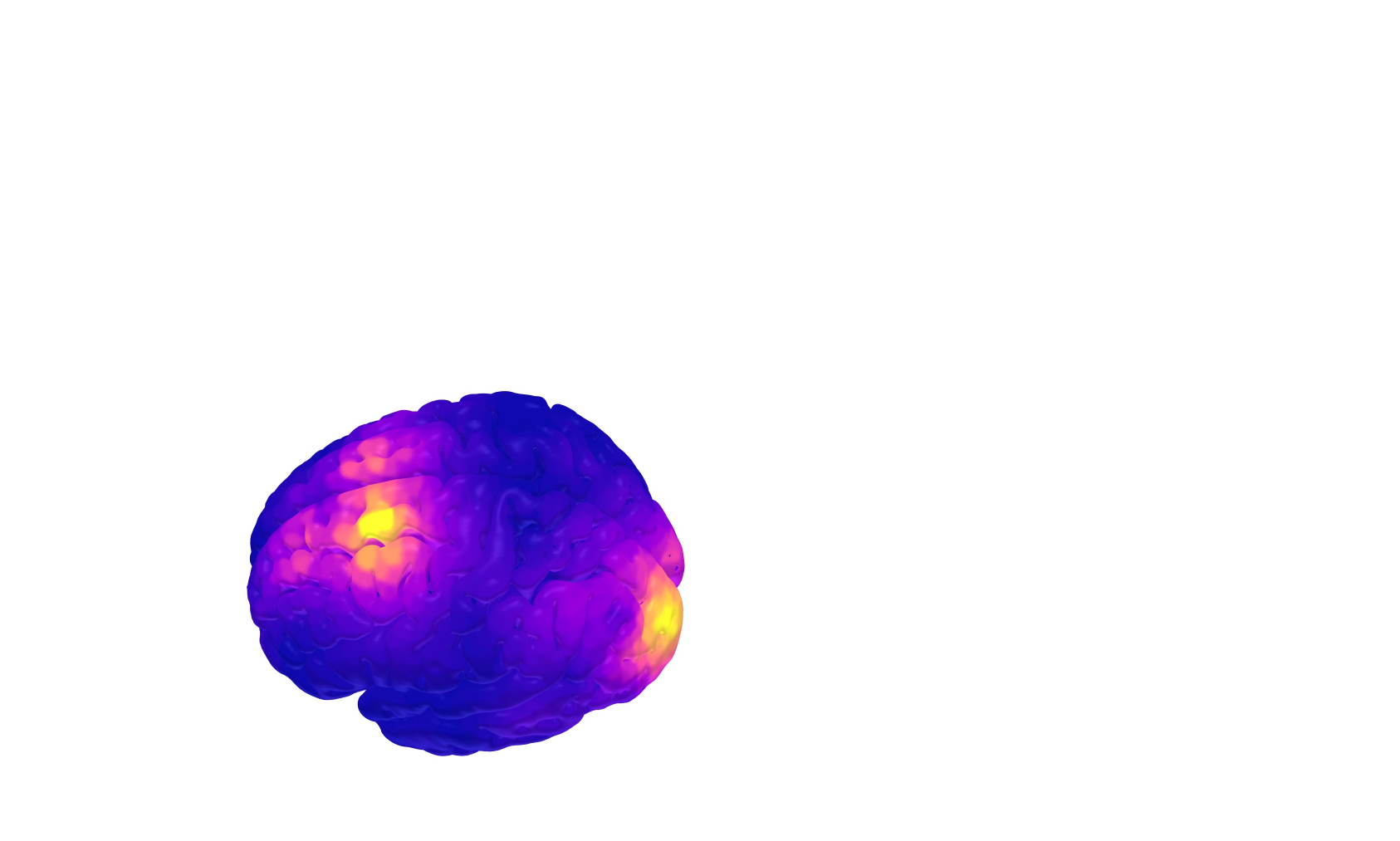}
        \end{minipage}\begin{minipage}{0.16\linewidth}
            \includegraphics[trim={5cm 1cm 14cm 8cm},clip,width=\linewidth]{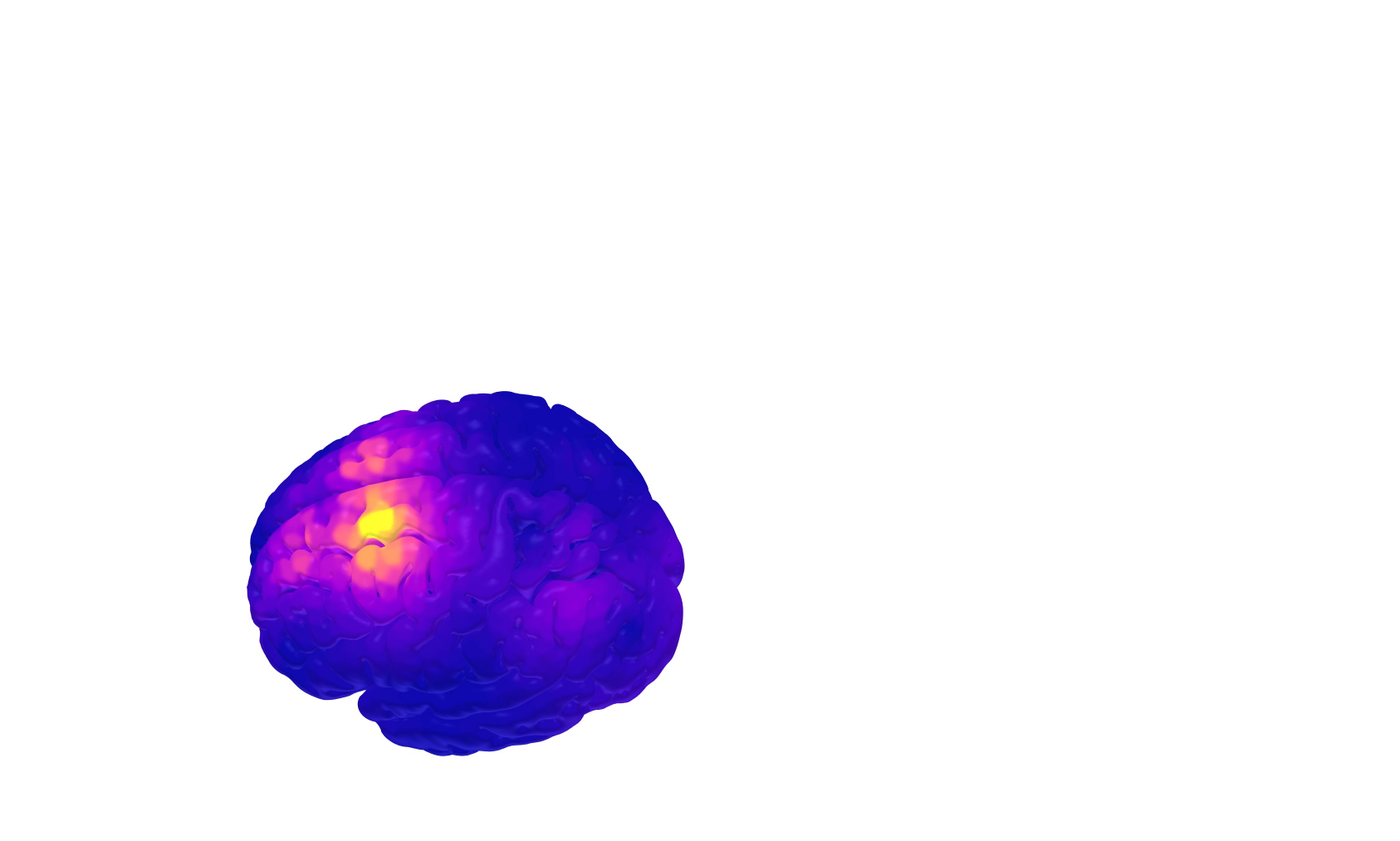}
        \end{minipage}\begin{minipage}{0.16\linewidth}
            \includegraphics[trim={5cm 1cm 14cm 8cm},clip,width=\linewidth]{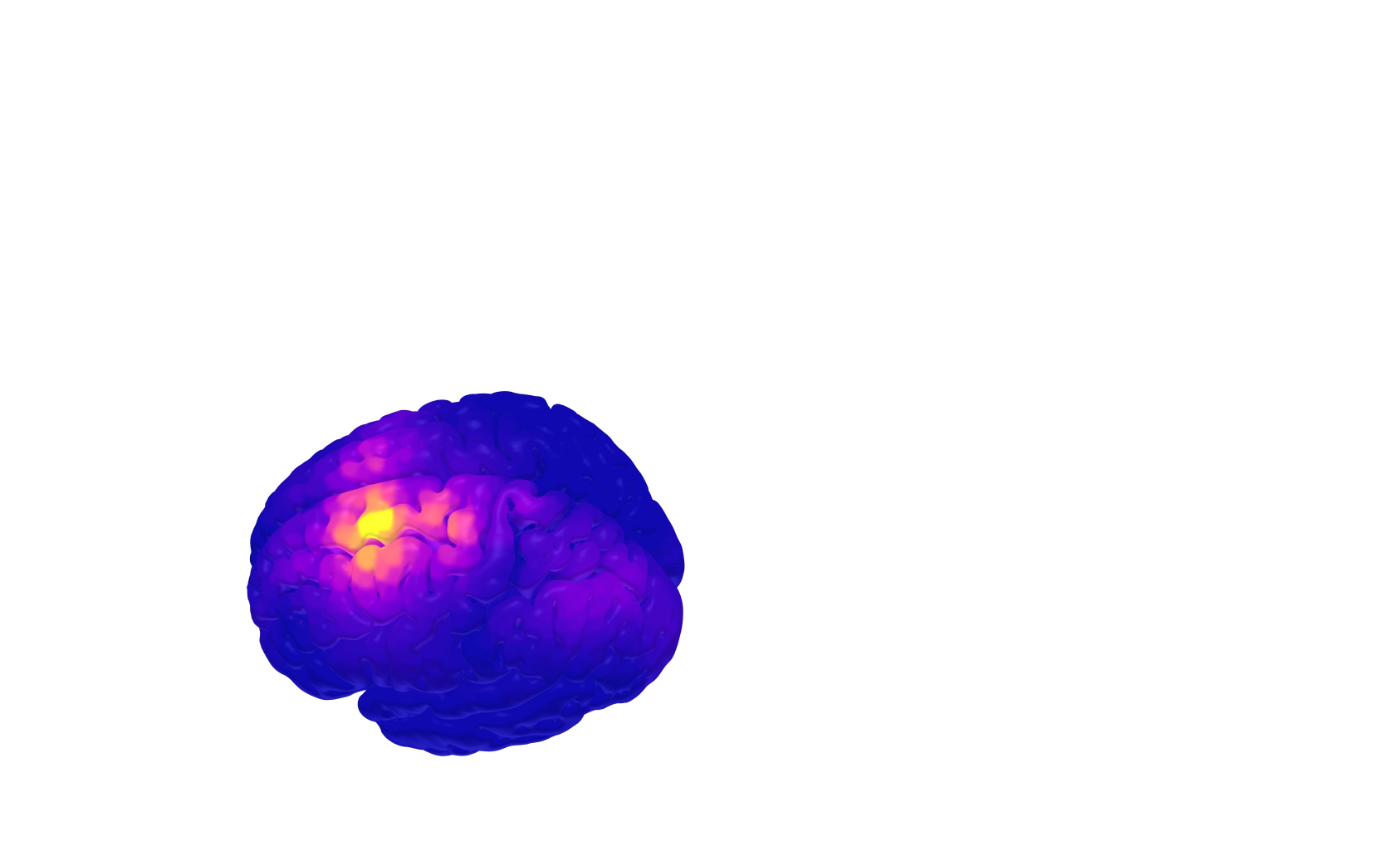}
        \end{minipage}
    \end{minipage}
    \noindent\rule[7pt]{\linewidth}{0.4pt}

    \begin{minipage}{0.02\linewidth}
        \rotatebox{90}{\scriptsize KF}
    \end{minipage}\begin{minipage}{0.98\linewidth}
        \begin{minipage}{0.16\linewidth}
            \includegraphics[trim={5cm 1cm 14cm 8cm},clip,width=\linewidth]{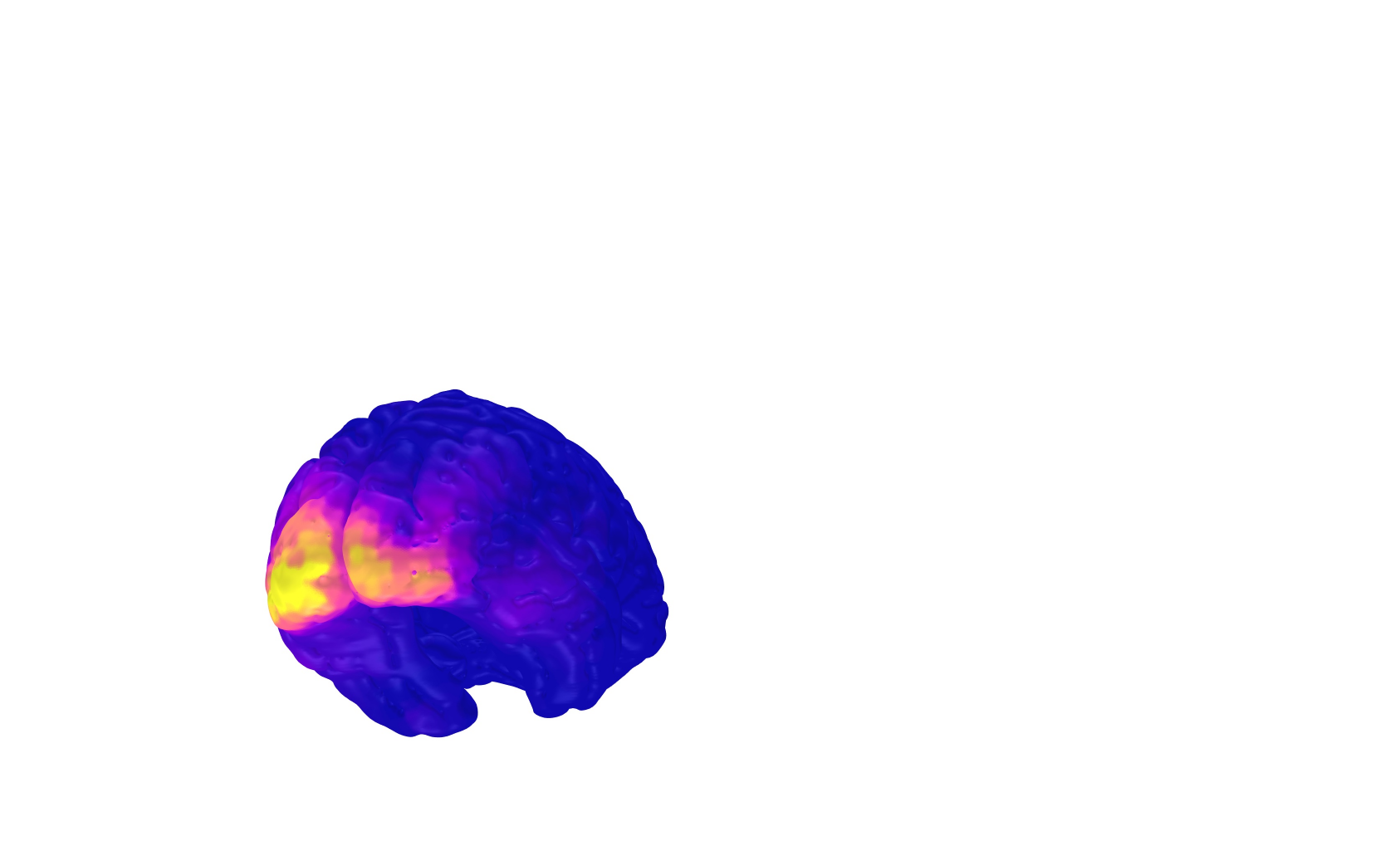}
        \end{minipage}\begin{minipage}{0.16\linewidth}
            \includegraphics[trim={5cm 1cm 14cm 8cm},clip,width=\linewidth]{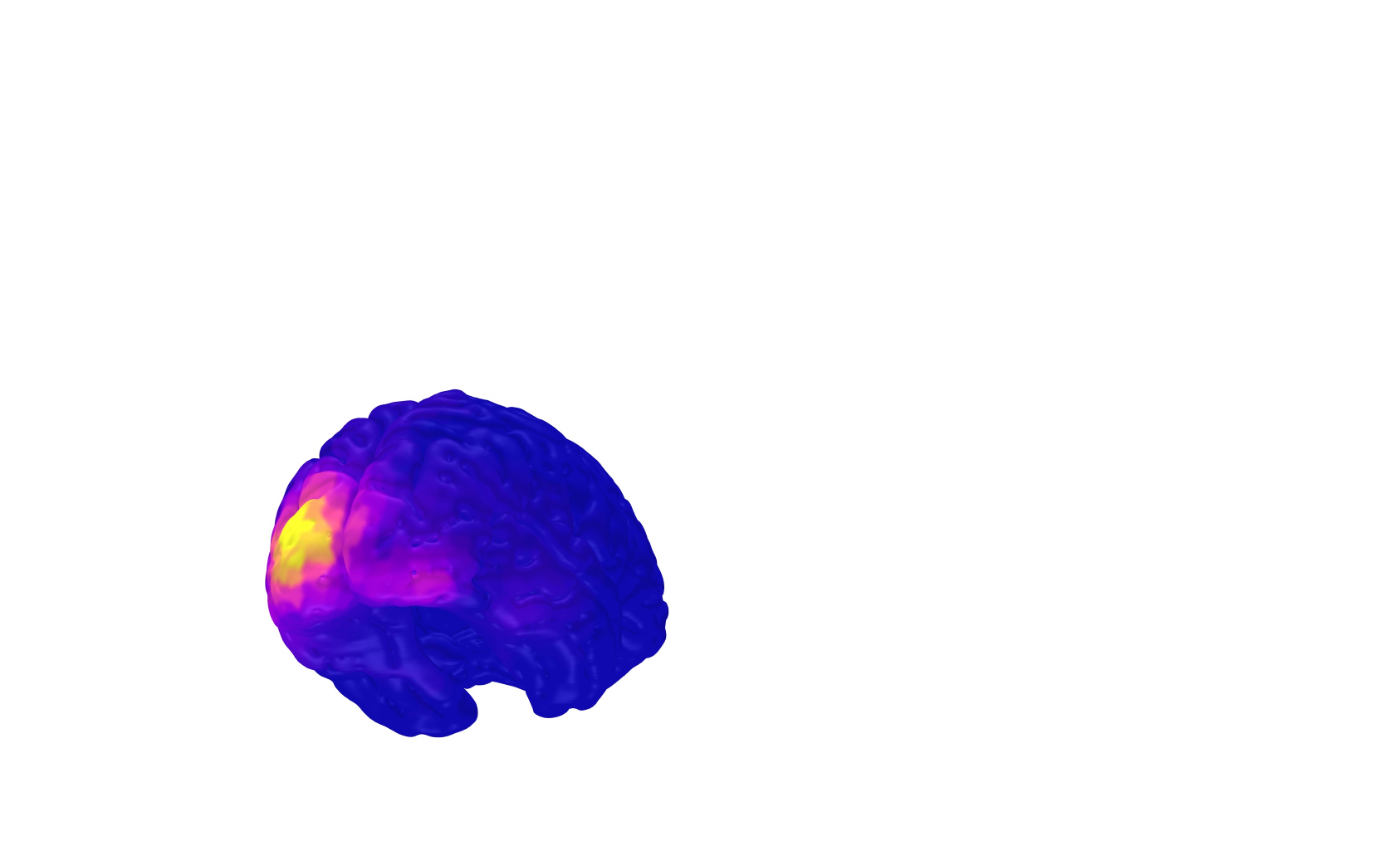}
        \end{minipage}\begin{minipage}{0.16\linewidth}
            \includegraphics[trim={5cm 1cm 14cm 8cm},clip,width=\linewidth]{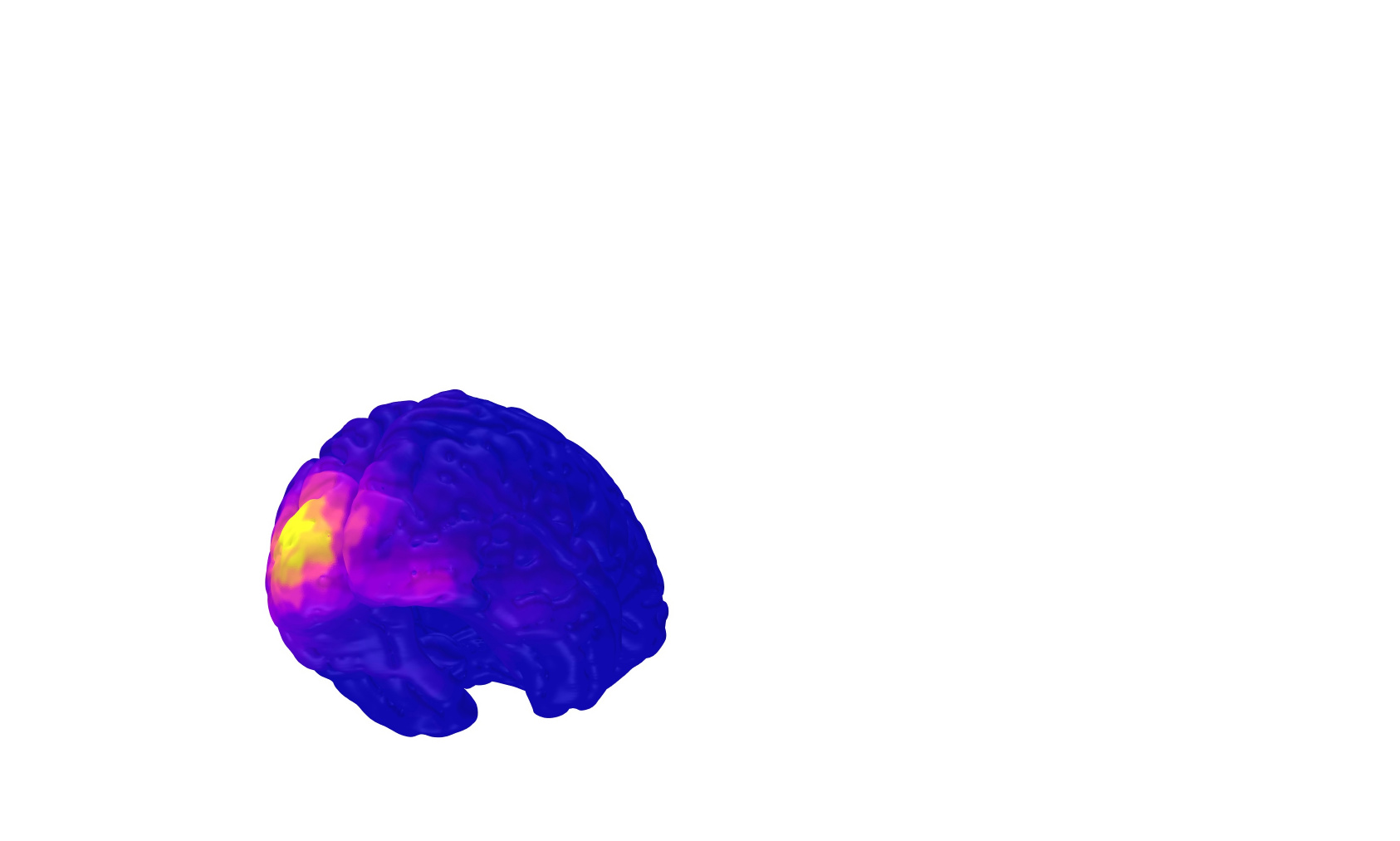}
        \end{minipage}\begin{minipage}{0.16\linewidth}
            \includegraphics[trim={5cm 1cm 14cm 8cm},clip,width=\linewidth]{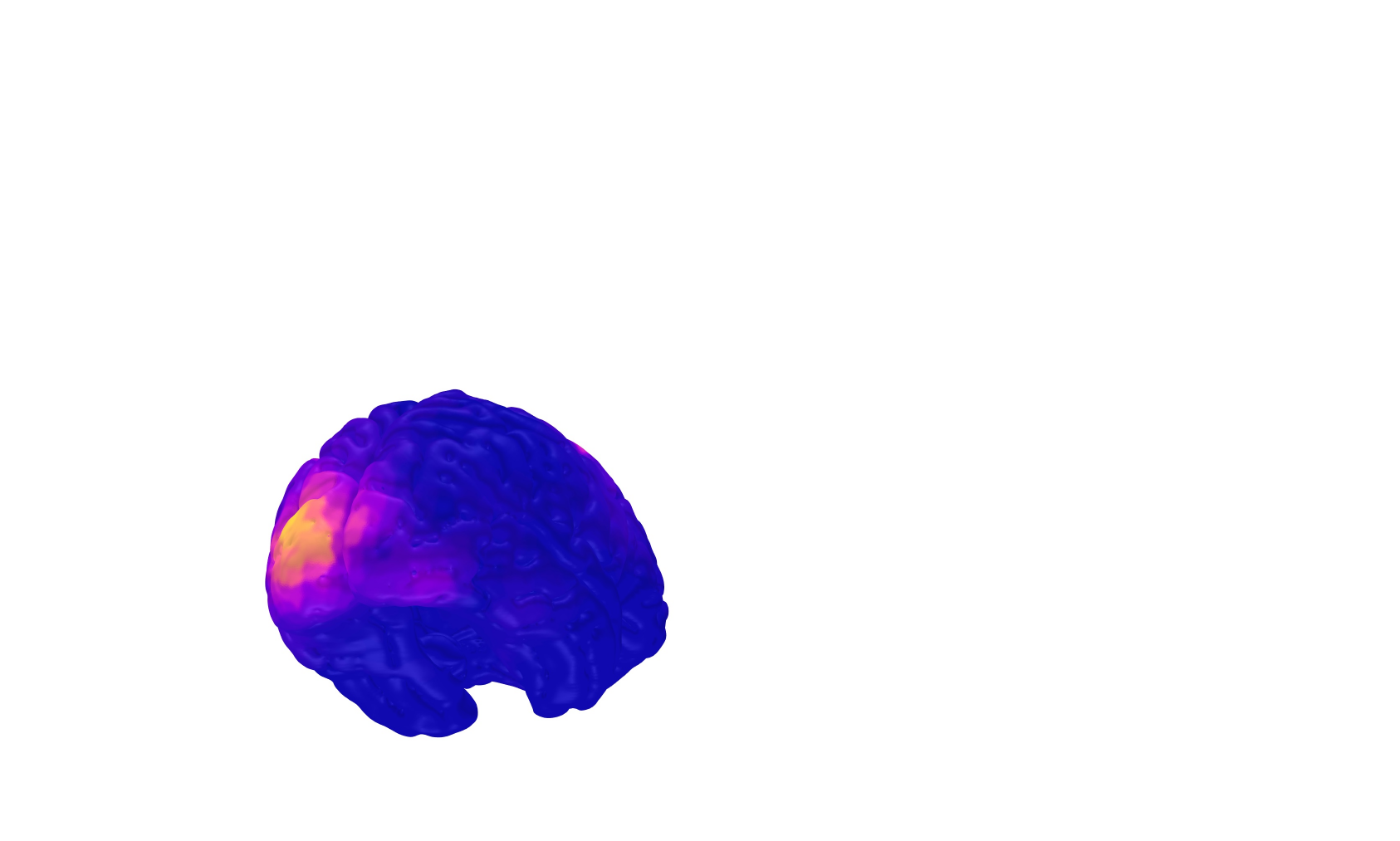}
        \end{minipage}\begin{minipage}{0.16\linewidth}
            \includegraphics[trim={5cm 1cm 14cm 8cm},clip,width=\linewidth]{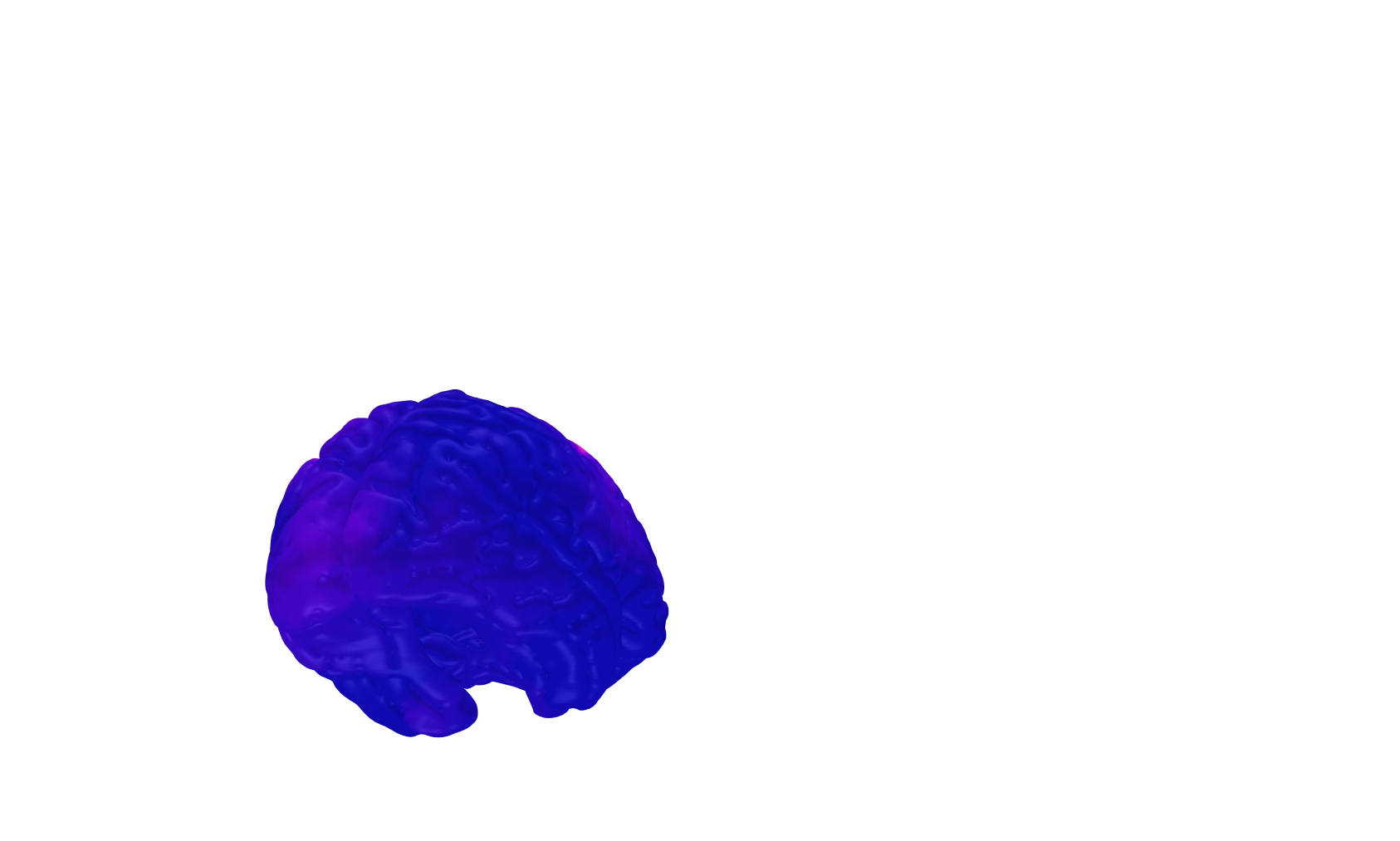}
        \end{minipage}\begin{minipage}{0.16\linewidth}
            \includegraphics[trim={5cm 1cm 14cm 8cm},clip,width=\linewidth]{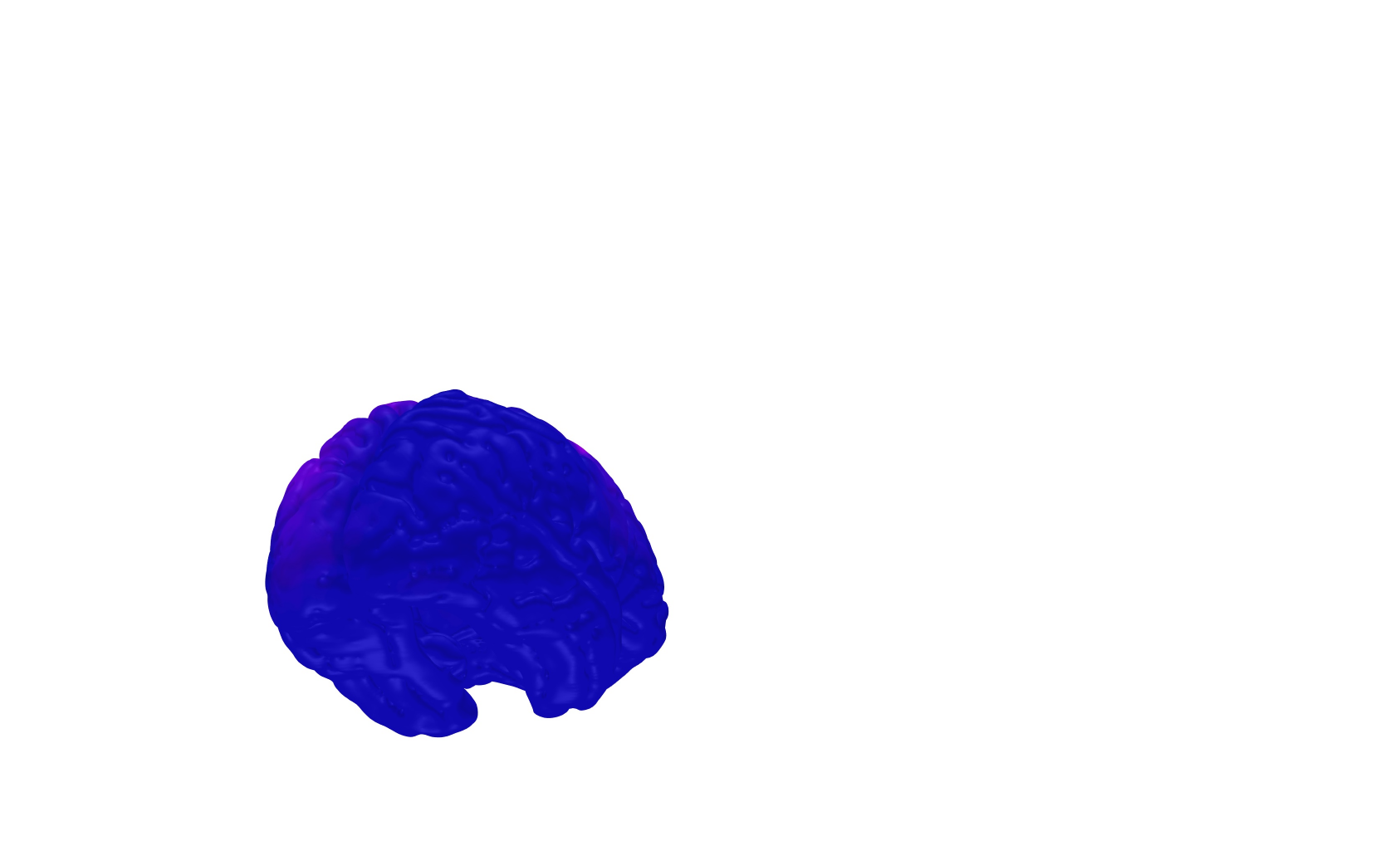}
        \end{minipage}

        \begin{minipage}{0.16\linewidth}
            \includegraphics[trim={5cm 1cm 14cm 8cm},clip,width=\linewidth]{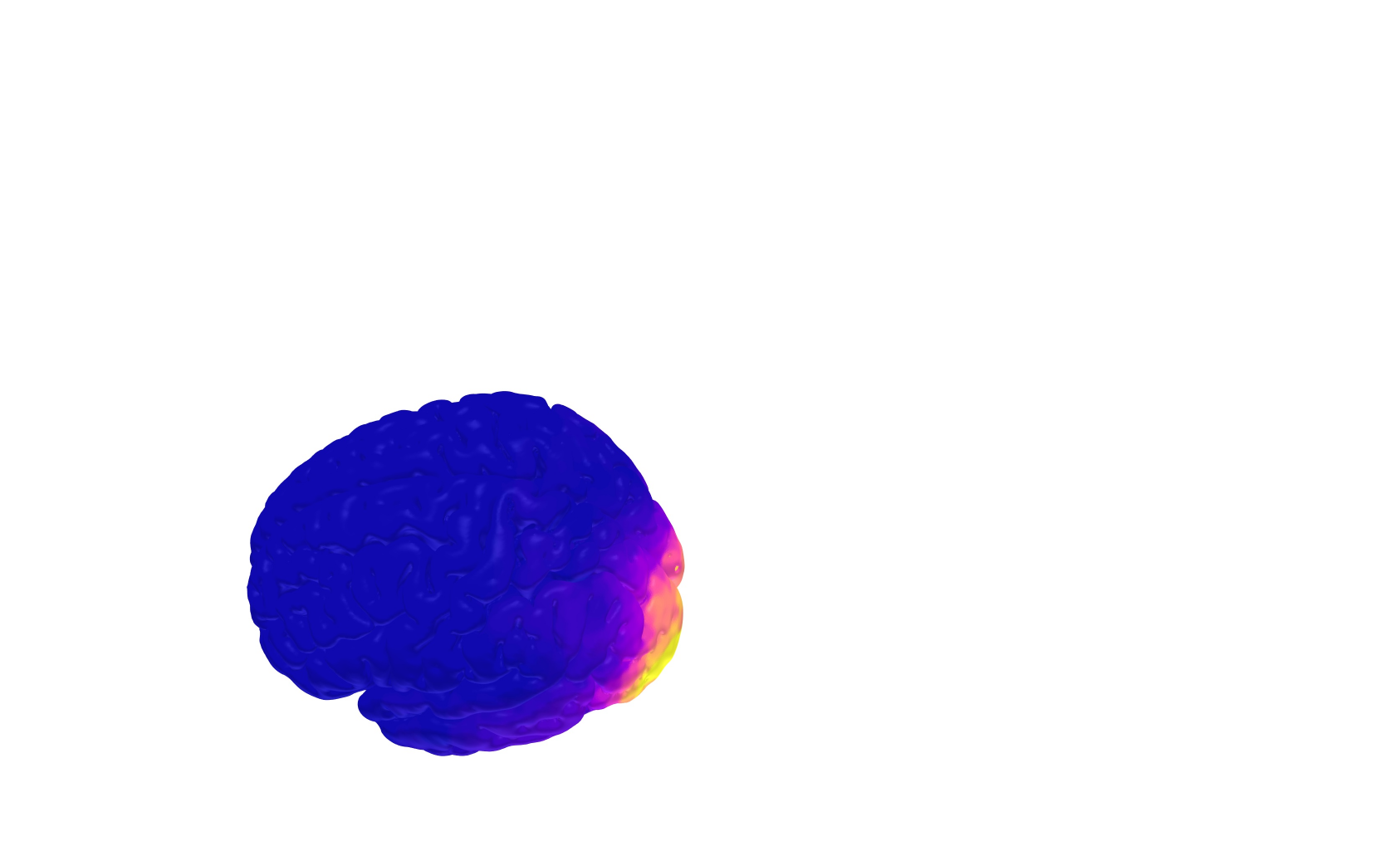}
        \end{minipage}\begin{minipage}{0.16\linewidth}
            \includegraphics[trim={5cm 1cm 14cm 8cm},clip,width=\linewidth]{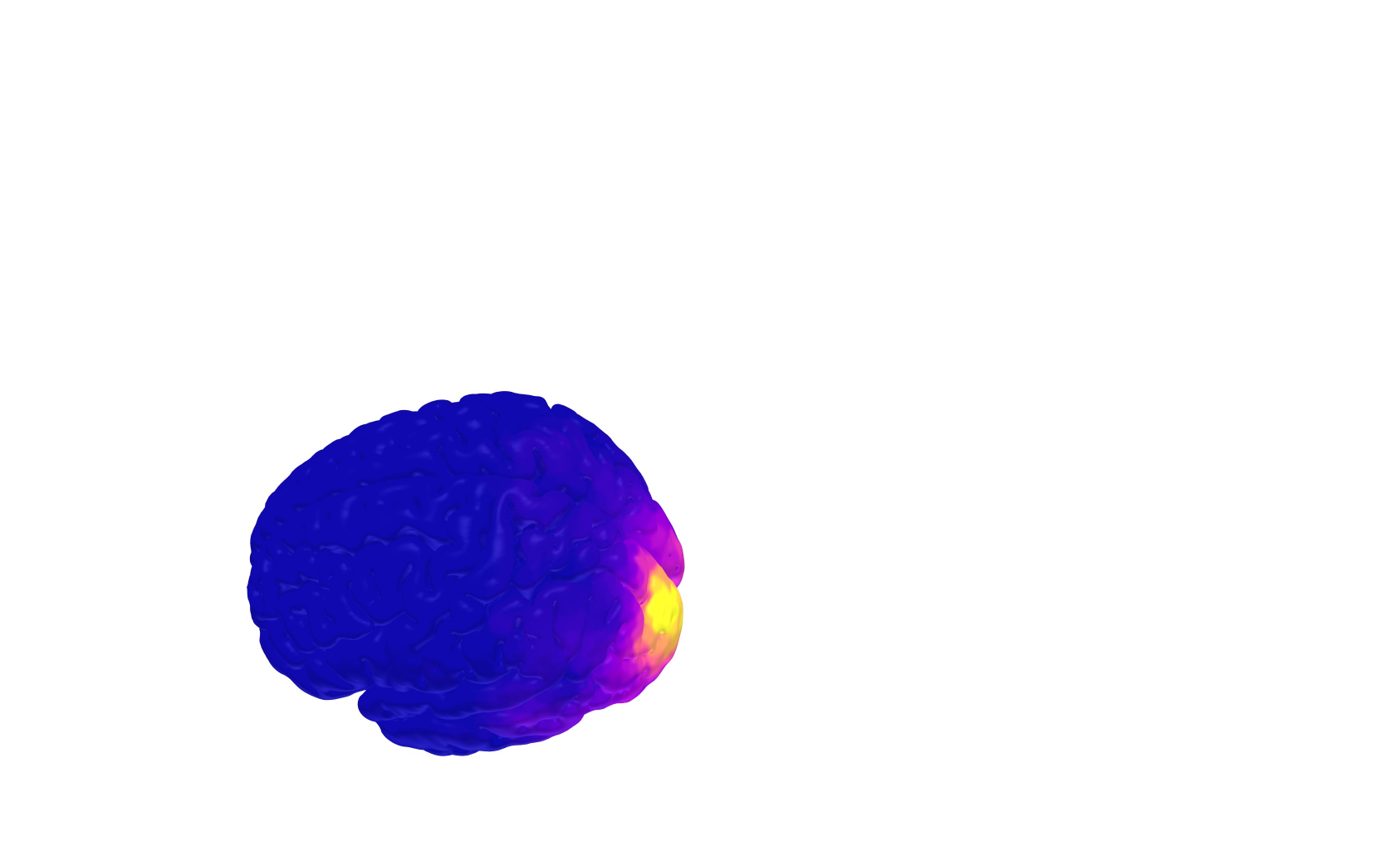}
        \end{minipage}\begin{minipage}{0.16\linewidth}
            \includegraphics[trim={5cm 1cm 14cm 8cm},clip,width=\linewidth]{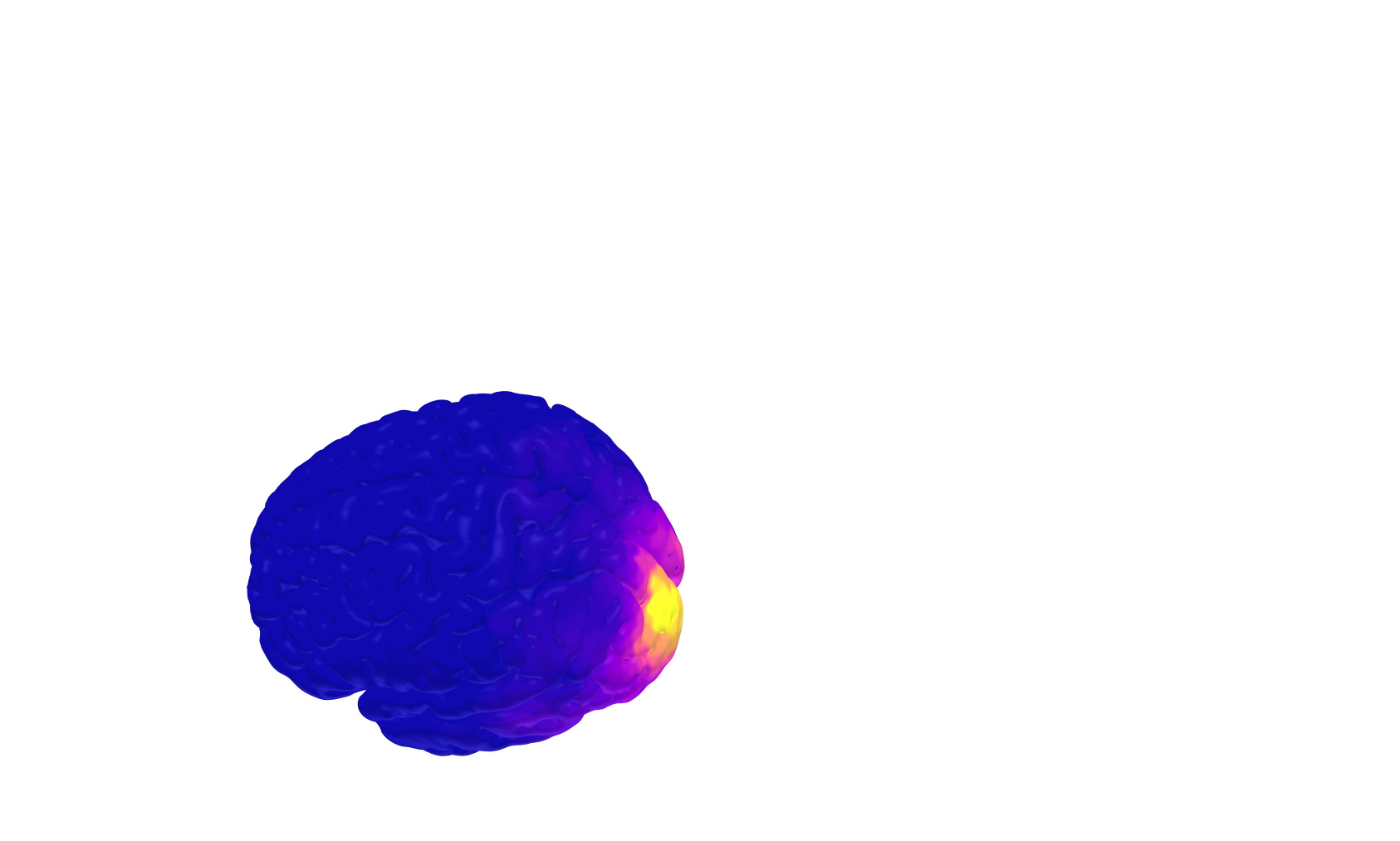}
        \end{minipage}\begin{minipage}{0.16\linewidth}
            \includegraphics[trim={5cm 1cm 14cm 8cm},clip,width=\linewidth]{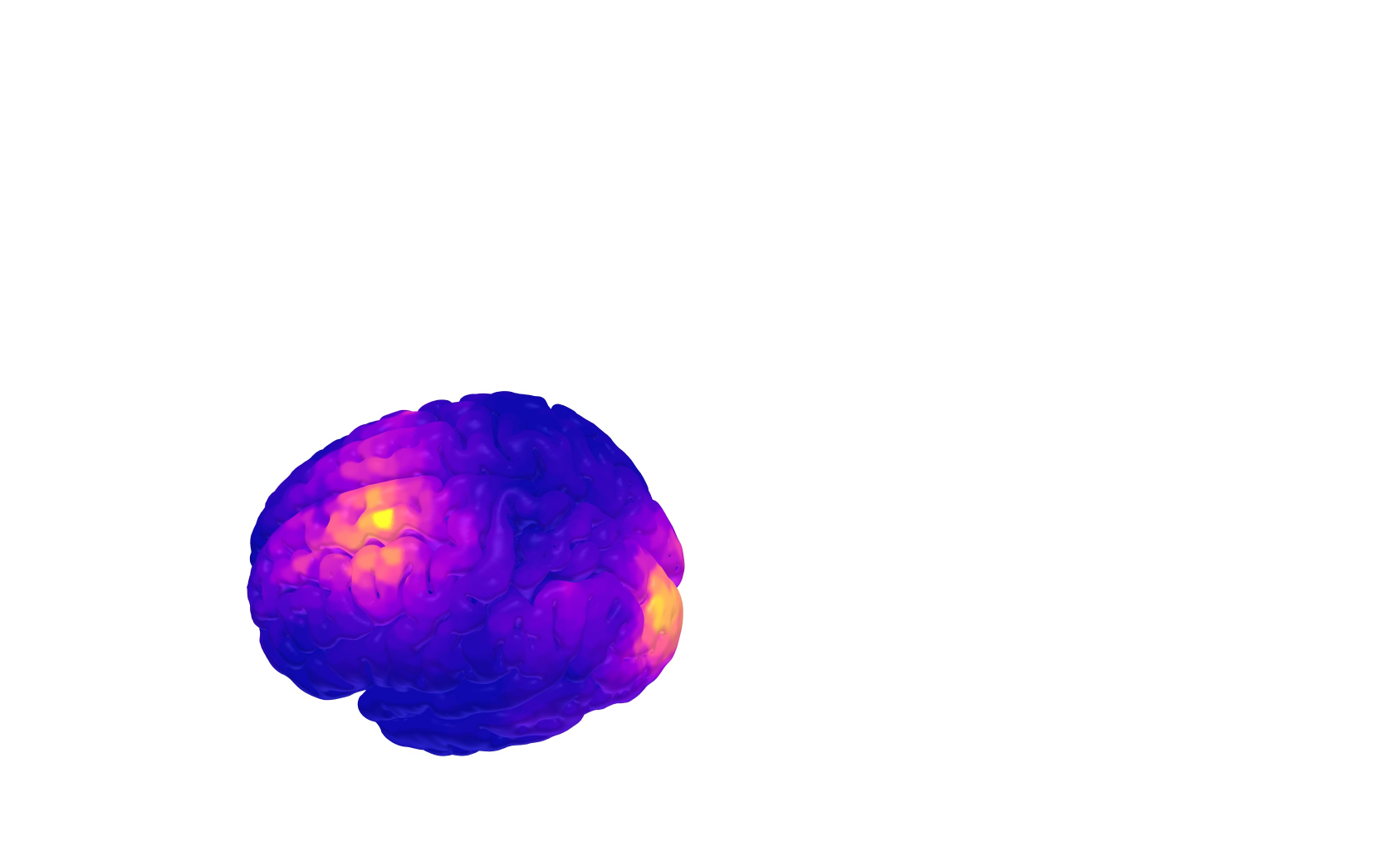}
        \end{minipage}\begin{minipage}{0.16\linewidth}
            \includegraphics[trim={5cm 1cm 14cm 8cm},clip,width=\linewidth]{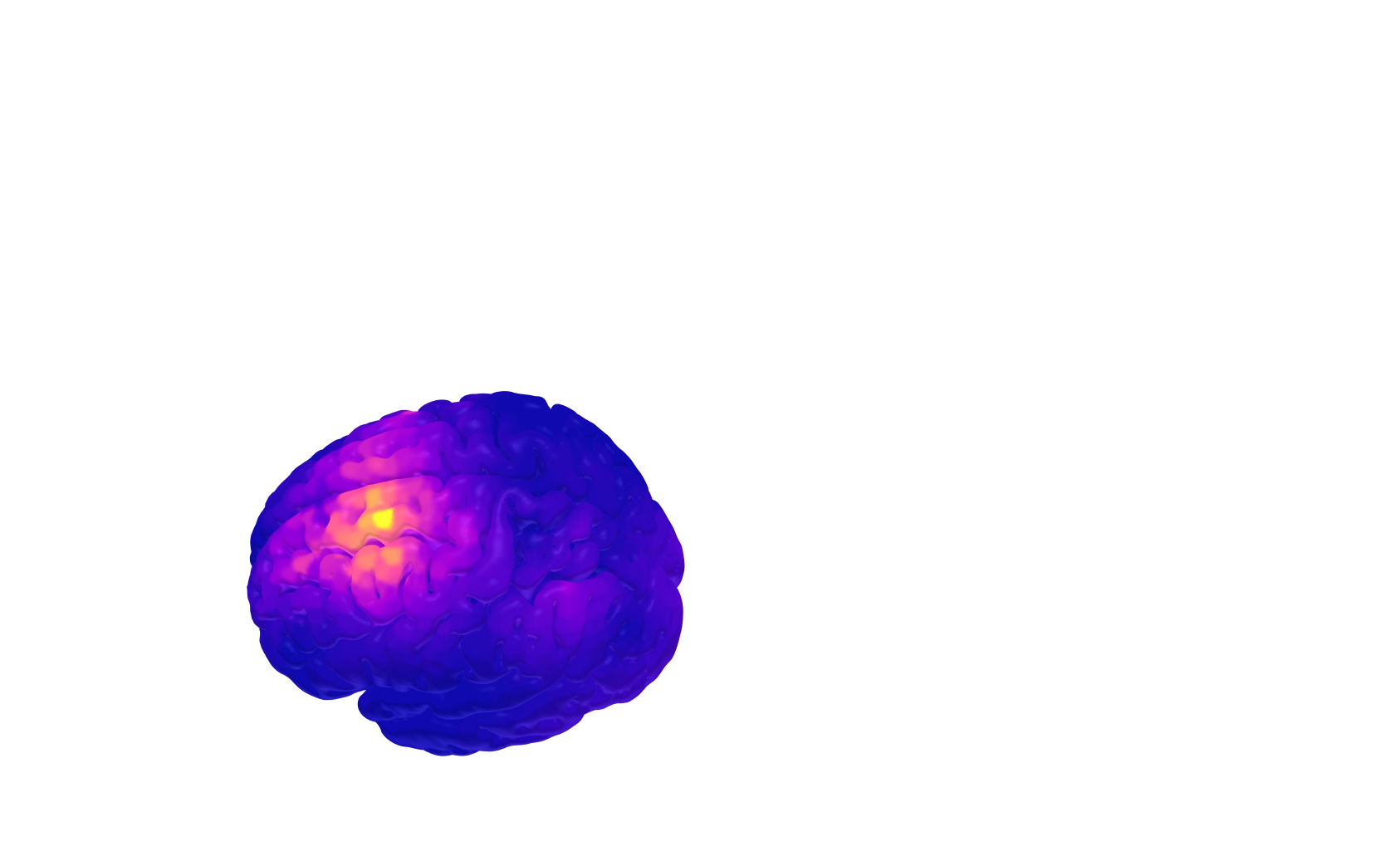}
        \end{minipage}\begin{minipage}{0.16\linewidth}
            \includegraphics[trim={5cm 1cm 14cm 8cm},clip,width=\linewidth]{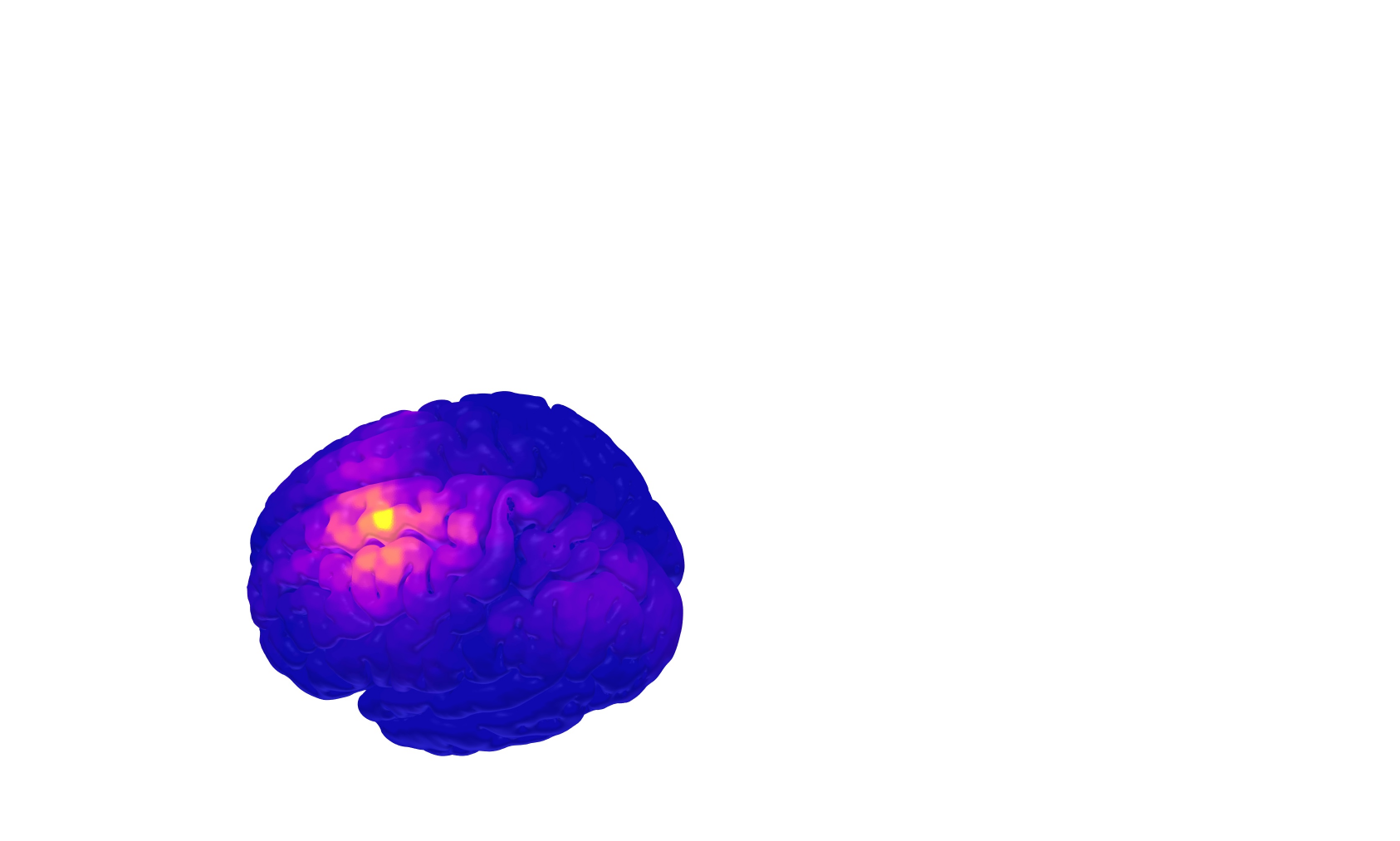}
        \end{minipage}
    \end{minipage}
    \noindent\rule[7pt]{\linewidth}{0.4pt}

    \begin{minipage}{0.02\linewidth}
        \rotatebox{90}{\scriptsize LCMV}
    \end{minipage}\begin{minipage}{0.98\linewidth}
        \begin{minipage}{0.16\linewidth}
            \includegraphics[trim={5cm 1cm 14cm 8cm},clip,width=\linewidth]{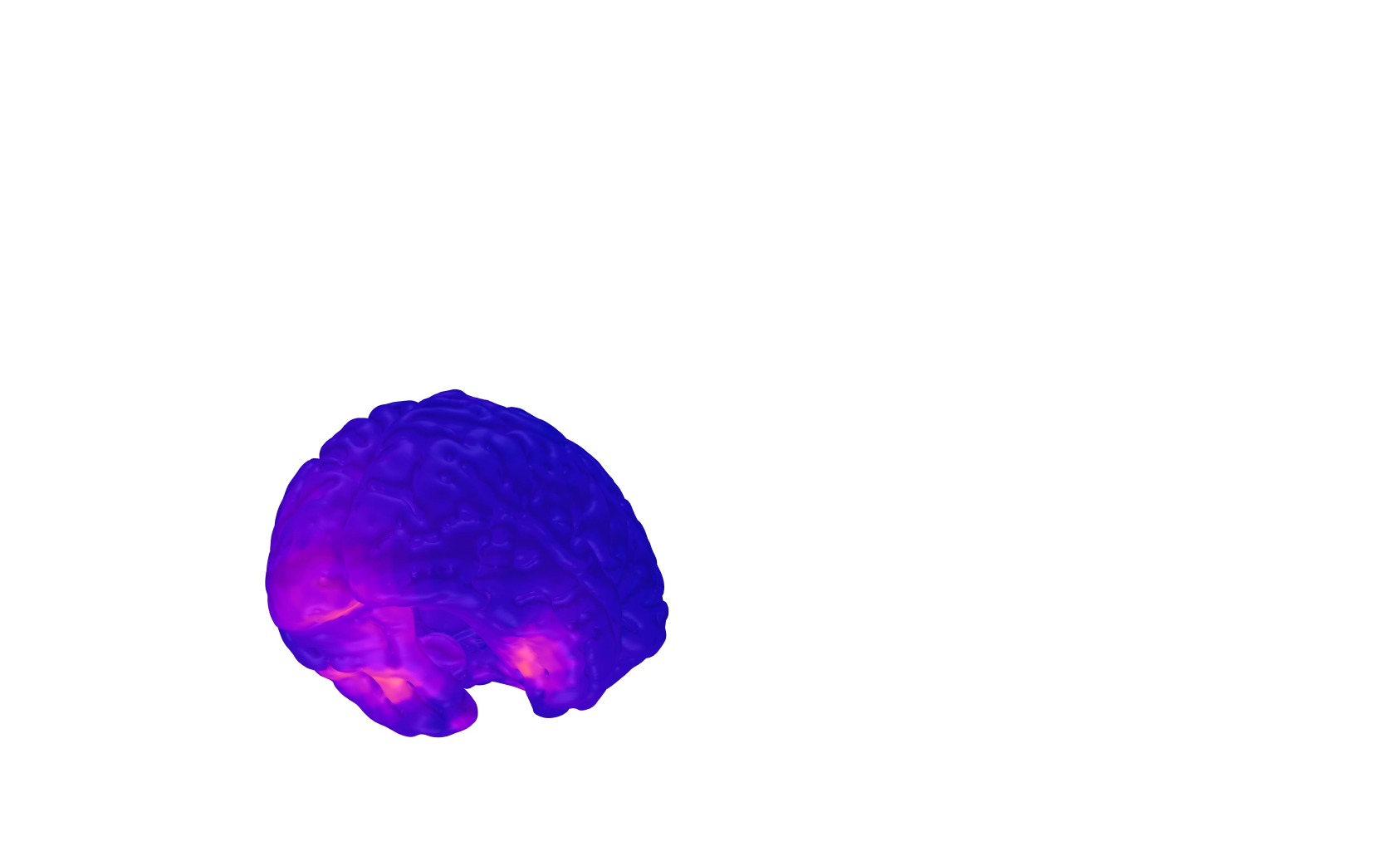}
        \end{minipage}\begin{minipage}{0.16\linewidth}
            \includegraphics[trim={5cm 1cm 14cm 8cm},clip,width=\linewidth]{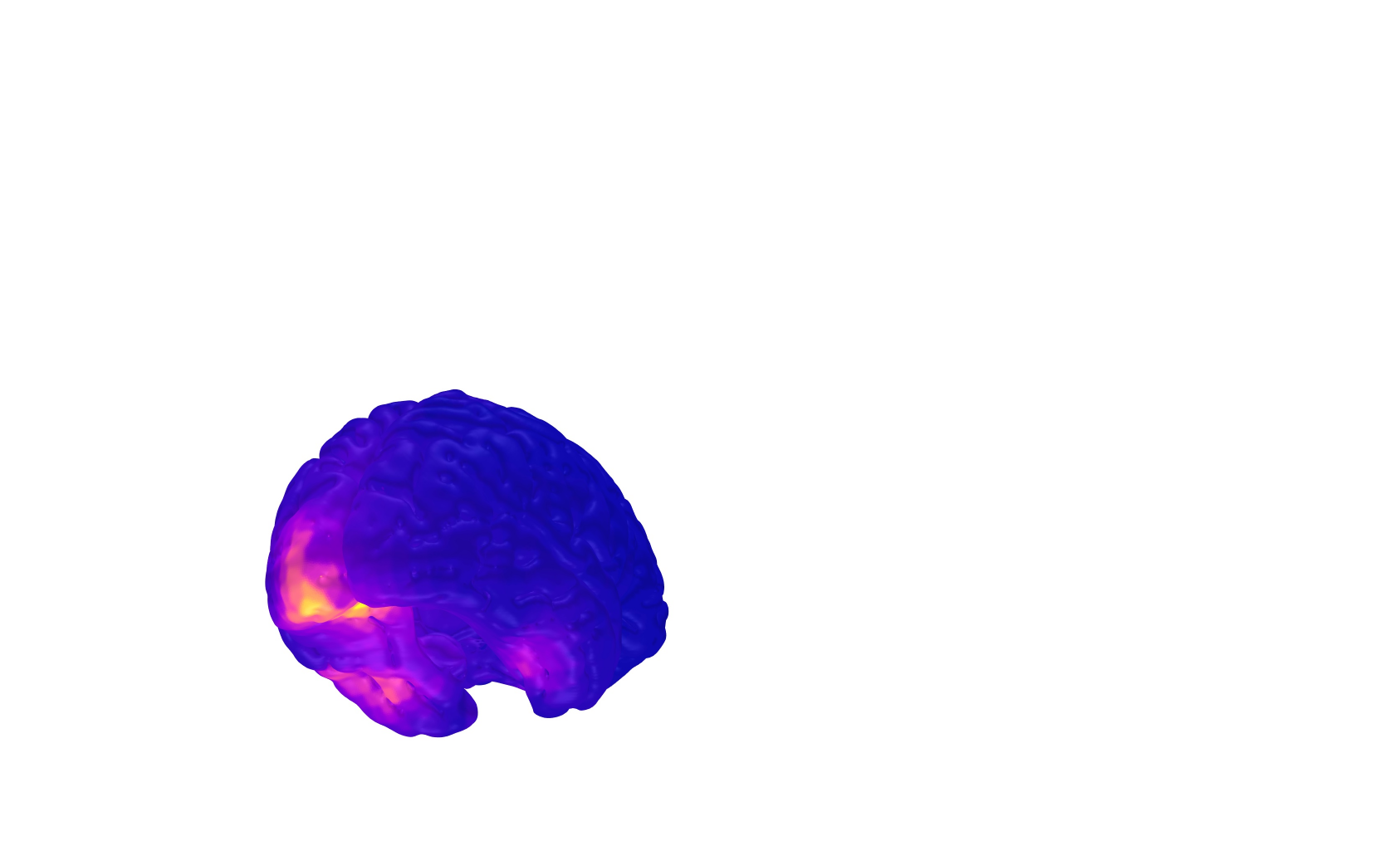}
        \end{minipage}\begin{minipage}{0.16\linewidth}
            \includegraphics[trim={5cm 1cm 14cm 8cm},clip,width=\linewidth]{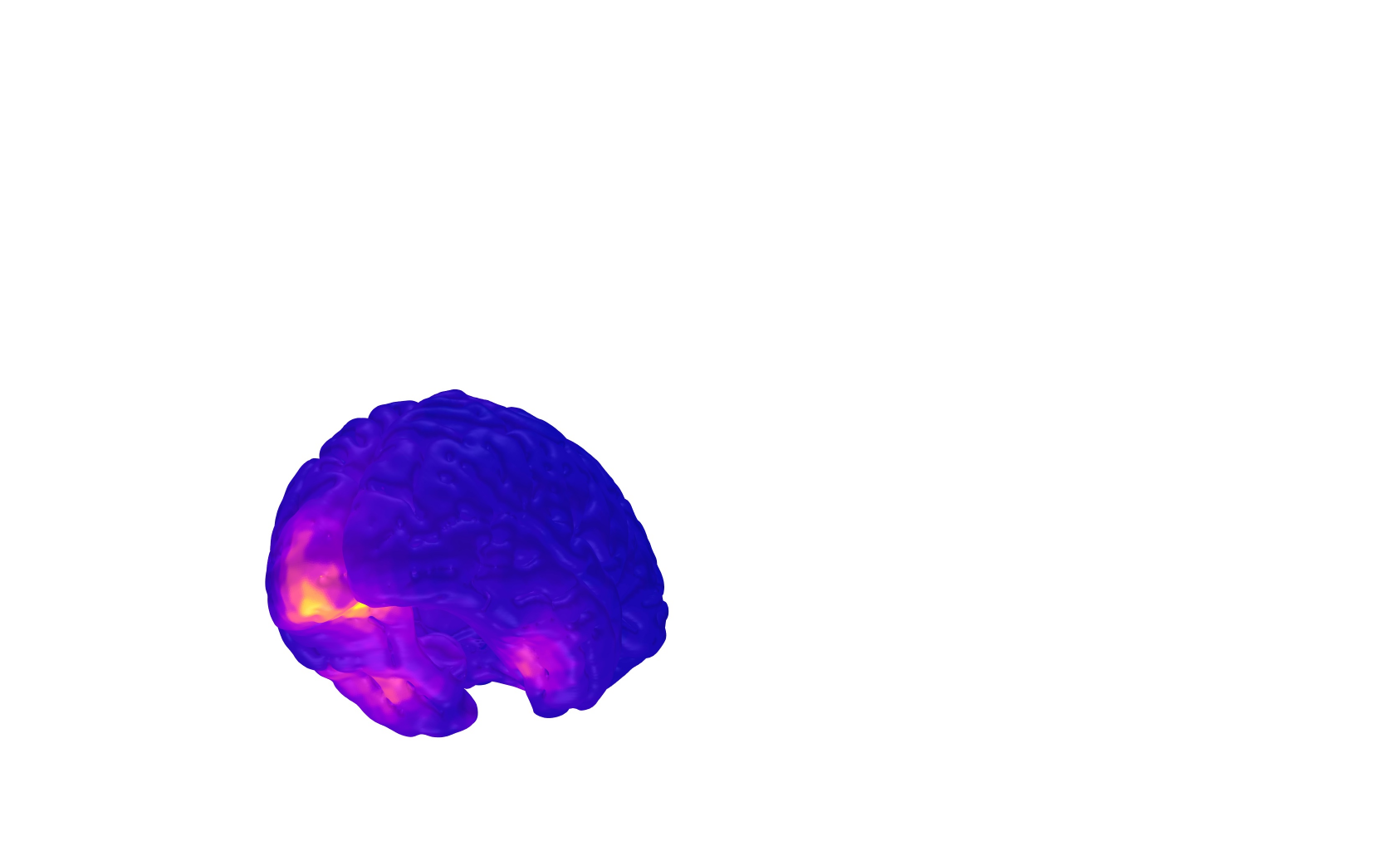}
        \end{minipage}\begin{minipage}{0.16\linewidth}
            \includegraphics[trim={5cm 1cm 14cm 8cm},clip,width=\linewidth]{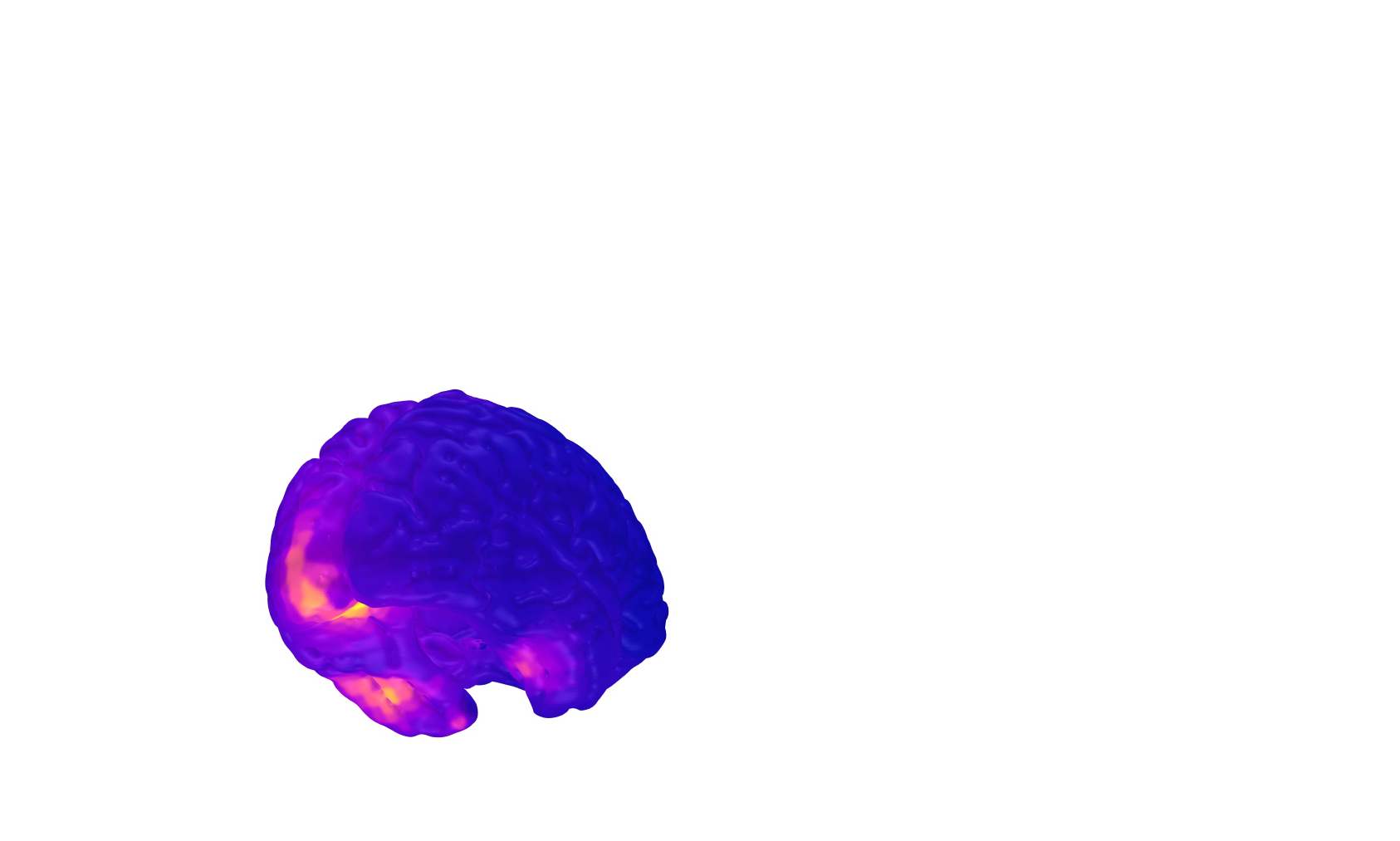}
        \end{minipage}\begin{minipage}{0.16\linewidth}
            \includegraphics[trim={5cm 1cm 14cm 8cm},clip,width=\linewidth]{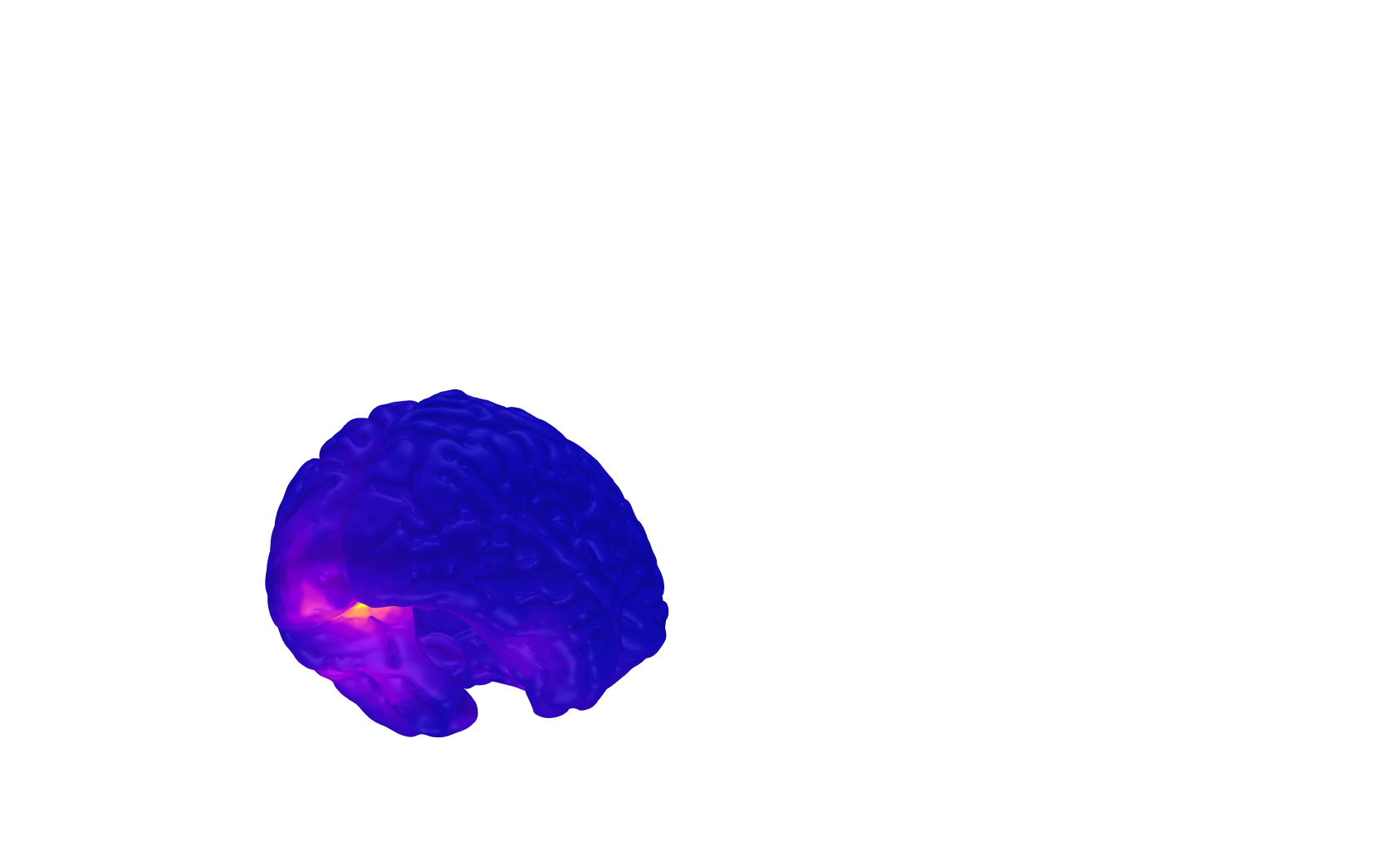}
        \end{minipage}\begin{minipage}{0.16\linewidth}
            \includegraphics[trim={5cm 1cm 14cm 8cm},clip,width=\linewidth]{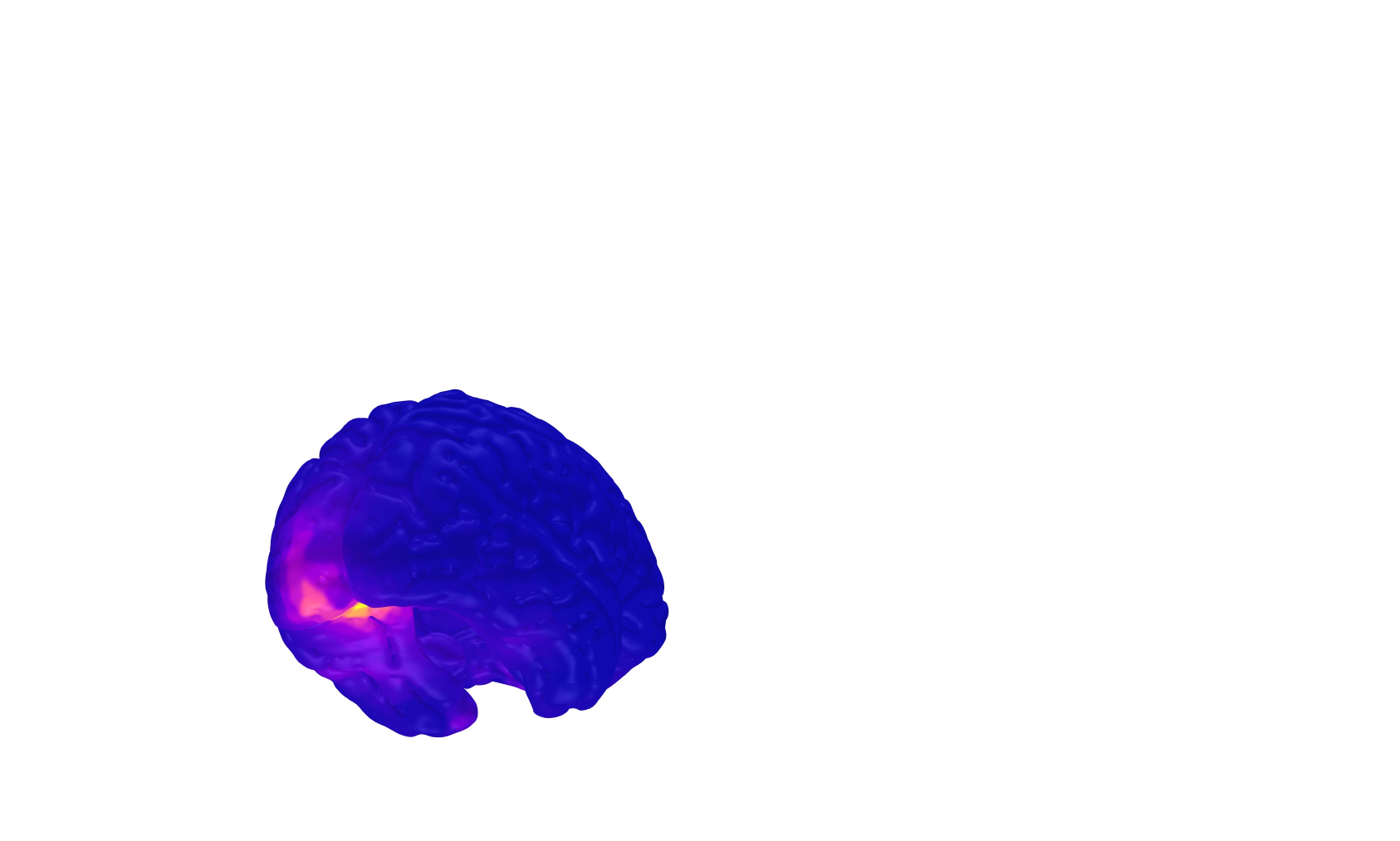}
        \end{minipage}

        \begin{minipage}{0.16\linewidth}
            \includegraphics[trim={5cm 1cm 14cm 8cm},clip,width=\linewidth]{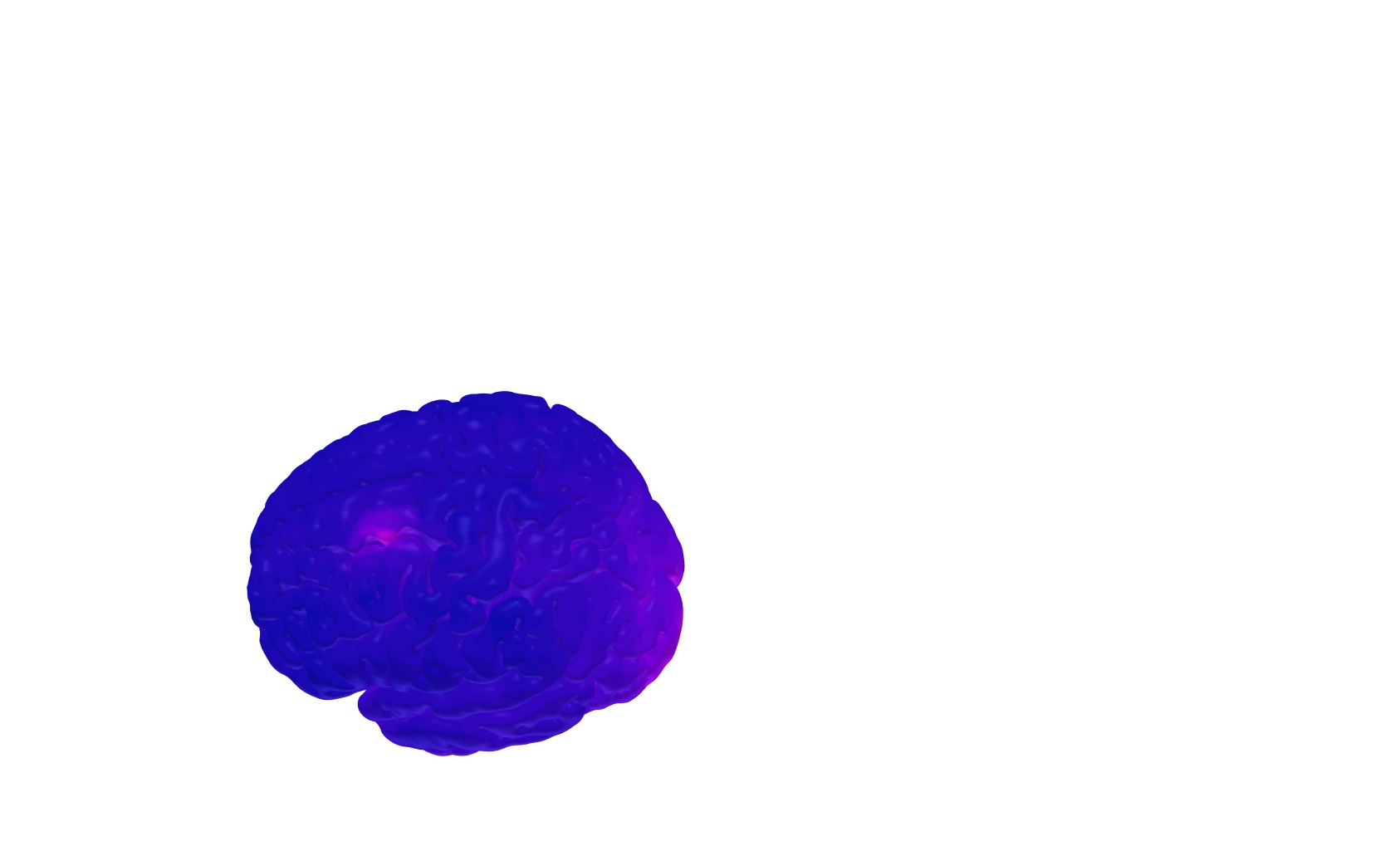}
        \end{minipage}\begin{minipage}{0.16\linewidth}
            \includegraphics[trim={5cm 1cm 14cm 8cm},clip,width=\linewidth]{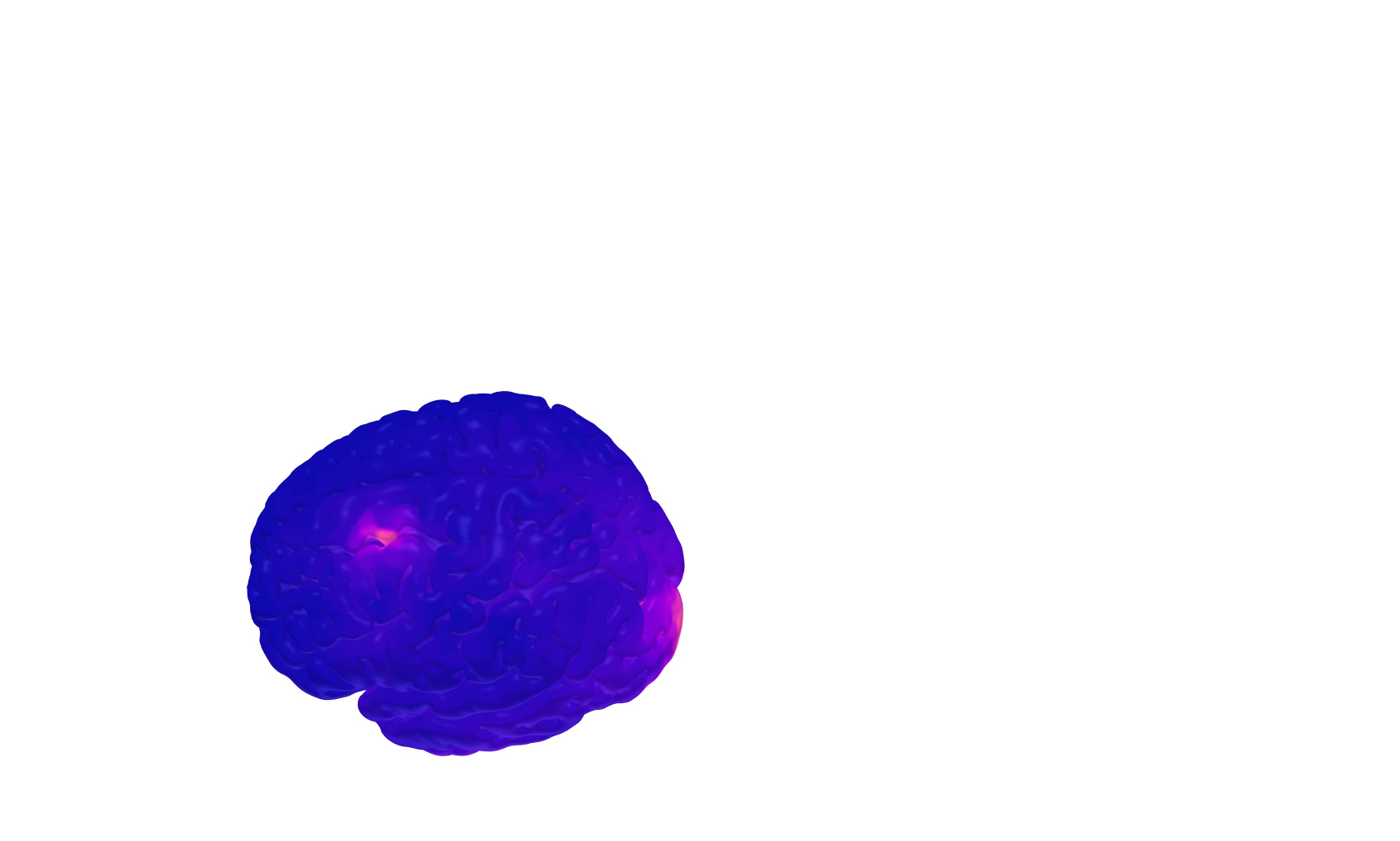}
        \end{minipage}\begin{minipage}{0.16\linewidth}
            \includegraphics[trim={5cm 1cm 14cm 8cm},clip,width=\linewidth]{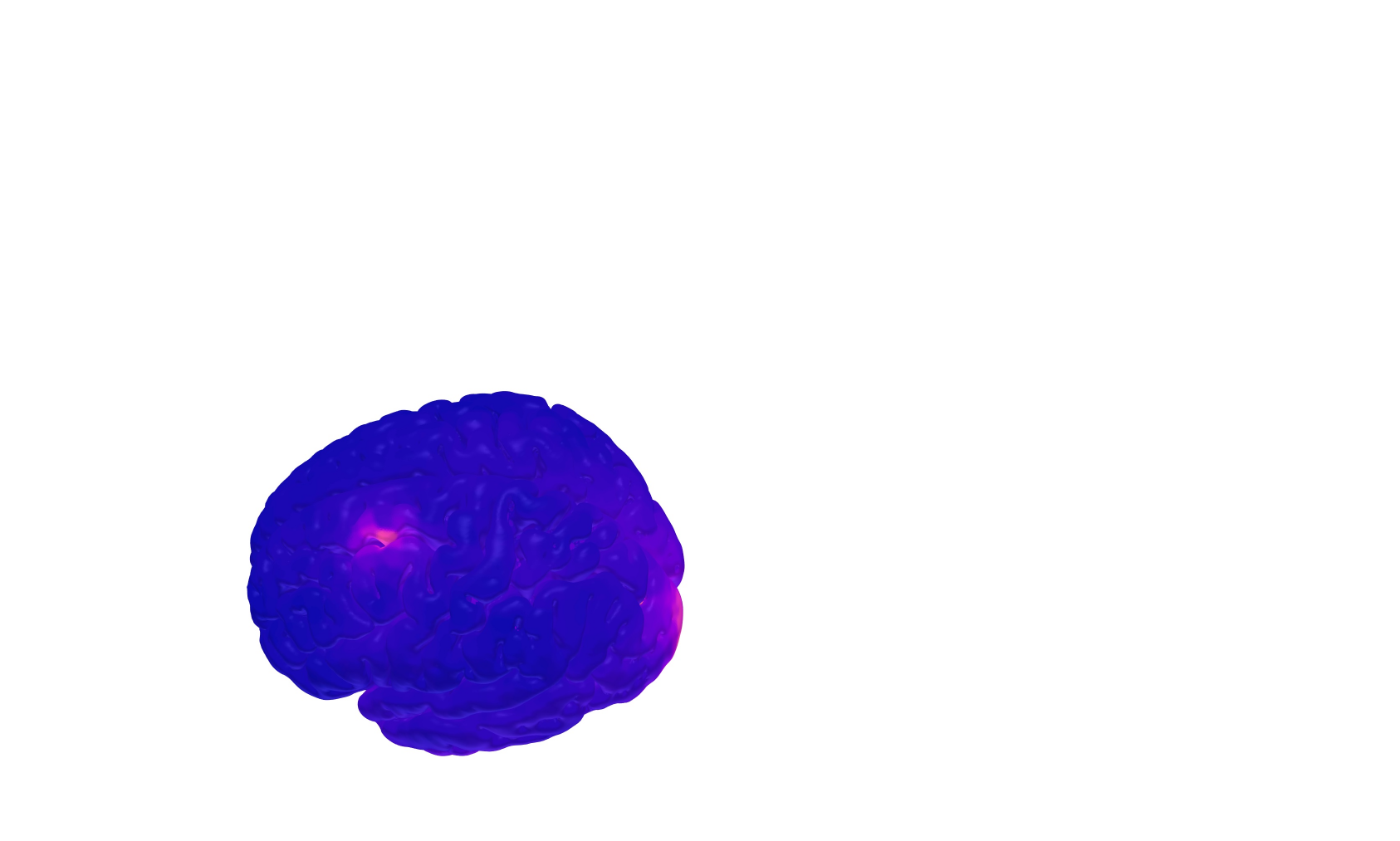}
        \end{minipage}\begin{minipage}{0.16\linewidth}
            \includegraphics[trim={5cm 1cm 14cm 8cm},clip,width=\linewidth]{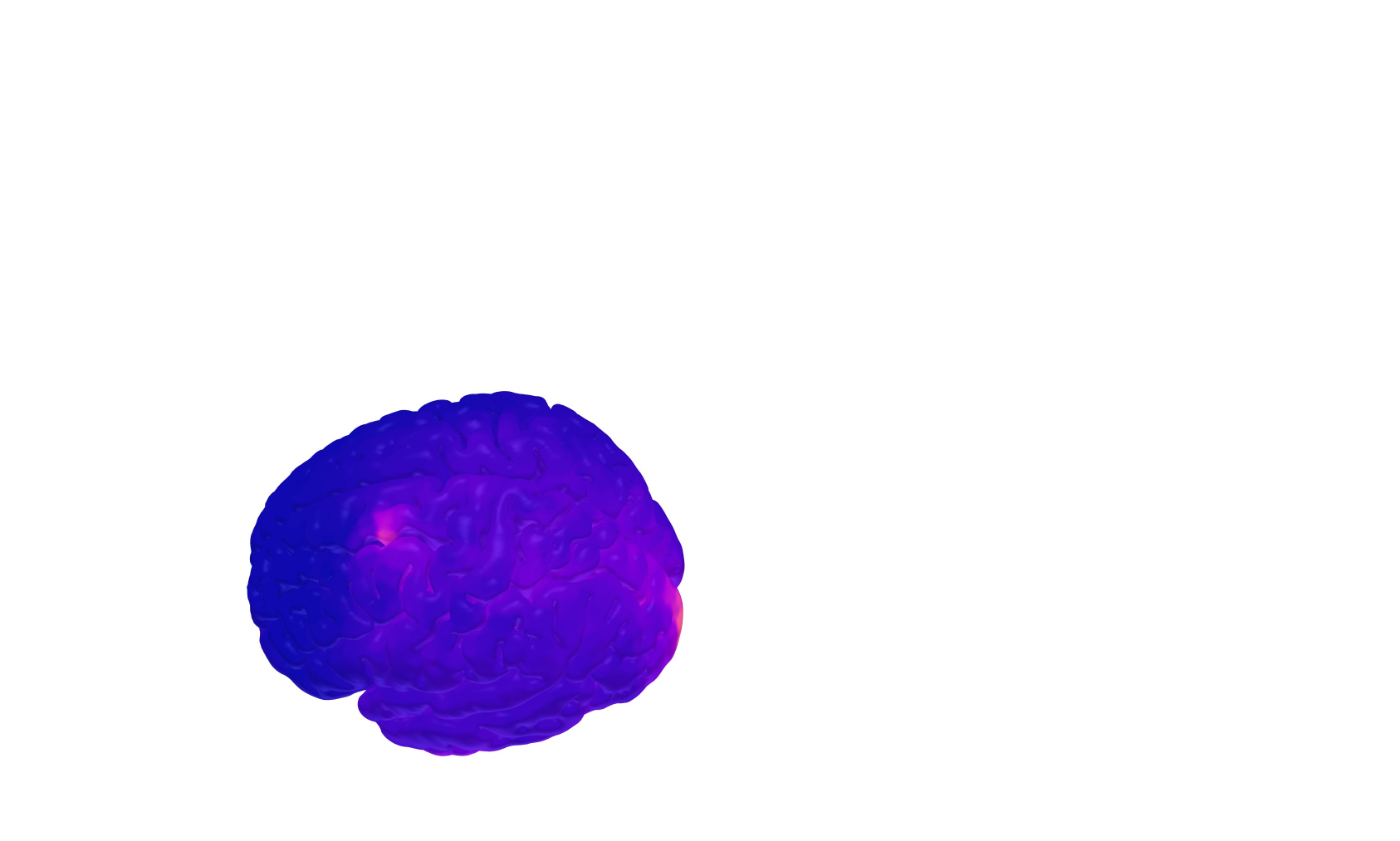}
        \end{minipage}\begin{minipage}{0.16\linewidth}
            \includegraphics[trim={5cm 1cm 14cm 8cm},clip,width=\linewidth]{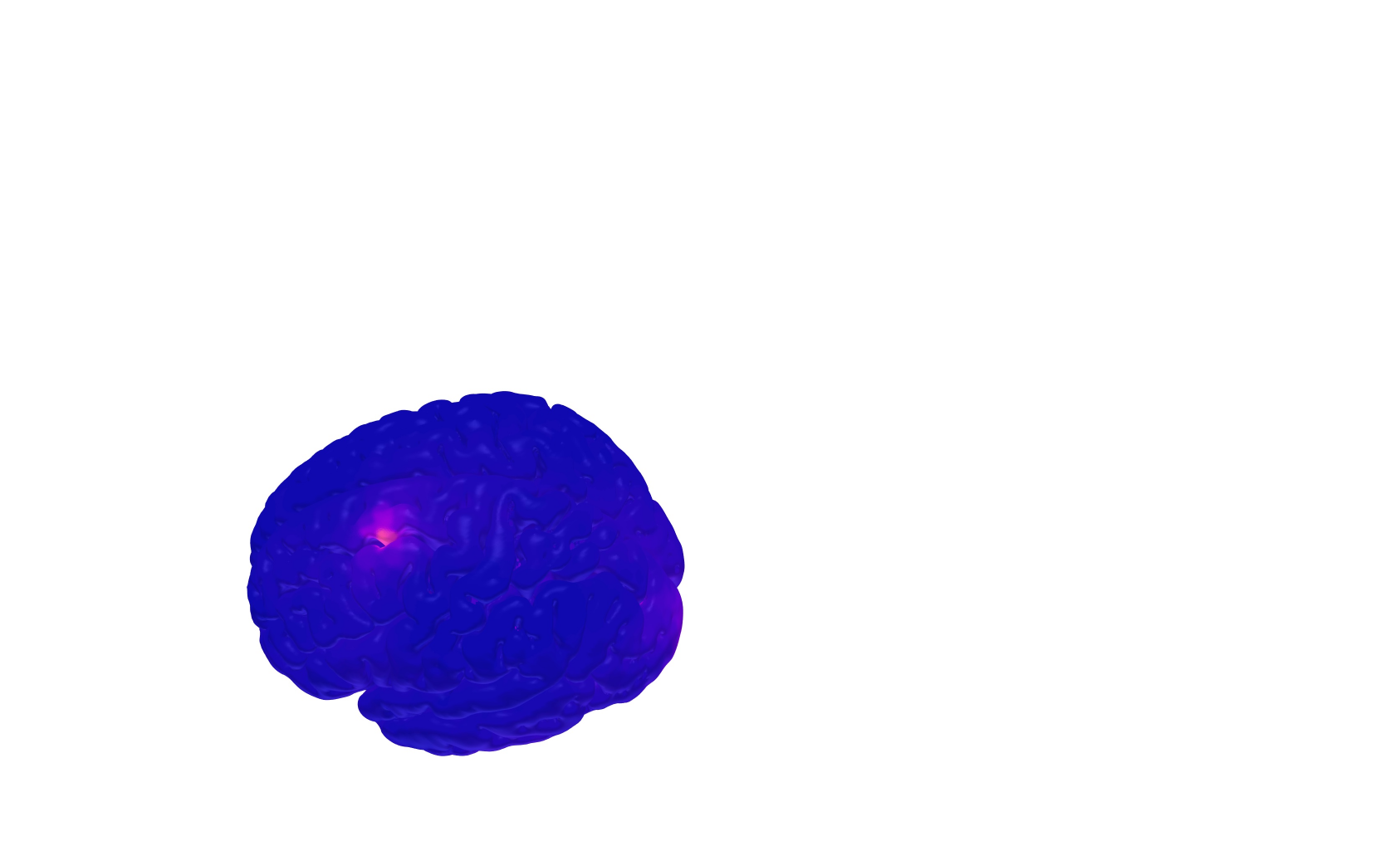}
        \end{minipage}\begin{minipage}{0.16\linewidth}
            \includegraphics[trim={5cm 1cm 14cm 8cm},clip,width=\linewidth]{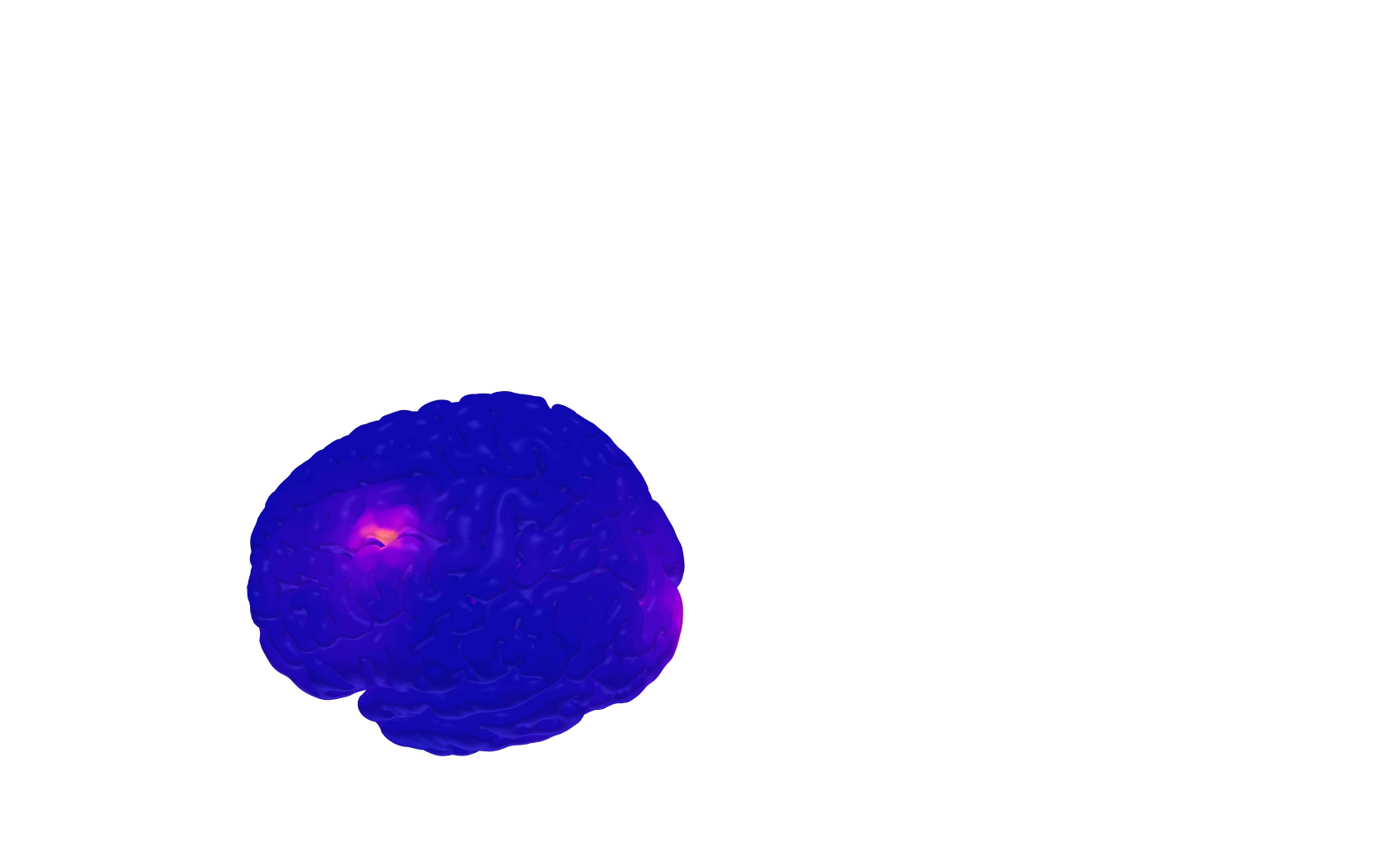}
        \end{minipage}
    \end{minipage}
    
    \caption{Spatial source reconstructions under the noiseless condition. Rows correspond to Change-rate Standardized Kalman filter (CrSKF), Standardized Kalman filter (SKF), Kalman filter (KF), and Linearly Constrained Minimum Variance (LCMV) beamformer; columns show the reconstructions at the modeled FVEP component times associated with V1, V2d, V2v+V3d, V3v, LIP, and FEF.}
    \label{fig:brain_noiseless}
\end{figure}

\begin{figure}[h!]
    \centering
    \begin{minipage}{0.02\linewidth}
        \hfill
    \end{minipage}\begin{minipage}{0.98\linewidth}
        \begin{minipage}{0.16\linewidth}
            \centering
            \scriptsize{V1}
        \end{minipage}\begin{minipage}{0.16\linewidth}
            \centering
            \scriptsize{V2d}
        \end{minipage}\begin{minipage}{0.16\linewidth}
            \centering
            \scriptsize{V2v+V3d}
        \end{minipage}\begin{minipage}{0.16\linewidth}
            \centering
            \scriptsize{V3v}
        \end{minipage}\begin{minipage}{0.16\linewidth}
        \centering
            \scriptsize{LIP}
        \end{minipage}\begin{minipage}{0.16\linewidth}
        \centering
            \scriptsize{FEF}
        \end{minipage}\end{minipage}
    
    \begin{minipage}{0.02\linewidth}
        \rotatebox{90}{\scriptsize CrSKF}
    \end{minipage}\begin{minipage}{0.98\linewidth}
        \begin{minipage}{0.16\linewidth}
            \includegraphics[trim={5cm 1cm 14cm 8cm},clip,width=\linewidth]{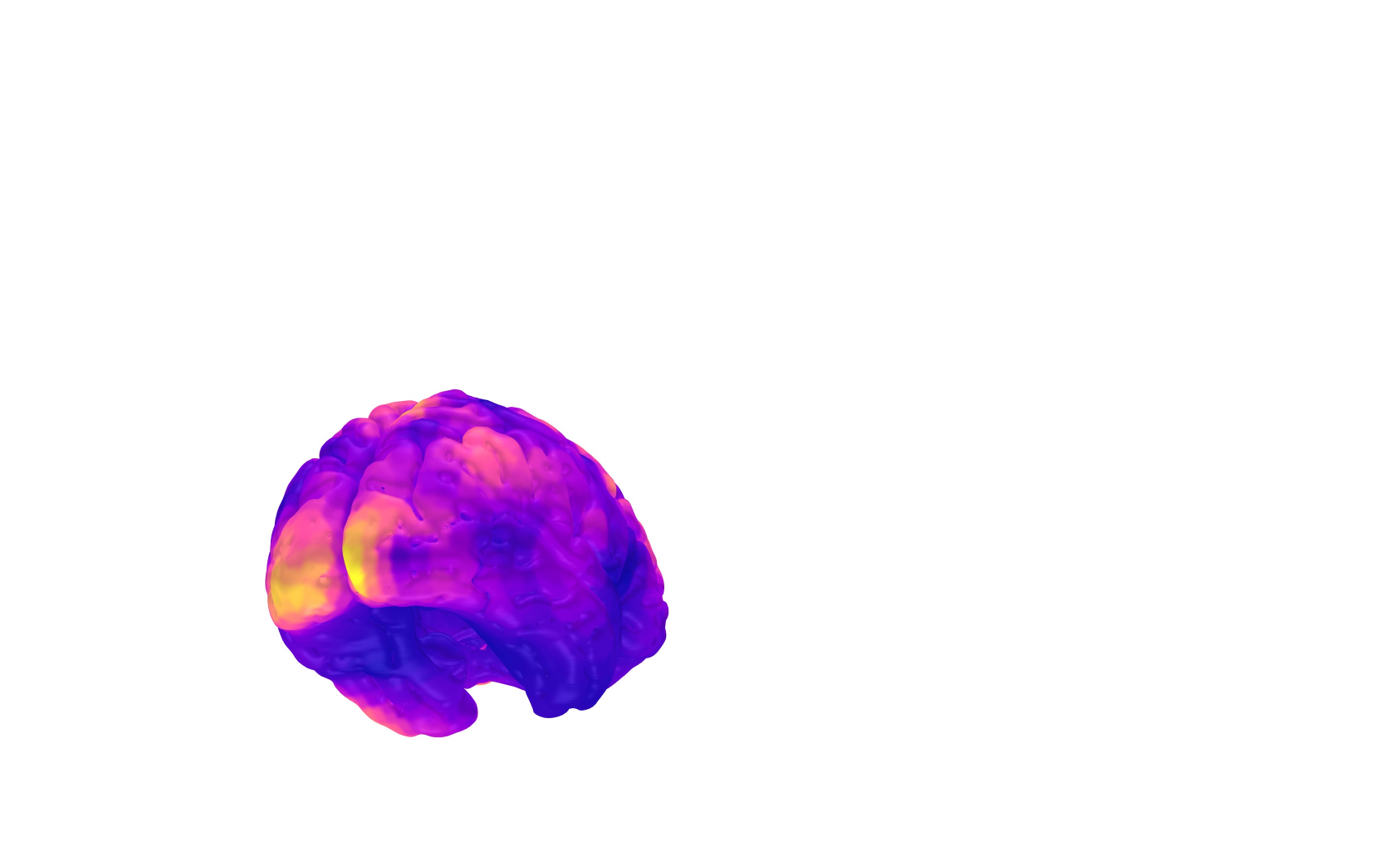}
        \end{minipage}\begin{minipage}{0.16\linewidth}
            \includegraphics[trim={5cm 1cm 14cm 8cm},clip,width=\linewidth]{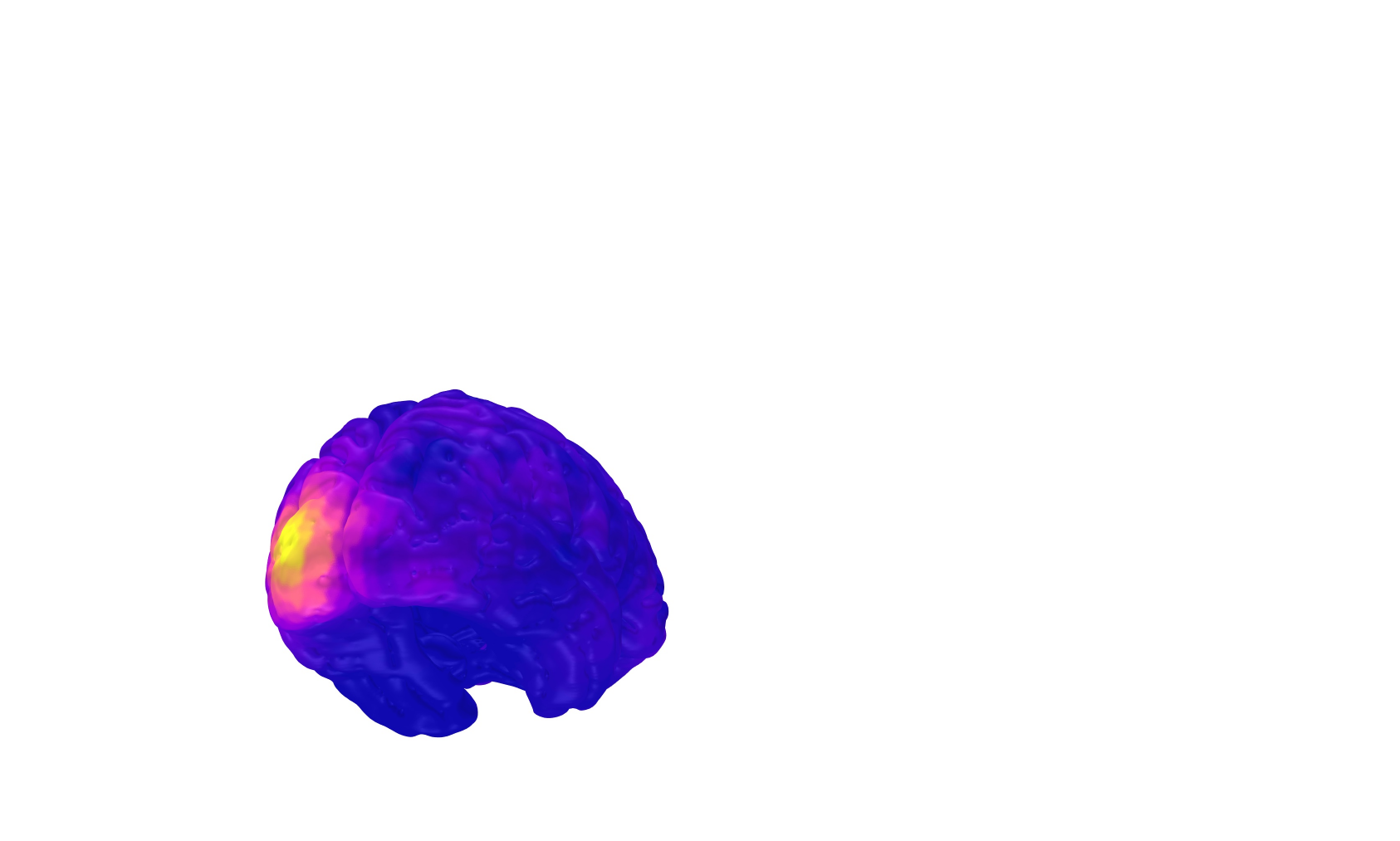}
        \end{minipage}\begin{minipage}{0.16\linewidth}
            \includegraphics[trim={5cm 1cm 14cm 8cm},clip,width=\linewidth]{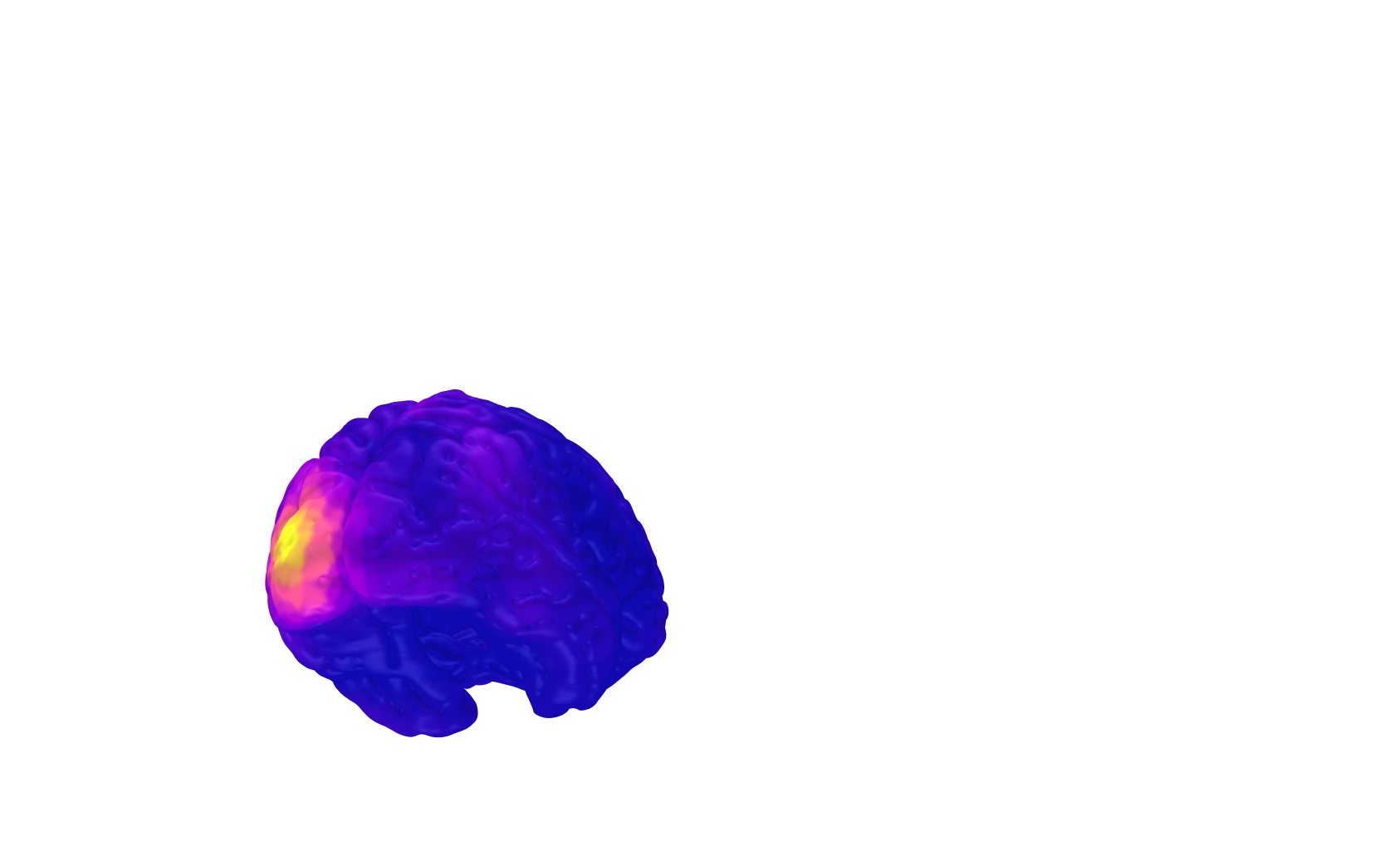}
        \end{minipage}\begin{minipage}{0.16\linewidth}
            \includegraphics[trim={5cm 1cm 14cm 8cm},clip,width=\linewidth]{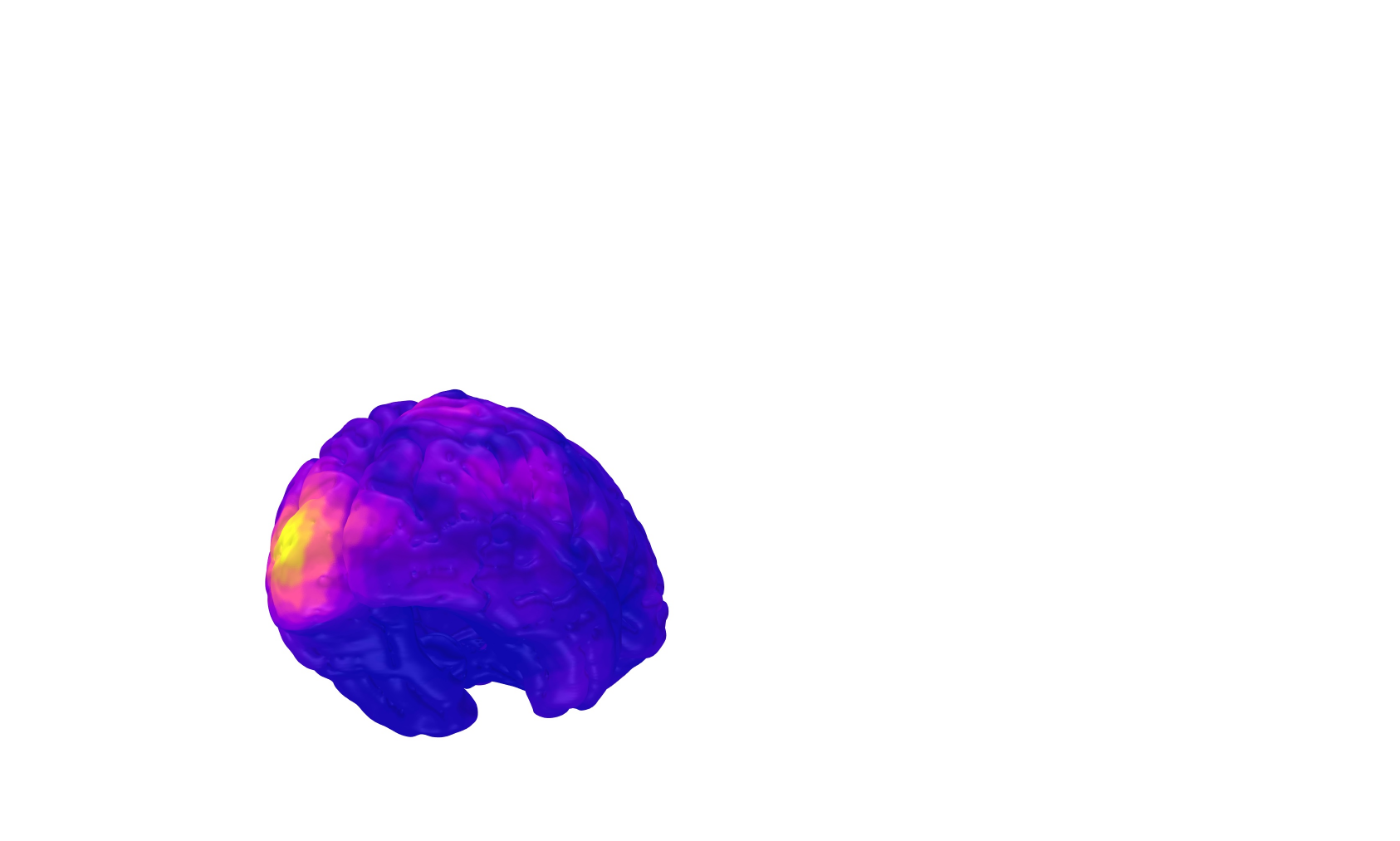}
        \end{minipage}\begin{minipage}{0.16\linewidth}
            \includegraphics[trim={5cm 1cm 14cm 8cm},clip,width=\linewidth]{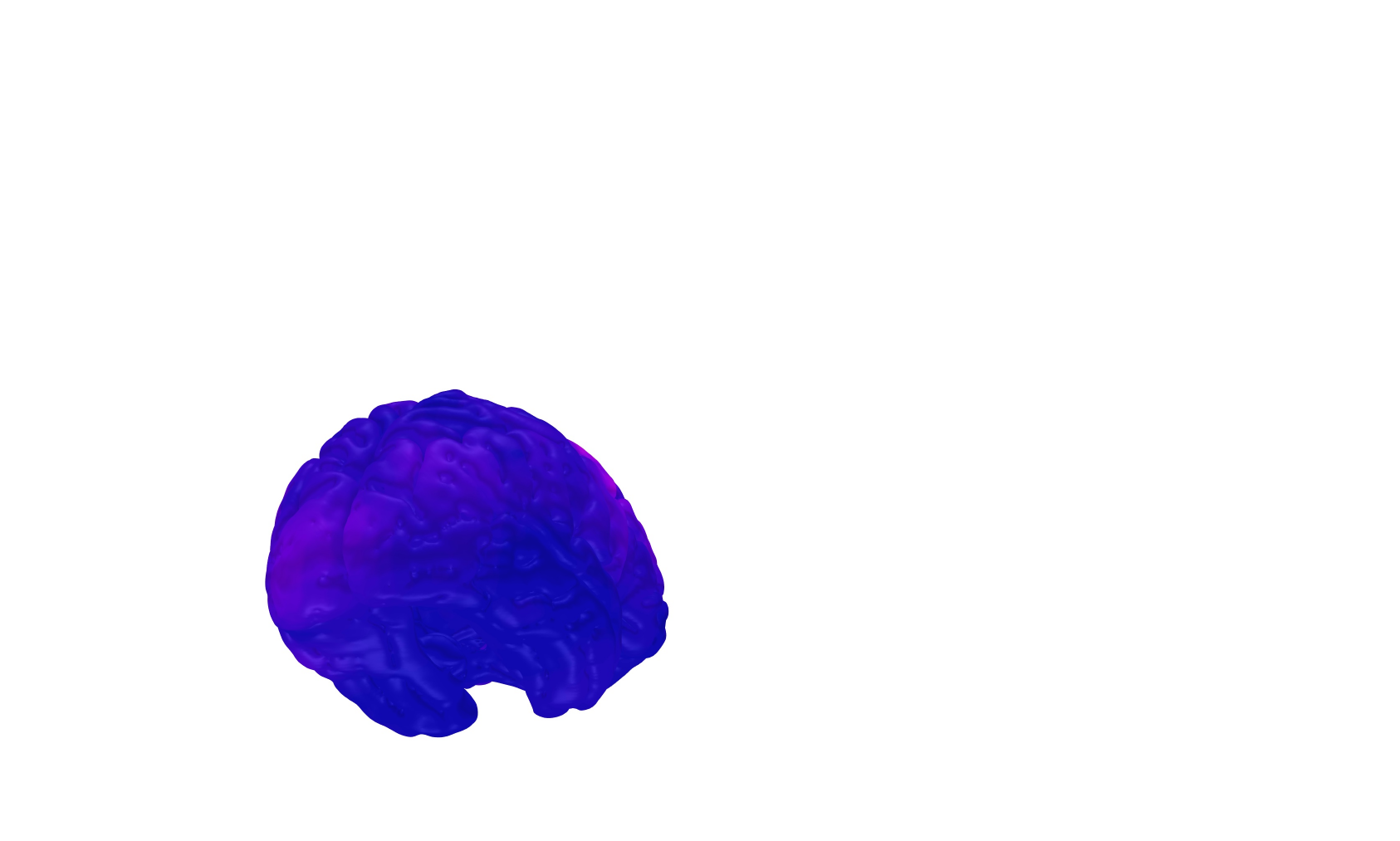}
        \end{minipage}\begin{minipage}{0.16\linewidth}
            \includegraphics[trim={5cm 1cm 14cm 8cm},clip,width=\linewidth]{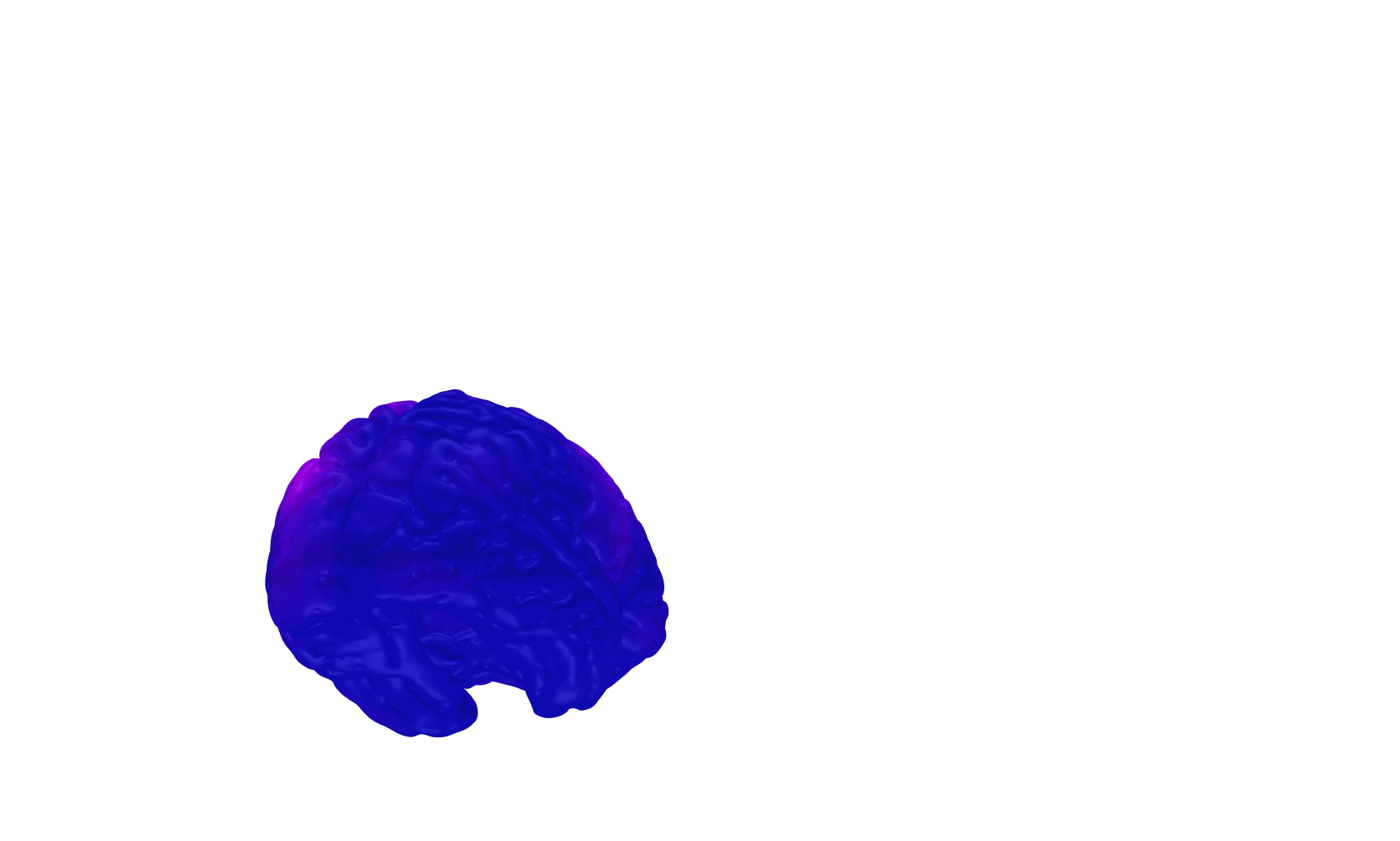}
        \end{minipage}

        \begin{minipage}{0.16\linewidth}
            \includegraphics[trim={5cm 1cm 14cm 8cm},clip,width=\linewidth]{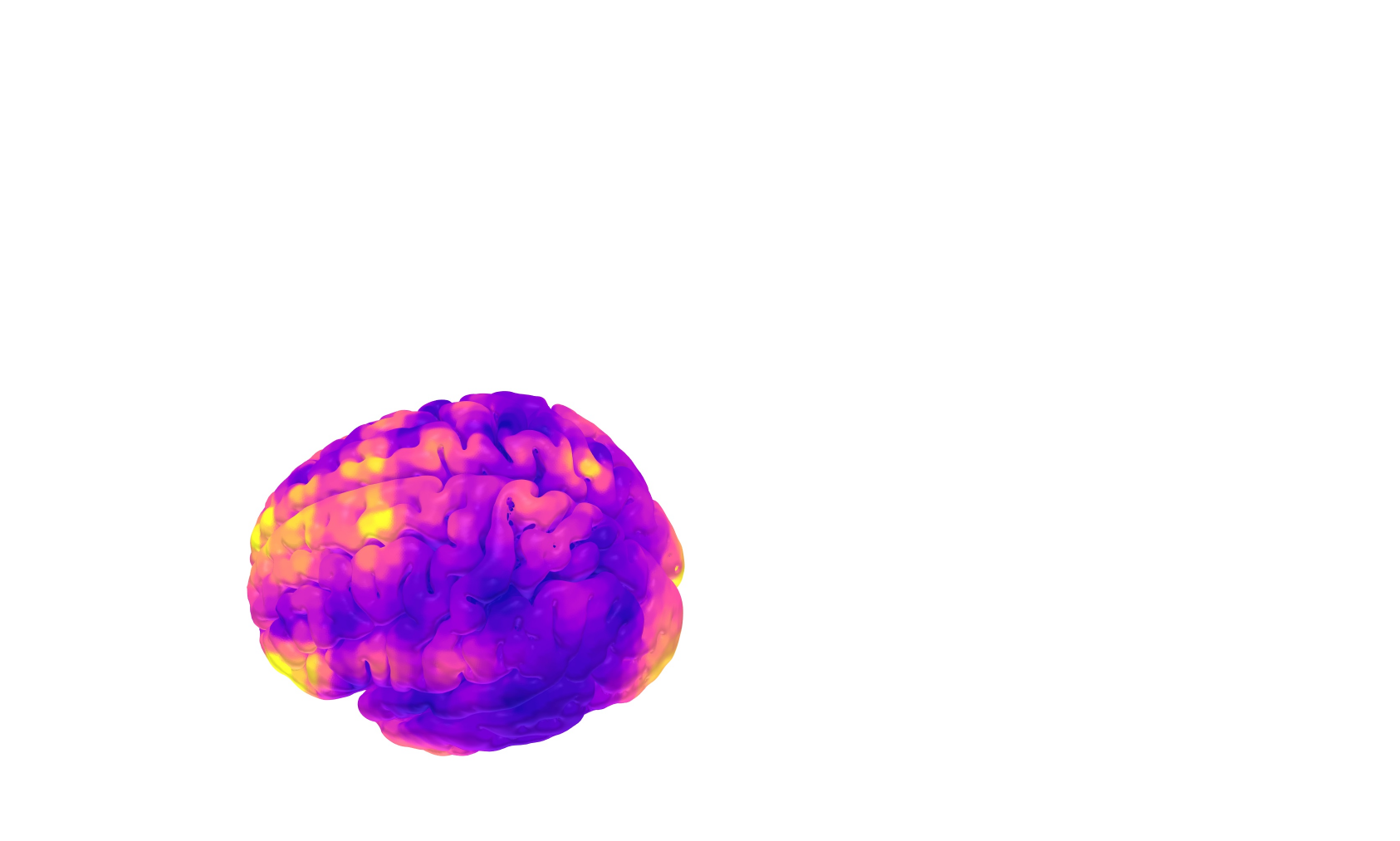}
        \end{minipage}\begin{minipage}{0.16\linewidth}
            \includegraphics[trim={5cm 1cm 14cm 8cm},clip,width=\linewidth]{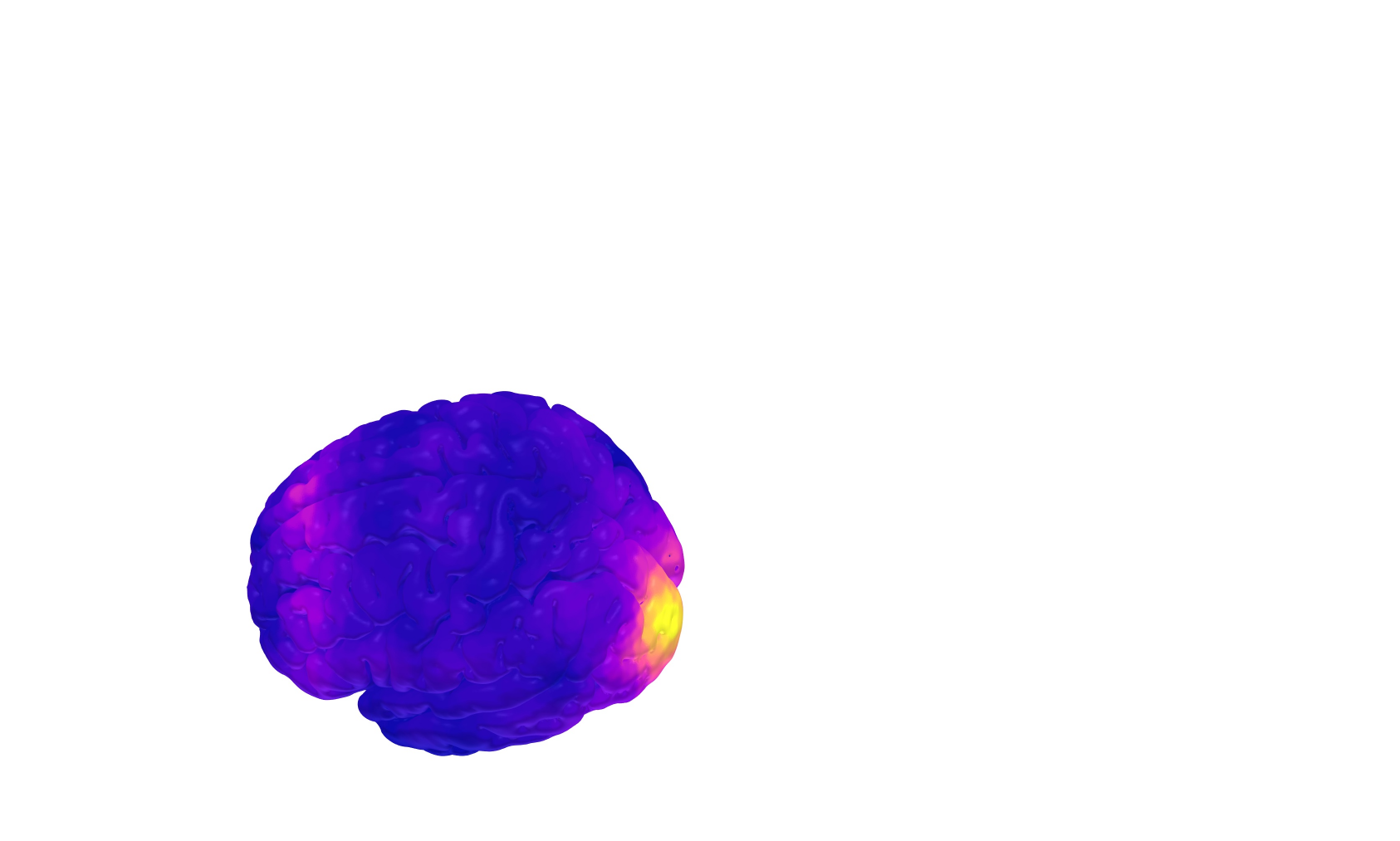}
        \end{minipage}\begin{minipage}{0.16\linewidth}
            \includegraphics[trim={5cm 1cm 14cm 8cm},clip,width=\linewidth]{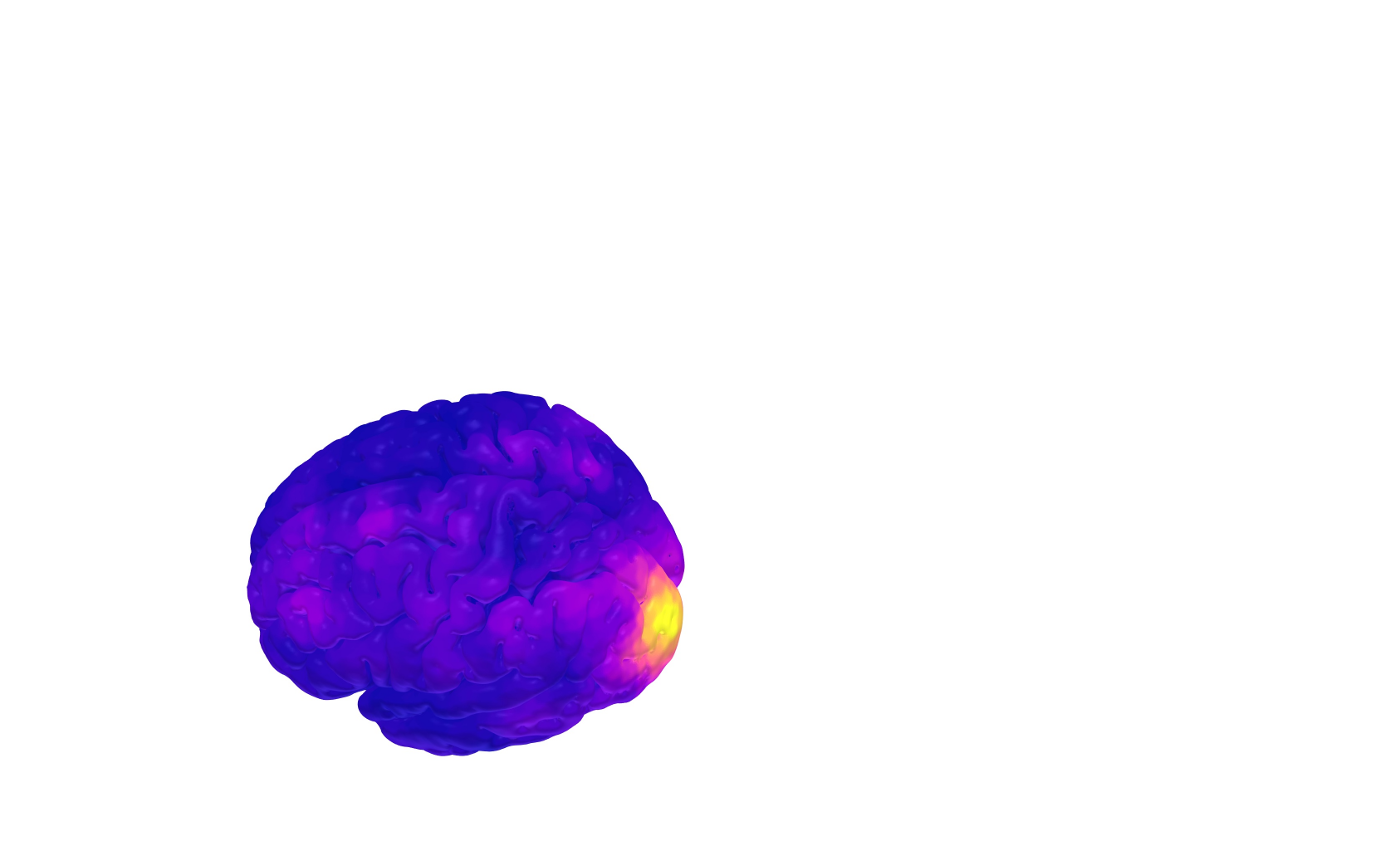}
        \end{minipage}\begin{minipage}{0.16\linewidth}
            \includegraphics[trim={5cm 1cm 14cm 8cm},clip,width=\linewidth]{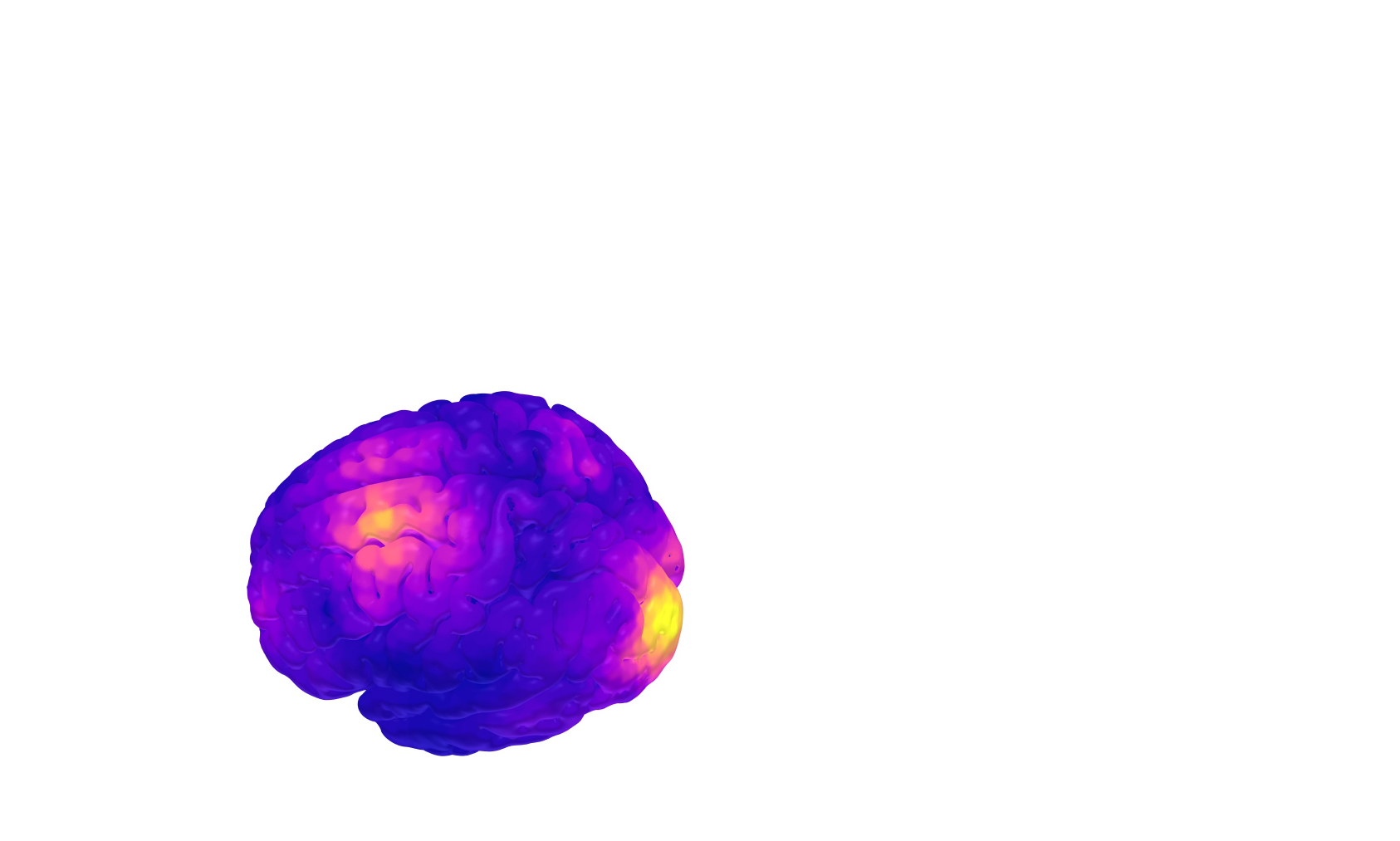}
        \end{minipage}\begin{minipage}{0.16\linewidth}
            \includegraphics[trim={5cm 1cm 14cm 8cm},clip,width=\linewidth]{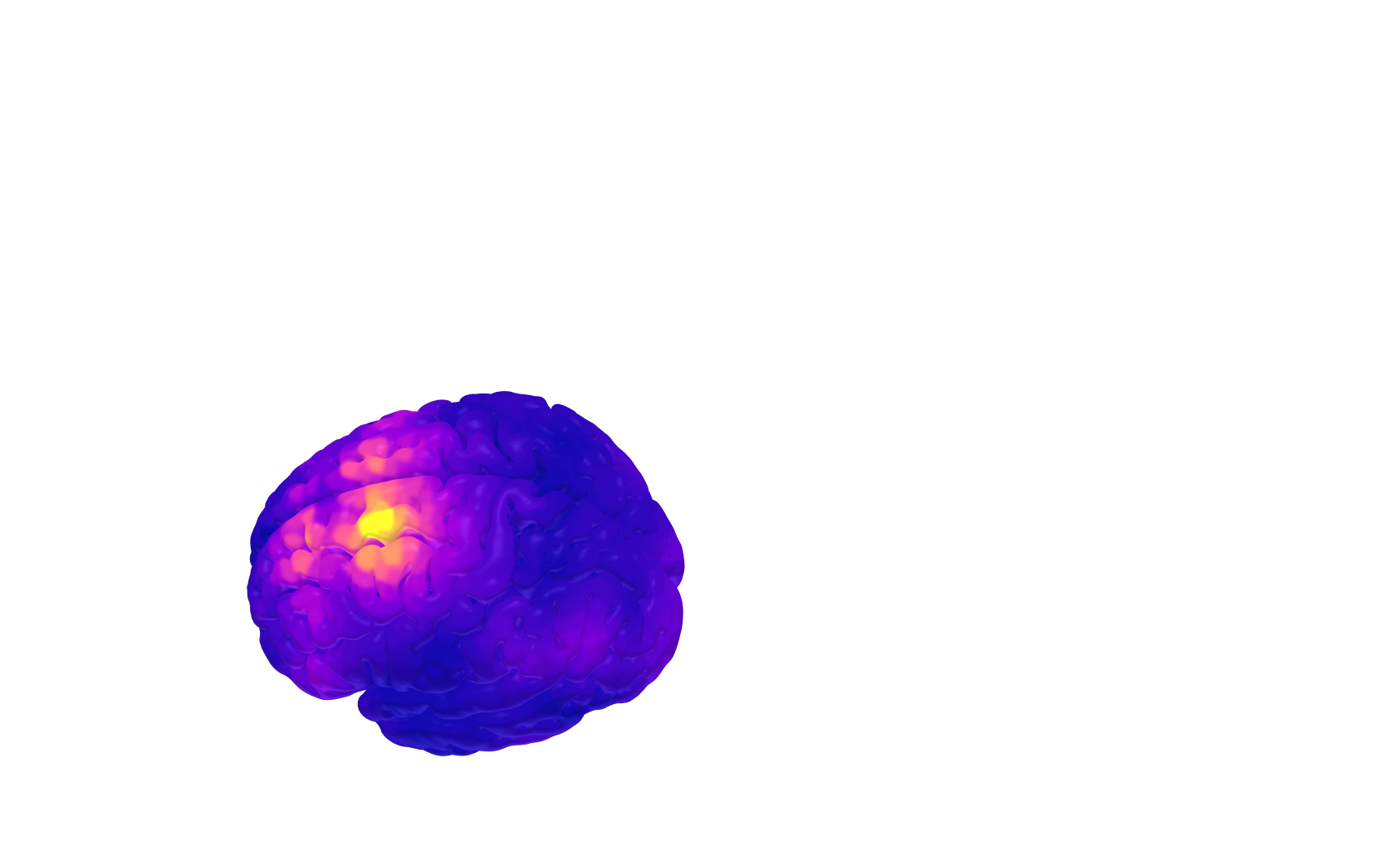}
        \end{minipage}\begin{minipage}{0.16\linewidth}
            \includegraphics[trim={5cm 1cm 14cm 8cm},clip,width=\linewidth]{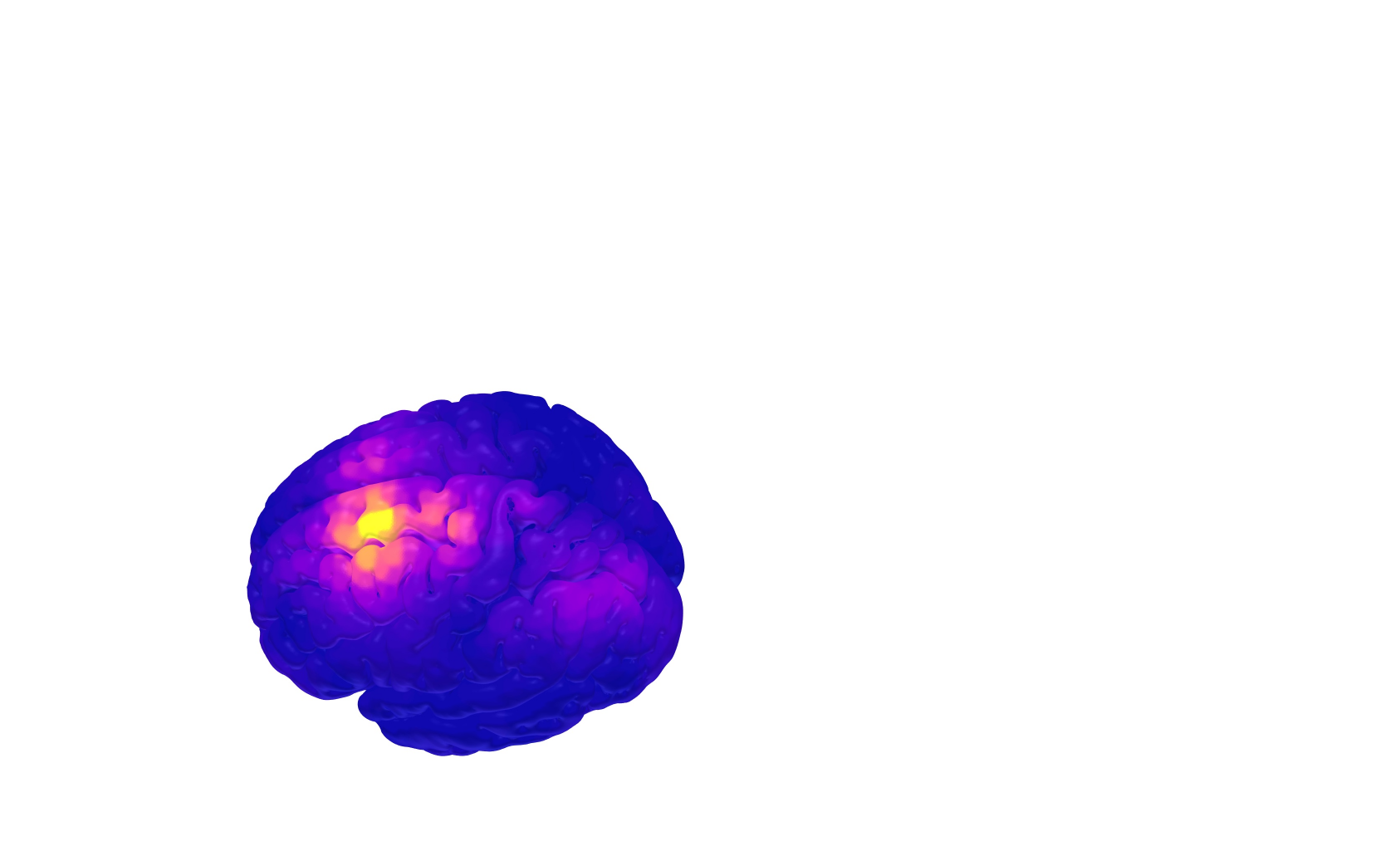}
        \end{minipage}
    \end{minipage}
    \noindent\rule[7pt]{\linewidth}{0.4pt}

    \begin{minipage}{0.02\linewidth}
        \rotatebox{90}{\scriptsize SKF}
    \end{minipage}\begin{minipage}{0.98\linewidth}
        \begin{minipage}{0.16\linewidth}
            \includegraphics[trim={5cm 1cm 14cm 8cm},clip,width=\linewidth]{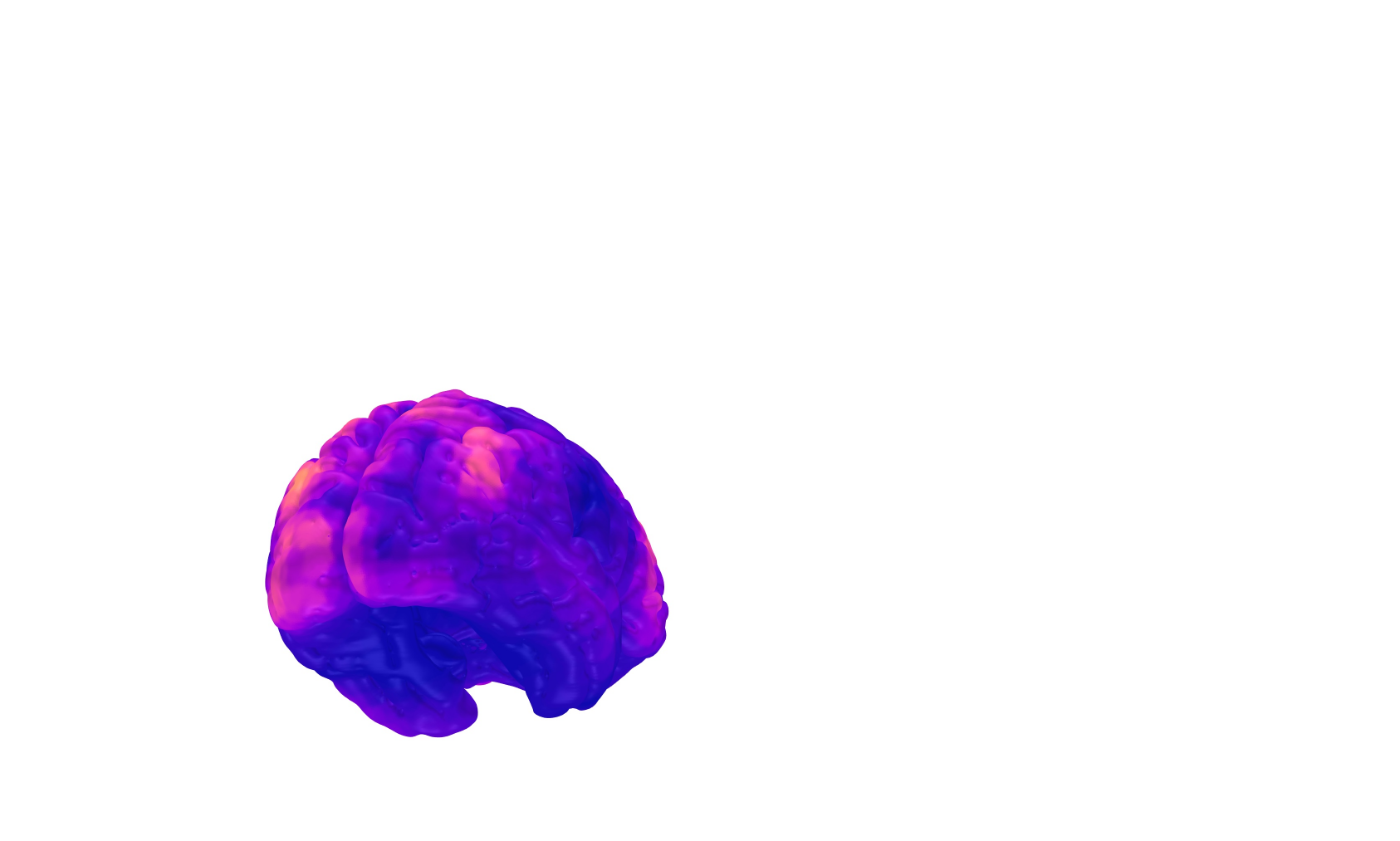}
        \end{minipage}\begin{minipage}{0.16\linewidth}
            \includegraphics[trim={5cm 1cm 14cm 8cm},clip,width=\linewidth]{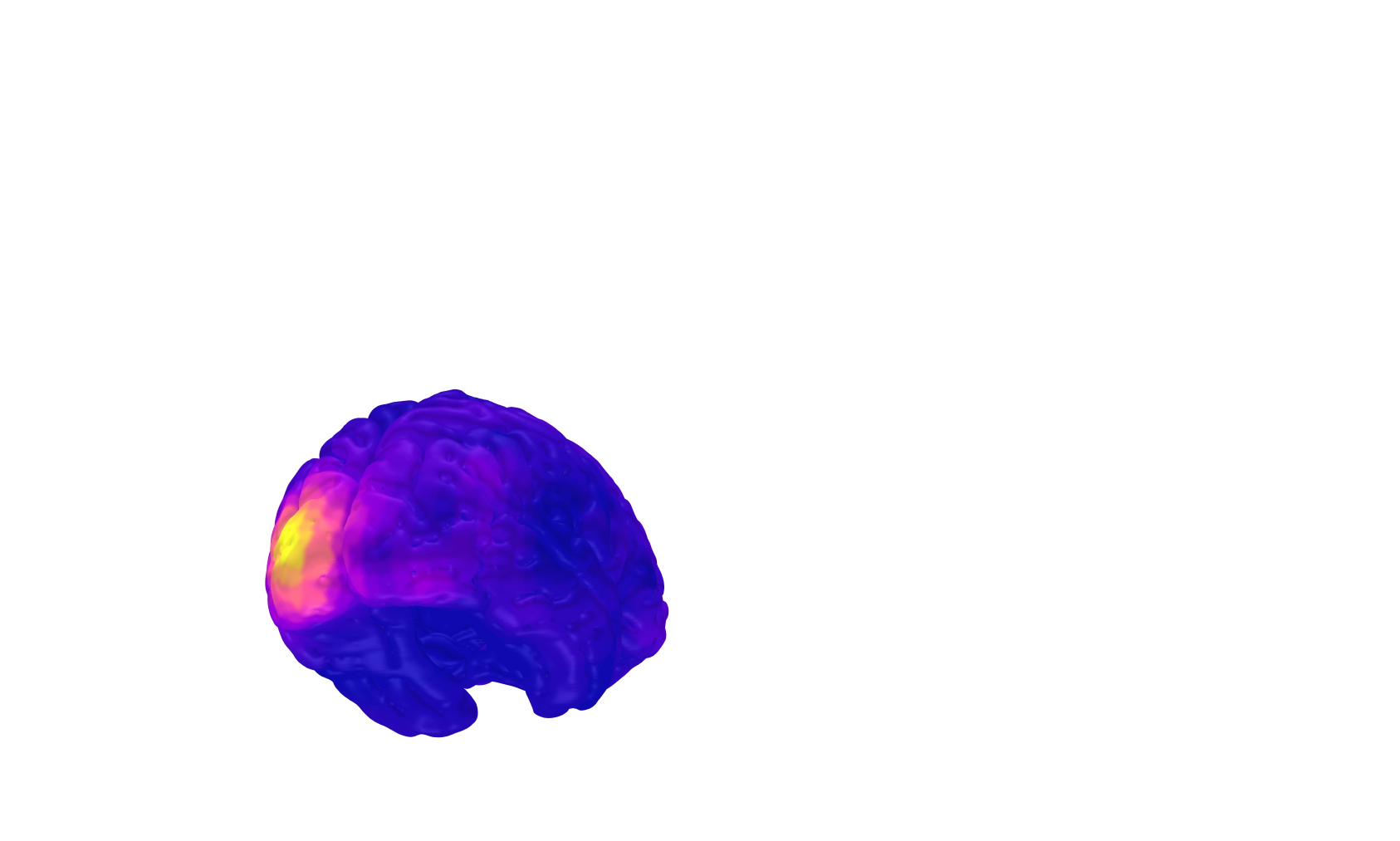}
        \end{minipage}\begin{minipage}{0.16\linewidth}
            \includegraphics[trim={5cm 1cm 14cm 8cm},clip,width=\linewidth]{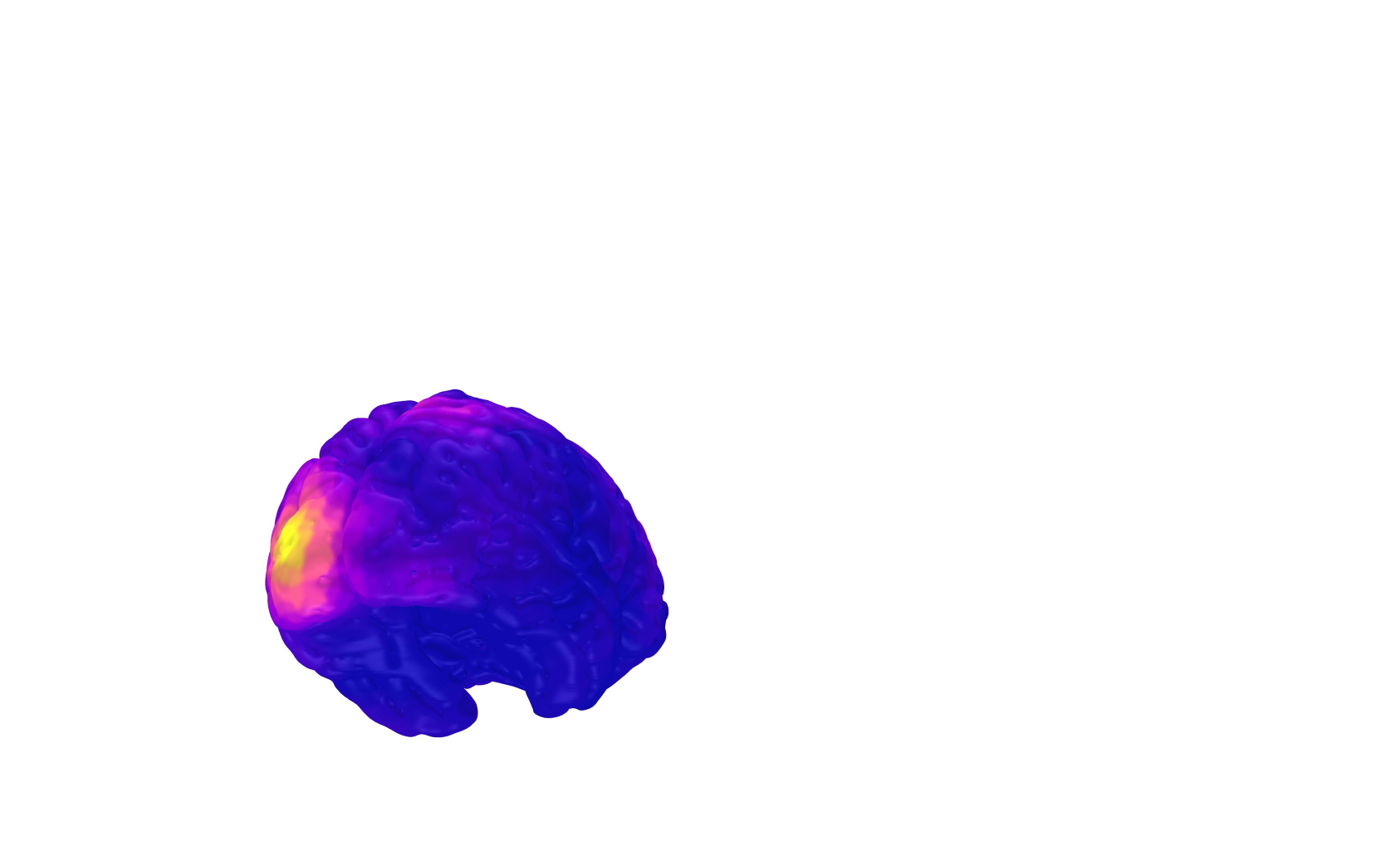}
        \end{minipage}\begin{minipage}{0.16\linewidth}
            \includegraphics[trim={5cm 1cm 14cm 8cm},clip,width=\linewidth]{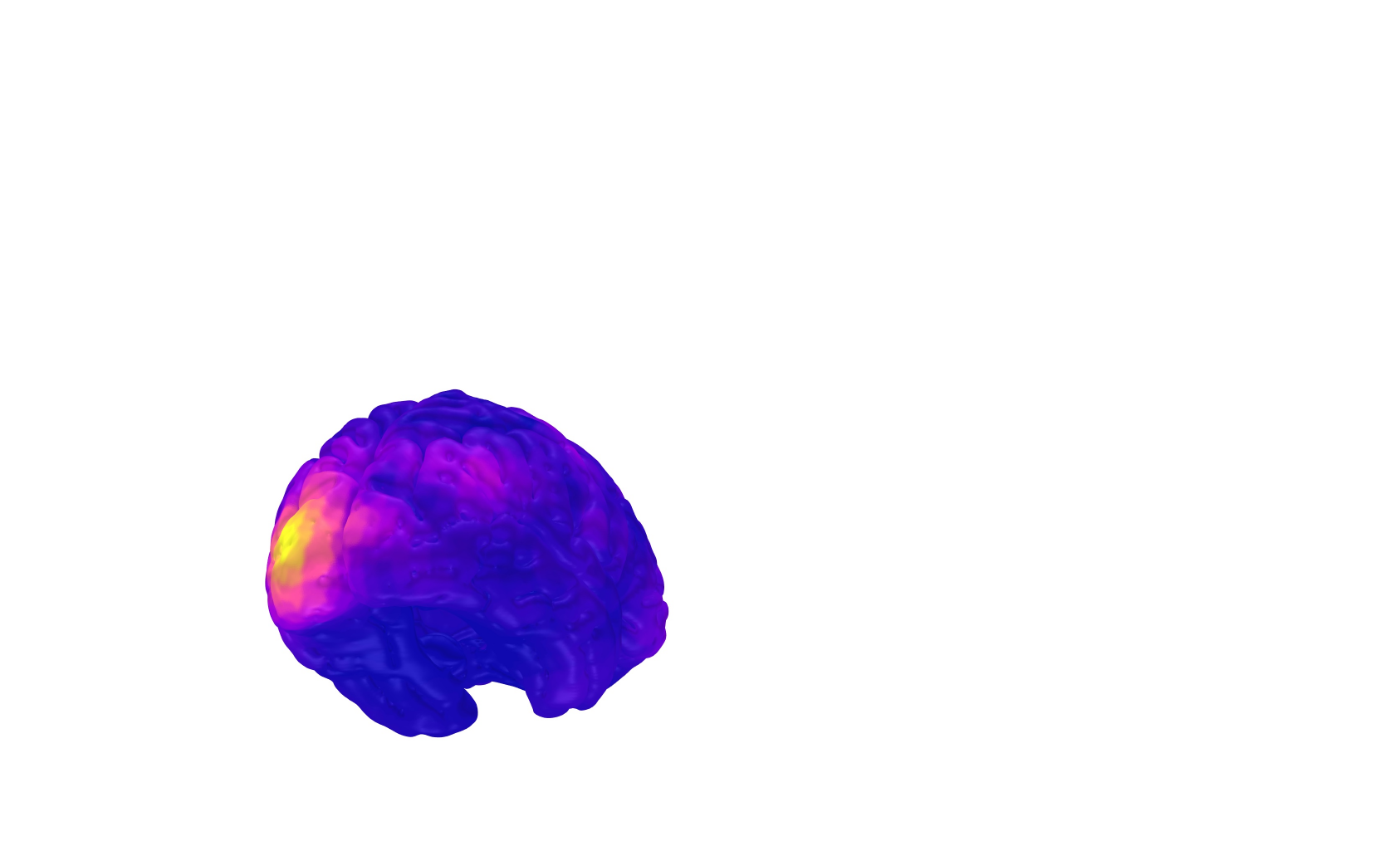}
        \end{minipage}\begin{minipage}{0.16\linewidth}
            \includegraphics[trim={5cm 1cm 14cm 8cm},clip,width=\linewidth]{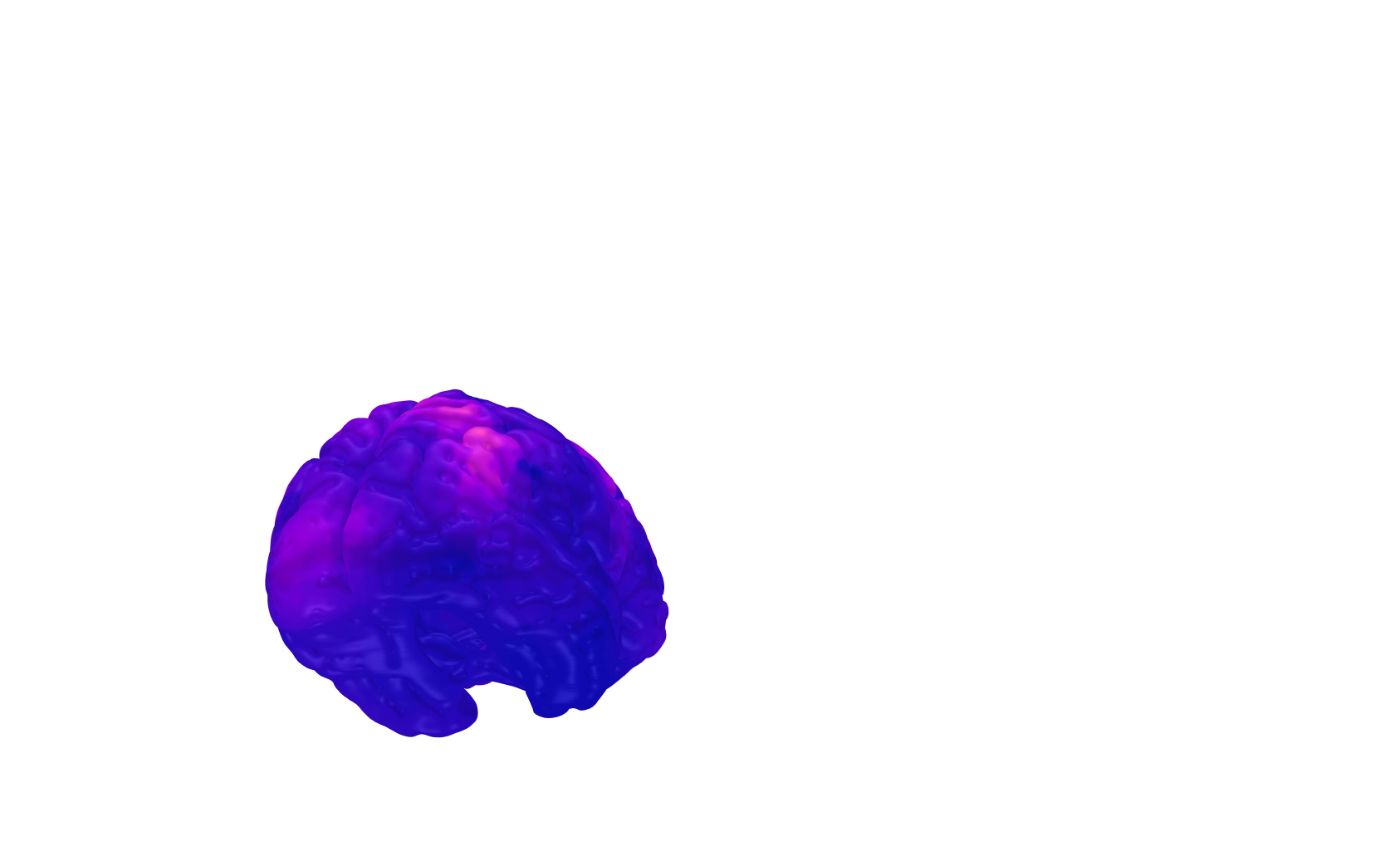}
        \end{minipage}\begin{minipage}{0.16\linewidth}
            \includegraphics[trim={5cm 1cm 14cm 8cm},clip,width=\linewidth]{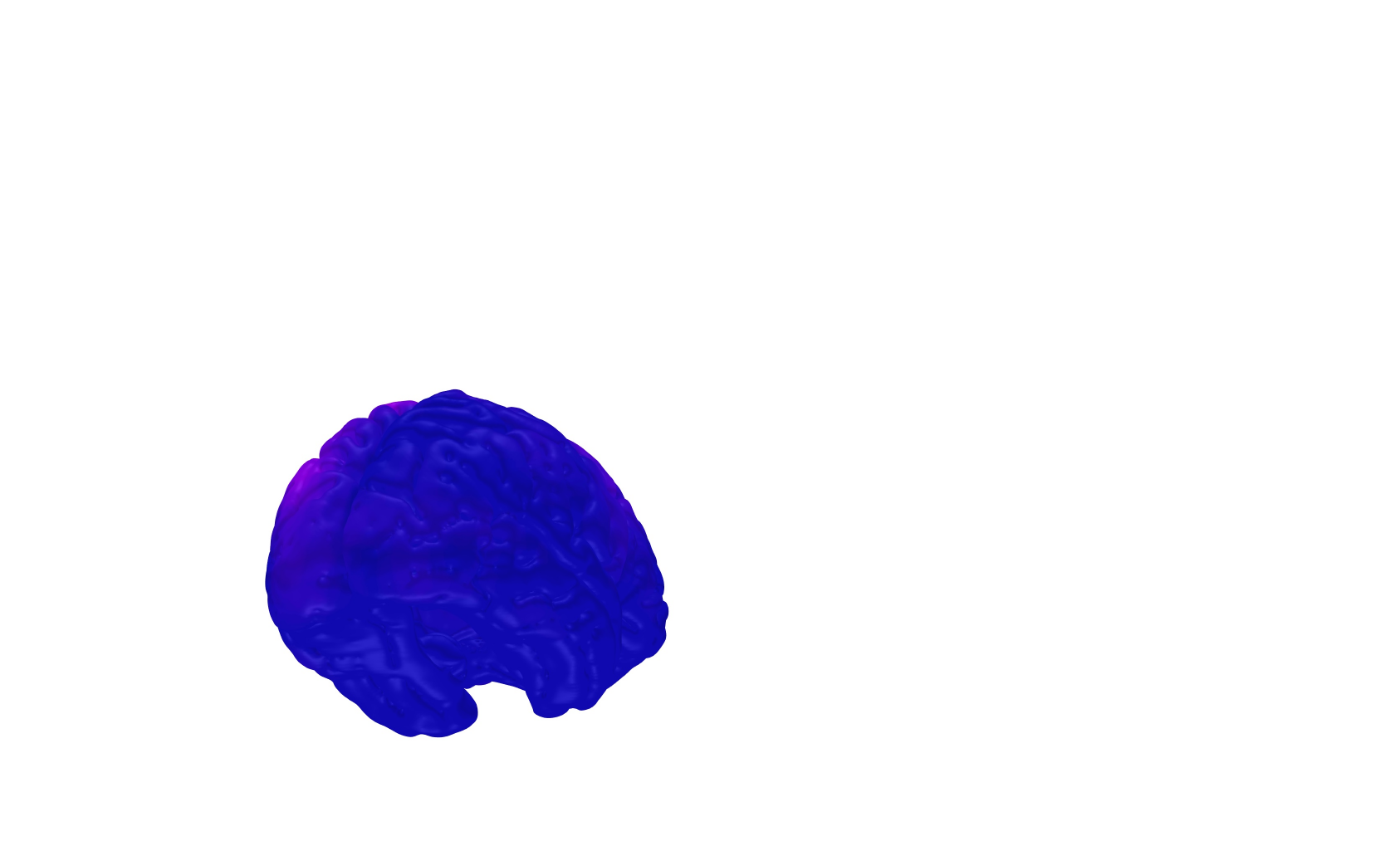}
        \end{minipage}

        \begin{minipage}{0.16\linewidth}
            \includegraphics[trim={5cm 1cm 14cm 8cm},clip,width=\linewidth]{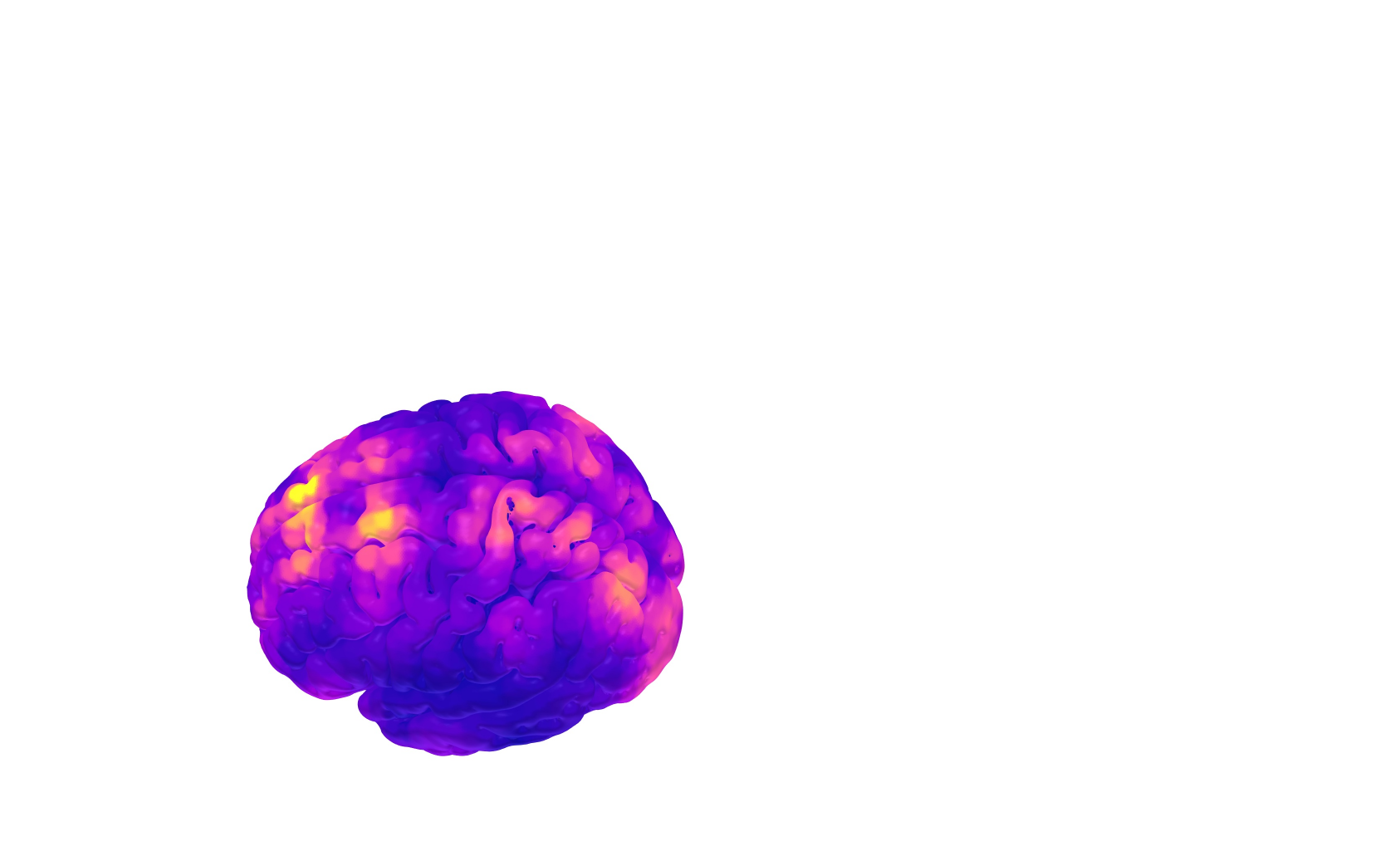}
        \end{minipage}\begin{minipage}{0.16\linewidth}
            \includegraphics[trim={5cm 1cm 14cm 8cm},clip,width=\linewidth]{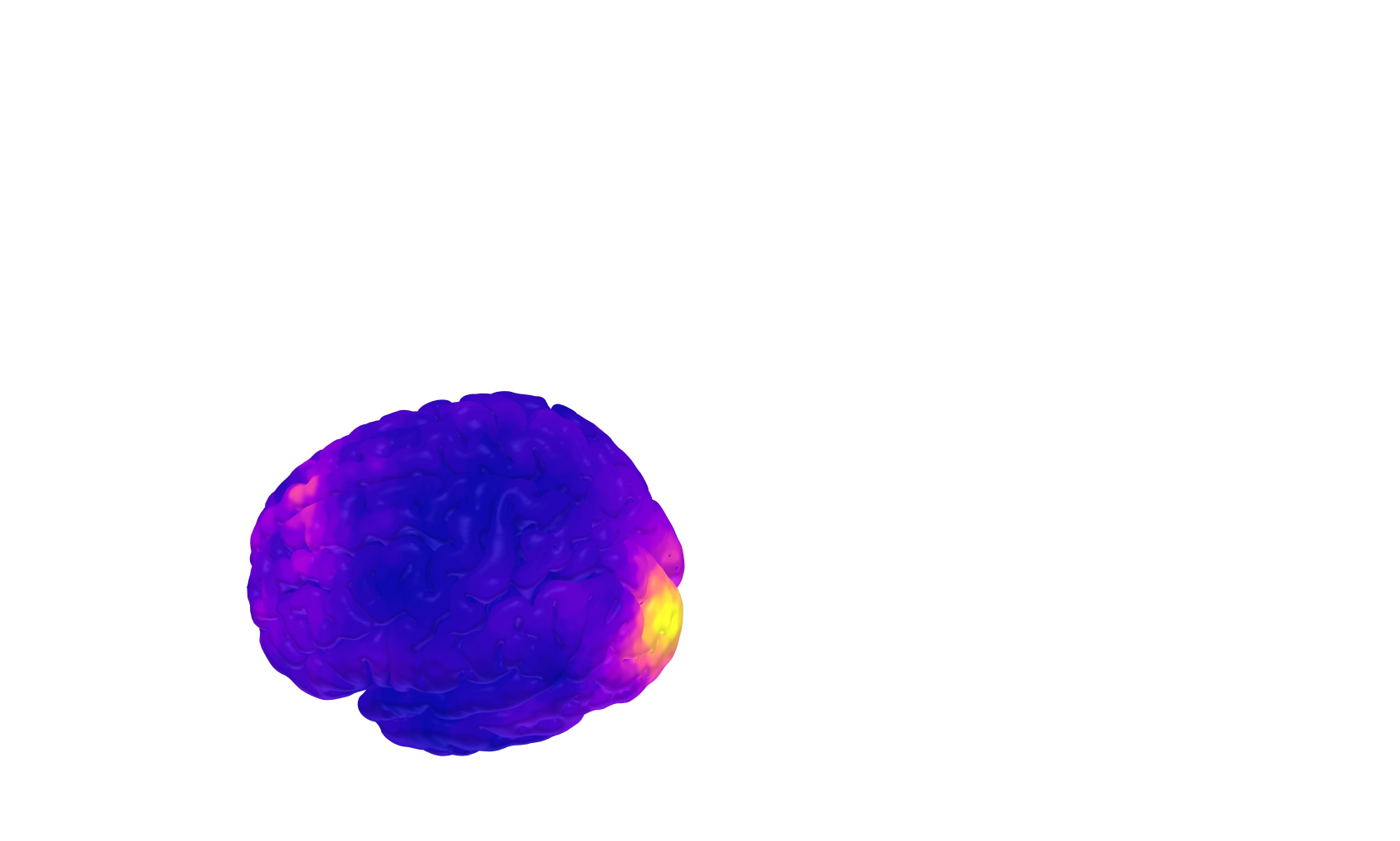}
        \end{minipage}\begin{minipage}{0.16\linewidth}
            \includegraphics[trim={5cm 1cm 14cm 8cm},clip,width=\linewidth]{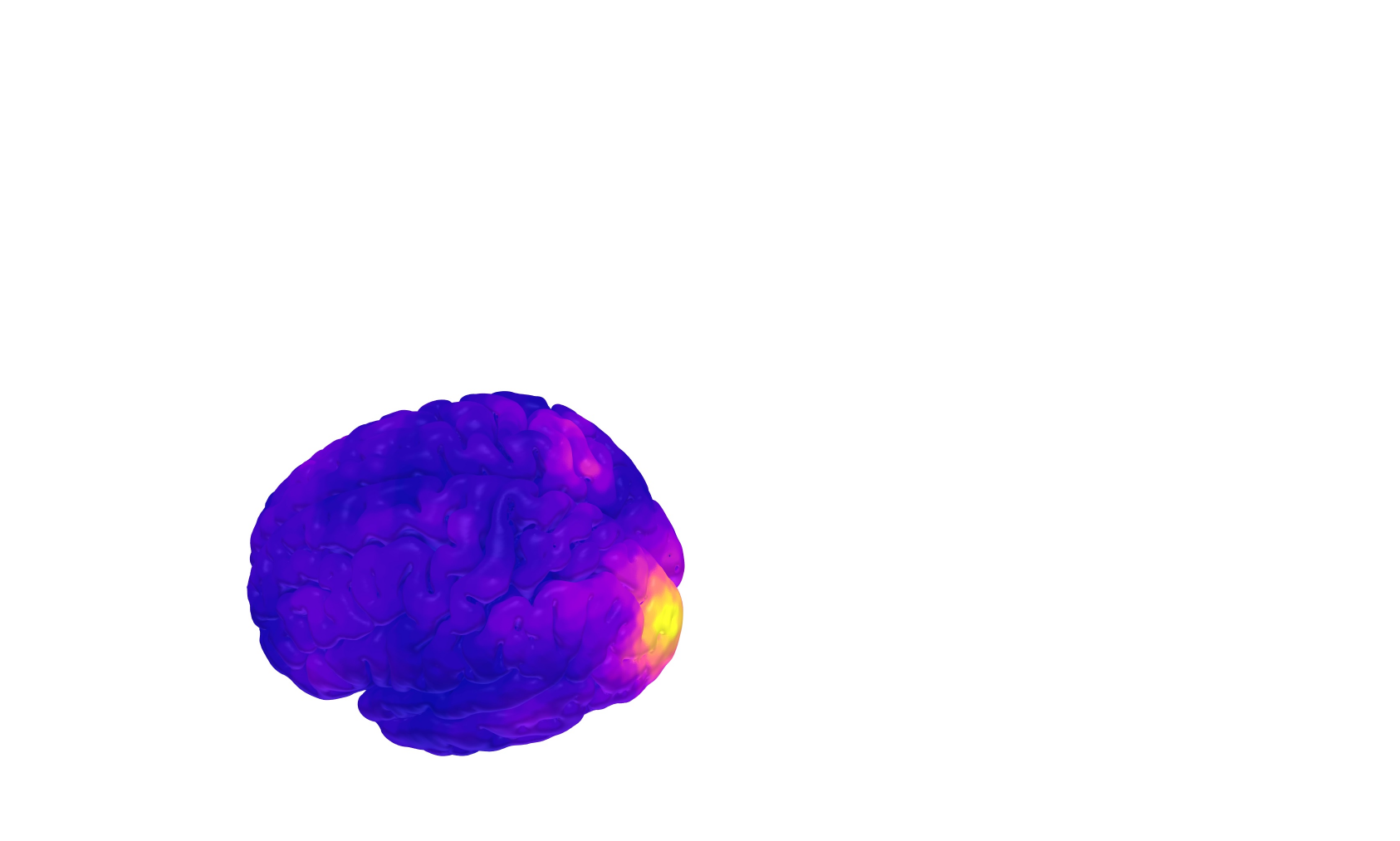}
        \end{minipage}\begin{minipage}{0.16\linewidth}
            \includegraphics[trim={5cm 1cm 14cm 8cm},clip,width=\linewidth]{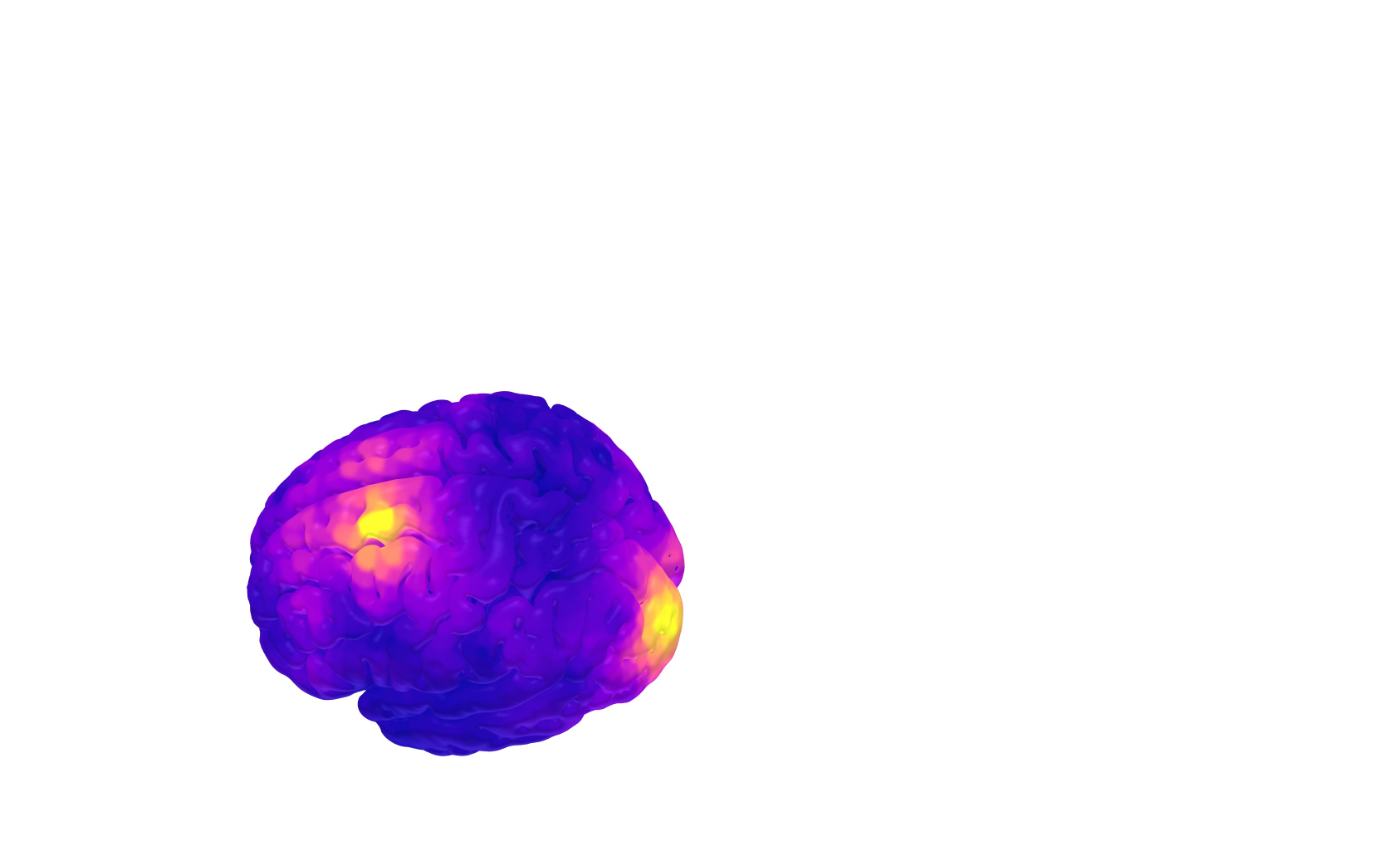}
        \end{minipage}\begin{minipage}{0.16\linewidth}
            \includegraphics[trim={5cm 1cm 14cm 8cm},clip,width=\linewidth]{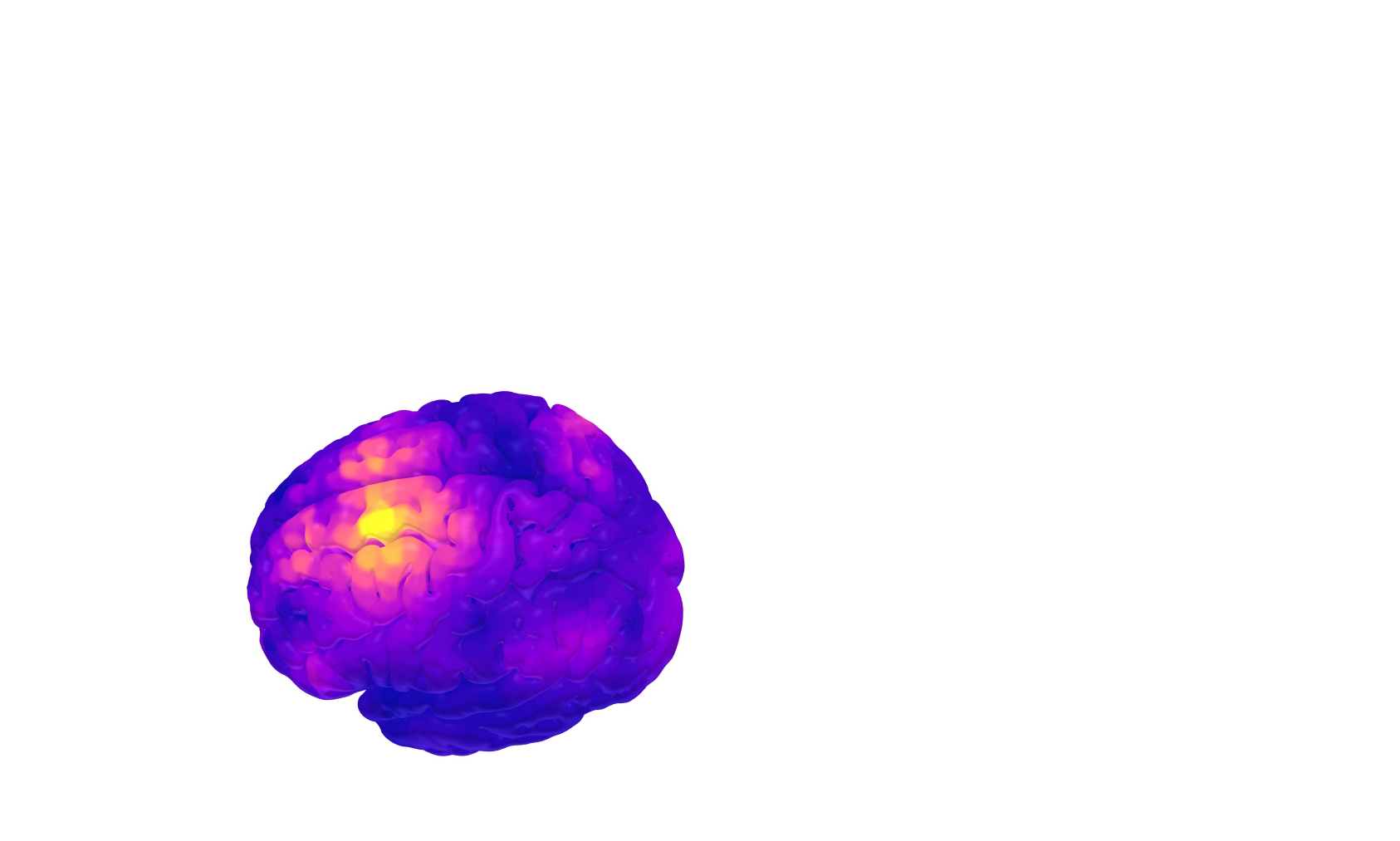}
        \end{minipage}\begin{minipage}{0.16\linewidth}
            \includegraphics[trim={5cm 1cm 14cm 8cm},clip,width=\linewidth]{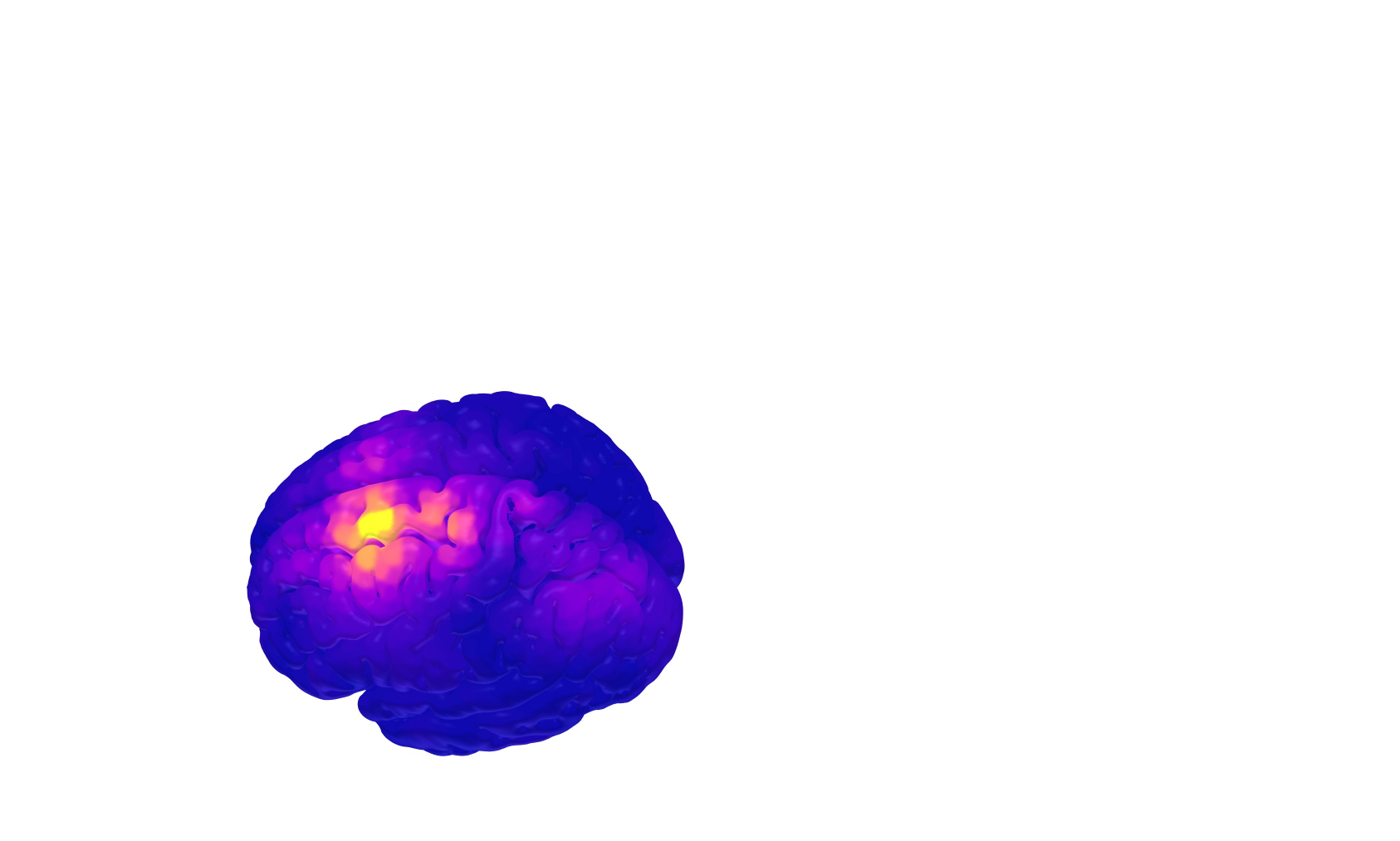}
        \end{minipage}
    \end{minipage}
    \noindent\rule[7pt]{\linewidth}{0.4pt}

    \begin{minipage}{0.02\linewidth}
        \rotatebox{90}{\scriptsize KF}
    \end{minipage}\begin{minipage}{0.98\linewidth}
        \begin{minipage}{0.16\linewidth}
            \includegraphics[trim={5cm 1cm 14cm 8cm},clip,width=\linewidth]{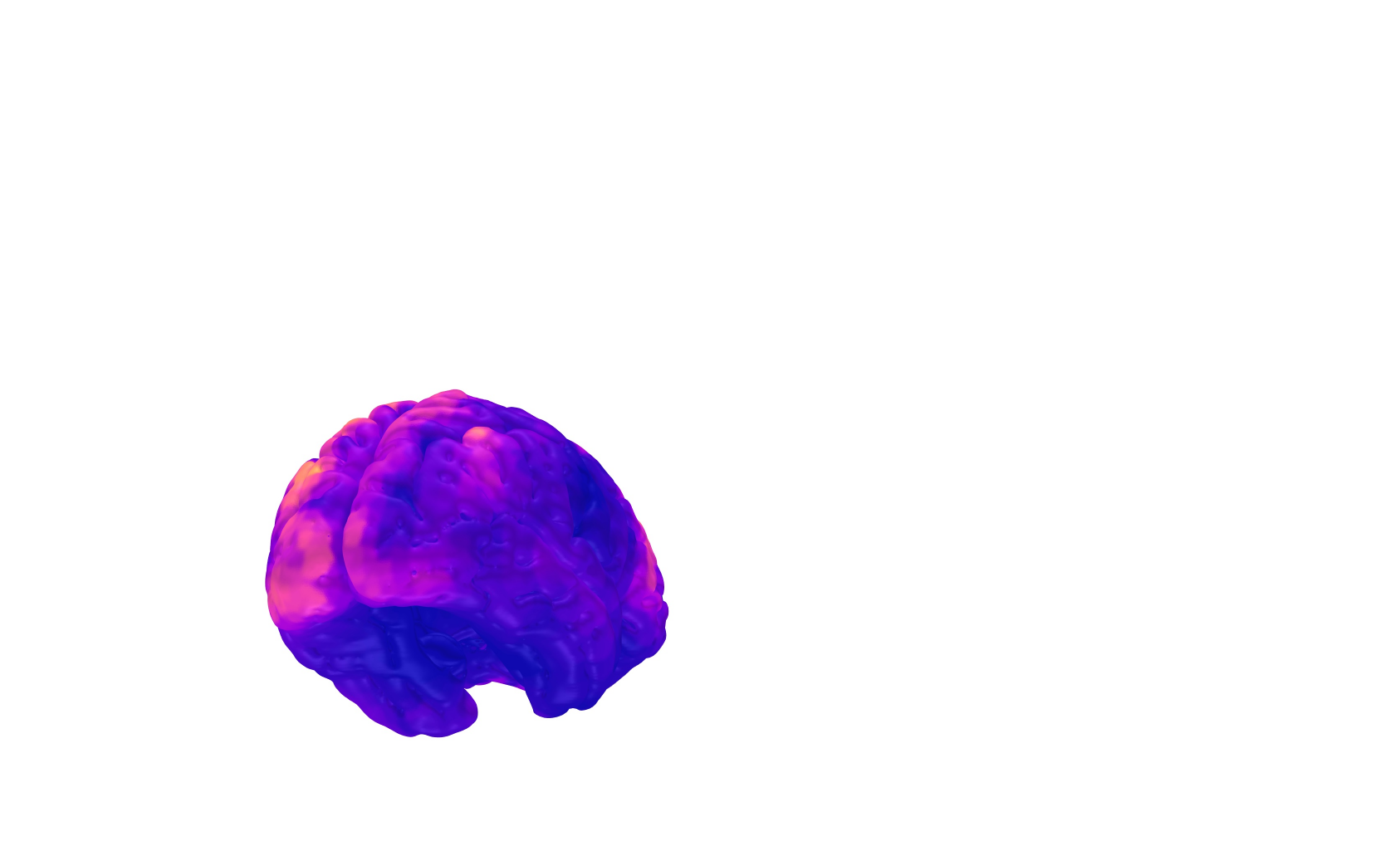}
        \end{minipage}\begin{minipage}{0.16\linewidth}
            \includegraphics[trim={5cm 1cm 14cm 8cm},clip,width=\linewidth]{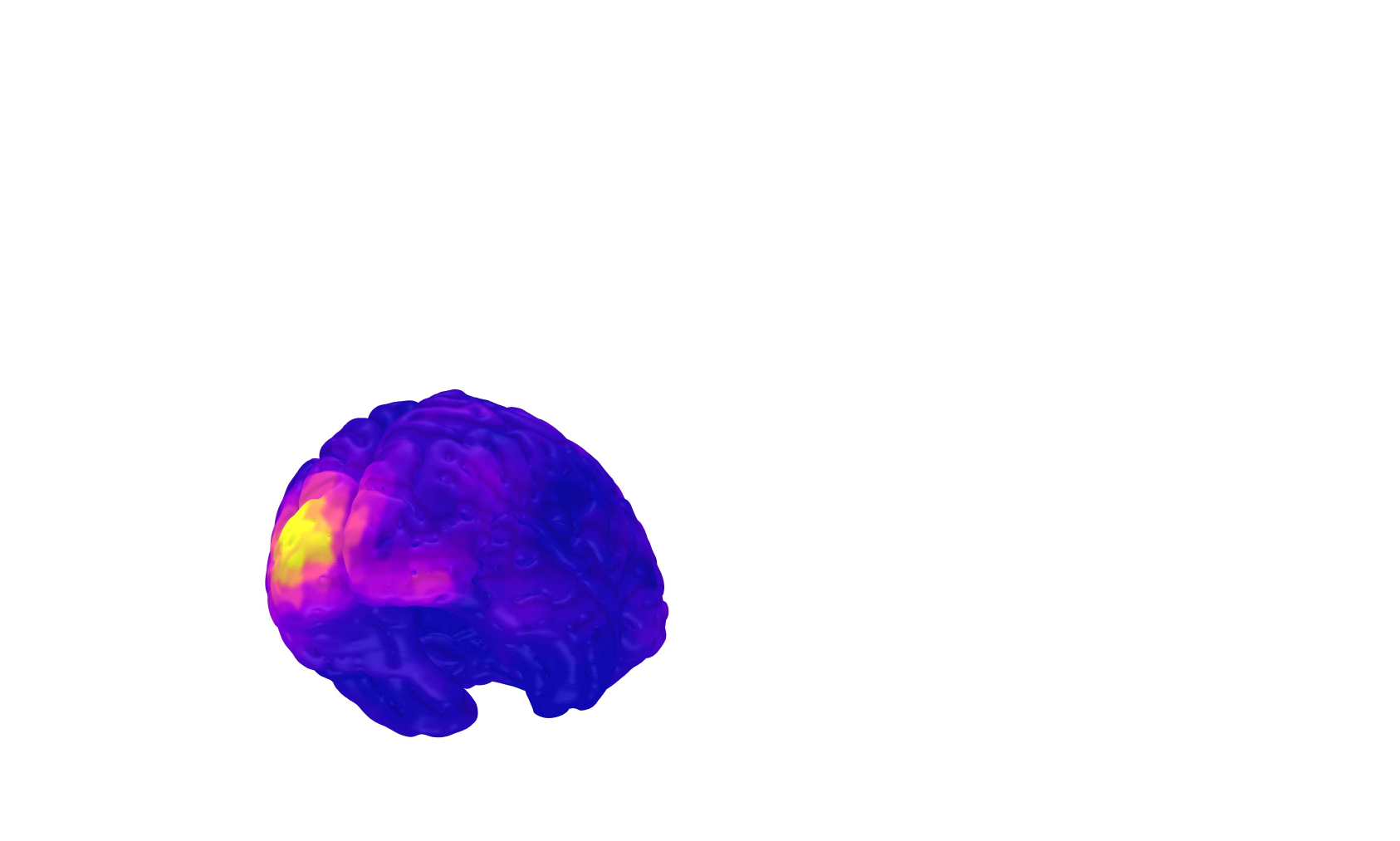}
        \end{minipage}\begin{minipage}{0.16\linewidth}
            \includegraphics[trim={5cm 1cm 14cm 8cm},clip,width=\linewidth]{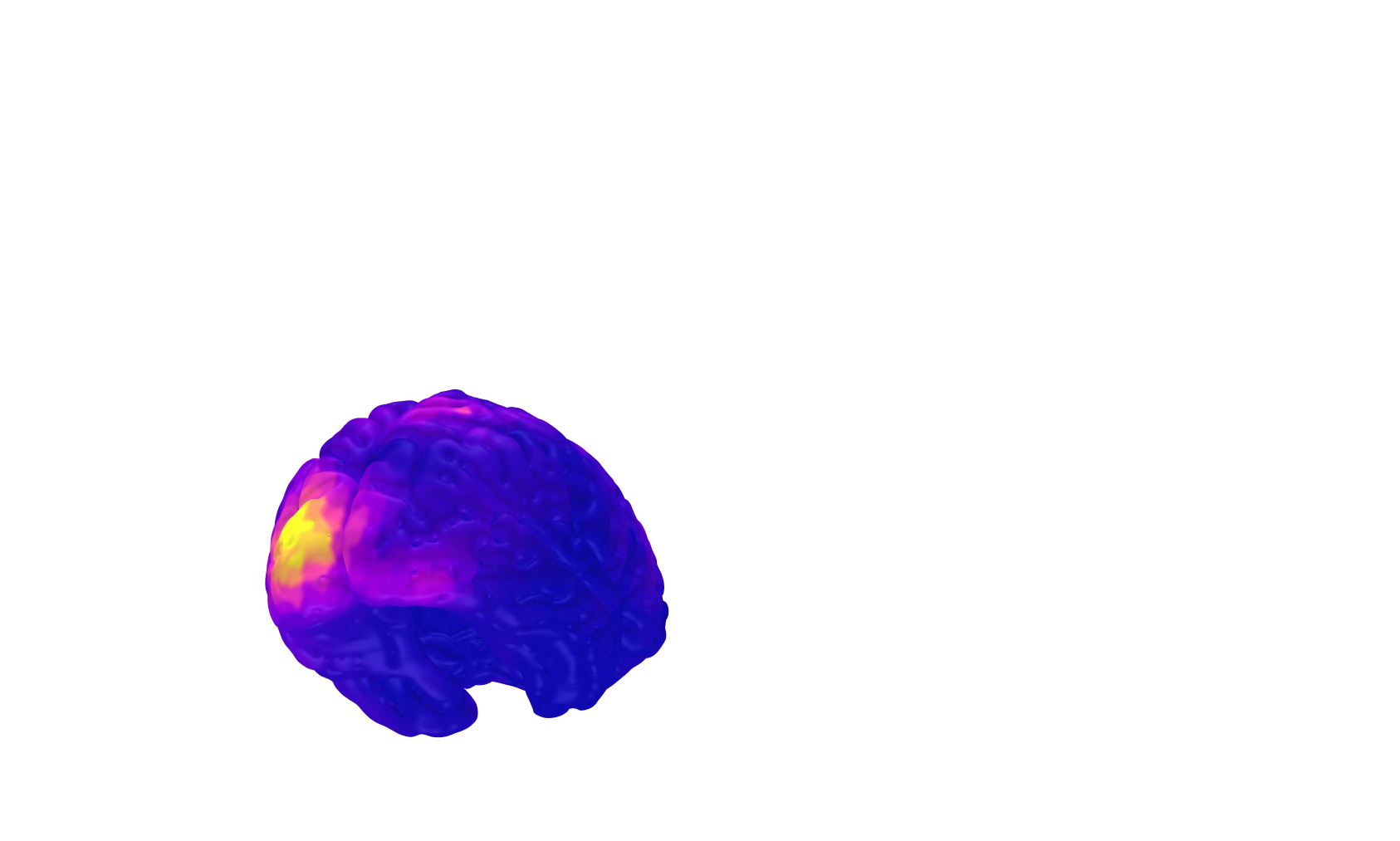}
        \end{minipage}\begin{minipage}{0.16\linewidth}
            \includegraphics[trim={5cm 1cm 14cm 8cm},clip,width=\linewidth]{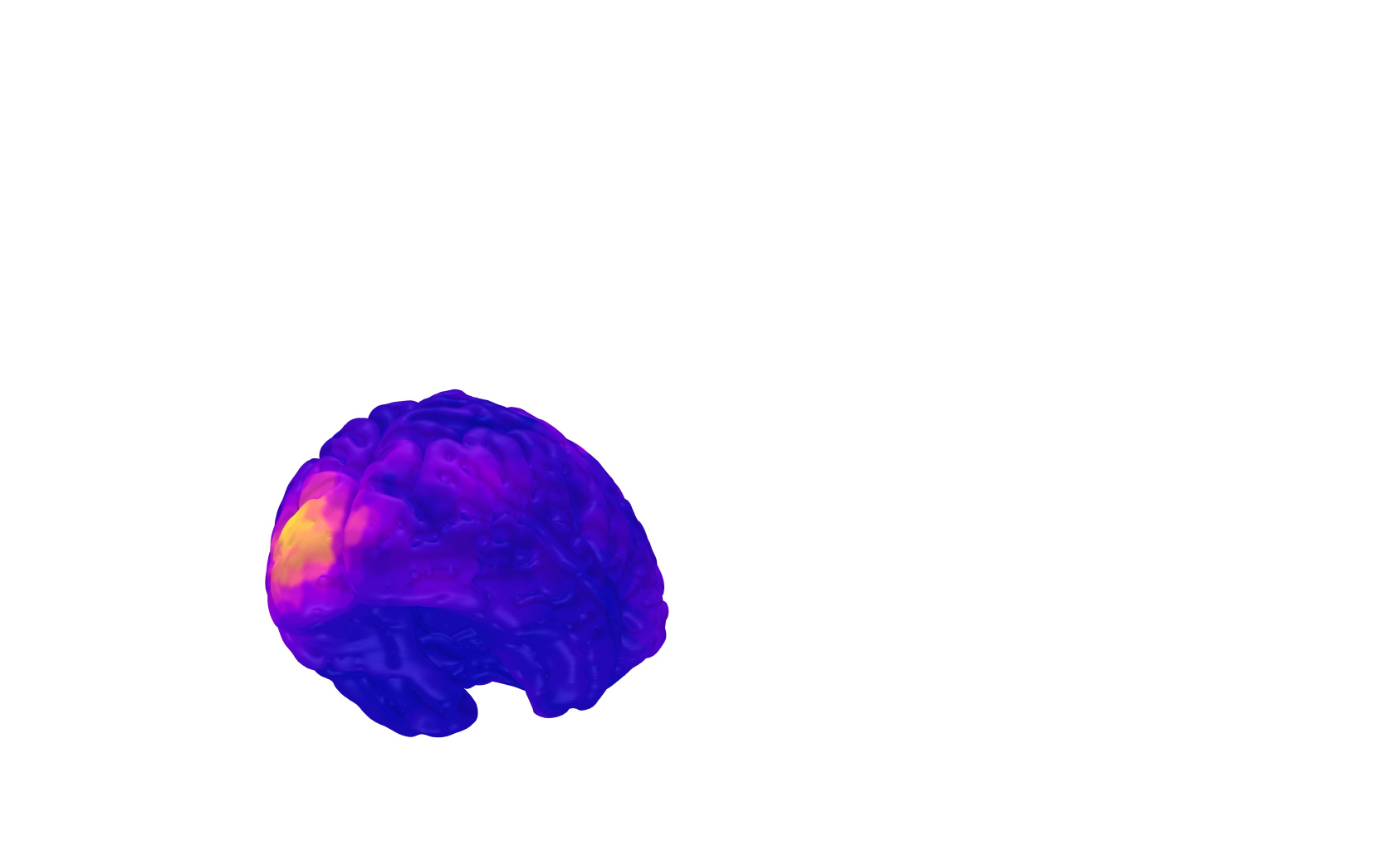}
        \end{minipage}\begin{minipage}{0.16\linewidth}
            \includegraphics[trim={5cm 1cm 14cm 8cm},clip,width=\linewidth]{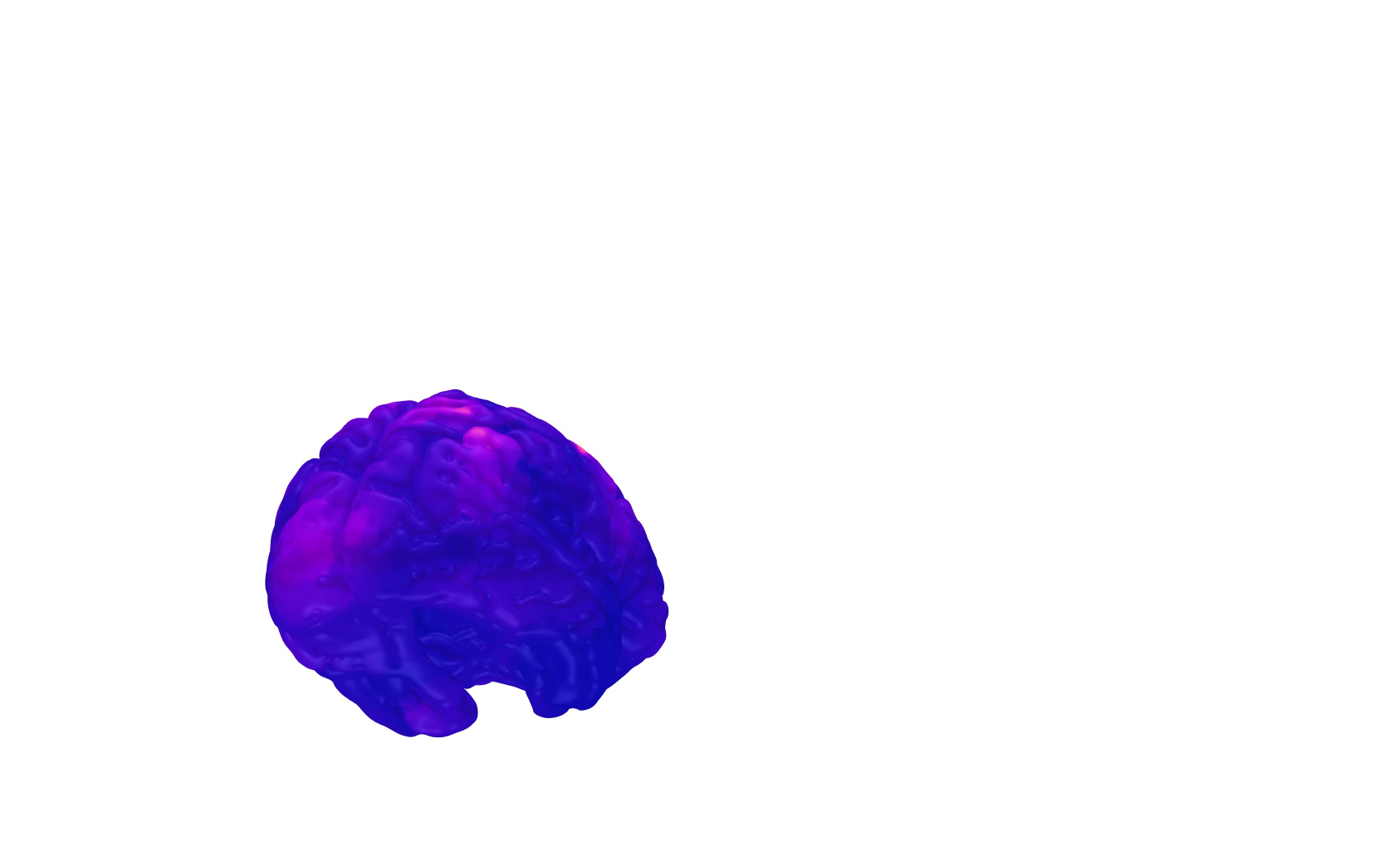}
        \end{minipage}\begin{minipage}{0.16\linewidth}
            \includegraphics[trim={5cm 1cm 14cm 8cm},clip,width=\linewidth]{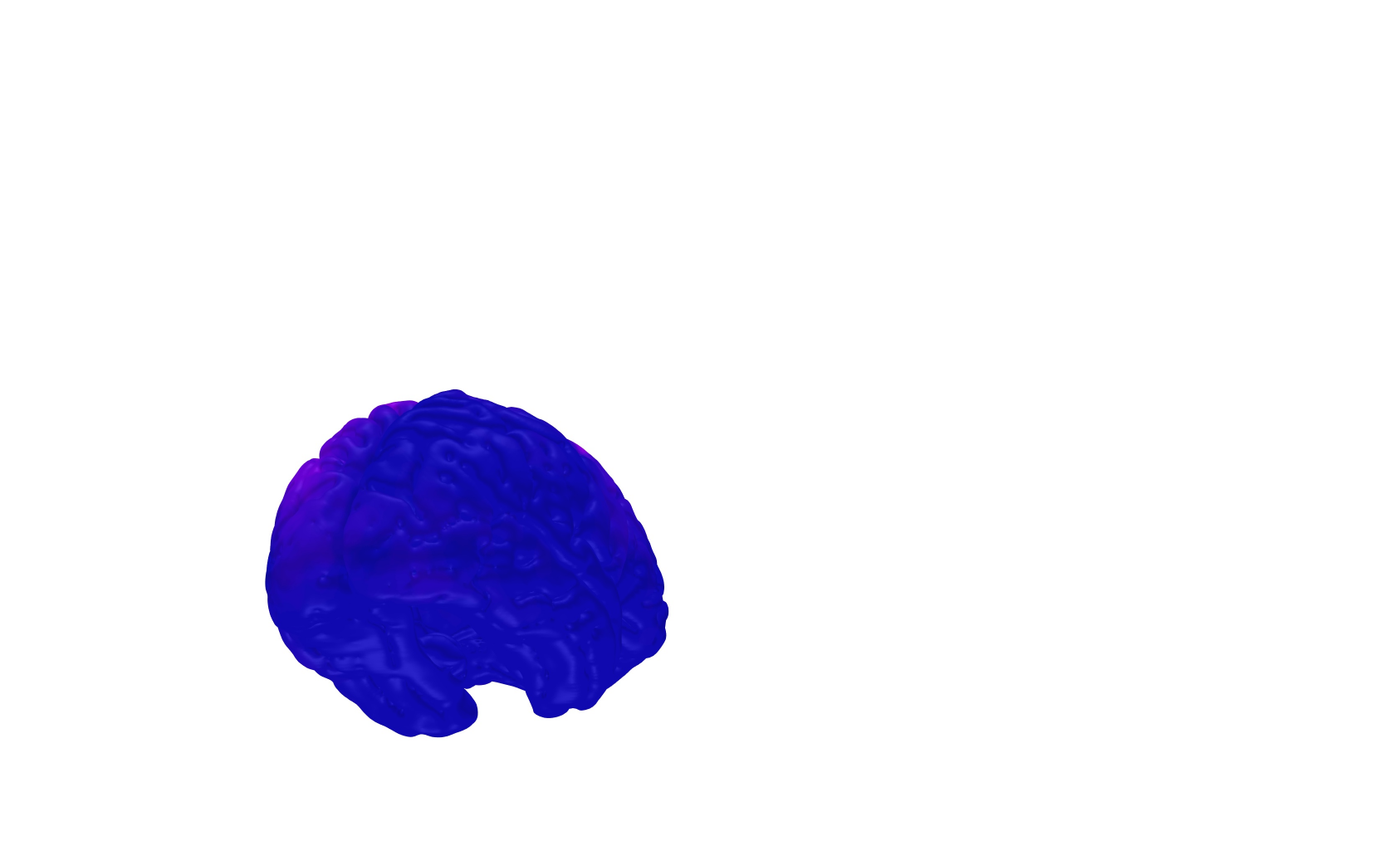}
        \end{minipage}

        \begin{minipage}{0.16\linewidth}
            \includegraphics[trim={5cm 1cm 14cm 8cm},clip,width=\linewidth]{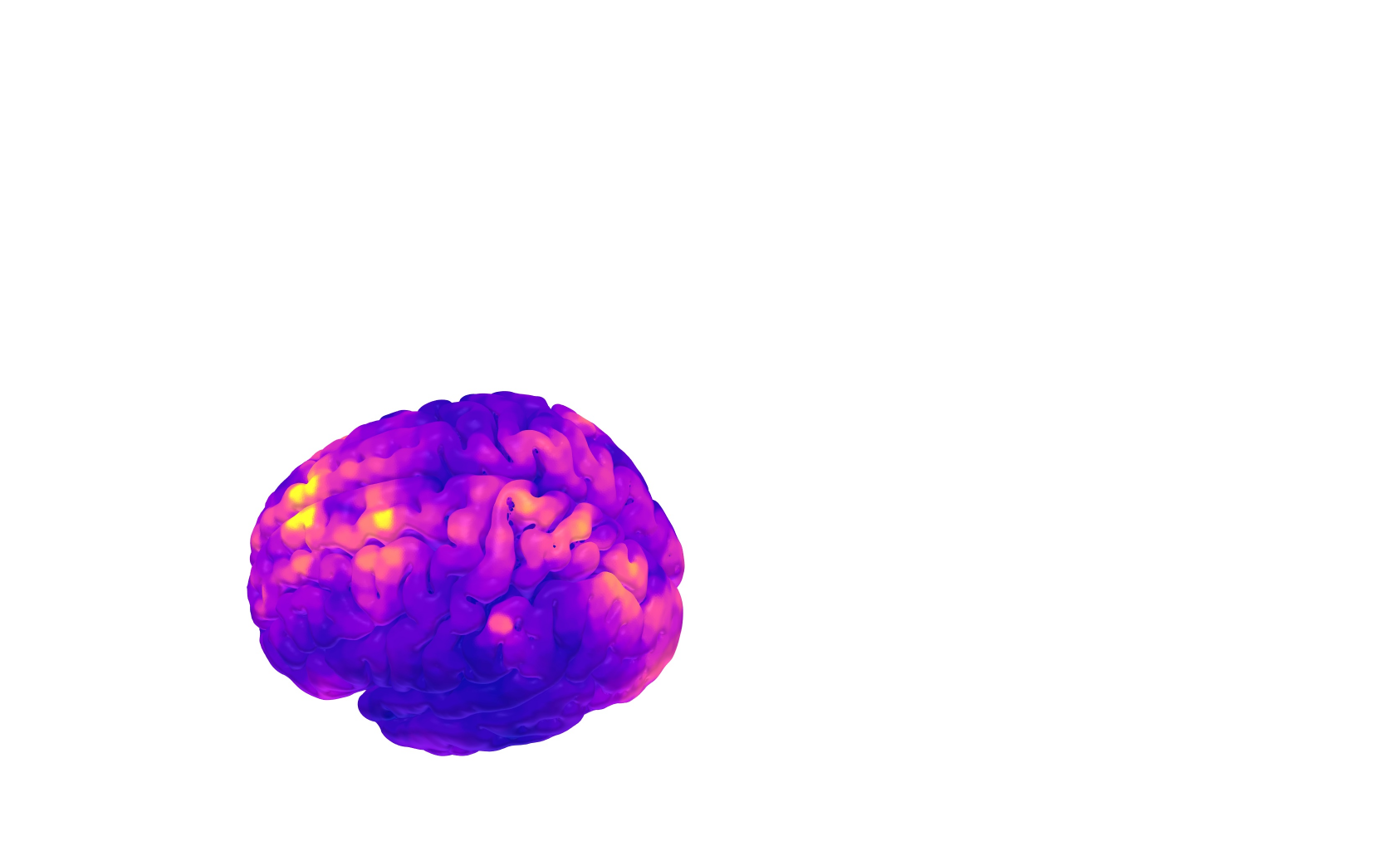}
        \end{minipage}\begin{minipage}{0.16\linewidth}
            \includegraphics[trim={5cm 1cm 14cm 8cm},clip,width=\linewidth]{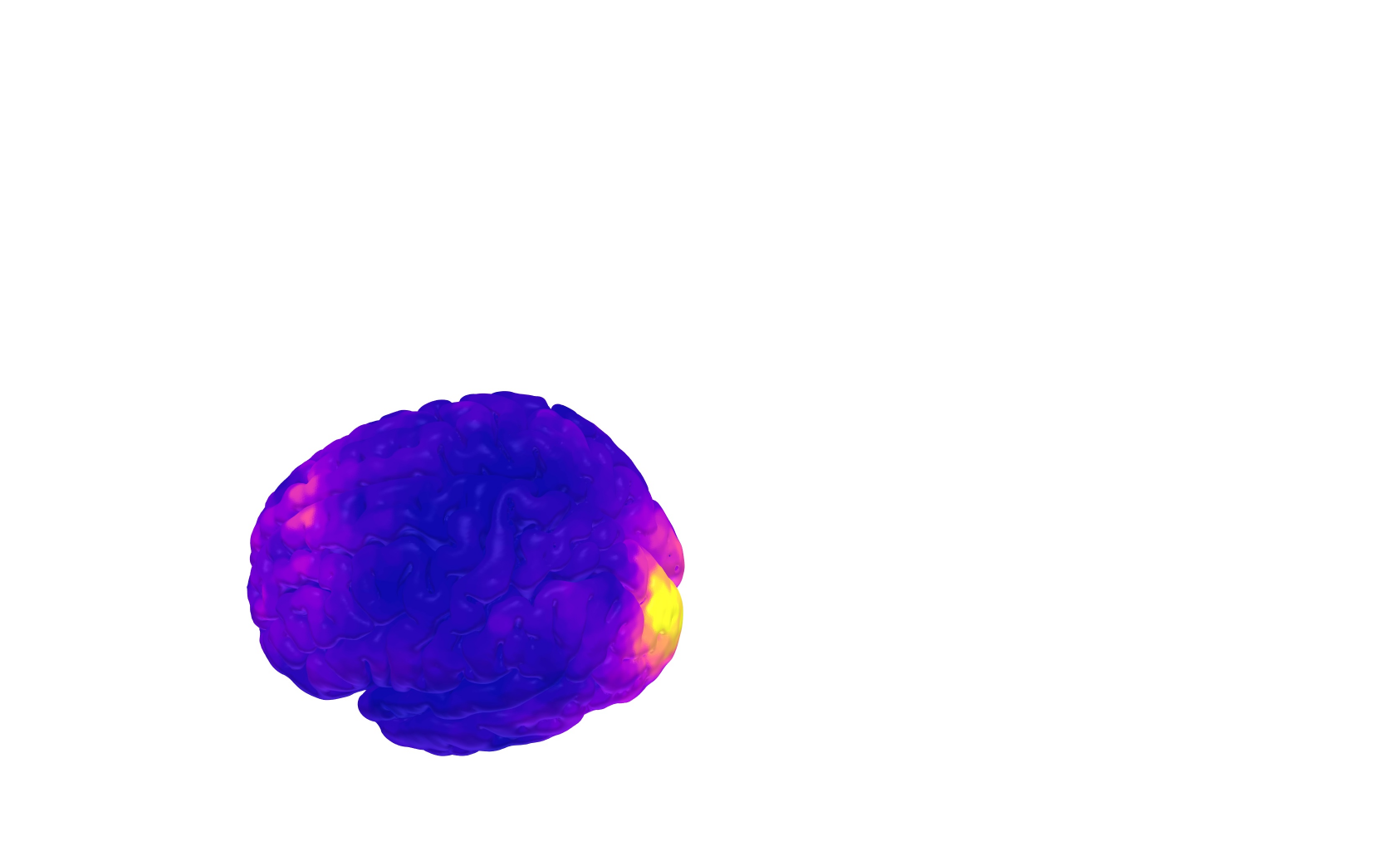}
        \end{minipage}\begin{minipage}{0.16\linewidth}
            \includegraphics[trim={5cm 1cm 14cm 8cm},clip,width=\linewidth]{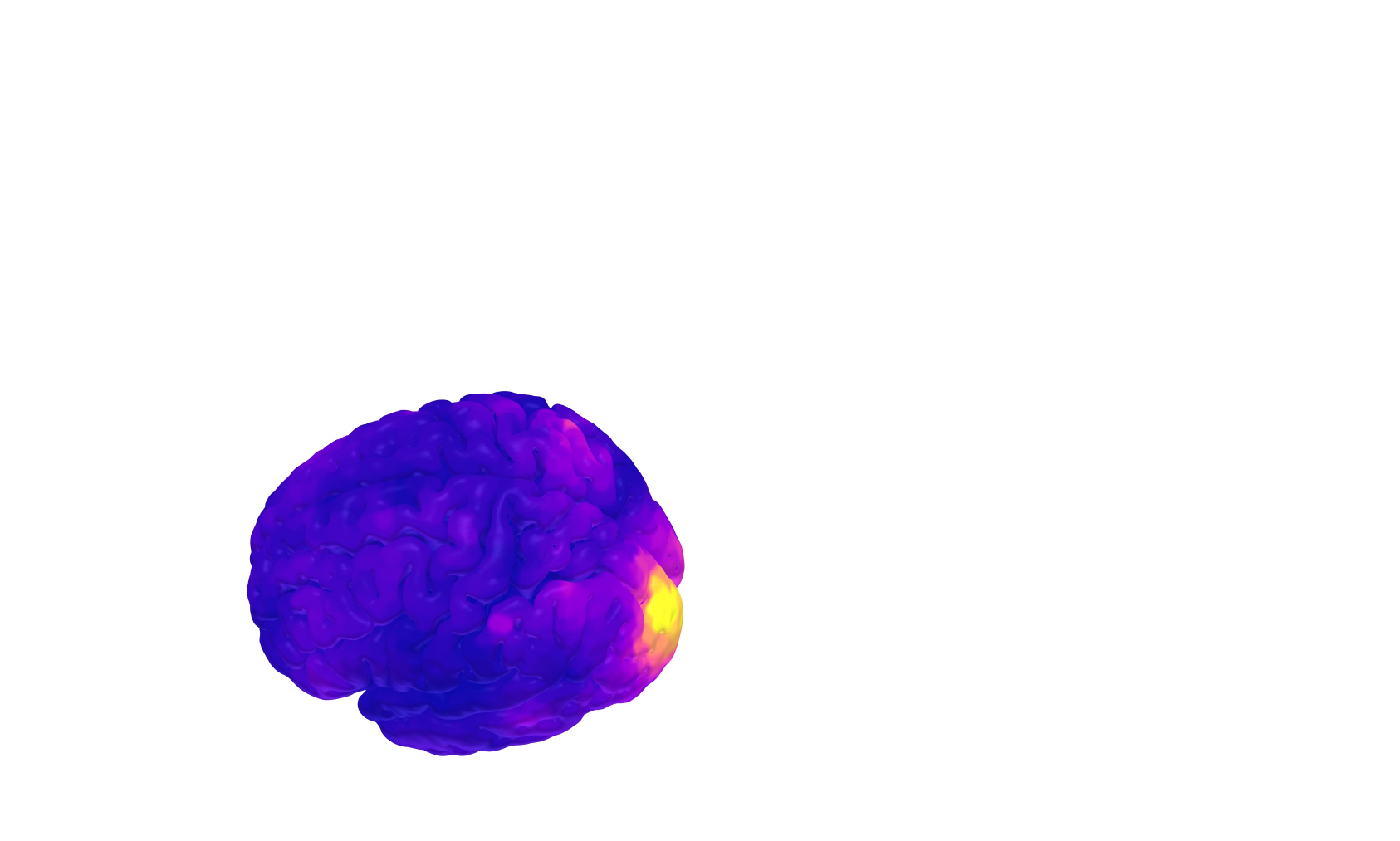}
        \end{minipage}\begin{minipage}{0.16\linewidth}
            \includegraphics[trim={5cm 1cm 14cm 8cm},clip,width=\linewidth]{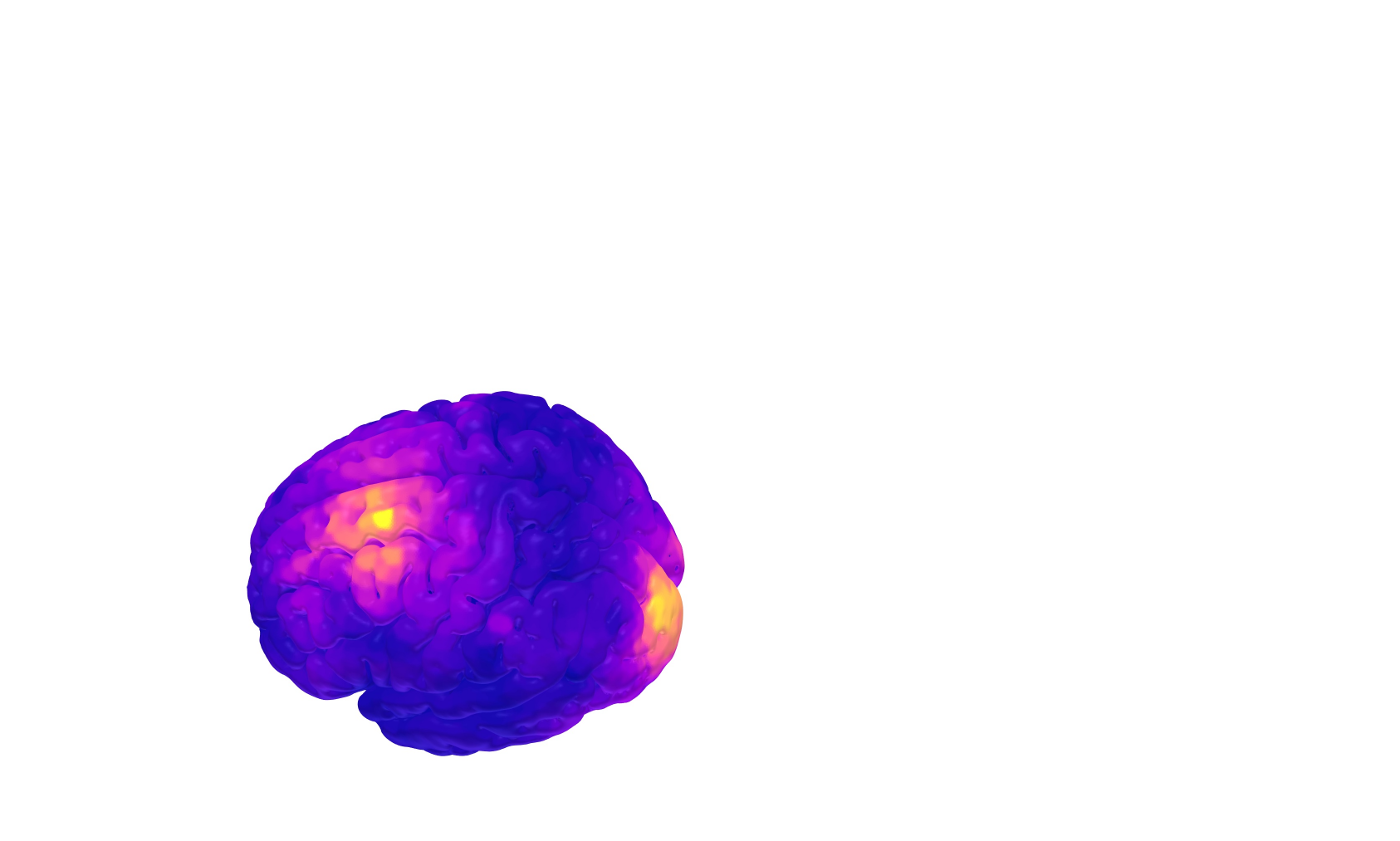}
        \end{minipage}\begin{minipage}{0.16\linewidth}
            \includegraphics[trim={5cm 1cm 14cm 8cm},clip,width=\linewidth]{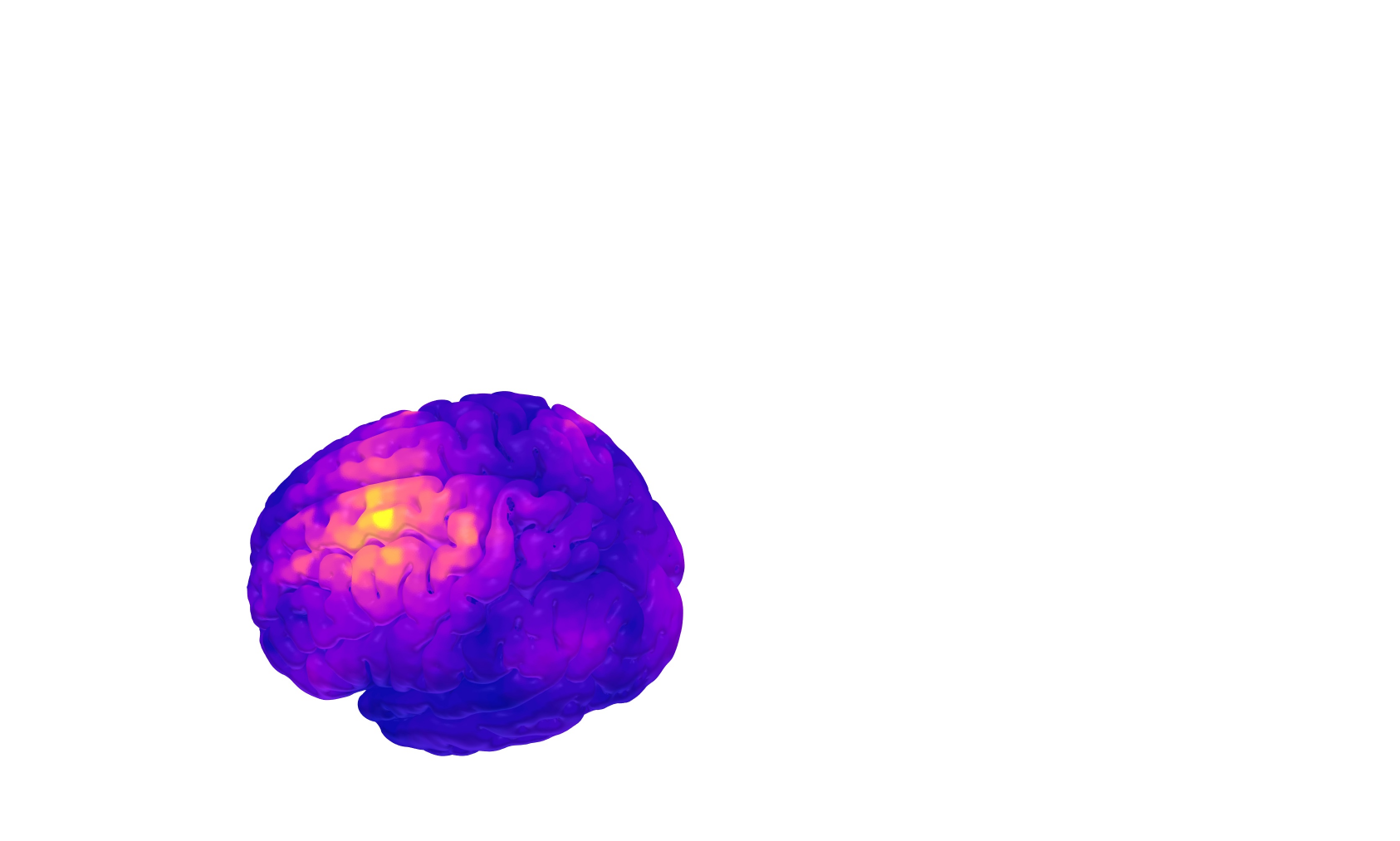}
        \end{minipage}\begin{minipage}{0.16\linewidth}
            \includegraphics[trim={5cm 1cm 14cm 8cm},clip,width=\linewidth]{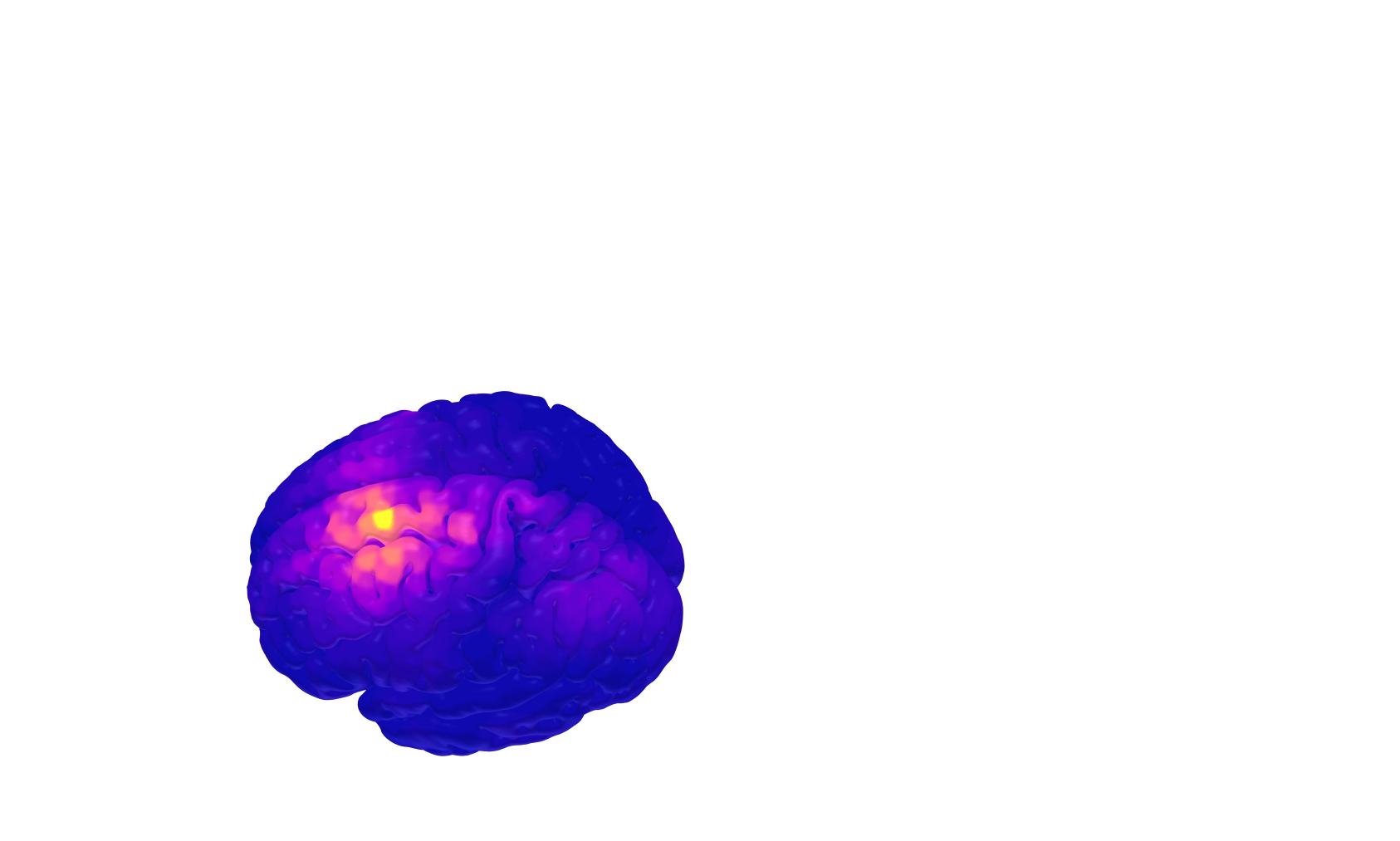}
        \end{minipage}
    \end{minipage}
    \noindent\rule[7pt]{\linewidth}{0.4pt}

    \begin{minipage}{0.02\linewidth}
        \rotatebox{90}{\scriptsize LCMV}
    \end{minipage}\begin{minipage}{0.98\linewidth}
        \begin{minipage}{0.16\linewidth}
            \includegraphics[trim={5cm 1cm 14cm 8cm},clip,width=\linewidth]{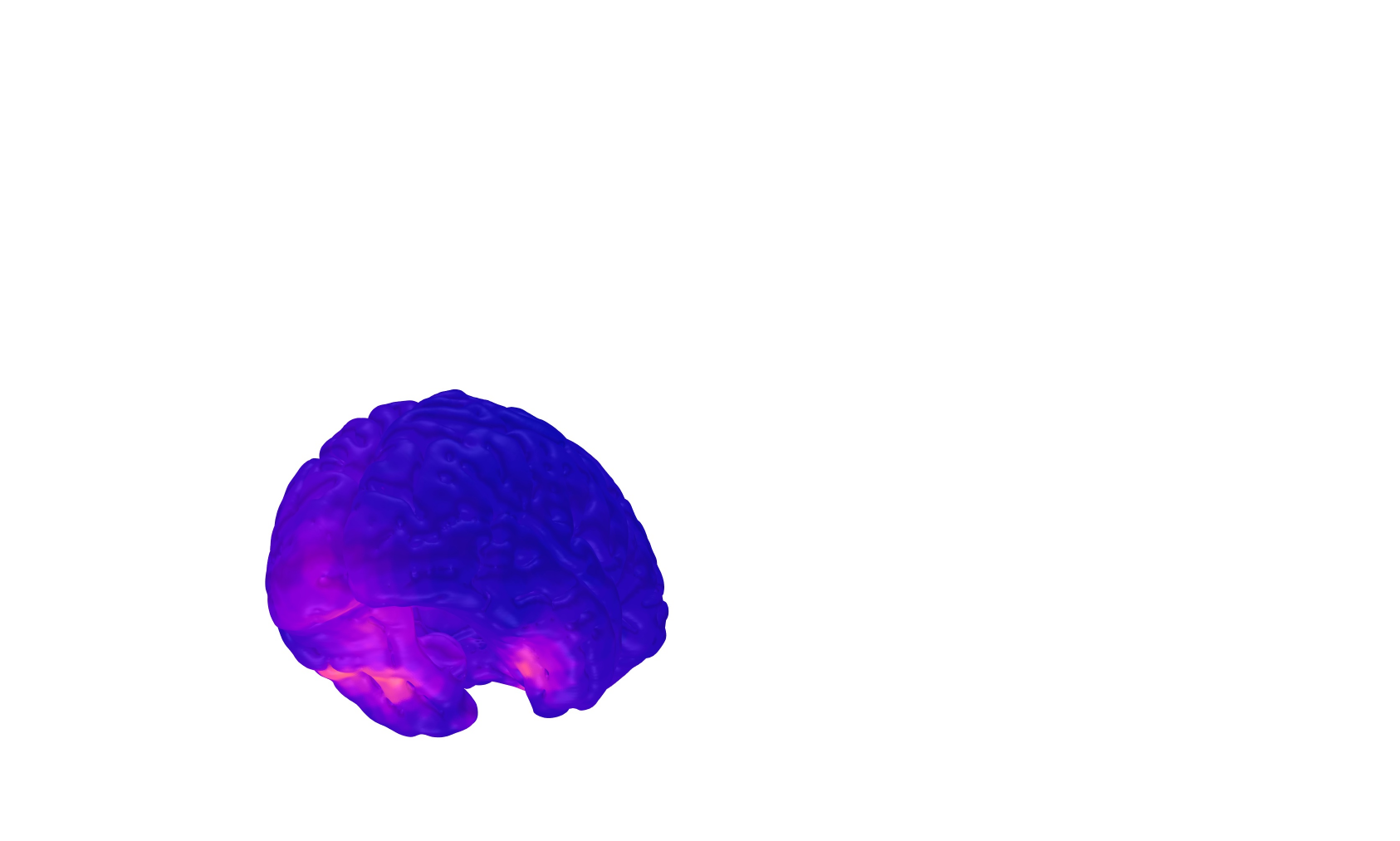}
        \end{minipage}\begin{minipage}{0.16\linewidth}
            \includegraphics[trim={5cm 1cm 14cm 8cm},clip,width=\linewidth]{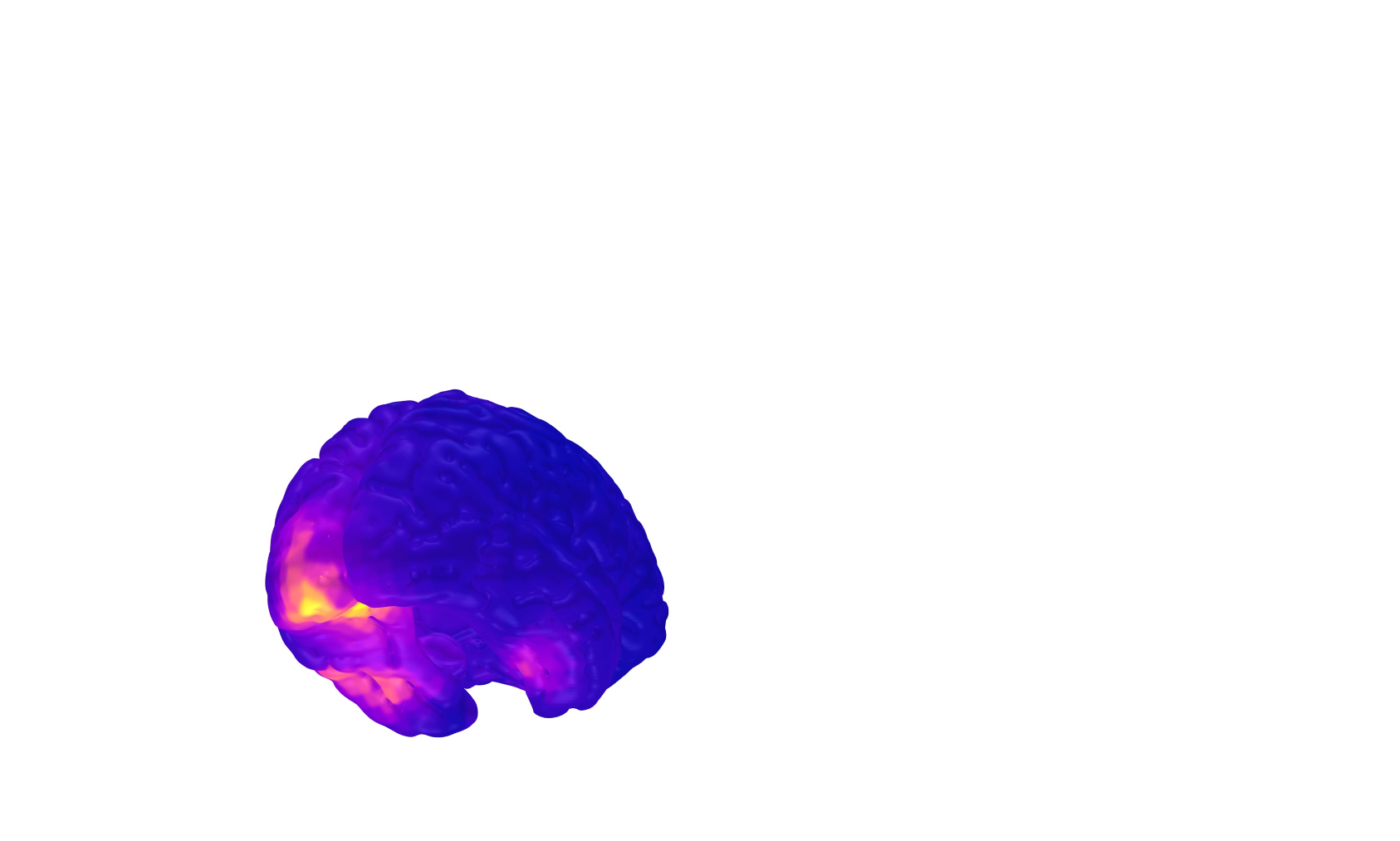}
        \end{minipage}\begin{minipage}{0.16\linewidth}
            \includegraphics[trim={5cm 1cm 14cm 8cm},clip,width=\linewidth]{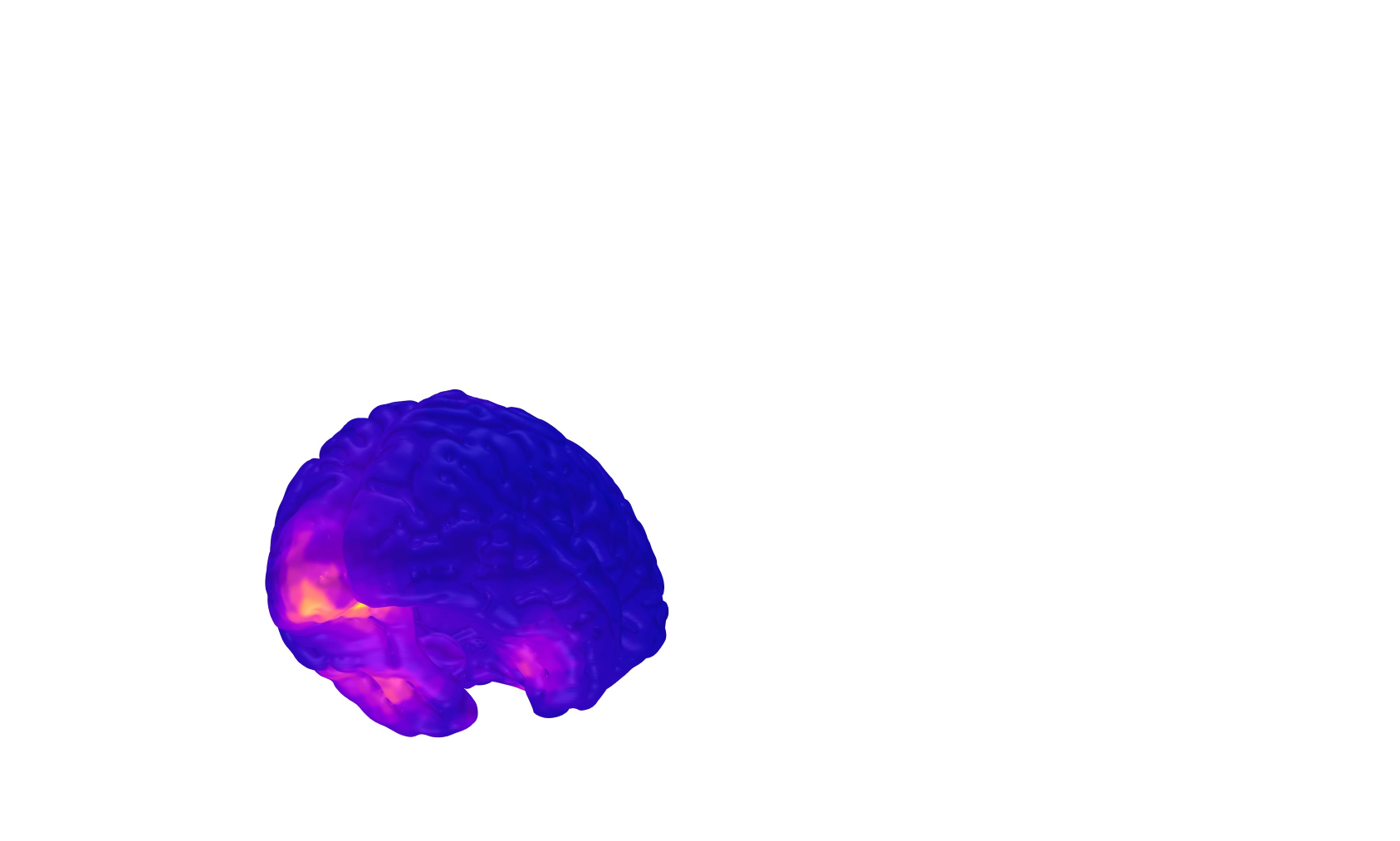}
        \end{minipage}\begin{minipage}{0.16\linewidth}
            \includegraphics[trim={5cm 1cm 14cm 8cm},clip,width=\linewidth]{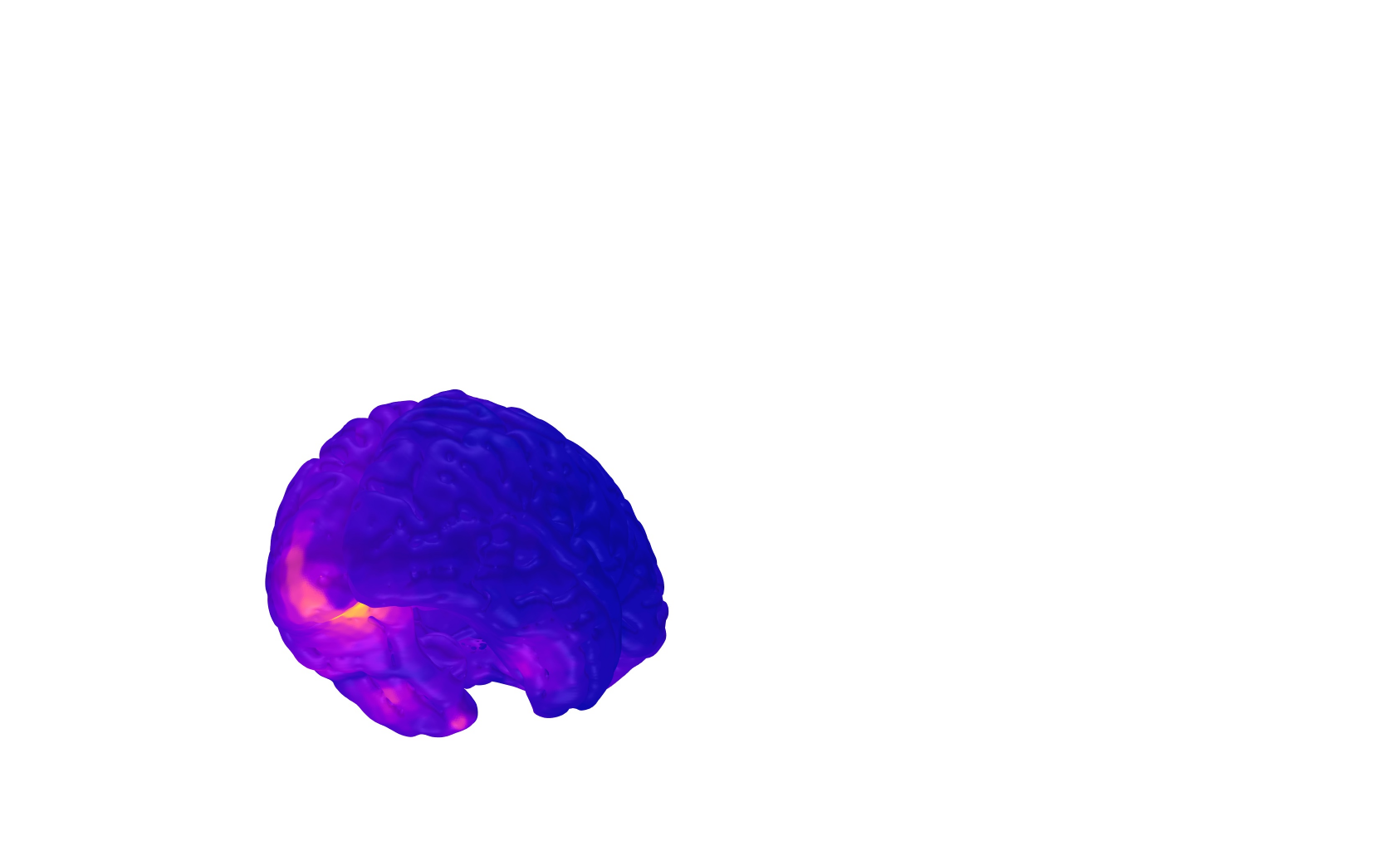}
        \end{minipage}\begin{minipage}{0.16\linewidth}
            \includegraphics[trim={5cm 1cm 14cm 8cm},clip,width=\linewidth]{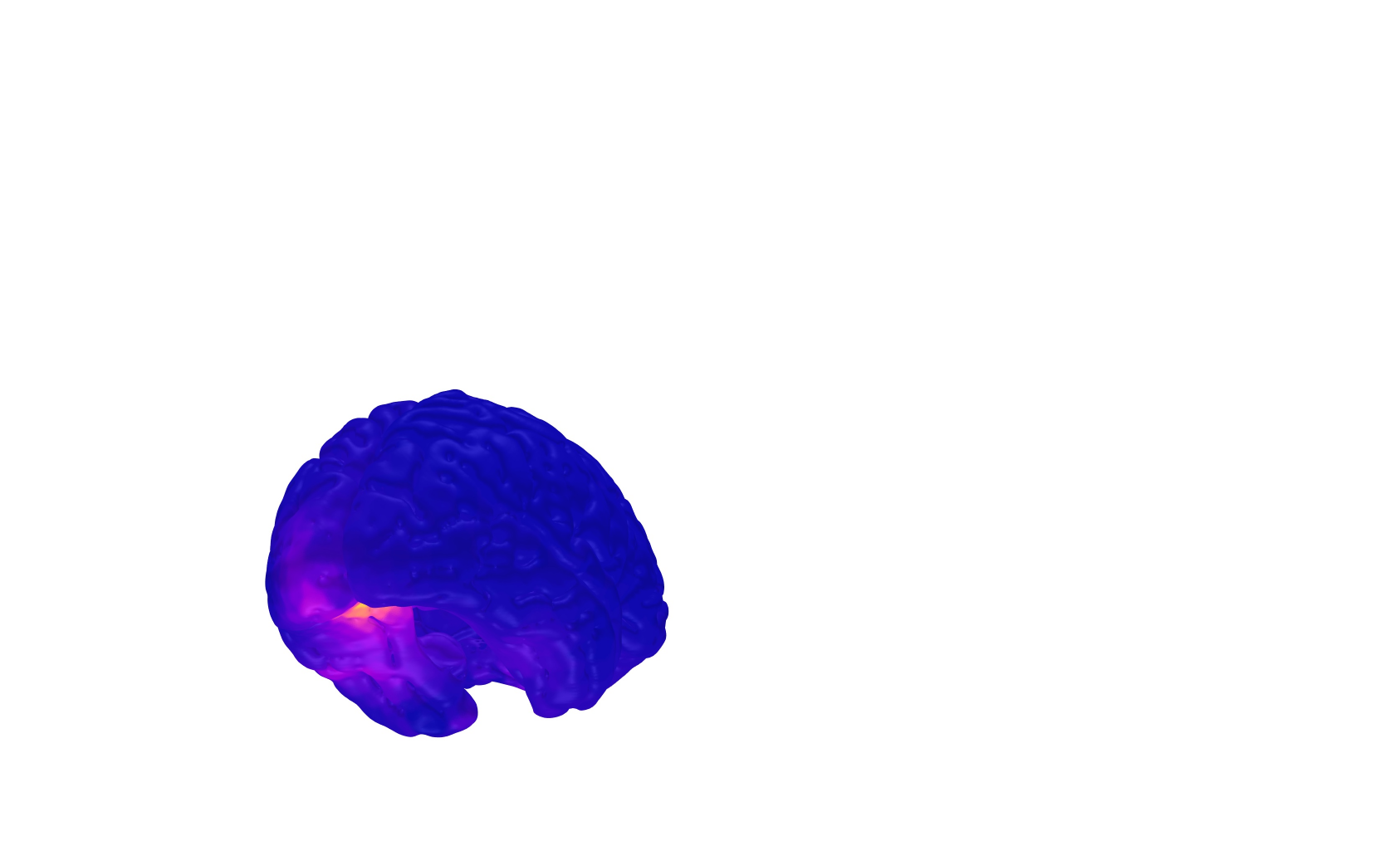}
        \end{minipage}\begin{minipage}{0.16\linewidth}
            \includegraphics[trim={5cm 1cm 14cm 8cm},clip,width=\linewidth]{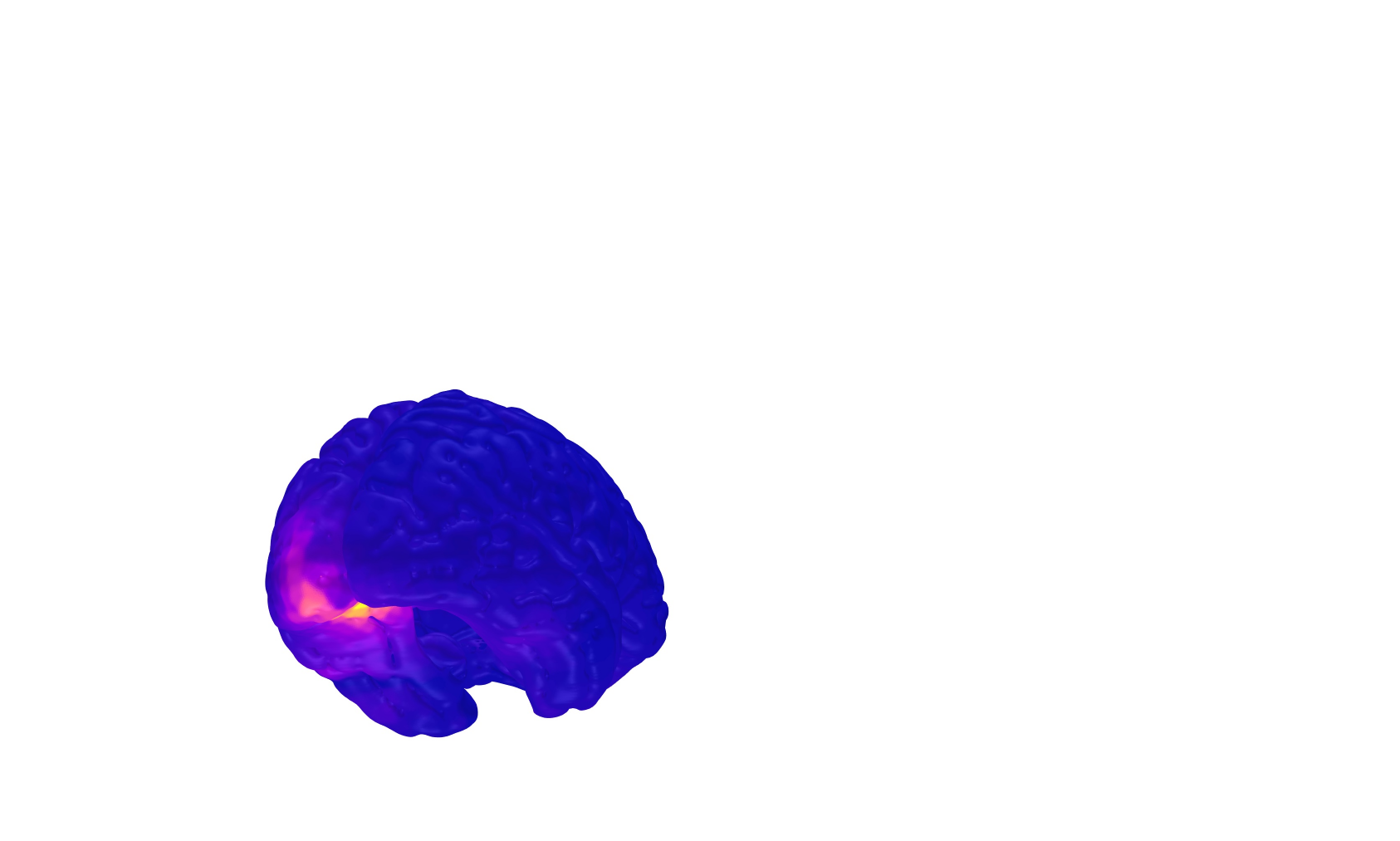}
        \end{minipage}

        \begin{minipage}{0.16\linewidth}
            \includegraphics[trim={5cm 1cm 14cm 8cm},clip,width=\linewidth]{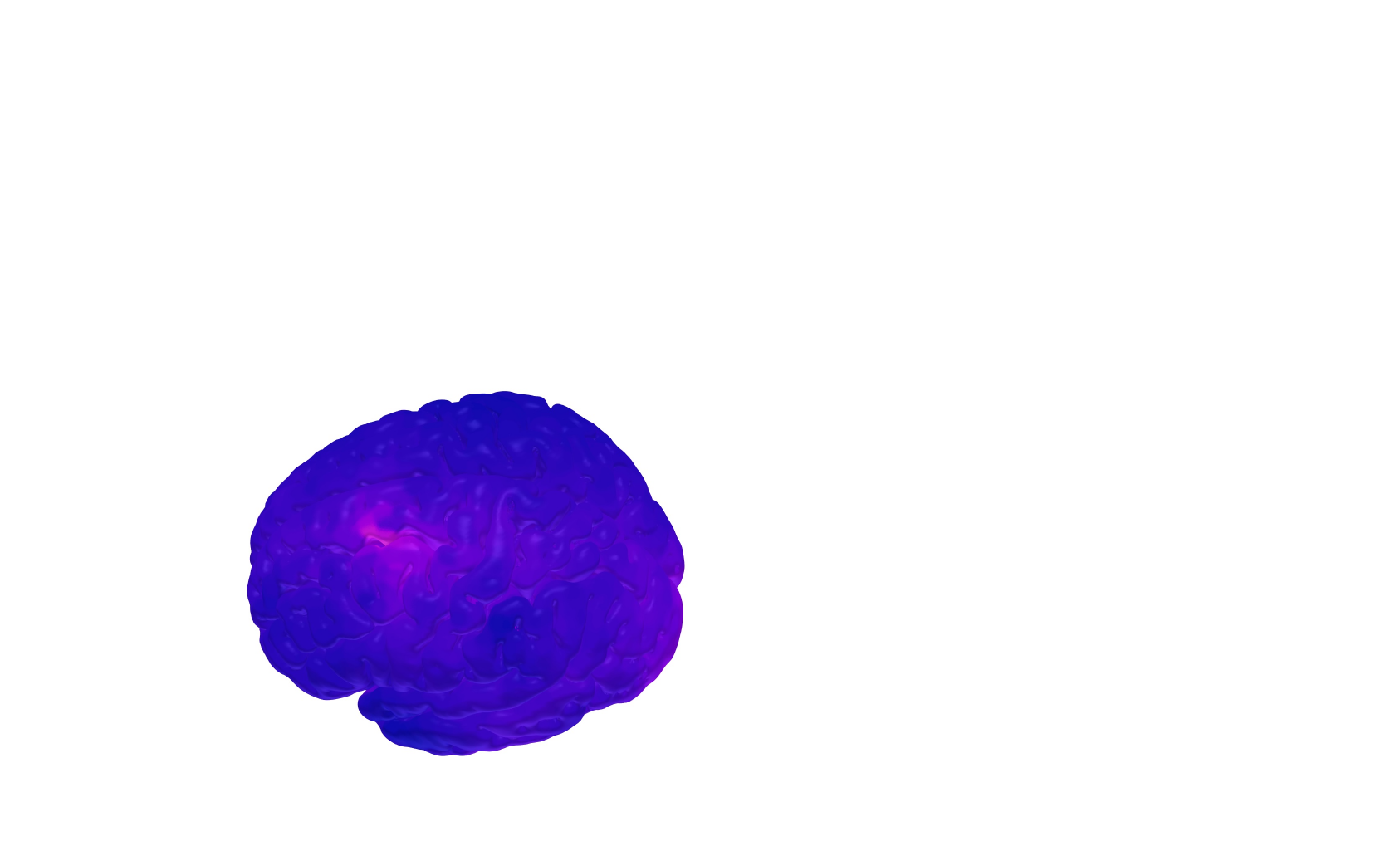}
        \end{minipage}\begin{minipage}{0.16\linewidth}
            \includegraphics[trim={5cm 1cm 14cm 8cm},clip,width=\linewidth]{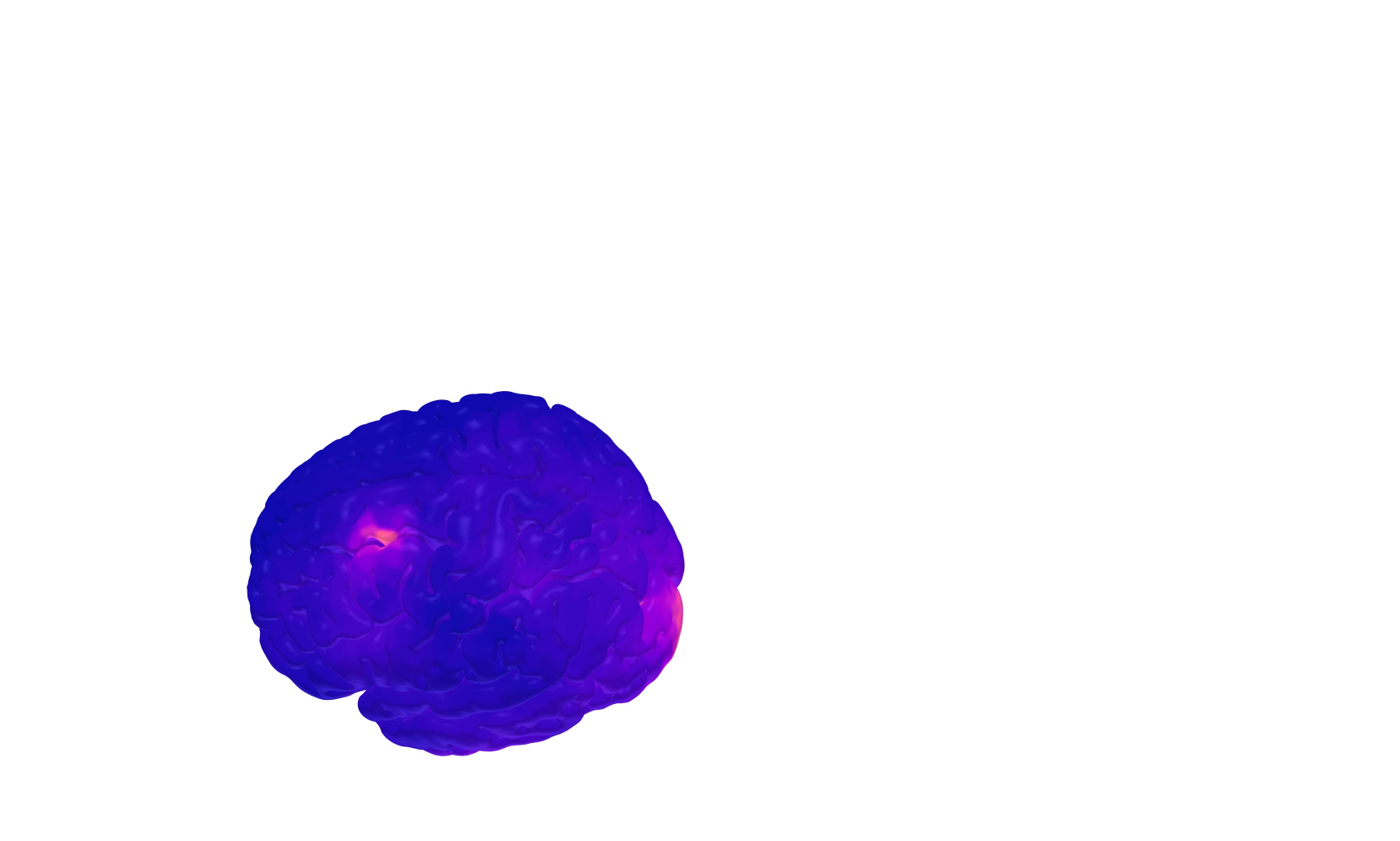}
        \end{minipage}\begin{minipage}{0.16\linewidth}
            \includegraphics[trim={5cm 1cm 14cm 8cm},clip,width=\linewidth]{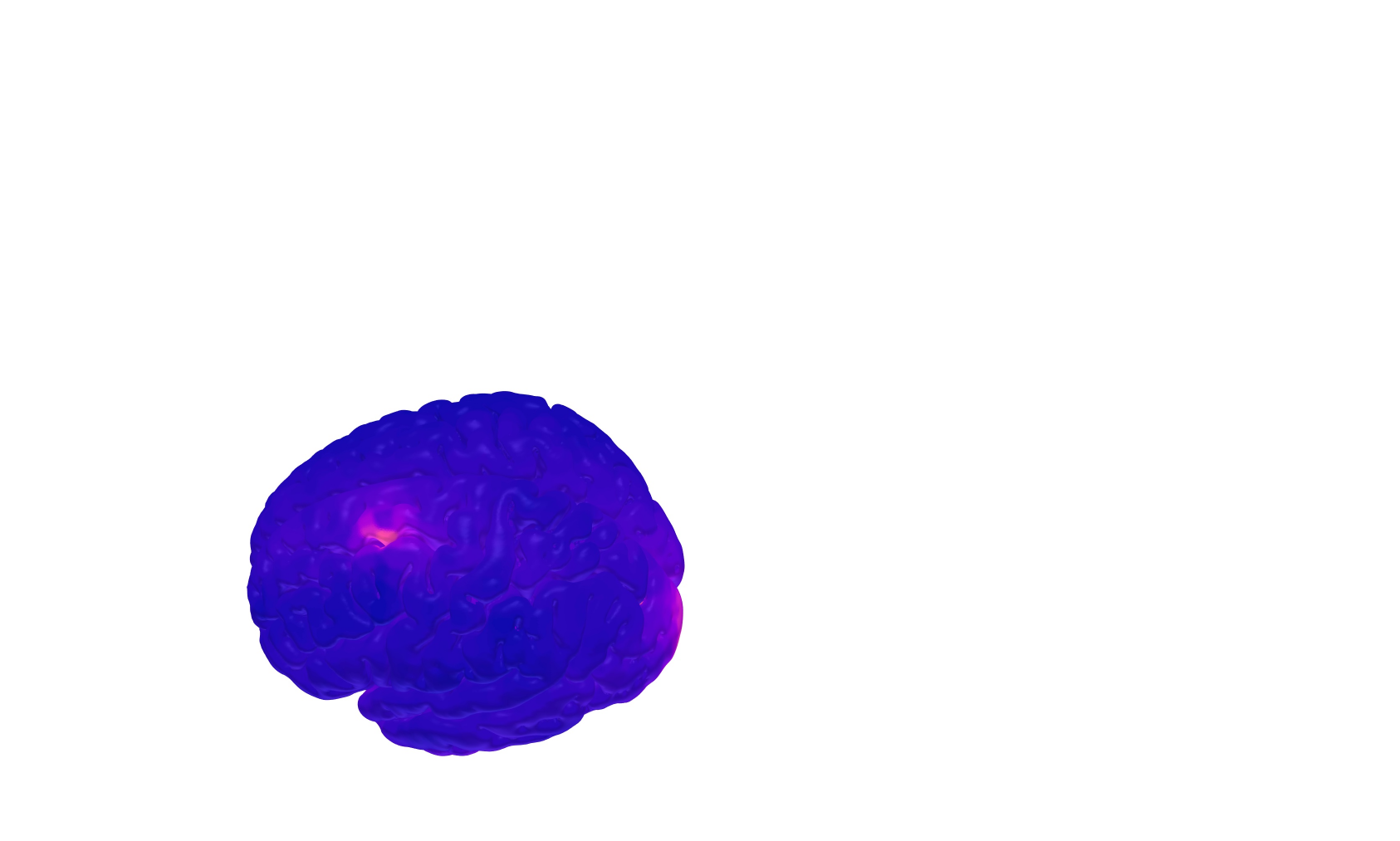}
        \end{minipage}\begin{minipage}{0.16\linewidth}
            \includegraphics[trim={5cm 1cm 14cm 8cm},clip,width=\linewidth]{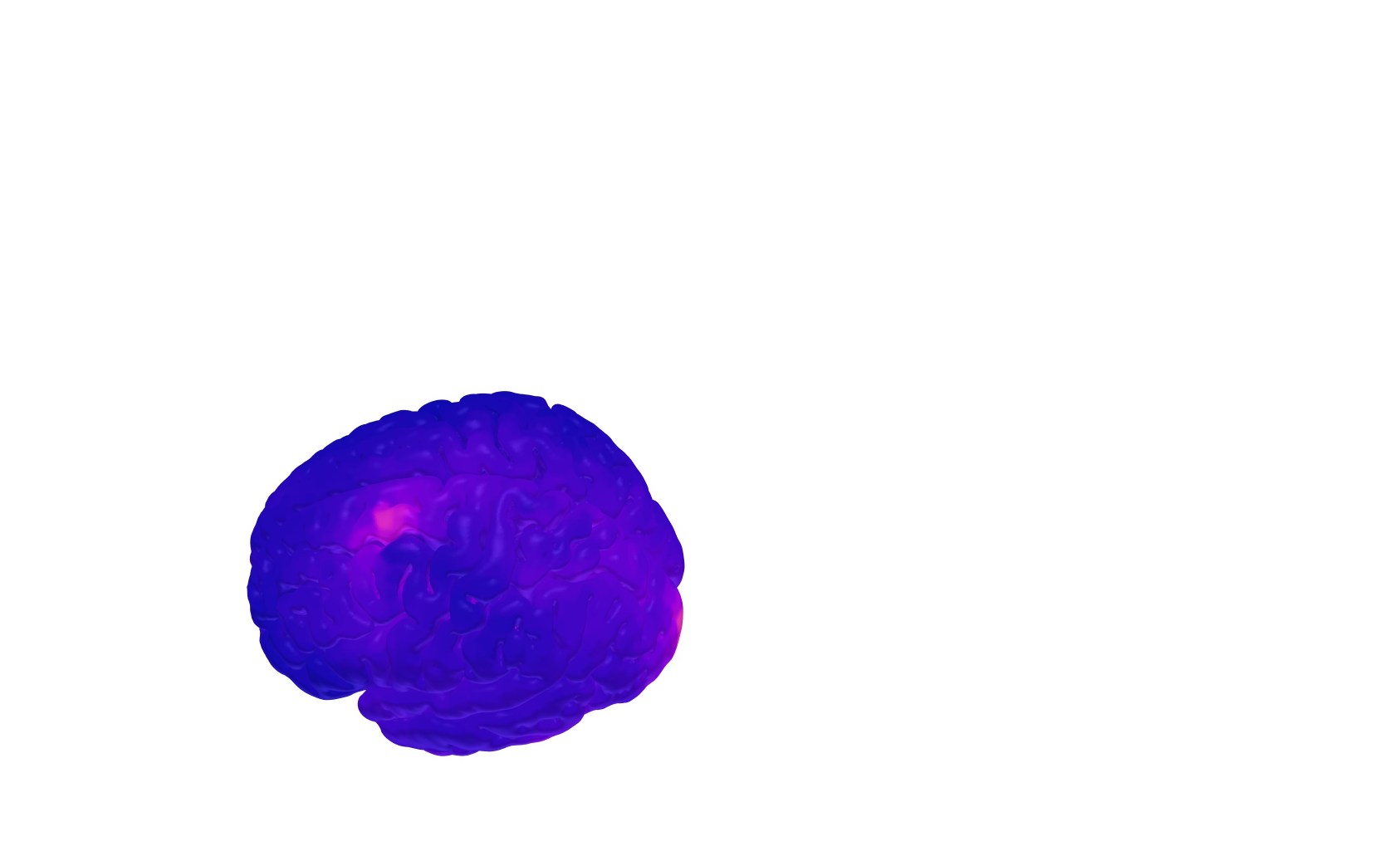}
        \end{minipage}\begin{minipage}{0.16\linewidth}
            \includegraphics[trim={5cm 1cm 14cm 8cm},clip,width=\linewidth]{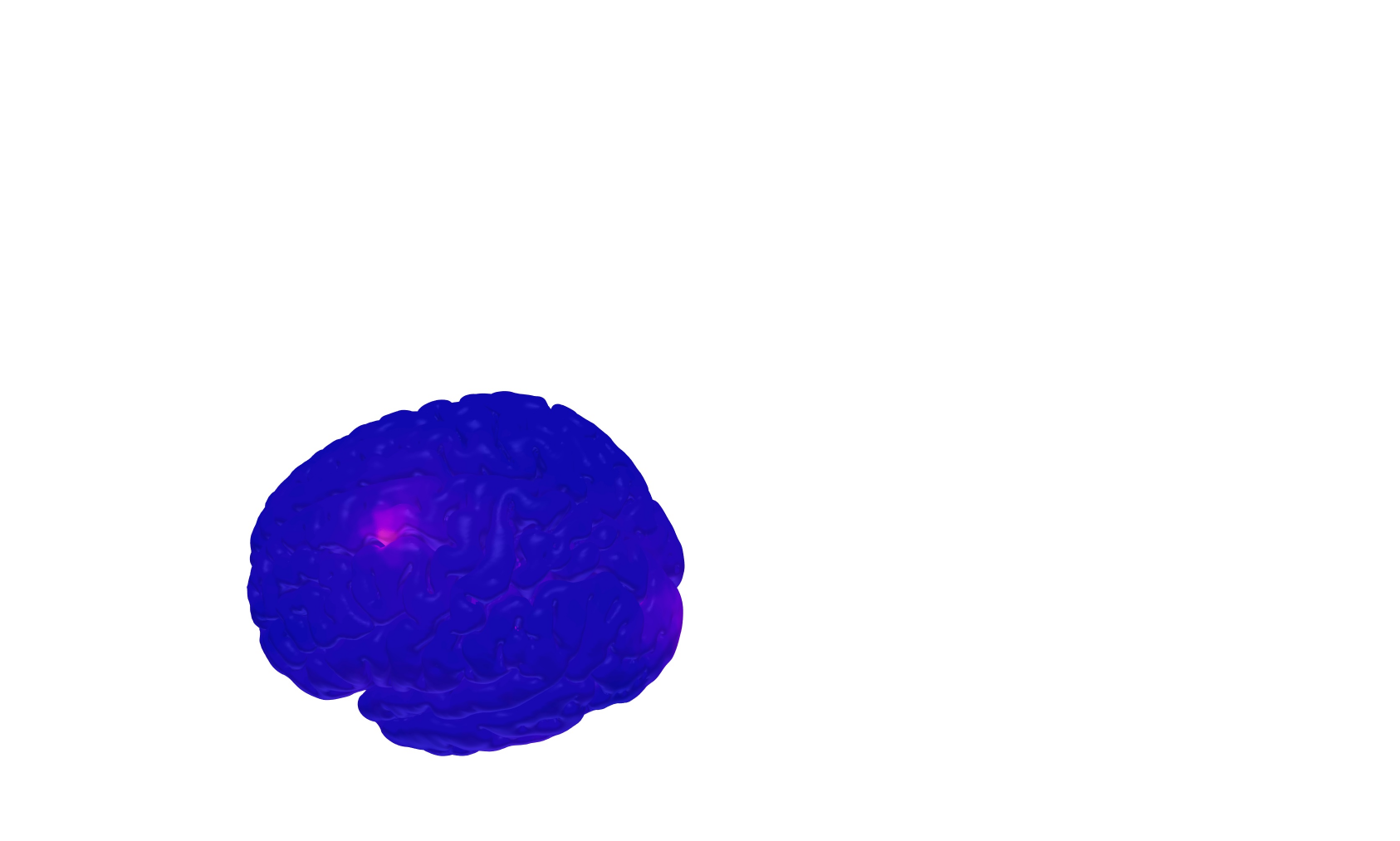}
        \end{minipage}\begin{minipage}{0.16\linewidth}
            \includegraphics[trim={5cm 1cm 14cm 8cm},clip,width=\linewidth]{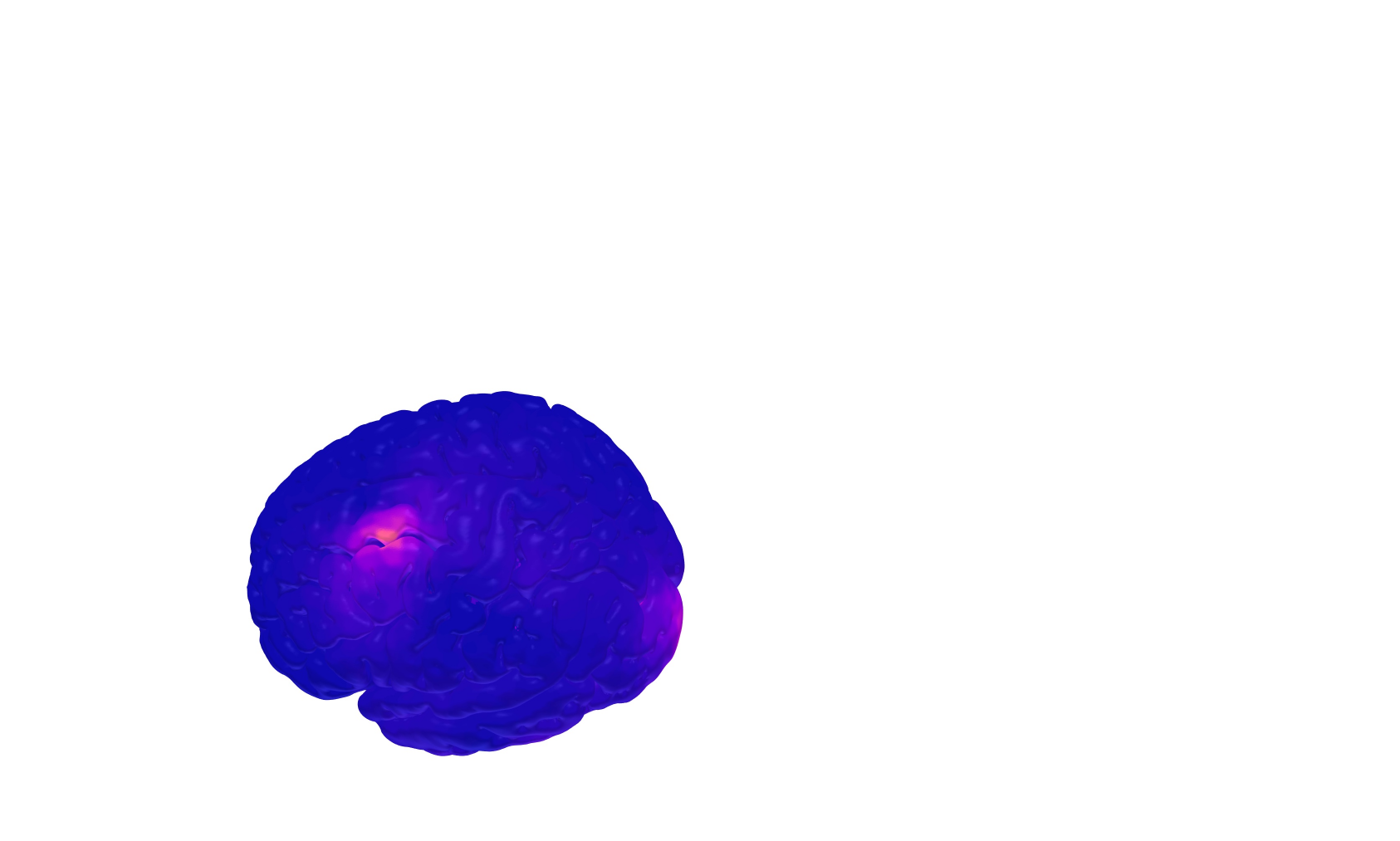}
        \end{minipage}
    \end{minipage}
    
    \caption{Spatial source reconstructions at 25 dB SNR. Rows correspond to Change-rate Standardized Kalman filter (CrSKF), Standardized Kalman filter (SKF), Kalman filter (KF), and Linearly Constrained Minimum Variance (LCMV) beamformer; columns show the reconstructions at the modeled FVEP component times associated with V1, V2d, V2v+V3d, V3v, LIP, and FEF.}
    \label{fig:brain_25dB}
\end{figure}

\begin{figure}[h!]
    \centering
    \begin{minipage}{0.02\linewidth}
        \hfill
    \end{minipage}\begin{minipage}{0.98\linewidth}
        \begin{minipage}{0.16\linewidth}
            \centering
            \scriptsize{V1}
        \end{minipage}\begin{minipage}{0.16\linewidth}
            \centering
            \scriptsize{V2d}
        \end{minipage}\begin{minipage}{0.16\linewidth}
            \centering
            \scriptsize{V2v+V3d}
        \end{minipage}\begin{minipage}{0.16\linewidth}
            \centering
            \scriptsize{V3v}
        \end{minipage}\begin{minipage}{0.16\linewidth}
        \centering
            \scriptsize{LIP}
        \end{minipage}\begin{minipage}{0.16\linewidth}
        \centering
            \scriptsize{FEF}
        \end{minipage}\end{minipage}
    
    \begin{minipage}{0.02\linewidth}
        \rotatebox{90}{\scriptsize CrSKF}
    \end{minipage}\begin{minipage}{0.98\linewidth}
        \begin{minipage}{0.16\linewidth}
            \includegraphics[trim={5cm 1cm 14cm 8cm},clip,width=\linewidth]{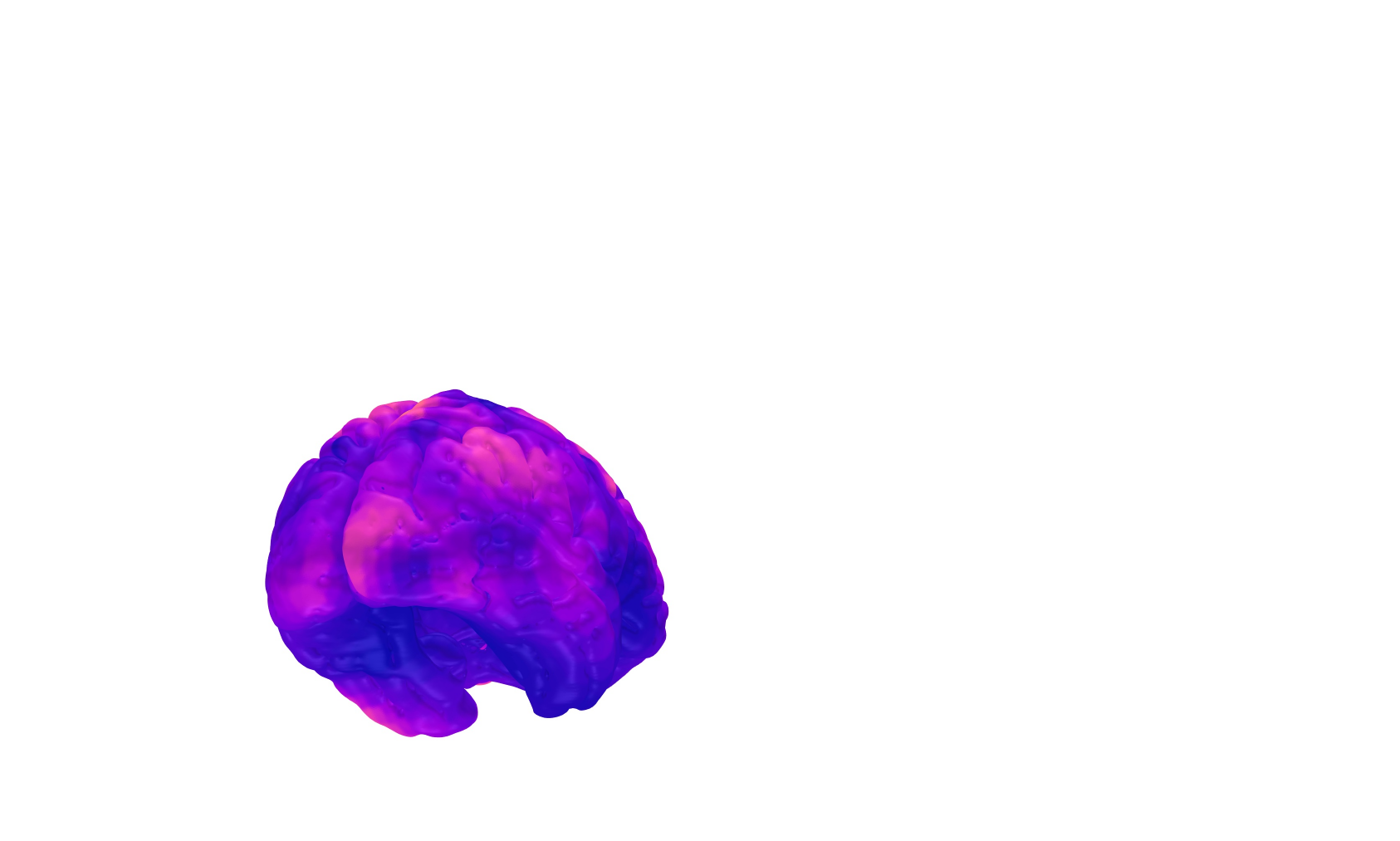}
        \end{minipage}\begin{minipage}{0.16\linewidth}
            \includegraphics[trim={5cm 1cm 14cm 8cm},clip,width=\linewidth]{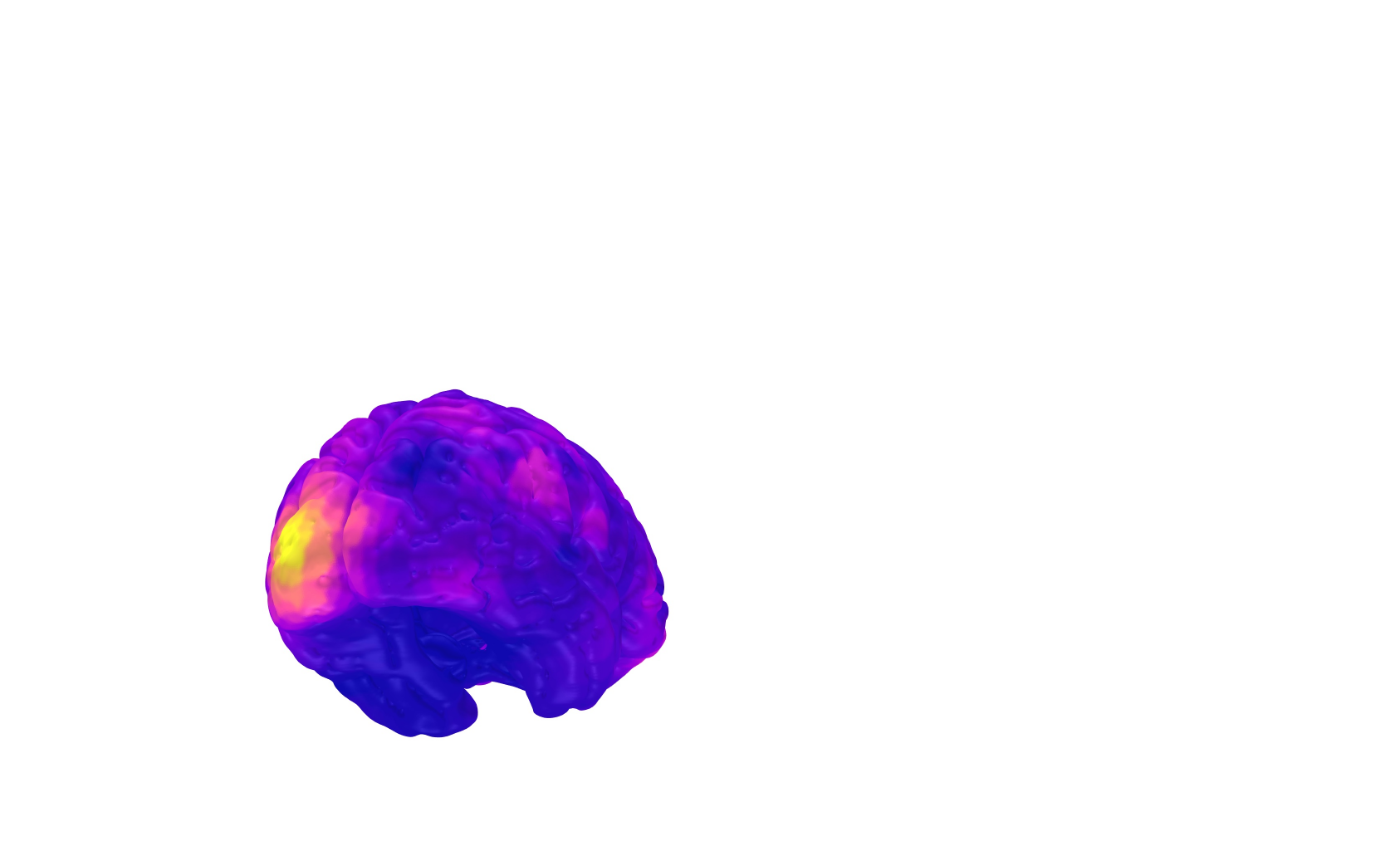}
        \end{minipage}\begin{minipage}{0.16\linewidth}
            \includegraphics[trim={5cm 1cm 14cm 8cm},clip,width=\linewidth]{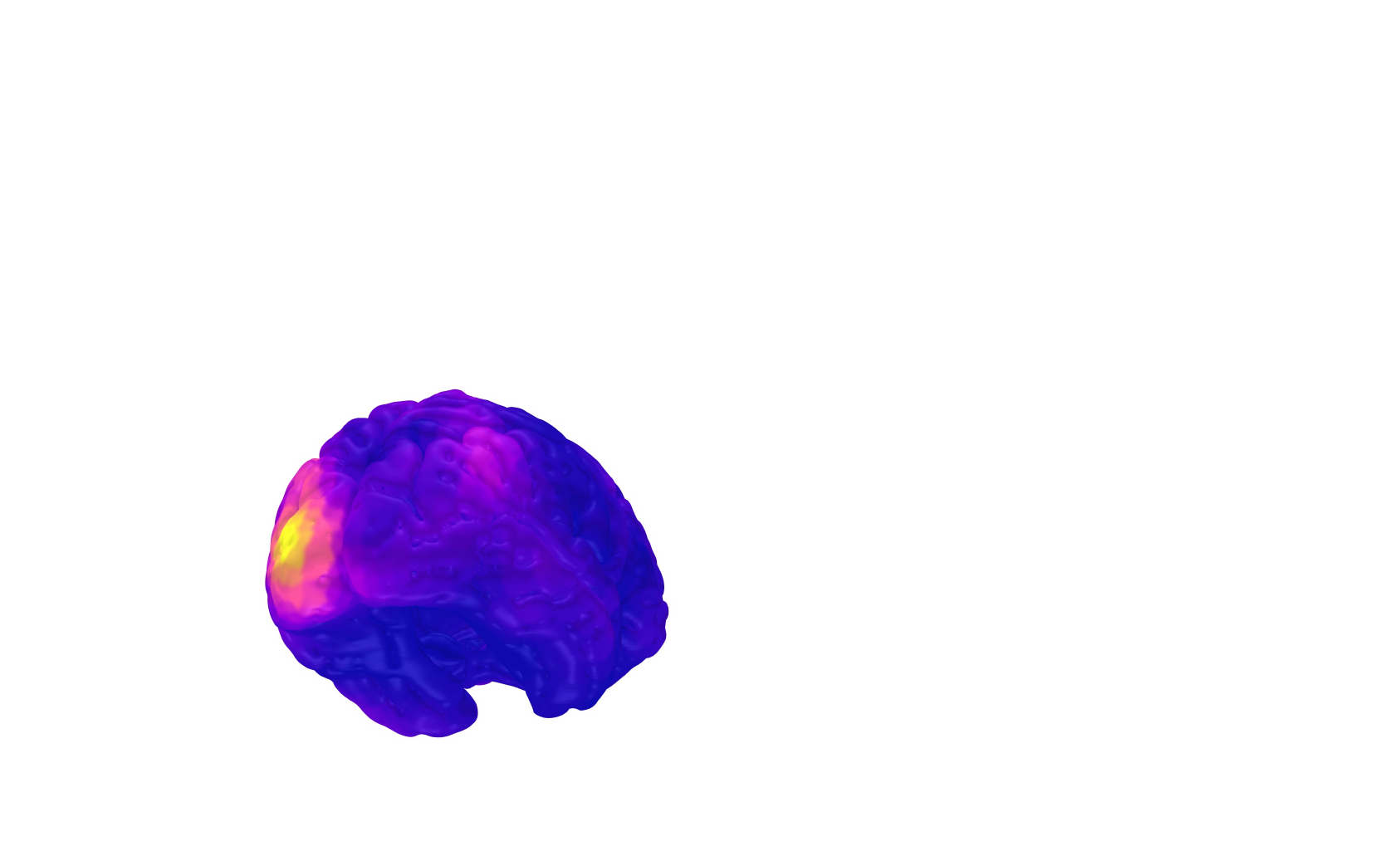}
        \end{minipage}\begin{minipage}{0.16\linewidth}
            \includegraphics[trim={5cm 1cm 14cm 8cm},clip,width=\linewidth]{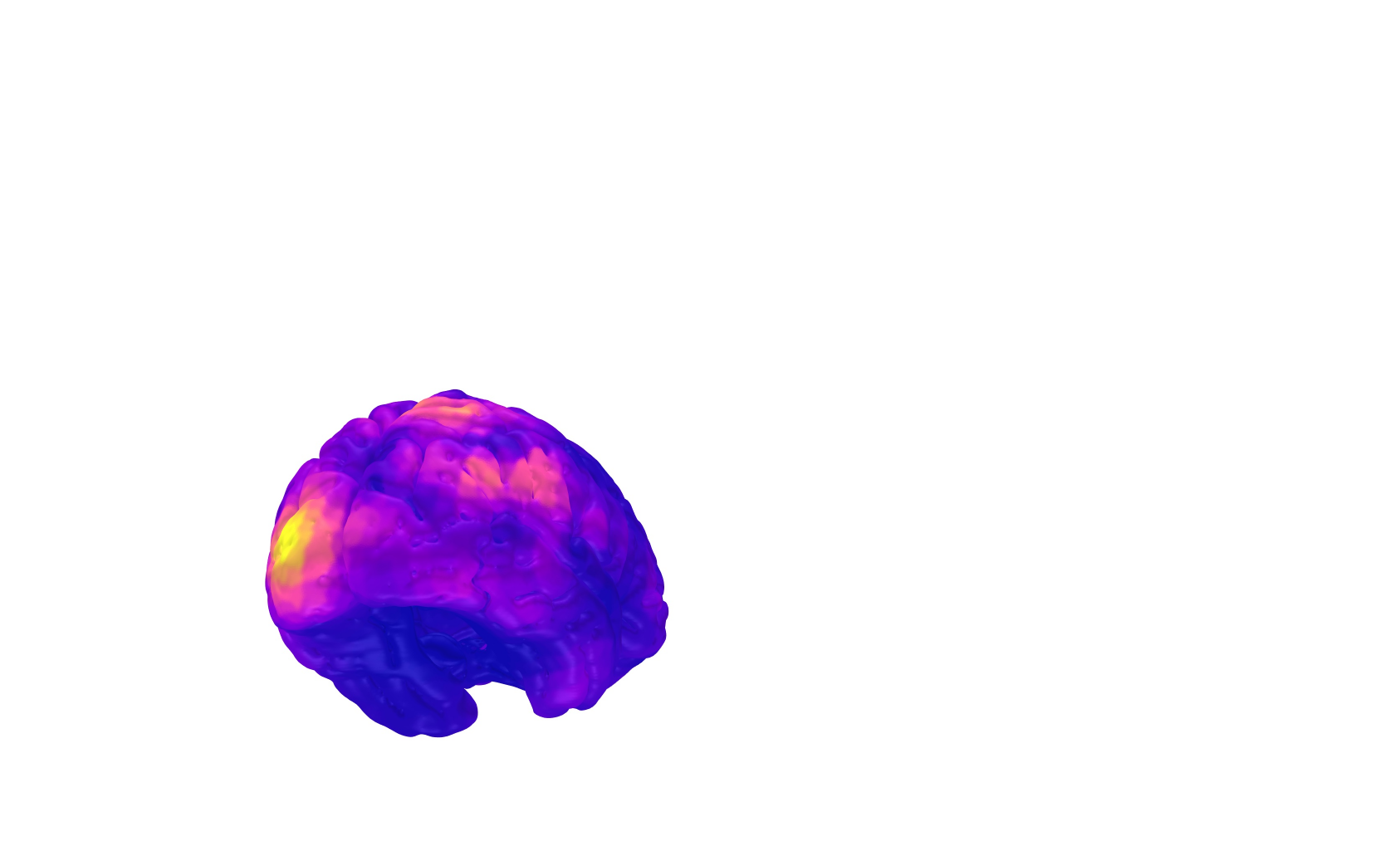}
        \end{minipage}\begin{minipage}{0.16\linewidth}
            \includegraphics[trim={5cm 1cm 14cm 8cm},clip,width=\linewidth]{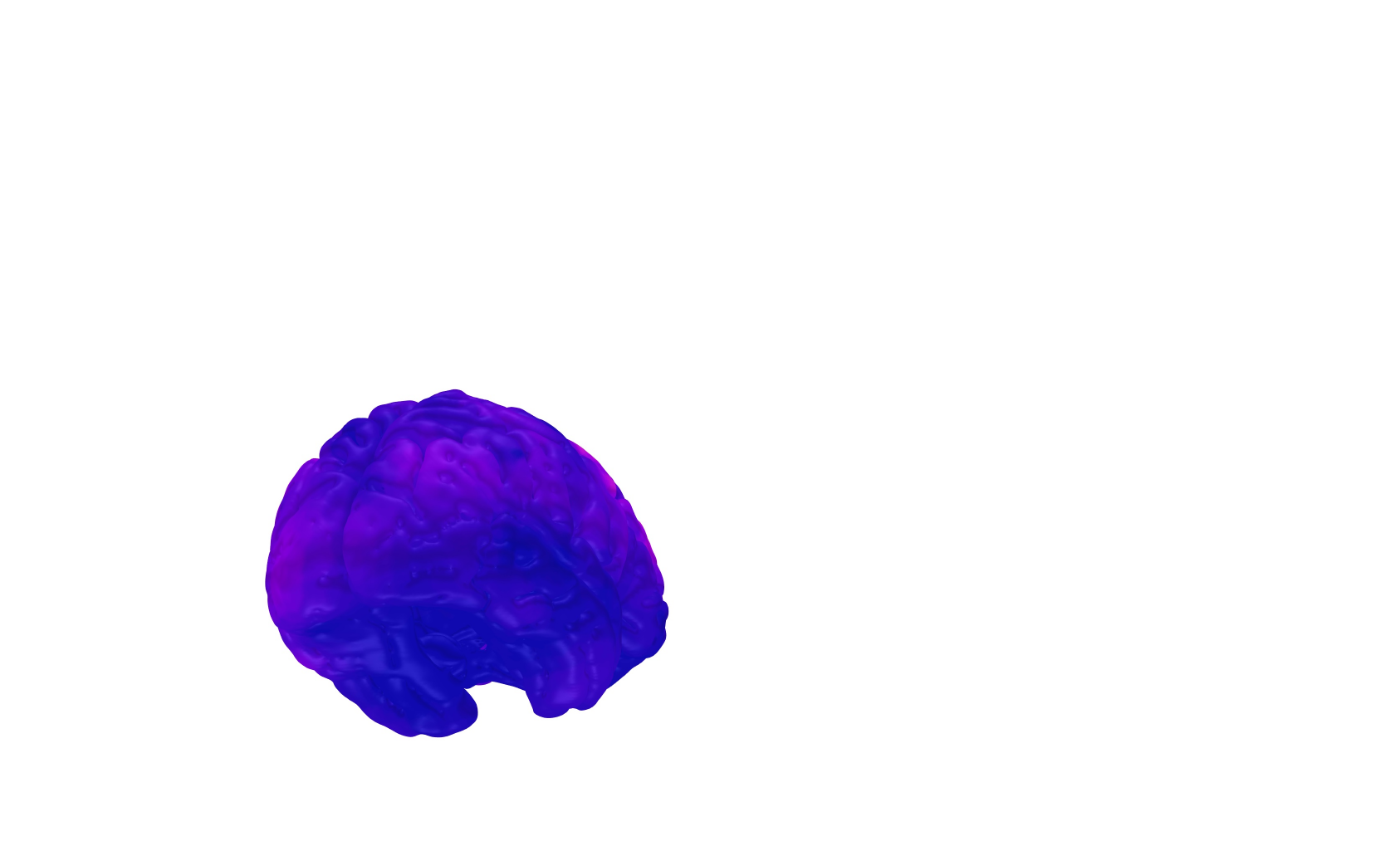}
        \end{minipage}\begin{minipage}{0.16\linewidth}
            \includegraphics[trim={5cm 1cm 14cm 8cm},clip,width=\linewidth]{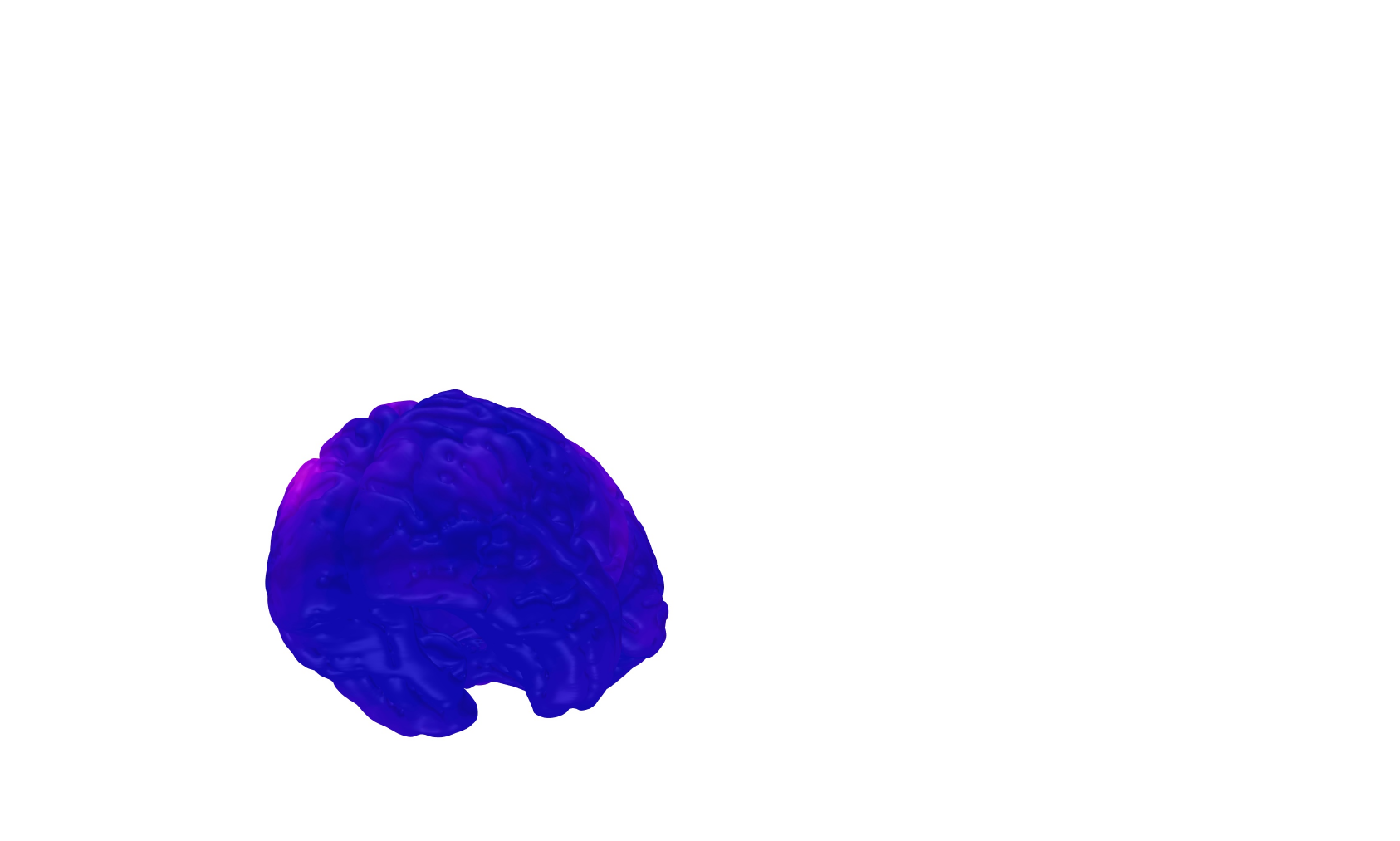}
        \end{minipage}

        \begin{minipage}{0.16\linewidth}
            \includegraphics[trim={5cm 1cm 14cm 8cm},clip,width=\linewidth]{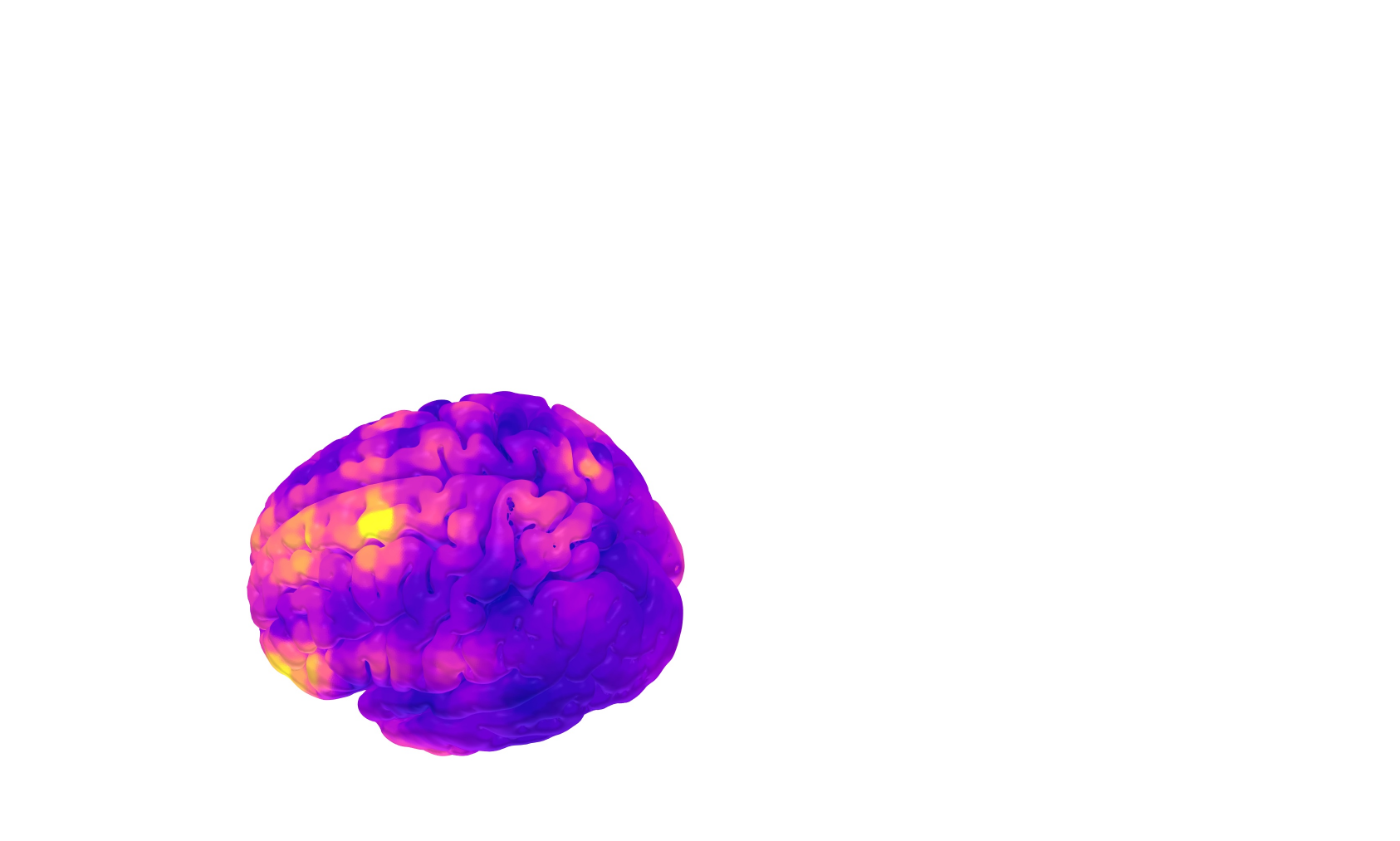}
        \end{minipage}\begin{minipage}{0.16\linewidth}
            \includegraphics[trim={5cm 1cm 14cm 8cm},clip,width=\linewidth]{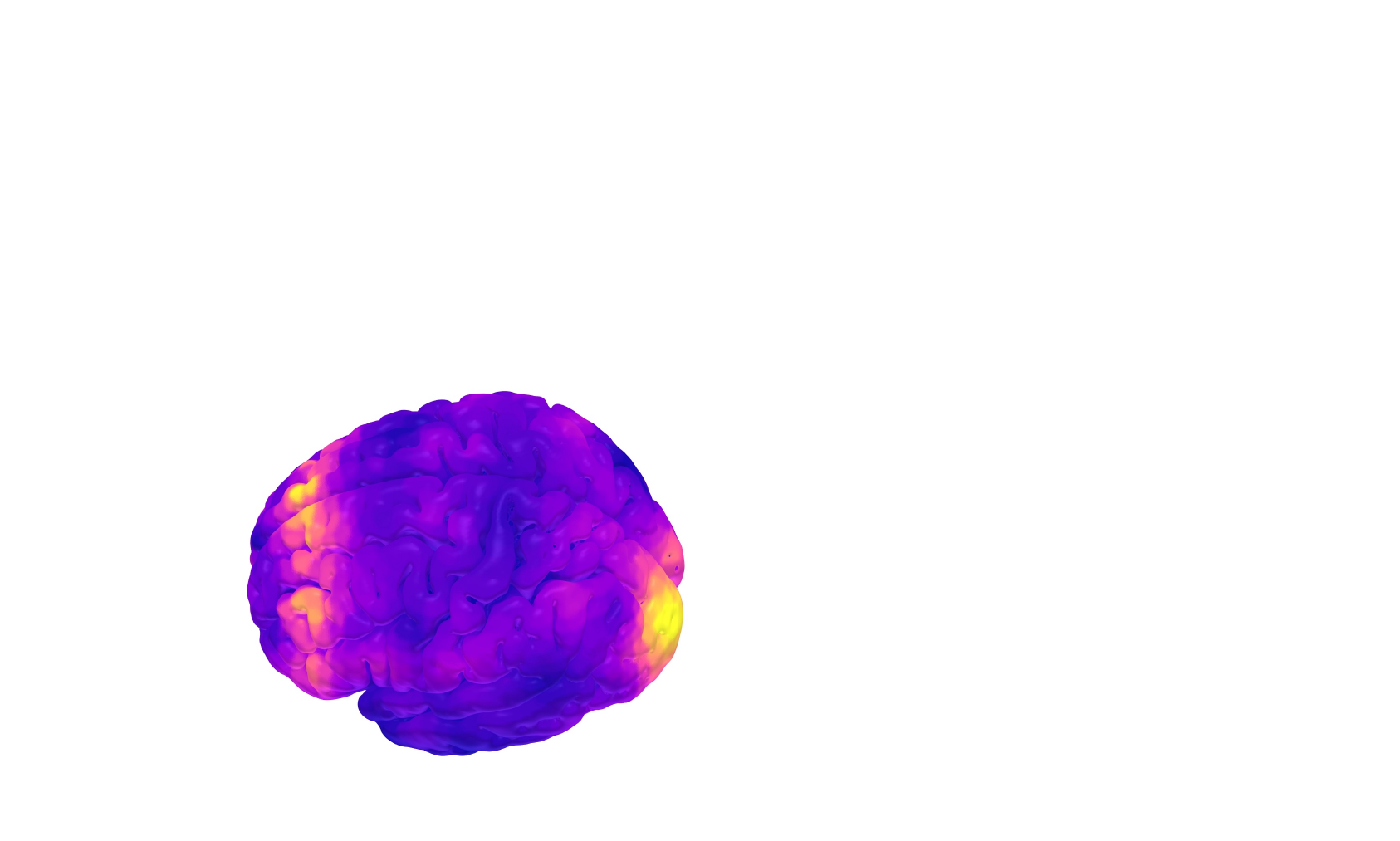}
        \end{minipage}\begin{minipage}{0.16\linewidth}
            \includegraphics[trim={5cm 1cm 14cm 8cm},clip,width=\linewidth]{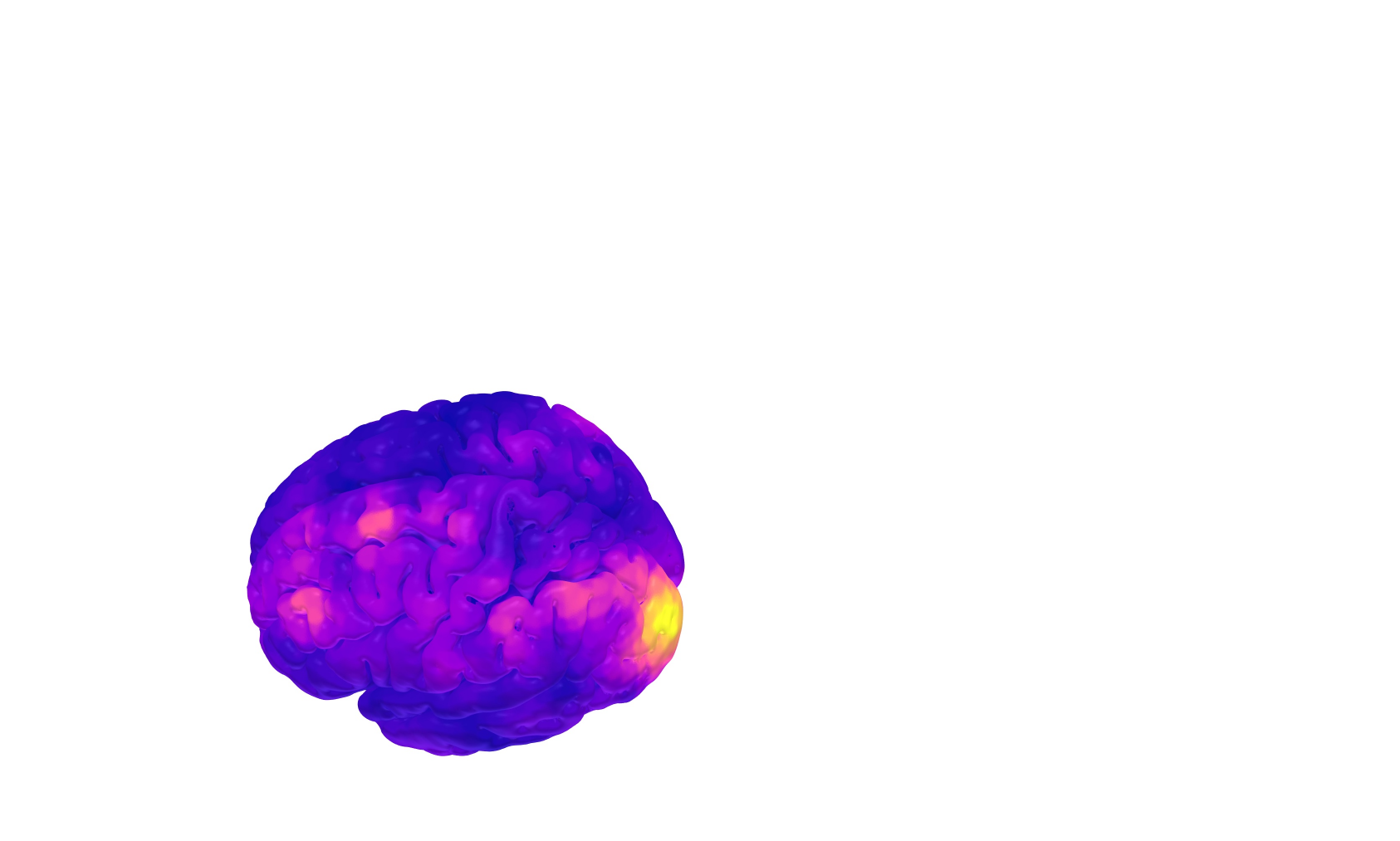}
        \end{minipage}\begin{minipage}{0.16\linewidth}
            \includegraphics[trim={5cm 1cm 14cm 8cm},clip,width=\linewidth]{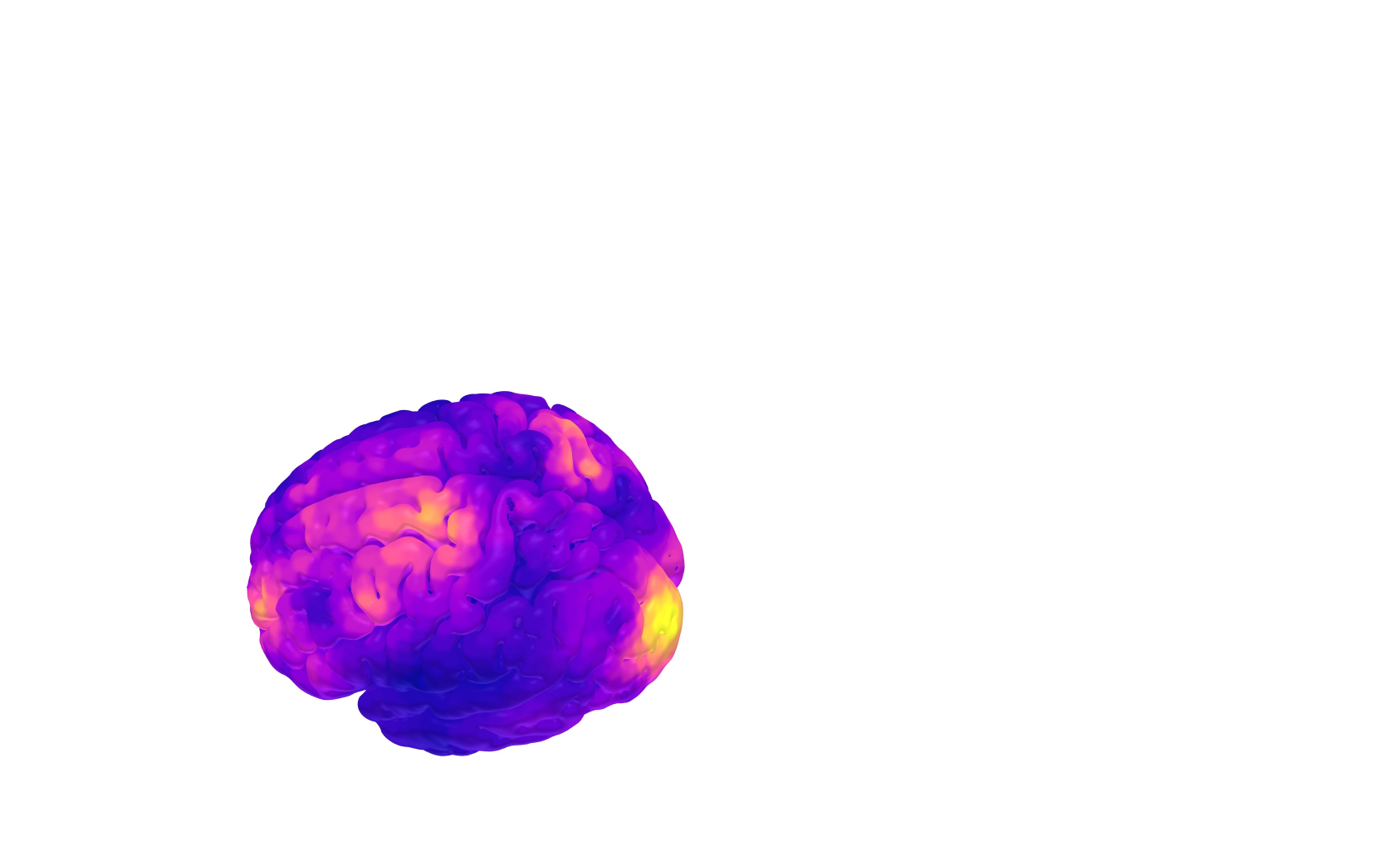}
        \end{minipage}\begin{minipage}{0.16\linewidth}
            \includegraphics[trim={5cm 1cm 14cm 8cm},clip,width=\linewidth]{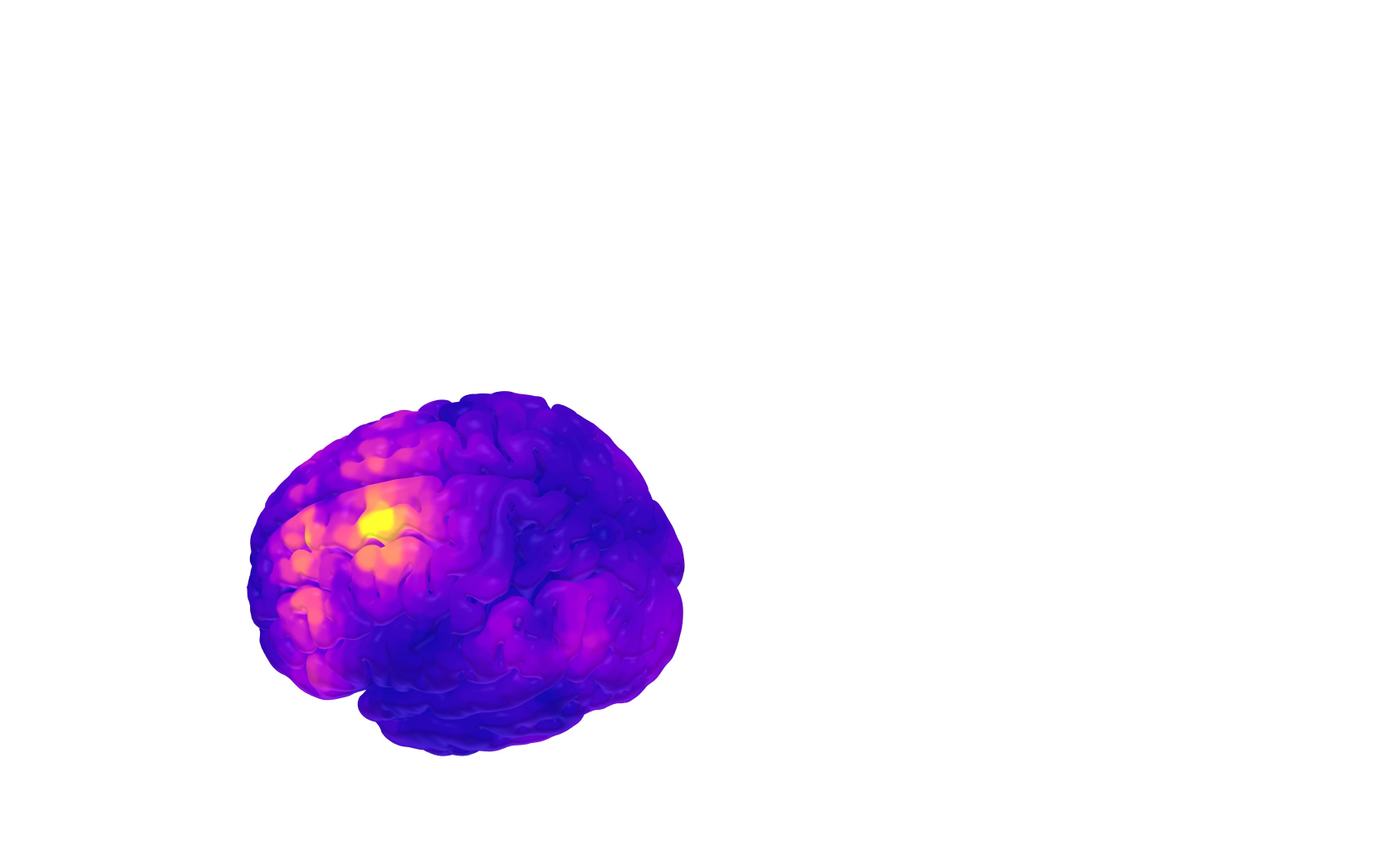}
        \end{minipage}\begin{minipage}{0.16\linewidth}
            \includegraphics[trim={5cm 1cm 14cm 8cm},clip,width=\linewidth]{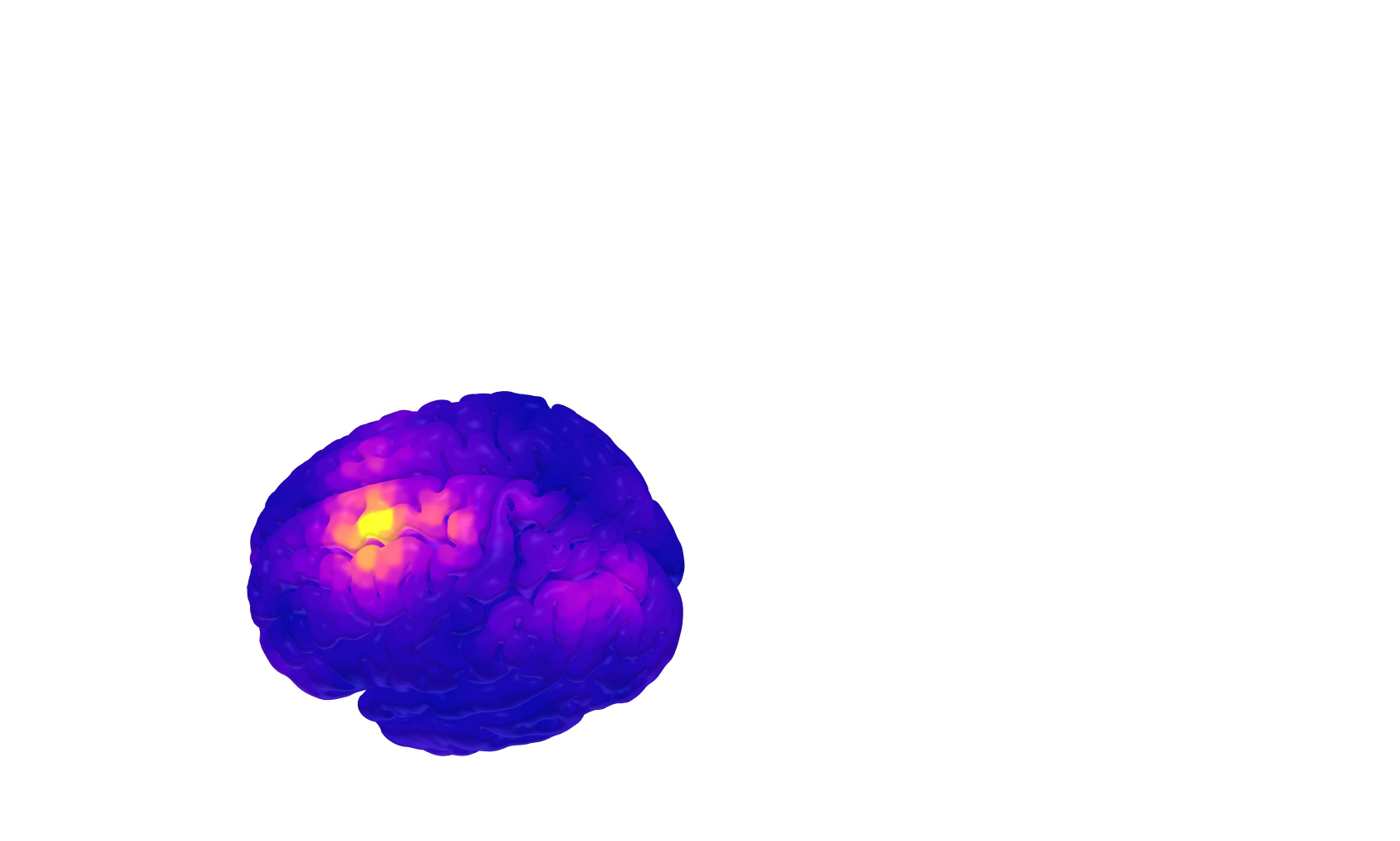}
        \end{minipage}
    \end{minipage}
    \noindent\rule[7pt]{\linewidth}{0.4pt}

    \begin{minipage}{0.02\linewidth}
        \rotatebox{90}{\scriptsize SKF}
    \end{minipage}\begin{minipage}{0.98\linewidth}
        \begin{minipage}{0.16\linewidth}
            \includegraphics[trim={5cm 1cm 14cm 8cm},clip,width=\linewidth]{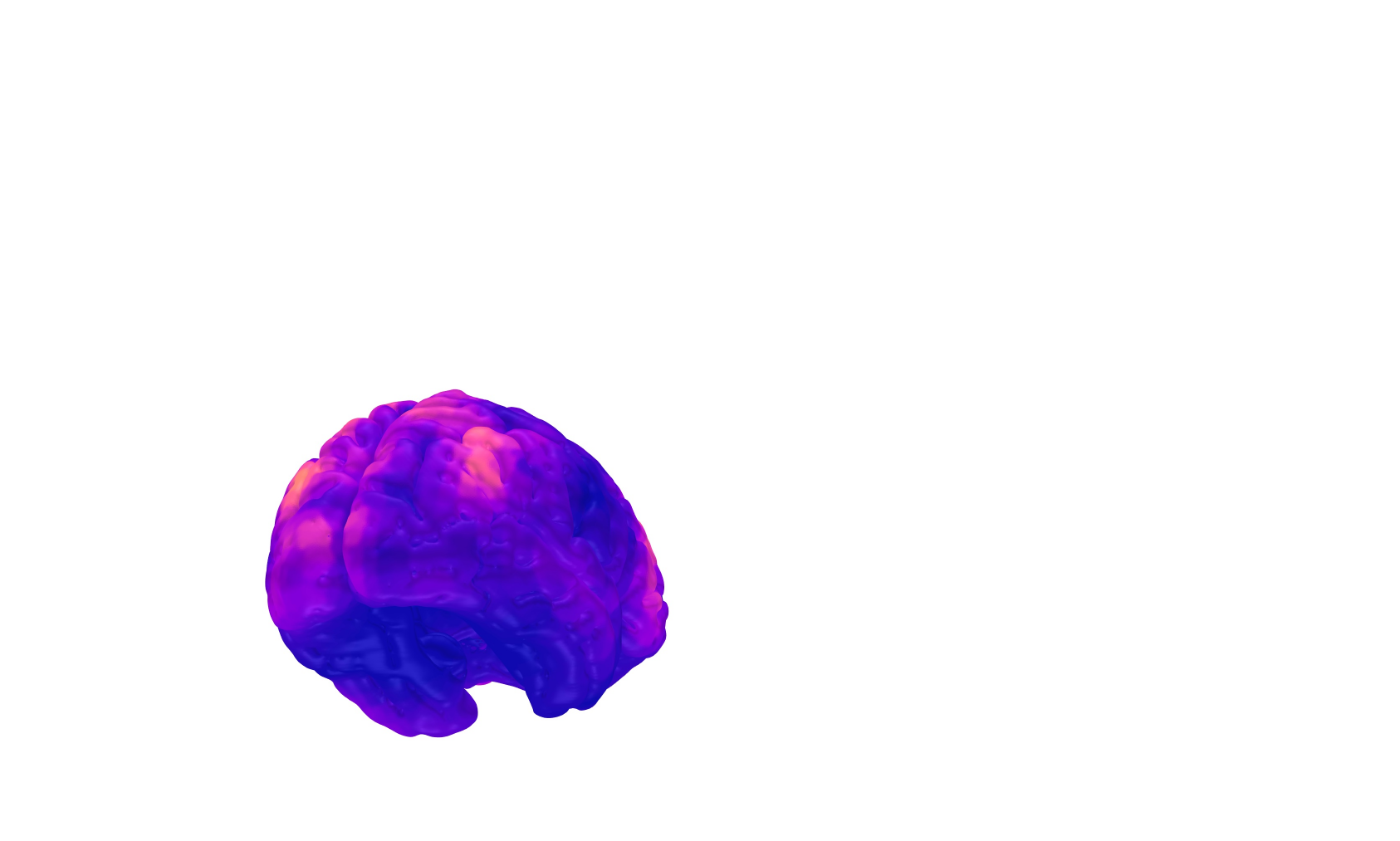}
        \end{minipage}\begin{minipage}{0.16\linewidth}
            \includegraphics[trim={5cm 1cm 14cm 8cm},clip,width=\linewidth]{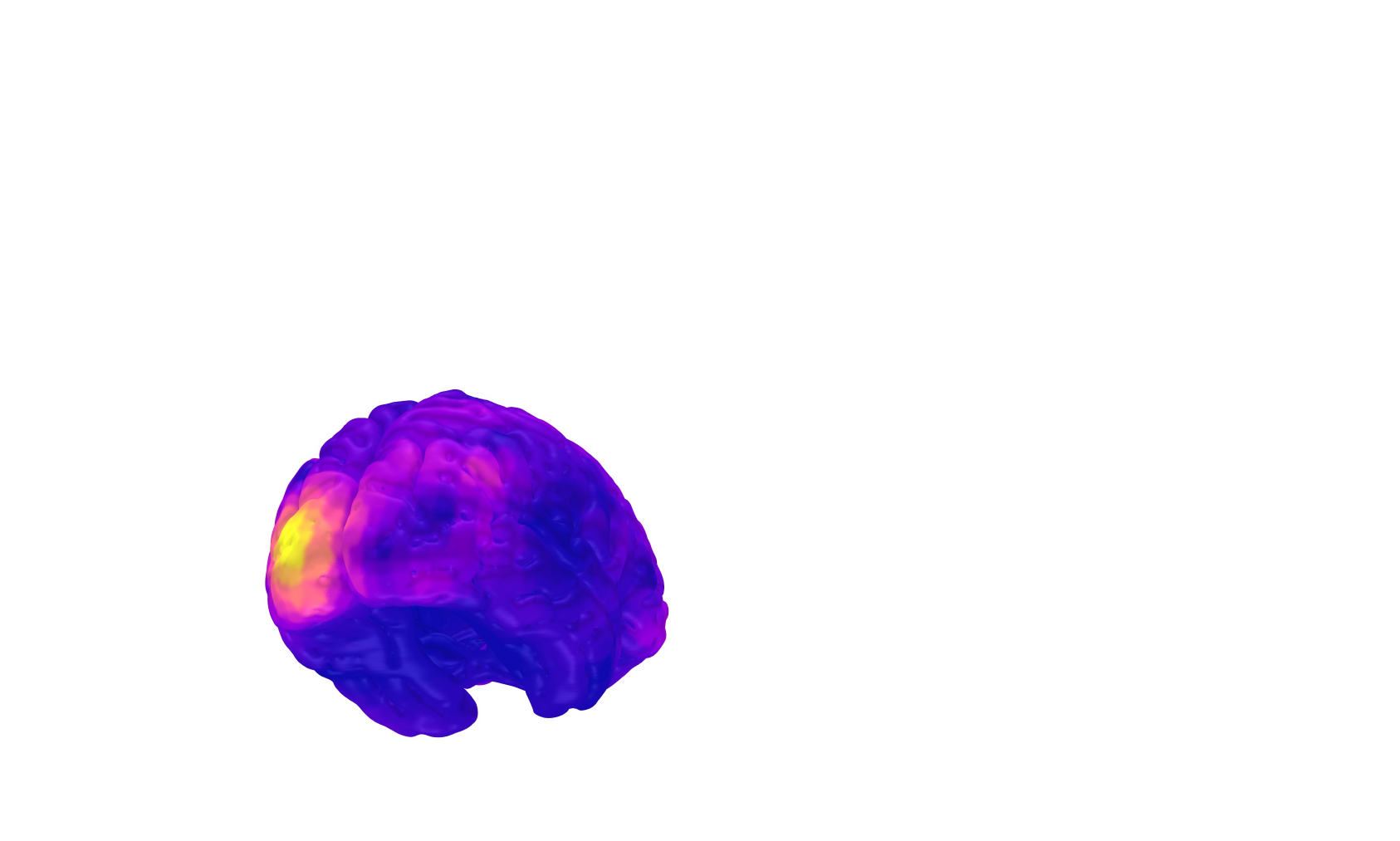}
        \end{minipage}\begin{minipage}{0.16\linewidth}
            \includegraphics[trim={5cm 1cm 14cm 8cm},clip,width=\linewidth]{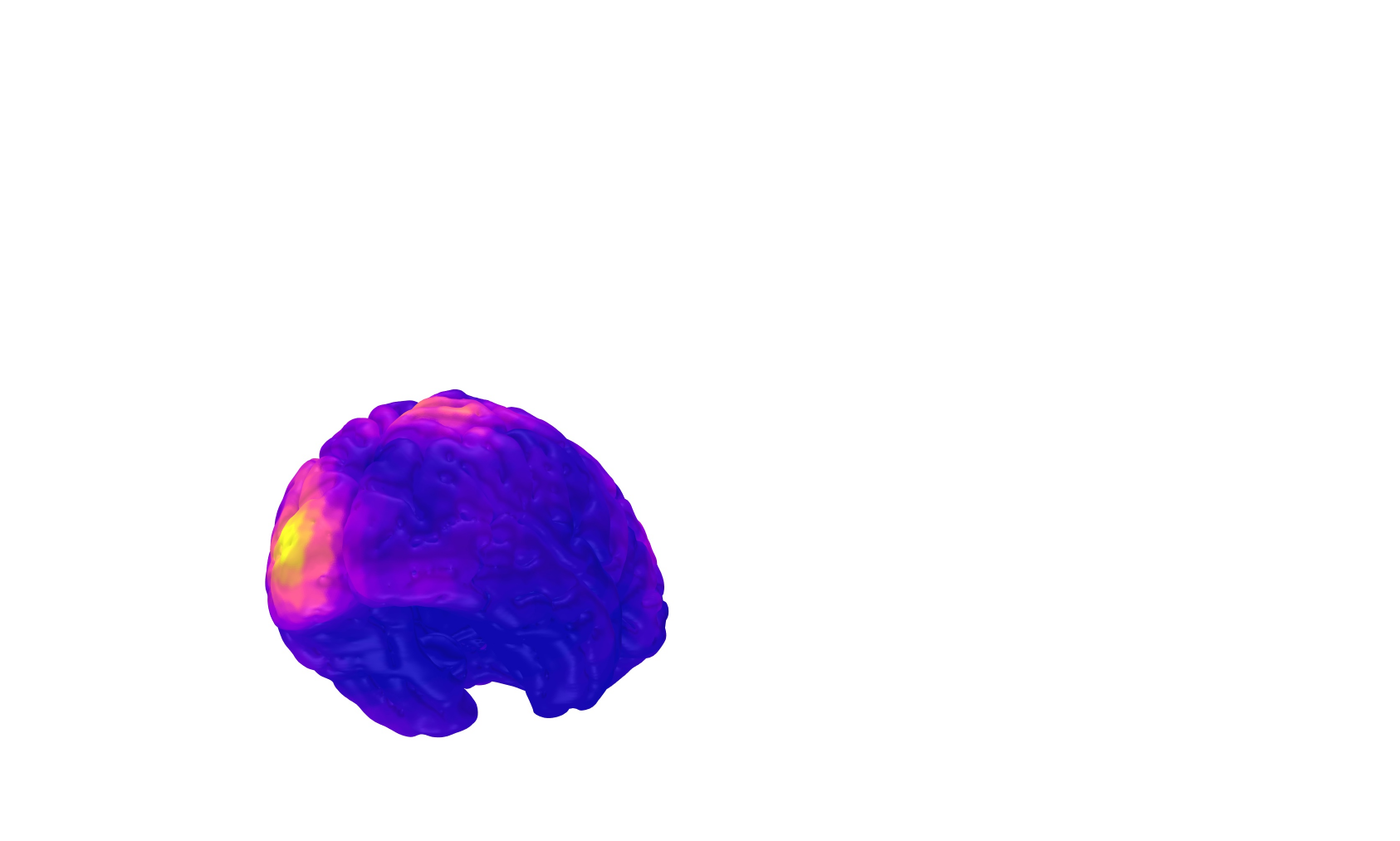}
        \end{minipage}\begin{minipage}{0.16\linewidth}
            \includegraphics[trim={5cm 1cm 14cm 8cm},clip,width=\linewidth]{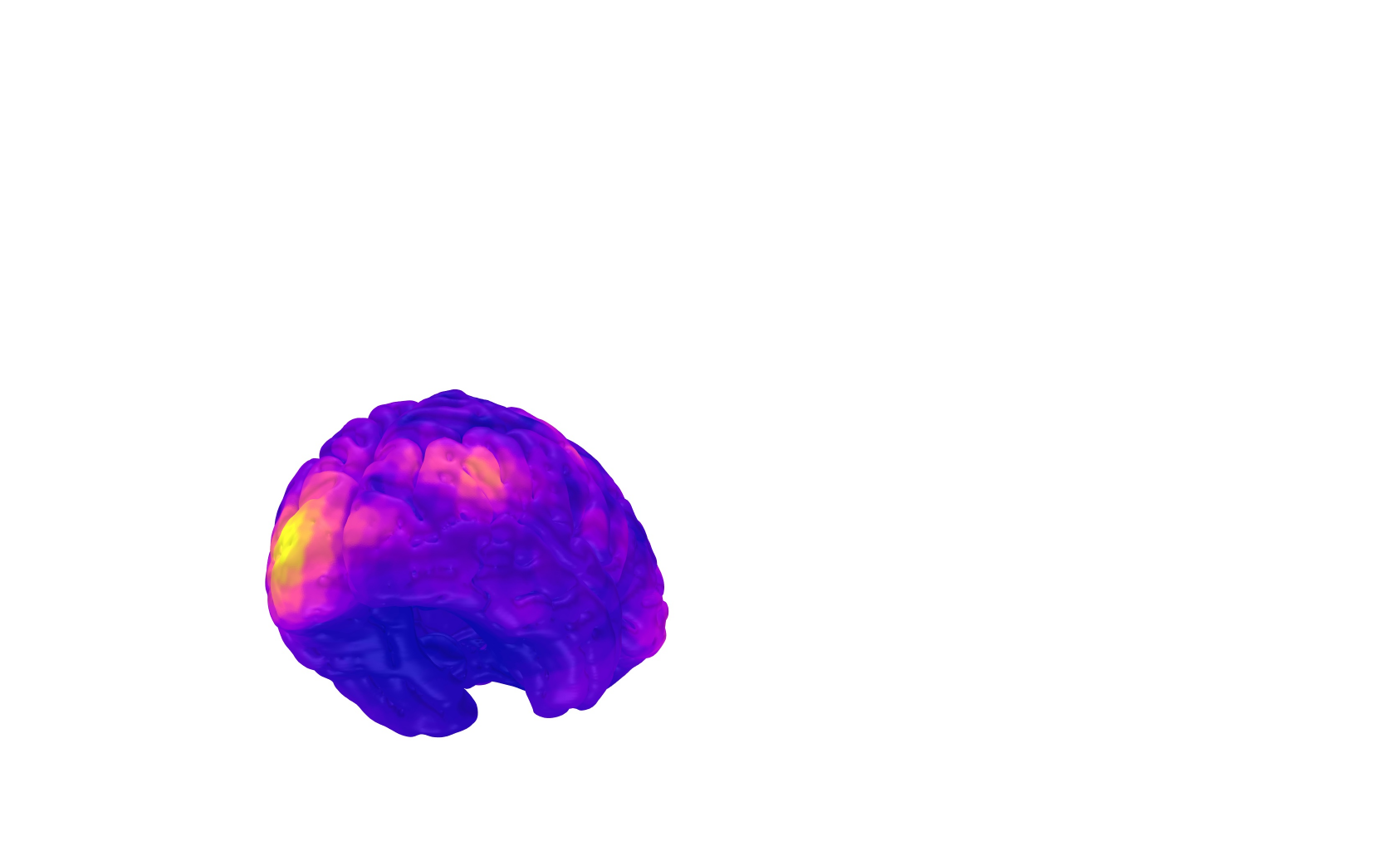}
        \end{minipage}\begin{minipage}{0.16\linewidth}
            \includegraphics[trim={5cm 1cm 14cm 8cm},clip,width=\linewidth]{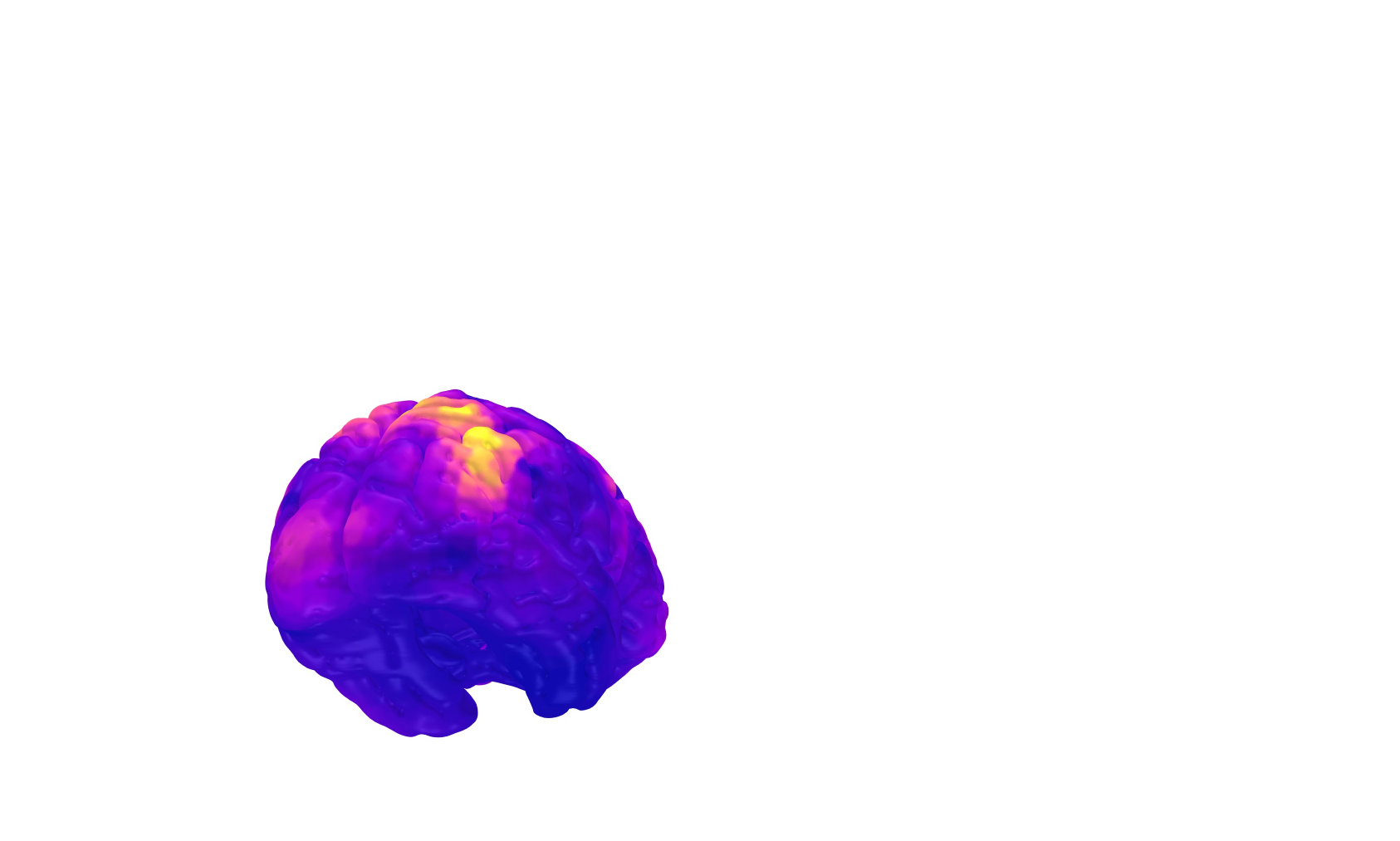}
        \end{minipage}\begin{minipage}{0.16\linewidth}
            \includegraphics[trim={5cm 1cm 14cm 8cm},clip,width=\linewidth]{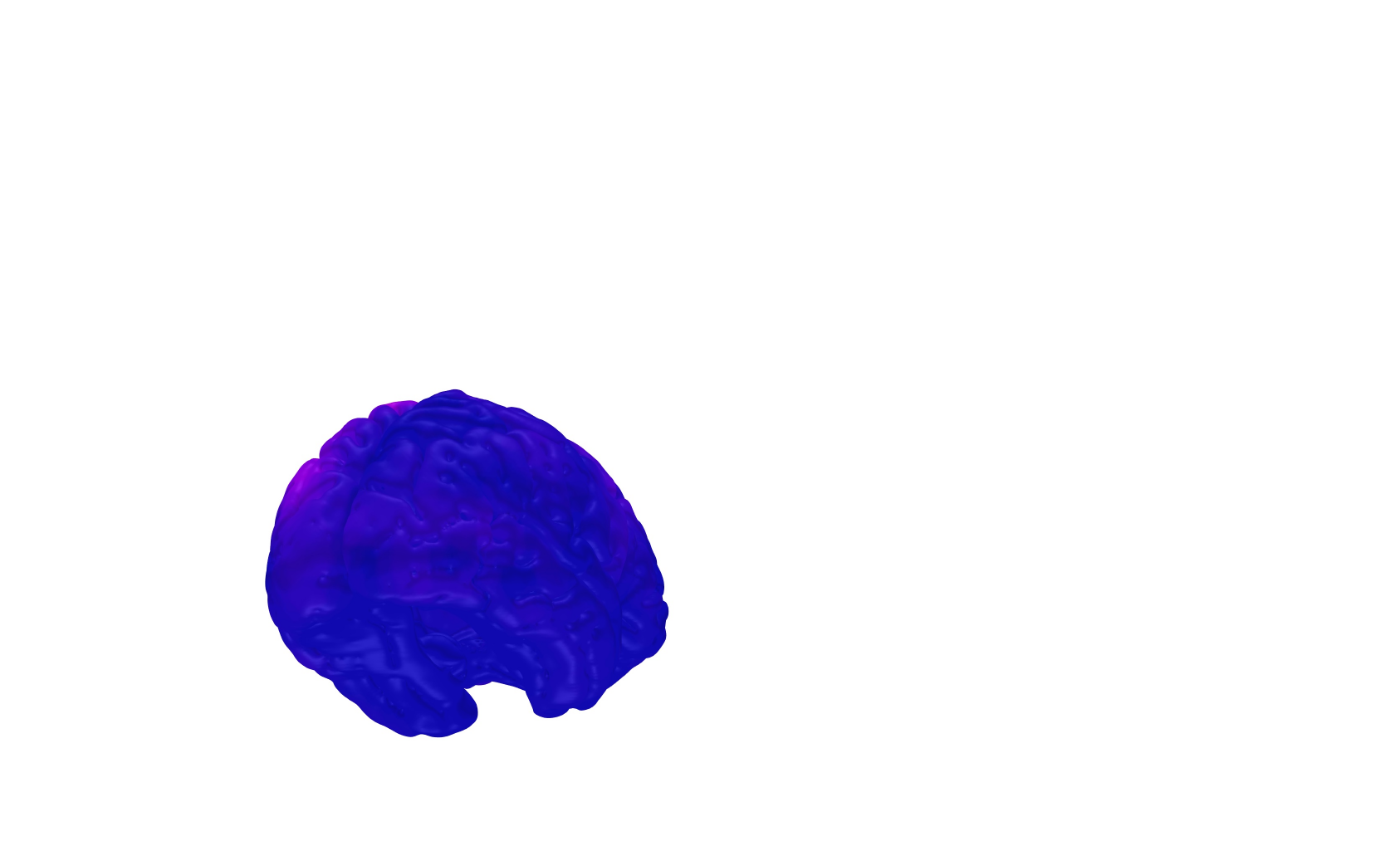}
        \end{minipage}

        \begin{minipage}{0.16\linewidth}
            \includegraphics[trim={5cm 1cm 14cm 8cm},clip,width=\linewidth]{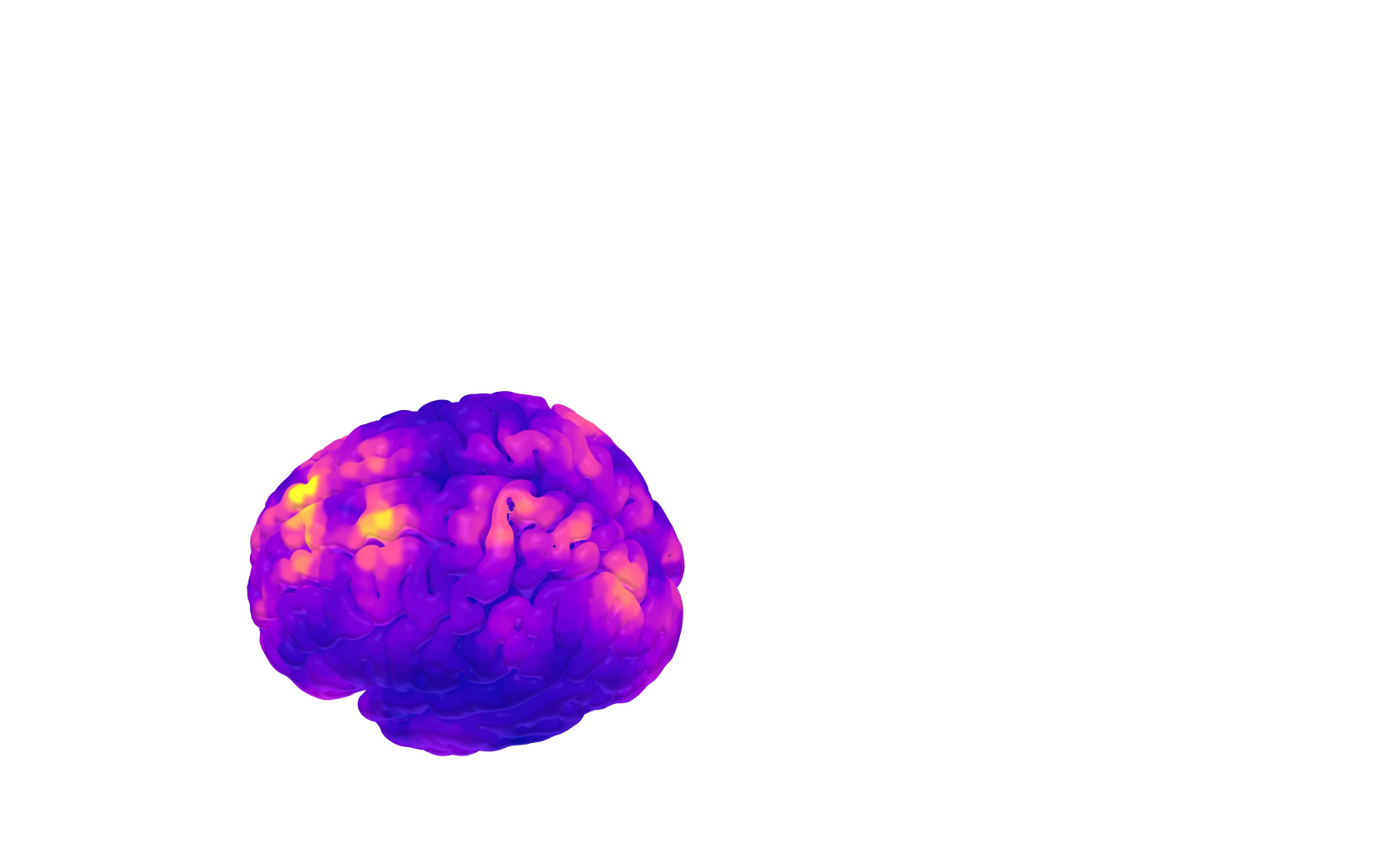}
        \end{minipage}\begin{minipage}{0.16\linewidth}
            \includegraphics[trim={5cm 1cm 14cm 8cm},clip,width=\linewidth]{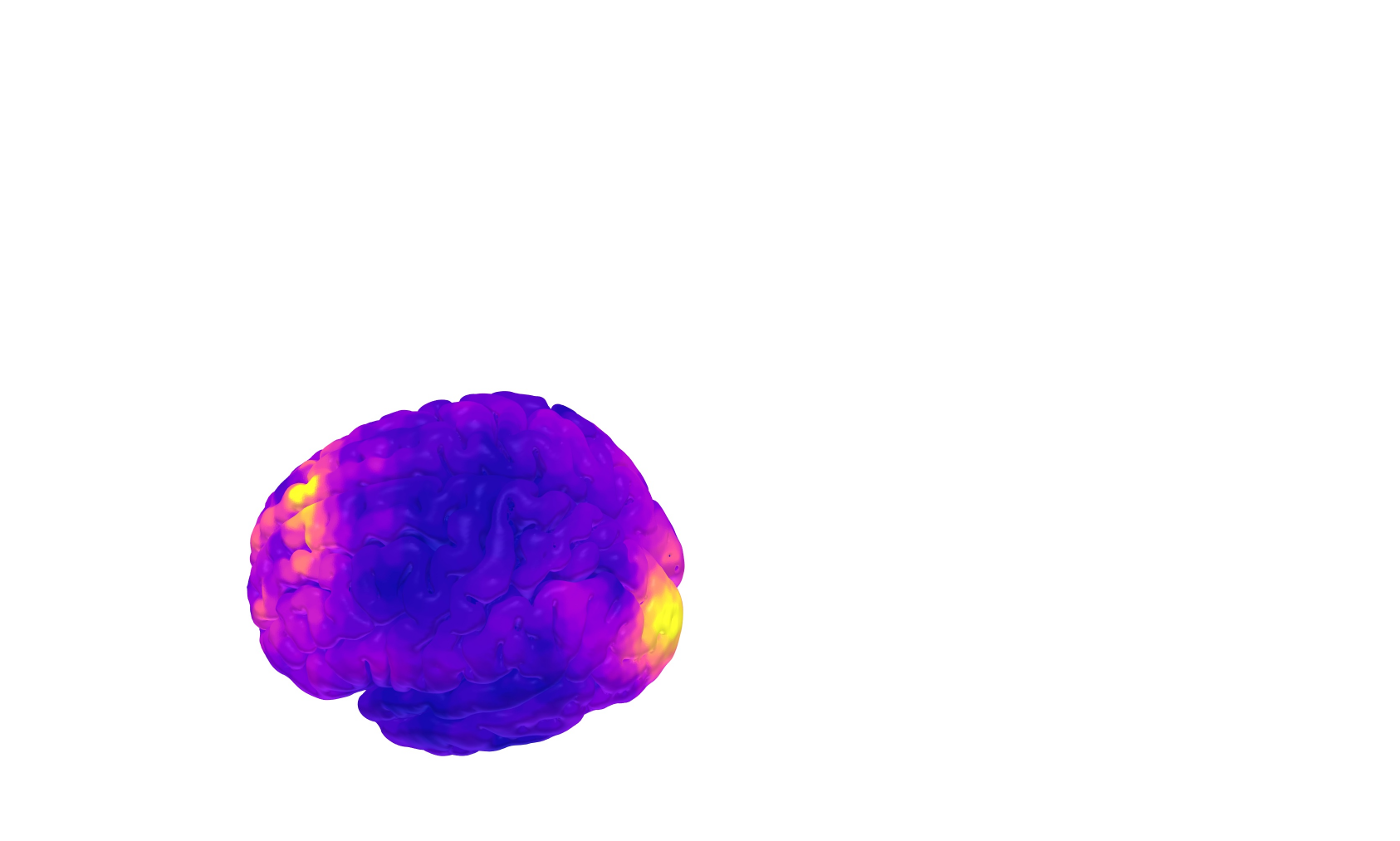}
        \end{minipage}\begin{minipage}{0.16\linewidth}
            \includegraphics[trim={5cm 1cm 14cm 8cm},clip,width=\linewidth]{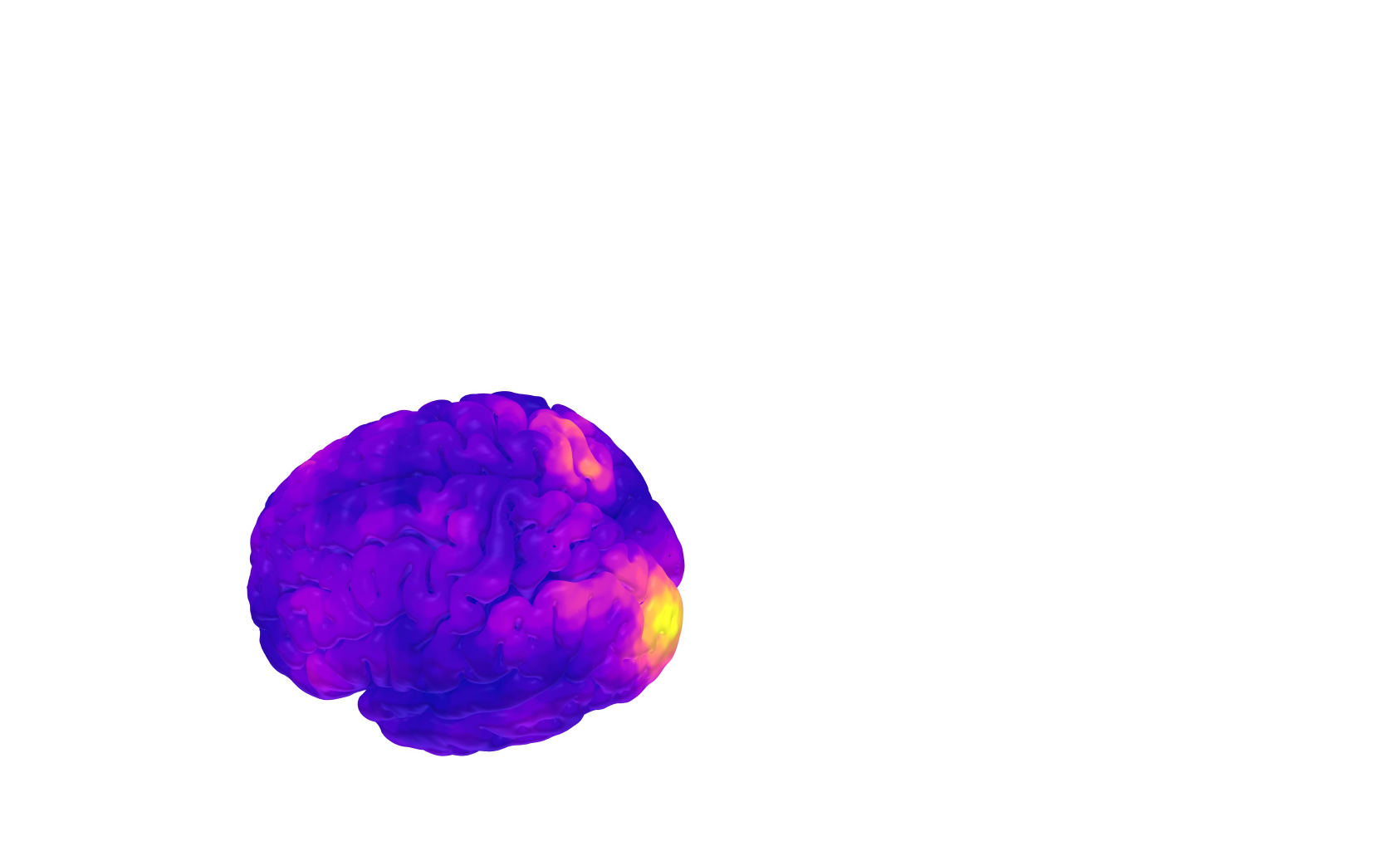}
        \end{minipage}\begin{minipage}{0.16\linewidth}
            \includegraphics[trim={5cm 1cm 14cm 8cm},clip,width=\linewidth]{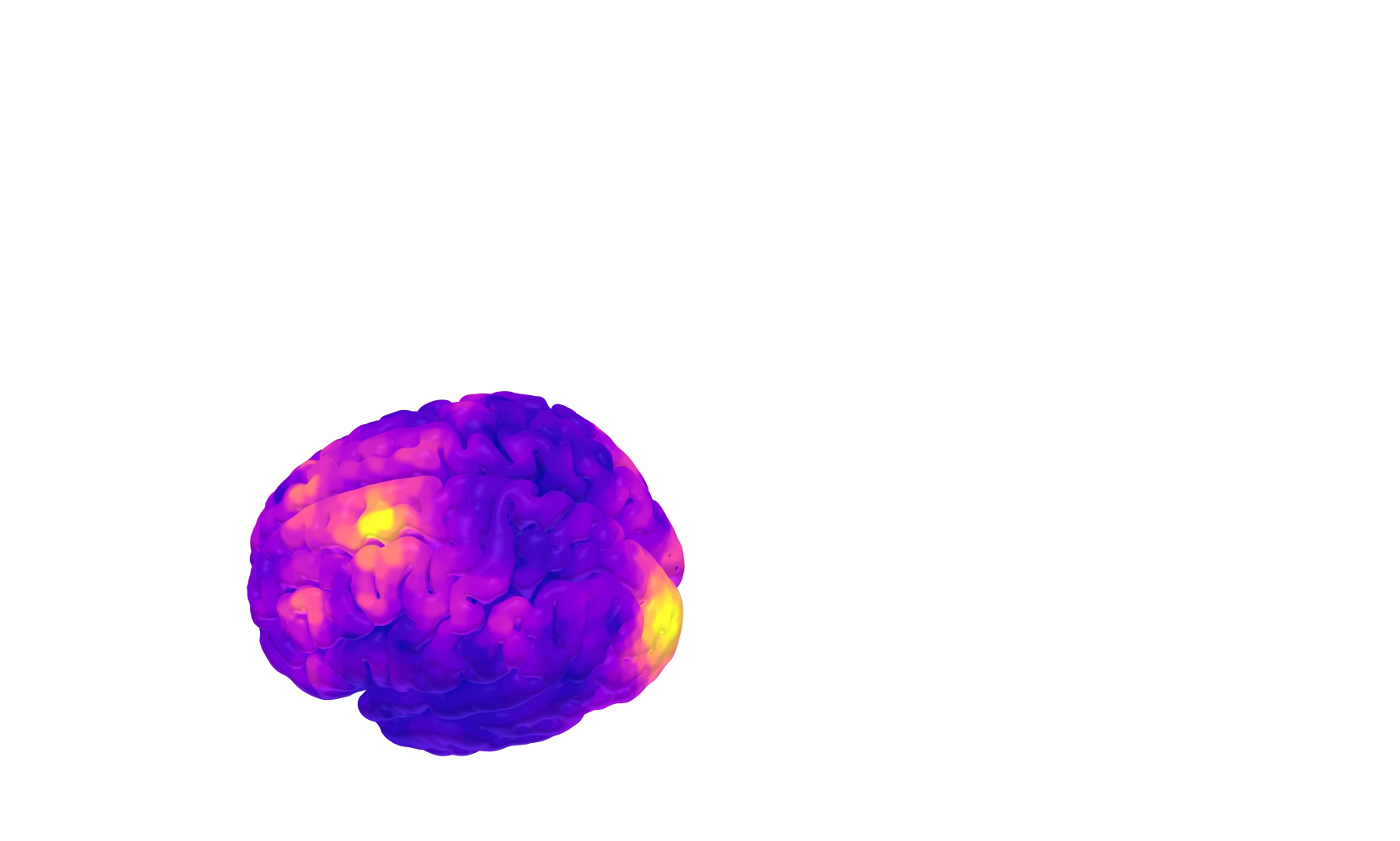}
        \end{minipage}\begin{minipage}{0.16\linewidth}
            \includegraphics[trim={5cm 1cm 14cm 8cm},clip,width=\linewidth]{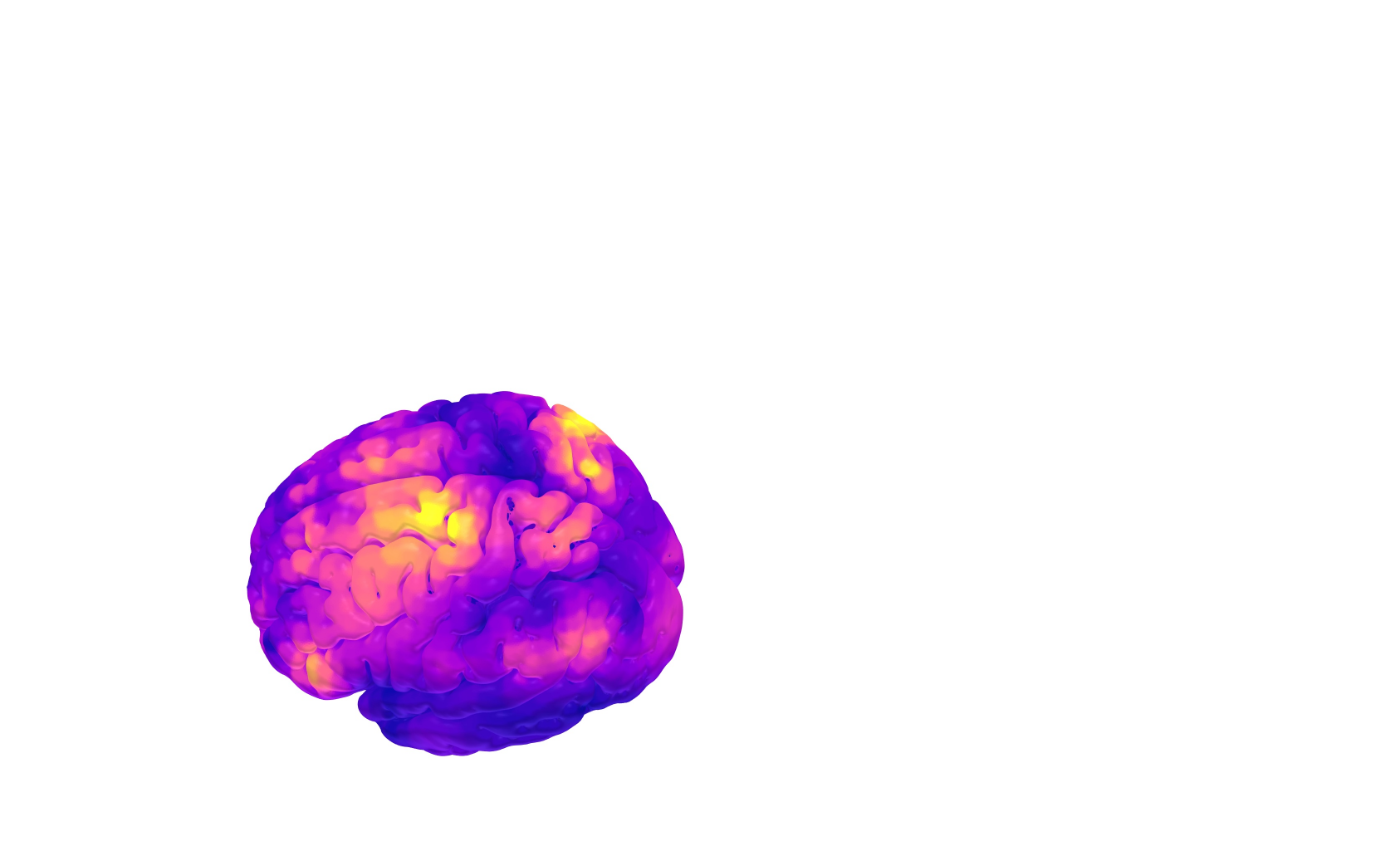}
        \end{minipage}\begin{minipage}{0.16\linewidth}
            \includegraphics[trim={5cm 1cm 14cm 8cm},clip,width=\linewidth]{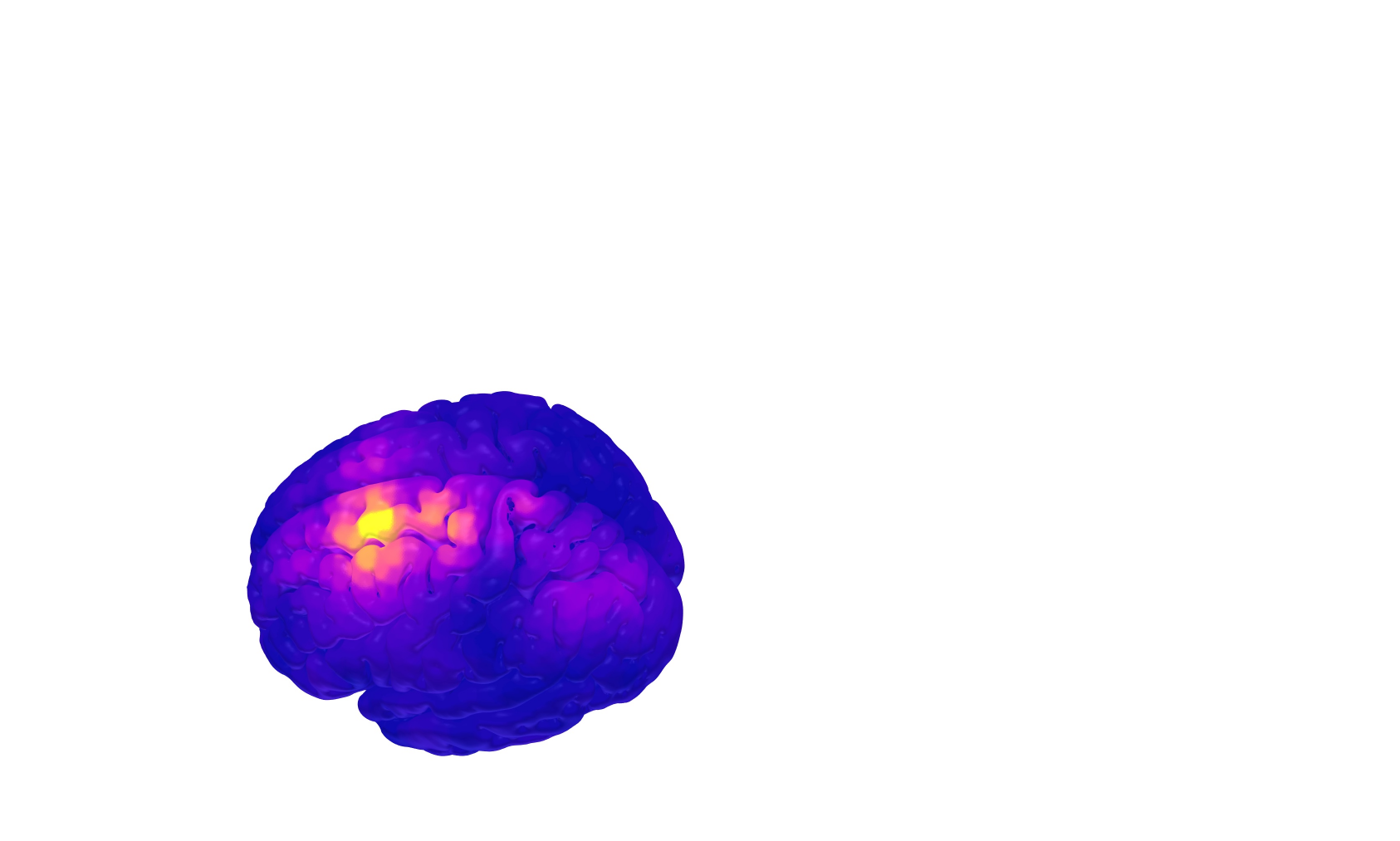}
        \end{minipage}
    \end{minipage}
    \noindent\rule[7pt]{\linewidth}{0.4pt}

    \begin{minipage}{0.02\linewidth}
        \rotatebox{90}{\scriptsize KF}
    \end{minipage}\begin{minipage}{0.98\linewidth}
        \begin{minipage}{0.16\linewidth}
            \includegraphics[trim={5cm 1cm 14cm 8cm},clip,width=\linewidth]{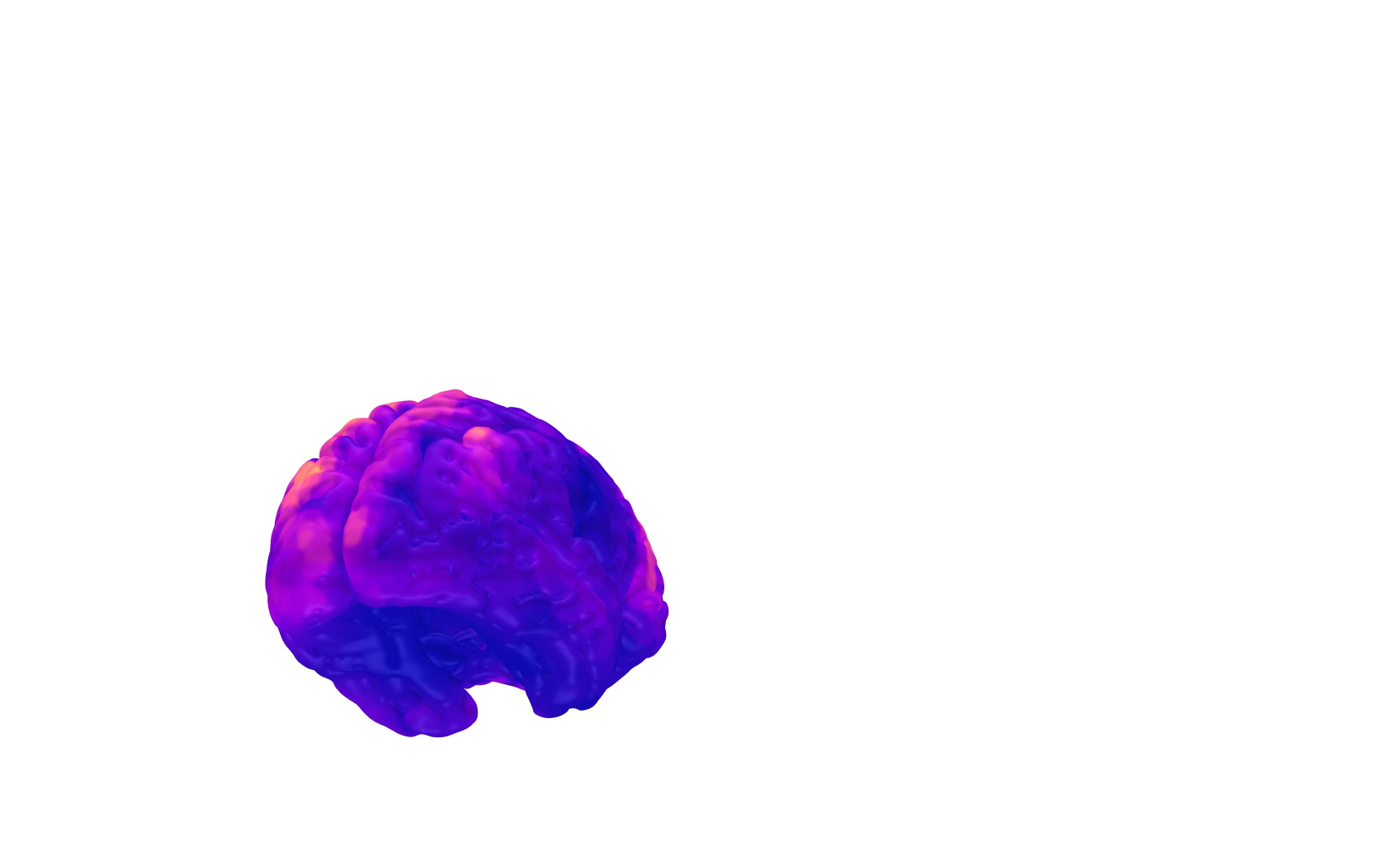}
        \end{minipage}\begin{minipage}{0.16\linewidth}
            \includegraphics[trim={5cm 1cm 14cm 8cm},clip,width=\linewidth]{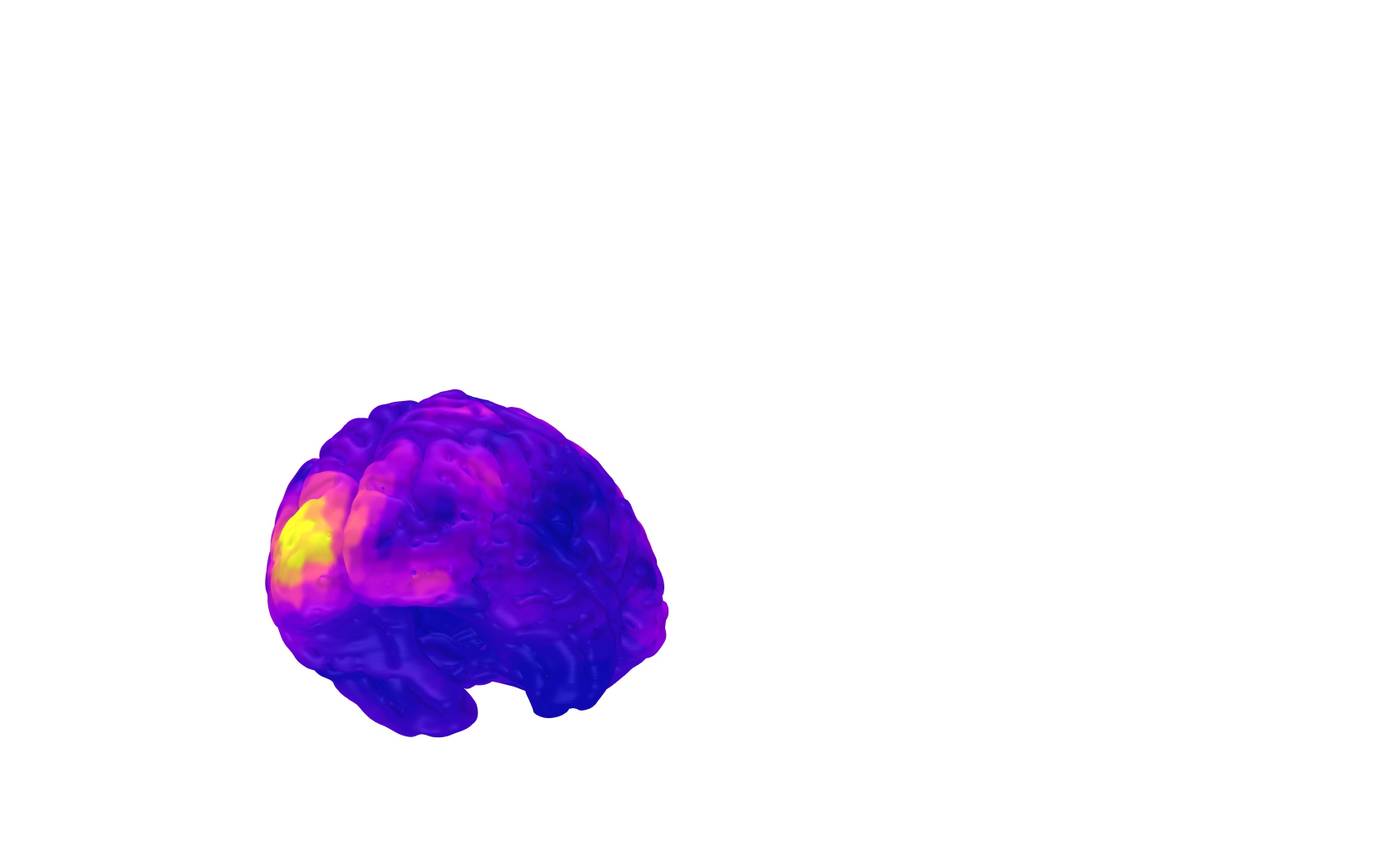}
        \end{minipage}\begin{minipage}{0.16\linewidth}
            \includegraphics[trim={5cm 1cm 14cm 8cm},clip,width=\linewidth]{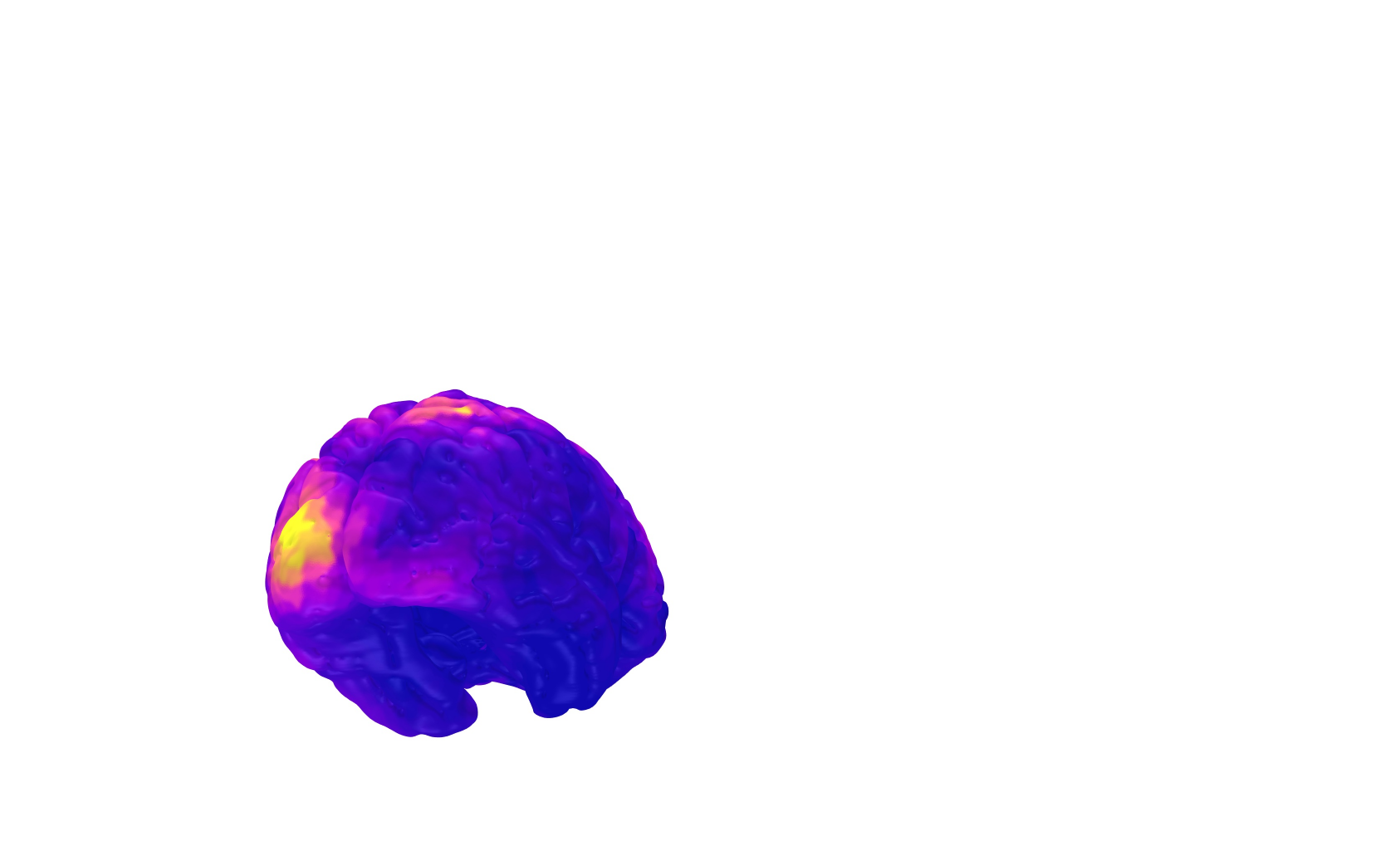}
        \end{minipage}\begin{minipage}{0.16\linewidth}
            \includegraphics[trim={5cm 1cm 14cm 8cm},clip,width=\linewidth]{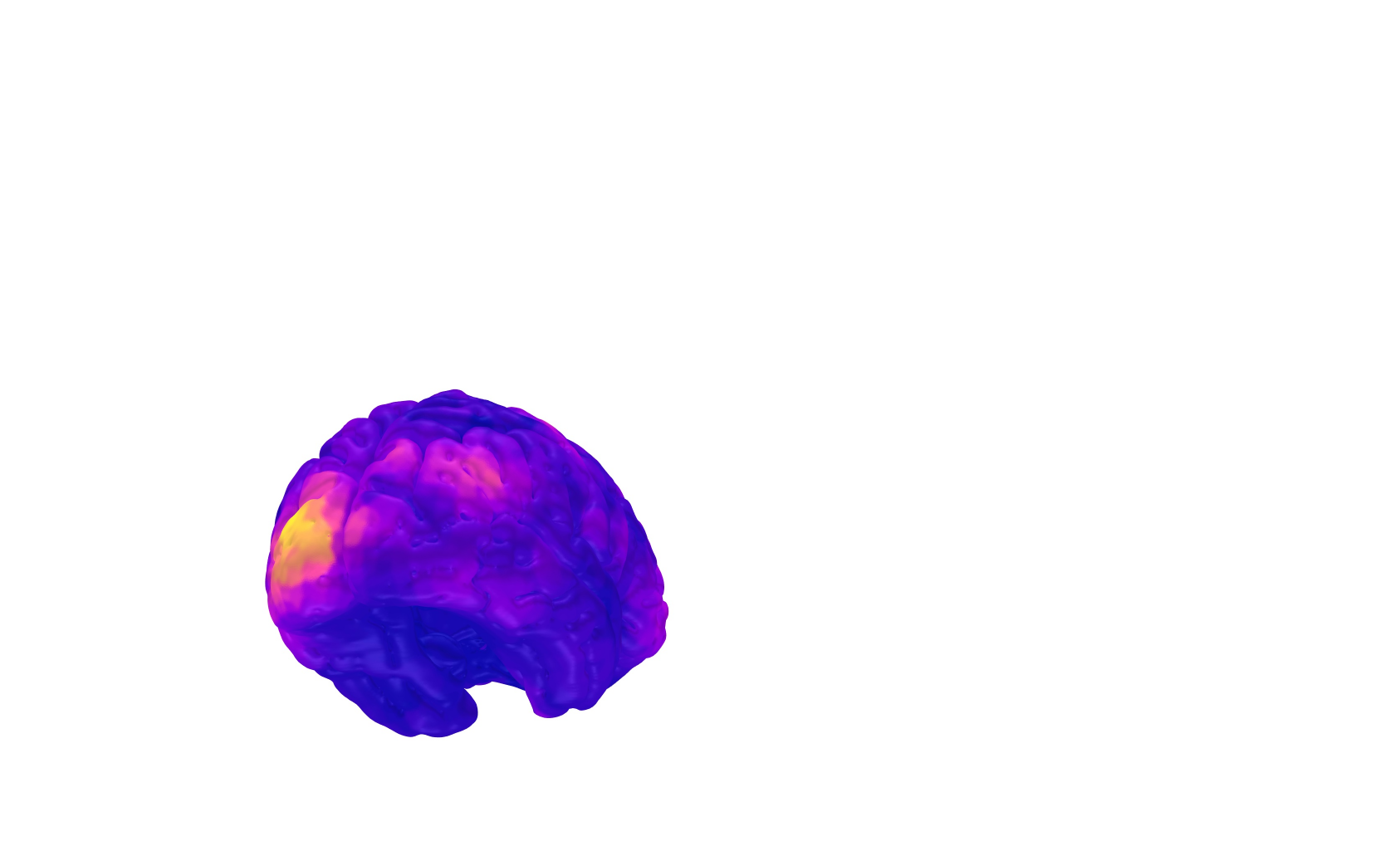}
        \end{minipage}\begin{minipage}{0.16\linewidth}
            \includegraphics[trim={5cm 1cm 14cm 8cm},clip,width=\linewidth]{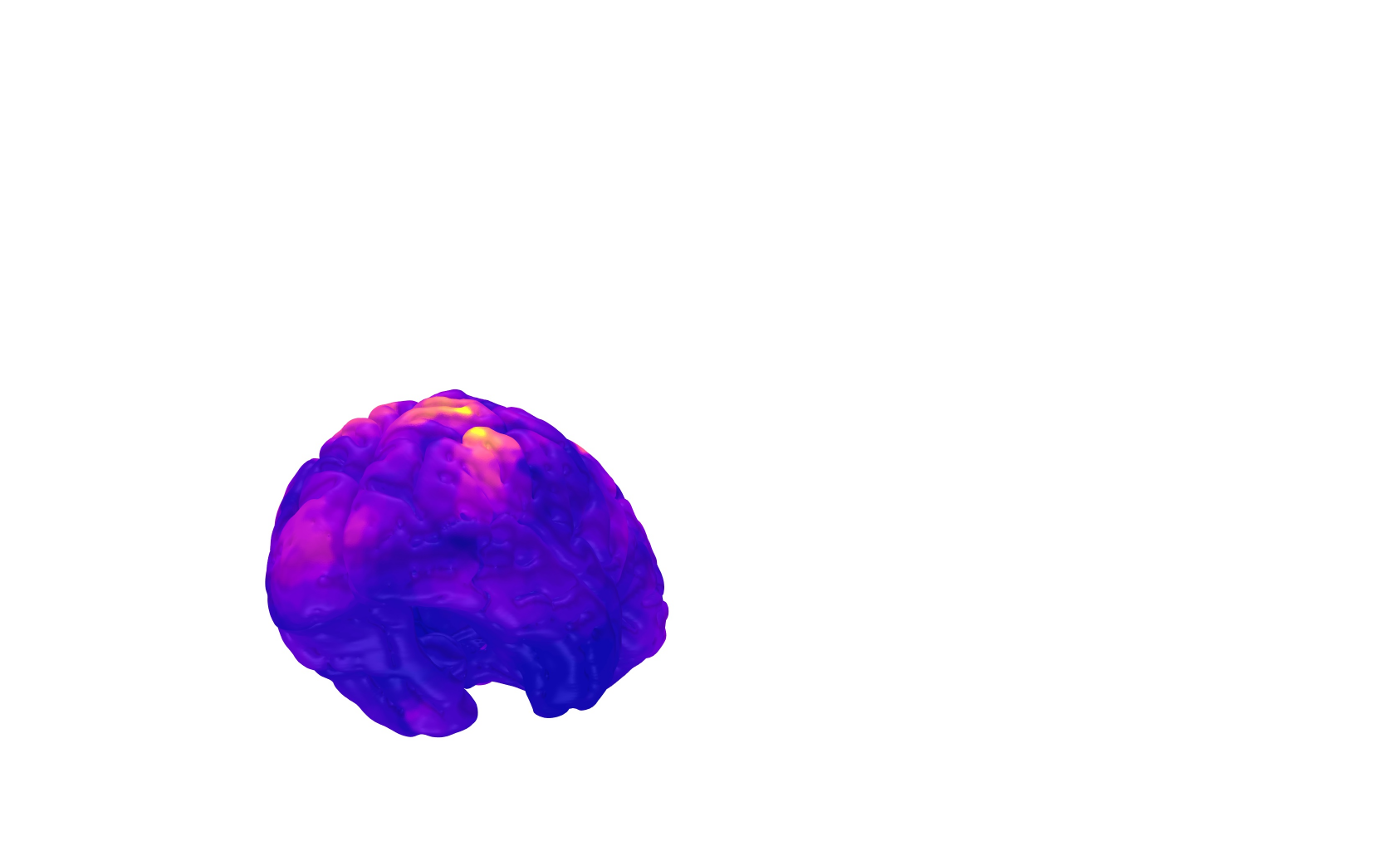}
        \end{minipage}\begin{minipage}{0.16\linewidth}
            \includegraphics[trim={5cm 1cm 14cm 8cm},clip,width=\linewidth]{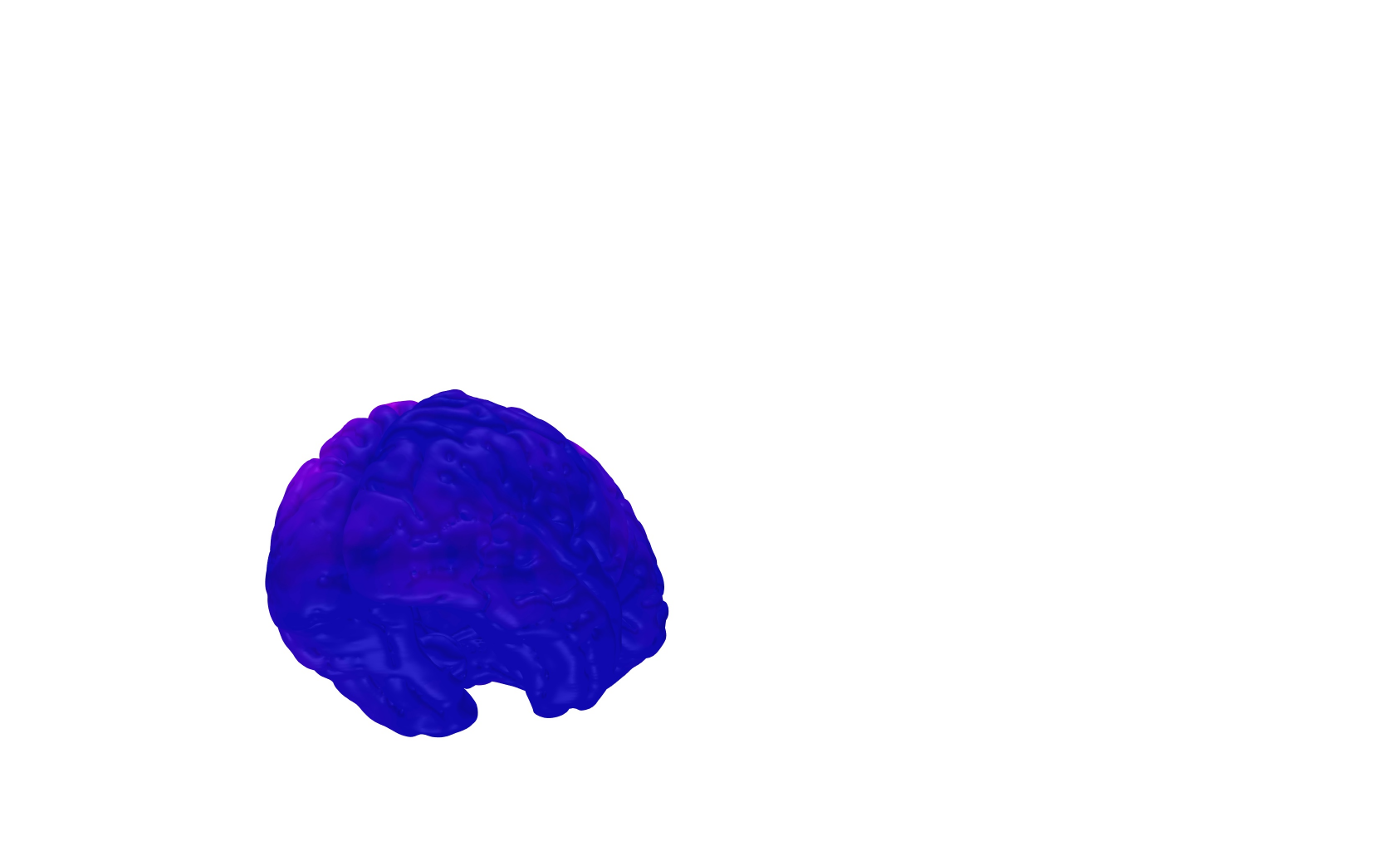}
        \end{minipage}

        \begin{minipage}{0.16\linewidth}
            \includegraphics[trim={5cm 1cm 14cm 8cm},clip,width=\linewidth]{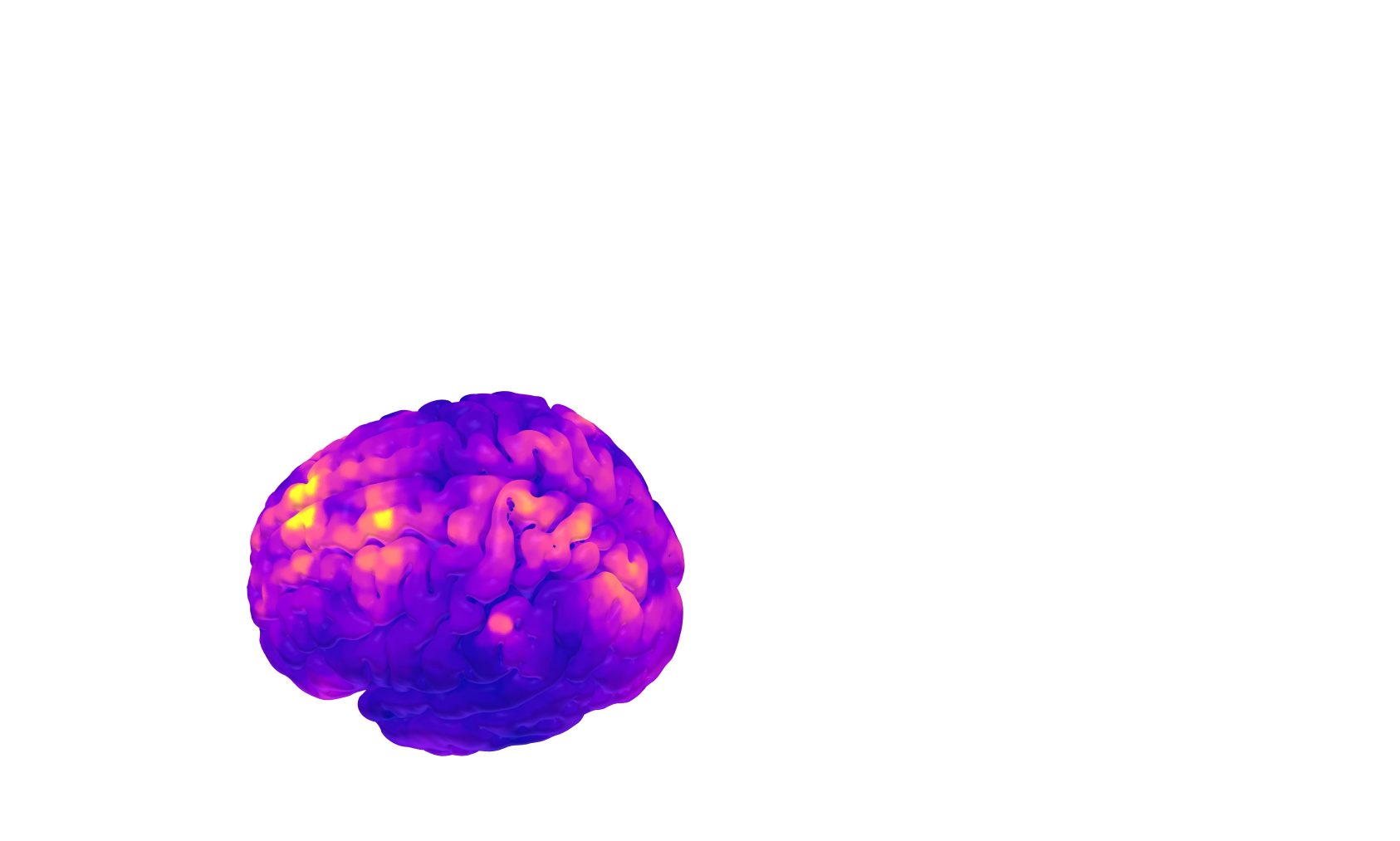}
        \end{minipage}\begin{minipage}{0.16\linewidth}
            \includegraphics[trim={5cm 1cm 14cm 8cm},clip,width=\linewidth]{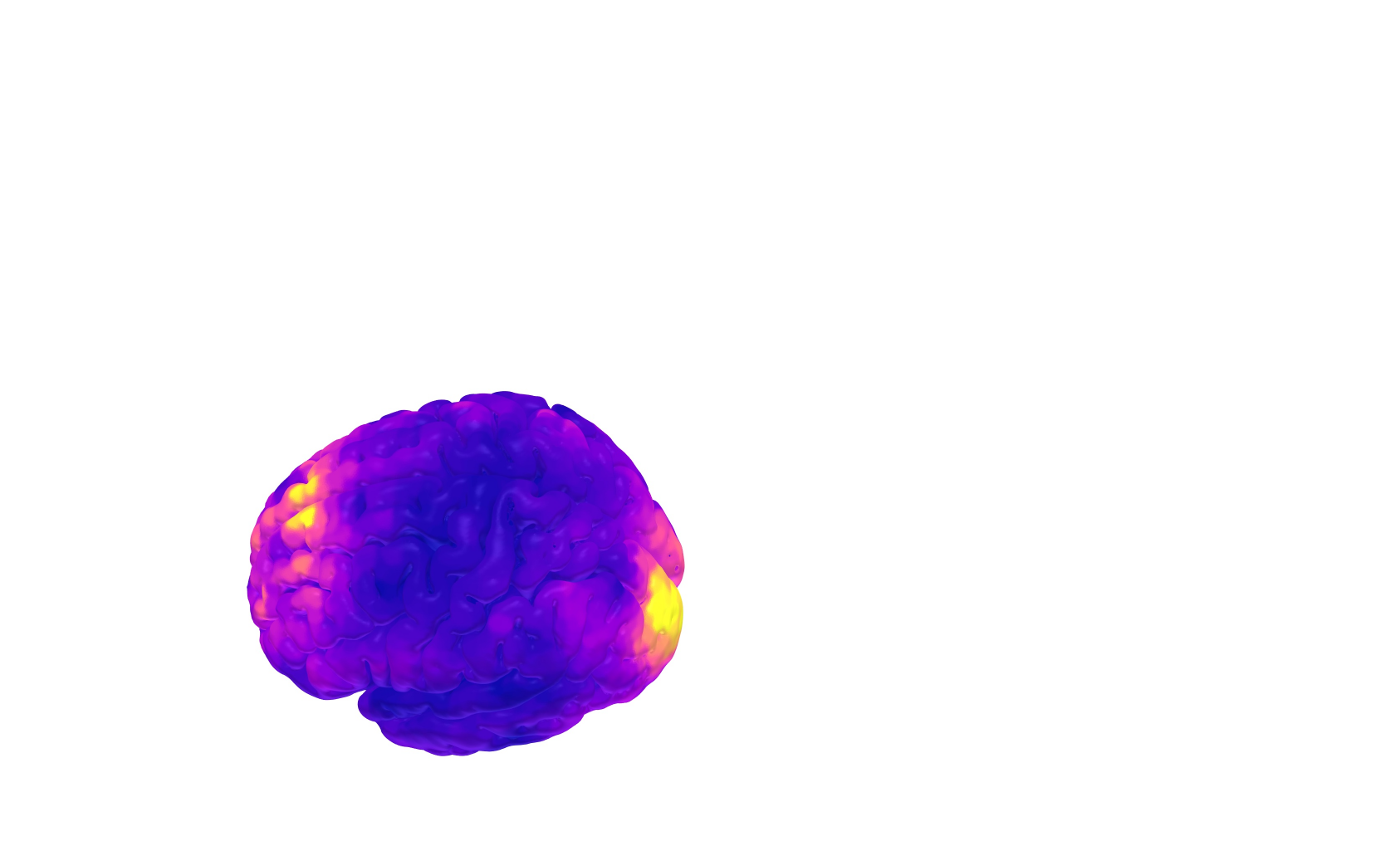}
        \end{minipage}\begin{minipage}{0.16\linewidth}
            \includegraphics[trim={5cm 1cm 14cm 8cm},clip,width=\linewidth]{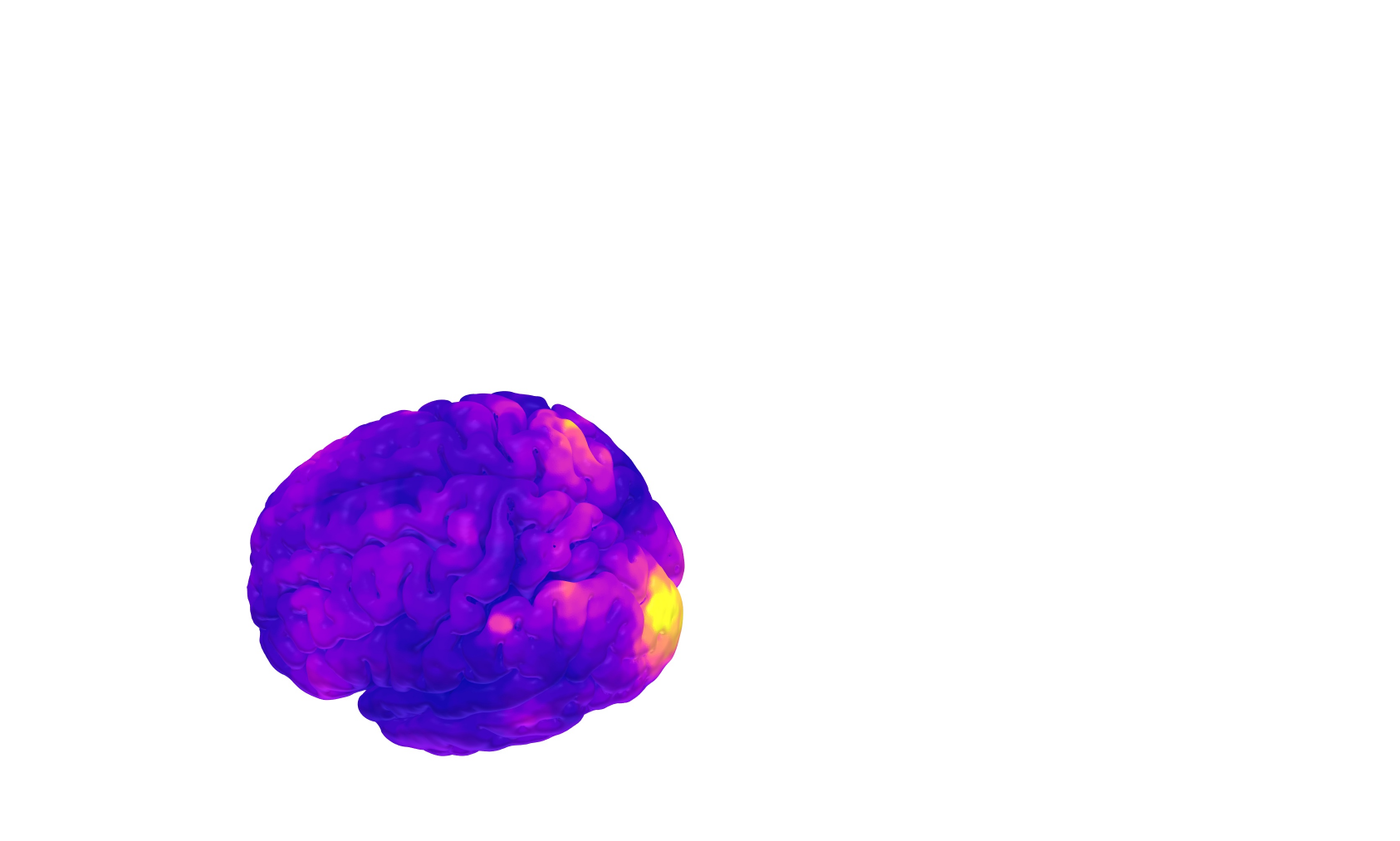}
        \end{minipage}\begin{minipage}{0.16\linewidth}
            \includegraphics[trim={5cm 1cm 14cm 8cm},clip,width=\linewidth]{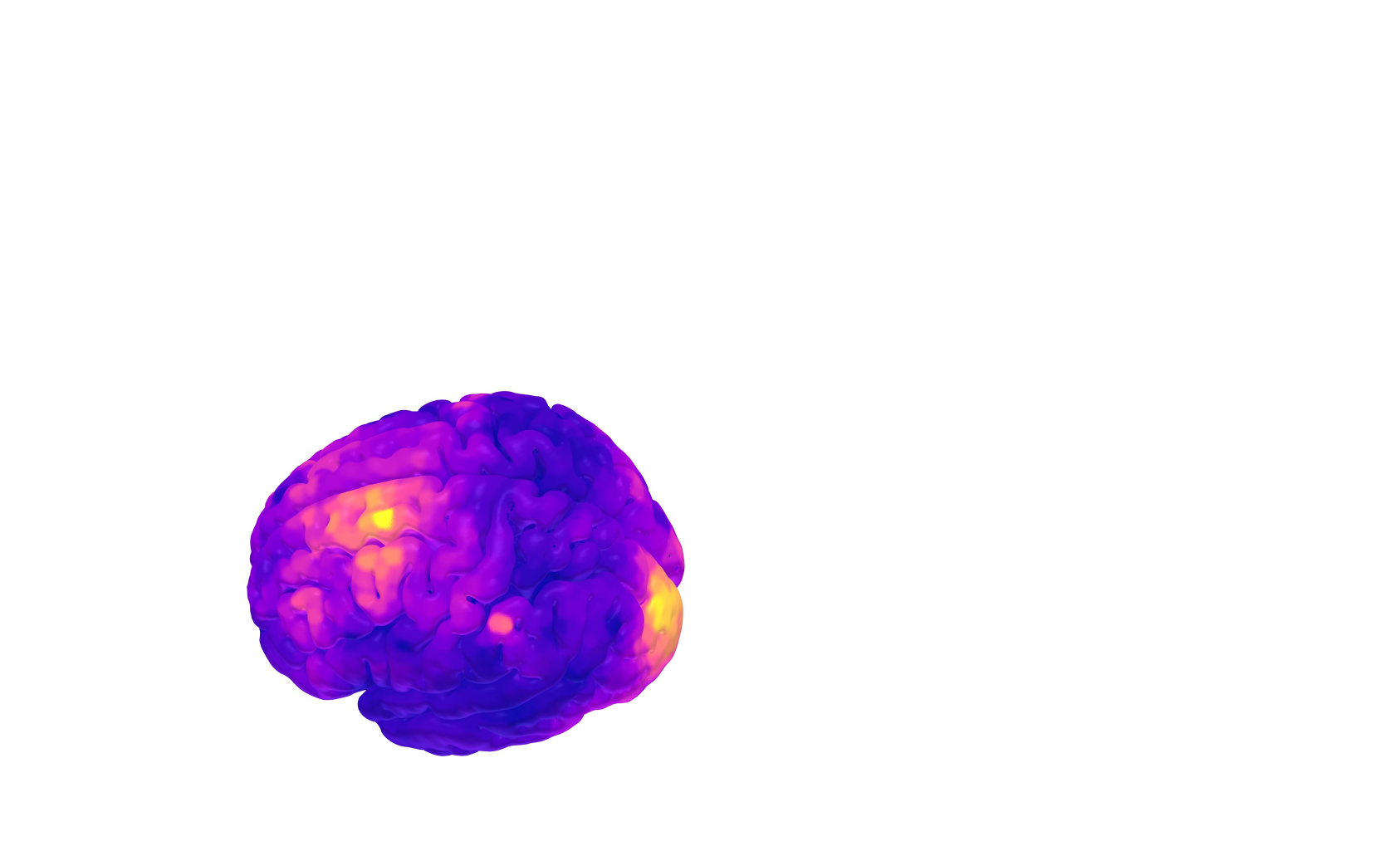}
        \end{minipage}\begin{minipage}{0.16\linewidth}
            \includegraphics[trim={5cm 1cm 14cm 8cm},clip,width=\linewidth]{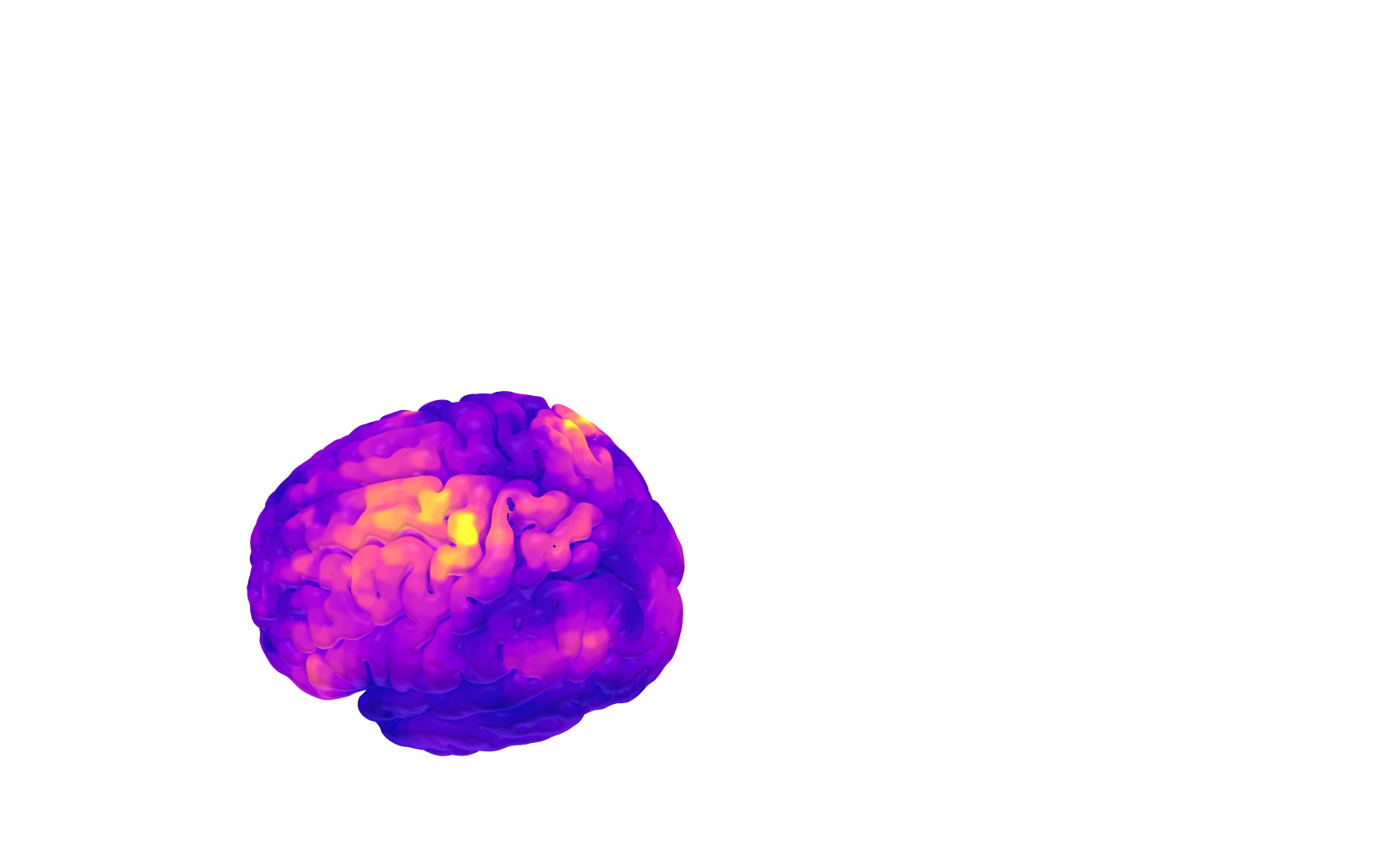}
        \end{minipage}\begin{minipage}{0.16\linewidth}
            \includegraphics[trim={5cm 1cm 14cm 8cm},clip,width=\linewidth]{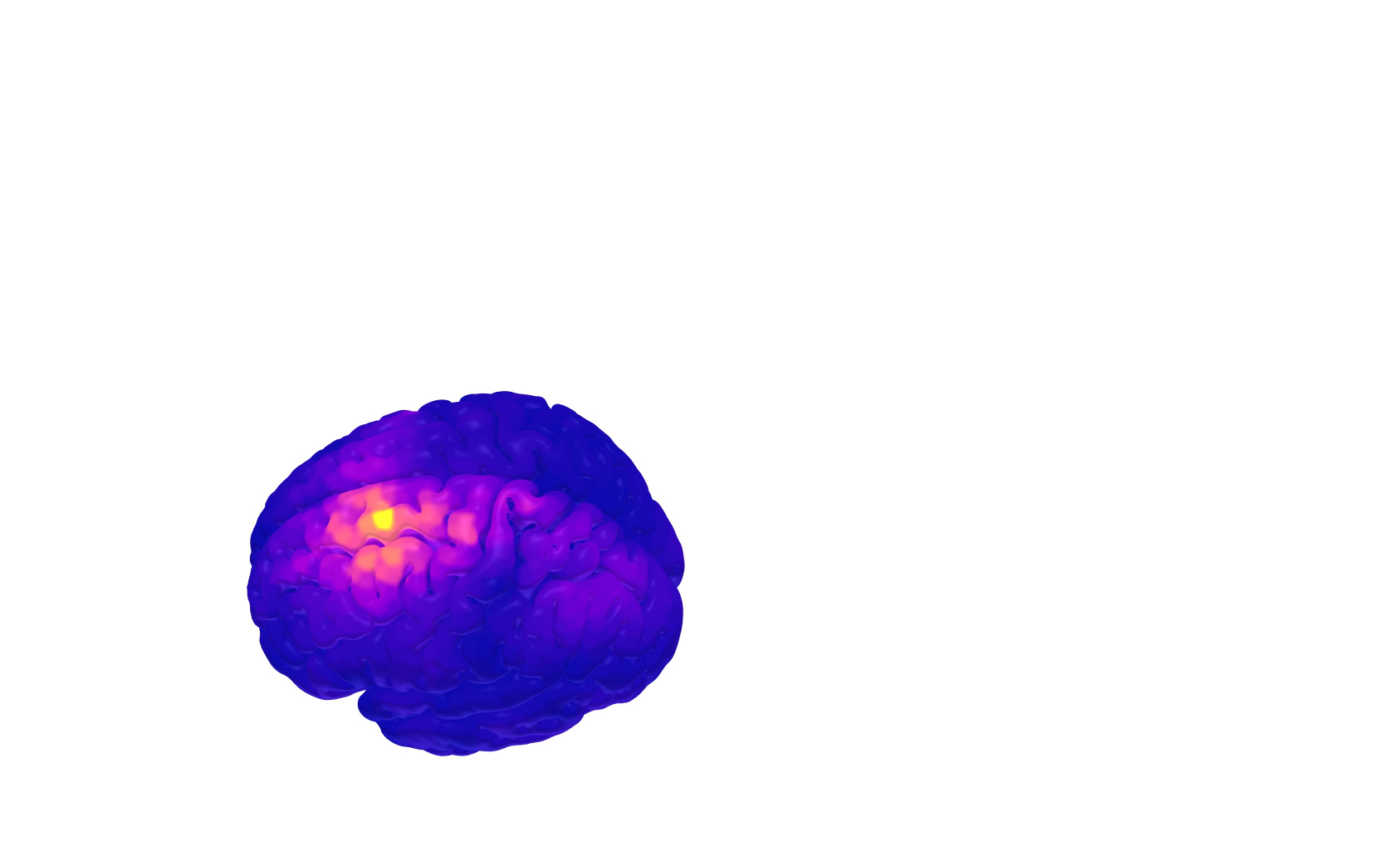}
        \end{minipage}
    \end{minipage}
    \noindent\rule[7pt]{\linewidth}{0.4pt}

    \begin{minipage}{0.02\linewidth}
        \rotatebox{90}{\scriptsize LCMV}
    \end{minipage}\begin{minipage}{0.98\linewidth}
        \begin{minipage}{0.16\linewidth}
            \includegraphics[trim={5cm 1cm 14cm 8cm},clip,width=\linewidth]{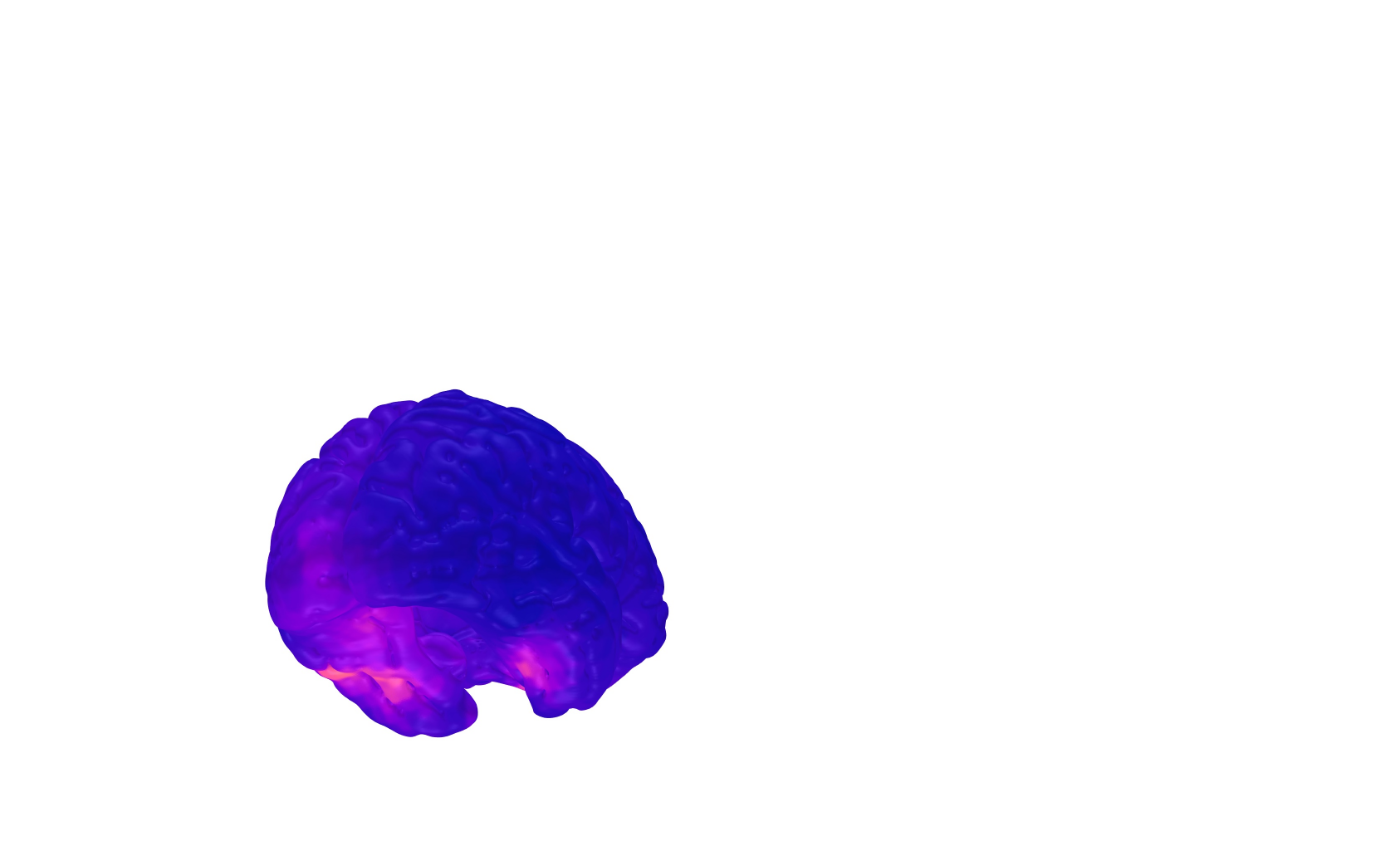}
        \end{minipage}\begin{minipage}{0.16\linewidth}
            \includegraphics[trim={5cm 1cm 14cm 8cm},clip,width=\linewidth]{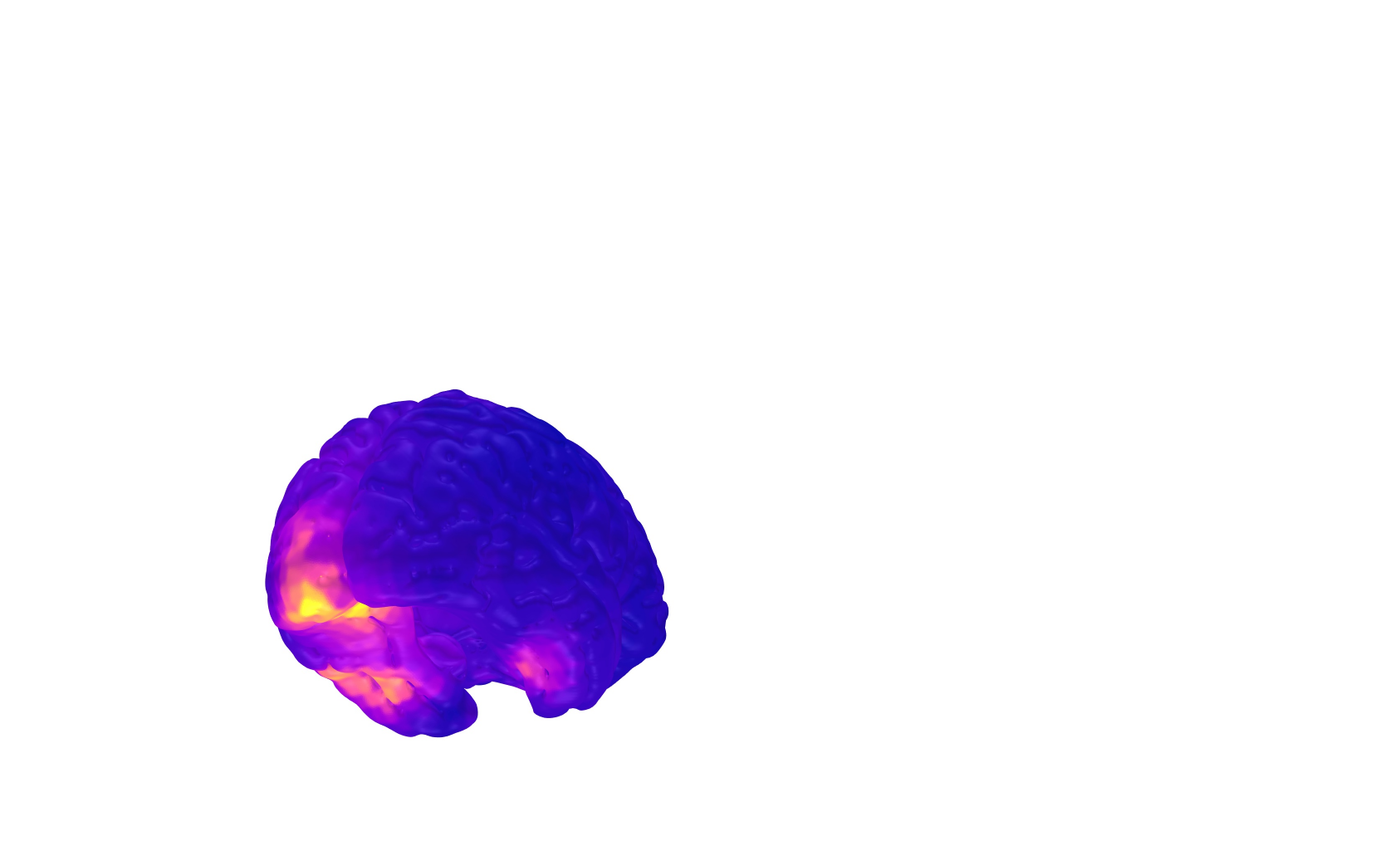}
        \end{minipage}\begin{minipage}{0.16\linewidth}
            \includegraphics[trim={5cm 1cm 14cm 8cm},clip,width=\linewidth]{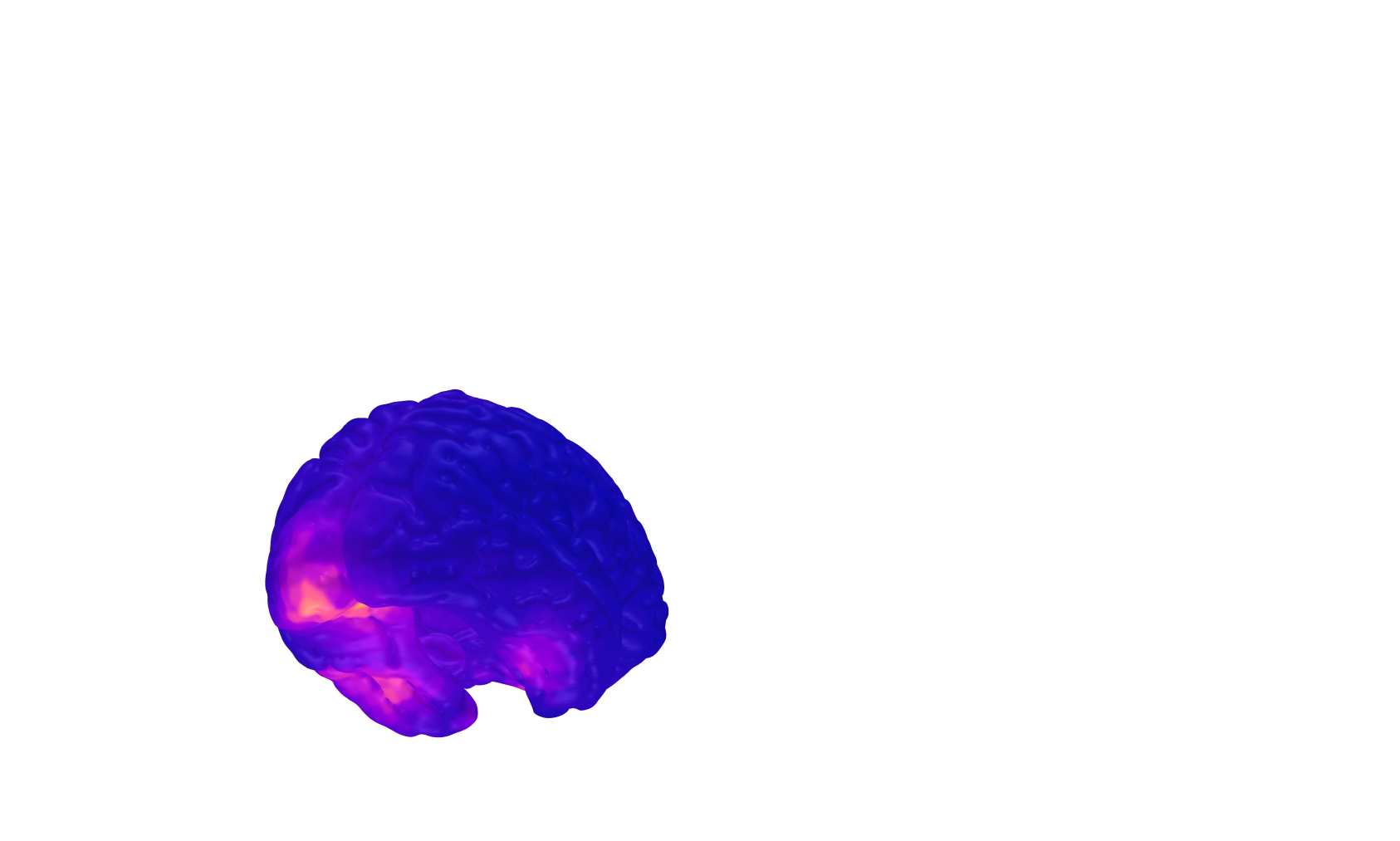}
        \end{minipage}\begin{minipage}{0.16\linewidth}
            \includegraphics[trim={5cm 1cm 14cm 8cm},clip,width=\linewidth]{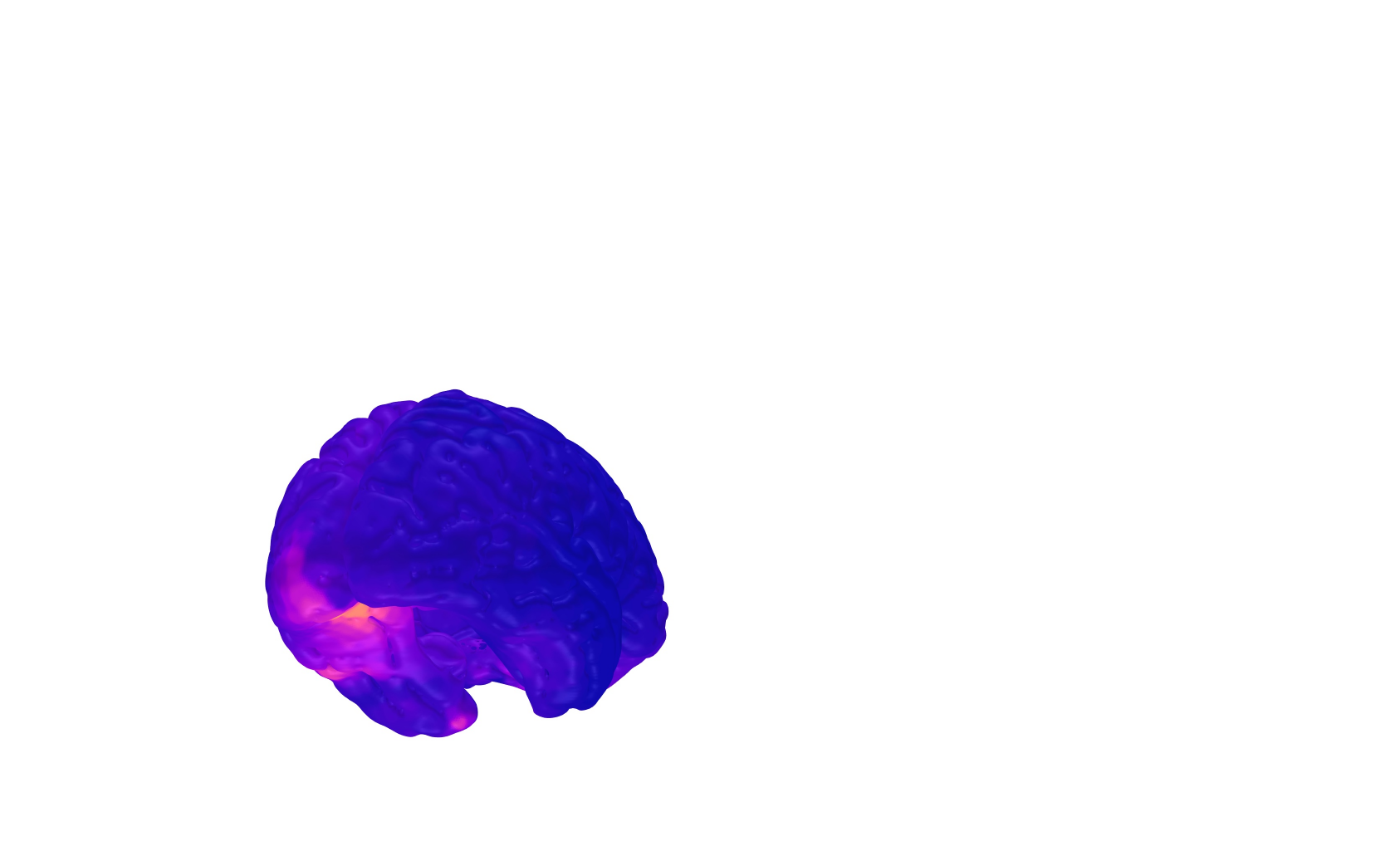}
        \end{minipage}\begin{minipage}{0.16\linewidth}
            \includegraphics[trim={5cm 1cm 14cm 8cm},clip,width=\linewidth]{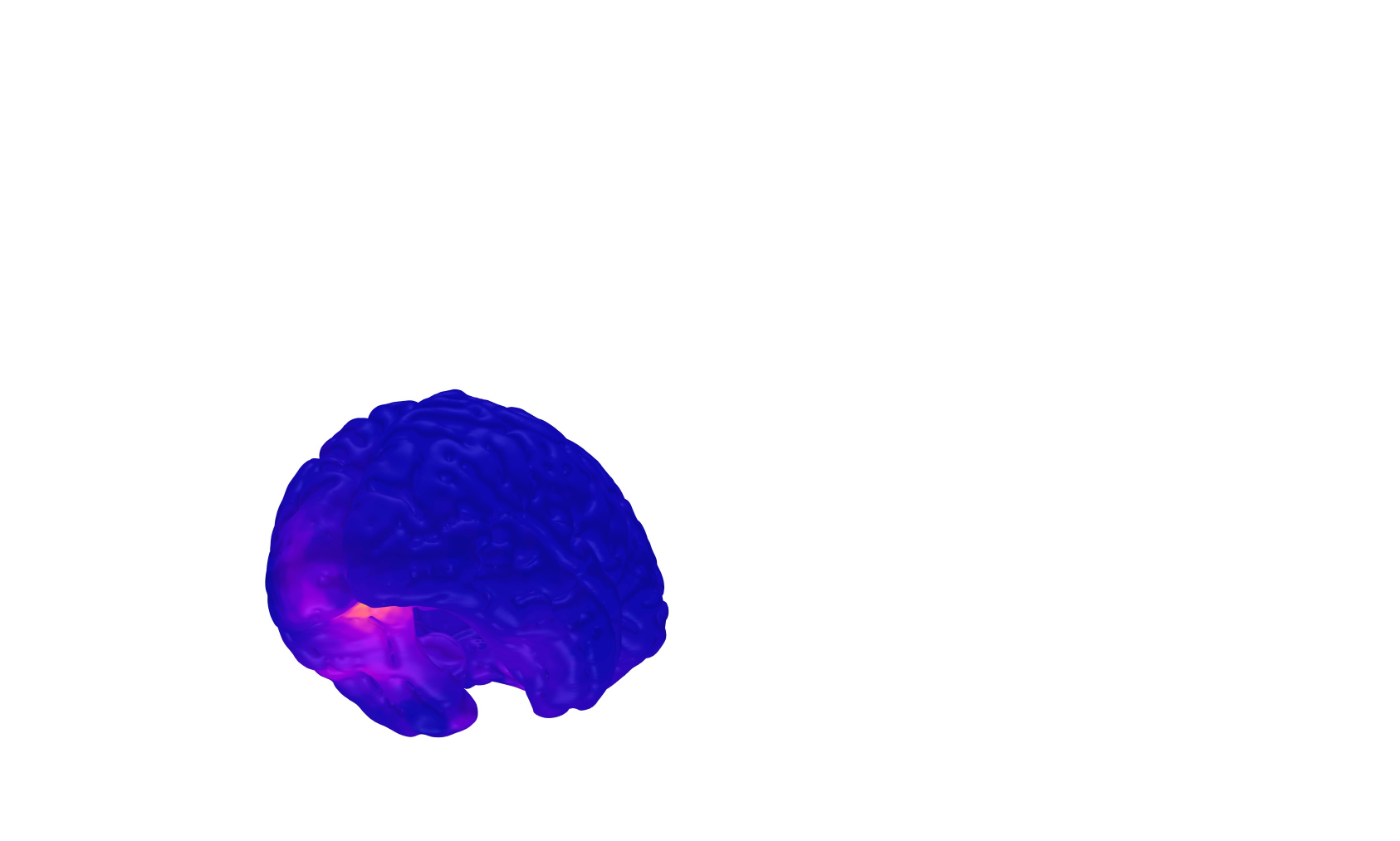}
        \end{minipage}\begin{minipage}{0.16\linewidth}
            \includegraphics[trim={5cm 1cm 14cm 8cm},clip,width=\linewidth]{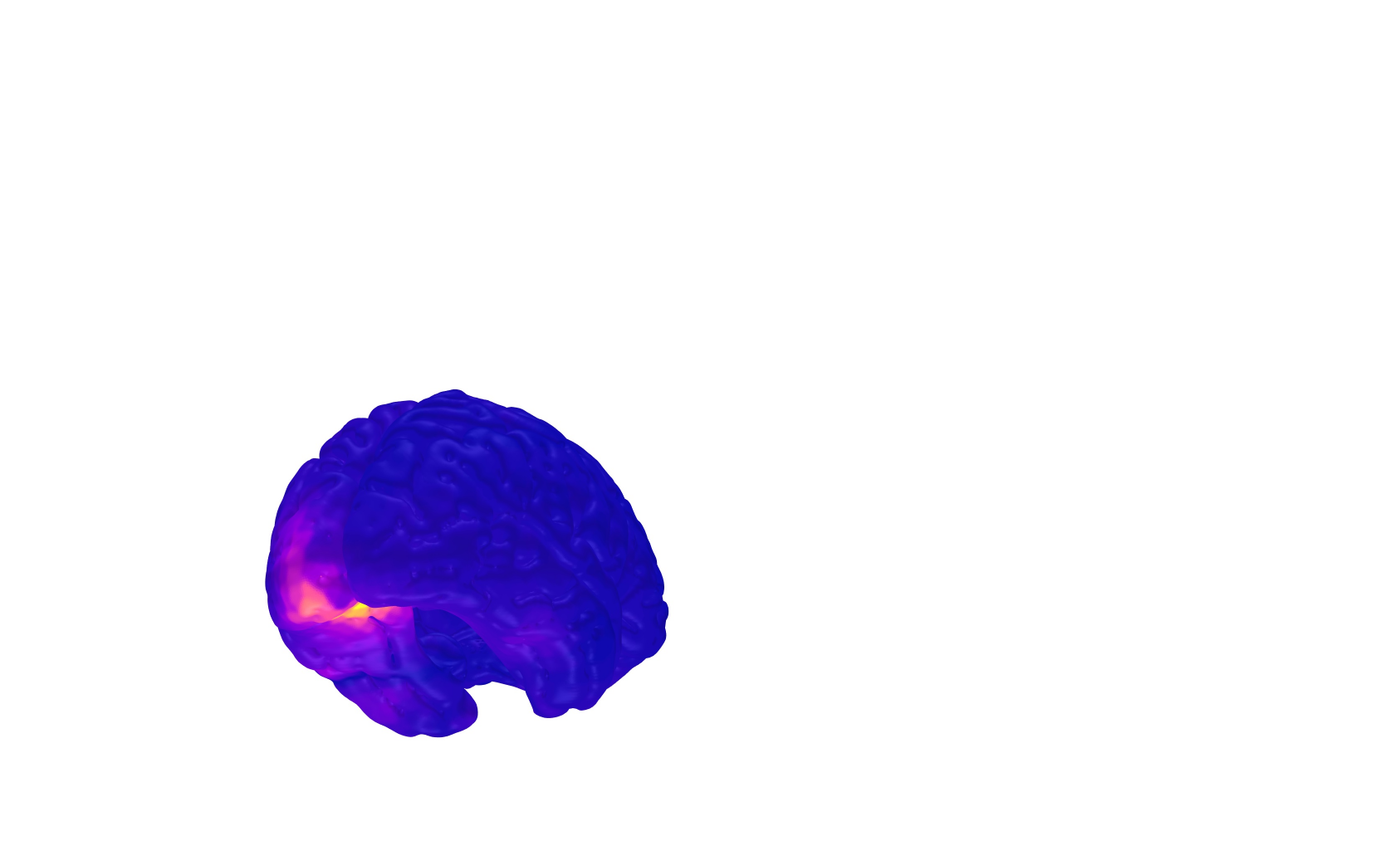}
        \end{minipage}

        \begin{minipage}{0.16\linewidth}
            \includegraphics[trim={5cm 1cm 14cm 8cm},clip,width=\linewidth]{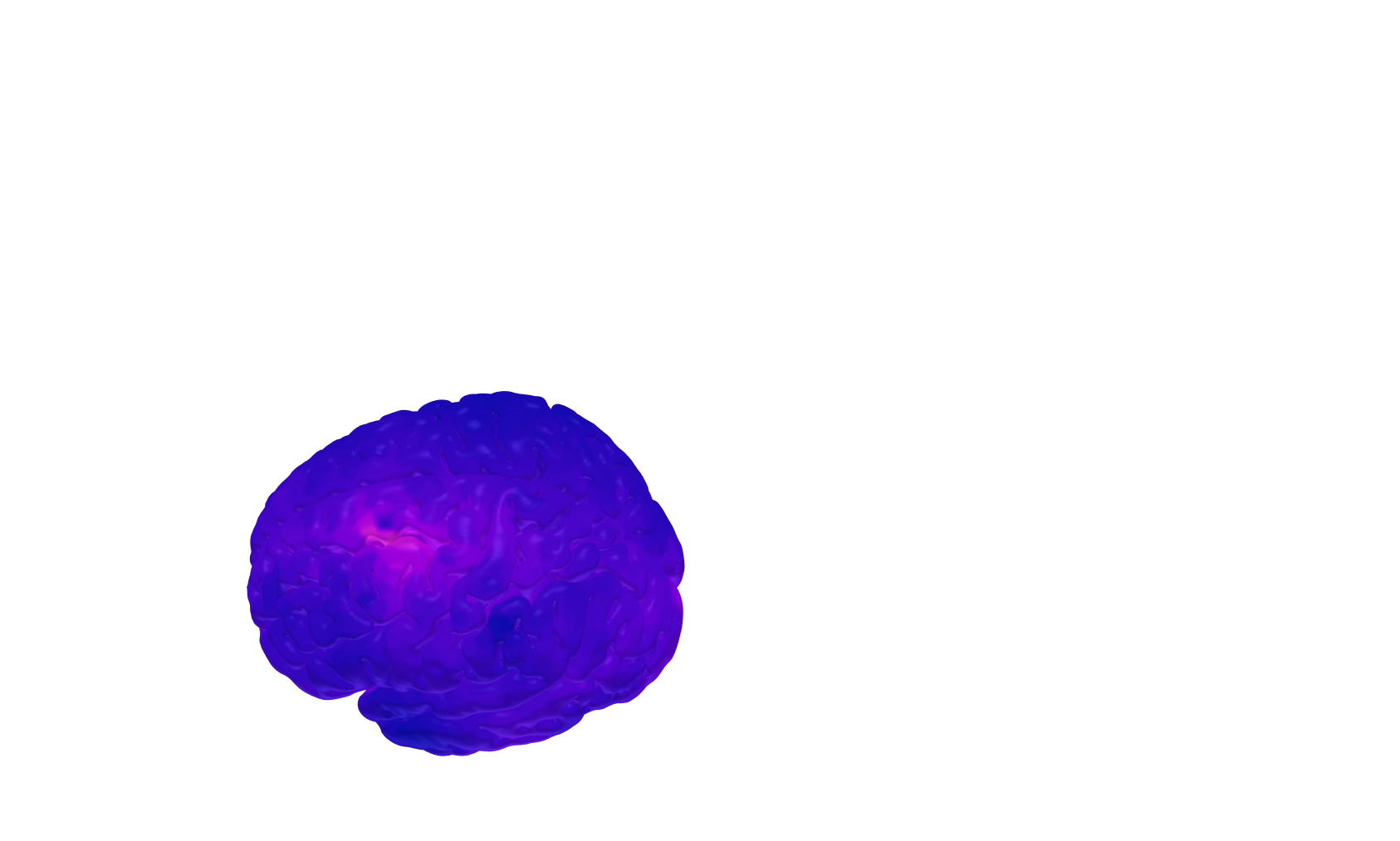}
        \end{minipage}\begin{minipage}{0.16\linewidth}
            \includegraphics[trim={5cm 1cm 14cm 8cm},clip,width=\linewidth]{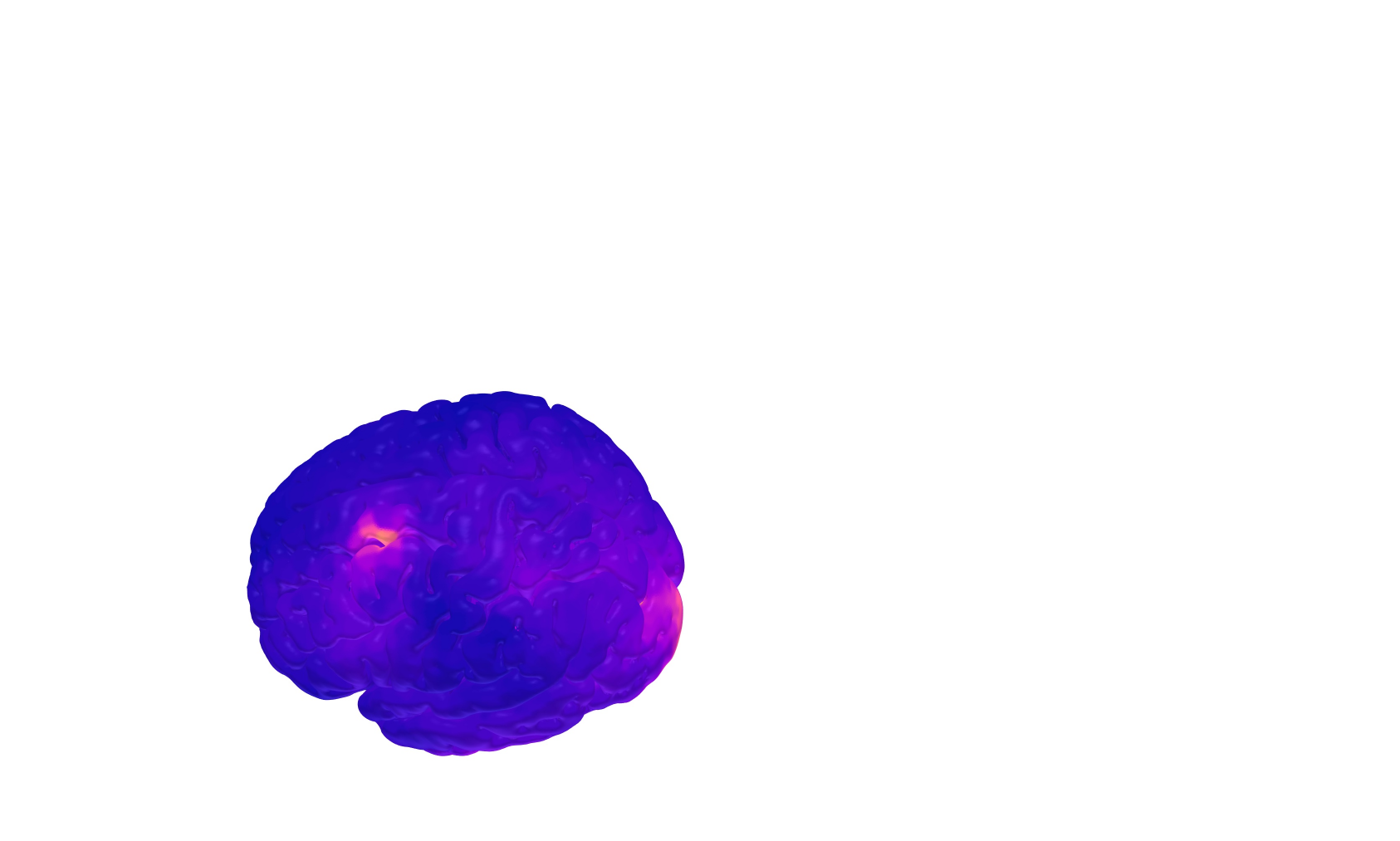}
        \end{minipage}\begin{minipage}{0.16\linewidth}
            \includegraphics[trim={5cm 1cm 14cm 8cm},clip,width=\linewidth]{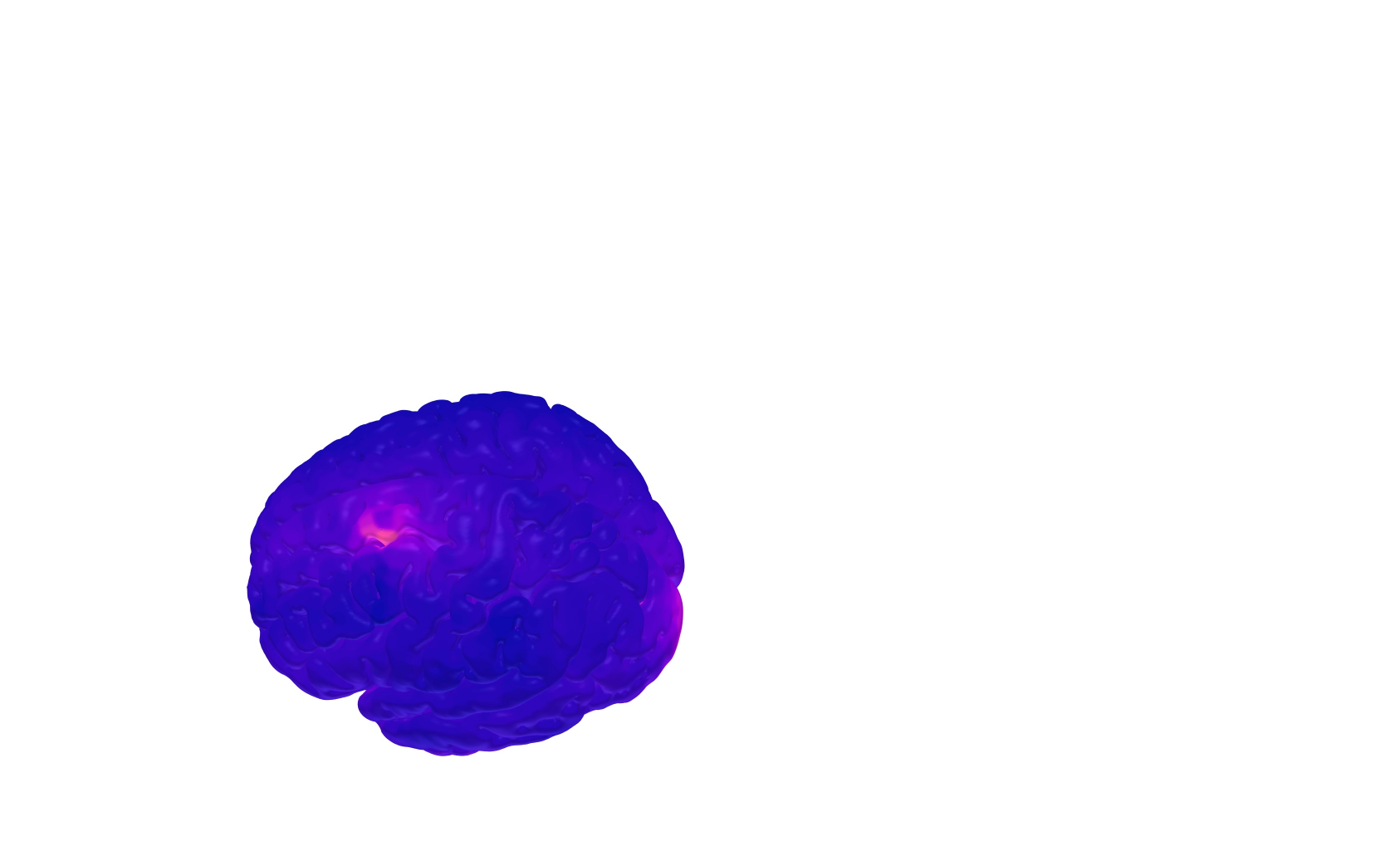}
        \end{minipage}\begin{minipage}{0.16\linewidth}
            \includegraphics[trim={5cm 1cm 14cm 8cm},clip,width=\linewidth]{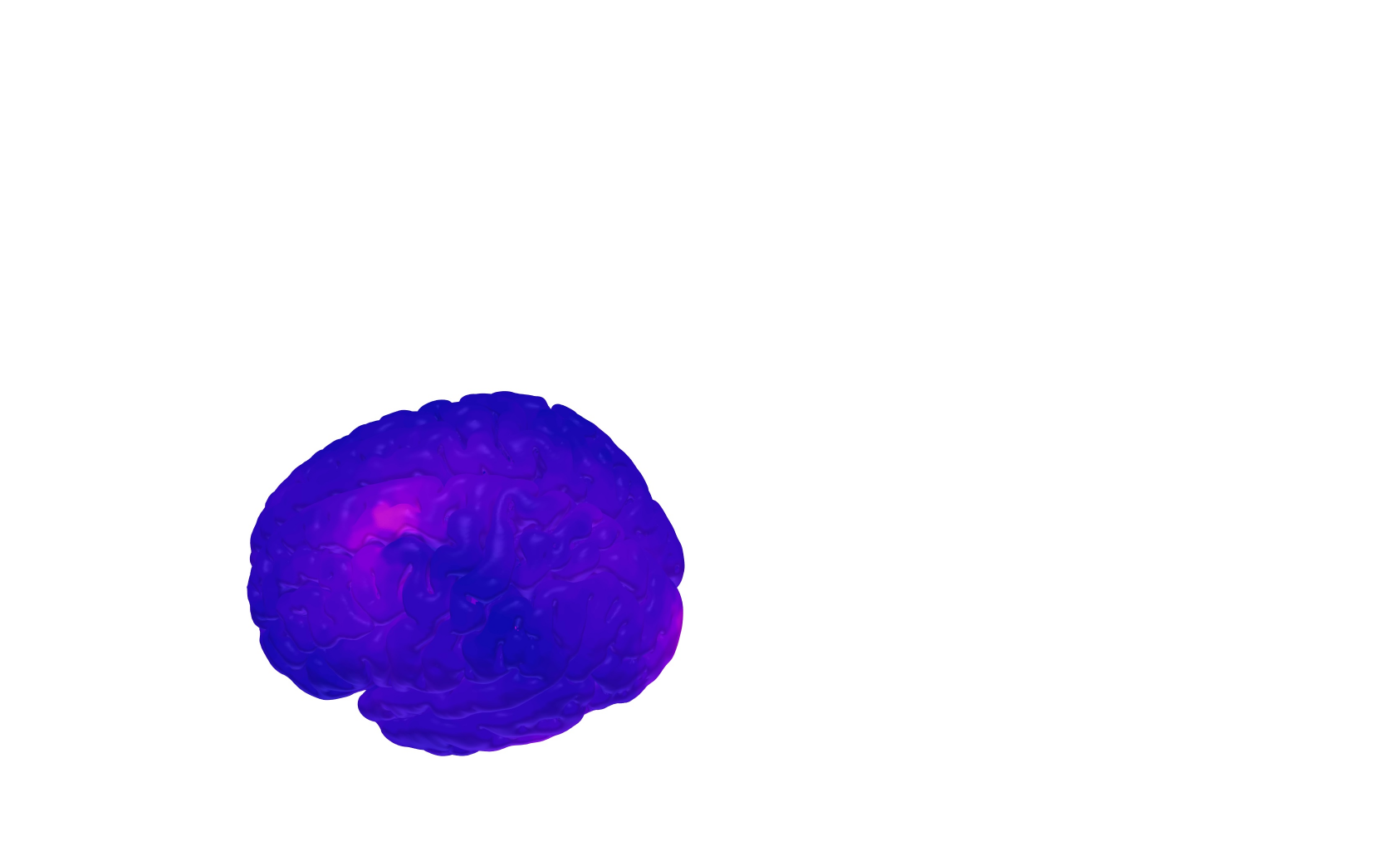}
        \end{minipage}\begin{minipage}{0.16\linewidth}
            \includegraphics[trim={5cm 1cm 14cm 8cm},clip,width=\linewidth]{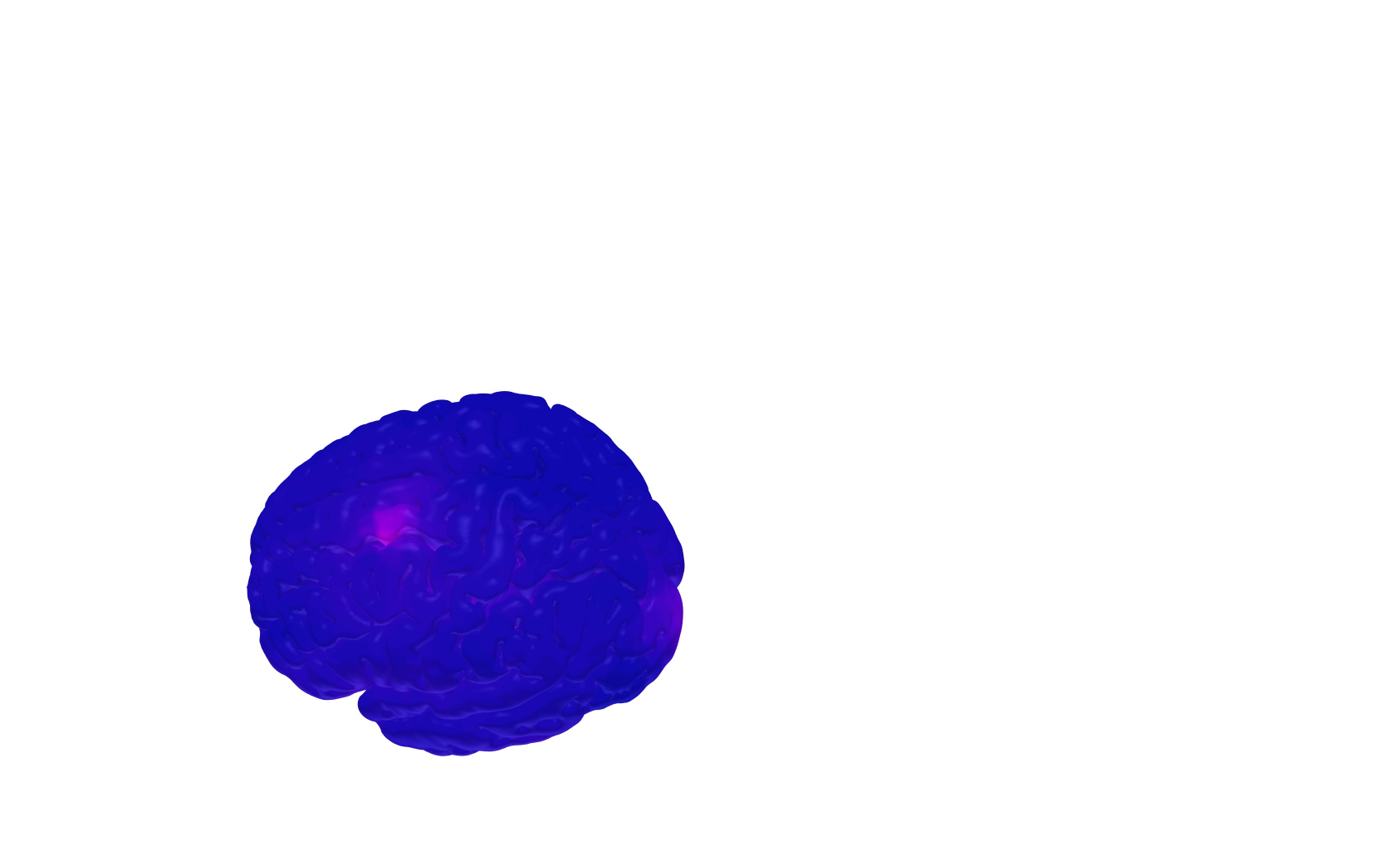}
        \end{minipage}\begin{minipage}{0.16\linewidth}
            \includegraphics[trim={5cm 1cm 14cm 8cm},clip,width=\linewidth]{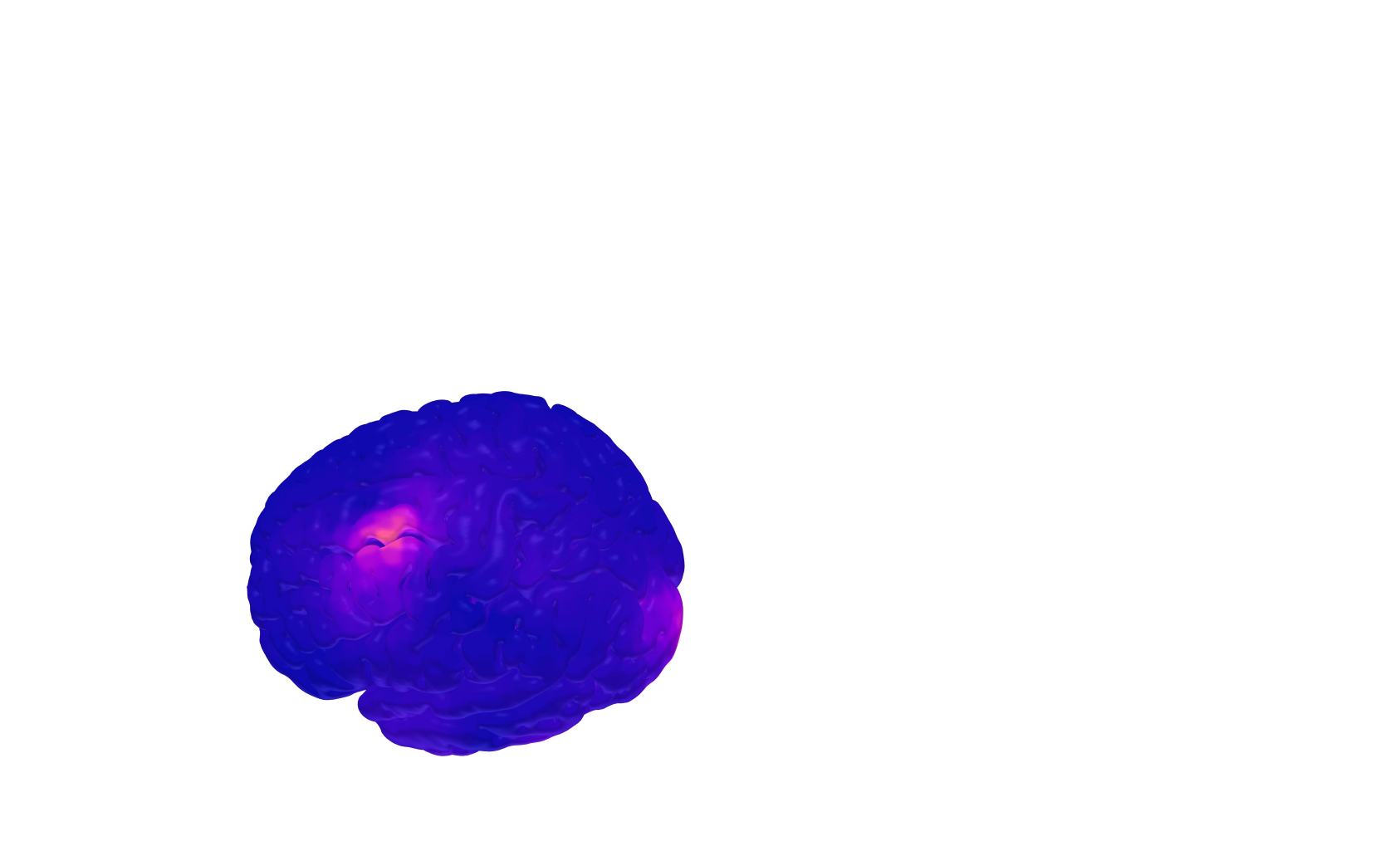}
        \end{minipage}
    \end{minipage}
    
    \caption{Spatial source reconstructions at 20 dB SNR. Rows correspond to Change-rate Standardized Kalman filter (CrSKF), Standardized Kalman filter (SKF), Kalman filter (KF), and Linearly Constrained Minimum Variance (LCMV) beamformer; columns show the reconstructions at the modeled FVEP component times associated with V1, V2d, V2v+V3d, V3v, LIP, and FEF.}
    \label{fig:brain_20dB}
\end{figure}

\subsection{Spatial reconstruction and distributional error}

Figures \ref{fig:brain_noiseless}--\ref{fig:brain_20dB} compare the spatial source distributions at the modeled FVEP component times. The standardized Kalman reconstructions remain more spatially concentrated than conventional KF across the displayed conditions, whereas LCMV shows stronger variation in spatial extent as sensor noise increases. The spatial maps also confirm that activity from nearby simulated sources is not confined to a single anatomical parcel, consistent with the parcel-wise mixing observed in the temporal curves.

To quantify spatial mismatch, the reconstructed and ground-truth source magnitudes are first restricted to nonnegative values and normalized to unit mass. Earth Mover's Distance (EMD) is then defined as the minimum transport cost required to transform the reconstructed distribution into the ground-truth distribution, using Euclidean distances between source coordinates as transport costs. Lower EMD therefore indicates closer spatial agreement independent of total amplitude. 

\begin{figure}[h!]
    \centering
    \begin{minipage}{0.02\linewidth}
        \rotatebox{90}{\scriptsize \hfill}
    \end{minipage}\begin{minipage}{0.245\linewidth}
        \centering
        {\bf CrSKF}
    \end{minipage}\begin{minipage}{0.245\linewidth}
        \centering
        {\bf SKF}
    \end{minipage}\begin{minipage}{0.245\linewidth}
        \centering
        {\bf KF}
    \end{minipage}\begin{minipage}{0.245\linewidth}
        \centering
        {\bf LCMV}
    \end{minipage}

    \begin{minipage}{0.02\linewidth}
        \rotatebox{90}{\scriptsize L.O.}
    \end{minipage}\begin{minipage}{0.245\linewidth}
        \centering
        \includegraphics[width=0.99\linewidth]{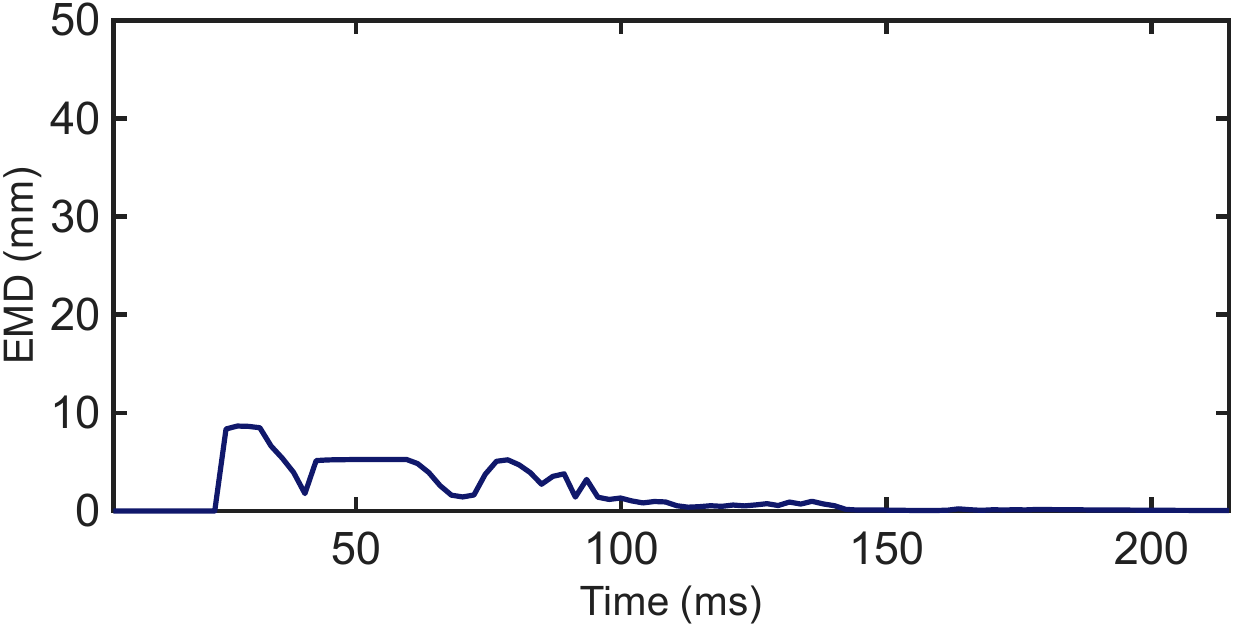}
    \end{minipage}\begin{minipage}{0.245\linewidth}
        \centering
        \includegraphics[width=0.99\linewidth]{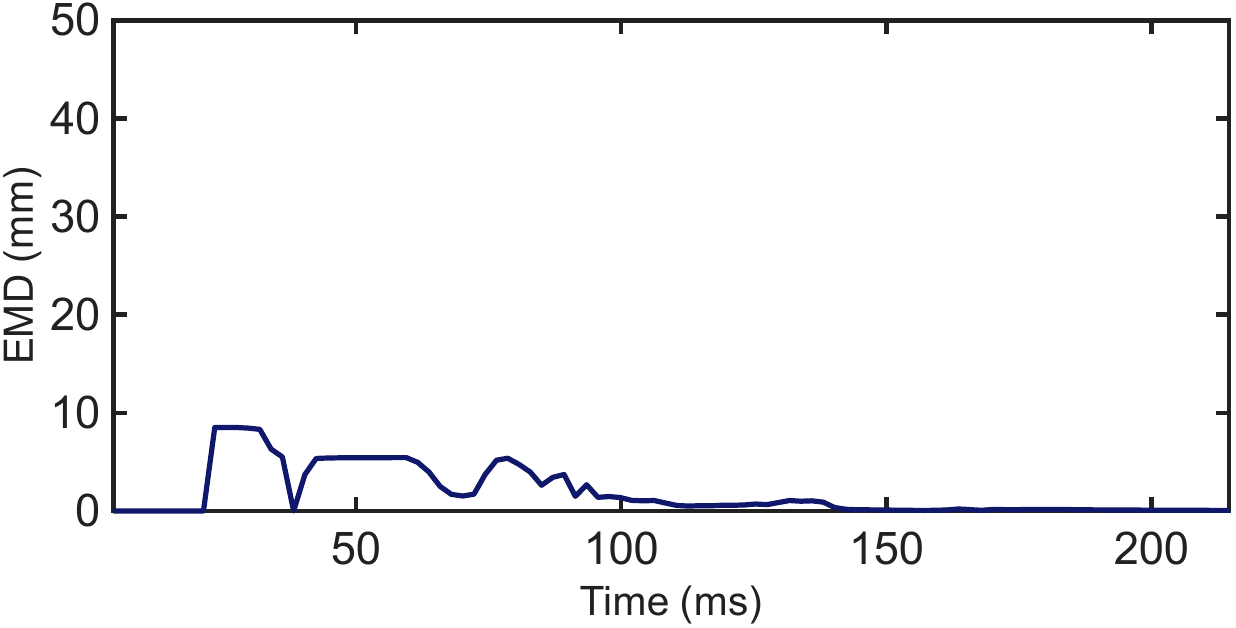}
    \end{minipage}\begin{minipage}{0.245\linewidth}
        \centering
        \includegraphics[width=0.99\linewidth]{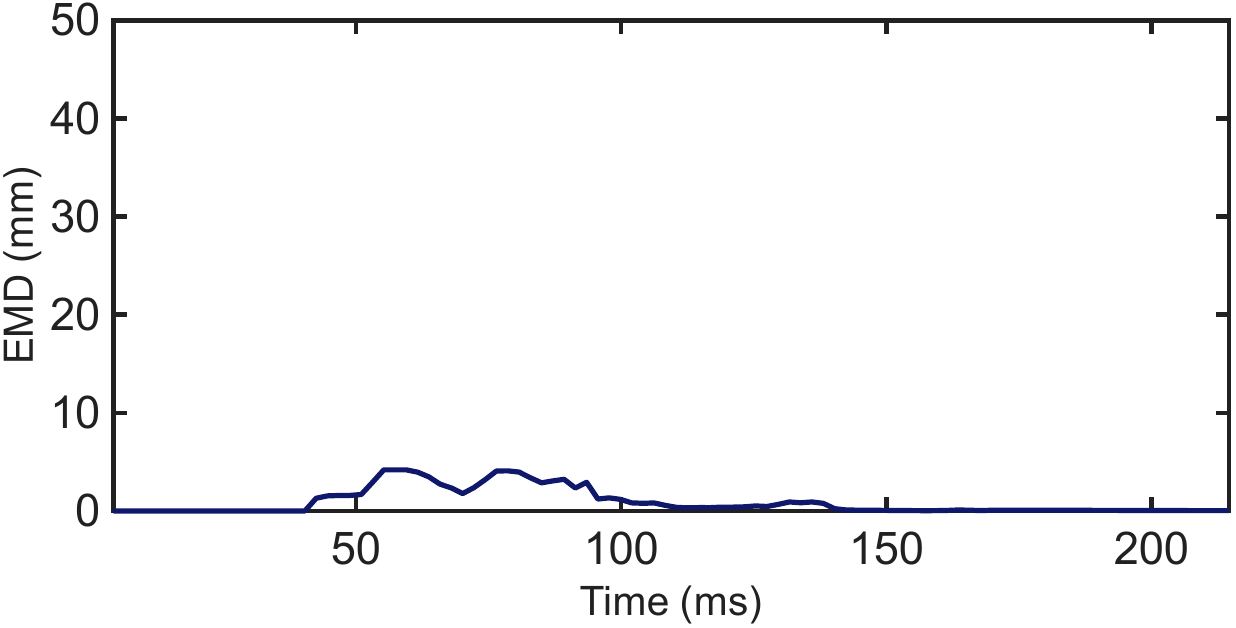}
    \end{minipage}\begin{minipage}{0.245\linewidth}
        \centering
        \includegraphics[width=0.99\linewidth]{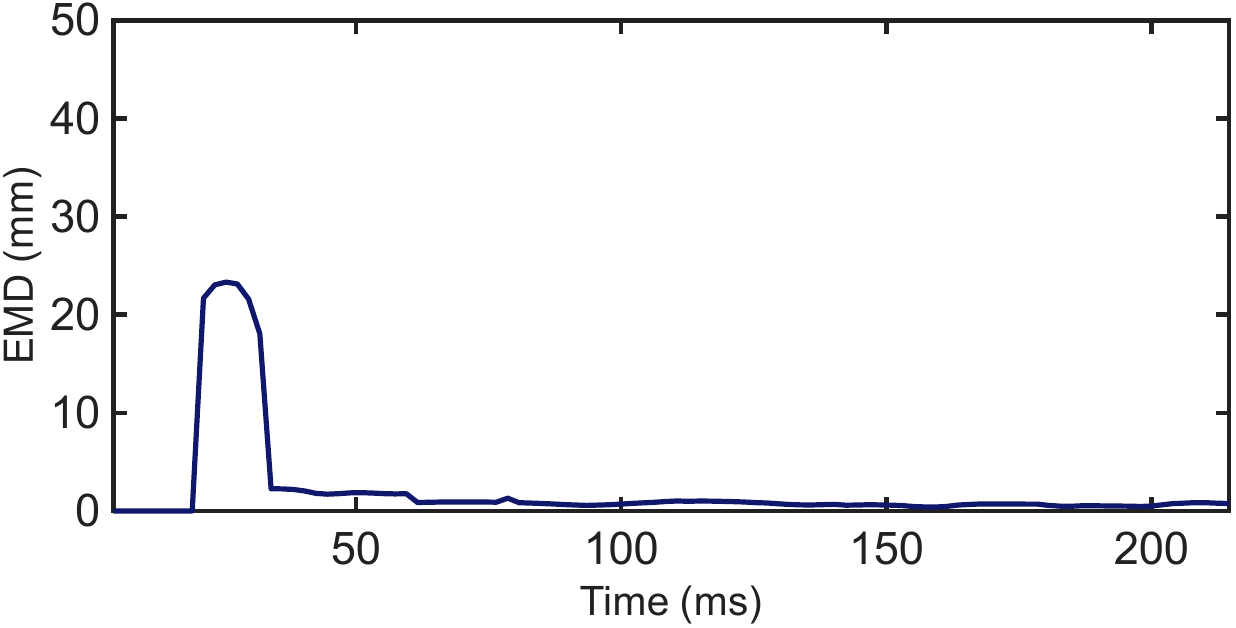}
    \end{minipage}

    \begin{minipage}{0.02\linewidth}
        \rotatebox{90}{\scriptsize Ling.}
    \end{minipage}\begin{minipage}{0.245\linewidth}
        \centering
        \includegraphics[width=0.99\linewidth]{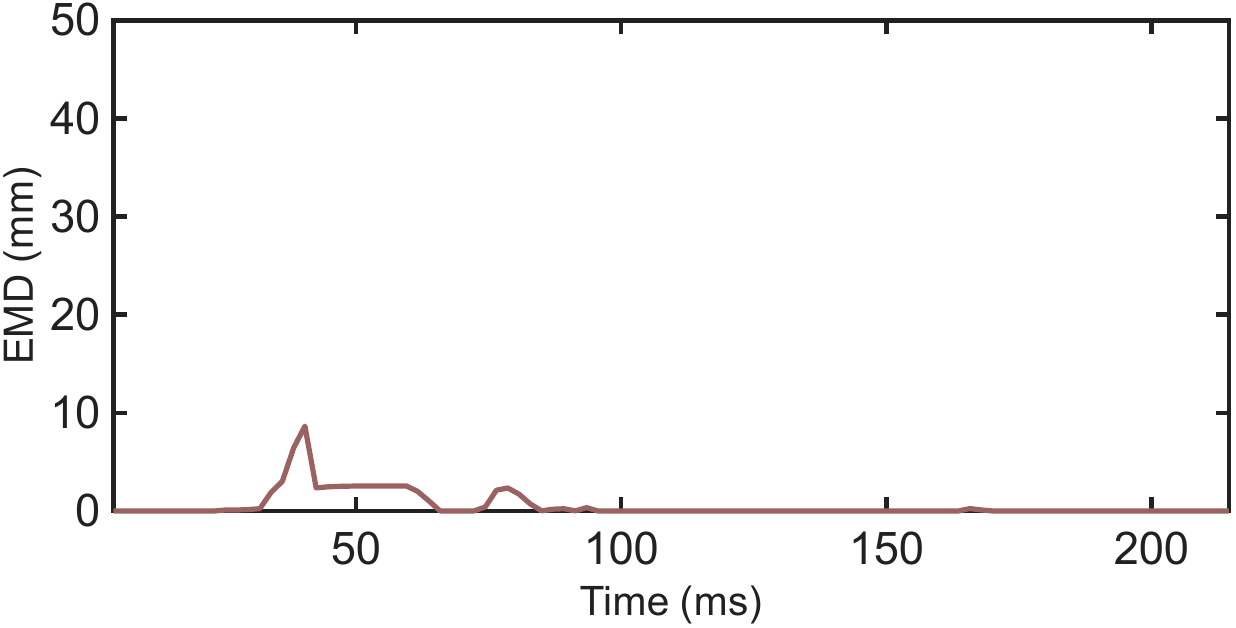}
    \end{minipage}\begin{minipage}{0.245\linewidth}
        \centering
        \includegraphics[width=0.99\linewidth]{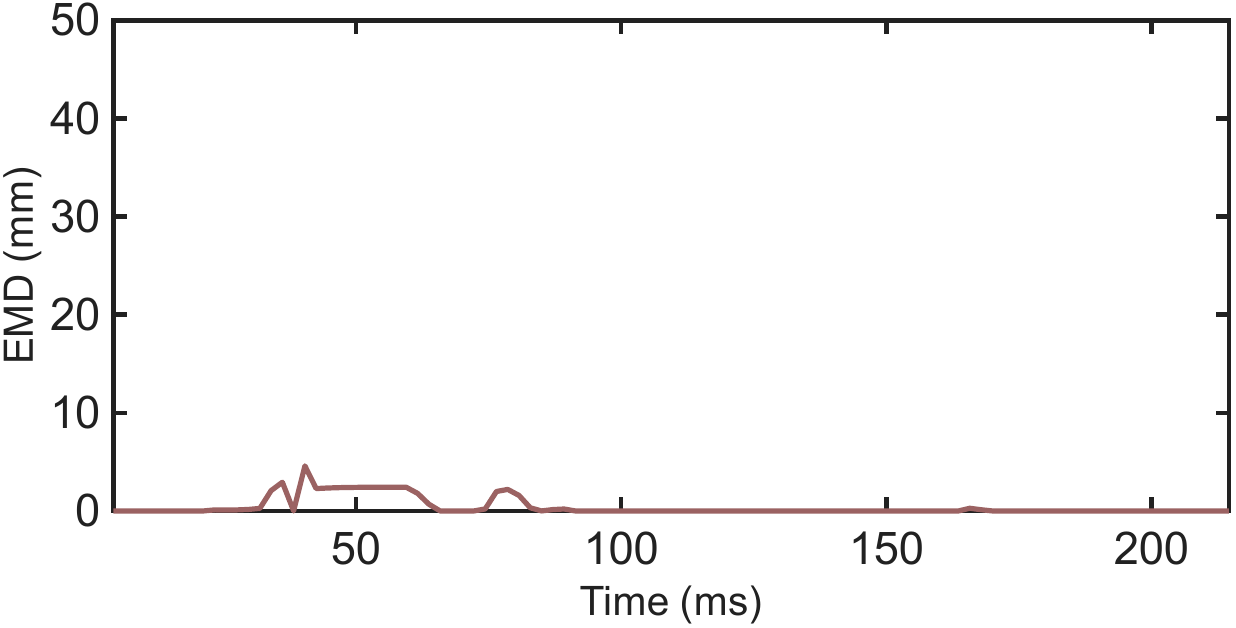}
    \end{minipage}\begin{minipage}{0.245\linewidth}
        \centering
        \includegraphics[width=0.99\linewidth]{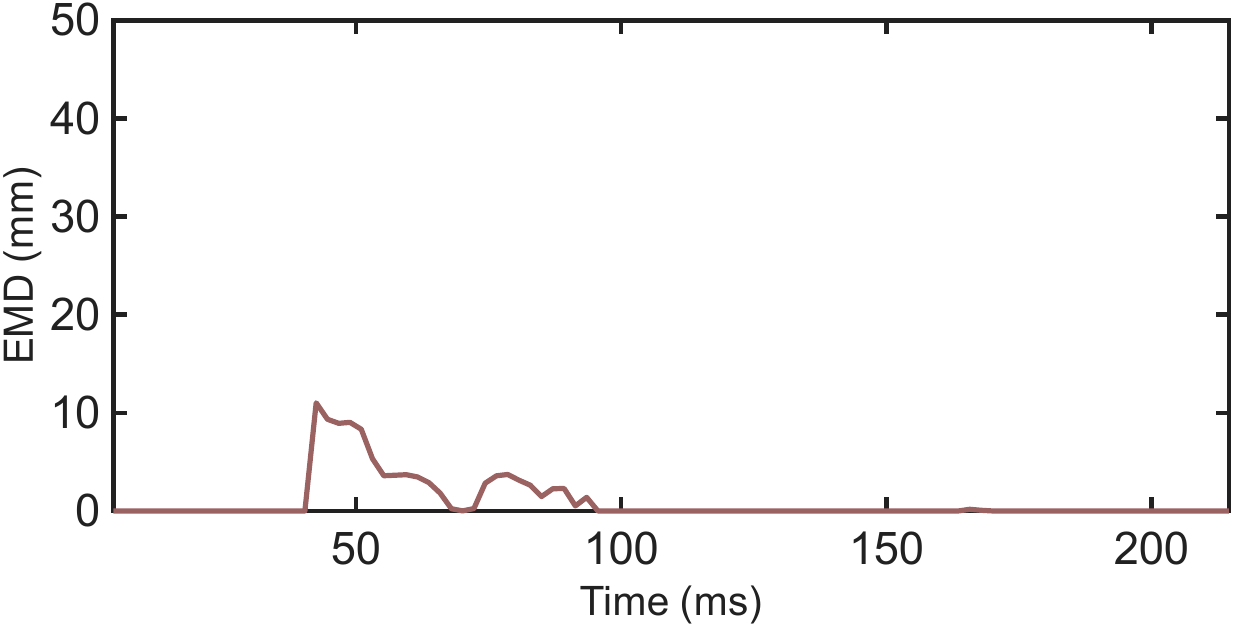}
    \end{minipage}\begin{minipage}{0.245\linewidth}
        \centering
        \includegraphics[width=0.99\linewidth]{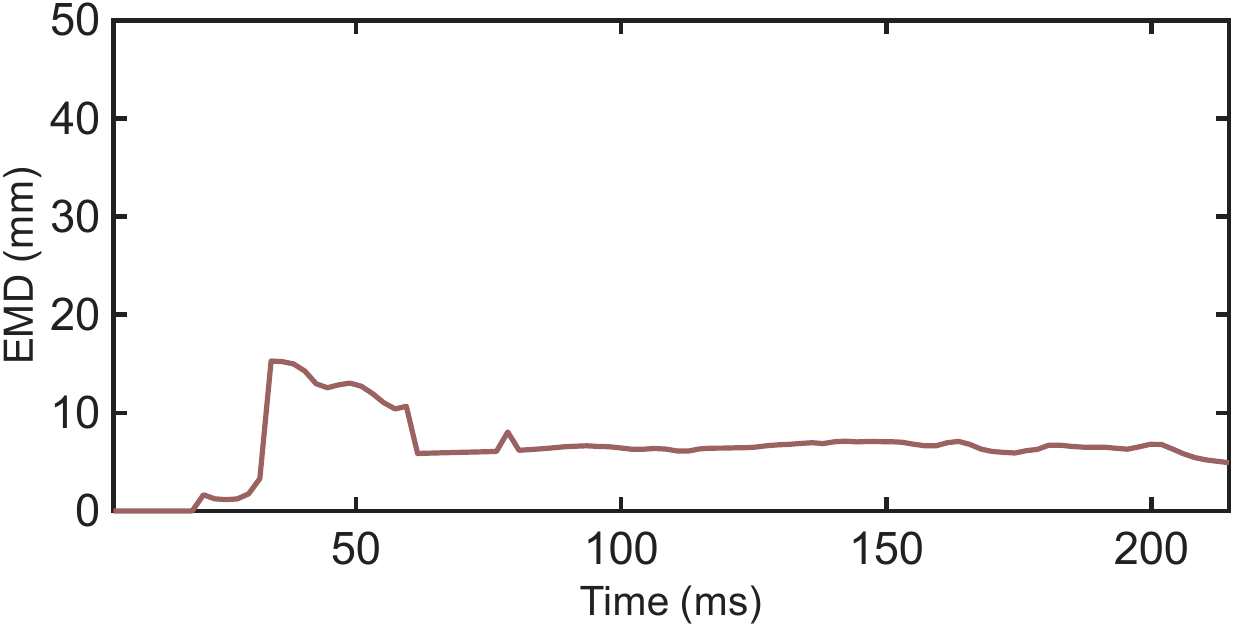}
    \end{minipage}

    \begin{minipage}{0.02\linewidth}
        \rotatebox{90}{\scriptsize S.P.}
    \end{minipage}\begin{minipage}{0.245\linewidth}
        \centering
        \includegraphics[width=0.99\linewidth]{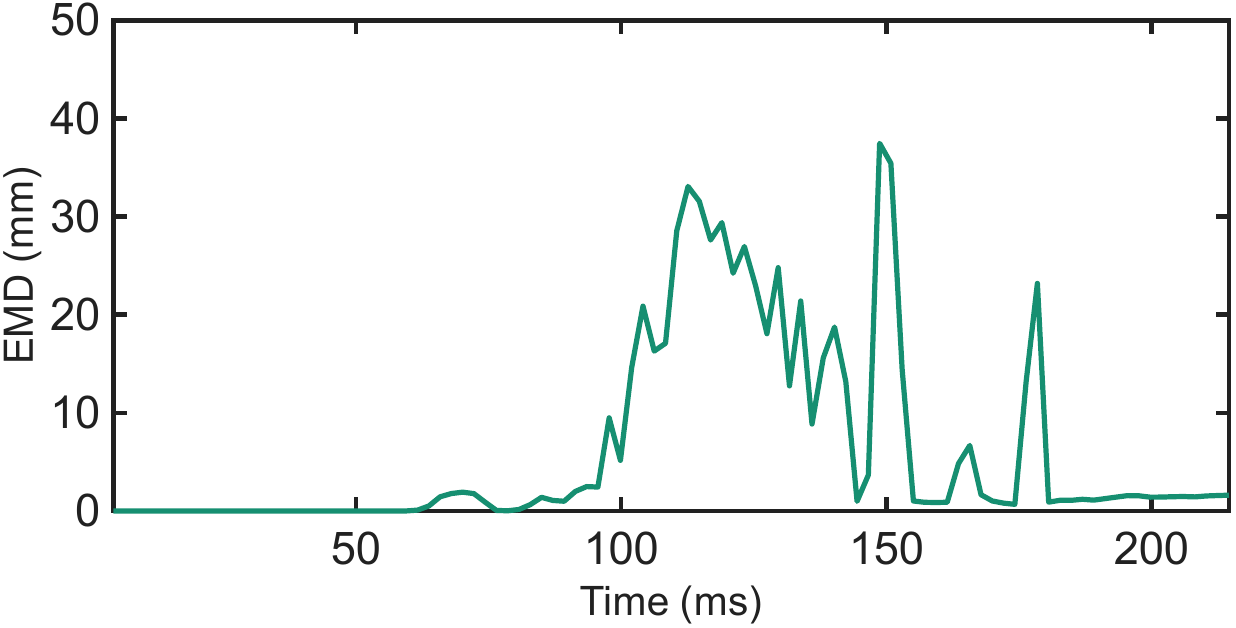}
    \end{minipage}\begin{minipage}{0.245\linewidth}
        \centering
        \includegraphics[width=0.99\linewidth]{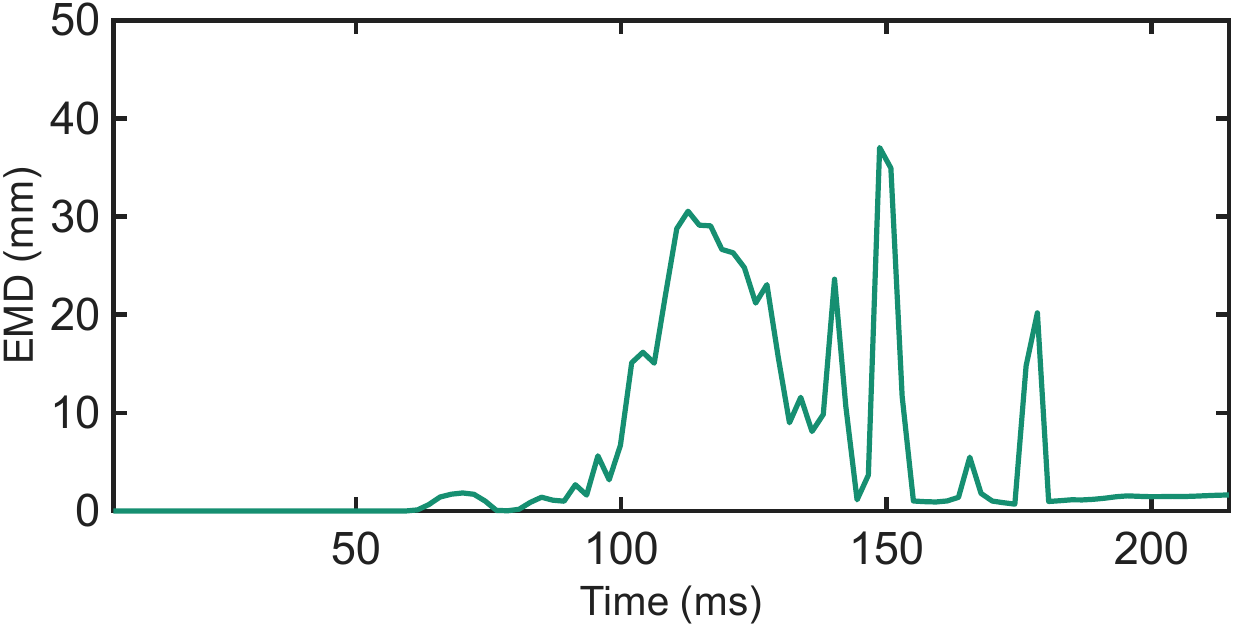}
    \end{minipage}\begin{minipage}{0.245\linewidth}
        \centering
        \includegraphics[width=0.99\linewidth]{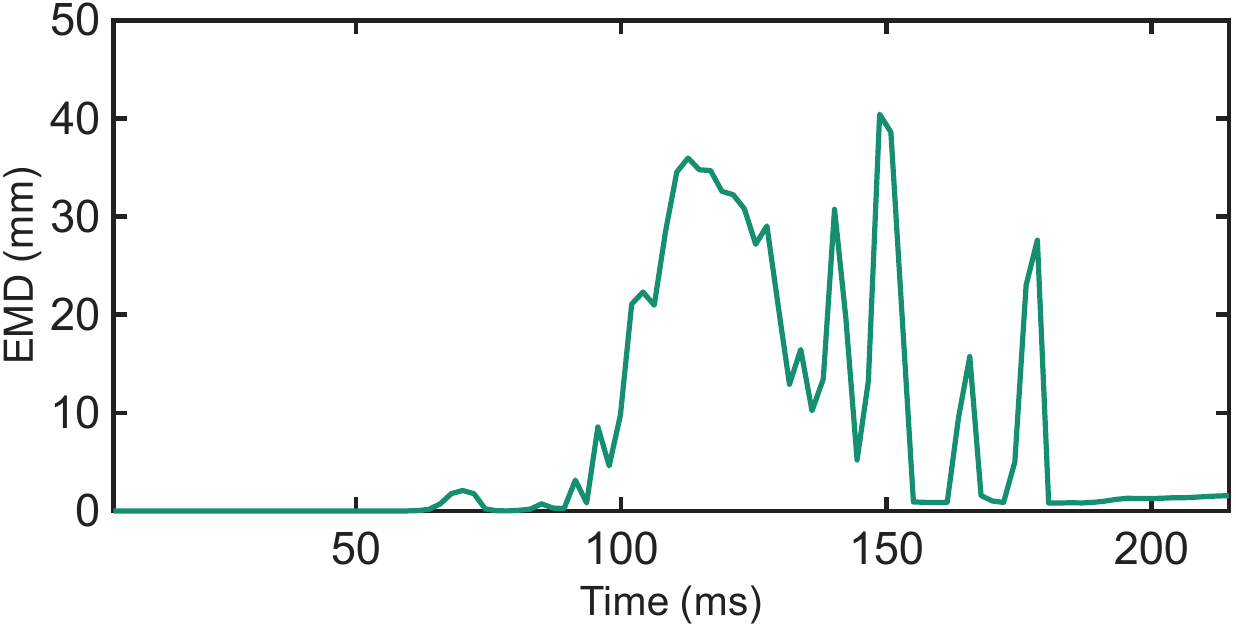}
    \end{minipage}\begin{minipage}{0.245\linewidth}
        \centering
        \includegraphics[width=0.99\linewidth]{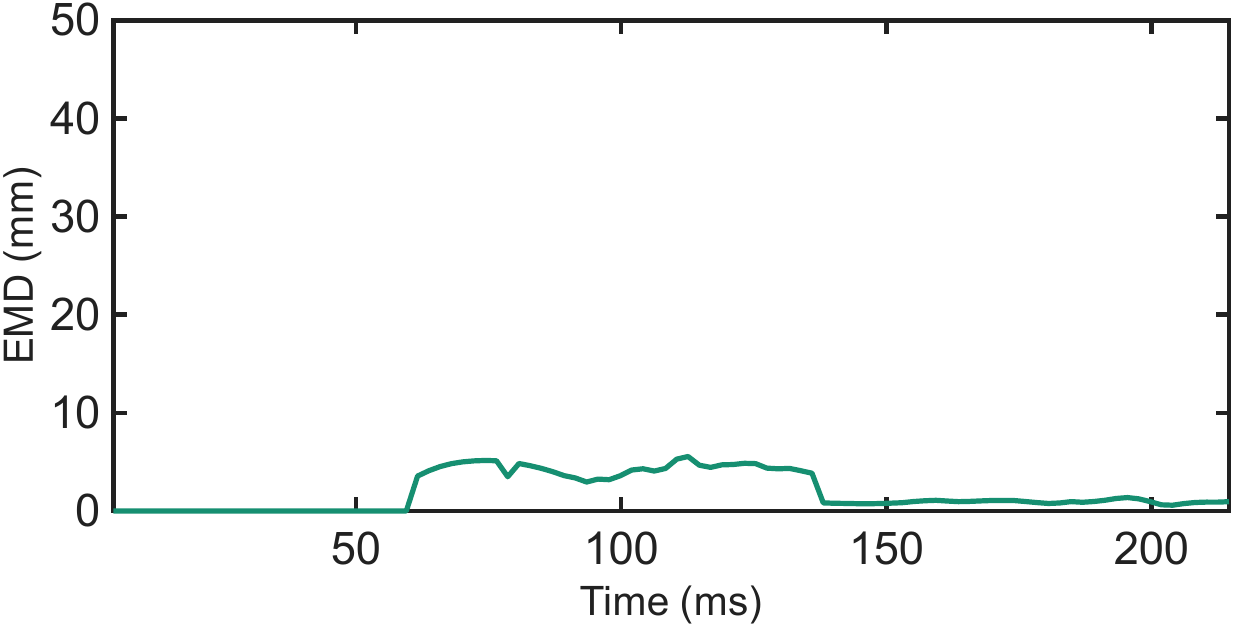}
    \end{minipage}

    \begin{minipage}{0.02\linewidth}
        \rotatebox{90}{\scriptsize C.M.F.}
    \end{minipage}\begin{minipage}{0.245\linewidth}
        \centering
        \includegraphics[width=0.99\linewidth]{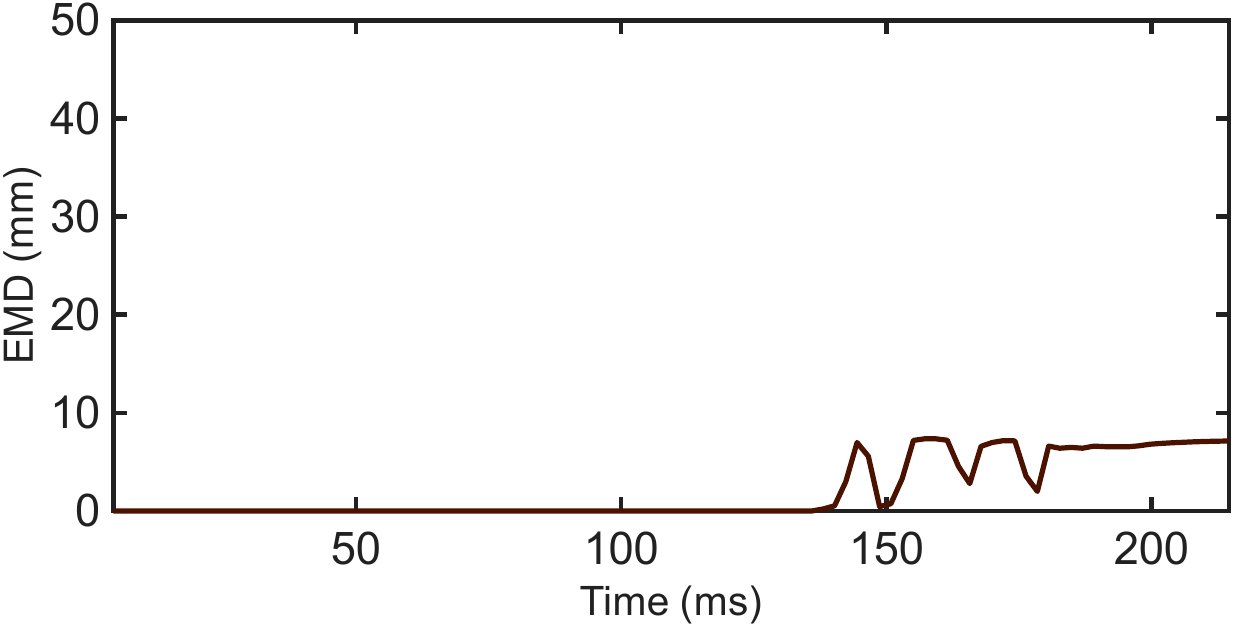}
    \end{minipage}\begin{minipage}{0.245\linewidth}
        \centering
        \includegraphics[width=0.99\linewidth]{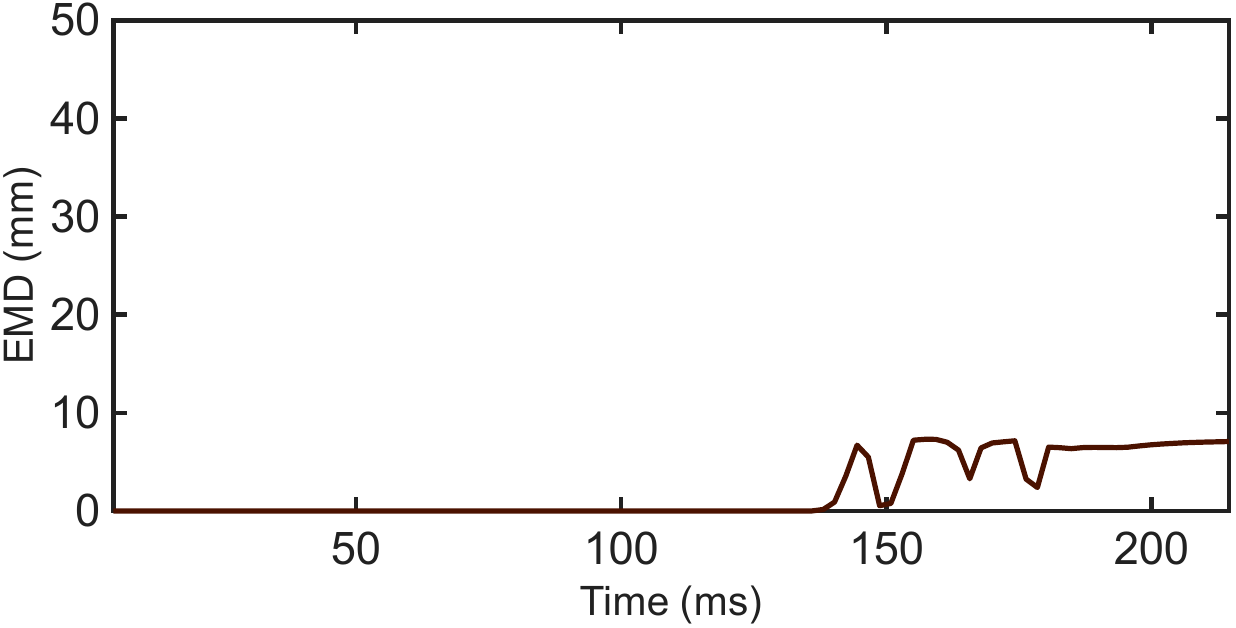}
    \end{minipage}\begin{minipage}{0.245\linewidth}
        \centering
        \includegraphics[width=0.99\linewidth]{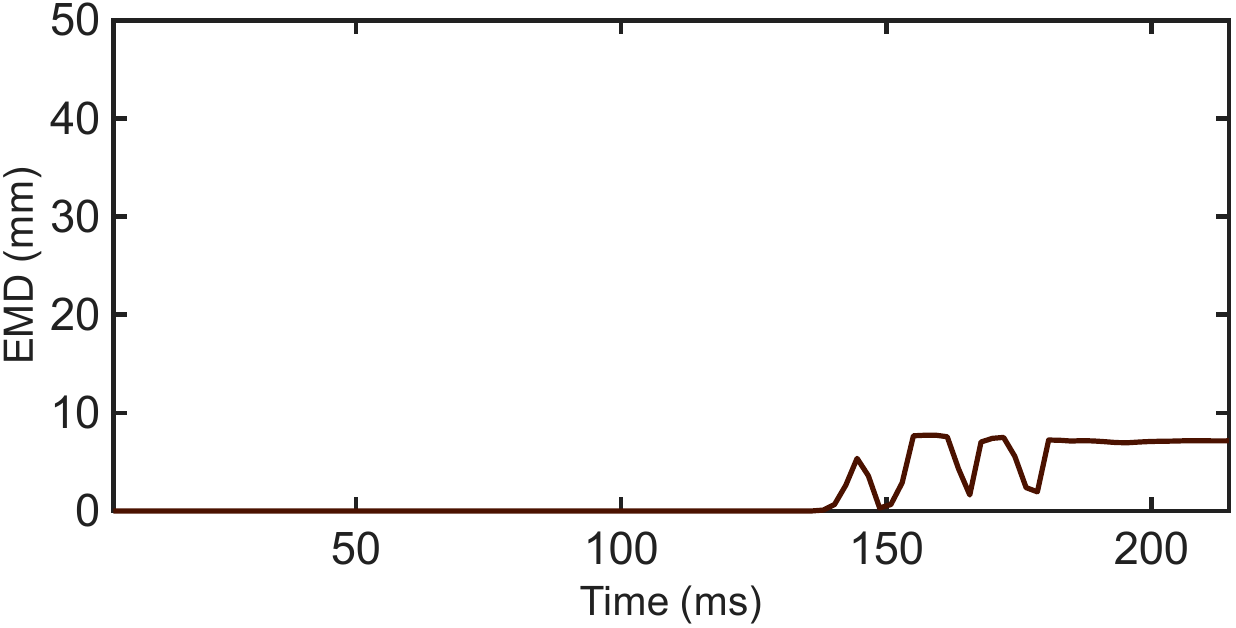}
    \end{minipage}\begin{minipage}{0.245\linewidth}
        \centering
        \includegraphics[width=0.99\linewidth]{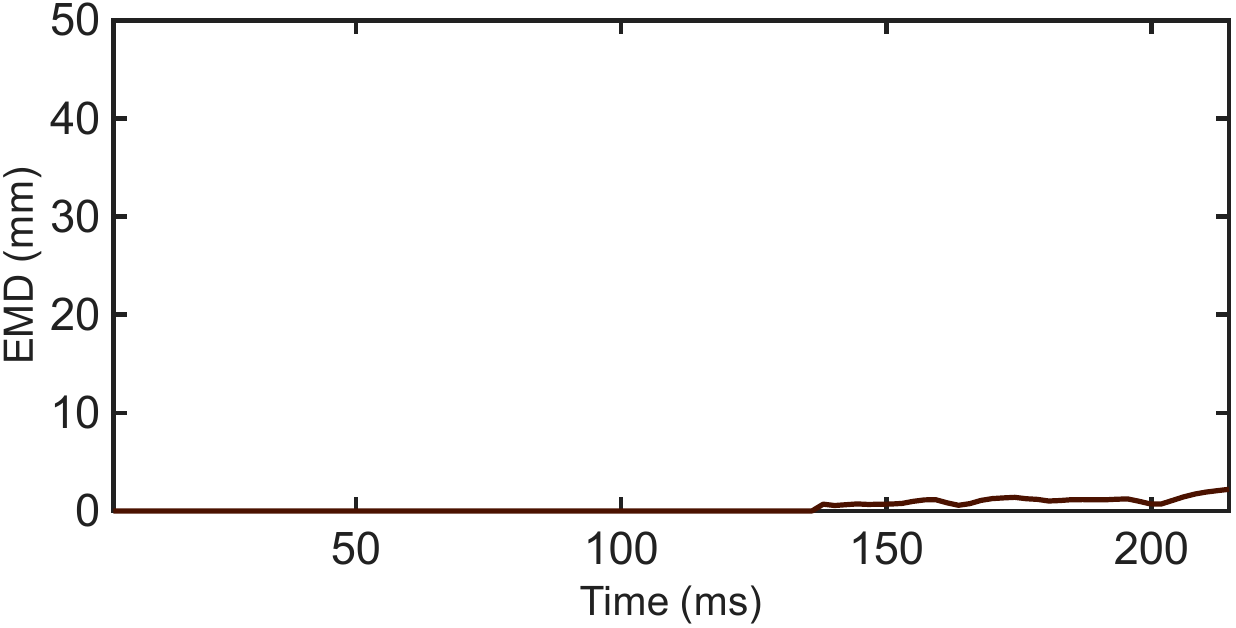}
    \end{minipage}

    \begin{minipage}{0.02\linewidth}
        \rotatebox{90}{\scriptsize Total}
    \end{minipage}\begin{minipage}{0.245\linewidth}
        \centering
        \includegraphics[width=0.99\linewidth]{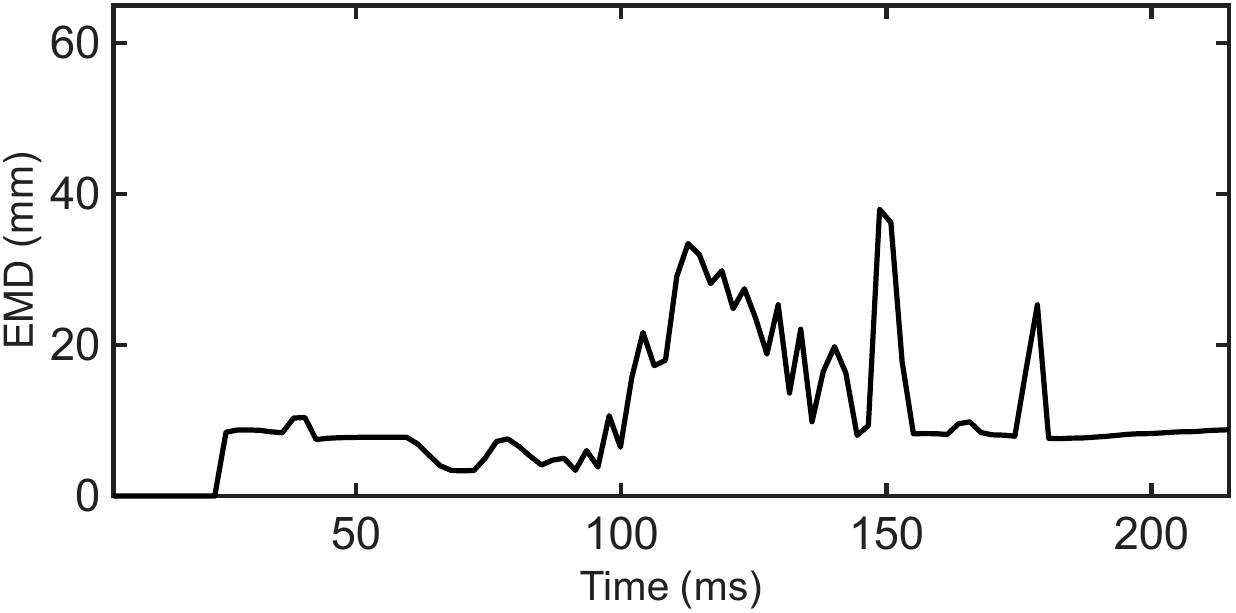}
    \end{minipage}\begin{minipage}{0.245\linewidth}
        \centering
        \includegraphics[width=0.99\linewidth]{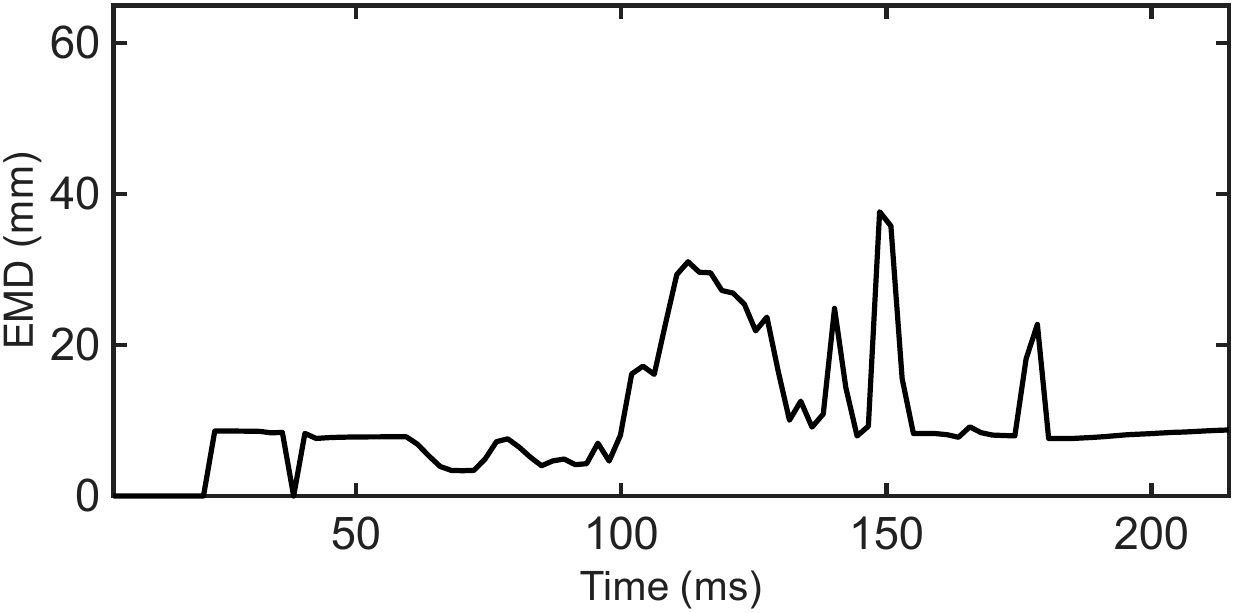}
    \end{minipage}\begin{minipage}{0.245\linewidth}
        \centering
        \includegraphics[width=0.99\linewidth]{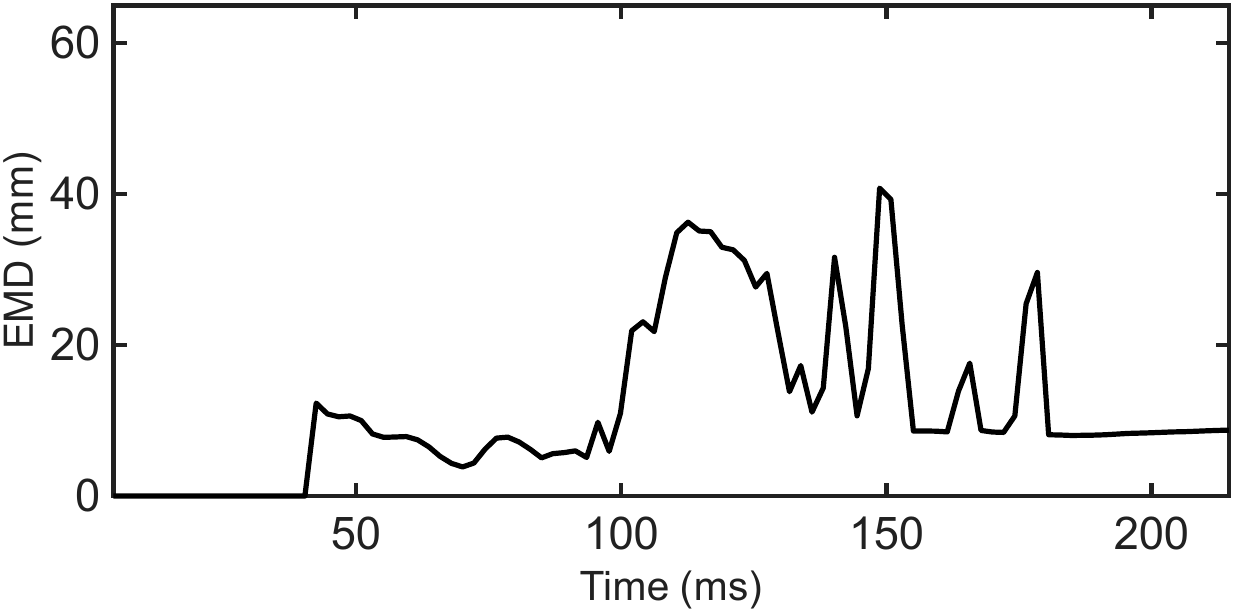}
    \end{minipage}\begin{minipage}{0.245\linewidth}
        \centering
        \includegraphics[width=0.99\linewidth]{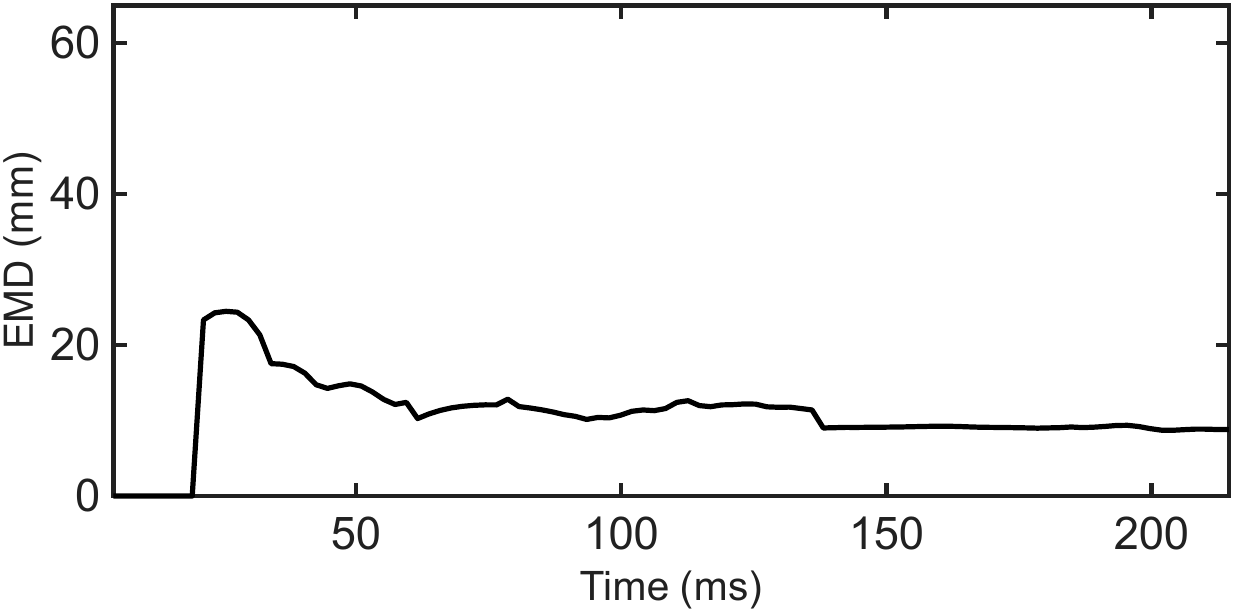}
    \end{minipage}

    \caption{Earth Mover's Distance (EMD) between the reconstructed and simulated ground-truth source distributions under the noiseless condition. Columns correspond to CrSKF, SKF, KF, and LCMV; rows show parcel-wise EMD for the lateral occipital (L.O.), lingual (Ling.), superior parietal (S.P.), and caudal middle frontal (C.M.F.) regions, together with the total EMD. Lower values indicate closer spatial agreement.}
    \label{fig:EMD0dB}
\end{figure}
\begin{figure}[h!]
    \centering
    \begin{minipage}{0.02\linewidth}
        \rotatebox{90}{\scriptsize \hfill}
    \end{minipage}\begin{minipage}{0.245\linewidth}
        \centering
        {\bf CrSKF}
    \end{minipage}\begin{minipage}{0.245\linewidth}
        \centering
        {\bf SKF}
    \end{minipage}\begin{minipage}{0.245\linewidth}
        \centering
        {\bf KF}
    \end{minipage}\begin{minipage}{0.245\linewidth}
        \centering
        {\bf LCMV}
    \end{minipage}

    \begin{minipage}{0.02\linewidth}
        \rotatebox{90}{\scriptsize L.O.}
    \end{minipage}\begin{minipage}{0.245\linewidth}
        \centering
        \includegraphics[width=0.99\linewidth]{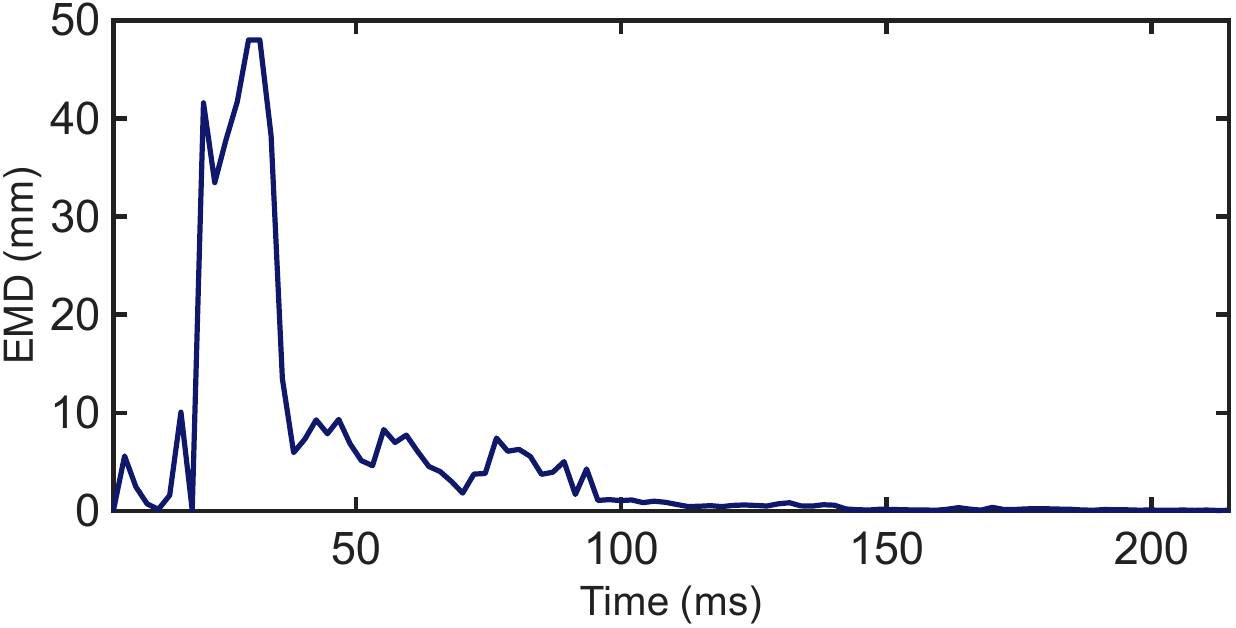}
    \end{minipage}\begin{minipage}{0.245\linewidth}
        \centering
        \includegraphics[width=0.99\linewidth]{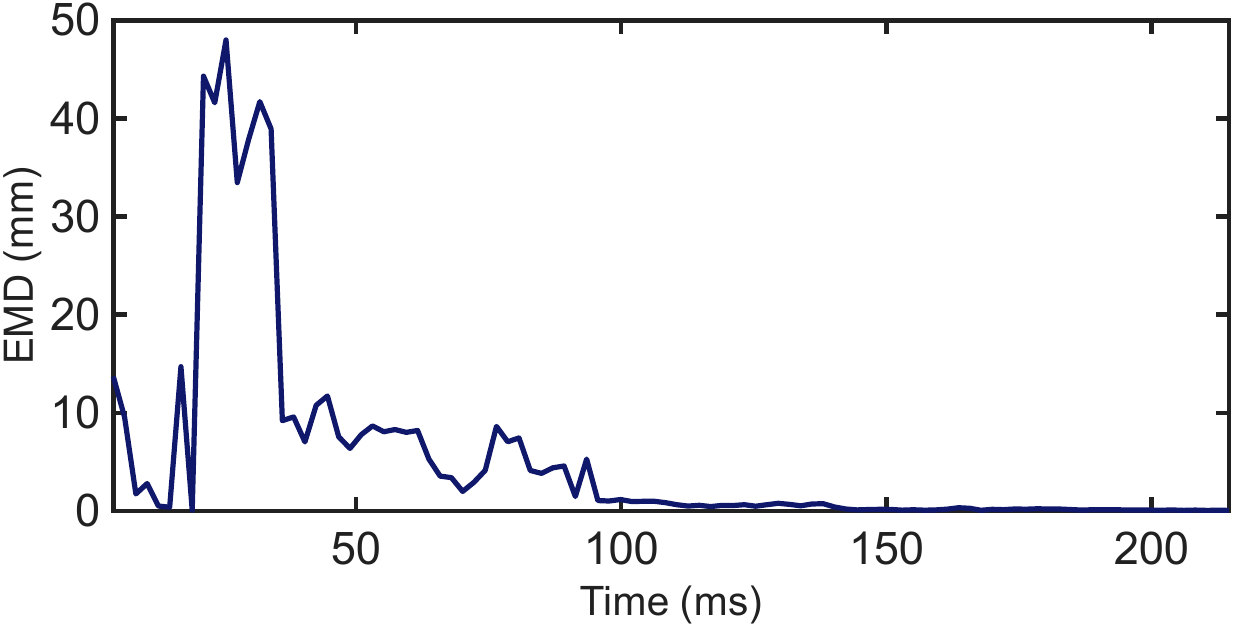}
    \end{minipage}\begin{minipage}{0.245\linewidth}
        \centering
        \includegraphics[width=0.99\linewidth]{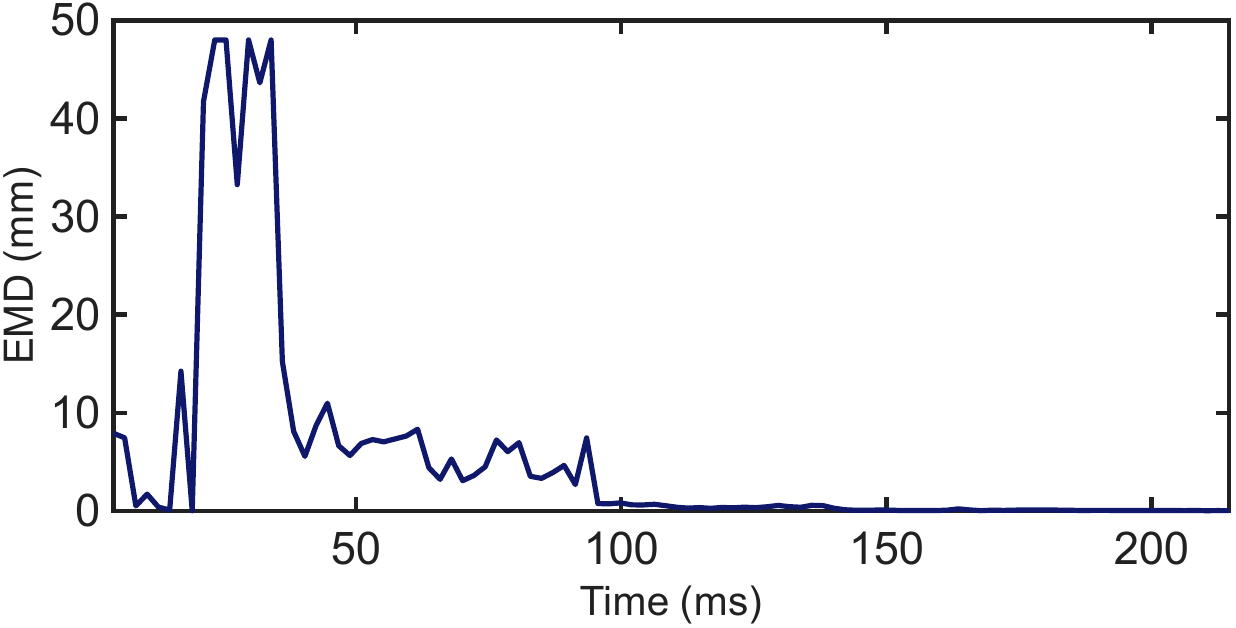}
    \end{minipage}\begin{minipage}{0.245\linewidth}
        \centering
        \includegraphics[width=0.99\linewidth]{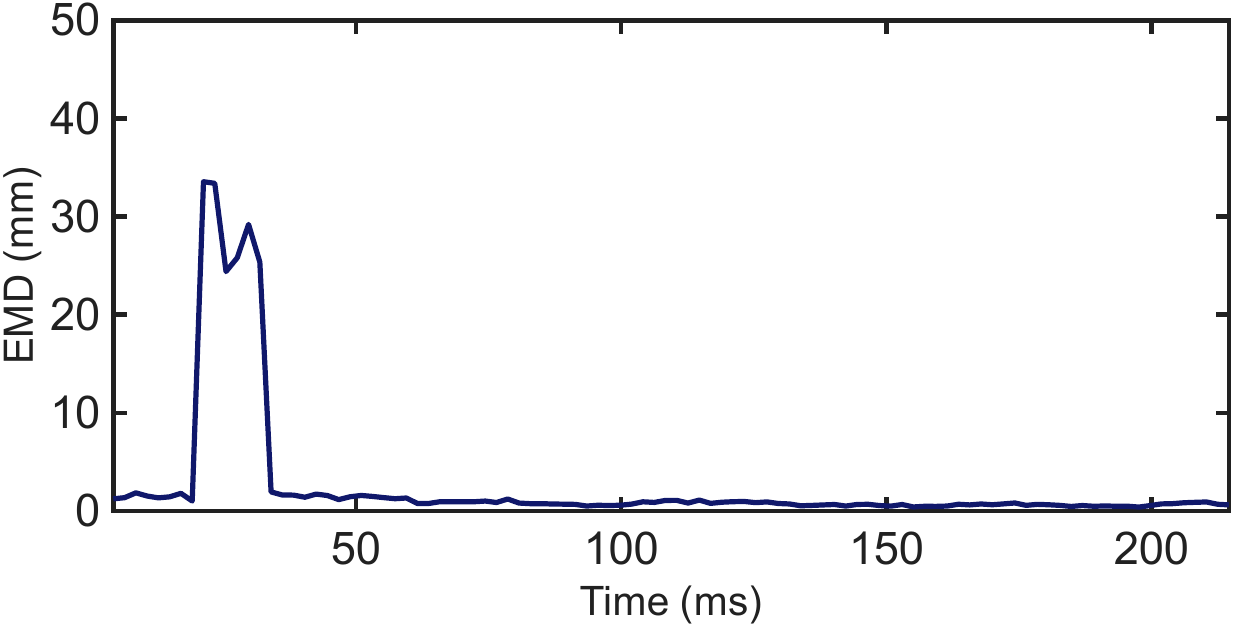}
    \end{minipage}

    \begin{minipage}{0.02\linewidth}
        \rotatebox{90}{\scriptsize Ling.}
    \end{minipage}\begin{minipage}{0.245\linewidth}
        \centering
        \includegraphics[width=0.99\linewidth]{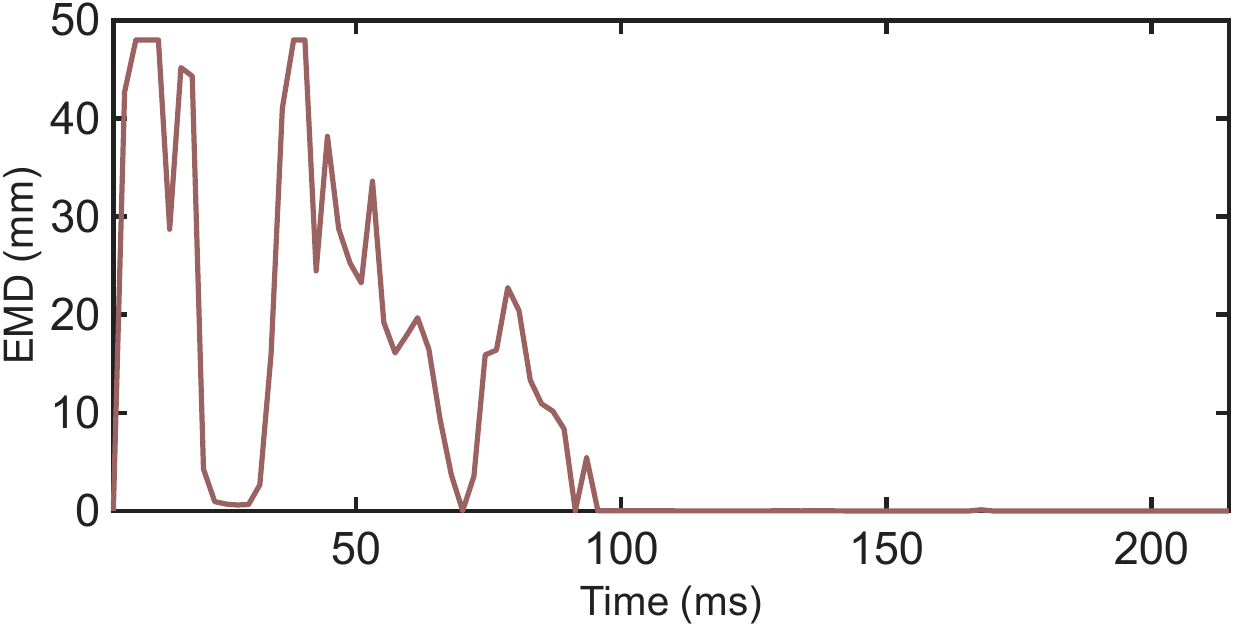}
    \end{minipage}\begin{minipage}{0.245\linewidth}
        \centering
        \includegraphics[width=0.99\linewidth]{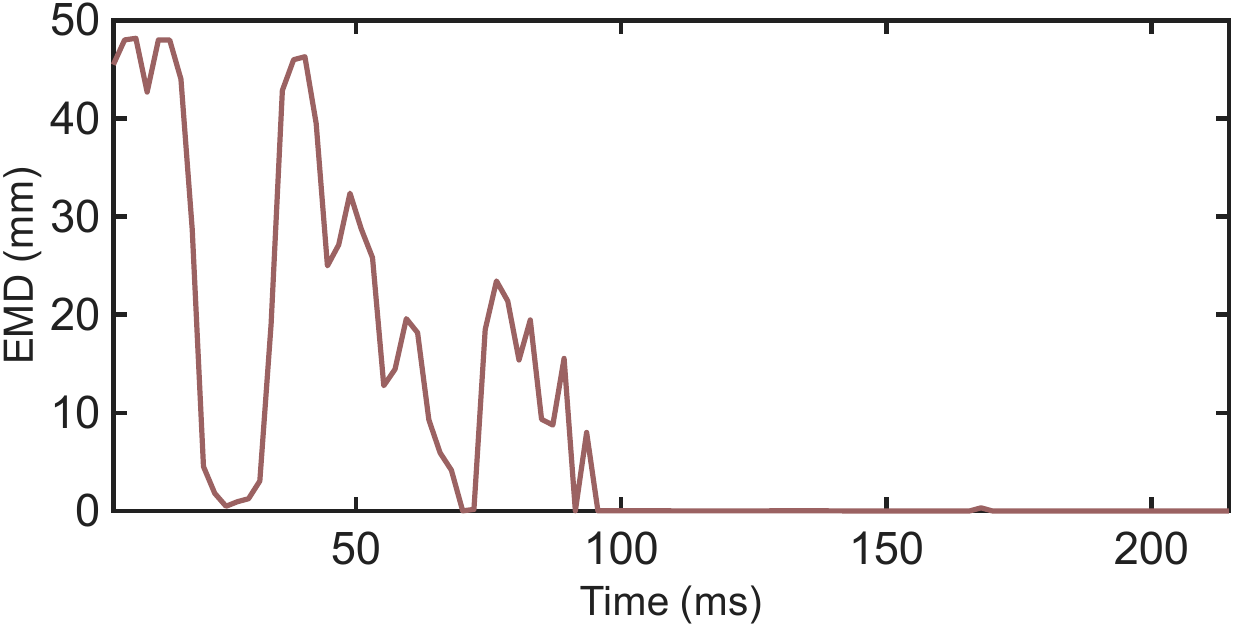}
    \end{minipage}\begin{minipage}{0.245\linewidth}
        \centering
        \includegraphics[width=0.99\linewidth]{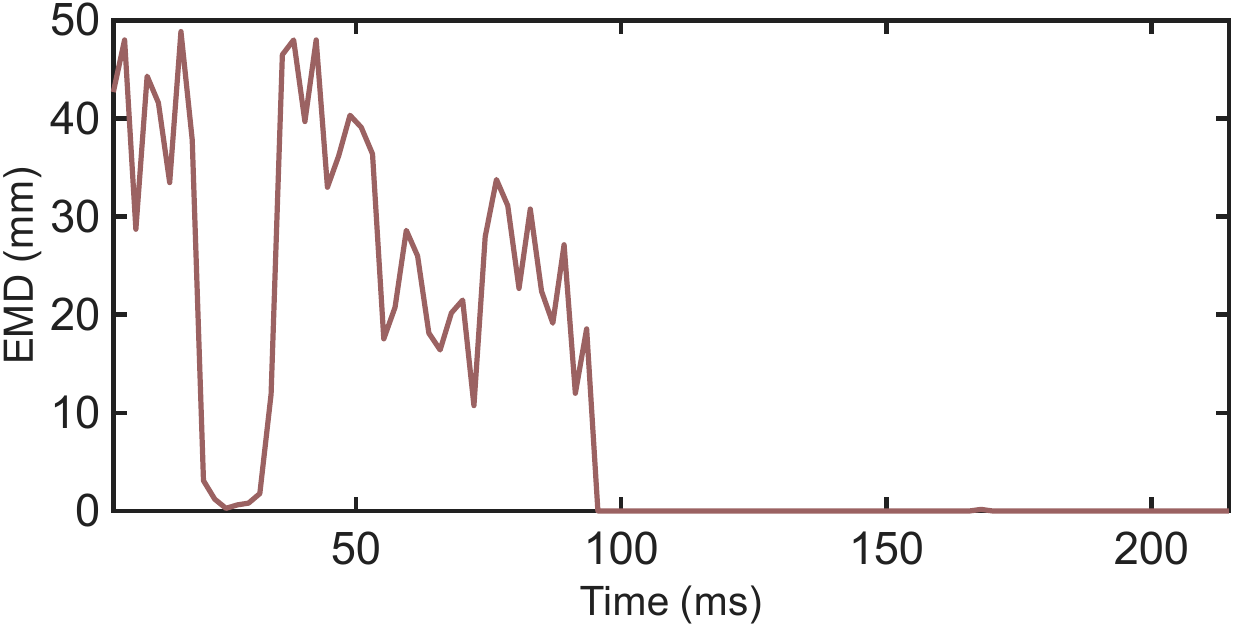}
    \end{minipage}\begin{minipage}{0.245\linewidth}
        \centering
        \includegraphics[width=0.99\linewidth]{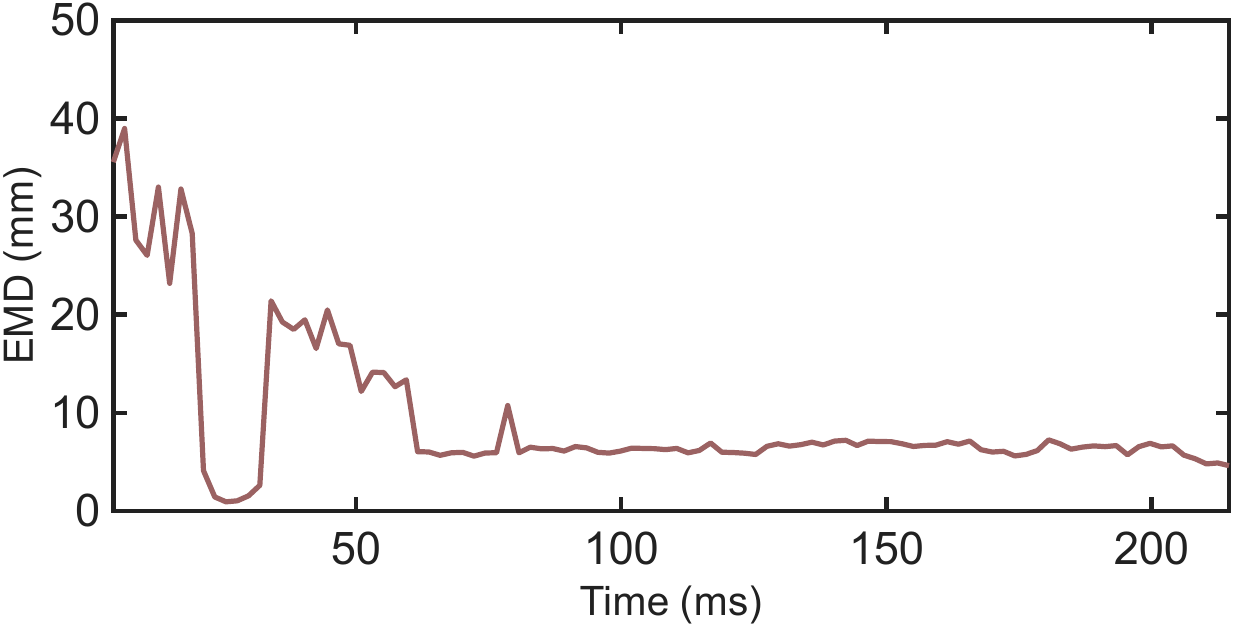}
    \end{minipage}

    \begin{minipage}{0.02\linewidth}
        \rotatebox{90}{\scriptsize S.P.}
    \end{minipage}\begin{minipage}{0.245\linewidth}
        \centering
        \includegraphics[width=0.99\linewidth]{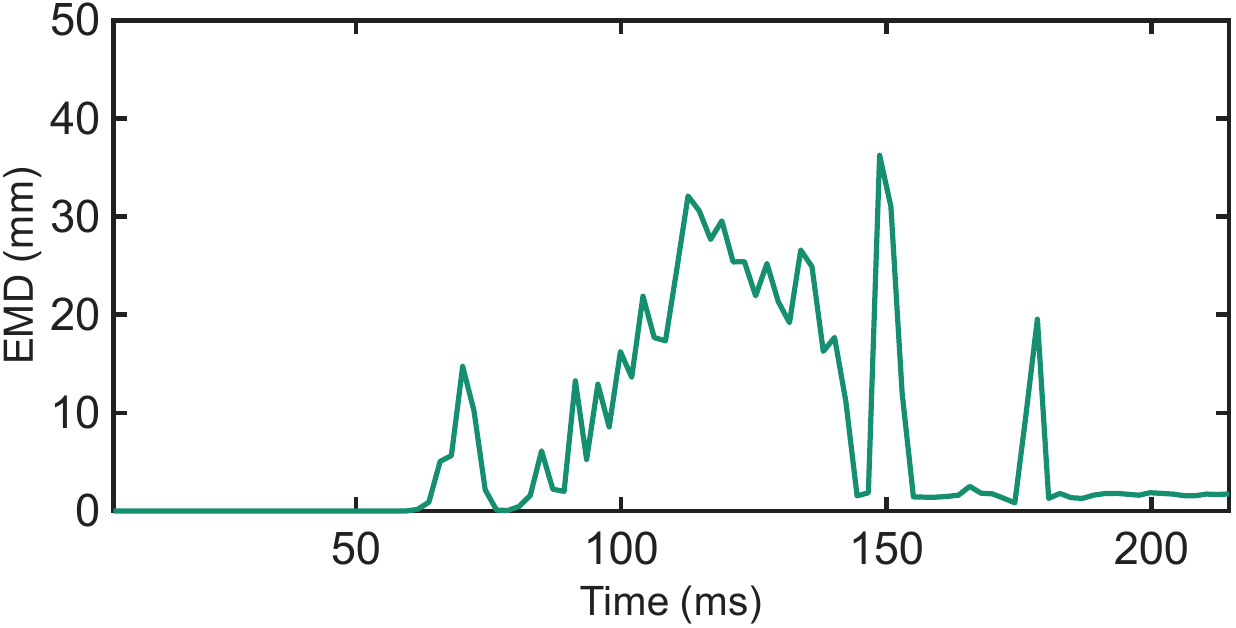}
    \end{minipage}\begin{minipage}{0.245\linewidth}
        \centering
        \includegraphics[width=0.99\linewidth]{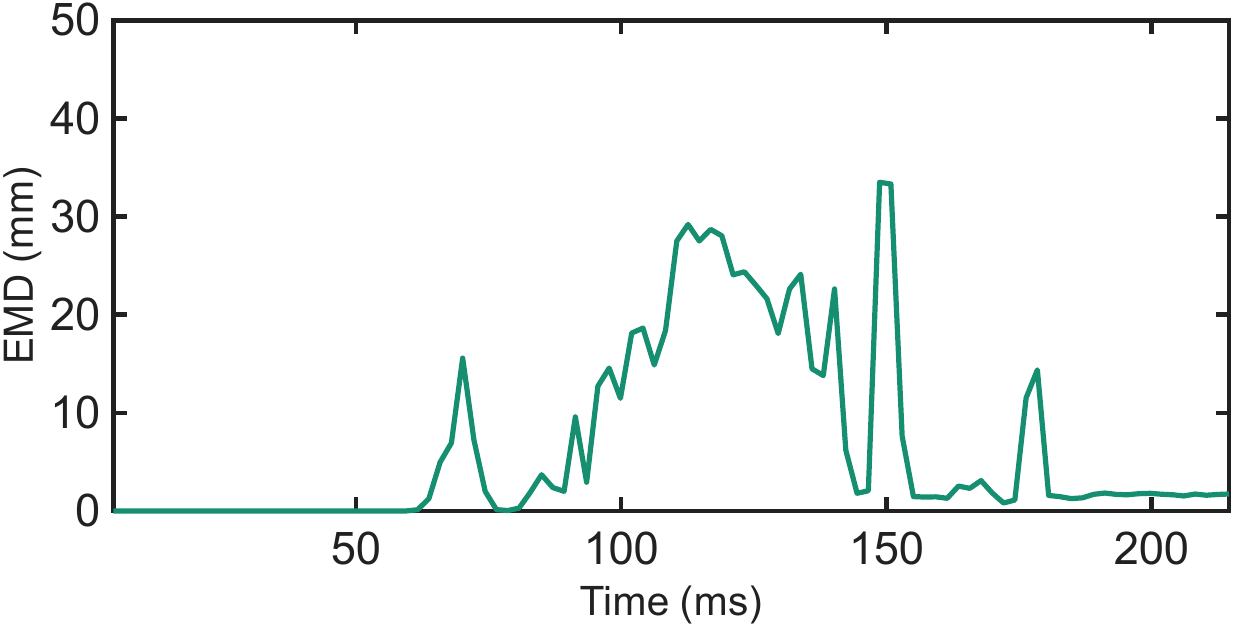}
    \end{minipage}\begin{minipage}{0.245\linewidth}
        \centering
        \includegraphics[width=0.99\linewidth]{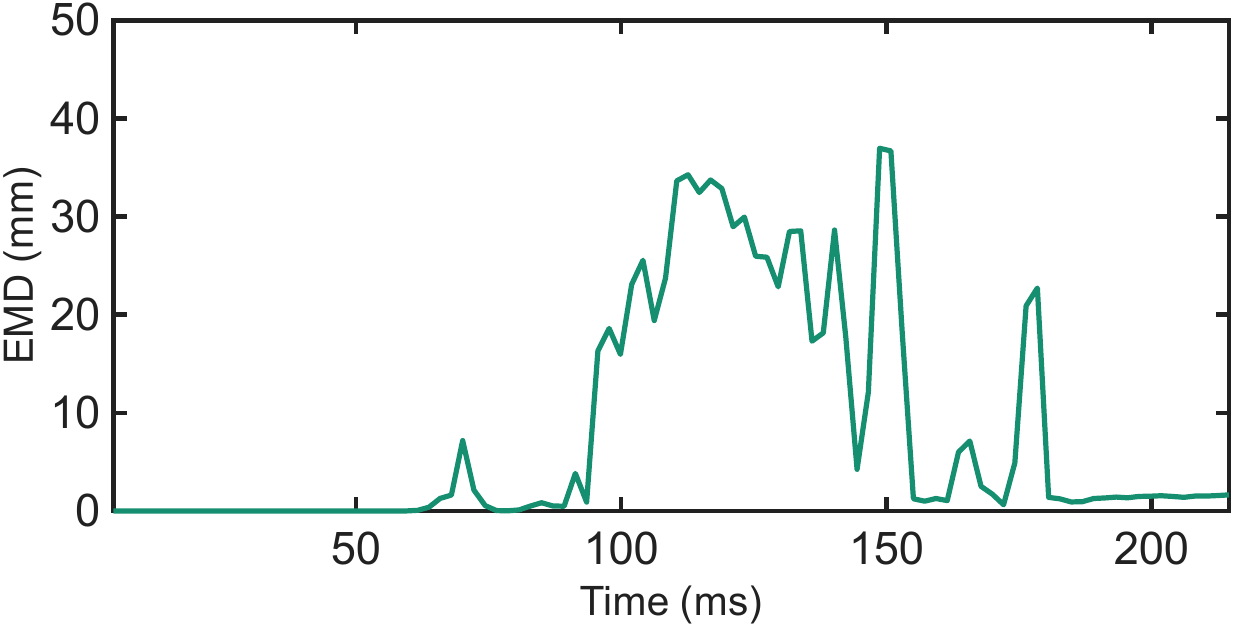}
    \end{minipage}\begin{minipage}{0.245\linewidth}
        \centering
        \includegraphics[width=0.99\linewidth]{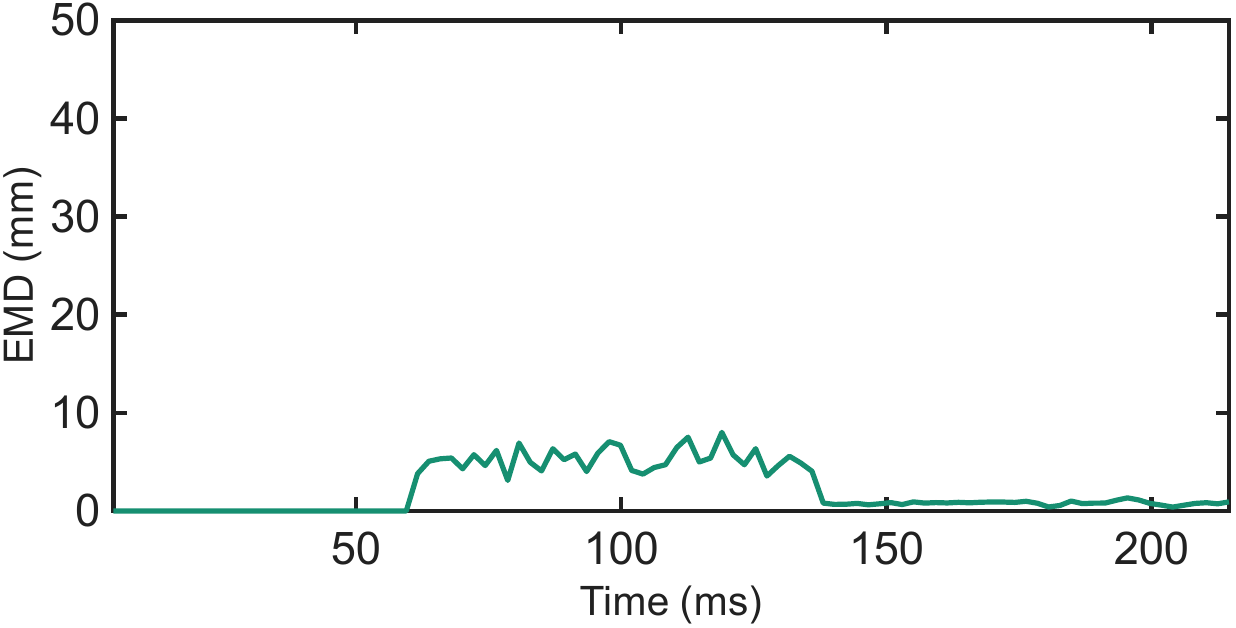}
    \end{minipage}

    \begin{minipage}{0.02\linewidth}
        \rotatebox{90}{\scriptsize C.M.F.}
    \end{minipage}\begin{minipage}{0.245\linewidth}
        \centering
        \includegraphics[width=0.99\linewidth]{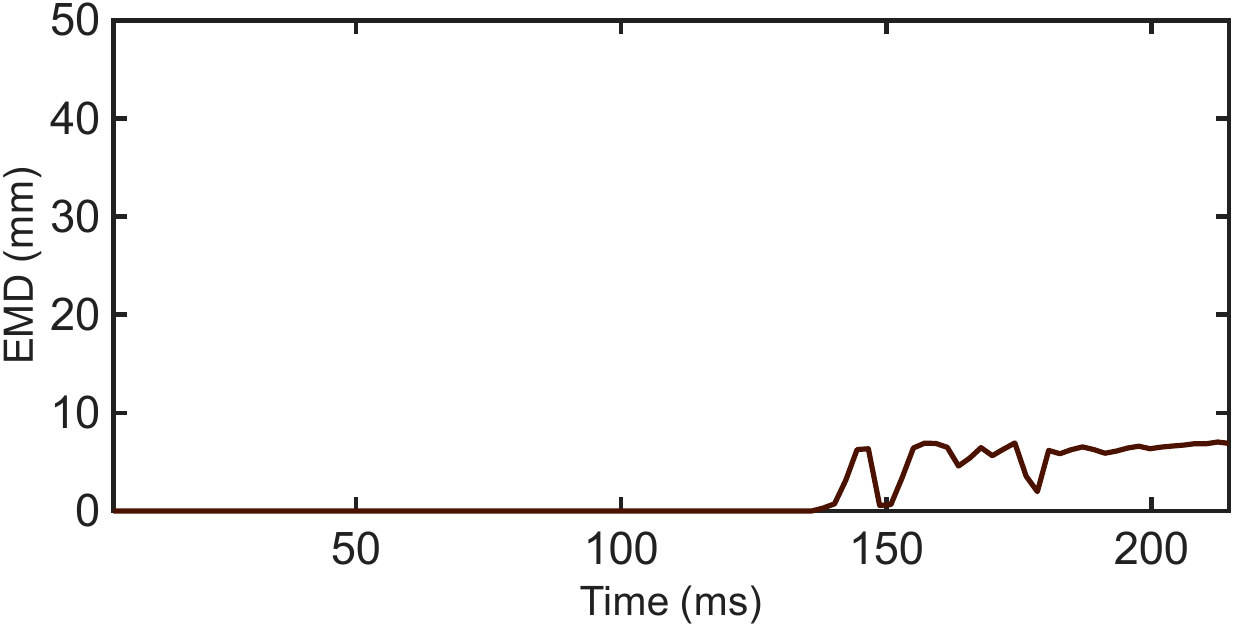}
    \end{minipage}\begin{minipage}{0.245\linewidth}
        \centering
        \includegraphics[width=0.99\linewidth]{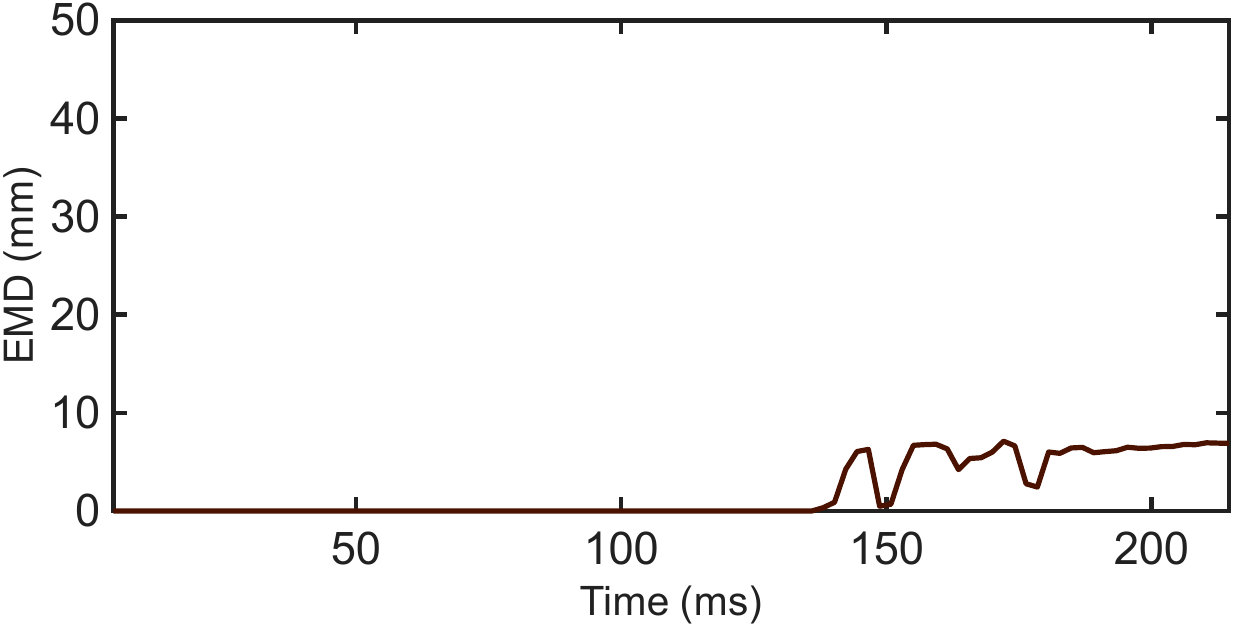}
    \end{minipage}\begin{minipage}{0.245\linewidth}
        \centering
        \includegraphics[width=0.99\linewidth]{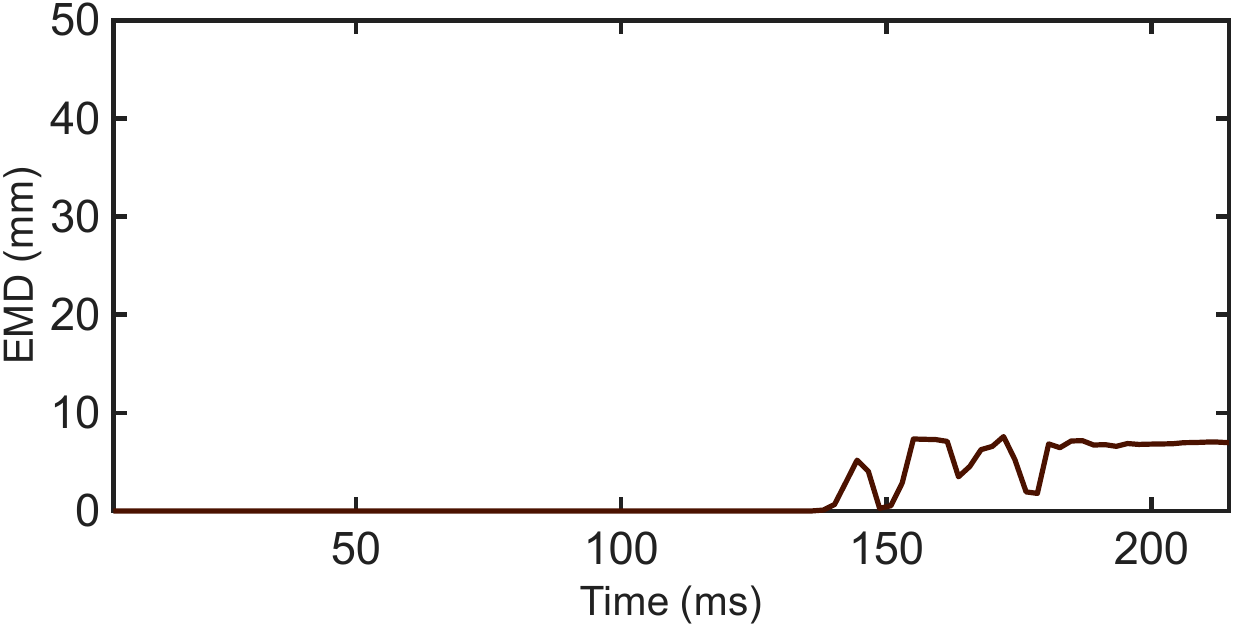}
    \end{minipage}\begin{minipage}{0.245\linewidth}
        \centering
        \includegraphics[width=0.99\linewidth]{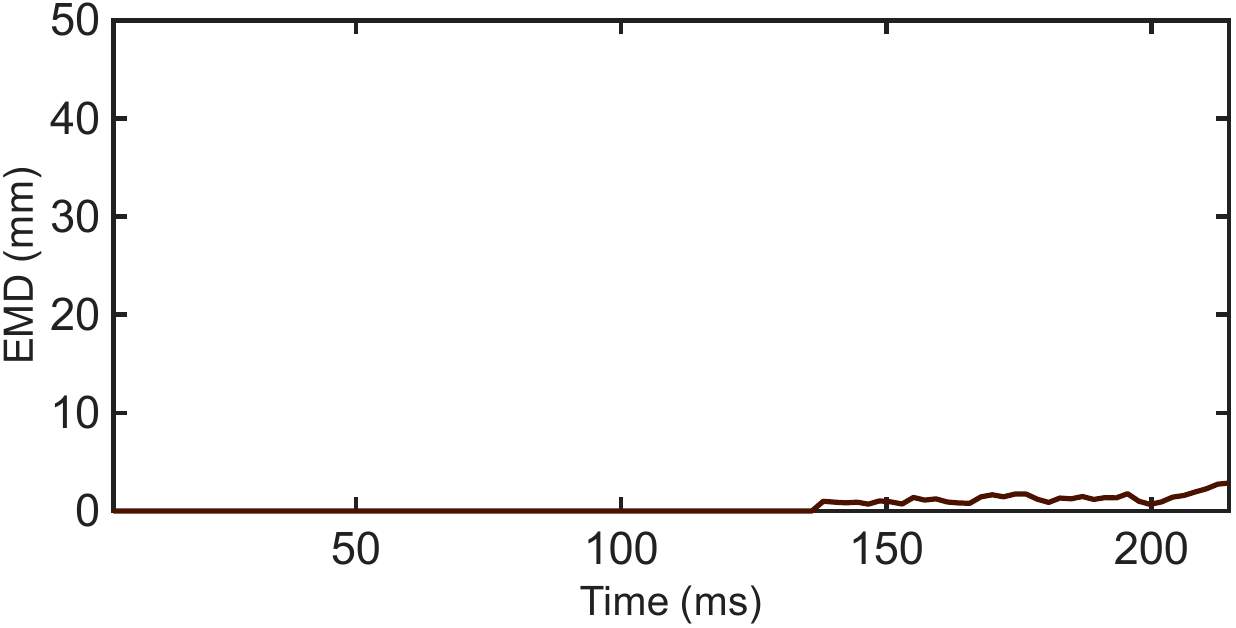}
    \end{minipage}

    \begin{minipage}{0.02\linewidth}
        \rotatebox{90}{\scriptsize Total}
    \end{minipage}\begin{minipage}{0.245\linewidth}
        \centering
        \includegraphics[width=0.99\linewidth]{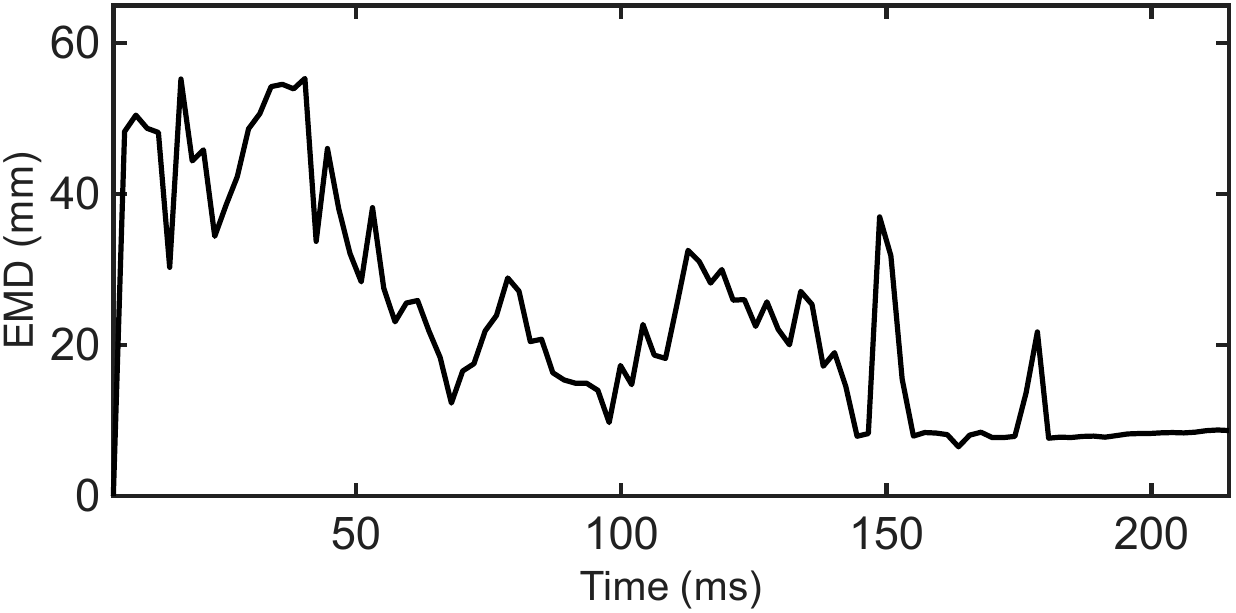}
    \end{minipage}\begin{minipage}{0.245\linewidth}
        \centering
        \includegraphics[width=0.99\linewidth]{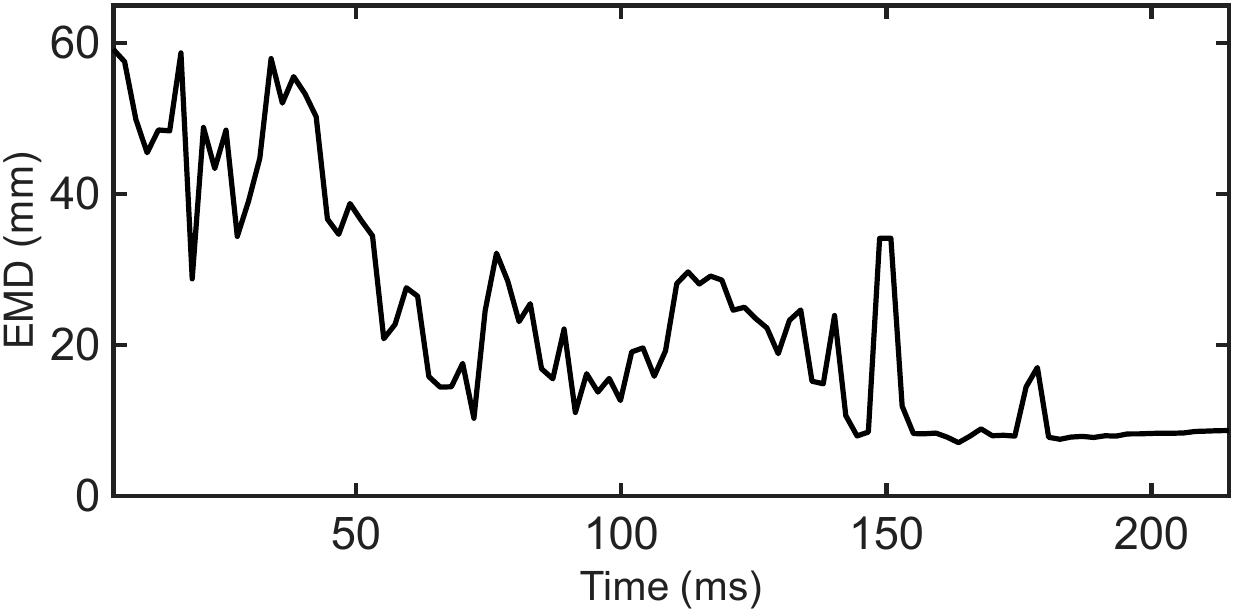}
    \end{minipage}\begin{minipage}{0.245\linewidth}
        \centering
        \includegraphics[width=0.99\linewidth]{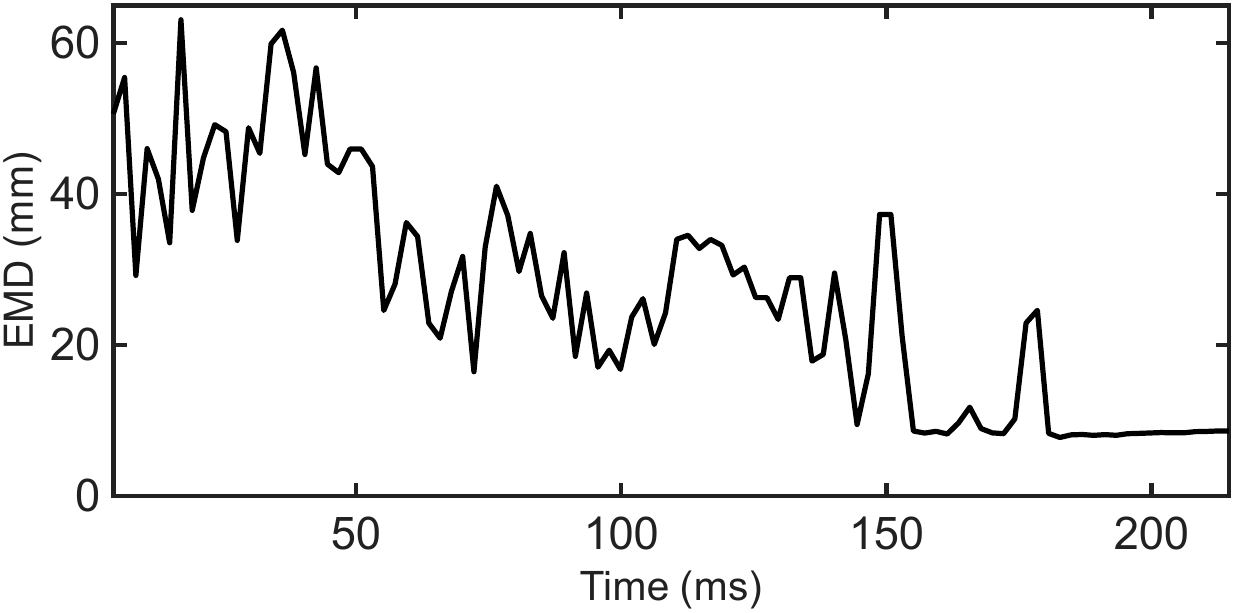}
    \end{minipage}\begin{minipage}{0.245\linewidth}
        \centering
        \includegraphics[width=0.99\linewidth]{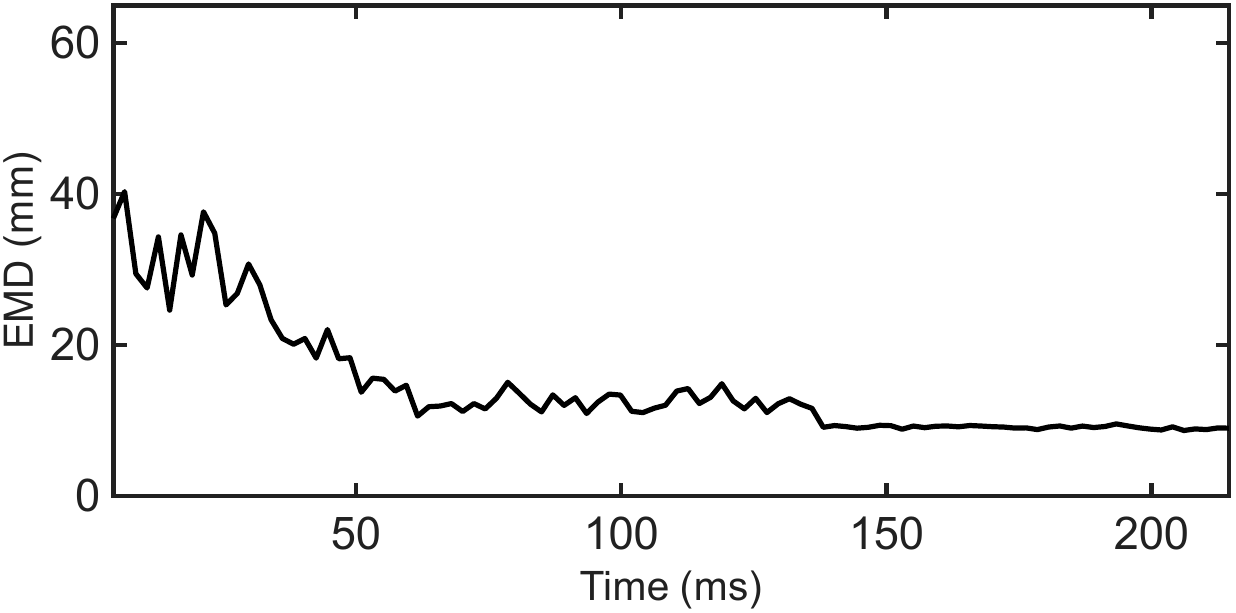}
    \end{minipage}

    \caption{Earth Mover's Distance (EMD) between the reconstructed and simulated ground-truth source distributions at 25 dB SNR. Columns correspond to CrSKF, SKF, KF, and LCMV; rows show parcel-wise EMD for the lateral occipital (L.O.), lingual (Ling.), superior parietal (S.P.), and caudal middle frontal (C.M.F.) regions, together with the total EMD. Lower values indicate closer spatial agreement.}
    \label{fig:EMD25dB}
\end{figure}
\begin{figure}[h!]
    \centering
    \begin{minipage}{0.02\linewidth}
        \rotatebox{90}{\scriptsize \hfill}
    \end{minipage}\begin{minipage}{0.245\linewidth}
        \centering
        {\bf CrSKF}
    \end{minipage}\begin{minipage}{0.245\linewidth}
        \centering
        {\bf SKF}
    \end{minipage}\begin{minipage}{0.245\linewidth}
        \centering
        {\bf KF}
    \end{minipage}\begin{minipage}{0.245\linewidth}
        \centering
        {\bf LCMV}
    \end{minipage}

    \begin{minipage}{0.02\linewidth}
        \rotatebox{90}{\scriptsize L.O.}
    \end{minipage}\begin{minipage}{0.245\linewidth}
        \centering
        \includegraphics[width=0.99\linewidth]{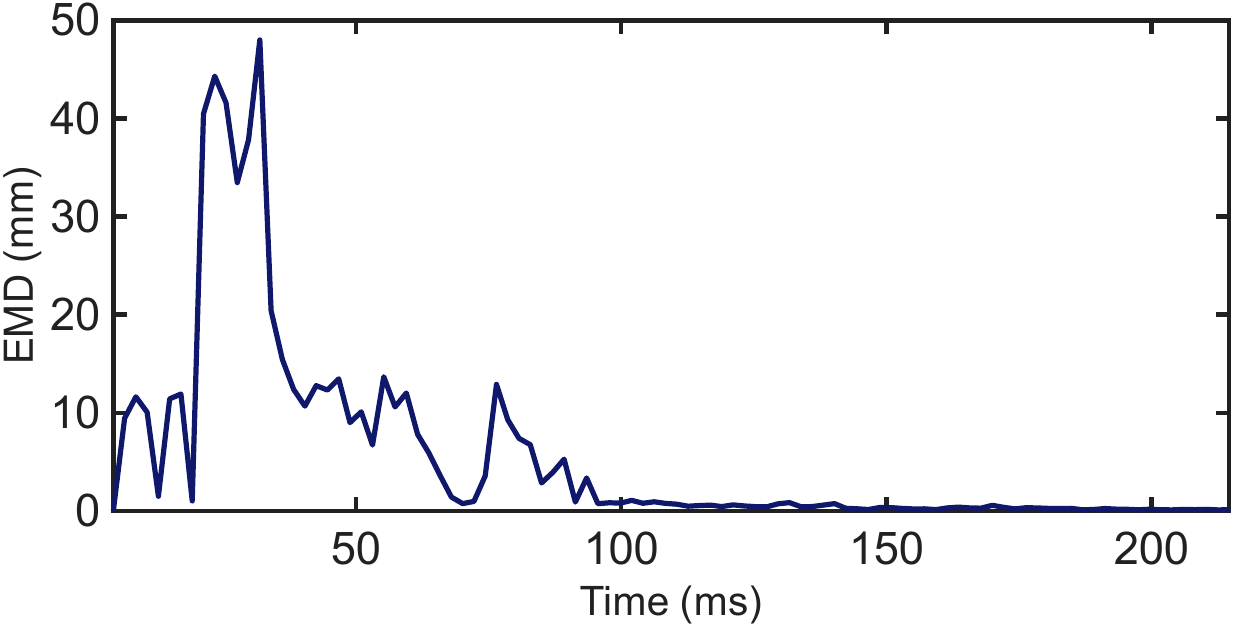}
    \end{minipage}\begin{minipage}{0.245\linewidth}
        \centering
        \includegraphics[width=0.99\linewidth]{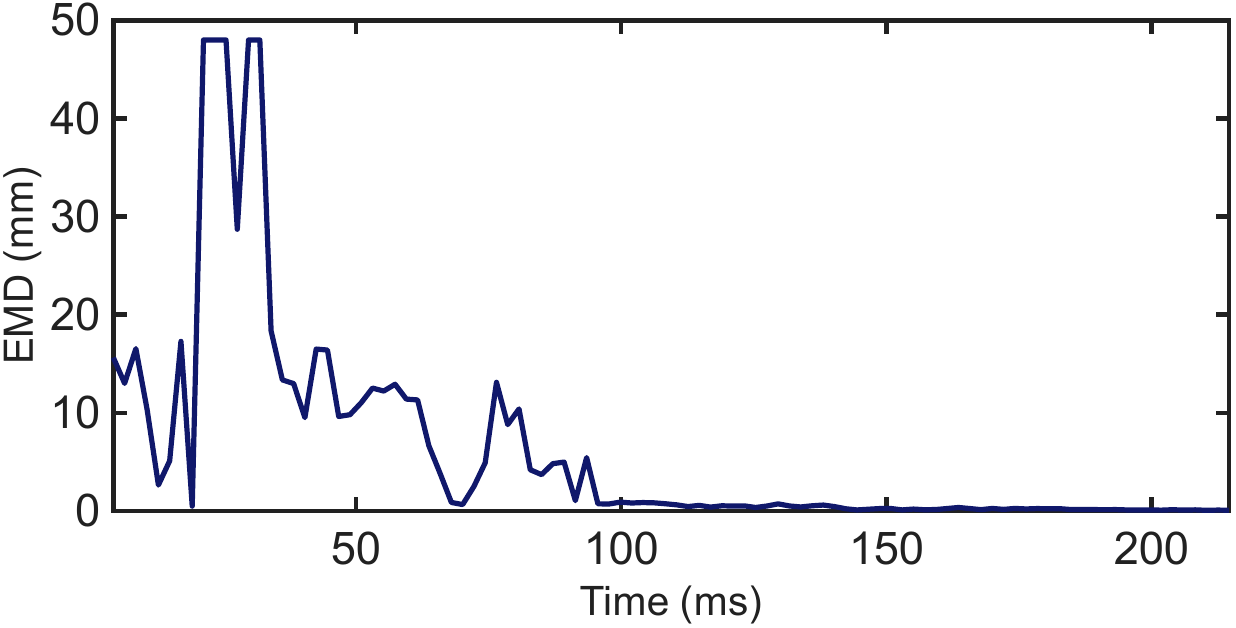}
    \end{minipage}\begin{minipage}{0.245\linewidth}
        \centering
        \includegraphics[width=0.99\linewidth]{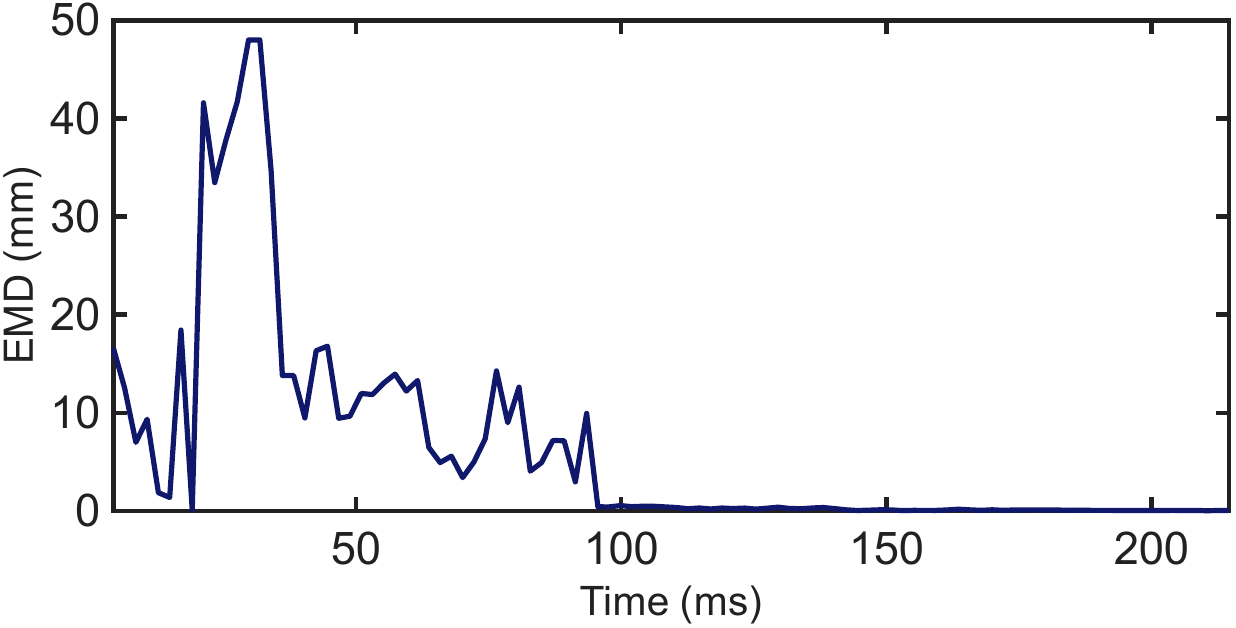}
    \end{minipage}\begin{minipage}{0.245\linewidth}
        \centering
        \includegraphics[width=0.99\linewidth]{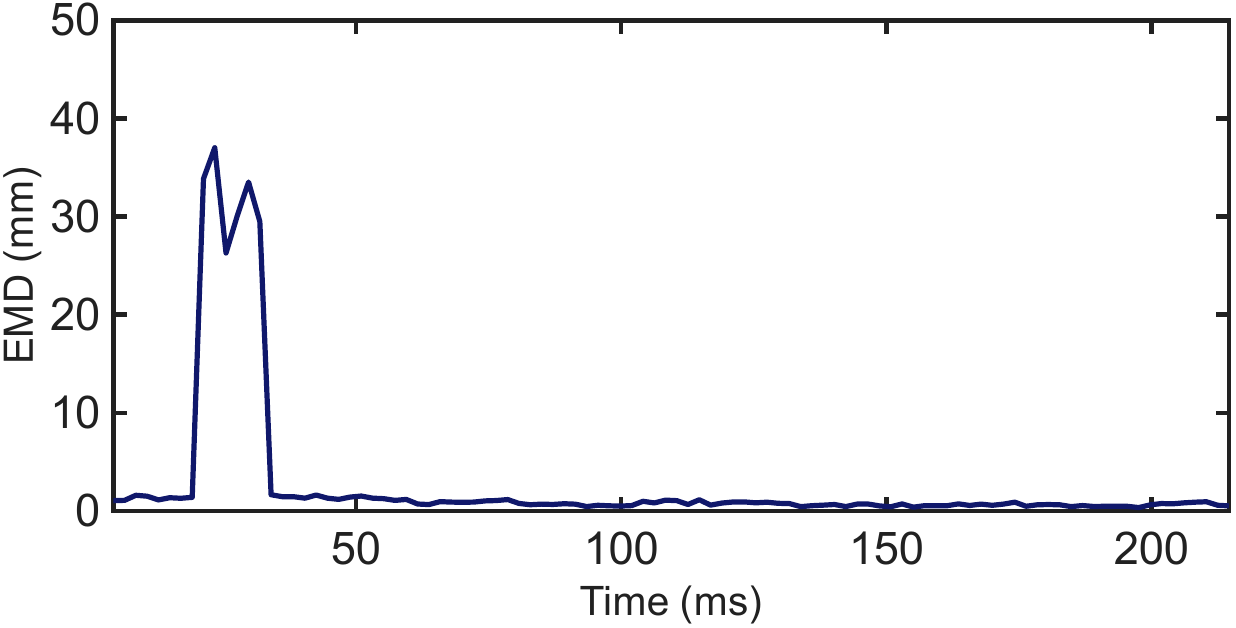}
    \end{minipage}

    \begin{minipage}{0.02\linewidth}
        \rotatebox{90}{\scriptsize Ling.}
    \end{minipage}\begin{minipage}{0.245\linewidth}
        \centering
        \includegraphics[width=0.99\linewidth]{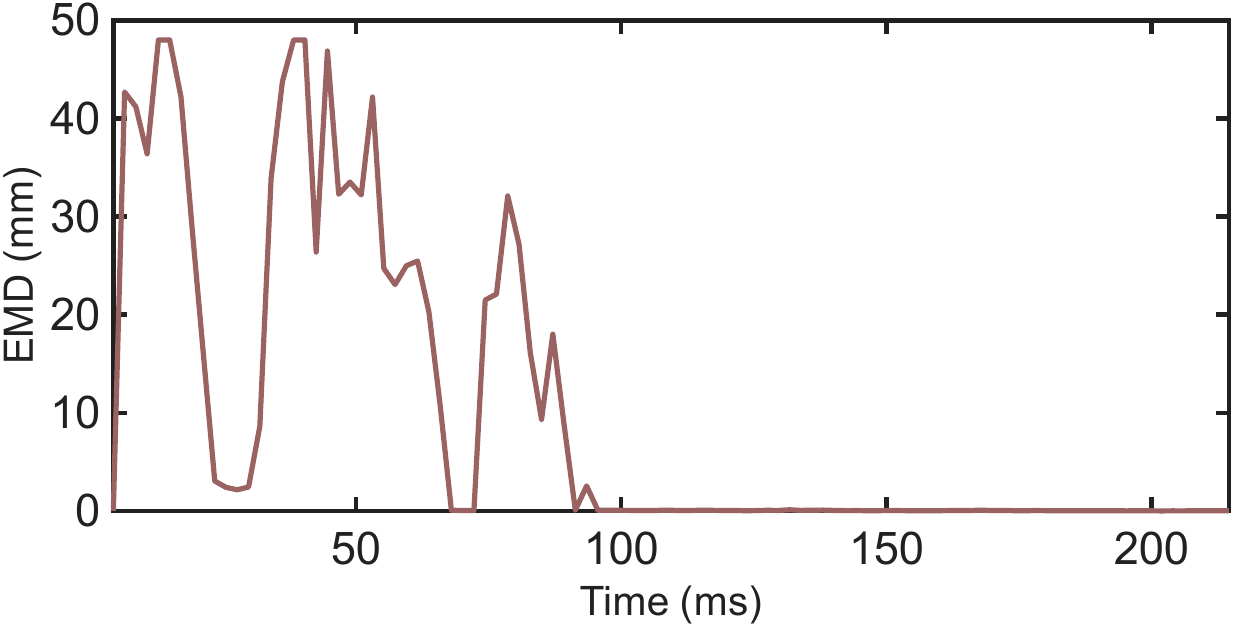}
    \end{minipage}\begin{minipage}{0.245\linewidth}
        \centering
        \includegraphics[width=0.99\linewidth]{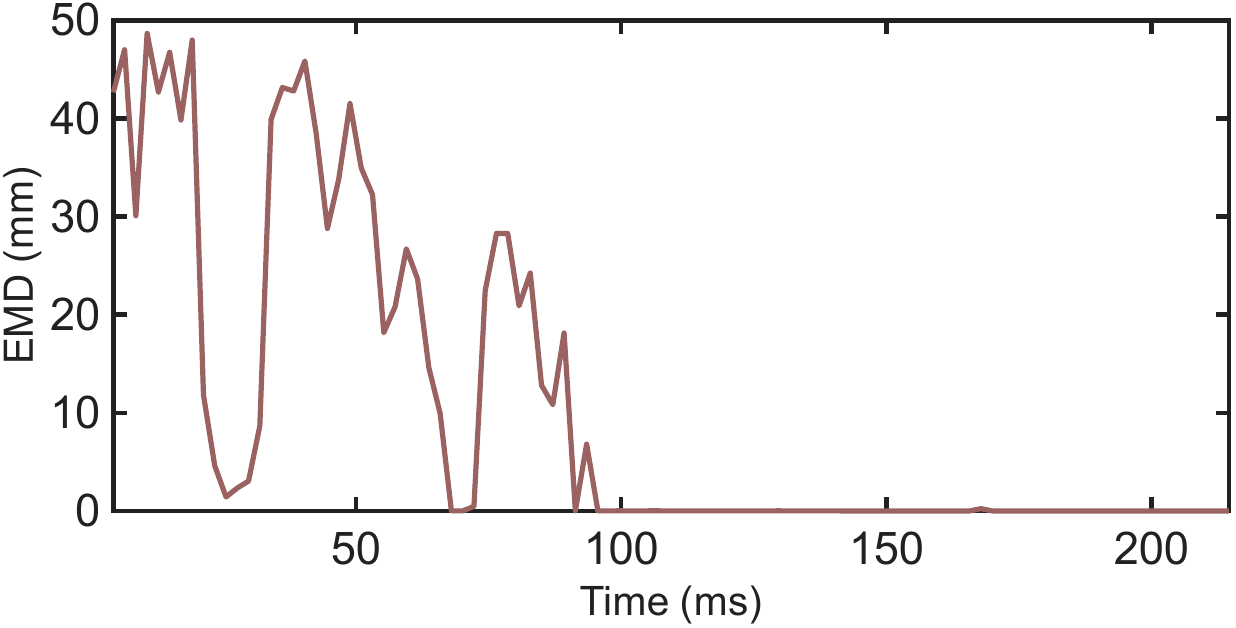}
    \end{minipage}\begin{minipage}{0.245\linewidth}
        \centering
        \includegraphics[width=0.99\linewidth]{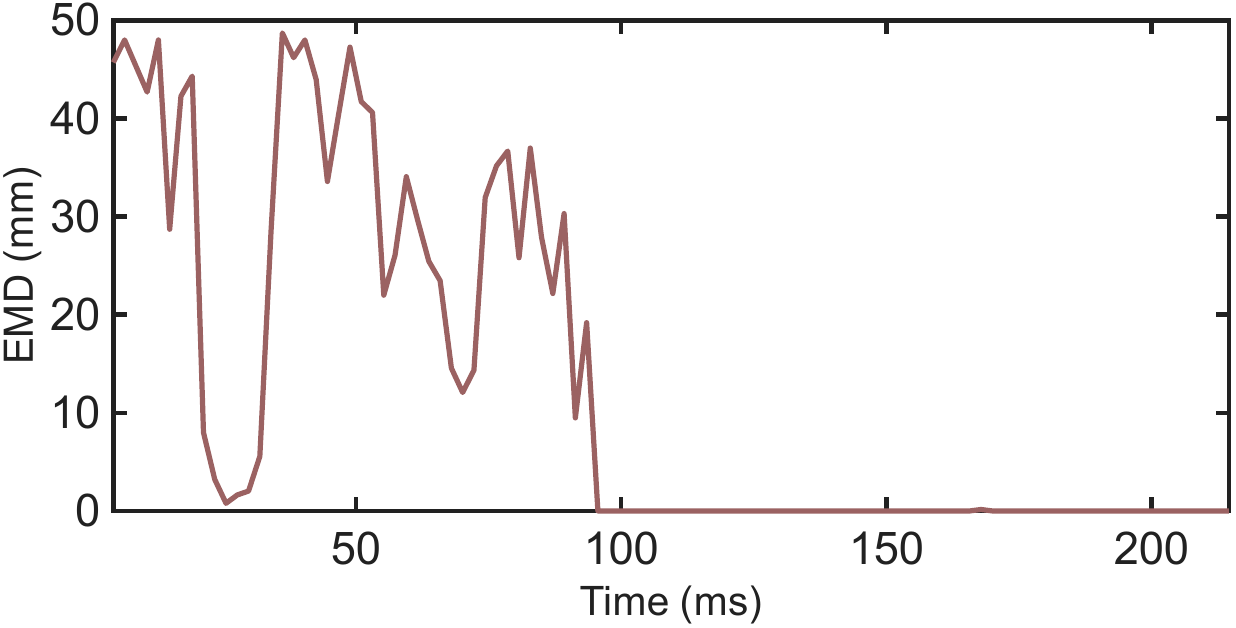}
    \end{minipage}\begin{minipage}{0.245\linewidth}
        \centering
        \includegraphics[width=0.99\linewidth]{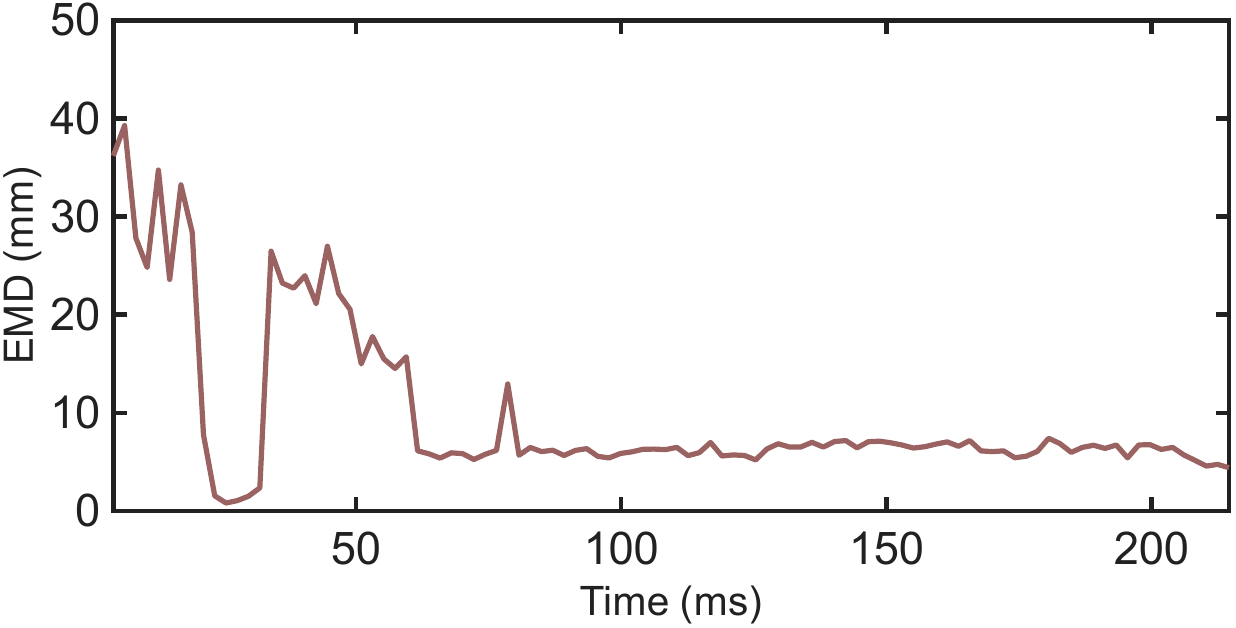}
    \end{minipage}

    \begin{minipage}{0.02\linewidth}
        \rotatebox{90}{\scriptsize S.P.}
    \end{minipage}\begin{minipage}{0.245\linewidth}
        \centering
        \includegraphics[width=0.99\linewidth]{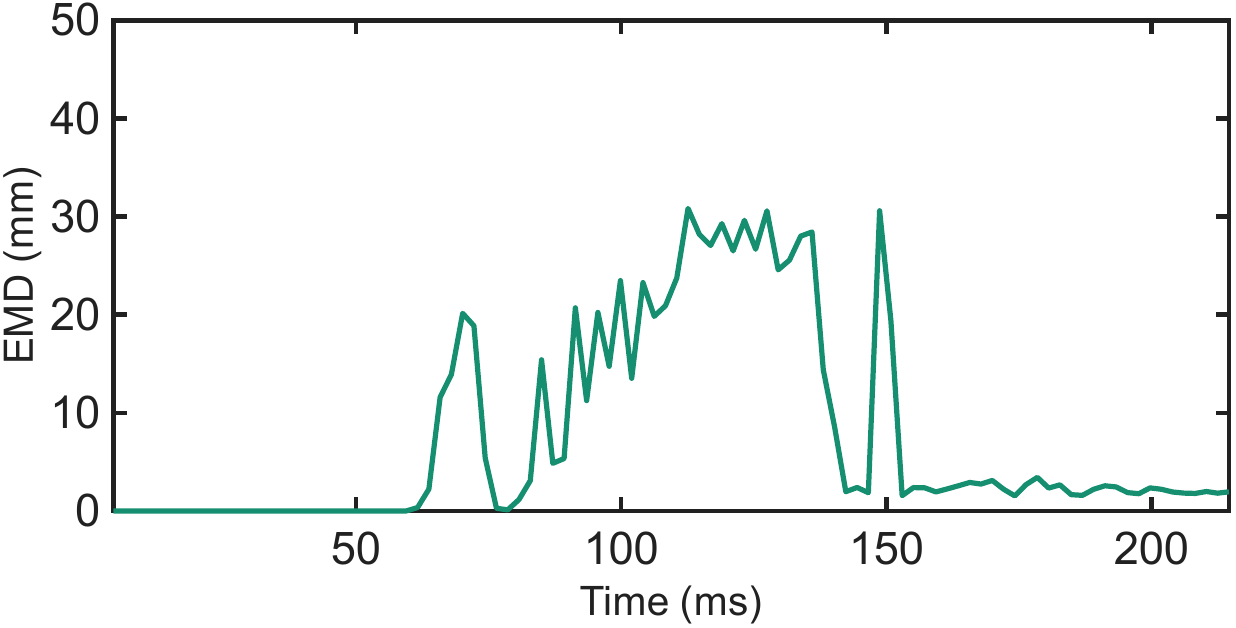}
    \end{minipage}\begin{minipage}{0.245\linewidth}
        \centering
        \includegraphics[width=0.99\linewidth]{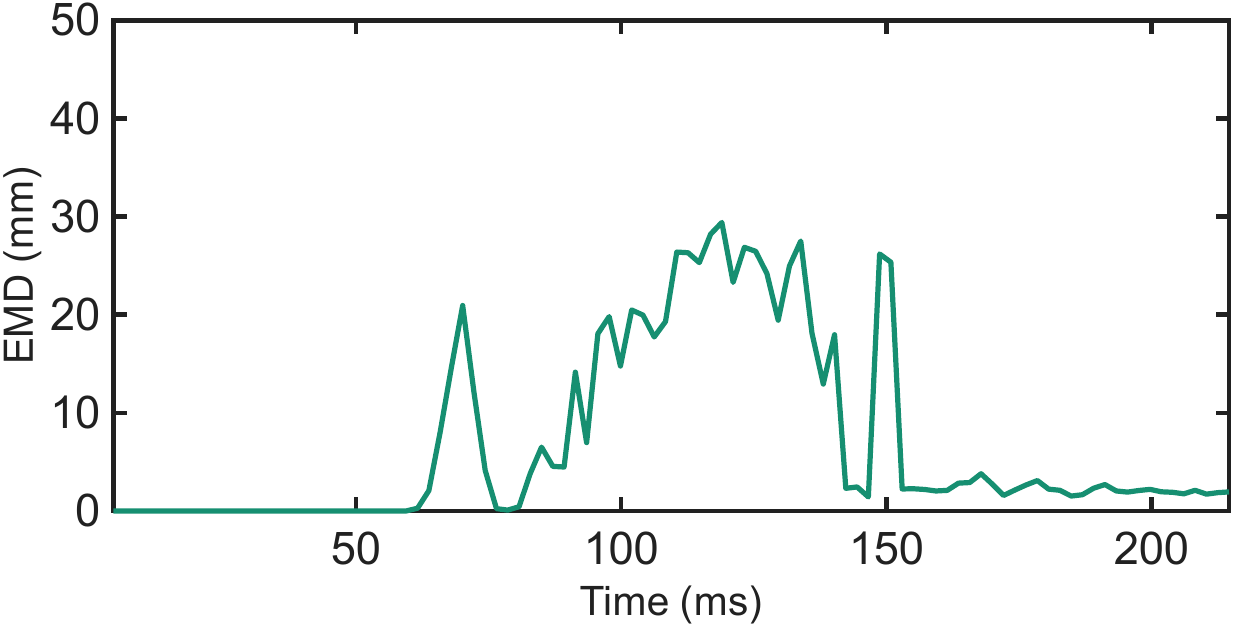}
    \end{minipage}\begin{minipage}{0.245\linewidth}
        \centering
        \includegraphics[width=0.99\linewidth]{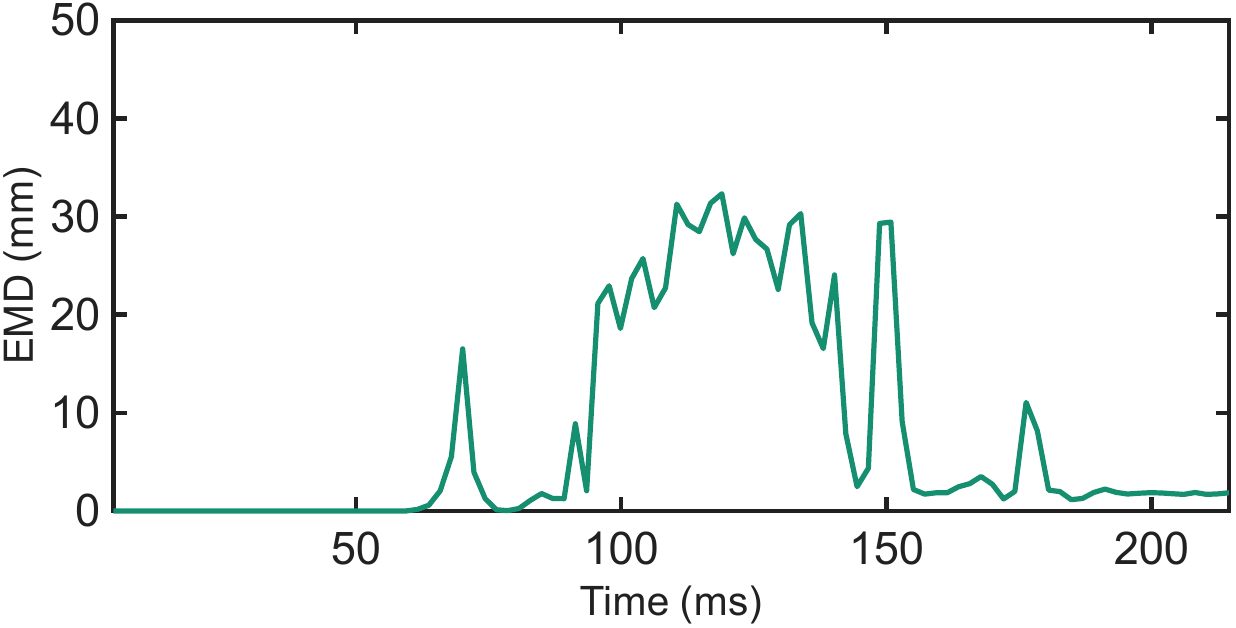}
    \end{minipage}\begin{minipage}{0.245\linewidth}
        \centering
        \includegraphics[width=0.99\linewidth]{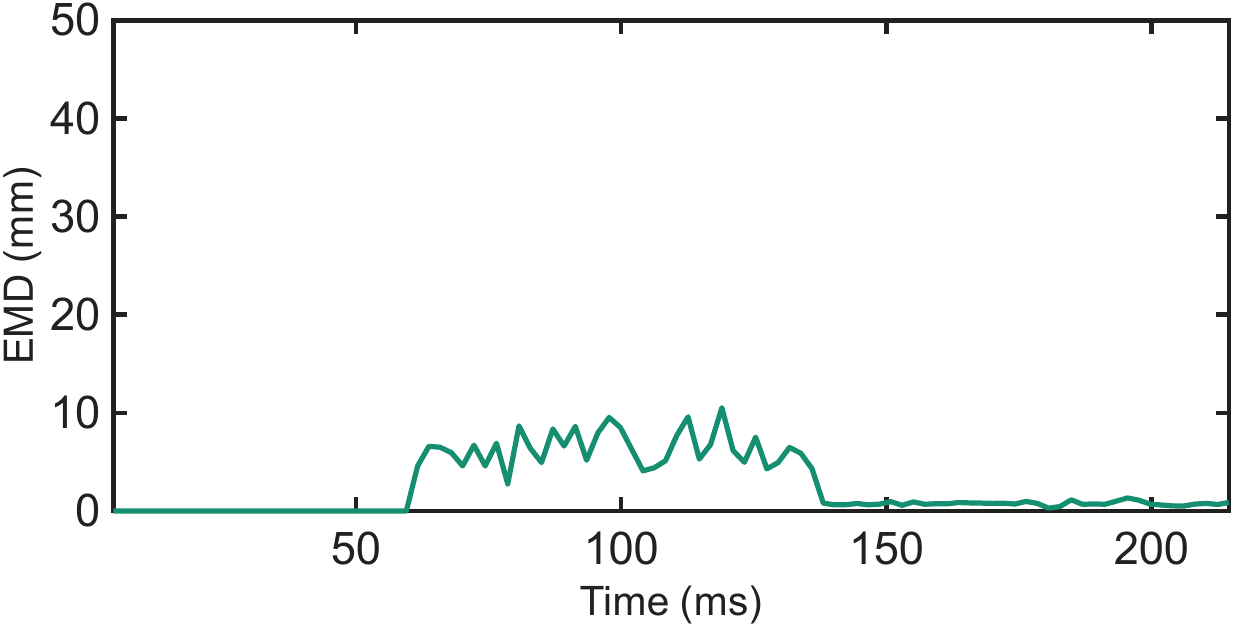}
    \end{minipage}

    \begin{minipage}{0.02\linewidth}
        \rotatebox{90}{\scriptsize C.M.F.}
    \end{minipage}\begin{minipage}{0.245\linewidth}
        \centering
        \includegraphics[width=0.99\linewidth]{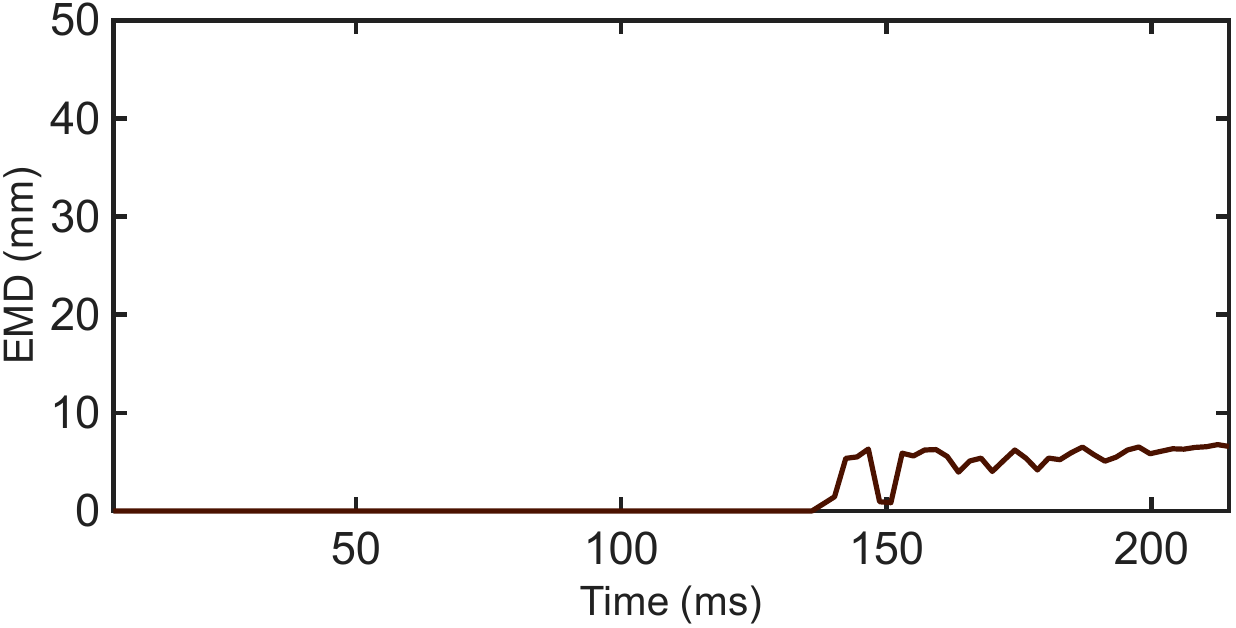}
    \end{minipage}\begin{minipage}{0.245\linewidth}
        \centering
        \includegraphics[width=0.99\linewidth]{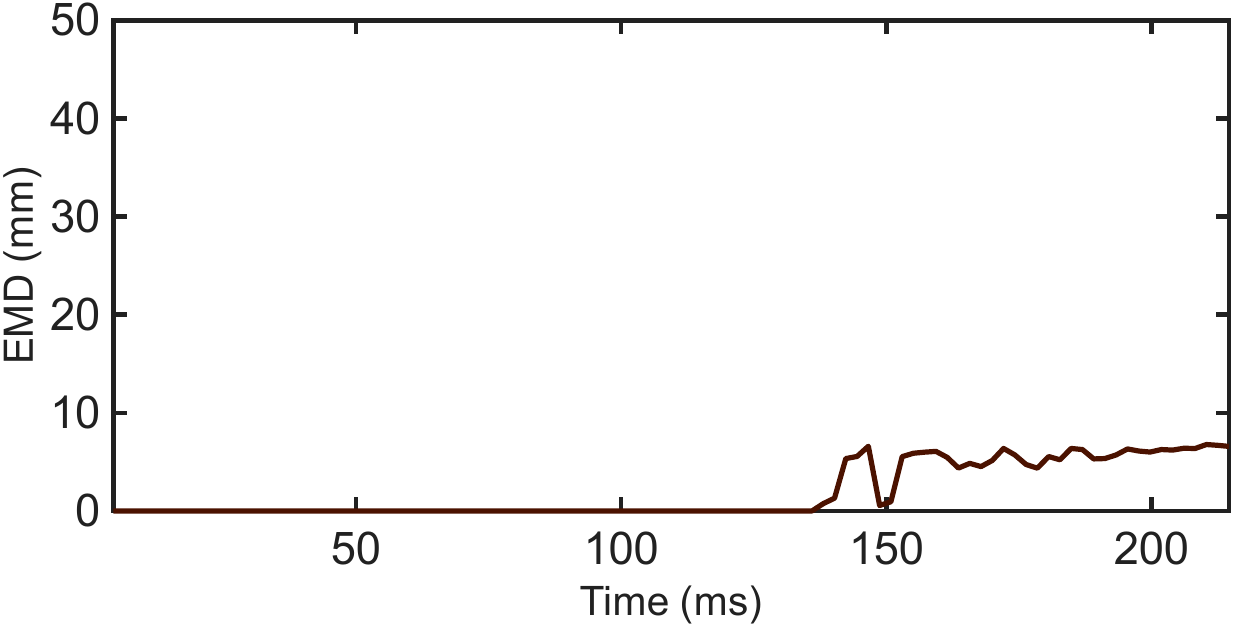}
    \end{minipage}\begin{minipage}{0.245\linewidth}
        \centering
        \includegraphics[width=0.99\linewidth]{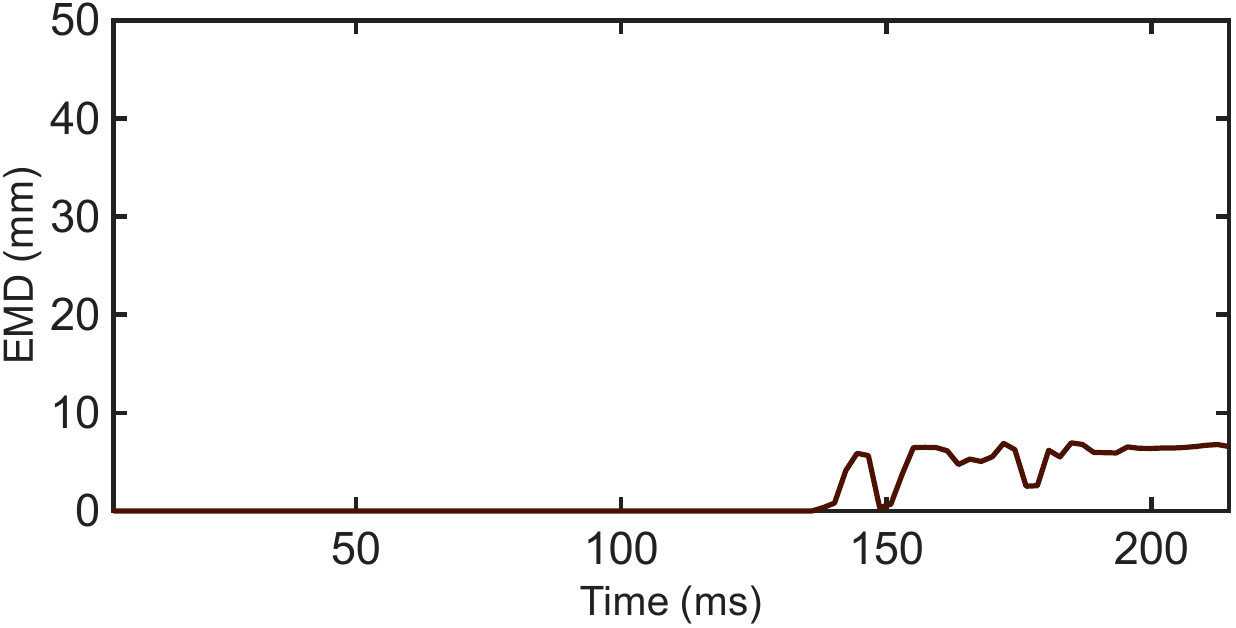}
    \end{minipage}\begin{minipage}{0.245\linewidth}
        \centering
        \includegraphics[width=0.99\linewidth]{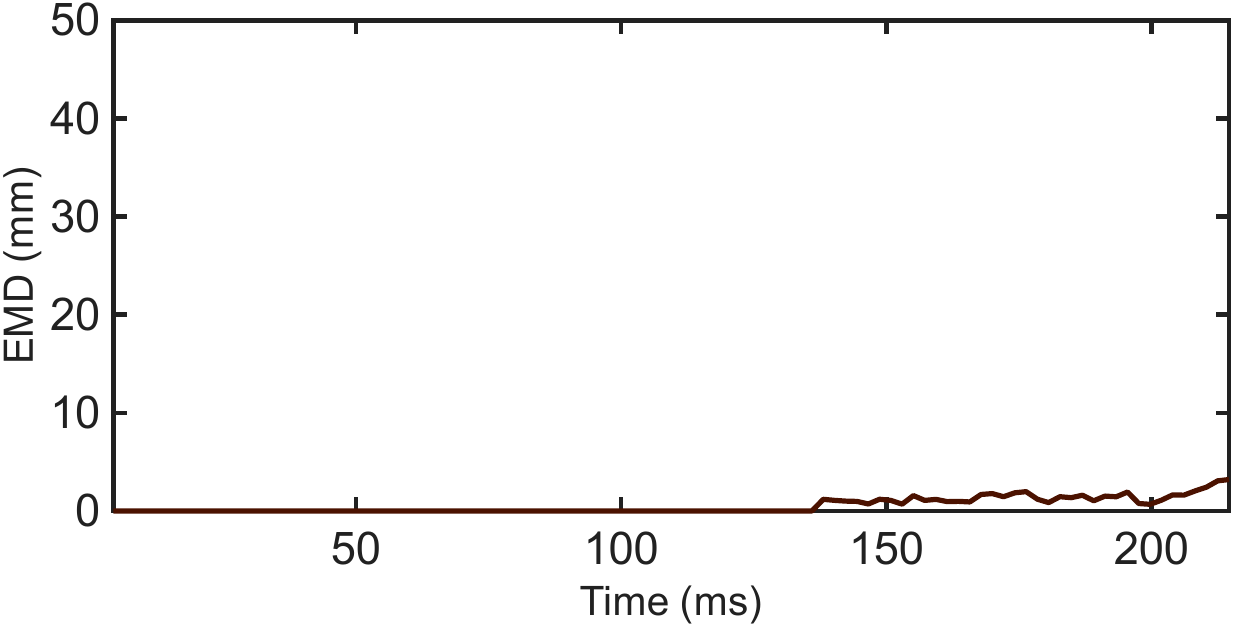}
    \end{minipage}

    \begin{minipage}{0.02\linewidth}
        \rotatebox{90}{\scriptsize Total}
    \end{minipage}\begin{minipage}{0.245\linewidth}
        \centering
        \includegraphics[width=0.99\linewidth]{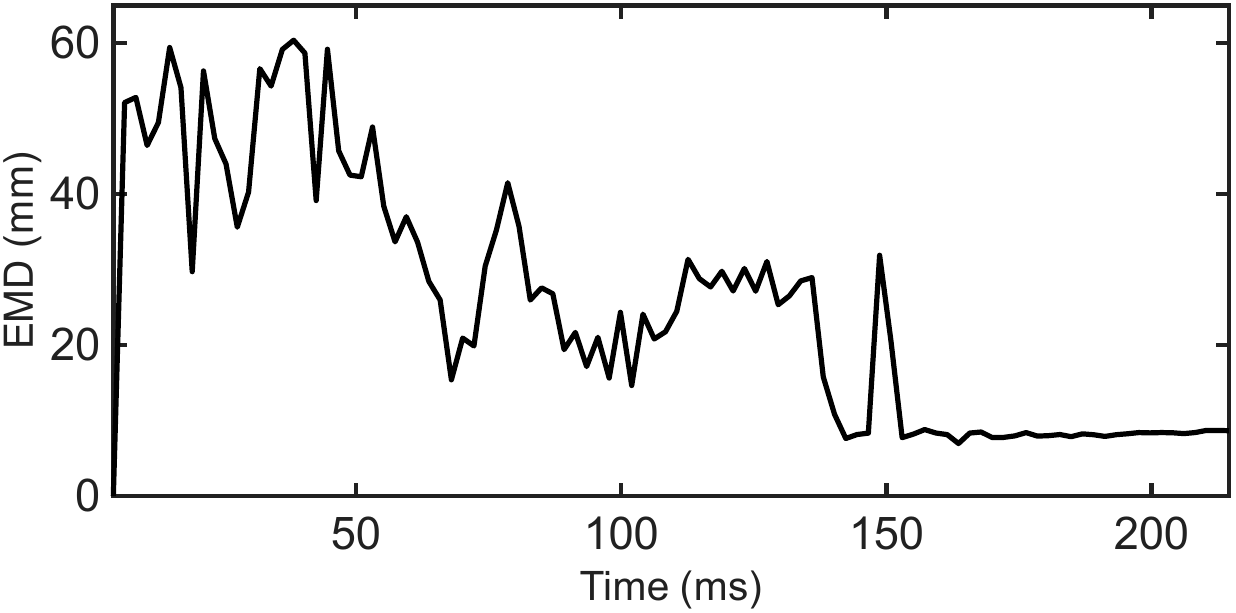}
    \end{minipage}\begin{minipage}{0.245\linewidth}
        \centering
        \includegraphics[width=0.99\linewidth]{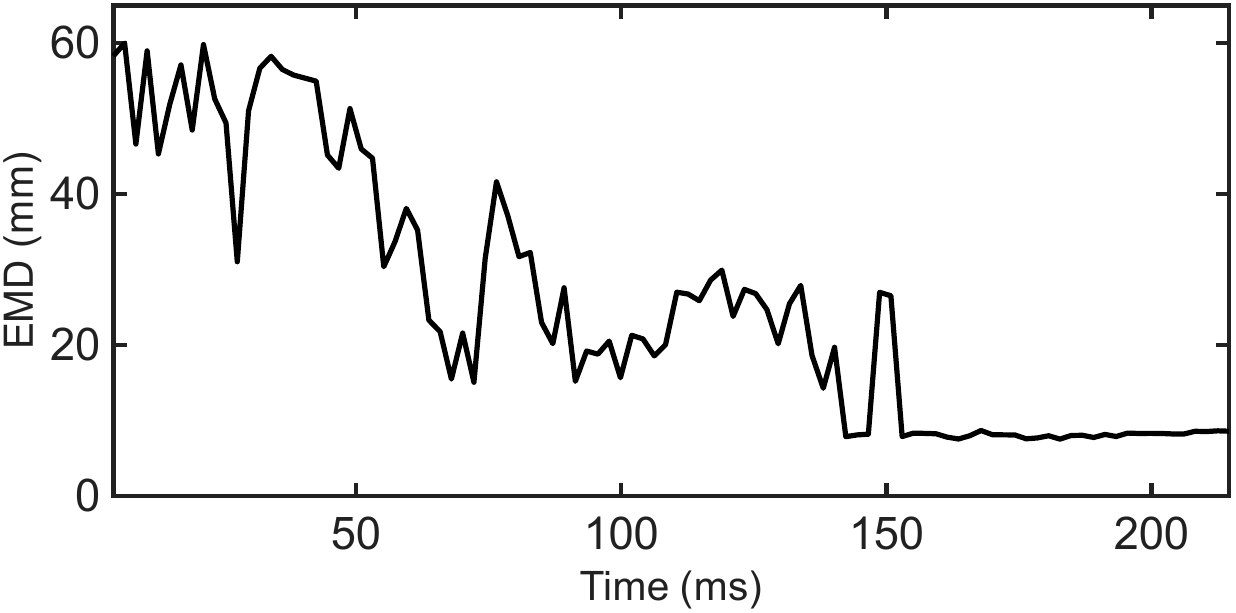}
    \end{minipage}\begin{minipage}{0.245\linewidth}
        \centering
        \includegraphics[width=0.99\linewidth]{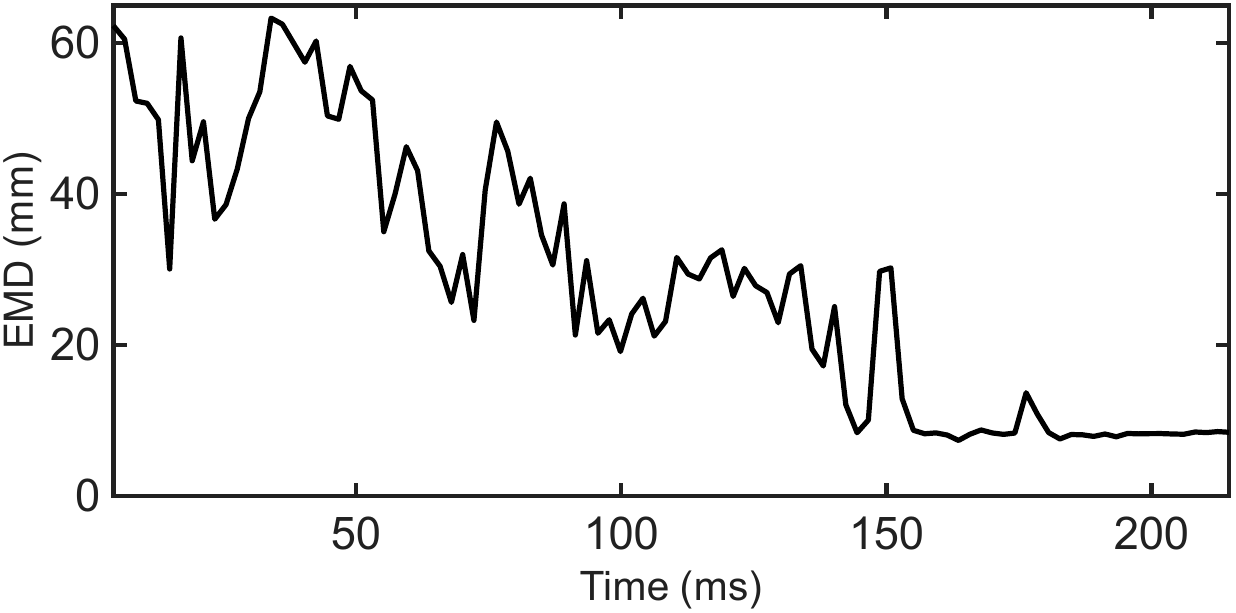}
    \end{minipage}\begin{minipage}{0.245\linewidth}
        \centering
        \includegraphics[width=0.99\linewidth]{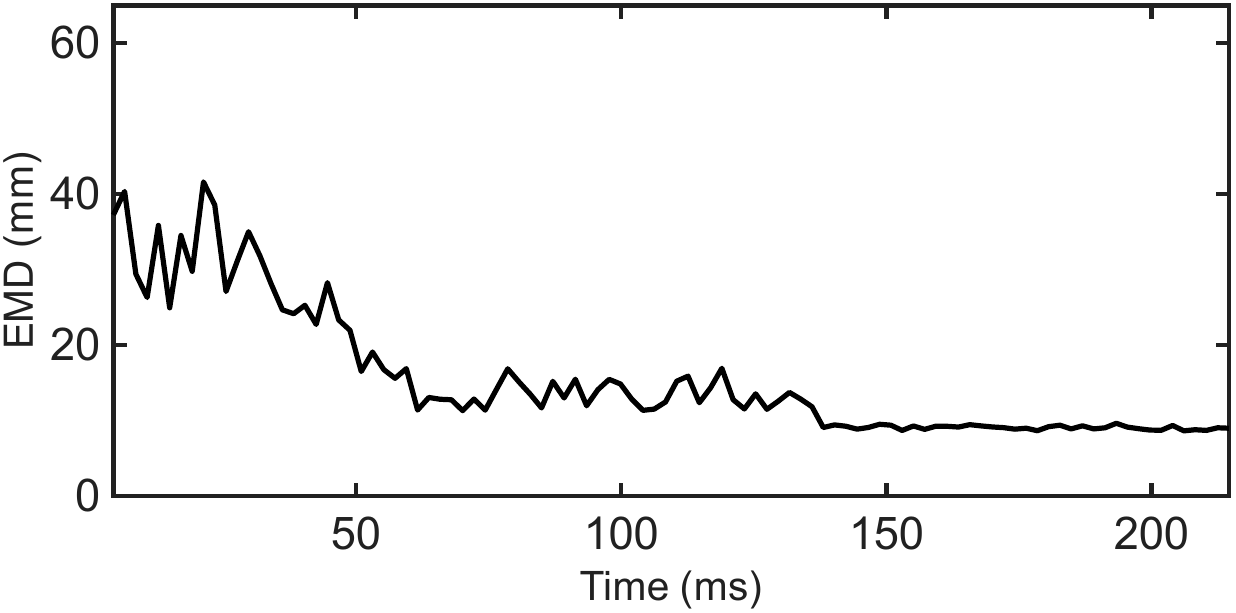}
    \end{minipage}

    \caption{Earth Mover's Distance (EMD) between the reconstructed and simulated ground-truth source distributions at 20 dB SNR. Columns correspond to CrSKF, SKF, KF, and LCMV; rows show parcel-wise EMD for the lateral occipital (L.O.), lingual (Ling.), superior parietal (S.P.), and caudal middle frontal (C.M.F.) regions, together with the total EMD. Lower values indicate closer spatial agreement.}
    \label{fig:EMD20dB}
\end{figure}

 Figures \ref{fig:EMD0dB}--\ref{fig:EMD20dB} show the resulting parcel-wise and total Earth Mover's Distance time courses. From the Earth Mover's Distance time course measured for the noiseless case (Figure \ref{fig:EMD0dB}), we can observe well-localized activity in the lateral occipital lobe for all the methods for most of the duration, as the EMD value is less than 1 \unit {\milli\metre}. However, after 25 \unit{\milli\second}. we observe an EMD spike for LCMV. Other well-localized activity is in the lingual and caudate middle frontal regions for LCMV. Total EMD for LCMV is nearly constant at 12 \unit {\milli\metre}, driven mainly by the summation of lingual and superior parietal EMDs. The total EMD of Kalman filters is mainly dictated by the superior parietal EMD spikes reaching up to 40 \unit{} 
With added Gaussian observation noise (25 and 20 \unit{\decibel}). We did not observe a drastic difference between the noisy-data experiments. However, compared with the noiseless case. We see a significant increase in EMD for the lingual and a suppression of EMD for the caudate middle frontal region for Kalman filter variants. The later occipital spike observed with LCMV without noise is present in noisy experiments for all the methods.  

\section{Discussion}

The present study combines a coupled Jansen--Rit neural mass model with a finite-element EEG forward model to create a controlled FVEP simulation framework for spatiotemporal source reconstruction. The causal source system branches V1 activity in the primary visual cortex into ventral and dorsal causal chains. The aim is to provide a structured and sufficiently realistic test case to evaluate, for example, inversion methods for spatiotemporal reconstruction or to recover the causal relations of individual dipolar sources. 

The realism of the generated data was tested by comparing the obtained PO7 electrode potential with that described in the literature \cite{odom2010TrueVEPPO7,AfrasiabiS2012TrueFVEPPO7}. Furthermore, we tested three variants of Kalman filtering (KF) and linearly constrained minimum variance (LCMV) for their ability to track activity across different parts of the brain atlas and localize in a focal manner. 

Similarly to the study by Hagler et al. \cite{Hagler2009}, we assumed the sources to be dipolar. We assume activity in the primary visual cortex (V1) branches into the ventral and dorsal parts of V2 and V3 areas. While the dorsal areas are associated with spatial information, motion-related processing, and visual information used to guide actions, the ventral areas are considered to be associated with visual form, object features, and shape processing; they should not be treated as completely separate pathways or as exclusively motion-processing versus form-processing regions, as both pathways process the visual information from the same stimulus \cite{Kravitz2011,Kravitz2013}. LIP and FEF are caused by salient visual stimuli such as a sudden flash as well as saccadic eye movements \cite{Thompson1999,Thompson2005}.

The reconstruction results reveal clear differences between the standardized Kalman approaches (SKF and Cr\-SKF), conventional KF, and LCMV. SKF and CrSKF recover the principal activity sequence more consistently in the parcel-wise time courses. This is consistent with the motivation of SKF as a spatiotemporal standardization of Kalman filtering \cite{Lahtinen2024SKF,Dilshanie2026SKF}. CrSKF produces the smoothest, most noise-robust trajectories because of its augmented state model and more flexible measurement-noise modeling. This smoothness helps suppress abrupt noise-driven fluctuations, but it can also reduce the apparent temporal separation of nearby sources; therefore, interpret it alongside the spatial maps and lag analysis rather than as an advantage in isolation.

Conventional KF tends to merge activity from neighboring parcels and to delay some reconstructed peaks. This behavior is consistent with the dependence of Kalman-based EEG inversion on the chosen evolution model \cite{GalkaAndreas2004KalmanEEG}. The random-walk model used here imposes temporal continuity but does not encode the known progression of the simulated visual network, so neighboring source contributions can remain strongly coupled in the reconstruction.

The examined spatial filter, LCMV, which is widely used for neuroelectromagnetic source analysis and connectivity studies \cite{VeenB1992Lois,HadjipapasAvgis2005LCMV,HuiHuaBrian2010LCMV,CoelliStefania2023LCMV}, can localize most of the sources accurately, but in the present simulation, its parcel-wise time courses are more strongly associated across regions and change more markedly as noise is added. This may be because the LCMV spatial filter is determined from the sensor covariance rather than from an explicit dynamical state model. 

The investigation of the lag between activity estimated at different parcels relevant to the modeled potentials reveals that KF variants also have difficulty separating lateral-occipital activity from lingual activity. Standardized Kal\-man filters suffer from an inability to distinguish lateral-occipital, lingual, and superior parietal activity from each other when noise is introduced to the observations. For a non-standard\-ized Kal\-man filter, the separation is weak in the aforementioned cases, but still observable. Activity in the caudal middle frontal region is separable from other source time evolutions. Overall, to accurately track simultaneous activity components, the observation data need to be as noiseless as possible, e.g., with multiple montages and montage-averaging.  

EMD complements the parcel-wise time courses by qu\-antifying spatial distributional mismatch. Its time dependence shows that inverse accuracy varies across the FVEP sequence, especially when multiple components overlap. This reinforces the value of evaluating temporal and spatial properties jointly rather than relying on a single parcel or peak-based criterion.

Overall, the previous studies of multifocal cortical sour\-ces demonstrate that such a source estimation is generally hard \cite{BelardinelliPaolo2012SMultiFocal,Bradley2016MultifocHard}, and ideal results cannot be expected even in two-source scenarios. Our results show that appropriate Kalman filtering can still yield valuable information. A possible benefit of a Kalman filter over other spatiotemporal methods \cite{OuHamalainen2009,LiuKe2016spatiotempBESTIES,LiuKe2019BESTIES2} is its recursive nature, which uses only past information to estimate the current state. 

Several limitations should be noted. First, the same lead field is used for forward simulation and inversion, so the analysis excludes errors caused by anatomical segmentation, tissue conductivity uncertainty, electrode-position errors, and other forms of forward-model mismatch. Second, the visual network is deliberately simplified: its source locations, delays, and relative gains are prescribed modeling choices rather than subject-specific physiological estimates. Second, the displayed noisy results are representative realizations rather than ensemble averages. Finally, the Jansen--Rit model is a reduced population model \cite{JansenRit1995}, so agreement with salient FVEP waveform features does not imply that all microscopic mechanisms of the visual response are reproduced.

Within these limitations, the framework extends previous studies of recovering complex spatiotemporal activity containing temporally overlapping components for SKF and CrSKF \cite{Lahtinen2024SKF,Lahtinen2026CRSKF} and LCMV \cite{HuiHuaBrian2010LCMV,PellegriniFranziska2023LCMV,CoelliStefania2023LCMV}.

A natural next step is to evaluate the synthetic model and inversion methods against experimental FVEP measurements. Once the data generation is evaluated against experimental data for its realism, it can serve as a benchmark for future spatiotemporal inversion method development.

\section{Conclusions}

A coupled Jansen--Rit neural mass model combined with a realistic finite-element EEG forward model provides a controlled means to synthesize FVEP data for evaluating spatiotemporal source reconstruction. In the present comparison, SKF and CrSKF preserve the principal temporal progression of the simulated visual response more consistently than conventional KF and LCMV, while CrSKF yields the smoothest temporal estimates. The analysis also shows that parcel-wise leakage, source overlap, and limited temporal separation can complicate interpretation when observation noise contaminates the data enough. The framework offers a basis for further evaluation of dynamic EEG inverse methods under progressively more realistic noise and forward-model mismatch.

\section*{Acknowledgements}
This study was supported by the Research Council of Finland (RCF) through the Center of Excellence in Inverse Modelling and Imaging 2018--2025 (353089), the  Flagship of Advanced Mathematics for Sensing, Imaging and Modelling (FAME) (359185). Sampsa Pursiainen and Qiong Tao were also supported by the Exploratory Study for Radar Tomography of Dimorphos project --- the Asteroid Moon of 65803 Didymos (RCF; decision 359198). We sincerely thank Prof.\ Carsten H.\ Wolters and the whole SIM-NEURO group (University of Münster, Münster, Germany) for fruitful discussions. We are also grateful to the German Academic Exchange Service (DAAD) and the Research Council of Finland (RCF; decisions 354976 and 367453) for supporting our travels to Münster. 

\section*{Author Contributions}
{\bf Qiong Tao}:\ Conceptualization (lead), Data Curation (equal), Software (Jansen--Rit model, lead), Methodology (supporting), Validation (equal), Writing – Review \& Editing (equal). {\bf Joonas Lahtinen}:\ Data Curation (equal), Formal Analysis (lead), Investigation (lead), Software (inversion \& quantitative measures, lead), Methodology (equ\-al), Validation (equal), Visualization (lead), Writing – Original Draft Preparation (lead).
{\bf Sampsa Pursiainen}:\ Conceptualization (supporting), Data Curation (equal), Funding Acquisition (lead), Methodology (equal), Project Administration (Lead), Software (visualization on brain, lead), Supervision (lead), Writing – Review \& Editing (equal)

\section*{Declaration of competing interests}
The authors declare no conflict of interest for this paper.



  \bibliographystyle{elsarticle-num} 
 \bibliography{bibliography}





\end{document}